\magnification=\magstep1

\input lect.def
\voffset=-1.0 truecm
%plain TeX with amssym files
%version of 12.03.96
\baselineskip=12.5pt plus 1.7pt

%\magnification=\magstep1
%\input komp.def
%\def\longpoints{\leaders\hbox to 0.5em{\hss.\hss}\hfill \hskip0pt}

\pageno=-1
\centerline{\bigbf COMPLEX CURVES IN ALMOST-COMPLEX }

\medskip
\centerline{\bigbf MANIFOLDS AND MEROMORPHIC HULLS}

\bigskip
\centerline{Sergei IVASHKOVICH -- Vsevolod SHEVCHISHIN}

\bigskip\bigskip\noindent{\bf Preface}

\bigskip\noindent
{\bigbf Chapter I. Local Properties of Complex Curves.}
\bigskip\noindent
\line{\bf Lecture 1. \bf Complex Curves in Almost-Complex Manifolds. \rm
\longpoints pp.\ 1--12}

\bgroup \baselineskip=11.5pt\parindent=0pt

\smallskip
\centerline{\vbox{\hsize=4.1truein  \noindent 
1.1.\ Almost-Complex Manifolds, Hermitian Metrics, Associated
(1,1)-Forms. 1.2. Existence of Calibrating and Tame Structures. 
1.3.\ Almost-Complex Submanifold,\break Complex Curves, Energy and Area. 
1.4.\ Symplectic Surfaces. 
1.5.\ Adjunction Formula for Immersed Symplectic Surfaces.
}}

\medskip\noindent
\line{\bf Appendix I. Chern Class and Riemann-Roch Formula.\rm
\longpoints pp.\ 13--14}

\smallskip \centerline{
\vbox{\hsize=4.1truein \noindent 
AI.1.\ First Chern Class of a Complex Bundle. 
AI.2.\ Rie\-mann-Roch Formula and Index of $\dbar $-Type Operators.
}}

\bigskip\noindent
\line{\bf Lecture 2. Local Existence of Curves.
\rm \longpoints pp.\ 15--21}

\smallskip
\centerline{
\vbox{\hsize=4.1truein  \noindent
2.1.\ Sobolev Imbeddings, Cauchy-Green Operators, Cal\-deron-Zygmund
    Inequality. 
2.2.\ Local Existence of Complex Curves in Nonintegrable Structures. 
2.3.\ Generalized Calderon-Zygmund Inequality. 
2.4.\ First A priori Estimate. 
2.5.\ Convergency outside of Finite Set of Points.
}}

\bigskip\noindent
\line{\bf Lecture 3. Positivity of Intersections of 
Complex Curves. \rm \longpoints pp.\ 22--37}

\smallskip\centerline{
\vbox{\hsize=4.1truein \noindent
3.1.\ Unique Continuation Lemma. 
3.2.\ Inversion of $\dbar_J+R$ and $L^{k,p}$-Topologies. 
3.3.\ Perturbing Cusps of Complex Curves. 
3.4.\ Primitivity. 
3.5.\ Positivity of Intersections.
\smallskip
}}

\noindent
\line{\bf Appendix II. The Bennequin Index and Genus Formula.
\rm \longpoints pp.\ 38--43}
\smallskip\rm
\centerline{
\vbox{\hsize=4.1truein \tenrm \noindent
AII.1.\ Bennequin Index of a Cusp.
AII.2.\ Proof of Adjunction Formula.
AII.3.\ A Matching Structures Lemma.
}}

\bigskip\noindent
{\bigbf Chapter II. Compactness Theorem.}

\medskip\noindent
\line{\bf Lecture 4. Space of Stable Curves and Maps. \rm 
\longpoints pp.\ 46--66}

\smallskip\rm\centerline{
\vbox{\hsize=4.1truein \tenrm\noindent
4.1.\ Stable Maps and Gromov Topology. 
4.2.\ Fenchel-Nielsen Coordinates on the Space of Nodal Curves. 
4.3.\ Complex Structure on the Space $\ttt_\Gamma$. 
4.4.\ Invariant Description of the Complex Structure on $\ttt_\Gamma$.
4.5.\ Example: Degeneration to a Half-cubic Parabola on the Language 
of Stable Curves.
}}

\bigskip\noindent
\line{\bf Lecture 5.  Gromov Compactness Theorem. \rm 
\longpoints pp.\ 67--84}

\nobreak
\smallskip\centerline{
\vbox{\hsize=4.1truein \noindent
5.1.\ Second A priori Estimate. 
5.2.\ Removal of Singularities. 
5.3.\ Compactness for the Curves with Free Boundary.
}}

\medskip\noindent
\line{\bf Appendix III. Compactness with Boundary Conditions. \rm 
\longpoints pp.\ 85--110}

\smallskip\rm\centerline{
\vbox{\hsize=4.1truein \tenrm\noindent
A3.1.\ Curves with Bondary on Totally Real Submanifolds. 
A3.2.\ A priori Estimates near a Totally Real \-Boundary.
A3.3.\ Long Strips and the Second A priori Estimate. 
A3.4.\ Gromov Compactness for Curves with Totally Real Boundary Conditions.
A3.5.\ Attaching an Analytic Disk to a Lagrangian Submanifold in $\cc^n$.
}}

\bigskip{\bigbf Chapter III. Global Properties and Moduli Spaces.}

\medskip{\bf Lecture 6. First Variation of $\bar\partial_J$-Equation. \rm 
\longpoints pp.\ 111--122}

\smallskip\centerline{
\vbox{\hsize=4.1truein
6.1.\ Symmetric Connections. 
6.2.\ Definition of $D_{u,J}$-Ope\-ra\-tor.
6.3.\ $\dbar $-type Operators. 
6.4.\ Holomorphic Structure on the Induced Bundle. 
}}

\bigskip
{\bf Lecture 7. Fredholm Properties of Gromov Operator. \rm 
\longpoints pp.\ 123--129}

\smallskip\centerline{
\vbox{\hsize=4.1truein \noindent
7.1.\ Generalized Normal Bundle. 
7.2.\  Surjectivity of \-$D_{u,J}$ - operator.
7.3.\ Tangent Space to the Moduli Space. 
7.4.\ Re\-para\-mete\-ri\-zations.
}}

\bigskip
{\bf Lecture 8. \bf Transversality. \rm 
\longpoints pp.\ 130--144}

\smallskip\centerline{
\vbox{\hsize=4.1truein 
8.1.\ Moduli Space of Nonparameterized Complex Cur\-ves. 
8.2.\ Transversal maps. 
8.3.\ Components of the Moduli Space. 
8.4.\ Moduli of Parameterized Curves. 
8.5.\ Gromov Non\-squee\-zing Theorem.
8.6.\ Exeptional Spheres in Symplectic $4$-Manifolds.
%8.7.\ Example of  Sikorav.
}}

\bigskip{\bigbf Chapter IV. Envelopes of Meromorphy.}

\bigskip{\bf Lecture 9. Deformation of Noncompact Curves. \rm 
\longpoints pp. 147--159}

\smallskip \centerline{
\vbox{\hsize=4.1truein
9.1.\ Banach Analytic Sets. 
9.2.\ Solution of a Cousin-type Problem.
9.3.\ Case of a Stein Curve. 
9.4.\ Degeneration to a Node. 
9.5.\ Banach Analytic Structure of the Stable Neighborhoods. 
9.6.\ Drawing Families of Curves.
}}

\bigskip\noindent
{\bf Lecture 10. Envelopes of Meromorphy of Two-Spheres. \rm 
\longpoints pp.\ 160--168}

\smallskip\centerline{
\vbox{\hsize=4.1truein
10.1.\ Continuity Principle relative to the K\"ahler Spaces. 
10.2.\ Proof of the Continuity Principle. 
10.3.\ Construction of Envelopes -- I.
10.4.\ Construction of Envelopes -- II. 
10.5.\ Examples.
10.6.\ Rigidity of Symplectic Imbeddings.
}}

\bigskip\noindent
{\bf Appendix IV. Complex Points and Stein Neighborhods. \rm 
\longpoints pp.\ 169--172}

\smallskip\centerline{
\vbox{\hsize=4.1truein
A4.1.\ Complex Points of Real Surfaces in Complex Surfaces.
A4.2.\ Cancellation of Complex Points.
A4.3.\ Nei\-ghborhoods of Real Surfaces.
}}

\newpage
\bigskip\noindent
{\bf Appendix V. Seiberg-Witten Invariants and Envelopes. \rm 
\longpoints pp.\ 173--181}

\smallskip\centerline{
\vbox{\hsize=4.1truein
A5.1.\ The Genus Estimate.
A5.2.\ The use of Seiberg-Witten Invariant.
A5.3.\ The Genus Estimate on Stein Surfaces and Envelopes.
A5.4.\ Two-spheres in $\cc^2$.
A5.5.\ Attaching Complex Disks to Strictly Pseudonconvex Domains.
}}

\bigskip\noindent
{\bf References. \rm 
\longpoints pp.\ 182--186}

\egroup

\bigskip
\newpage
\pageno=0

%\footline=\toks11
%\nopagenumbers\footline {\hfil }
\centerline{\bf P r e f a c e}

\bigskip\bigskip
This are the notes of a course, given by the first author for the
Graduiertenkollegs (= graduate students) at the Ruhr-University of Bochum, in
December, 1997. Previously, shorter courses were given in G\"oteborg, Warsaw
and Bonn.

\bigskip
These lectures pursued two main tasks:

\bigskip\noindent
{\bf First} --- to give a systematic and self-contained introduction to the
Gromov theory of so-called pseudoholomorphic curves (a term which will be
completely eliminated in these notes). This is handled in Chapters I,II,III.

\medskip\noindent
{\bf Second} --- to explain our joint results on envelopes of meromorphy of
real two-spheres in complex two-dimensional manifolds, which we obtained using 
Gromov theory. We do this in Chapter IV.

\bigskip
It was impossible, of course, to omit the original motivation of Gromov in the
development of the subject:

\medskip
\noindent
 - in Appendix III, \S A3.5, we explain the Gromov's idea how 
to attach an analytic disk to a Lagrangian submanifold in $\cc^n$;

\noindent
 - in \S 8.5 we give the proof of his ``non-squeezing theorem"; 

\noindent
 - in \S 8.6 we study exeptional spheres in symplectic $4$-manifolds.  

\smallskip\smallskip In Appendix V, following S. Nemirovski, we discuss one 
more approach 
to the problem of constructig envelopes of meromorphy of spheres. This 
approach uses Seiberg-Witten theory. As an example if his results 

\smallskip\noindent
 - we prove that a complex disk cannot be attached from the 
 outside to the strictly pseudoconvex domain in $\cc^2$, which is 
 diffeomorphic to the  ball.

\bigskip
The authors would like to give their thanks to Alan Huckleberry, the head of
Graduiertenkolleg program at the Ruhr-University and to SFB-236 for the
constant support of our joint research. We would like also to give our 
thanks to Barbara Huckleberry for the uncountable number of grammatical 
corrections.

\bigskip
\bigskip
\bigskip
\hfill{Sergei  and  Seva, November 1999.}

\smallskip
\hfill{Bochum --- Lille --- Lviv.}
\vfil

%%%%%%%%%%%%%%%%%%%%%%%%%%%%%%%%%%%%%%%%%%%%%%%%%%%%%%%%%%%%%%
\newpage

\noindent
{\bigbf Chapter I. Local Properties of Complex Curves.}

\bigskip
In this chapter we shall concentrate on local properties of complex
curves in almost-complex manifolds.The first lecture is very basic and mainly
defining. Its aim is to introduce a couple of notions and recall a few standard
facts. In the second one we prove, following Vekua and Sikorav, the existence
of $J$-complex curves (locally through any point and any direction) for any
almost-complex structure $J$ of smoothness $C^{1,\alpha}$, $\alpha >0$. One
more principal result here is the so-called {\it first a priori estimate} in
\S 2.4. The third lecture, as well as {\sl Appendix 2}, is devoted to the 
deeper study of complex curves in real 4-dimensional manifolds. We shall
prove the positivity of intersections and the Adjunction Formula, showing
that local properties of complex curves in nonintegrable almost-complex 
structures are similar to that in integrable ones.

\bigskip\noindent
{\bigbf Lecture 1}

\smallskip\noindent
{\bigbf Complex Curves in Almost-Complex Manifolds.}

\medskip\noindent
{\bigsl 1.1. Almost-Complex Manifolds, Hermitian Metrics, Associated
(1,1)-forms.}

\smallskip
Let $X$ be a real manifold, and $TX$ its tangent bundle.

\state Definition 1.1.1. {\it A continuous endomorphism $J:TX\to TX$ such
that $J^2=-\id $ is called an almost-complex structure on $X$. A pair $(X,J)$
is called an almost-complex manifold.
}

\state Exercise. Prove that locally one can always find vector fields
$e_1,e_2,...,e_{2n-1},e_{2n}$ such that $Je_{2k-1}=e_{2k}, k=1,...,n$.
In particular, $X$ is even-dimensional.

\smallskip\noindent
\state Examples.  $1^0.$ $(\rr^{2n},J\st )$. Denote by $e_1,...,e_{2n}$
the standard basis in $\rr^{2n}$. Define a {\it  standard} complex
structure $J\st $ in $\rr^{2n}$ by $J\st (e_{2k-1})=e_{2k}$.

If $(V,J)$ is an $\rr $-linear space  endowed with an endomorphism $J$,
such that $J^2=-\id $, then a structure of a $\cc $-linear space on $V$
is defined by $z\cdot v=(x+iy)\cdot v:=x\cdot v +y\cdot Jv$. In the case
$(\rr^{2n},J\st )$ this gives $\cc^n$.

\smallskip
\noindent $2^0.$ Let $x_0\in (X,J)$. By an affine change of coordinates in
some chart $U\ni x_0$ we can always suppose that $x_0=0$ and $J(0)=J\st $
in basis ${\d \over \d x_1},...,{\d \over \d x_{2n}}$. One can now put
$y_j=x_{2j}$ in order to have $J(0) {\d \over \d x_j}={\d \over \d y_j}$.

\smallskip\noindent 
\sl $3^0.$ Complex manifolds. \rm By definition, a complex manifold is an 
almost-complex manifold $(X,J)$ such that in the neighborhood of 
any point $x_0\in X$ a
coordinate system $x_1,y_1,...,x_n,y_n$ can be found such that
$J {\d \over \d x_j}={\d \over \d y_j}$ everywhere in this neighborhood.

Equivalently there is in a neighborhood of each point a $C^1$-diffeomorphism
$\phi $ into $(\rr^{2n},J\st )$ such that its differential is
$(J,J\st )$-linear. The latter means that $d\phi \scirc J=J\st \scirc d\phi $
as mappings $TX\to \rr^{2n}$.

\smallskip\noindent
\state Definition 1.1.2. {\it Such an almost-complex structure is called
complex.
}

\smallskip\noindent
\sl $4^0.$ Riemann surfaces. \rm Riemann surfaces are oriented, connected
two-dimensional manifolds. Fix some Riemannian metric $g$ on such $S$.
For $v\in TS$ define $Jv=\{$ rotation of $v$ in positive direction onto
$90^\circ \} $. $J$ is an almost-complex structure on $S$. If $g_1=\lambda g$
is another, conformally equivalent to $g$ metric, then this construction
leads to the same almost-complex structure.

Vice-versa, given $J$ on $S$ , denote by $J^{\cc }$ the $\cc $-linear
extension of $J$ onto the complexified tangent bundle $T^{\cc }S:=TS\otimes
\cc $. Here $J^{\cc }$ is defined by $J^{\cc }(v\otimes \lambda )=Jv\otimes
\lambda $ for $v\in TS$ and $\lambda \in \cc $. Take an
eigenvector $w$ of $J^{\cc }$ in $T^{\cc }S$ with an eigenvalue $-i$.
Write $w=u+iv$, with
$u,v\in TS$ (more accurately $w=u\otimes 1 + v\otimes i$), and put $g(u,v)=0,
g(u)=g(v)=1$. This defines a Riemannian metric, which corresponds to $J$.
The statement $J^{\cc }w=-iw=v-iu$ means that $Ju=v, Jv=-u$.

The freedom in choosing $w$ here permits us to define $g$ up to a real factor,
\ie, $J$ corresponds to a class of conformally equivalent metrics.

\state Exercise. Prove the last assertion. More precisely: if $\lambda w$ is
another eigenvector of $J^{\cc }$ with an eigenvalue $-i$ and a
metric $g_{\lambda }$ is constructed from $\lambda w$ as $g$ from $w$ above,
then $g=\vert \lambda \vert^2g_{\lambda }$.

\smallskip
We shall see later in {\sl Corollary 2.2.2} that any almost-complex 
structure on the Riemann surface is complex.

\smallskip
\state Definition 1.1.3. {\it $J$-Hermitian metric on $X$ is a $\rr $-bilinear
form $h:TX\times TX\to \cc $ such that

(a) $h(Ju,v)=\overline h(v,Ju)=ih(u,v)$;

(b) $h(u,u)>0$ for $u\not= 0$.
}

\smallskip
Decompose $h(u,v)=g_h(u,v) -i\omega_h(u,v)$. Then :

\state Exercise.
(a) $g_h$ is a Riemannian metric on $X$ and $\omega_h$ is an exterior 2-form.

(b) $g_h(u,v)=\omega_h(u,Jv)$ and thus $h(u,v)=\omega_h(u,Jv) -
i\omega_h(u,v)$. More over, one has $\omega_h (Ju,Jv)=\omega_h (u,v)$.

(c) $\omega_h(u,v)=-g_h(u,Jv)$.

\smallskip
\state Definition 1.1.4. The form {\it $\omega_h$ is called a (1,1)-form 
associated to $h$.}

\smallskip
Conversely,

\state Definition 1.1.5. {\it An exterior $2$-form $\omega $ on $X$ is called
$J$-calibrated if

(a) $\omega (Ju,Jv)=\omega (u,v)$;

(b) $\omega (u,Ju)>0$ for $u\not= 0$.}

\smallskip Thus,  $\omega $ clearly defines a $J$-Hermitian metric
$h^{\omega }(u,v)=\omega (u,Jv) - i\omega (u,v)$.

\smallskip A triple $(X,J,h)$, where $h$ is $J$-Hermitian is called
an almost-Hermitian manifold. We shall not use this term in these notes.

\state Example. To explain the coefficient $-i$ chosen in the decomposition
$h(u,v)=g(u,Jv) - i\omega (u,v)$ consider $(\rr^{2n}, \omega\st, J\st)$.
Here $\omega\st = \sum_{j=1}^ndx_j\wedge dy_j$
is a standard symplectic form in $\rr^{2n}$. The Hermitian metric here
is $h\st = \sum_{j=1}^ndz_j\otimes d\bar z_j = \sum_{j=1}^n (dx_j\otimes
dx_j + dy_j\otimes dy_j) - i\sum_{j=1}^n(dx_j\otimes dy_j - dy_j\otimes dx_j)
= g\st -i\omega\st $.

Thus $\omega\st = -\im h\st $ is a canonical expression for a K\"ahler
form of the Hermitian metric.

\state Definition 1.1.6. {\it We say that an exterior $2$-form $\omega $ on $X$
tames an almost-complex structure  $J$ if $\omega (v,Jv)>0$ for $v\not= 0$.
}

\smallskip
\state Remark. {\rm Here and later we will use the convention
that for forms $f$ and $h$ on the vector space $V$ the exterior
product is defined as $f\wedge g: = f\otimes g - g\otimes f$, \ie,
without coefficient ${1\over 2}$.
}

\medskip\noindent\bigsl 1.2. Existence of Calibrating and Tame Structures. \rm

\nobreak
\smallskip
Here we shall prove here the following simple but very important statement:

\state Proposition 1.2.1. {\it
Let $X$ be a real manifold and $\omega $ a nowhere degenerate
exterior $2$-form on $X$, \ie, $\omega^n\not=0$, where $2n=\dim_{\rr } X$.
On
the open subset $U\subset X$ a $\omega $-calibrating almost-complex structure
$J$ is given. Then for any relatively compact open $U_1\comp U$ there exists
a $\omega $-calibrating almost-complex structure $J_1$ on the whole $X$ such that
$J_1\mid_{U_1}=J\mid_{U_1}$.
}

\state Proof. Consider a Riemannian metric $g(u,v):=\omega (u,Jv)$ on
$U$. Find a Riemannian metric $g_1$ on $X$ such that $g_1\mid_{U_1}=g\mid_{
U_1}$. Since $\omega $ is not degenerate, there exists a (unique!)
isomorphism
$A_1:TX\to TX$ such that
$$
\omega (u,v) = g_1(A_1u,v)\eqno(1.2.1)
$$
for all $u,v\in TX$. Further, $\omega $ is antisymmetric, so
$$
g_1(A_1u,v) = \omega (u,v) = -\omega (v,u) = -g_1(A_1v,u) = -g_1(u,A_1v).
$$
Thus $A_1^*=-A_1$. So $A_1^*A_1=A_1A_1^*=-A_1^2$ is positively
definite and self-adjoint \wrt $g_1$. Let $Q_1:=\sqrt{-A_1^2}$ be a positive
square root of $-A_1^2$. Put
$$
J_1=A_1Q_1^{-1}.\eqno(1.2.2)
$$
$J_1$ is an almost-complex structure. Indeed
$$
J_1^2=A_1Q_1^{-1}A_1Q_1^{-1}=A_1^2(Q_1^{-1})^2=A_1^2(-A_1^{-2})=-\id .
$$
Let us check that on $U_1$ one has $J_1=J$. In fact, on $U_1$ we have
$g=g_1$ and so
$$
\omega (Ju,Jv) = \omega (u,v) = g(A_1u,v) = \omega (A_1u,Jv),
$$
and thus $A_1=J$ on $U_1$. So $Q_1=\sqrt{-A_1^2}=\id $. From here
we have $J_1=A_1Q_1^{-1}=A_1=J$ on $U_1$.

\state Remark. Note that at this point of the proof  we
construct a correspondence $P:g\to J$, which maps a Riemannian metric
$g$ to an $\omega $-calibrated a.-c.\ structure $J$ with $g(u,v)=\omega (u,Jv)
$. This correspondence is obviously a continuous map from the space of
metrics to the space of structures. Moreover, note that if $g$ appears as
$g=g_J(u,v)=\omega (u,Jv)$ for some $\omega $-calibrated $J$, then $P(g_J)=J$.

\smallskip
Let us now check that $\omega $ is $J_1$-calibrated. First:
$$
\omega (J_1u,J_1v)=g_1(A_1J_1u,J_1v)=g_1(A_1^2Q_1^{-1}u,J_1v)=
g_1(Q_1u,-J_1v)=
$$
$$
=-g_1(Q_1u,J_1v)=-g_1(u,Q_1J_1v)=-g_1(u,A_1v)=-g_1(A_1v,u)=-\omega (v,u)=
$$
$$
=\omega (u,v).
$$
Second:
$$
\omega (u,J_1u)=g_1(A_1u,J_1u)=g_1(A_1u,A_1Q_1^{-1}u)=
$$
$$
=g_1(u,u)>0
$$
for nonzero $u$.\qed

This remark leads to the following

\state Corollary 1.2.2. {\it Let $(X,\omega )$ be as above. Then:

1. There exists a $\omega $-calibrating almost-complex structure on $X$.

2. The space of $\omega $-calibrating almost-complex structures is
contractible.}

\state Proof. The first statement of this corollary is a particular case
of the proposition when $U=\emptyset $.

To prove the second one, we fix some calibrating structure $J_0$. Denote
by $g_{J_0}$ the corresponding Riemannian metric on $X$, \ie, 
$g_{J_0}(u,v)=\omega (u,J_0v)$. The space  of Riemannian metrics is a
convex cone $\calc$ in $\hom_{\rr }(TX,T^*X)$. Therefore there exists 
a contraction
$\Psi : \calc \times [0,1]\to \calc $ to $g_{J_0}$,
\ie, $\Psi (\cdot , 0)=\id $, $ \Psi (\cdot ,1 )\equiv g_{J_0}$ and
$\Psi (g_{J_0},t)\equiv g_{J_0}$.

Consider the following map $\Phi :\calj_{\omega }^c\to \calc $ from the space
of $\omega $-calibrating structures  to the space of metrics:
$$
\Phi (J)(u,v)=g_{J}(u,v)=\omega (u,Jv).
$$
In the proof of Proposition 1.2.1 (see Remark there) we showed that
$$
P\scirc \Phi =\id :\calj_{\omega }^c\to \calj_{\omega }^c.
$$

Now $P\scirc \Psi (\cdot ,t)\scirc \Phi $ will be a contraction of
$\calj_{\omega }^c$ to $J_0$.
\qed

\def\bfx{{\bold x}}
\def\bfy{{\bold y}}
\def\sfL{\mathop{\sf L}}
\def\sfK{\mathop{\sf K}}
\def\bfone{\mathop{1\mkern-5mu{\rom l}}}
\def\sft{{\sf t}}

Let us prove the following 

\state Proposition 1.2.3. {\it For any symplectic manifold $(X,\omega)$, the 
set $\calj_\omega$ of $\omega$-tame almost complex structures on $X$ is a
non-empty contractible manifold}.

\medskip
The proposition follow immediately from the following result from linear 
algebra.

\state Lemma 1.2.4. {\it Let $V$ be a (finite-dimensional) real vector space
and $\omega$ a linear symplectic form on $V$. Then the set $\jj_\omega$ of 
$\omega$-tame linear complex structures on $V$ is a non-empty open contractible
subset in the set $\jj$ of all linear complex structures on $V$}.

\state Proof. {\sl Step 1}. The linear version of the Darboux theorem states 
that there exists a basis $(\bfx_1, \ldots, \bfx_n, \bfy_1, \ldots, \bfy_n)$ 
of $V$ such that for corresponding linear coordinates $(x_1, \ldots, x_n, 
y_1, \ldots, y_n)$ on $V$ 
$$
%\textstyle
\omega=\sum_{i=1}^n dx_i \wedge dy_i.
\eqno\rom(1.2.3)
$$
Set $J_0 (\bfx_i) \deff \bfy_i$, $J_0 (\bfy_i) \deff -\bfx_i$. It is easy to
see that $J_0$ is a linear complex structure on $V$. Moreover, the $\cc$-valued
linear functionals $z_i \deff x_i + \isl\, y_i$ induce the isomorphism $\phi
:(V,J_0) \cong (\cc^n, J\st)$, $J\st \deff \isl$, of {\sl complex} linear 
spaces, such that 
the given symplectic form $\omega$ on $V$ is mapped into the standard 
sympletic form $\omega\st$ on $\cc^n$ given by the same formula (A.T.1), \ie 
$\phi^*\omega\st = \omega$. Thus we can conclude that there exists a 
$J_0$-Hermitian metric $h$ on $V$ with the following properties:

\sli $g(v,w) \deff \re h(v,w)$ is a symmetric positively definite form on
$V$ such that $h(v,w)= g(v,w) + \isl\,\omega(v,w)$ for any $v,w \in V$; 

\slii $g(J_0 v, J_0w) = g(v,w)$ and $\omega(J_0 v, J_0w) = \omega(v,w)$.

\sliii $g(v,w) = \omega(v, J_0w)$ and $\omega(v,w)=g(J_0v,w)$ for any $v,w 
\in V$. 

\noindent
In particular, $\omega(v, J_0w)$ is positive definite, \ie $J_0 \in 
\jj_\omega$. This shows that $\jj$ is non-empty.

\medskip\noindent
{\sl Step 2. For any $J \in \jj_\omega$ the operator $J + J_0$ is invertible}.
Really, otherwise we have $Jv = - J_0v$ for some non-zero $v\in V$ and then
$0< \omega(v, Jv) = -\omega(v, J_0v) <0$, a contradiction.

For $J \in \jj$ with $J + J_0$ invertible we set $\sfL(J) \deff -(J - J_0) (J+
J_0)\inv$. The equivalent definitions are 
$$
\matrix\format \r&\,\c\, & \l \\
\sfL (J) &=&  -(J - J_0) (J+J_0)\inv = 
-\bigl((J - J_0)J_0\bigr) \cdot \bigl( (J+J_0) J_0 \bigr)\inv \cr
\vph&=&(\bfone + JJ_0 )(\bfone - JJ_0  )\inv 
      = (\bfone - JJ_0)\inv (\bfone + JJ_0).
\endmatrix
\eqno\rom(1.2.4)
$$

\medskip\noindent
{\sl Step 3. If $J + J_0$ is invertible, then $W= \sfL(J)$ is 
$J_0$-antilinear, \ie $W J_0= - J_0 W$}. Really, 
$$
\matrix \format \c &\,\c\, &\r &\,\c \,&\r &\,\c\, &\c
\vbox to 11pt{} \\
W J_0 &=& -(J - J_0) (J+J_0)\inv J_0 &=& 
(\bfone - JJ_0)\inv (\bfone + JJ_0) J_0 &=&
\cr 
&&
(\bfone -JJ_0)\inv (J_0- J)&=&
- (\bfone -JJ_0)\inv J(\bfone + JJ_0) &=& 
\cr
&&
(\bfone -JJ_0)\inv J\inv (\bfone +JJ_0) &=& 
(J - J^2 J_0)\inv (\bfone +JJ_0) &=&
\cr
&&
(J_0 + J)\inv (\bfone +JJ_0) &=& 
\bigl( (\bfone -JJ_0) J_0 \bigr)\inv (\bfone +JJ_0) &=& 
\cr
&&
J_0\inv(\bfone -JJ_0)\inv (\bfone +JJ_0) &=& 
- J_0(\bfone -JJ_0)\inv (\bfone +JJ_0) &=& - J_0 W
\endmatrix
$$
 
\medskip\noindent
{\sl Step 4. $J \in \jj_\omega$ implies $\bfone - W^\sft W \gg 0$ for $W= 
\sfL(J)$}, \ie the operator $\bfone - W^\sft W$ is positively definite, 
$g(v,v) > g(Wv,Wv)$ for any non-zero $v \in V$. The conjugation of $W \in 
\endo(V)$ is done using the metric $g$.

First we note that $J_0^\sft = J_0\inv =-J_0$. Since $\bfone -JJ_0$ is 
invertible, the positivity of $\bfone - W^\sft W$ is equivalent to the 
positivity of $(\bfone -JJ_0)^\sft(\bfone - W^\sft W) (\bfone -JJ_0)$. 
Computing this operator we obtain
$$
\matrix\format \c \vbox to 11.5pt{} \\
(\bfone -JJ_0)^\sft(\bfone - W^\sft W) (\bfone -JJ_0) = 
(\bfone -JJ_0)^\sft(\bfone -JJ_0) - 
\cr
(\bfone -JJ_0)^\sft\bigl( (\bfone +JJ_0) (\bfone -JJ_0)\inv \bigr)^\sft
\bigl( (\bfone +JJ_0) (\bfone -JJ_0)\inv \bigr)(\bfone -JJ_0)=
\cr
(\bfone -JJ_0)^\sft(\bfone -JJ_0) -
(\bfone +JJ_0)^\sft (\bfone +JJ_0) =
\cr
(\bfone +J_0 J^\sft)(\bfone - JJ_0) - (\bfone - J_0 J^\sft)(\bfone +JJ_0)=
\cr
\bfone +J_0 J^\sft - JJ_0  -  J_0 J^\sft JJ_0 -
\bfone + J_0 J^\sft - JJ_0 +  J_0 J^\sft JJ_0 =
\cr
2(J_0 J^\sft - JJ_0).
\endmatrix
$$

Now, writing any non-zero $v\in V$ as $v= J_0w$, we obtain
$$
\matrix\format \l &\,\c\,\;& \l \vbox to 12.5pt{} \\
g(v, (J_0 J^\sft - JJ_0)v) &=&  g(J_0w, (J_0 J^\sft - JJ_0)J_0w)= 
\cr
g(J_0w, J_0J^\sft J_0w) + g(J_0w, Jw) &=& 
g(w, J^\sft J_0w) + g(J_0w, Jw) = 
\cr
g(Jw, J_0w) + g(J_0w, Jw) &=& 2\omega(w, Jw)>0.
\endmatrix
$$

Thus $\sfL$ maps $\jj_\omega$ into the set 
$$
\calw \deff 
\{\, W \in \endo(V)\;:\; W J_0 =- J_0W,\; \bfone -W^\sft W \gg0 \,\}.
\eqno\rom(1.2.5)
$$

\medskip
It is easy to see that $\calw$ is contractible. Thus it is sufficient to
show that $\sfL: \jj_\omega \to \calw$ is a diffeomorphism. Explicit 
calculations show that the inverse map should be given by
$$
J= \sfK(W) \deff J_0 (\bfone + W) (\bfone - W)\inv= 
J_0 (\bfone - W)\inv (\bfone + W) = 
J_0 {\bfone + W \over \bfone - W}
\eqno\rom(1.2.6)
$$

\medskip\noindent
{\sl Step 5. Let $W\in \calw$. Then $J= \sfK(W) = J_0 {\bfone + W \over 
\bfone - W}$ is well-defined and $J^2 = -\bfone$}. Really, the condition
$\bfone -W^\sft W \gg0$ implies that $\bfone \pm W$ is invertible, and then 
$$
\matrix\format \c &\,\c\,&\l  \vbox to 11.5pt{} \\
J^2 &=& 
-J_0 (\bfone + W) (\bfone - W)\inv J_0\inv (\bfone + W) (\bfone - W)\inv
 \cr
&=& - (\bfone + J_0 W J_0\inv) (\bfone - J_0 W J_0\inv)\inv
(\bfone + W) (\bfone - W)\inv
 \cr
&=& - (\bfone - W ) (\bfone + W )\inv
(\bfone + W) (\bfone - W)\inv \qquad  =- \bfone 
\endmatrix
$$

\medskip\noindent
{\sl Step 6. Let $W\in \calw$ and $J= \sfK(W)$. Then $J+ J_0$ is invertible
and $J \in \jj_\omega$}. In fact, $J + J_0 =  J_0\left( \bfone + {\bfone + 
W \over \bfone - W} \right) = {2 J_0 \over \bfone - W}$ is invertible. Now, 
repeating argumantations from {\sl Step 4} the positivity of $\bfone - 
W^\sft W$ is equivalent to the positivity of $J_0 J^\sft - JJ_0$, which means 
the tameness of $J$.

The lemma is proved.\qed

\state Remark. One can consider $\sfK:\calw \to \jj_\omega \subset \jj$ 
a generalized Cayley transformation defined on a bounded domain $\calw$ in 
the complex linear space $\barr \endo_\cc(V,J_0)$ of $J_0$-antilinear
endomorphisms of $V$.

\medskip\noindent
{\bigsl 1.3. Almost-Complex Submanifolds, Complex Curves, Energy and Area.}

\smallskip

Let $J$ be some linear complex structure on $\rr^{2n}$, \ie, $J\in
\endo(\rr^{2n}), J^2=-\id $. Let $h$ be some $J$-Hermitian metric on
$\rr^{2n}$, $\omega =\omega_h$ its associated $(1,1)$-form and $g=g_h$
corresponding Riemannian metric. Put $\sigma_k =  {1\over k!}\omega^k$.

\smallskip
\state Wirtinger Inequality. {\it For any $g$-orthonormal system
$v_1,...,v_{2k}$ in $\rr^{2n}$ one has
$$
\vert \sigma_k(v_1,...,v_{2k})\vert \le 1, \eqno(1.3.1)
$$
with equality taking place iff the subspace $\<v_1,...,v_{2k}\>$ is $J$ -
invariant.
}

\smallskip
\state Remark. \rm We shall also call $J$-invariant subspaces  $J$-complex.

\state Proof. We shall prove (1.3.1) by induction on $k$.

Let $k=1$. Then $\vert \omega (v_1,v_2)\vert =\vert g(v_1,Jv_2)\vert
\le \Vert v_1\Vert_g\cdot \Vert Jv_2\Vert_g=$ $\Vert v_1\Vert_g\cdot \Vert v_2
\Vert_g$, with equality iff $v_1$ and $Jv_2$ are collinear. Thus the
subspace $\<v_1,v_2\>$ is $J$-invariant.

Now let $V$ be a subspace of dimension $2k$. Put $\omega'=\omega\mid_V$.
Find a $g$-orthonormal base $e_1,...,e_{2k}$ of $V$, s.t. $\omega'=
\lambda_1e^1\wedge e^2 + ... +\lambda_ke^{2k-1}\wedge e^{2k}$. By the
case $k=1$ we have $\vert \omega' (e_{2p-1},e_{2p})\vert=
\vert \omega (e_{2p-1},e_{2p})\vert =\vert \lambda_p\vert \le 1$, with
equality taking place iff $e_{2p}={+\over -}Je_{2p-1}$. So, for $\sigma'_k
=\sigma_k\mid_V$ we obtain
$$
\vert \sigma_k'(e_1,...,e_{2k})\vert = \vert {1\over k!}\omega^k(e_1,...,
e_{2k})\vert = \vert \lambda_1...\lambda_k\vert \le 1,
$$
with equality iff $e_{2p}={+\over -}Je_{2p-1}$  for all $1\le p\le k$.
\qed

\smallskip
A submanifold $Y$ of an almost-complex manifold $(X,J)$ is called an
almost-complex (or $J$-complex) if the tangent spaces to $Y$ are
$J$-invariant.

\smallskip
\state Definition 1.3.1. {\it Nowhere degenerate exterior two-form $\omega$ 
on $X$ is called symplectic if $d\omega =0$.}

\smallskip 
A pair $(X,\omega )$ is called a symplectic manifold.

\state Corollary 1.3.1. {\it Let $(X,w,J)$ be a symplectic manifold with $w$
being $J$-calibrated. Then $J$-complex submanifolds are minimal and
their volume is given by
$$
vol_{2k}(Y) = {1\over k!}\int_Yw^k. \eqno(1.3.2)
$$}

\state Proof. Let $Y$ be a $J$-complex submanifold of $X$, and $Y_1$
some submanifold of dimension $2k$ as well as $Y$. Suppose $\d Y=\d Y_1$ and
$Y\cup (-Y_1) \sim 0$. Denote by $dy$ and $dy_1$ the volume forms of $Y$ and
$Y_1$ with respect to metric $g(u,v)=\omega (u,Jv)$. Then we have
$$
0=    \int_{Y\cup (-Y_1)}\sigma_k = \int_{Y}\sigma_k  - \int_{Y_1}\sigma_k  =
\int_{Y}\sigma_k (TY)dy - \int_{Y_1}\sigma_k(TY_1)dy_1\ge
$$
$$
\ge \int_{Y}dy  - \int_{Y_1}dy_1=vol(Y)-vol(Y_1),
$$
with equality taking place if and only if $\sigma_k(T_pY_1)=1$
for all $p\in Y_1$, \ie, when $Y_1$ is also $J$-complex.
\qed

\state Definition 1.3.2. {\it A $C^0\cap L^{1,2}$-map $u :(Y,J_Y)
\to (X,J_X)$ is called holomorphic if for a.-a. $x\in Y$
$$
du_x\scirc J_Y(x)= J_X(u(x))\scirc du_x\eqno(1.3.3)
$$
\noindent as mappings  $T_xY\to T_{u(x)}X$.
}

\smallskip
\state Exercise. {\rm Check that for complex valued functions, \ie,
$u :(\rr^2,i)\to (\rr^2,i)$, (1.3.3) is a Cauchy-Riemann equation.
}

\smallskip
Sometimes one calls such $u$ a $(J_Y,J_X)$-holomorphic map. In the
special case when $(Y,J_Y)$ is a Riemann surface $(S,J_S)$, one calls
$(S,J_S,u)$ a parameterized $J_X$-complex curve. Its image $u(S)$
one simply calls a $J_X$-complex curve.

If some $J$-calibrated exterior $2$-form $\omega $ on an almost-complex
manifold $(X,J)$ is chosen, one defines the $\omega $-area of a $J$-complex
curve $u:(S,J_S)\to (X,J)$ as
$$
area[u(S)] = \int_Su^*\omega .\eqno(1.3.4)
$$
Remark that if $\omega $ is symplectic, then by {\sl Corollary 1.3.1} this 
area coincides  with $g$-area  for the metric $g(\cdot ,\cdot ) =
\omega (\cdot ,J\cdot )$ and $J$-complex curves are $g$-minimal surfaces.

Remark also that for a $J$-complex curve $u:(\Delta ,J\st )\to (X,J)$,
parameterized by a standard unit disk, one has
$$
\int_{\Delta }u^*\omega = \int_{\Delta }\omega  (du({\d \over \d x}),
du({\d \over \d y}))dx\wedge dy =  -\int_{\Delta }{1\over 2}g(du({\d \over
\d x}),Jdu(i{\d \over \d x}))dx\wedge dy =
$$
$$
=-\int_{\Delta }{1\over 2}g(du({\d \over \d x}),J^2du({\d \over \d x}))dx
\wedge dy = \int_{\Delta }{1\over 2}g({\d u\over \d x}, {\d u \over \d x}))
dx\wedge dy = {1\over 2}\Vert {\d u\over \d x}\Vert^2_{L^2(\Delta ,X)}.
$$
Thus
$$
\int_{\Delta }u^*\omega  = \Vert {\d u\over \d x}
\Vert^2_{L^2(\Delta ,X)}+\Vert {\d u\over \d y}\Vert^2_{L^2(\Delta ,X)}.
\eqno(1.3.5)
$$
The right hand side of (1.3.5) is called an energy of a $C^0\cap L^{1,2}$-map
(not necessary holomorphic!) from $\Delta $ to a Riemannian manifold $(X,g)$. 
Thus, (1.3.5) tells us that, for holomorphic maps, energy and the area of
the image, taken with multiplicities, are equal.

\state Remark. %Our definition of the area uses the following fact. 

Let $g$ be a Riemannian
metric on $C$ compatible with $j_C$, $h$ a Riemannian metric on $X$, and $u:C
\to X$ a $J$-holomorphic immersion. Then $\norm{du}^2_{L^2(C)}$ is independent
of the choice of $g$ and coincides with the area of the image $u(C)$ \wrt
the metric $h_J(\cdot, \cdot) \deff \half(h(\cdot, \cdot) + h(J\cdot, J\cdot)
)$. The metric $h_J$ here can be seen as a ``Hermitization'' of $h$ \wrt
$J$. The independence of $\norm{du}^2_{L^2(C)}$ of the choice of a metric $g$
on $C$ in the same conformal class is a well-known fact, see, \eg, [S-U]. Thus
we can use the flat metric $dx^2 + dy^2$ to compare the area and the energy.
For a $J$-holomorphic map we get
$$
\norm{du}^2_{L^2(C)}= \int_C |\d_x u|_h^2 + |\d_y u|_h^2 =
\int_C |\d_x u|_h^2 + |J \d_x u|_h^2 = \int_C |du|_{h_J}^2 = 
\area_{h_J}(u(C)),
$$
where the last equality is another well-known result, see, \eg, [G]. Since we
consider changing almost complex structures on $X$, it is useful to know that
we can use any Riemannian metric on $X$ having a reasonable notion of area.

\bigskip\noindent{\bigsl 1.4. Symplectic Surfaces.}

\smallskip

\rm First  we recall some elementary
facts about orthogonal complex structures in $\rr^4$.

In $\rr^4$ with coordinates $x_1, y_1, x_2, y_2$ consider the standard
symplectic form $\omega\st  = dx_1\wedge dy_1 + dx_2\wedge dy_2 $ and the
standard
complex structure $J\st$ defined by the operator
$$
J\st =
\left( \matrix
0&-1&0&0\cr
1&0&0&0\cr
0&0&0&-1\cr
0&0&1&0\cr
\endmatrix \right).
$$

\medskip
Let  $\jj$ denote the set of all orthogonal complex structures in $\rr^4$
giving $\rr^4$ the same orientation as $J\st$. Orthogonality here means
just that $J$ is an orthogonal matrix. The same orientaion means that for any 
pair $x_1, x_2$ of $J$-independent vectors the basis $x_1, Jx_1, x_2, Jx_2$ 
gives the same orientation of $\rr^4$ as $J\st$. 

\smallskip
\state Exercise. {\rm Prove that this orientation does not depend on the 
particular choice of $x_1, x_2$ and coincides with the orientation given 
by $\omega\st^2$.
}

\smallskip
 The following lemma summarizes the
elementary facts which we need for the sequel.

\state Lemma 1.4.1. \it The elements of $\jj$ have the form
$$
J_{s} =
\left(
\matrix
0&-s &c_1&c_2\cr
s &0&c_2&-c_1\cr
-c_1&-c_2&0&-s \cr
-c_2&c_1&s &0\cr
\endmatrix \right),\eqno(1.4.1)
$$
with $c_1^2 + c_2^2 + s^2 = 1$. One also has for $x\in \rr^4$
$$
\omega\st  (x, J_s x) = s \Vert x\Vert^2 . \eqno(1.4.2)
$$

\smallskip
\rm
We remark that the set $\jj$ is a unit
two-dimensional sphere $S^2$ in $\rr^3$ with coordinates $c_1, c_2,s $. We
note also that the number $\omega\st (x, J_s x)$ does not depend on the
choice of a \sl unit \rm vector $x$. One also remarks that the standard
structure corresponds to the north pole of $S^2$ and structures tamed by
$\omega\st $ constitute the upper half-sphere.

\smallskip
\state Exercise. 1. Note that $\jj = \{ J\in O(4): J^2=-E\} = \{ J: 
J^t=J^{-1}=-J\} $. Prove that condition $J^t=-J$ implies that 
$$
J = \left(
\matrix
0&-s &c_1&c_2\cr
s &0&c_3&c_4\cr
-c_1&-c_3&0&-t \cr
-c_2&-c_4&t &0\cr
\endmatrix \right), 
$$
for some $s,t,c_i, i=1,...,4.$

\noindent
{\bf 2.} Denote by $C$ the matrix $\matrix c_1& c_2 \cr c_3& c_4\endmatrix$.
Prove that condition $J^2=-E$ implies that:

\noindent
a) 
$$
C = \left(
\matrix
\pm c_1&c_2\cr
\mp c_2&c_1\cr
\endmatrix \right)
$$
with $c_1^2+c_2^2=1-s^2$;

\noindent
b) $\vert t\vert =\vert s\vert $;

\noindent
c) and, finally, that $s=t$ and for $s\not= \pm 1$ the matrix
${1\over \sqrt{1-s^2}}\cdot C$ lies in $O_-(2)$ (the set of orthogonal 
$2\times 2$-matrices with determinant $-1$) iff the structure 
$$
J_s := 
\left(
\matrix
0&-s &c_1&c_2\cr
s &0&c_3&c_4\cr
-c_1&-c_3&0&-s \cr
-c_2&-c_4&s &0\cr
\endmatrix \right) 
$$
defines the same orientation of $\rr^4$ as $J\st = J_1$.

\noindent
{\bf 3.} Check finally that $\omega\st (x,J_sx)=s\cdot \Vert x\Vert^2$.

\smallskip  
Let $M_1$ and $M_2$ be two smooth oriented surfaces in the unit ball $B\subset
\rr^4$ with zero as a common point. Let $v_1, w_1$ and $v_2, w_2$ be oriented
bases of $T_0M_1$ and $T_0M_2$, respectively. Suppose that $M_1$ and $M_2$
intersect transversally at zero, {\sl \ie}, $v_1, w_1, v_2, w_2$ is the basis
of $\rr^4$. We say that they intersect positively if this basis gives
the same orientation of $\rr ^4$ as the standard one.

\smallskip
Let $(X,\omega )$ be a manifold with nowhere degenerate exterior $2$-form
$\omega $.

\smallskip\noindent
\state Definition 1.4.1. {\it An immersion $u:S\to X$ of a real surface
$S$ into $X$ is called $\omega $-positive if $u^*\omega $ never vanishes.
}

\smallskip 
If $\omega $ is symplectic we call such immersions {\it symplectic}.

\smallskip\noindent
\state Definition 1.4.2. {\it An almost-complex structure $J$ is said to be
tamed by an exterior $2$-form $\omega $ if $\omega (u,Ju)>0$ for any
nonzero $u\in TX$}.

\smallskip
In the following lemma we suppose for simplicity of proof that $\dim_{\rr }X
=4$.

\state Lemma 1.4.2. {\it Let $M$ be a $\omega $-positive
compact surface immersed into $(X,\omega )$ with only double positive local
self-intersections, and let $U_1\subset \subset U$ be  neighborhoods of $M$.
Then for any given $\omega $-tamed a.-c.\ structure $J$ there
exists a smooth family $\{ J_t\} _{t\in [0, 1]} $ of almost-complex
structures on $X$ such that:

\smallskip
a) $J_0$ is the  given structure $J$ on $X$;

b) for each $t\in [0, 1]$ the set $\{ x\in X: J_t(x) \not= J_0(x) \} $ is
contained in $U_1$;

c) $M$ is $J_1$-holomorphic;

d) all $\{J_t\}$ are tamed by the given form $\omega $, \ie,
$\omega(v,J_tv) > 0$ for every nonzero $v\in TX$.
}

\smallskip
\state Proof. Let $N$ be a normal bundle to $M$ in $X$ and $V_1$ a
neighbourhood of the zero section in $N$. Shrinking $V_1$ and $U_1$ we can
assume
that $U_1$ is an image of $V_1$ under an $\exp$-map in $J$-Hermitian metric
$h$ associated to $\omega $. More precisely, one should take the  Riemannian
metric associated to $h$, \ie $g = \re h$. Shrinking $V_1$ and $U_1$
once more, if necessary, we can extend the distribution of tangent planes to
$M =$ (zero section of $N$) to the distribution $\{ L_x\} _{x\in V_1}$ of
$\omega $-positive planes on $V_1$. Here we do not distinguish between
$\omega$ and its lift onto $V_1$ by $\exp$. Denote by $N_x$ the subspace
in $T_xV_1$ which is $g$-orthogonal to $L_x$. Perturbing $g$ we can choose
the distribution $\{ L_x\}$
in such a way that if $\exp(x)= \exp(y)$ for some $x\not= y$ from $V_1$ then
$d\exp_x(L_x)=d\exp_y(N_y)$ and $d\exp_y(L_y)=d\exp_x(N_x)$. In particular,
we suppose that $M$ intersects itself $g$-orthogonally. For every $x\in V_1$
choose an orthonormal basis $e_1(x), e_2(x)$ of $L_x$ such that $\omega_x(
e_1(x), e_2(x)) > 0$, and an orthonormal basis $e_3(x), e_4(x)$ of $N_x$ such
that the basis $e_1(x), e_2(x), e_3(x), e_4(x)$ gives the same orientation of
$T_xV_1$ as $\omega^2$.

Define an almost complex structure $J^1$ on $U_1$ by $J^1e_1 = e_2$, $J^1 e_3 =
e_4$. Then
$J^1$ depends smoothly on $x$, even when $e_j(x)$ are not smooth
in $x$. Furthermore, we have  $\omega _x(e_1(x), J_x^1e_1(x)) > 0$ 
and thus
from {\sl Lemma 1.1.1} we see that the relation
$\omega _x(e_3(x), J_x^1e_3(x)) > 0$ is also satisfied.
This means
that our $J^1$ is tamed by $\omega $. Note also that $M$ is $J^1$-holomorphic.

Denote by $\jj_x$ the sphere of $g$-orthogonal complex structures on $T_x
X$ as
in {\sl Lemma 1.4.1}. Let $\gamma_x$ be the shortest geodesic on $\jj_x$
joining $J(x)$ - our given integrable structure, with $ J_x^1$. Put $J_t(x)
= J^1_{\gamma _x(t\cdot \Vert \gamma _x\Vert \cdot \phi (x ))}$. Here $\Vert
\gamma_x \Vert $ denotes the length of $\gamma _x$, $\phi$ is smooth on $X$
with support in $U_1$ and identically one in the neighborhood of $M$. The
curve $\{ J_t \}$ satisfies all the~conditions of our lemma.
\qed

\smallskip\noindent
\state Lemma 1.4.3. {\it Under the conditions of Lemma 1.4.2 the structure
$J_1$ can be made complex in the neighborhood of $M$. 
}

\smallskip
The proof is left to the reader.

\bigskip\noindent
{\bigsl 1.5. Adjunction Formula for Immersed Symplectic Surfaces.}

\nobreak
\smallskip
Let us now prove the Adjunction Formula for \sl immersed symplectic \rm
surfaces.
Let $u:\bigsqcup_{j=1}^dS_j\to (X,\omega)$ be a reduced compact symplectic
surface (see Definition 1.4.1) immersed into a symplectic four-dimensional
manifold. Let $g_j$ denote
the genus of $S_j$ and $M_j=u(S_j)$. Put $M:=\bigcup_{j=1}^dM_j$ and denote
by $[M]^2$ the homological
self-intersection number of $M$. Define a geometrical self-intersection number
$\delta $ of $M$ in the following way. Perturb $M$ to obtain a symplectic
surface $\widetilde M$ with only transversal double points. Then $\delta$
will be
the sum of indices of the intersection over those double points. Those indices
can be equal to $1$ or $-1$.

\state Exercise. Let $P$ be an oriented plane in $\rr^4$. Call $P$ symplectic 
if $\omega\st (P)>0$. Find two symplectic planes $P_1,P_2\in \rr^4$ which
intersect transversally at origin with intersection index $-1$.

\medskip
Let $J$ be some  almost complex structure which is tamed by
$\omega$, {\sl \ie,} $\omega(\xi, J\xi)>0$ for any nonzero $\xi \in TX$.
Denote by $c_1(X, J)$ the first Chern class of $X$ with respect to
$J$. Since, in fact, $c_1(X, J)$ does not depend on continuous changes of $J$
and since  the set of $\omega$-tamed almost complex structures is
contractible, we usually omit the dependence of $c_1(X)$ on $J$.

\state Lemma 1.5.1. {\it Let $M=\bigcup_{j=1}^d M_j$ be a compact immersed
symplectic surface in four-dimensional symplectic manifold $X$. Then
$$
\sum_{j=1}^d g_j = {[M]^2 - c_1(X)[M]\over2} + d - \delta .\eqno(1.5.1)
$$}

\state  Proof.
By replacing every $M_j$ by its small perturbation, we can suppose that $M_j$
has only transversal double self-intersection points. Let $N_j$ be a normal
bundle to $M_j$ and let $\widetilde M_j$ denote the zero section
of $N_j$. Also let $\exp_j $ be the exponential map from a~neighborhood $V_j$
of $\widetilde M_j\subset N_j$ onto the neighborhood $W_j$ of $M_j$.
Lift $\omega$ and $J$ onto $V_j$. Since $\widetilde M_j$ is embedded to $V_j$,
we can apply {\sl Lemma 1.4.2} to obtain the $\omega$-tame
almost complex structure $J_j$ on $V_j$ such that  $\widetilde M_j$ is
$J_j$-holomorphic.

For every $j$ we now have the following exact sequence of complex bundles:

\smallskip
$$
0\longrightarrow TS_j {\buildrel du\over \longrightarrow }
E_j \buildrel \pr \over \longrightarrow
N_j \longrightarrow 0 \eqno(1.5.2)
$$

\smallskip
Here $E_j=(u^*TX)\ogran_{S_j}$ is endowed with complex structure given by
$J_j$.
Since $du$ is nowhere degenerate {\sl complex} linear morphism, $N_j \deff
E_j/du(TM_j)$ is a complex line bundle over $S_j$. From (1.5.2) we get

\smallskip
$$
c_1(E_j) = c_1(TS_j) + c_1(N_j). \eqno(1.5.3)
$$

\smallskip
Observe now that $c_1(E_j) = c_1(X)[M_j]$ and that $c_1(TS) = \sum_{j=1}^d
c_1(TS_j) = \sum_{j=1}^d(2-2g_j) = 2d - 2\sum_{j=1}^d g_j$. Furthermore,
$c_1(N_j)$ is the algebraic number of zeros of a generic smooth section
of $N_j$. To compare this
number with the self-intersection of $M_j$ in $X$, note that if we move $M_j$
generically to obtain $M_j'$, then the~intersection number $int(M_j,
M_j')$ is equal to the algebraic number of zeros of generic section
of $N_j$ plus two times the
sum of intersection numbers of $M_j$ in self-intersection points, {\sl \ie,}
$[M_j]^2 = c_1(N_j) + 2\delta_j$. So
$$
c_1(X)[M_j] = 2 - 2 g_j + [M_j]^2 - 2\delta _j . \eqno(1.5.4)
$$
Now it only remains to take the sum over $j=1,\ldots, d$ and to
remark that the intersection points of $M_i$ with $M_j$ for $i\not=j$ 
make a double
contribution to $[M]^2$.
\qed

\newpage

\smallskip\noindent
{\bigbf Appendix 1}

\smallskip\noindent
{\bigbf Chern Class and Riemann-Roch Formula}

\medskip\noindent
{\bigsl A1.1. First Chern Class.}

\smallskip\rm
Let $L\to M$ be a complex line bundle over a real manifold $M$. One of
the possible definitions consists of taking a real rank two bundle over $M$
with an operator $J\in \End (L)$, satisfying $J^2=-\id_L$. One can than
locally find a frame $e_1(x), e_2(x)$ with $Je_1(x)=e_2(x)$. This gives
a covering $\{ U_\alpha \} $ of $M$ together with isomorphisms of
complex line bundles $\phi_{\alpha }:L\mid_{U_{\alpha }}\to U_{\alpha }
\times \cc $, \ie, a standard definition of a complex line bundle. Sometimes
we shall mark as $e_1^J, e_2^J, \phi_{\alpha }^J$ the corresponding objects
to underline their dependence on $J$.

Denoting by $\cala=\cala_M$ and  by $\cala^*=\cala_M^*$ the
sheaves of complex valued (resp.\ complex valued nonvanishing) functions on
$M$, we observe the following exact sequence

$$
0\to \zz \buildrel{i}\over {\longrightarrow} \cala \buildrel{exp(2\pi i
\cdot )}\over {\longrightarrow} \cala^* \to 0. \eqno(A1.1.1)
$$

Here $i$ is an imbedding of the sheaf of locally constant integer valued
functions into $\cala$, end $exp(2\isl\pi \cdot ):f\to e^{2\isl\pi f}$.
The sequence(A1.1.1) gives rise to the following long exact sequence 
of \v Cech
cohomologies
$$
0=\sfh^1(M,\cala)\buildrel{\exp(2\isl\pi  \cdot)}\over {\longrightarrow}
\sfh^1(M,\cala^*\buildrel{\delta }\over {\longrightarrow} \sfh^2(M,\zz )
\to \sfh^2(M,\cala)=0. \eqno(A1.1.2)
$$
Equalities $\sfh^1(M,\cala)=\sfh^2(M,\cala)=0$ follow from the fact that
the sheaf $\cala$ admits a partition of unity.

Classes from $\sfh^1(M,\cala^*)$ are the defining cocycles of complex line
bundles - one more possible definition of the line bundle. In terms of
local trivializations $\{ \phi_{\alpha }\} $ such cocycles can be obtained
as $\phi_{\alpha ,\beta }=\phi_{\alpha }\scirc \phi_{\beta }^{-1}
:(U_{\alpha }\cap U_{\beta })\times \cc \to (U_{\alpha }\cap U_{\beta })
\times \cc $, \ie, $\phi_{\alpha ,\beta }\in \cala^*_{U_{\alpha }\cap 
U_{\beta }}$.

\state Definition A1.1.1. {\it If  $\{ \phi_{\alpha ,\beta }\}
\in \sfh^1(M, \cala^*)$ is a defining cocycle of a complex line bundle $L$,
then $\delta (\{ \phi_{\alpha ,\beta }\} )\in \sfh^2(M, \zz)$ is called the
first Chern class of $L$ and is usually denoted as $c_1(L)$.
}

\smallskip\rm
For the complex bundle $E$ of complex rank $r$ the first Chern class is
defined as $c_1(\Lambda^rE)$. If $E=TX$,  the tangent bundle to an
almost-complex manifold $X$, then one simply writes $c_1(X)$ or $c_1(X,J)$ if
an almost-complex structure is needed to be specified. 

If an almost-complex structure $J$ on the real bundle $E$ varies continuously,
then the corresponding trivializations $\{ \phi_{\alpha }\} $ above
(on $\Lambda^r E$) can be obviously chosen to also vary continuously.
Thus $c_1(E,J)$ varies continuously. But $c_1(E,J)\in \sfh^2(M, \zz )$, \ie,
takes values in a discrete group. So, it does not change at all. This
simple but important observation together with {\sl Corollary 1.2.2} leads to
the following

\state Corollary A1.1.1. {\it Let $\omega $ be a nondegenerate exterior two-form
on the even-dimensional real manifold $X$. Then $c_1(X,J)$ does not depend
on the choice of $\omega $-calibrating (and even $\omega $-compatible) 
almost-complex structure $J$.
}

\medskip\noindent\sl
A1.2. Riemann-Roch Formula and index of $\dbar $-type operators.

\smallskip\rm
For the holomorphic bundle $E$ of complex rank $r$ over a compact Riemann
surface $S$ denote by
$\calo_E$ the sheaf of its holomorphic sections. In a usual way one
defines the cohomology groups $\sfh^0(S, \calo_E$ and $\sfh^1(S, \calo_E)$.
Denote $h^i=dim_{\cc }\sfh^i$. By  $g$ we denote the genus of $S$ and
$$
c_1 = \int_Sc_1(E).
$$
These numbers are related by the classical

\state Riemann-Roch Formula. {\it For a holomorphic bundle $E$ over a Riemann
surface $S$ one has

$$
h^0-h^1 =  c_1 + r\cdot (1 - g). \eqno(A1.2.1)
$$
}

\smallskip\rm This formula can be interpreted as the formula for the
index of some operators acting on the spaces of smooth sections of $E$.
On the sheaf of smooth sections of a holomorphic bundle over a Riemann
surface, or more generally over a complex manifold, one can define 
$\dbar $-operators. Those are  $\cc $-linear operators
$\dbar :\Gamma^{1,p}(S,E)\to \Gamma^p_{0,1}(S,E)$ satisfying

$$
\dbar (f\cdot \sigma ) = \dbar_Sf\otimes \sigma + f\cdot \dbar \sigma .
\eqno(A1.2.2)
$$

\smallskip
Here by  $\Gamma^{1,p}(S,E)$ we denote the Sobolev space of $(1,p)$-
smooth sections of $E$, and by $\Gamma^p_{0,1}(S,E)$ the space of $(0,1)$
$L^p$-integrable forms with coefficients in $E$. $\dbar_S$ is a canonical
$\dbar $-operator on $S$.

If one additionally fixes some Hermitian metric on  $E$, then such an
operator is determined uniquely if one imposes the additional condition to
preserve the scalar product, see [G-H] Ch.0 for details.

The operator $\dbar $ being elliptic is Fredholm and its index is defined as
$\ind \dbar := \dim \ker \dbar \allowbreak - \dim \coker \dbar $. Remark that
$\ind \dbar = h^0-h^1$ and so by the Riemann-Roch formula
$$
\ind \dbar = c_1 + r\cdot (1 - g). \eqno(A1.2.3)
$$

\state Definition A1.2.1. {\it An $\rr $-linear operator
$D:\dbar :\Gamma^{1,p}(S,E)\to \Gamma^p_{0,1}(S,E)$ which can be
represented as $\dbar + R$ with $R\in C^0(S, \hom _{\rr }(E,\Lambda^{0,1}
\otimes E))$ shall be called a $\dbar $-type operator.
}

\smallskip\rm This is again an elliptic (Fredholm) operator, which is
homotopic
to $\dbar $. So by the homotopy invariance of the  index we have that for
any $\dbar $-type operator $D$
$$
\ind_\rr D = 2\cdot (c_1 + r\cdot (1 - g)). \eqno(A1.2.4)
$$
The reader should take into account that since $D$ is real, the real
dimensions in the last formula are considered. That is why the number $2$
appears.

\newpage\noindent
{\bigbf Lecture 2}

\smallskip\noindent
{\bigbf Local Existence of Curves}

\medskip\noindent
{\bigsl 2.1. Sobolev Imbeddings, Cauchy-Green Operators, Calderon-Zygmund
Inequality.}

\smallskip\rm
For a natural $k$ and a real $p\ge 1$ the Sobolev space $L^{k,p}
(\Delta ,\cc^n)$
consists of functions $f\in L^p(\Delta ,\cc^n)$ such that their derivatives
up to the order $k$ are also in $L^p(\Delta ,\cc^n)$. One puts
$\Vert f\Vert_{k,p}:= \Sigma_{0\le \vert i\vert \le k}\Vert D^if\Vert_p$.

For $0<\alpha <1$ one considers also the H\"older spaces
$C^{k,\alpha }(\Delta , \cc^n)$. $C^{k,\alpha }(\Delta , \cc^n)$  is the
space of $f\in C^k(\Delta )$ such that
$$
\Vert f\Vert_{k,\alpha }:=
\Vert f(x)\Vert_{C^k}  + \sup_{x\not= y}{\Vert D^kf(x)-D^kf(y)\Vert
\over \vert x-y\vert^{\alpha }}<\infty.
$$
Where $D^kf$ denotes the vector of derivatives of $f$ of order $k$.

One has the following important

\state Sobolev Imbeddings. {\it

(a) For $1\le p<2$ and $1\le q\le {2p\over 2-p}$

$$
L^{k,p}(\Delta , \cc^n) \subset  L^{k-1,q}(\Delta , \cc^n); \eqno(2.1.1)
$$

(b) for $2<p\le \infty $ and $0\le \alpha \le 1-{2\over p}$

$$
L^{k,p}(\Delta , \cc^n) \subset C^{k-1,\alpha }(\Delta , \cc^n).
\eqno(2.1.2)
$$
 Moreover, the imbedding (2.1.1) is a bounded operator and imbedding (2.1.2)
 is  compact.
}
\smallskip
\rm Existence of such imbeddings means, in particular, the existence of 
universal constants $C_q$ and $C_{\alpha }$ s.t. $\Vert f\Vert_{L^{k-1,q}(
\Delta )}\le C_q\cdot \Vert f\Vert_{L^{k,p}(\Delta )}$ and 
$\Vert f\Vert_{C^{k-1,\alpha }(\Delta )}\le C_{\alpha }\cdot 
\Vert f\Vert_{L^{k,p}(\Delta )}$. This is the most frequently used form
of Sobolev Imbedding Theorem.

In the Schwartz spaces $\cals (\cc ,\cc^n)$ and $\cals'(\cc ,\cc^n)$
consider
the operators $\d ={\d \over \d z}$, $\dbar = {\d \over \d \bar z}$ and
$T=T_{CG}={1 \over 2\pi iz}*(\cdot )$, $\bar T={1 \over 2\pi i\bar z}*
(\cdot )$. Note that Cauchy-Green operators
$T$ and $\bar T$ act from $\cals $ only to $\cals'$.
Nevertheless one has the following identities on $\cals $, on $L^p(\cc )$
and on $\cals'$:
$$
\dbar \scirc T=T\scirc \dbar = \id , \eqno(2.1.3)
$$
and
$$
\d \scirc \bar T=\bar T\scirc  \d = \id .
$$

\rm Recall also the

\smallskip
\state Calderon-Zygmund Inequality. {\it For all $1<p<\infty $ there is
a constant $C_p$ such that for all $f\in L^p(\cc ,\cc^n )$  
$$
\Vert (\d \scirc T)(f)\Vert_{L^p}\le C_p\cdot \Vert f\Vert_{L^p}
\eqno(2.1.4)
$$
and
$$
\Vert (\dbar \scirc \bar T)(f)\Vert_{L^p}\le C_p\cdot
\Vert f\Vert_{L^p}.
$$
}

\smallskip
This implies that for $f\in L^p(\cc ,\cc^n )$ and $g=Tf$ (or $g=\bar Tf$) one
has
$g\in L^p_\loc(\cc ,\cc^n )$ and $\Vert dg\Vert_{L^p}(\cc ,\cc^n )\le (1+C_p)
\Vert f\Vert_{L^p(\cc ,\cc^n)}$. For the proof of (2.1.4) we refer to
[Mc-Sa].

Properties (2.1.3) and (2.1.4) also imply that the Cauchy-Green operator $T$ is
a bounded linear operator from $L^{k,p}(\Omega , \cc^{n})$ to
$L^{k+1,p}(\Omega , \cc^{n})$ if $\Omega \comp \cc $. The same is true in
H\"older spaces $C^{k,\alpha }$.

\smallskip We shall repeteadly use 

\state H\"older Inequality. {\it Let $p,q>1$ and let ${1\over p}+{1\over q}=
{1\over r}$. Then for all $f\in L^p(\Delta )$, $g\in L^q(\Delta )$ we have 
that $fg\in L^r(\Delta )$ and 

$$
\Vert fg\Vert_{L^r(\Delta )}\le \Vert f\Vert_{L^p(\Delta )}\cdot 
\Vert g\Vert_{L^q(\Delta )}.\eqno(2.1.5)
$$
}
\smallskip
The behavior of $L^p$-norms under dilatations is also frequently used. 

\state Lemma 2.1.1. {\it Let $h\in L^p(\Delta )$, $\tau >0$ and $\pi_{\tau }:
\Delta\to \Delta $ denotes the contraction $\pi_{\tau }:z\to \tau z$. Put 
$\pi_{\tau }^*h(z)=h(\pi_{\tau }(z))=h(\tau z)$ and $\Delta_{\tau }=
\pi_{\tau }(\Delta )=\{ z\in \cc : \vert z\vert <\tau \} $. Then 

$$
\Vert \pi_{\tau }^*h\Vert _{L^p(\Delta )} = \tau^{-{2\over p }}\cdot 
\Vert h\Vert_{L^p(\Delta_{\tau })}.\eqno(2.1.6)
$$
}
\state Proof. {\rm 
$$
\Vert \pi_{\tau }^*h\Vert _{L^p(\Delta )} = \bigl(\int_{\Delta }\vert h(\tau z)
\vert^pdz\wedge d\bar z\bigr)^{1\over p} =  \bigl(\int_{\Delta_{\tau } }
\vert h(w)\vert^p{1\over \tau^{2}}dw\wedge d\bar w\bigr)^{1\over p} =
\tau^{-{2\over p}}\cdot \Vert h\Vert_{L^p(\Delta_{\tau })},
$$
where $w=\tau z$.
}
\qed  
\smallskip
\smallskip\noindent
{\bigsl 2.2. Local Existence of Curves}

\state Proposition 2.2.1. {\sl (Local existence of curves)} {\it Let $(X,J)$
be
an almost-complex manifold with $J$ of class $C^{k,\alpha }, k\ge 1, \alpha
>0$, and $x_0\in
X$. Then for every $v\in T_{x_0}X$ small enough there exists a $J$-complex
curve $u:(\Delta ,0)\to (X,x_0)$ such that $du(0)({\d \over \d x})=v$.
}

\state Proof. {\rm Take a chart $B\ni x_0$ with $x_0=0$ and $J(0)=J\st$.
Denote by $z=x+iy$ the coordinate in $\Delta $, and by $u_1,...,u_{2n}$ the
coordinates in $B$. The Cauchy-Riemann equation for $u:(\Delta ,J_{\Delta })
\to (B,J)$ has, in our local coordinates, the form
$$
{\d u\over \d y} + J(u)du(J_{\Delta }{\d \over \d y}) = 0,
$$
where $J_{\Delta }$ denotes the canonical a.-c. structure in
$\Delta $, \ie, a multiplication by $i$. Using the fact that 
$J_{\Delta }({\d \over \d y})=-{\d \over
\d x}$, one obtains
$$
{\d u\over \d y} - J(u){\d u\over \d x} = 0.
$$
This is equivalent to
$$
{\d u\over \d \bar z} - Q(u){\d u\over \d z} = 0, \eqno(2.2.1)
$$
where $Q(u)=[J\st + J(u)]^{-1}\scirc [J\st - J(u)]$.

After rescaling, \ie, considering $Q_t(u):=Q(tu)$ and $u_{t,\tau }(z)
:= t^{-1}u(\tau z)$, we can assume that $\Vert Q\Vert_{C^1}$ is
sufficiently small. Note also, that $(2.2.1)$ is equivalent to the
holomorphicity of $h=(\id - T[Q(u){\d \over \d z}])u$.

Consider a $C^1$-mapping
$$
\Phi :]-1,1[ \times C^{1,\alpha }(\Delta ,B)\to C^{1,\alpha }(\Delta ,\cc^n)
$$
given by
$$
\Phi (t,u) = (\id - T[Q(tu){\d \over \d z}])u. \eqno(2.2.2)
$$
Note that $\Phi (0,u) = \id (u) = u$. Thus, the Implicit Function Theorem
tells us that there exists a $t_0>0$, such that for $\vert t\vert <t_0$
the map $\Phi (t,\cdot )$ is a $C^1$-diffeomorphism of the neighborhood
of  zero in $C^{1,\alpha }(\Delta ,B)$ onto the neighborhood of zero $V$
in $C^{1,\alpha }(\Delta , \cc^n)$.

For $w\in \cc^n$ small enough, a holomorphic function $h_w(z):=z\cdot w$
belongs to $V$. Put $u_{t,w}=\Phi (t,\cdot )^{-1}(h_w)$. Then $tu_{t,w}$
is $J$-holomorphic. In fact, $h_w=\Phi (t,\cdot )[{1\over t}u_{t,w}]=
{1\over t}[\id - TQ(u_{t,w}){\d \over \d z}]u_{t,w}$. Moreover $u_{0,w}=
h_w$, so than $du_{0,w}(0)=w$-linear map from $\cc $ to $\cc^n$. This shows
that for $t>0$ small enough $w\to du_{t,w}(0)({\d \over \d x})$ is a
diffeomorphism between the neighborhoods of zero in $\cc^n$.
\qed

\smallskip
We immediately obtain the following

\smallskip\noindent
\state Corollary 2.2.2. {\it Every almost complex structure on a Riemann
surface is complex.}

\medskip\noindent
{\bigsl 2.3. Generalized Calderon-Zygmund Inequality. }

\smallskip\rm
Consider  a continuous linear complex structure $J(z)$ on the trivial
bundle $\cc \times \rr^{2n}\to \cc $, \ie, $J(z)$ is a continuous family
of endomorphisms $\rr^{2n}\to \rr^{2n}$ with $J(z)^2=-Id$. Define an
operator $\dbar_{J}:\cals' (\cc ,\rr^{2n})\to \cals'(\cc ,\rr^{2n})$
by the formula
$$
(\dbar_{J}f)(z)={1\over 2}[\d_xf(z) + J(z)\d_yf(z)].
$$
If $J_0$ is another continuous complex structure on $\cc \times \rr^{2n}$,
then
for $f\in L^p(\cc ,\rr^{2n})$ it holds that
$$
\Vert (\dbar_J\scirc T - \dbar_{J_0}\scirc T)f\Vert_{L^p(\cc )}\le
\Vert J - J_0\Vert_{L^{\infty }(\cc )}\cdot \Vert d(Tf)\Vert_{L^p(\cc )}
\le
$$
$$
\le \Vert J - J_0\Vert_{L^{\infty }(\cc )}(1+C_p)\Vert f\Vert_{L^p(\cc )}.
\eqno(2.3.1)
$$
If we take $J_0(z)\equiv J\st $, the standard structure in $\cc^n$, then
as was remarked above, $\dbar_{J_0}\scirc T: L^p(\cc ,\cc^n)\to L^p(\cc,\cc^n)$
is an identity. So from (2.3.1) we see that there exists 
$\eps_p={1\over 1+C_p}$ 
such that if $\Vert J - J\st \Vert <\eps_p$, then $\dbar_{J}\scirc T: L^p(\cc,
\cc^n)\to L^p(\cc , \cc^n)$ is an isomorphism. Moreover, since 
$\dbar_J\scirc T = \dbar_{J\st}\scirc T + (\dbar_J-\dbar_{J\st})\scirc T$, 
we have 
$$
(\dbar_J\scirc T)\inv = (\id + (\dbar_J-\dbar_{J\st})\scirc T)^{-1} =
\Sigma_{n=0}^{\infty }(-1)^n[(\dbar_J -
\dbar_{J\st})\scirc T]^n.\eqno(2.3.2)
$$
This shows, in particular, that $(\dbar_J\scirc T)\inv $ does not depend
on $p>1$. Now we shall prove the following statement, which can be viewed
as a generalization of the Calderon-Zygmund estimate.
}

\smallskip
\state Lemma 2.3.1. {\it For any $u\in L^{1,2}(\cc ,\rr^{2n})$ with compact
support, any
continuous $J$ with $\Vert J-J\st \Vert_{L^{\infty }(\cc ,End(\rr^{2n}))}<
\eps_p$ the condition $\dbar_Ju\in L^p(\cc ,\rr^{2n})$ implies
$$
\Vert du \Vert_{L^p(\cc )}\le C\cdot \Vert \dbar_Ju\Vert_{L^p(\cc )}
\eqno(2.3.3)
$$
with some $C=C(p,\Vert J-J\st \Vert_{L^{\infty }(\cc )})$.
}

\state Proof. \rm Put $v=u-T\scirc \dbar_{J\st}u $. Then $\dbar_{J\st}v=0$.
So $v$ is holomorphic and descends at infinity. Thus $v=0$, which implies
$u=(T\scirc \dbar_{J\st})u$. By the Calderon-Zygmund inequality, to estimate
$\Vert du\Vert_{L^p(\cc )}$ it is sufficient to estimate $\Vert \dbar_{J\st}u
\Vert_{L^p(\cc )}$. We have  $(\dbar_J\scirc T)\scirc \dbar_{J\st}u =
\dbar_Ju
\subset L^p(\cc )\cap L^2(\cc )$. From (2.3.2) we obtain  $\dbar_{J\st}u\in
L^p(\cc )\cap L^2(\cc )$ with an estimate
$$
\Vert \dbar_{J\st }u\Vert_{L^p(\cc )}\le \sum_{n=0}^{\infty }\Vert
(\dbar_J-\dbar_{J\st})\scirc T\Vert_p^n\cdot \Vert \dbar_Ju\Vert_{L^p(\cc )}
\le
C\cdot \Vert \dbar_Ju\Vert_{L^p(\cc )},
$$
where $C=C(p,\Vert J-J\st\Vert_{L^{\infty }})=\sum_{n=1}^{\infty } 
\left(\vert J-J\st\vert_{L^{\infty }}(1+C_p)\right)^n$, provided that
$\Vert J-J\st\Vert_{L^{\infty }}<\eps_p={1\over 1+C_p}$.

\smallskip
This yields (2.3.3).
\qed

\smallskip
\state Corollary 2.3.2. {\it If $u\in (L^{1,2}\cap C^0)(\Delta ,\rr^n)$ is 
$J$-holomorphic for a continuous a.-c. structure $J$ belongs to $L^{1,p}_\loc
(\Delta ,\rr^2)$ for all $2<p<\infty $. 
}

\bigskip\noindent
{\bigsl 2.4. First A priori Estimate.}

\medskip\noindent
\rm
Let $(X,J_0)$ be an almost complex manifold, and the tensor $J_0$ 
is supposed to be
of
class $C^0$ (\ie continuous) only. Let $K\comp X$ be some compact
in $X$. We fix
some Riemannian metric $h$ on $X$ and denote
by $\mu (J_0)$ the module of continuity of $J_0$. All norms and distances
are taken with respect to $h$.

Recall that for a map $u:X\to Y$ into a metric space an oscillation of $u$
on a subset $D\subset X$ is defined as $\osc (D,u):=\sup \{ d_Y(u(x),u(y)):
x,y\in D\} $.

\smallskip
\state Lemma 2.4.1. {\it (First a priori estimate).  For every $2< p<\infty $
there exists an $\eps_1 =\eps_1(p,\mu (J_0),K,h)$ and
$C_p=C(p,\mu (J_0),K,h)$ such that for any
$J\in C^0$, $\vert\vert J - J_0\vert\vert_{C^0(K)}<\eps_1 $ and
every $J$-holomorphic map $u\in C^0\cap L^{1,2}(\Delta ,X)$ with
$u(\Delta )\subset K$, satisfying $\Vert du\Vert_{L^2(
\Delta )}<\eps_1 $ the following holds
$$
\Vert du\Vert_{L^p({1\over 2}\Delta )}\le C_p\cdot
\Vert du\Vert_{L^2(\Delta )}.\eqno(2.4.1)
$$

}

\smallskip
\state Proof. {\rm This will be done in several steps.

\noindent
\sl Step 1. We prove first the inequality (2.4.1)
for the case when $K\comp U\subset \cc^n$, $h$ Euclidean metric,
$J_0$ the standard complex structure in $\cc^n=\rr^{2n}$, and $\Vert J-
J_0\Vert_{L^{\infty }}< \eps_p$ for $\eps_p$ from {\sl Lemma 2.3.1}.

\smallskip
\rm Note that in this case the condition that $\Vert du\Vert_{L^2(\Delta )}$ 
should be small is not needed.

To prove the first step consider a $J$-holomorphic map
$u:\Delta \to (\rr^{2n},J)$, $u(\Delta )\subset K$ and $\Vert J-J\st \Vert
<\eps_p$. Define on $\Delta \times \rr^{2n}$ a linear complex structure
$J(z)=(u^*J)(z)$. Then $u$ defines a $J$-holomorphic section of
$(\Delta \times \rr^{2n}, J)$ with $\Vert J-J\st \Vert_{L^{\infty } (\Delta )}
< \eps_p$. Extend $J$ onto $\cc \times \rr^{2n}$ with the same estimate.

Let $\psi $ be a non-negative cut-off function supported in $\Delta_{3/4}$
and equal to one on $\Delta_{1/2}$. Put $u_1=u\psi $. Then $u_1\in 
L^{1,2}(\Delta)$ (because $u\in L^{1,2}(\Delta)$) and $\dbar_Ju_1=u\dbar_J\psi
\in L^p(\cc )$ with $\Vert \dbar_Ju_1\Vert_{L^p(\Delta )}=\Vert u\dbar_J\psi 
\Vert_{L^p(\Delta )} \le C\Vert du\Vert_{L^2(\Delta )}$ for any $p$ by the 
Sobolev imbedding $L^{1,2}(\Delta ,\cc )\to L^p(\Delta ,\cc )$ for all 
$p<\infty $. Now the
generalized Calderon-Zygmund estimate (2.3.3) applies in obtaining 
an estimate of \sl Step 1.

\rm Using the Sobolev imbedding $L^{1,p}\subset C^{1-{2\over p}}$ and obvious
properties of $L^p$-norms with respect to dilatations, one 
easily derives from {\sl Step 1} the following

\medskip\noindent\sl
Step 2. \sl Fix $2<p<\infty $. There exists $\eps_2=\eps_2(p,\mu (J_0),K,h)>
0$ such that for every $J$ with $\Vert J-J_0\Vert_{L^{\infty }}<\eps_p$
and any $J$-holomorphic $u$ with
$\osc(u,\Delta (x,r))\le \eps_2$ on any $\Delta (x,r)\subset \Delta$
it holds that
$$
\osc(u,\Delta (x,{r\over 2}))\le  Cr^{{2\over p}-1}\cdot
\vert \vert du\vert \vert_{L^2(\Delta (x,r))},\eqno(2.4.2)
$$

\noindent
and 
$$
\Vert du\Vert_{L^p(\Delta (x,r/2))}\le C_p\cdot 
\Vert du\Vert_{L^2(\Delta (x,r))}.
$$

\smallskip\noindent\rm
In fact,
$$
\osc(u,\Delta (x,{r\over 2}))\le C_1\vert \vert du\vert \vert_
{L^p(\Delta (x,r/2))}=C(\int_{\Delta (0,{r\over 2})}\Vert du\Vert^pdxdy)^{1
\over p}=
$$
$$
=C(\int_{\Delta (0,{1\over 2})}\Vert du(rw){1\over r}\Vert^p)^{1\over p}=
d(rw)d\overline{(rw)})^{1\over p}=r^{2-p\over p}\Vert du_r\Vert_{L^p(\Delta
(0,{1\over 2}))}
$$
for $u_r(w)=u(rw)$.

\smallskip Furthermore, if $\eps_2>0$ is such that $\Vert J(x) - J(y)\Vert 
<\eps_p$ for $h(x,y)<\eps_2$, $x,y\in K$, then applying Step 1 to $u_r$ we 
get 

$$
r^{2-p\over p}\Vert du_r\Vert_{L^p(\Delta (0,{1\over 2}))}\le Cr^{2-p\over p}
\Vert du_r\Vert_{L^2(\Delta (0,1))} = 
Cr^{2-p\over p}\Vert du\Vert_{L^2(\Delta (0,r))}.
$$

Let $\alpha \ge 0$ be the following continuous function: $\alpha \equiv 1$ for
$t\le {1\over 2}$ and $\alpha (t)\equiv 0$ for $t\ge 3/4$. On the interval
${1\over 2}\le t\le {3\over 4}, \alpha (t)=3-4t$. Put for $x\in \Delta $

$$
f(x):=\max \{ t: 0\le t \le {1\over 8}, \osc(u,\bar\Delta(x,t\cdot \alpha
(\vert x\vert )))\le \eps_2\}.
$$

\noindent
Clearly $f$ is continuous and $f\equiv {1\over 8}$ for $ {3\over 4}\le
\vert x \vert <1$.

\smallskip\noindent\sl
Step 3. \it There exists an $\eps_1=\eps_1(p,J_0,K,h)>0$ such that
$f(x)\equiv {1\over 8}$.

\smallskip\rm
Suppose that there is an $x_0$ with $f(x_0)=\min \{ f(x):x\in \Delta \}
< {1\over 8} $. It is clear that $f(x_0)>0$.

Take the disc $\Delta (x_0,a), a:=f(x_0)\alpha (\vert x_0\vert )$. Note that
$$
\osc(u,\Delta (x_0,a))=\eps_2\eqno(2.4.3)
$$
Condition (2.4.2) together with the Sobolev
embedding $L^{1,4}(\Delta )\subset C^{0,{1\over 2}}(\Delta )$ tells us
(because
$\osc(u,\Delta (x_0,a))=\eps_2$) that $\osc(u,\Delta (x_0,{a\over 2}))$
$\le C\cdot \vert \vert du\vert \vert_{L^2(\Delta (x_0,a))}$. Take a point
$x_1 \in \Delta (x_0,a)$ on the distance not more then ${3\over 4}a$ from
$x_0$. We have from
$f(x_0)={a\over \alpha (\vert x_0\vert )}$ that 
$f(x_1)\ge {a\over \alpha (\vert x_0\vert )}$ and
thus $f(x_1)\alpha (\vert x_0\vert )\ge a$. At the same time $\alpha (\vert
x_1\vert )\ge \alpha (\vert x_0\vert ) - 3a$, so $f(x_1)\alpha (\vert x_1
\vert )\ge a-3a\cdot f(x_1)\ge {a\over 2}$, because $f(x_1)\le {1\over 8}$.
That means that $\osc(u,\Delta(x_1,{a\over 2}))\le \eps_2$ and thus
$\osc(u,\Delta(x_1,{a\over 4}))\le C\cdot \vert \vert du\vert \vert_{L^2
(\Delta (x_1,{a\over 2}))}$. Thus,
$\osc(u,\Delta(x_0,a))\le 4C\cdot \vert \vert du\vert \vert_{L^2(\Delta )}$.
If $\eps_1 $ were taken smaller than ${\eps_2\over 4C}$, then we would obtain a
contradiction with (2.4.3). Step 3 is proved.

This means that $\osc(u,\Delta(x,{1\over 8}))\le \eps_2$ for any $x\in
{1\over 2}\Delta $. Therefore, Step 2 with $r={1\over 8}$ gives us 
the conclusion
of the Lemma. \qed

This statement can be used to prove once more the following statement, which 
was already proved in Corollary 2.3.2:

\smallskip
\state Corollary 2.4.2. {\it A $J$-holomorphic map $u:\Delta
\to (X,J)$ is $L^{1,p}$-continuous for any $p<\infty $, provided $J$
is of class $C^0$.
}

\state Proof. Note that $u_{\eps }(z):=u(\eps z)$ is
also $J$-holomorphic and $\Vert du_{\eps }\Vert_{L^2(\Delta (0,1))} =
\Vert du\Vert_{L^2(\Delta (0,\eps ))}$. For $\eps $ small enough we
shall have
$\Vert du_{\eps }\Vert_{L^2(\Delta (0,1))} =
\Vert du\Vert_{L^2(\Delta (0,\eps ))}<\eps_1 $ - from Lemma 2.4.1. Now
estimate  (2.4.1) gives us $L^{1,p}$-continuity of $u_{\eps }$ and thus of
$u$ in the neighborhood of zero.

\smallskip\noindent\bigsl
2.5. Convergence outside of a finite set of points.

\smallskip\rm
Another immediate consequence of the a priori estimate (2.4.1) is the
following

\smallskip
\state Corollary 2.5.1. {\it Let $\{ J_n\}$ be a sequence of continuous
almost complex structures on $X$ such that $J_n\to J$ in $C^0$-topology
on $X$. Let $u_n\in C^0\cap L^{1,2}_{loc }(\Delta ,X)$ be a sequence of
$J_n$-holomorphic maps
such that $\Vert du_n\Vert_{L^2(\Delta )}\le \eps_1$. Then there exists a
subsequence $u_{n_k}$ which $L^{1,p}_\loc$-converge to a $J$-holomorphic
map $u_\infty$ for all $p\ge 2$.}

\state Proof. The main estimate together with the Sobolev imbedding $L^{1,p}
\subset C^{1-{2\over p}}$ gives us a subsequence $\{ u_{n_k}\} $ which
converges to $u_\infty$ in $C^{\alpha }(\Delta )$ for any $\alpha
<1$. Take $\phi \in C_0^{\infty }(\Delta )$, $\phi\mid_{\Delta (r)}
\equiv 1$ and consider $\dbar_{J_{n_k}}(\phi u_{n_k})=\d_x(\phi u_{n_k})
+ J_{n_k}(u_{n_k})\d_y(\phi u_{n_k})=$ $(\d_x\phi +J_{n_k}(u_{n_k})
\d_y\phi)u_{n_k}$, which is $C^0$ and thus $L^p$-convergent for any
$p\ge 2$. The a priori estimate for $\dbar_J$-operator (\ie Calderon-Zygmund
inequality 2.3.1) gives us $L^{1,p}_\loc$-convergency of $\{ u_{n_k}\} $.
\qed

\smallskip This corollary implies another one:

\smallskip
\state Corollary 2.5.2. {\it Let $\{ J_n\}$ be a sequence of continuous
almost-complex structures on $X$ such that $J_n\to J\in C^0$ in $C^0$-topology
on $X$. Let $u_n\in C^0\cap L^{1,2}_{loc }(\Delta ,X)$ be a sequence of
$J_n$- holomorphic maps
such that $\Vert du_n\Vert_{L^2(\Delta )}\le C$. Then there exists a
subsequence $\{ u_{n_k}\} $ and a finite set of points $\{ x_1,...,x_l\} $ in
$\Delta $, such that $\{ u_{n_k}\} $  is $L^{1,p}_\loc$-convergent on
compacts in $\Delta \setminus \{ x_1,...,x_l\} $ for all $p\ge 2$. 
Moreover $l\le {3C\over \eps_1}$.
}

\state Proof. \rm For any $n\in \nn $, cover our disk $\Delta $ by a finite
number of disks of radii ${1\over n}$ in such a way that no point in
$\Delta $ belongs at the same time to more than  $3$ of them. Then there is
at most $[{3C\over \eps_1}$ disks, on which the energy of our maps
(after going to a subsequence) is more then $\eps_1$. Thus, on the complement
to those disks the {\sl Corollary 2.5.1} applies.

Taking $n$ bigger and bigger (and passing to a subsequence) we obtain the
desired result.
\qed

\smallskip One can prove also the following regularity statements about 
$j$-holomorphic maps:

\smallskip
\state Lemma 2.5.3. {\it a) If the structure $J$ is Lipschitz then 
$J$-holomorphic maps are of class $C^1$;

\noindent
2) If $J\in C^k$ then $J$-holomorphic maps are in $C^{k,\alpha }$ for 
all $\alpha <1$.
}

\newpage
\noindent{\bigbf Lecture 3}

\smallskip\noindent
{\bigbf Positivity of Intersections of Complex Curves.}

\smallskip\noindent \sl 3.1. Unique Continuation Lemma.

\nobreak\smallskip\rm

We start with a unique continuation-type lemma for $\dbar$-unequalities, 
namely with {\sl Lemma 3.1.1} (compare with [Ar] and [Hr-W]). In the proof we
shall use the following special version
of the theorem of Harvey and Polking [Ha-Po]: 

\smallskip
\state Theorem of Harvey-Polking. {\sl Let $f:\Delta
\to \cc^n$ be locally $L^2$-integrable.  Assume that for some $g \in
L^1_\loc (\Delta, \cc^n)$  the equation $\dbar f =g$ holds (in the weak
sence) in the punctured disc $\check\Delta$. Then $\dbar f =g$ holds in
the whole disc $\Delta$.
}

\medskip
\state Lemma 3.1.1. {\it Suppose that the function $u\in L^2_\loc(\Delta,
\cc^n)$ is not identically $0$, $\dbar u\in L^1_\loc(\Delta,\cc^n)$
and satisfies {\sl a.e.} the inequality
$$
\vert \dbar u \vert \le h\cdot \vert u \vert\eqno(3.1.1)
$$
for some nonnegative $h\in L^p_\loc(\Delta)$ with $2<p<\infty$. Then

\sli $u\in L^{1, p}_\loc(\Delta)$, in particular $u\in C^{0,\alpha}_\loc
(\Delta)$ with $\alpha\deff1-{2\over p}$;

\slii for any $z_0\in\Delta$ such that $u(z_0)=0$ there exists
$\mu\in\nn$---the~multiplicity of zero of $u$ in $z_0$---such that
$u(z)=(z-z_0)^\mu \cdot g(z)$ for some $g\in L^{1, p}_\loc (\Delta)$ with
$g(z_0)\not=0$.
}

\bigskip
\state Proof. Statement \sli is easily obtained by increasing the 
smoothness argument. Let us for the methodological reasons give it in 
full details.  

\smallskip\noindent\sl
Step 0. $u$ is in $L^{p_0}_\loc$ for $p_0:=p>2$. 

\smallskip \rm From (3.1.1) and the H\"older inequality we see that 
$h\vert u\vert \in L_\loc^{{1\over {1\over 2}+{1\over p}}}=
L^{{2p\over p+2}}_\loc$. Therefore one 
obtains $\bar\partial u\in L^{2p\over p+2}_\loc(\Delta)$. Consequently, 
$u\in L^{1, {2p\over p+2}}_\loc(\Delta)$ due to ellipticity of $\dbar$ and 
because ${2p\over p+2} >1$. By the Sobolev imbedding we have  that 
$u\in L_\loc^{{2{2p\over p+2}\over 2-{2p\over p+2}}}=L^p_\loc(\Delta)$.

Thus we proved that $u\in L^{p_0}_\loc$ for $p_0:=p>2$.

\smallskip\noindent\sl
Step 1. $u$ is in $L^{p_1}_\loc$ for $p_1:={2p\over 4-p}>p_0$.

\smallskip\rm Again, $\vert \dbar u\vert \le h\vert u\vert\in 
L_\loc^{{1\over {1\over p}+{1\over p}}}=L_\loc^{{p\over 2}}$. Therefore 
$u\in L_\loc^{1,{p\over 2}}\subset L_\loc^{{p\over 2-{p\over 2}}}= 
L_\loc^{{2p\over 4-p}}$.

\smallskip\noindent\sl
Step n. There exists $r>1$ such that $u\in L_\loc^{p_n}$, where 
$p_n\ge rp_{n-1}$ for all $n$.

\smallskip\rm $u\in L^{p_{n-1}}$ implies $\vert \dbar u\vert \le h\vert u\vert 
\in L^{{1\over {1\over p}+{1\over p_{n-1}}}}=L^{{pp_{n-1}\over p+p_{n-1}}}$. 
Therefore $u\in  L^{1,{pp_{n-1}\over p+p_{n-1}}}\subset 
L^{{2pp_{n-1}\over 2p+2p_{n-1}-pp_{n-1}}}=L^{p_n}$ with $p_n={1\over {1\over p}
+{1\over p_{n-1}}-{1\over 2}}$. One easily checks that ${p_n\over p_{n-1}}\ge 
r>1$, which does not depend on $n$.

\smallskip\noindent\sl
Step $\infty $. $u$ is in $C^{0,\alpha }$ with $\alpha =1-{2\over p}$.

\smallskip\rm When $p_n\to \infty $ $\vert \dbar u\vert \le h\vert u\vert \in 
L^{{1\over {1\over p}+{1\over p_n}}}=L^q$, where $q\sim p$, because ${1\over 
p_n}\sim 0$. Therefore $u\in L^{1,q}$ with $q>2$ close to $p$. Sobolev 
Imbedding tells in this case that $L^{1,q}\subset C^{0,\beta }$ with $\beta 
=1-{1\over q}$. Finally we see that $\dbar u\in L^p$, therefore $u\in C^{0,
\alpha }$.

\smallskip
\slii  Now suppose that $u(z_0)=0$. Then, due to the~H\"older continuity, we
have $\vert u(z) \vert \le C\vert z-z_0 \vert^\alpha$ for $z$ close enough to
$z_0$ and consequently $u_1(z) \deff u(z) / (z-z_0)$ is from $L^2_\loc(\Delta)
$. The~theorem of Harvey-Polking provides that $\dbar u_1\in
L^1_\loc(\Delta)$ and $u_1$ also satisfies
inequality (3.1.1). In particular, $u_1$ is also continuous. Iteration of this
procedure gives the possibility of defining the~multiplicity of zero of $u$ in
$z_0$ provided we show that after a finite number of steps we obtain the
function $u_N$ with $u_N(z_0)\not= 0$. To do this we may assume that $z_0=0$.
Let $\pi_\tau (z)\deff\tau\cdot z$ for $0<\tau<1$. Then $u_\tau\deff
\pi_\tau^*(u)$ satisfies the inequality $\vert \dbar u_\tau \vert \le
\tau\pi_\tau^*h \cdot \vert u_\tau \vert$. Since
$$
\Vert \tau\pi_\tau^*h \Vert_{L^p(\Delta)} =
\tau^{1-2/p}\cdot \Vert h \Vert_{L^p(\pi_\tau(\Delta)),}
$$
we can also assume that $\Vert h \Vert_{L^p(\Delta)}$ is small enough.
Fix a~cut-off function $\varphi\in C^{\infty }_0(\Delta)$ which is identically 
$1$
in ${1\over2}\Delta$, the~disk of radius $1\over2$. Then
$$
\Vert \varphi u \Vert_{L^p(\Delta)}\le
C_1\cdot \Vert \varphi u \Vert_{L^{1, {2p\over2+p}}(\Delta)}\le
C_2\cdot \Vert \dbar(\varphi u) \Vert_{L^{2p\over2+p}(\Delta)}\le
$$
$$
\le C_2\cdot (\Vert \varphi \dbar u \Vert_{L^{2p\over2+p}(\Delta)}+
\Vert \dbar\varphi  u \Vert_{L^{2p\over2+p}(\Delta)})\le
%%$$ $$ \le
C_2\cdot ( \Vert \varphi h u  \Vert_{L^{2p\over2+p}(\Delta)}
+  \Vert \dbar\varphi u \Vert_{L^{2p\over2+p}(\Delta)} \le
$$
$$
\le  C_3(\Vert \varphi u \Vert_{L^p(\Delta)}\cdot \Vert h \Vert _
{L^p(\Delta )} + \Vert u\Vert _{L^p(\Delta \setminus {1\over2}\Delta )}
\cdot \Vert \dbar\varphi \Vert _{L^2(\Delta \setminus {1\over2}\Delta )}).
$$

\smallskip
Here we used the~fact that
the~support of $\dbar\varphi$ lies in $\Delta\bss{1\over2}\Delta$. Since
$\Vert h  \Vert_{L^p(\Delta)}$ is small enough we obtain the~estimate
$$
\Vert u \Vert_{L^p({1\over3}\Delta)}
\le C\cdot \Vert u \Vert_{L^p(\Delta\setminus {1\over2}\Delta)}
$$
with the~constant $C$ independent of $u$. Thus if the~multiplicity of
zero of $u$ in $z_0=0$ is at least $\mu$, then
$$
\Vert z^{-\mu}u \Vert_{L^p({1\over3}\Delta)}
\le C\cdot \Vert z^{-\mu}u \Vert_{L^p(\Delta\bss{1\over2}\Delta))}
$$
which easily gives
$$
\Vert u \Vert_{L^p({1\over3}\Delta)} \le C
\left(\hbox{$2\over3$}\right)^{-\mu}
\cdot \Vert u \Vert_{L^p(\Delta\bss{1\over2}\Delta))}.
$$
Now one can easily see that either $u$ has isolated zeros of finite
multiplicity, or $u$ is identically zero. Yet the~last case is excluded by
the hypothesis of the~lemma.\qed

\medskip
\state Lemma 3.1.2. {\it
Under the~hypothesis of {\sl Lemma 3.1.1} suppose
additionally that $u$ satisfies {\sl a.e.} the inequality
$$
\vert \dbar u(z) \vert \le \vert z-z_0 \vert^\nu
h(z)\cdot \vert u(z) \vert,\eqno(3.1.2)
$$
with $z_0\in \Delta$, $\nu\in\nn$, and $h\in L^p_\loc(\Delta)$, $2<p<\infty$.
Then
$$
u(z)=(z-z_0)^\mu\bigl(P^{(\nu)}(z) + (z-z_0)^\nu g(z)\bigr),
\eqno(3.1.3)
$$
where $\mu\in\nn$ is the multiplicity of $u$ in $z_0$, defined above,
$P^{(\nu)}$ is a polynomial in $z$ of degree $\le \nu$ with $P^{(\nu)}(z_0)
\not=0$, $g\in L^{1, p}_\loc(\Delta,\cc^n)\hookrightarrow C^{0,\alpha}$,
$\alpha=1-{2\over p}$, and $g(z)=O(\vert z-z_0\vert^\alpha)$.
}

\medskip\noindent
\bf Proof. \rm The proof uses the~same idea as in the~previous lemma. 
Define $u_0(z)
\deff \msmall{u(z)\over (z-z_0)^\mu}$ and $h_1(z) \deff h(z)\cdot \vert u_0(z)
\vert$. Due to {\sl Lemma 3.1.1}, $u_0\in C^{0,\alpha}$, $u_0(z_0)\not=0$,
$h_1\in L^p_\loc$, and $u_0$ satisfies {\sl a.e.} the~inequality $\vert \dbar
u_0(z) \vert \le \vert z-z_0 \vert^\nu h_1(z)$.

Set $a_0=u_0(z_0)$. Since $u_0(z)-a_0 = O(\vert z-z_0 \vert^\alpha)$,
we have $u_1\deff\msmall{u_0(z) - a_0 \over z-z_0}\in L^2_\loc$. Applying
the~theorem of Harvey-Polking once more, we obtain  $\vert \dbar u_1 \vert
\le \vert z-z_0 \vert^{\nu-1} h_1$, and consequently $u_1\in C^{0,\alpha}$,
$u_1(z)-u_1(z_0)=O(\vert z-z_0 \vert^\alpha)$.
Repeating this procedure $\nu$ times, we obtain the~polynomial
$$
P^{(\nu)}(z)=a_0 + (z-z_0)a_1 + \cdots + (z-z_0)^\nu a_\nu
$$
with
$$
a_k\deff\lim_{z\to z_0}
\msmall{ u(z)-\sum_{i=0}^{k-1} (z-z_0)^i a_i \over (z-z_0)^k },
\qquad 0\le k\le\nu,
$$
and the~function
$$
g(z)\deff \msmall{{u(z)-P^{(\nu)}(z) \over (z-z_0)^\nu }},
$$
which satisfies the~conclusion of the~lemma.
\qed

\medskip

\state Corollary 3.1.3. {\it Let $J$ be a Lipschitz-continuous almost complex
structure in a neighborhood $U$ of $0\in \cc^n$ such that $J(0)=J\st$,
the~standard complex structure in $\cc^n$. Suppose that $u:\Delta\to U$
is a $J$-holomorphic $C^1$-map with $u(0)=0$. Then there exist uniquely
defined $\mu\in\nn$ and a (holomorphic) polynomial $P^{(\mu-1)}$ of degree
$\le\mu-1$ such that $u(z)=z^\mu\cdot P^{(\mu-1)} +z^{2\mu -1}v(z)$
with $v(0)=0$ and $v\in L^{1,p}(\Delta,\cc^n)$ for any $p<\infty$.
}

\medskip
\state Proof. {\rm In fact, by the hypothesis of the lemma $du + J(u)
\scirc du \scirc J_\Delta =0$ and hence
$$
\vert \dbar u \vert  =
\bigl\vert \msmall{1\over2}( du + J\st \scirc du \scirc J_\Delta)\bigr\vert =
\bigl\vert \msmall{1\over2} (J\st - J(u))\scirc du \scirc J_\Delta \bigr\vert
\le \Vert J \Vert_{C^{0,1}(B)} \cdot \vert u \vert  \cdot
\vert du \vert.
\eqno(3.1.4)
$$

\smallskip\noindent
Thus by {\sl Lemma 3.1.1}, $u(z) = z^\mu w(z)$ with some $\mu \in\nn$
and $w\in L^{1,p}_\loc(\Delta, \cc^n)$ for any $p<\infty$.
If $\mu =1$ we are done. Otherwise $ du(z)/z^{\mu-1} =
\mu\,w\,dz + z\,dw \in L^p_\loc(\Delta, \cc^n)$ for any $p<\infty$.
Thus, the corollary follows now from {\sl Lemma 3.1.2}.
}
\qed

\bigskip
\state Corollary 3.1.4. {\it Let $u:S\to (X, J)$ be a~$J$-holomorphic map.
Then for any $p\in X$ the~set $u\inv(p)$ is discrete in $S$, provided that
$J$ is Lipschitz-continuous.
}

\smallskip
\state Proof. {\rm Take an~{\sl integrable} complex structure $J_1$ in some
neighborhood $U$ of $p\in X$, such that $J(p)=J_1(p)$. Take also
$J_1$-holomorphic coordinates $w_1,\ldots, w_n$ in $U$, such that $w_i(p)=0$.
Then the~statement follows easily from {\sl Corollary 3.1.3.
}
\qed

\bigskip
\noindent \sl 3.2. Inversion of $\dbar_{J}+R$ and $L^{k,p}$-topologies.

\smallskip\nobreak
\rm
Now suppose that $J\in C^0(\Delta, \endo_\rr(\rr^{2n}))$ satisfies
the~identity $J^2\equiv -1$,\ie, $J$ is a continuous almost complex structure
in
the trivial $\rr^{2n}$-bundle over $\Delta$. Define the~$\rr$-linear
differential
first order operator $\dbarj:L^{1, p}(\Delta,\rr^{2n}) \to L^p(\Delta,\rr^{2n})$
by setting
$$
\dbarj(f)=\msmall{1\over2}\left(
\msmall{\d f \over \d x} + J \msmall{\d f\over \d y} \right).\eqno(3.2.1)
$$
For example, for $J\equiv J\st$, the~standard complex structure in $\rr^{2n}=
\cc^n$, the~operator $\dbarj$ is a usual Cauchy-Riemann operator $\dbar$.
The~operator $\dbarj$ is elliptic and possesses nice regularity properties in
Sobolev spaces $L^{k, p}$ with $1<p<\infty$ and H\"older spaces $C^{k,\alpha}$
with $0<\alpha<1$. The~following two statements are typical for (nonlinear)
elliptic PDE and produces a result which we need for the purpose of this paper.

\smallskip
\state Lemma 3.2.1. {\it Let $J$ be $C^k$-continuous with $k\ge0$ and $\dbarj$
be as above. Also let $R$ be an~$\endo (\rr^{2n})$-valued function in
$\Delta$ of class $L^{k,p}$ with $1<p<\infty$. If $k=0$, we also assume that
$p>2$. Suppose that for $\dbarj f +Rf \in L^{k,p}$ for some $f\in L^{1,
1}(\Delta,\rr^{2n})$. Then $f\in L^{k+1,p}_\loc (\Delta,\rr^{2n})$ and
for $r<1$
$$
i) \hbox{ } \Vert f \Vert_{L^{k+1, p}(\Delta(r))}
\le C_1(J, \Vert R\Vert_{L^{k,p}}, k, p, r)
\bigl( \Vert \dbarj f + Rf \Vert_{L^{k, p}(\Delta)} +
\Vert f \Vert_{L^1(\Delta)}\bigr)
%\leqno{\sli }}
$$

\smallskip\noindent
If, in addition, $J$ and $R$ are $C^{k,\alpha}$-smooth with $0<\alpha<1$,
then $f\in C^{k+1,\alpha}_{\loc }(\Delta,\rr^{2n})$ and
$$
ii) \hbox{ } \Vert f \Vert_{C^{k+1,\alpha}(\Delta (r))}
\le C_2( J , \Vert R \Vert_{C^{k,\alpha}},k, \alpha, r)\cdot
\bigl( \Vert \dbarj f + Rf \Vert_{C^{k,\alpha}(\Delta)}+
\Vert f \Vert_{L^1(\Delta).}\bigr)
%\leqno{\slii  }}
$$
If, additionally,  
$\norm{ J-J\st }_{C^k(\Delta)} + \norm{R }\Vert_{L^{k,p}(\Delta)}$ 
(resp.,\ $\Vert J-J\st \Vert_{C^{k,\alpha}(\Delta)}+ 
\Vert R \Vert_{C^{k,\alpha}(\Delta)}$) 
is small enough, then there exists a~linear bounded operator
$T_{J, R}:L^{k,p}(\Delta,\rr^{2n}) \to L^{k+1, p}(\Delta,\rr^{2n})$
$($resp.,\ $T_{J, R}:C^{k, \alpha}(\Delta,\rr^{2n}) \longrightarrow
C^{k+1,\alpha}(\Delta,\rr^{2n})$\thinspace{\rm)}
such that $(\dbarj + R)\scirc T_{J, R}\equiv \id$ and $T_{J,
R}(f)\vert_{z=0}=0$.
}

\medskip
\state Proof. The~estimates \sli and \slii of the~theorem are obtained from
the~ellipticity of $\dbarj$ by standard methods. See [Mo] for the~case of
general elliptic systems or [Sk] for the~direct proof. We note only that
for $J$ being only continuous the~constant $C_1$ in \sli crucially depends on
the modulus of continuity of $J$, whereas in all the~other cases of the~theorem
the~dependence is essentially only on the~corresponding norm of $J$.

\smallskip
Let $T$ be a composition ${\d\over\d z} \scirc G$ where $G(f)$ is
the~solution of the~Poisson equation $\Delta u=f$ with the~Dirichlet boundary
condition on $\d\Delta$. Then for $J\equiv J\st$ and $R\equiv0$ we can define
$T_{J\st,0}(f) \deff T(f)- T(f)(0)$. In the general case we set
$$
T_{J, R}=\sum_{n=0}^\infty (-1)^nT_{J\st, 0}\scirc
\bigl((\dbarj-\dbar_{J\st}+R)\scirc T_{J\st, 0}\bigr)^n.
$$
Due to the Sobolev imbedding theorem, the~series converges in an~appropriated
norm, provided $\Vert J-J\st \Vert + \Vert R \Vert$ is small enough.
\qed

\smallskip
\state Remark. {\rm The last statement of this lemma Implies that for 
$\Vert J-J\st \Vert + \Vert R\Vert $ small enough $L^{k+1,p}(\Delta ,\rr^{2n})
=H_{\dbar_J+R}\oplus T_{J,R}(L^{k,p}(\Delta ,\rr^{2n}))$ where  
$H_{\dbar_J+R}=\ker (\dbar_J+R)$. This will be used later in the proof of 
Lemma 6.3.4.

\bigskip
\state Corollary 3.2.2. {\it Let $k\in \nn$, $q>2$, and $J$
a~$C^k$-smooth almost complex
structure in $X$. Also let $(S, J_S)$ be a~complex curve. Suppose that
$L^{1,q}$-map $u:S\to X$ satisfies the~equation
$$
du + J \scirc du \scirc J_S = 0.
$$
Then $u$ is $L^{k+1,p}$-smooth for any $p<\infty$. If, in addition, $J$ is
$C^{k,\alpha}$-smooth with $0<\alpha<1$, then $u$ is $C^{k+1,\alpha}$-smooth.

Let $J^{(n)}$ (resp.\ $J^{(n)}_S$) be a sequence of almost complex
structures on $X$ (resp.\ on $S$), which $C^k$-converges to $J$
(resp.\ to $J^{(n)}_S$) and let $u_n:S \to X$ be a sequence of
$(J^{(n)}_S, J^{(n)})$-holomorphic maps. Then the $C^0$-convergence
$u_n \longrightarrow u$ implies  the $L^{k+1,p}$-convergence for
any $p<\infty$, and $C^{k+1, \alpha}$-convergence if
$(J^{(n)}_S, J^{(n)})$ converge to $(J_S, J)$ in $C^{k,\alpha}$,
$0<\alpha<1$.
}

\smallskip
\state Proof. The map $u$ is continuous and in local coordinates
$w_1,\ldots, w_{2n}$ on $X$ and $z=x+iy$ on $S$ the~equation has the~form
$$
du(z) + J(u(z)) \scirc du \scirc J_M = 0,
$$
which is equivalent to $\dbar_{J\scirc u} u =0$.  Using {\sl Lemma 3.2.1}
and induction in $k'=0\ldots k$, one can obtain the~regularity of $u$.

Similarly, for $J^{(n)}$ and $u^{(n)}$ satisfying the hypothesis of the
corollary one gets
$$
\dbar_{J\scirc u}(u^{(n)}-u)=
(\dbar_{J\scirc u} - \dbar_{J^{(n)} \scirc u^{(n)}})u^{(n)}
\longrightarrow 0 \quad\hbox{in $L^{k',p}_\loc$
(resp.\ in $C^{k',\alpha}_\loc$)}
$$
by induction in $k'=0\ldots k$.
\qed

\smallskip
\state Remark. The corollary implies that for a compact Riemann surface
$S$ the topology in the space $\calp$ of $(J_S, J)$-holomorphic maps $u:S\to
X$ is independent of the particular choice of the functional space
$L^{k', p}(S,X)\supset \calp$ with $1\le k' \le k+1$, $1< p<\infty$, and
$kp>2$, provided $J_S$ and $J$ are changing $C^k$-smoothly. In the same
way, {\sl Lemma 3.2.1} implies that for $J\in C^k$ with $k\ge1$ the
differential structure on $\calp$ is also independent of the particular
choice of a functional space $L^{k',p}(S,X)$ with $1\le k' \le k$,
$1<p<\infty$ and $k'p>2$.

%%%%%%%%% vstavka 2  %%%%%%%%%%%%%%%%%%%

\medskip
\state Lemma 3.2.3. {\it For any natural numbers $\mu>\nu\ge0$ and real
numbers $p>2$, $\alpha<{2\over p}$ and $\gamma>0$ there exists $C =
C(\mu, \nu, \alpha, p, \gamma)>0$ with the following property.
Let $J$ be an almost complex structure in
$B \subset \cc^n$ with $J(0)=J\st$ and let $u:\Delta \to B$ be a
$J$-holomorphic map of the form $u(z)=z^\mu(P(z) +z^\nu v(z))$, where
$P(z)$ is some (holomorphic) polynomial of the degree $\nu$ and
$v\in L^{1,p}(\Delta, \cc^n)$ with $v(0)=0$. Suppose that $\Vert J-J\st
\Vert_{C^1(B)} \le \gamma$ and $\Vert u\Vert_{L^{1,p}(\Delta)} \le \gamma$.
Then for any $0<r<{1\over2}$
$$
\bigl\Vert z dv \bigr\Vert_{C^0(\Delta(r))} +
\bigl\Vert z dv \bigr\Vert_{L^{1,p}(\Delta(r))}
\le
C \cdot r^\alpha\cdot \Vert u \Vert_{L^{1,p}(\Delta)},
\eqno(3.2.2)
$$

In particular, $du(z) = d(z^\mu P(z)) + o(|z|^{\mu+ \nu -1 +\alpha})$
for any $\alpha<1$.
}

\state Proof. For $0<r<{3\over4}$ we define the map $\pi_r:B\to B$, setting
$\pi_r(w) \deff r^\mu w$. We also set $J^{(r)} \deff \pi_r^*J$, $u^{(r)}(z)
\deff \pi_r^{-1} \scirc u(rz)$, $P^{(r)}(z)\deff P(rz)$ and $v^{(r)}(z)=
r^\nu v(rz)$. Then $J^{(r)}$ is an almost complex
structure in $B$ with $\Vert J^{(r)} -J\st \Vert_{C^1(B)} \le r^\mu
\Vert J-J\st\Vert_{C^1(B)}$, and $u^{(r)} \equiv z^\mu (P^{(r)}(z) +
z^\nu v^{(r)}(z))$ is a $J^{(r)}$-holomorphic.

Without losing generality, we can suppose that $\alpha>0$. Set $\beta \deff
1+\alpha - {2\over p}$ and $q\deff {2\over1-\beta}$. Then $\alpha<\beta<1$,
$\beta + {2\over p} -1=\alpha$, and $q>2$. {\sl Lemma 3.1.3} implies that
$\Vert v \Vert_{C^{0,\beta}(\Delta(2/3))} + \Vert dv \Vert_{L^q(\Delta(2/3))}
\le C_1 \cdot \Vert u \Vert_{L^{1,p}(\Delta)}$. Here the constant $C_1$, as
well as the  constants $C_2,\ldots,C_6$ below, depend only 
on $\mu, \nu, p, \alpha$,
and $\gamma$, but are independent of $r$. Consequently,
$$
\Vert u^{(r)} - z^\mu P^{(r)}(z) \Vert_{L^{1,q}(\Delta)}
\le C_2\cdot r^{\nu + \beta} \cdot
\Vert u \Vert_{L^{1,p}(\Delta)}.
$$
Furthermore, due to {\sl Corollary 3.2.2}, we have $\Vert u^{(r)}
\Vert_{L^{2,p}(\Delta)} \le C_3\cdot \Vert u \Vert_{L^{1,p}(\Delta)}$. Thus
$$
\Vert \dbar_{\!J\st}( u^{(r)} - z^\mu P^{(r)}(z) ) \Vert_{L^{1,p}(\Delta)} =
\Vert (\dbar_{\!J\st} -\dbar_{\!J^{(r)} \scirc u^{(r)} } )
( u^{(r)}) \Vert_{L^{1,p}(\Delta)}
\le
$$
$$
\le C_4\cdot r^\mu \cdot
\Vert u \Vert_{L^{1,p}(\Delta)}.
$$
Applying {\sl Lemma 3.2.1}, we obtain
$$
\bigl\Vert z dv^{(r)} \bigr\Vert_{L^{1,p} (\Delta(2/3))} \le
C_5 \cdot r^{\nu + \beta} \cdot \Vert u \Vert_{L^{1,p}(\Delta)},
\eqno(3.2.3)
$$
which is equivalent to
$$
\bigl\Vert z dv \bigr\Vert_{L^{1,p}(\Delta(2r/3))} \le
C_5 \cdot r^{\beta + 2/p -1} \cdot \Vert u \Vert_{L^{1,p}(\Delta)}.
$$
On the other hand, $(3.2.3)$ implies that
$$
\bigl\Vert z dv^{(r)} \bigr\Vert_{C^0 (\Delta(2/3))} \le
C_6 \cdot r^{\nu + \beta} \cdot \Vert u \Vert_{L^{1,p}(\Delta)},
$$
and consequently
$$
\bigl\Vert z dv \bigr\Vert_{C^0(\Delta(2r/3))} \le
C_6 \cdot r^\beta \cdot \Vert u \Vert_{L^{1,p}(\Delta)}.
$$
\qed

%%%%%%%%% kinec' vstavky 2  %%%%%%%%%%%%%%%%%%%

\medskip
\noindent \sl 3.3. Perturbing Cusps of Complex Curves.

\medskip\rm
For understanding the rest of this lecture the reader should consult
\S 6.2 and \S 6.3 for the definition and some elementary properties of
Gromov's $\dbar $-operator.

\smallskip
\state Lemma 3.3.1. {\it
For a given $p$, $2<p<\infty$, and $\gamma>0$ there exist constants
$\epsi=\epsi(p,\gamma)$ and $C=C(p,\gamma)$ with the following
property. Suppose that $J$ is a $C^1$-smooth almost complex
structure in the unit ball $B \subset \cc^n$
with $\Vert J-J\st \Vert_{C^1(B)}
\le \epsi$ and $u\in  L^{1,p}(\Delta, B(0,{1\over2}))$ is
a $J$-holomorphic map with $\Vert u \Vert_{L^{1,p}(\Delta)}
\le \gamma$. Then for every almost complex structure $\tilde J$ in
$B$ with $\Vert \tilde J-J \Vert_{C^1(B)} \le \epsi$
there exists a $\tilde J$-holomorphic map $\tilde u: \Delta \to
B$ with
$$
\Vert \tilde u-u \Vert_{L^{1,p}(\Delta)} \le
C\cdot \Vert \tilde J-J \Vert_{C^1(B)}
$$
such that $\tilde u(0) = u(0)$.
}

\smallskip\rm
\state Proof. Let $J_t$ be a curve of $C^1$-smooth almost
complex structures in $B$ starting at $J_0=J$ and depending
$C^1$-smoothly on $t\in [0,1]$. Consider an ordinary differential
equation for $u_t \in L^{1,p} (\Delta, B)$ with the initial
condition $u_0=u$ and

$$
\msmall{du_t\over dt} = - T_{J_t \scirc u_t,\, R_t}
\biggl(\msmall{dJ \over dt} \scirc du_t \scirc J_\Delta \biggr),
$$

where $R_t$ is defined by the relation $D_{J_t, u_t}=\dbar_{J_t, u_t}
+ R_t$. (See {\sl paragraph 6.2}). Since $J_0$ and $R_0$ satisfy
the hypothesis of  {\sl Lemma 3.2.1} and $R_t$ depends
$L^p$-continuously on $J_t\in C^1$ and $u_t \in
L^{1,p}$, the solution exists in some small interval $0\le t \le
t_0$. For such a solution using (6.3.1) one has
$$
\left\Vert \msmall{du_t\over dt} \right\Vert_{L^{1,p}} +
\left\Vert \msmall{dR_t\over dt} \right\Vert_{L^p}
\le C\cdot (\Vert u_t\Vert_{L^{1,p}} +
\Vert J_t  - J\st \Vert_{C^1} )
\cdot \left\Vert \msmall{dJ_t\over dt} \right\Vert_{C^1}.
$$
This implies the existence of the solution of our equation for all
$t\in [0,1]$, provided
$$
\int_{t=0}^1
\left\Vert \msmall{dJ_t\over dt} \right\Vert_{C^1}
dt \le \epsi.
$$
\qed

\smallskip
In this paragraph we suppose that some $p$ with $2<p<\infty$ is fixed.

\smallskip\nobreak
\state Lemma 3.3.2. {\it
For every  $\gamma >0$ and every pair of integers $\nu\ge1$, $\mu\ge1$
there exists an $\epsi=\epsi(\mu,\nu,\gamma)>0$
such that for an almost complex structure $J_0$ in $B$ with
$\Vert J_0-J\st \Vert_{C^1(B)}\le \epsi$, $J_0(0)=J\st$, and for
a $J_0$-holomorphic map $u_0: \Delta \to B(0, {1\over2})$ with $\Vert u_0
\Vert_{L^{1,p}(\Delta)} \le\gamma $ and with the~multiplicity $\mu$ at
$0\in \Delta$ the~following holds:

\smallskip
\sli If $\nu\le2\mu-1$, then for every almost complex structure $J$
in $B$ with $J(0)=J\st$ and for every $v\in L^{1,p}(\Delta,\cc^n)$ with
$\Vert v \Vert_{L^{1,p}(\Delta)} + \Vert J-J\st \Vert_{C^1(\Delta)}
\le\epsi$ there exists $w\in L^{1,p}
(\Delta,\cc^n)$ with $w(0)=0$, satisfying
$$
\dbarj(u_0+z^\nu(v+w))=0\eqno(3.3.1)
$$
and
$$
\Vert w \Vert_{L^{1,p} (\Delta)}\le
C\cdot \left(
\Vert v \Vert_{L^{1,p}(\Delta)} +
\Vert J-J_0 \Vert_{C^1(\Delta)}
\right).\eqno(3.3.2)
$$

\smallskip
\slii If $\nu\ge2\mu$, then for every $v\in L^{1,p}(\Delta,\cc^n)$
with $\Vert v \Vert_{L^{1,p}(\Delta)}\le\epsi$ there exists $w\in
L^{1,p}(\Delta,\cc^n)$ with $w(0)=0$, satisfying
$$
\bar\partial_{J_0}(u_0+z^\nu(v+w))=0\eqno(3.3.3)
$$
and
$$
\Vert w \Vert_{L^{1,p} (\Delta)}\le
C\cdot \Vert v \Vert_{L^{1,p}(\Delta)}.
\eqno(3.3.4)
$$
}

\smallskip
\state Proof. If $\nu\le 2\mu -1$, we fix a curve $J_t$, $t\in [0,1]$, of
$C^1$-smooth almost complex structures in $B$, which starts in $J_0$,
finishes in $J$ and satisfies the conditions $J_t(0)=J\st$
and $\int_{t=0}^1\Vert \dot J_t \Vert_{C^1(B)} dt \le 2\epsi$,
where $\dot J_t \deff dJ_t/dt$. If $\nu\ge2\mu$, we simply set
$J_t\equiv J_0$.

As in the previous lemma, we want to find a needed $w$ by solving
for $t\in [0,1]$ and $w_t\in L^{1,p}(\Delta,\cc^n)$ the equation
$$
z^{-\nu}\dbar_{J_t}(u_0 + z^\nu(t\cdot v + w_t)) =0.
$$
However, this time we need to consider the fact that
now we are dealing with different (almost) complex structures on $B$, namely
$J\st$ and $J_t$ for any fixed $t\in [0,1]$. Thus, we write the
last equation in the more correct form:
$$
(x + y J\st)^{-\nu}\dbar_{J_t}
(u_0 + (x + y J\st)^\nu(t\cdot v + w_t)) =0.
\eqno(3.3.5)
$$

\smallskip
The differentiation of (3.3.5) with respect to $t$ gives
$$
(x + y J\st)^{-\nu}D_{u_t, J_t}
((x + y J\st)^\nu(v + \dot w_t))
+(x + y J\st)^{-\nu}\dot J_t
\scirc du_t \scirc J_\Delta
=0,
\eqno(3.3.6)
$$
where $u_t\deff u_0+ (x + y J\st)^\nu(t\cdot v + w_t)$
and $J\st$ also denotes the pull-back of the standard
complex structure on $E\deff u_t^*TB\cong \cc^n$.

\smallskip
First we show that the operator $(x + y J\st)^{-\nu} \scirc
D_{u_t, J_t} \scirc (x + y J\st)^\nu$ has the form
$\dbar_{J^{(\nu)}_t} + R^{(\nu)}_t$ for an appropriate (almost)
complex structure $J^{(\nu)}_t$ in $E\cong \cc^n$ and
$\rr$-linear operator $R^{(\nu)}_t$. In fact, the explicit formula
(6.2.5) for $D_{u, J}$ shows that taking the standard connection
$\nabla$ in $TB\cong \cc^n$ and identifying $\Lambda^{(0,1)}\Delta
\cong \cc$ one gets
$$
D_{u_t, J_t}(v)=\msmall{1\over2}\left( \msmall{\d \over \d x}v +
J_t \msmall{\d \over \d x}v \right) + R^{(0)}_t(v)
$$
with $R^{(0)}_t\in C^0(\Delta, \endo_\rr(\cc^n))$. Moreover,
the formula (6.3.1) implies that
$$
\Vert R^{(0)}_t \Vert_{L^p(\Delta)} \le  \Vert J_t \Vert_{C^1(B)}
\cdot \Vert du_t \Vert_{L^p(\Delta)}.
$$
Hence
$$
2(x + y J\st)^{-\nu} \scirc D_{u_t, J_t}
\scirc (x + y J\st)^\nu   =  \Bigl[\msmall{ \d \over \d x} +
(x + y J\st)^{-\nu} \scirc J_t \scirc (x + y J\st)^\nu
\msmall{\d \over \d y}  \Bigr]+
$$
$$
+ \Bigl[ \nu \cdot(x + y J\st)^{-\nu} \scirc
(1 + J_t \scirc J\st) \scirc (x + y J\st)^{\nu -1}
+ 2(x + y J\st)^{-\nu} \scirc R^{(0)}_t \scirc (x + y J\st)^\nu
\Bigr]=
$$
$$
=2\dbar_{J^{(\nu)}_t} + 2\,R^{(\nu)}_t.
$$
One has the obvious identities $(J^{(\nu)}_t)^2 \equiv -1$,
$$
\Vert J^{(\nu)}_t -J\st \Vert_{C^0(\Delta)}=
\Vert (x + y J\st)^{-\nu} \scirc (J_t -J\st)
\scirc (x + y J\st)^\nu \Vert_{C^0(\Delta)} =
\Vert J_t -J\st \Vert_{C^0(B)},
$$
and
$$
\Vert (x + y J\st)^{-\nu} \scirc R^{(0)}_t \scirc (x + y J\st)^\nu
\Vert_{L^p(\Delta)} =\Vert R^{(0)}_t \Vert_{L^p(\Delta)} .
$$
Furthermore, $1 + J_t(0) J\st=0$; hence $\bigl\Vert \> |z|^{-1}(1 +J_t J\st)
\bigr\Vert_{C^0(\Delta)} \le \Vert J_t -J\st \Vert_{C^1(B)}$. This gives
the estimate
$$
\Vert R^{(\nu)}_t \Vert_{L^p(\Delta)} \le \bigl(C_1 \cdot \nu +
C_2 \cdot \Vert du_t \Vert_{L^p(\Delta)} \bigr) \cdot
\Vert J_t -J\st \Vert_{C^1(B)}.
$$

\smallskip
To show the existence of the solution of (3.3.5) for all
$t\in [0,1]$, we need the estimate
$$
\left\Vert z^{-\nu}\cdot \dot J_t (u_0 + z^\nu w)
\scirc d(u_0 + z^\nu w) \right\Vert_{L^p(\Delta)}
\le C_1(\mu, \nu) \cdot \Vert u_0  \Vert_{L^{1,p}(\Delta)}^2
\cdot \Vert \dot J_t \Vert_{C^1(B)}
\eqno(3.3.7)
$$
for all sufficiently small $w\in L^{1,p}(\Delta, \cc^n)$.
The estimate trivially holds for $\nu\ge 2\mu$, because
in this case $\dot J_t\equiv 0$. Otherwise for
$\lambda \deff \min\{\mu, \nu\}$ from
{\sl Corollary 3.1.3.} one obtains
$$
\left\Vert z^{-\lambda}(u_0 + z^\nu w)
\right\Vert_{L^\infty (\Delta)} +
\left\Vert z^{-\lambda +1} d(u_0 + z^\nu w)
\right\Vert_{L^p(\Delta)}
\le C_2(\mu, \nu)\cdot
\left\Vert u_0 \right\Vert_{L^{1,p}(\Delta)}
\eqno(3.3.8)
$$
which gives the estimate (3.3.7).

\smallskip
Now consider the ordinary differential equation for $w_t \in L^{1,p}
(\Delta, B)$ with the initial condition $w_0\equiv 0$ and
$$
\msmall{dw_t\over dt} = - T_{J^{(\nu)}_t,\, R^{(\nu)}_t}
\Bigl(z^{-\nu}\cdot \dot J_t \scirc du_t \scirc J_\Delta  + v \Bigr).
\eqno(3.3.9)
$$
As in {\sl Lemma 3.3.1} one has the estimate
$$
\biggl\Vert \msmall{dw_t\over dt} \biggr\Vert_{L^{1,p}(\Delta)} +
\biggl\Vert \msmall{dR_t^{(\nu)}\over dt} \biggr\Vert_{L^p(\Delta)}
\le
$$
$$
\le C\cdot (\Vert u_t\Vert_{L^{1,p}(\Delta)} +
\Vert J_t - J\st \Vert_{C^1(B)} )
\cdot \left( \bigl\Vert \msmall{dJ_t\over dt} \bigr\Vert_{C^1(B)}
+ \bigl\Vert v \bigr\Vert_{L^{1,p}(\Delta)} \right),
$$
which implies the existence of the solution of (3.3.9) for all
$t\in [0,1]$. \qed

%%%%%%%%%%%%%%%%%%%%%%%%%%%%%%%%%%%%%%%%%%%%%%%%%%%%%%%%%%%%%%%%%%

\bigskip\noindent\sl
3.4. Primitivity.

\medskip\rm
We need to study first the notions of primitivity and distinctness of
complex curves in $C^1$-smooth nonintegrable structures. Let us start with
the following

\smallskip\noindent
\state Definition 3.4.1.  Two  $J$-holomorphic maps $u_1: (S_1,j_1) \to X$
and $u_2: (S_2,j_2) \to X$ with $u_1(a_1)=u_2(a_2)$ for some $a_i\in S_i$ are
called {\sl distinct in $(a_1,a_2)$}
if there are no neighborhoods $U_i\subset S_i$ of $a_i$ with $u_1(U_1) =
u_2(U_2)$.
We call $u_i: (S_i,j_i) \to X$ {\sl distinct} if they are distinct for all
pairs $(a_1,a_2)\in S_1 \times S_2$ with $u_1(a_1)=u_2(a_2)$.

\state Definition 3.4.2. A $J$-holomorphic map $u: (S,j) \to X$ is called {\sl
primitive} if there are no disjoint non-empty open sets $U_1, U_2\subset S$
with $u(U_1) =u(U_2)$. Note that  $u$ must be nonconstant.

An important regularity property of a $J$ - complex  curve with $J\in C^1$
is contained in the following

\smallskip\noindent
\bf Theorem 3.4.1. \it Let $(S_1, j_1)$ and $(S_2, j_2)$ be smooth connected
complex curves and $u_i : (S_i,j_i) \to (X,J)$ non-constant $J$-holomorphic
maps with $J \in C^1$. If there are non-empty open sets $U_i \subset S_i$ with
$u_1(U_1)= u_2(U_2)$, then there exists a smooth {\sl connected} complex curve
$(S, j)$  and a $J$-holomorphic map $u: (S,j) \to (X,J)$ such that
$u_1(S_1) \cup u_2(S_2) =u(S)$ and $u:S \to X$ is primitive.

Moreover, maps $u_i :S_i \to X$ factorize through $u: S \to
X$, \i.e., there exist holomorphic maps $g_i : (S_i,j_i) \to (S,j)$ such that
$u_i = u \scirc g_i$.

\smallskip\noindent
\bf Proof. \rm

\smallskip\noindent
{\sl Step 1.} \it For every $J$-holomorphic imbedding $u:\Delta \to (X,J)$
with $J\in C^1$ and for every $p\in \Delta$ there exist a neighborhood 
$U \subset X$ of
$u(p)$, a $C^1$-smooth integrable complex structure $J_1$ in $U$ and
$J_1$-holomorphic coordinates $(w_1, \ldots, w_n)=(w_1,w')$ (\ie
$w'=(w_2, \ldots, w_n))$ in $U$ such that
$u(\Delta) \cap U = \{ (w_1, w')\in U: w'=0\}$ 
and $J|_{u(\Delta)} = J_1|_{u(\Delta)}$. In particular  
$w_1|_{u(\Delta)}$ is a holomorphic coordinate in $u\inv(U) \subset \Delta$ 
and 

$$
|J(w_1, w')-J_1(w_1, w')| \le C(|w'|).
\eqno(3.4.1)
$$
\rm

\smallskip For the proof take any $J$-Hermitian metric
$h$ on $X$. Let $N_1$ denote the normal bundle $u_1^*TX/ du(TS_1)$ with fiber
$N_{1,z}$ over $z\in S_1$. Then the exponential map $\phi: \xi \in N_{1,z} \mapsto
\exp_{u_1(z)} (\xi)$ is a local diffeomorphism in some neighborhood $V$
of the zero section of $N_1$ over $S_1$. Note that along the zero section 
we can
identify $\phi$ with $u_1$.

Equip $N_1$ with a canonical holomorphic structure of the quotient $u_1^*TX/
du(TS_1)$. Since $h$ was chosen $J$-Hermitian, the operator of the complex
structure $J_{N_1}$ coincides along zero section with the structure $\phi^*J
\equiv u_1^*J$. Note that $J_{N_1}$ is {\sl integrable} in $V$.

\smallskip\noindent
{\sl Step 2. (``Unique Continuation Lemma'').} \it
Let the open subset $U \subset X$, $C^1$-smooth complex structure
$J_1$ in $U$, and $J_1$-holomorphic coordinates $(w_1, \ldots, w_n)=(w_1,w')$
be as in the previous step. Set $C \deff \{ (w_1, w')\in U: w'=0\}$.
Let $v: \Delta \to U$ be a $J$-holomorphic map. Then either $v(\Delta) \subset
C$ or $v\inv(C)$ is discrete. \rm

To show this, we use $(w_1, w')$ as complex coordinates and consider $U$ as
an open subset of $\cc^n$, so that $J_1= J\st$. Set $v' \deff w' \scirc v$,
\i.e., $v'$ is obtained from $v$ by forgetting the first component. Then
$v\inv(C)=v^{\prime\,-1}(0)$. Furthermore,
$$
|\dbar_{J\st}v'(z)| \le |\dbar_{J\st}v(z)| = |\dbar_{J_1}v(z)| =
|\dbar_Jv(z)-\dbar_{J_1}v(z)| \le
$$
$$
\le \norm{dv}_{L^\infty} \cdot
|J(v'(z))- J_1(v'(z))|,
$$
which means that locally $|\dbar_{J\st}v'| \le C |v'|$. By {\sl Lemma 3.1.1}
either $v^{\prime\,-1}(0)$ is discrete, or $v'$ is identically zero.
\qed

\state Remarks.\ 1. Note that, in the latter case, $v:\Delta \to C$ 
is holomorphic.
\par\noindent \bf
2. \rm The ``Unique Continuation Lemma'' remains true for Lipschitz $J$
(\i.e., $J\in C^{0,1}$)
but not for H\"older and continuous $J$ (\i.e., $J\in C^{0,\alpha}$ with
$0\le\alpha<1$). To see this, we note that if $J\in C^{0,\alpha}$ with
$0<\alpha\le1$, then any $J$-holomorphic map $u: \Delta \to (X,J)$ is
$C^1$-smooth by standard regularity for the  elliptic equation $\dbar_J u=0$.
If $u$ is an immersion, then $du: T\Delta \to u^*TX$ is an imbedding of complex
(not holomorphic!) bundles, and $N \deff u^*TX/ du(T\Delta)$ has a canonical
complex structure. Equip $N$ with any holomorphic structure compatible with
the complex one. Then {\sl Step 1} from above goes through and we can
locally construct $J_1$, which coincides with $J$ along some neighborhood
of a given point $p\in \Delta$. Instead of (1) we obtain
$$
|J(w_1, w')-J_1(w_1, w')| \le C(|w'|^\alpha).
\eqno(3.4.2)
$$
Thus the ``similarity principle'' and {\sl Step 2} are applicable if $\alpha=1$
(Lipschitz case) but fail if $\alpha<1$ (H\"older case).

\state Exercise. {\rm Construct a counter example to 
``unique continuation lemma''
for the H\"ol\-der $J$ in the following way. Take two different smoothly
imbedded disks
which coincide along some open subset, and try to construct a complex
structure $J$ to make both $J$ - complex . It appears that for
appropriately chosen disks one succeeds with $J \in C^{0,\alpha}$ for at least
some $0<\alpha<1$.
}

\medskip\noindent
{\sl Step 3}. {\it Proof of the Theorem in the case when $u_i: (S_i, j_i) 
\to
(X,J)$  are immersions.}

%Set $M\deffu_1(S_1) \cup u_2(S_2)$.
\rm
Consider set $\cala$ of those $V$ such that $V$ is an open subset in one of
$S_i$ and $u_i: V \to X$ is injective. Equip every
$V \in \cala$ with the complex structure and with a $J$-holomorphic map $\ti u:
V\to X$ induced from $S_i$. Define the equivalence relation $\sim$ on the
disjoint union $\sqcup_{V\in \cala}V$, identifying points $p'\in V'$ and
$p'' \in V''$ if there exist $V\in \cala$ and $p\in V$ such that
$\ti u(V) \subset \ti u(V') \cap \ti u(V'')$ and $\ti u(p) = \ti u(p') =
\ti u(p'')$.

Set $S \deff \bigl(\sqcup_{V\in \cala}V\bigr) / \sim$ with projection $\pi :
\sqcup_{V\in
\cala}V \to S$. Induce the quotient topology on $S$ declaring
$\{ \pi(V) :V\in \cala\}$ as the basis of the topology. Then $\ti u$ induces
the map $u: S \to X$ which is holomorphic. We shall show that $S$ and $u$
obey the desired properties. The main point at this step is to show that the
quotient topology on $S$ is Hausdorff.

Let $p'$ and $p''$ be two distinct points on $S$. If $u(p') \not= u(p'')$,
then there exist disjoint neighborhoods $u(p') \in U' \subset X$ and
$u(p'') \in U'' \subset X$ and their pre-images $V' \deff u\inv(U')$ and
$V'' \deff u\inv(U'')$ are disjoint neighborhoods of $p'$ and $p''$.
It remains to consider the case when $u(p') = u(p'')$.

Let $\ti p' \in \ti V'$ and $\ti p''\in \ti V''$ be representatives, \i.e.,
$\ti V', \ti V'' \in \cala$, $\pi(\ti p')=p'$, and $\pi(\ti p'')=p''$.
Shrinking $\ti V'$, if needed, we may assume that for $\ti u(\ti V')$
there exists a neighborhood $U \subset X$, an (integrable)
complex structure $J_1$ in $U$ and $J_1$-holomorphic  coordinates $(w_1, w')$
in $U$ with the properties of {\sl Step 1}. Shrinking $\ti V''$, we 
may assume that
$\ti u(\ti V'')$ is also contained in $U$. Using {\sl Step 2}, we may
additionally assume that $u(p')=u(p'')$ is the only intersection point of
$\ti u(\ti V')$ and $\ti u(\ti V'')$. Then $\pi(\ti V')$ and $\pi(\ti V'')$
are disjoint neighborhoods of $p'$ and $p''$. This shows that $S$ is
Hausdorff.

By assumption, maps $u_i: S_i \to X$ are immersions. Thus every $S_i$ can be
covered by open sets $V \subset S_i$ such that restrictions $u_i|_V$ are
imbeddings. Every such $V$ belongs to the atlas $\cala$, and the canonical
projection $\pi: V \to S$ is a holomorphic imbedding. By the definition of
the equivalence relation $\sim$, these local maps $\pi: V\subset S_i \to S$
can be glued together to holomorphic maps $g_i: S_i \to S$. One can see that
the constructed $(S,j)$, $u:S \to X$, and $g_i: S_i \to S$ have the desired
properties.
%%%%%%%%%%%%%%%%%%%%%page 28%%%%%%%%%%%%%%%%%%%%%%%%%%%

\smallskip\noindent\bf Non immersed case. 
 {\it Here we consider the case where $u_i: (S_i,j_i) \to (X,J)$
are not necessarily immersions.}

\smallskip\rm  From {\sl Proposition 6.3.2} of {\sl Lecture 6},
maps $u_i$ define bundle homomorphisms $du_i : TS_i \to E_i :=
u_i^*(TX)$ which are holomorphic w.r.t. the canonical holomorphic structure on
$E_i$. This implies that the zero set $A_i:=\{y \in S_i:  du_i(y)=0 \}$ is
discrete in $S_i$. In particular, every $u_i$ is an immersion outside $A_i$.

Set $\check S_i \deff S_i \bs A_i$. Then $u_i: \check S_i \to X$ are 
immersions,
and the previous arguments remain valid. Let $(\check S, \check \jmath) $,
$\check u:\check S \to X$, and $\check g_i: \check S_i \to \check S$ be
the corresponding objects.

Numerate points in $A_1$ in any order, $A_1 = \{\, a_1, a_2, \ldots \,\}$. 

\smallskip\noindent {\sl Step 4.} {\it We
state that for $k=1,2,\ldots\,$ there exist Riemann surfaces $S^{(k)}$ and
holomorphic maps $g_1^{(k)}: S_1^{(k)} \to S^{(k)}$ such that $S_1^{(k)} \deff
\check S_1 \cup \{a_1, \ldots a_k\}$, $g_1^{(k)}|_{\check S} =\check g_1$, and
$S^{(k)} = \check S \cup g_1^{(k)}(\{a_1, \ldots a_k\})$. }

\smallskip In other words,
at the $k$-th step we add to $S_1^{(k-1)}$ the point $a_k$, and 
extend $g_1^{(k-1)}
: S_1^{(k-1)} \to S^{(k-1)}$ through $a_k$ holomorphically, adding, 
if needed, some new point to $S^{(k-1)}$.

The proof of this statement proceeds inductively.
For a given $k$ take a sufficiently small neighborhood $V_k \subset S_1$ of
$a_k$ which contains no other points of $A_1$. Set $\check V_k \deff V_k
\bs \{a_k\}$. If the restriction $g_1^{(k-1)}|_{\check V_k}: \check V_k \to
S^{(k-1)}$ extends to a holomorphic map $g_1^{(k)}|_{V_k}: \check V_k \to
S^{(k-1)}$, one has nothing to do but set $S^{(k)}\deff S^{(k-1)}$.

Otherwise we set $S^{(k)}\deff S^{(k-1)} \sqcup \{b_k\}$ and $g_1^{(k)}(a_k)
= b_k$. We must prove  existence of appropriate topology and complex
structure on the constructed $S^{(k)}$. We give a proof of these properties
in full detail.

\smallskip
Choose $C^2$-smooth local complex (not
holomorphic) coordinates $(w_1, \ldots, w_n)$ in a neighborhood of
$u_1(a_k) \in X$ such that a complex structure $J'$ defined by $(w_1, \ldots,
w_n)$ coincides with $J$ in $u_1(a_k)$. Let $z$ be a local complex coordinate
on $S_1$ in a neighborhood of $a_k$. We may assume that $(w_1, \ldots, w_n)$
(resp.\ $z$) are centered in  $u_1(a_k)$ (resp.\ in $a_k$).
{\sl Lemma 3.1.1} provides that in these coordinates the map $u_1 : S_1 \to X$ has
the form $u_1(z) = z^l \cdot (v_1(z), \ldots,v_n(z))$ with $\vec v=(v_1,
\ldots, v_n) \in L^{1,p}$ and $\vec v(0) \not=0$.

Choose a coordinate $w_p$ such that the corresponding component $v_p$ does not
vanish at $z=0$, $v_p(0) \not=0$. Then there exists the $l$-th root of 
this $v_p$,
\i.e., $v_p(z) = (f(z))^l$ for some continuous non-vanishing function $f(z)$.
Since $v_p(z)$ is $L^{1,p}$-smooth, so is $f(z)$. Shrinking our
neighborhood $V_k$, if needed, we may assume that $f(z)$ is not vanishing in
$V_k$. Thus for an appropriate neighborhood $U$ of $u_1(a_k)$ there exists
a $C^2$-smooth map $w_p:U \to \cc$ such that $h_p \deff w_p\scirc u_1 (z)=
(zf(z))^l$ for some $L^{1,p}$-smooth non-vanishing function $f(z)$ defined
in a neighborhood of $a_k$. We may additionally assume that
$$
|f(z)-f(0)|\le \eps <\!\!<1 \qquad\hbox{for all $z\in V_k$}.
\eqno(3.4.3)
$$
Note that for every sufficiently small $b \in \cc$ we have at least one
solution $z\in V_k$ of the equation $h_p(z)=b$ with $|z  f(0)|^l \le
2|b|$. In order to see this, one considers $h_p$-images of 
circles $\gamma_t \deff
\{\, |zf(0)|^l = t|b|\,  \}$. Then due to (3.4.3), the image 
$h_p(\gamma_{1\over2})$
lies in the disk $\Delta(|b|)$, whereas $h_p(\gamma_2)$ lies outside of
$\Delta(|b|)$. Homotopy argument shows that $b\in \psi(\gamma_t)$ for some
$t\in ]{1\over2}, 2[$. Similarly, one obtains the estimate
$$
c\cdot |h_p(z)|^{1\over l} \le |z| \le C\cdot |h_p(z)|^{1\over l}
\qquad\hbox{for $z\in V_k$}
\eqno(3.4.4)
$$
with some constants $c$ and $C$.
Take a sufficiently small disk $\Delta(r)$ in $\cc$ (range of $h_p(z)= \omega_p
\scirc
u_1(z)$) and consider the set $V'_k \deff \{ z\in V_k: |z|^l < {2r \over
|f(0)|^l}, \, |h_p(z)|< r\}$. Consider a sequence $z_i \in V'_k$ such that
$h_p(z_i)$ converges to some $b\in \Delta(r)$. Then some subsequence of $z_i$
converges to $z^*$ with $|z^*| \le {2r \over |f(0)|^l}$. But then $|h_p(z^*)|=
|z^*f(z^*)|^l =|b|$ which implies $|z^*|^l = {|b| \over  |f(z^*)|^l} <{2r
\over |f(0)|^l}$, so that $z^* \in V'_k$. This shows that the map $h_p: V_k
\to \Delta(r)$ is proper. By (4), $h_p(z)=0$ for $z\in V'_k$ implies $z=0$.
Thus, for $\check V'_k \deff \{z\in V'_k : z\not=0\} = V'_k \bs\{a_k\}$ the map
$h_p: \check V'_k \to \check\Delta(r)$ is also proper.

Recall that the bundle homomorphism $du_1: TS_1 \to E_1\deff u_1^*TX$
is holomorphic with respect to the canonical holomorphic structure on $E_1$.
Thus, in $V'_k$ we can represent $du_1$ in the form $du_1(z)= z^{l-1}\cdot
s(z)$ with $s(z)$ also holomorphic and $s(0)\not=0$. This implies that for
$z\not=0 \in V'_k$ sufficiently close to $0$ the image $du_1(T_zS_1)$ of
the tangent space $T_zS_1$ is close to a real 2-dimentional space
$\cc s(0)\subset E_1|_0=T_{u_1(a_k)}X$. This is a complex line in
$T_{u_1(a_k)}X$ generated by $s(0)$.
We may assume that the coordinates $(w_1, \ldots, w_n)$ were chosen in such
a way that the $p$-th component of the vector $s(0)$ is non-zero. This means
that the linear map $dw_p|_{du_1(a_k)}: T_{u_1(a_k)}X \to \cc$ is not
degenerated $\cc s(0) \subset T_{u_1(a_k)}X$. Thus, for any $b$ sufficiently
close to $u_1(a_k)$, and for any real 2-dimensional space $L \subset T_bX$
sufficiently close to $\cc s(0)$ the linear map $dw_l|_b: T_bX \to \cc$ is
not degenerated on $L$.

Hence we can conclude that $h_p: \check V'_k \to \check \Delta(r)$ is not
degenerated provided $r$ was chosen sufficiently small. This means that
$h_p: \check V'_k \to \check \Delta(r)$ is a covering. Note that $\pi_1(\check
\Delta(r)) = \zz$ and that $\pi_1(\check V'_k)$ contains $\zz$ as a subgroup
generated by a small circle $\{ |z| = \rho\}$. From inclusion $(u_1)_*:
\pi_1(\check V'_k) \hookrightarrow \pi_1(\check \Delta(r))$ we conclude
that $\pi_1(\check V'_k) =\zz$, and that $h_p: \check V'_k \to \check
\Delta(r)$ is a finite covering. Since a small circle $\{ |z| = \rho\}$
around $a_k \in V'_k$ is a generator of $\pi_1(\check V'_k)$, we conclude
that $V'_k$ is a disk and $\check V'_k$ is a punctured disk.

Now recall that the map $h_p:\check V'_k \to \check \Delta(r)$ can be seen as
a composition of $g^{(k-1)}: \check V'_k \to S^{(k-1)}$ with a coordinate
function $w_p$ restricted on $\check u_1(\check V'_k)$. Set $\check W_k \deff
g^{(k-1)} (\check V'_k)$ and consider a holomorphic map $g^{(k-1)}: \check
V'_k \to \check W_k$. Then this map is again proper and non-degenerate
and therefore is also a finite covering. This implies that $\check W_k$ is also
a punctured disk.

Let $\check\psi: \check W_k \to \check \Delta$ be a biholomorphism which can
be seen as a local chart on $\check W_k \subset S^{(k-1)}$.
Define holomorphic structure on $S^{(k)}= S^{(k-1)} \sqcup \{b_k\}$ in the
following way. Interpret $b_k$ as a puncture point of $\check W_k$, set
$W_k \deff \check W_k\sqcup \{b_k\}$, extend $\check \psi$ to the map
$\psi: W_k \to \Delta$ setting $\psi(b_k) \deff 0$. Extend the topology on
$S^{(k)}$ in such a way that $\psi: W_k \to \Delta$ becomes a homeomorphism.
Declare $\psi$ to be a holomorphic coordinate on $S^{(k)}$. Note that the
composition $\check\psi \scirc g^{(k-1)}: \check V'_k \to \check \Delta$ is
bounded; thus it extends to a holomorphic map from $V'_k$ to $\Delta$.
Consequently, $\check\psi\inv: \check \Delta \to \check W_k \subset S^{(k-1)}$
cannot be holomorphically extended to $\psi\inv: \Delta \to S^{(k-1)}$, 
because
otherwise this would mean holomorphic extensibility of $g^{(k-1)}: S_1^{(k-1)}
\to S^{(k-1)}$ through $a_k$. This property insures that $S^{(k)}$ is
Hausdorff. Thus constructed, $S^{(k)}$ and $g^{(k)}: S_1^{(k)} \to S^{(k)}$
obey  the desired properties.

\medskip
Setting $S' \deff \cup_k S^{(k)}$ we obtain the minimal extension
of $\check S$ such that the map $\check g_1: \check S_1 \to \check S$ extends
holomorphically to $g'_1: S_1 \to S'$. A similar construction should be done
to extend $\check g_2 : \check S_2 \to S'$ to a holomorphic map
$g_2: S_2 \to S$. This finishes the proof of the {\sl theorem }.
\qed

\smallskip\rm
In the spirit of {\sl Theorem 3.4.1} one obtains the following result.

\state Proposition 3.4.2. \it Let $(X,J)$ be an almost complex manifold with
$J \in C^1$. Then any primitive $J$-holomorphic map $u: (S,j) \to X$
with connected irreducible $(S,j)$ is injective outside some countable
subset $A\subset S$. Moreover, $A$ can be chosen in such a way that
for any compact $K\Subset S$ the intersection $A \cap K$ is finite. \rm

\state Proof. We have shown that $A_1 \deff \{a\in S: du(a)=0\}$ is discrete;
therefore it is countable. Set $\check S:= S\setminus A_1$.
Consider
the set $A_2 \deff \{ (a,b) \in \check S\times \check S : u(a)=u(b) \}$.
We state that
$A_2$ is discrete in $\check S\times \check S\setminus \Delta $, where 
$\Delta $ is a diagonal. Otherwise there must exist a
sequence $(a_n,b_n)$ converging to some pair $(a^*,b^*) \in \check S\times
\check S$. Take a sufficiently
small neighborhood $V \subset \check S$ of $a^*$ and find a neighborhood
$U \subset X$ of $u(a^*)$, which obeys the properties of {\sl Step 1} from
the previous proof.
Note that obviously $u(b^*) = u(a^*)$. Now we obtain a contradiction
with argumentations from {\sl Step 2} above. This shows that $A_2$ must be
discrete in $\check S\times \check S\Delta $.

Thus it is also countable. Let $A_3$ be the
projection of $A_2$ onto the first $\check S$. Then $A:=A_1\cup A_2$ obeys
the desired properties.
\qed

\state Remark. In fact, positivity of the intersection of $J$ - complex  
curves
with $J\in C^1$ provides that there exists a {\it discrete} subset $A\subset
S$ fulfilling the condition of {\sl Proposition 3.4.2}. We will prove 
this in next section.

\bigskip\noindent\sl
3.5. Positivity of Intersection.

\smallskip\rm
Let us first recall the notion of the intersection number (index) of two
surfaces in $\rr^4$.  Let $M_1$ and $M_2$ be two-dimensional, oriented,
smooth surfaces in $\rr^4$ passing through the origin.
%, \ie $M_i$ is a $C^1$-image
We suppose further that both $M_1,M_2$ intersect the unit sphere $S^3$
transversally by curves $\gamma_1$ and $\gamma_2$, respectively, and that
$\gamma_1$ and $\gamma_2$ do not meet one another. Let $\tilde M_i$ be
small perturbations of $M_i$ making them intersect transversally.

\state Definition 3.5.1. {\it The intersection number of $M_1$ and $M_2$ is defined
to be the algebraic intersection number of $\tilde M_1$ and $\tilde M_2$. If
$M_1,M_2$ intersect only at zero, we also say that the number just
defined is the intersection index of $M_1$ and $M_2$ at zero. It will be
denoted by $\delta_0(M_1,M_2)$ or $\delta_0$.
}

\state Remark. {\rm This number is independent of the particular choice of
perturbations $\tilde M_i$. We shall use the fact that the intersection
number of $M_1$ and $M_2$ is equal to the {\it linking number} $l(\gamma_1,
\gamma_2)$ of the (reducible in general) curves $\gamma_i$ on $S^3$, see
[Rf].
}

\smallskip\rm
In the following theorem the structure $J$ is supposed to be of class
$C^1$.

\smallskip\noindent
\state Theorem 3.5.1. \it Let $u_i:\Delta\to(\rr^4, J)$, $i=1, 2$ be two primitive
distinct $J$ - complex  disks  such that $u_1(0)=u_2(0)$. Set $M_i \deff u_i
(\Delta)$. Let $Q = M_1\cap  M_2$ be the intersection set of the disks.
If $J$ is $C^1$-smooth, then the following is true.

\smallskip
\sli The set $\{\, (z_1,z_2) \in \Delta \times \Delta: 
u_1(z_1)= u_2(z_2)\,\}$ is
a discrete subset of $\Delta\times \Delta $. In particular, $u_1(\Delta) \cap
u_2(\Delta)$ is a countable set;

\smallskip
\slii The~intersection index in any such point of $Q$ is strictly positive.
Moreover, if $\mu_1$ and $\mu_2$ are the~multiplicities of $u_1$ and $u_2$ in
$z_1$ and $z_2$, respectively, with $u_1(z_1)=u_2(z_2)=p$, then the
intersection number $\delta_p$  of branches of $M_1$ and $M_2$ at $z_1$ and
$z_2$ is at least $\mu_1\cdot \mu_2$;

\smallskip
\sliii $\delta_p=1$ iff $M_1$ and $M_2$ intersect at $p$ transversally.
\rm

\state Proof.
{\sl Case 1. The map $u_1$ has no critical points}. 
Thus $u_1: \Delta \to \rr^4$ is an
immersion. Take any $(a_1,a_2) \in \Delta \times \Delta$ with $u_1(a_1)=u_2
(a_2)$. In a neighborhood of $u_1(a_1)\in \rr^4$ we find local coordinates
$(w_1, w_2)$ with the properties from {\sl Step 1} of the proof of {\sl
Theorem 3.4.1}. We may assume that in these coordinates the map $u_1$ has
the form $u_1(z)=(z,0)$ and that $a_1$ is the point with $z=0$. Choose the 
local
coordinate $z$ on the second disk such that $z=0$ is the point $a_2$.
Consider representation of the map $u_2: \Delta \to \rr^4$ in coordinates
$(w_1,w_2)$ such that $u_2(z)= (w_1(z),w_2(z))$. From the inequality
(3.4.1) and {\sl Lemma 3.1.1}  we conclude that $w_2(z)=
z^\nu b(z)$ for some $b\in L^{1,p}$ with $b(0)\not =0$. This implies that the
only intersection point of images of small neighborhoods of $a_1$ and $a_2$
is $u_1(a_1)=u_2(a_2)=p$, and its intersection index is $\nu$. This provides
the properties \sli -- \sliii for the case when one of the maps 
$u_1$ or $u_2$ is
an immersion.

\medskip\noindent
{\sl Case 2. Both $u_1$ and $u_2$ are immersions outside $0\in \Delta$, but
$du_1(0)=0=du_2(0)$. Moreover, $u_1(0)=u_2(0)$.}

First we collect some information about the local behavior of $u_i$ at $0\in
\Delta$. We consider only the map $u_1$; the same procedure should be done
for $u_2$.

Without loss of generality we may assume that the point $u_1(0)=u_2(0)$ is
the center of coordinates in $\rr^4$. We may also assume that the standard
complex structure $J\st$ in $\rr^4$ coincides with the given structure $J$ in
the point $P$. Let $(w_1,w_2)$ be the standard complex coordinate in $(\rr^4,
J\st) =\cc^2$. To avoid possible confusion, denote by $\ti u=\ti u(z)$ the
representation of $u$ in coordinates $(w_1,w_2)$. Then $\ti u_1(z)= z^{\mu_1}
v_1(z)$ with $v_1 \in L^{1,p}$ and $v_1(0) \not=0$ by the ``similarity
principle''.

Now consider the bundle $E_1:= u_1^*T\rr^4$ over $\Delta$. Two complex
structures in $\rr^4$, $J$ and $J\st$, induce two complex structures in
$E_1$, $u_1^*J$ and $u_1^*J\st$. We always equip $E_1$ with the co-structure
$u_1^*J$; otherwise the contrary is noted explicitly. In particular, $E_1$
equipped with $u_1^*J$ has the canonical holomorphic structure such that the
homomorphism $du_1: T\Delta \to E_1$ is holomorphic. Consequently, $du_1(z)=
z^{\nu_1} s_1(z)$ for some holomorphic section $s_1$ of $E_1$ with $s_1(0)
\not =0$.

On the other hand, since the tangent bundle $(T\cc^2, J\st)$ is trivial, the
complex bundle $(E_1, u_1^*J\st)$ has  natural trivialization. Denote this
trivialization by $\psi: (E_1, u_1^*J\st) \to \cc^2$, where $\cc^2$ stands
for the trivial bundle over $\Delta$. Then the homomorphism $\psi(z): E_1|_z
\to \cc^2$ is an $\rr$-linear isomorphism in general, and complex linear
exactly for those $z\in\Delta$, for which $J(u_1(z))= J\st(u_1(z))$.
In particular, this holds for $z=0$. This implies $\psi(z^\nu s(z))= z^\nu
\psi(s(z)) + O(|z|^{\nu+1})$. From our construction of $psi$ follows the
equality $\psi \scirc du_1(z) = d\ti u_1(z)$. Thus $d\ti u_1(z)= z^{\nu_1}
\cdot \psi \scirc s_1(z) + O(z^{\nu_1+1})$. Comparing this asymptotic
expansion with $\ti u_1(z)= z^{\mu_1} v_1(z)$ we conclude that $\nu_1=
\mu_1-1$ and $s_1(0)= \mu_1 v_1(0)$.

From now on we do not distinguish the map $u_1$ from its representation
$\ti u_1$ in coordinates $(w_1, w_2)$. Note that we can use the same
coordinates $(w_1, w_2)$ and the structure $J\st$ considering the map $u_2$.
Thus we get the asymptotic relations
$$
\eqalignno{
u_i(z) &= z^{\mu_i} v_i(0) + O(|z|^{\mu_i+ \alpha}) &(3.5.1)
\cr
du_i(z) &= z^{\mu_i-1} v_i(0) + O(|z|^{\mu_i-1+\alpha}) &(3.5.2)
}
$$
with some H\"older exponent $\alpha>0$.
This implies the transversality of small $J$ - complex  disks $u_i(\Delta
(\rho))$ to small spheres $S^3_r$. More precisely, there exist radii $\rho>0$
and $R>0$ such that for any $0<r<R$ the $J$-curves $u_i(\Delta (\rho))$
intersect the sphere $S^3_r\deff \{ |w_1|^2+ |w_2|^2=r^2 \}$ transversally
along smooth immersed circles $\gamma_i(r)$. In fact, the asymptotic relation
(4) provides that for any $\theta \in [0,2\pi]$ there exists at least one
solution of the equation $\bigl|u_i(\rho_{i} e^{\isl\theta}) \bigr| =r$ with
$\rho_i< \rho$, and that any such solution $\rho_i$ must be close to $\left(
{r \over |v_1(0)|}\right)^{1/\mu_i}$. Then one uses (5) to show that the set
$\ti\gamma_i(r) \deff \{ z\in\Delta: \bigl|u_i(z) \bigr| =r \}$ is, in fact,
a smooth imbedded curve in $\Delta$, parameterized by $\theta\in [0,2\pi]$, and
that $u_i: \ti\gamma_i(r) \to S^3_r$ is an immersion with the image
$\gamma_i(r)$.

Taking an appropriate small subdisk and rescaling, we may assume that $\rho=1
=R$. Note that the points of the self- (resp.\ mutual) intersection of
$\gamma_i(r)$ are self-( or resp.\ mutual) intersection points of $u_i
(\Delta)$. Let us call $r \in ]0,1[$ non-exceptional if curves $\gamma_i(r)
\subset S^3_r$ are imbedded and disjoint. Thus $r^* \in ]0,1[$ is exceptional
if $S^3_{r^*}$ contains intersection points of $u_i (\Delta)$.

The ``unique continuation lemma'' of Step 2 of the proof of {\sl Theorem
3.4.1} provides that any such
intersection point is isolated in the punctured ball $\check B \deff \{ 0<
|w_1|^2 + |w_2|^2 <1 \}$. This implies that either there exist finitely many
exceptional radii $r^*\in ]0,1[$, or that they form a 
sequence $r^*_n$ converging
to $0$.

Denote $M_i(r) \deff u_i (\Delta) \cap B(r)$. For non-exceptional $r$ we can
correctly define the intersection index of $M_1(r)$ with $M_2(r)$.

By {\sl Corollary 3.1.3}, the~maps $u_i$ can be represented in the~form
$u_i(z)=z^{\mu_i}\cdot v_i(z)$ with $v_i\in L^{1,p}(\Delta,\cc^2)$ such
that $v_1(0) \not= 0 \not= v_2(0)$. As above, we consider two cases.

\smallskip\noindent
{\sl Case 1. The~vectors $v_1(0)$ and $v_2(0)$ are not collinear.}

It is easy
to see that, in this case, $0\in\cc^2$ is an~isolated intersection point of
$u_1(\Delta)$ and $u_2(\Delta)$ with multiplicity exactly $\mu_1\cdot \mu_2$.
The asymptotic formula (3.4.4) provides that if vectors $v_1(0)$ and $v_2(0)$ are
complex linear independent, then $0\in \rr^4$ is an isolated intersection
point of $Q =u_1(\Delta) \cap u_2(\Delta)$ with the index $\mu_1 {\cdot}
\mu_2$. Thus we have only finitely many intersection points $p\in Q$.
Since all other points $p\in Q$ are smooth, the intersection index 
in every such
point is positive. Thus in the case of non-collinear $v_1(0)$ and $v_2(0)$
for any non-exceptional $r>0$ the intersection index of $M_1(r)$ and $M_2(r)$
is positive.

{\sl Case 2. The~second case is when the~vectors $v_1(0)$ and $v_2(0)$ are
collinear.}

The idea is to ``turn'' the~map $u_2$ a little and to reduce
the~case to the~previous one. So let $T\in \endc(\cc^2)$ be a~linear
unitary map which is close enough to identity such that $T(v_2(0))$ is not
collinear to $v_1(0)$. Define $J^T\deff T\inv \scirc J \scirc T$ so that
$\Vert J^T -J \Vert_{C^1(B)} \le \Vert T - \id \Vert$. Applying {\sl
Lemma 3.3.1} with $u_0=z^{\mu _2}v_2$ and $v=0$, we find $w\in L^{1,p}(
\Delta,\cc^2)$ with $w(0)=0$ such that $\dbar_{J^T} (z^{\mu_2}(v_2+w)) =0$.
The~map $\tilde u_2\deff T(z^{\mu_2}(v_2+w))$ is the required ``turned''
$J$-holomorphic map. Since such a ``turn'' can be made small enough,
the~intersection number of $u_i(\Delta)$ in the~ball $B_{r^-}$ does not
change.

This implies that the~intersection number of $u_1(\Delta)$ and
$u_2(\Delta)$ in any ball $B_r$ is not less than $\mu_1\cdot \mu_2$. Another
conclusion is that $0\in\cc^2$ is an isolated intersection point of
$u_j(\Delta)$. Otherwise we could find a sequence $r_i\msmall{\searrow}0$ with
at least one intersection point of $u_j(\Delta)$ in every spherical shell
$B_{r_i} \!\backslash B_{r_{i+1}}$; therefore the~intersection number of
$u_j(\Delta)$ in the~balls $B_{r_i}$ would be strictly decreasing in $i$.

Thus the statements (i) and (ii) are proved. The proof of (iii) is now obvious
and follows from the observation that $\mu _1 = \mu _2 =1$ in this case.
\qed

\smallskip
\state Corollary 3.5.2. {\it Let $u_i:S_i\to(X, J)$, $i=1, 2$ be compact
irreducible $J$ - complex  curves such that $u_1(S_1)=M_1\not=u_2(S_2)=M_2$.
Then
they have finitely many intersection points and the~intersection index in any
such point is strictly positive. Moreover, if $\mu_1$ and $\mu_2$ are the~
multiplicities of $u_1$ and $u_2$ in such a~point $p$, then the~intersection
number of $M_1$ and $M_2$ in $p$ is at least $\mu_1\cdot \mu_2$.
}

%%%%%%%%%%%%%%%%%%%%%%%%%%%%page 34%%%%%%%%%%%%%%%%%%%%%%%%%%%%

\newpage
\break

\bigskip
\noindent
{\bigbf Appendix 2}

\smallskip\noindent
{\bigbf The Bennequin Index and Genus Formula.}

\smallskip\noindent
\sl A2.1. Bennequin Index of a Cusp.

\smallskip
\rm Let $u : (\Delta , 0)\to (\cc^2, J, 0) $ be a germ of a nonconstant
 complex  curve in zero. Without loss of generality we always suppose
that $J(0)=J\st$. Taking into account that zeros of $du$ are isolated, we
can suppose that $du$ vanishes only at zero. Furthermore, let $w_1, w_2$ be
the standard complex coordinates in $(\cc^2, J\st)$. According to 
{\sl Lemma 3.1.3} we can write our curve in the form

\smallskip
$$
u(z) = z^\mu\cdot a + o(\vert z\vert^{\mu +\alpha}),
\qquad a\in \cc^2 \setminus \{ 0\},\> 0<\alpha<1. \eqno(A2.1.1)
$$

\smallskip

For $r>0$ define $F_r\deff TS^3_r \cap J(TS^3_r)$ to be the distribution
of $J$-complex planes in the tangent bundle $TS^3_r$ to the sphere of radius
$r$. $F_r$ is trivial, because $J$ is homotopic to $J\st = J(0)$. By $F$
we denote the distribution $\cup_{r>0}F_r\subset \cup_{r>0}TS^3_r \subset 
TB^*$, where  $TB^*$ is the tangent bundle to the punctured ball in $\cc^2$.

\smallskip
\noindent
\bf Lemma A2.1.1. \it The (possibly reducible) curve $\gamma_r = M\cap S^3_r$ is
transversal  to $F_r$ for all sufficiently small $r>0$.

\smallskip\noindent
\bf Proof. \rm In fact, due to  {\sl Lemma 3.2.3} one has $du(z) = \mu
z^{\mu-1} \cdot a\,dz + o(\vert z\vert^{\mu -1+\alpha})$. Since
$J\approx J\st$
for $r$ sufficiently small, $T\gamma_r$ is close to $J\st  n_r$, \rm
where $ n_r$ \rm  is the field of normal vectors to $S^3_r$.

On the other hand, for sufficiently small $r$, the distribution $F_r$ is
close to the one of $J\st$ complex planes in $TS^3_r$, which is orthogonal
to $J\st n_r$.
\qed

\smallskip
This fact permits us to define the Bennequin index of $\gamma_r$. Namely, take
any nonvanishing section $\vec v$ of $F_r$ and move $\gamma_r$ along the vector
field $v$ to obtain a curve $\gamma'_r$. We can make this move for a small
enough time, so that $\gamma'_r $ does not intersect $\gamma_r$. Following
Bennequin [Bn], we have

\state Definition A2.1.1. Define the \sl Bennequin index $b(\gamma_r)$ \rm to
be the linking number of $\gamma_r$ and $\gamma'_r$.

\smallskip
\rm
This number does not depend on $r>0$, taken sufficiently small, because
$\gamma_r$ is homotopic to $\gamma_{r_1}$ for $r_1<r$ within the curves
transversal to $F$, see [Bn]. It is also independent of the particular choice
of the field $\vec v$. Thus for the standard complex structure $J\st$ in
$B\subset \cc^2$ we use $\vec v\st(w_1,w_2)\deff  (-\bar w_2, \bar w_1)$
for calculating the Bennequin index of the curves
on sufficiently small spheres. For an arbitrary almost complex structure $J$
with $J(0)=J\st$ we can find the vector field $\vec v_J$, which is defined
in a small punctured neighborhood of $0\in B$, is a small perturbation
of $\vec v\st$, and lies in the~distribution $F$ defined by $J$.

Denote by $B_{r_1, r_2}$ the spherical shell $B_{r_2}\setminus B_{r_1}$ for
$r_1<r_2$.

\state Lemma A2.1.2. \it Let $\Gamma $ be an immersed $J$ - complex  curve 
in the neighborhood of $\overline B_{r_1, r_2}$ such that all self 
intersection points of $\Gamma $  are contained in $B_{r_1, r_2}$ and for 
every $r_1\le r\le r_2$ all components of the curve $\gamma_r = \Gamma \cap 
S^3_r$ are transversal to $F_r$. Then
$$
b(\gamma_{r_2}) = b(\gamma_{r_1}) + 2\cdot \sum\nolimits
_{x\in Sing(\Gamma) }\delta_x,
\eqno(A2.1.2)
$$
where the sum is taken over self-intersection points of $\Gamma $.

\state Proof. \rm Move $\Gamma $ a little along $v_J$ to obtain
$\Gamma^\varepsilon$. By $\gamma^\varepsilon_{r_1},
\gamma^\varepsilon_{r_2}$ denote the intersections $\Gamma^\varepsilon \cap
S^3_{r_1},\Gamma^\epsi\cap S^3_{r_2}$, which are of course the moves
of $\gamma_{r_j}$ along $v_J$. We have  $l(\gamma_{r_2},\gamma
_{r_1}^\epsi) - l(\gamma_{r_1},\gamma_{r_1}^ \epsi) =
{\sf int}(\Gamma ,\Gamma^\epsi)$, where $l(\cdot ,\cdot )$ is the~linking
number and ${\sf int}(\cdot ,\cdot )$ is the~intersection number.

Now let us calculate $int(\Gamma ,\Gamma^\epsi)$. From {\sl Theorem
3.5.1} we know that there are only a finite number $\{ p_1,\ldots, p_N\} $ of
self-intersection points of $\Gamma $. Take one of them, say $p_1$. Let $M_1,
\ldots, M_d$ be the disks on $\Gamma $ with a common point $p_1$ and otherwise
mutually disjoint. More precisely we take $M_j$ to be irreducible
components of $\Gamma \cap B_{p_1}(\rho)$ for $\rho >0$ small enough.
Remark that $M_j$ are transversal to $v_J$, so their moves $M_j^\varepsilon$
don't intersect them, {\sl \i.e.,} $M_j\cap M_j^\epsi = \emptyset $.
Note also that ${\sf int}(M_k, M_j) = {\sf int}(M_k, M^\epsi_j)$
for $\varepsilon >0$ sufficiently small. So $int(M_k, M_j) =
int(M_k, M^\epsi_j) + int(M^\epsi_k, M_j)$.
Finally $\delta_{p_1} = \sum_{1\le k<j\le d} int(M_k, M_j) =
int(\Gamma \cap B_{p_1}(\rho ),\Gamma^\epsi\cap
B_{p_1}(\rho ))$. This means that $int(\Gamma ,\Gamma^\epsi)
= 2\cdot \sum_{j=1}^N\delta_{p_j}$.
\qed

\medskip\noindent
\sl A2.2. Proof of Adjunction Formula.

\nobreak\rm
In $\S \S 3.4, 3.5$ we have proved that compact $J$ - complex  curve 
with a finite number
of irreducible components $M = \bigcup_{i=1}^d M_i$ has only 
a finite number of
nodes ({\sl \i.e.,} self-intersection points) points, provided $J$ is of class
$C^1$.

For each such point $p$ we can introduce, according to {\sl Definition 3.5.1},
\rm the~self-in\-ter\-sec\-tion
number $\delta_p(M)$ of $M$ at $p$. Namely, let
$S_j$ be a parameter curve for $M_j$, {\sl \ie} $M_j$ is given as an image of
the $J$-holomorphic map $u_j : S_j\to M_j$. We always suppose that the
parameterization $u_j$ is primitive, {\sl \i.e.,} they cannot be 
decomposed like
$u_j=v_j\scirc r$ where $r$ is a nontrivial covering of $S_j$ by another
Riemann surface. Denote by $\{ x_1,\ldots, x_N\}$ the set of all pre-images
of $p$ under $u : \bigsqcup_{i=1}^d S_j\to X$, and take mutually disjoint
disks $\{ D_1,\ldots, D_N\}$ with centers $x_1,\ldots, x_N$ such that their
images have no other common points different from $p$. For each pair $D_i,
D_j$, $i\not= j$, define an intersection number as in {\sl Definition 3.5.1}
and take the sum over all different pairs to obtain $\delta_p(M)$.

Now put $\delta = \sum_{p\in D(M)}\delta_p(M)$, where the sum is taken over the
set $D(M)$ of all nodes of $M$, {\sl \i.e.,} points which have at least two
pre-images.

Consider now the set $\{ p_1,\ldots, p_L\} \subset \bigcup_{j=1}^d S_j = S$
of all cusps of $M$, {\sl \i.e.,} points where the differential of the
appropriate parameterization vanishes. Take a small ball $B_r$ around 
$u(p_i)$ and a small disk $D_{p_i}$ around $p_i$.
Let $\gamma^i_r\deff u(\Delta_{p_i}) \cap \partial B_r$ and $b_i$ be
the~Bennequin index of $\gamma^i_r$, defined in {\sl Definition A2.1.1}.

\state Definition A2.2.1. The~number $\varkappa_i\deff (b_i +1)/2$
is called the~{\sl conductor of the~cusp $a_i=u(p_i)$}.

\smallskip Also let $\varkappa\deff \sum_{i=1}^L \varkappa_i$. We are now
able to state the genus formula.

\medskip
\noindent
\bf Theorem A2.2.1. \it Let $M =\cup_{j=1}^d M_j$ be a compact $J$ - complex 
curve in an almost complex surface $(X, J)$ with the~distinct irreducible
components $\{M_j\}$, where $J$ is of class $L^{1,p}$. Then

\smallskip
$$
\sum_{j=1}^d g_j = {[M]^2 - c_1(X, J)[M]\over2} + d -\delta - \varkappa ,
\eqno(A2.2.1)
$$

\smallskip
\noindent
where $g_j$ are the genera of parameter curves $S_j$.

\state Proof.
\rm The main line of the proof is the reduction of a general case to
the case where $M$ is immersed, which was proved in $\S 1$.

Let $u: \bigsqcup_{j=1}^d S_j \longrightarrow X$ be a $J$-holomorphic map,
which parameterizes the curve $M$. Also let $\{x_1,\ldots,x_n\}$ be the
set of cusp-points of $M$, {\sl \i.e.,} the images of critical points of $u$.
Our reduction procedure is local, and we make our constructions in a
neighborhood of every point $x_j$ separately. Therefore, from now on 
we fix such
a point $x$. Due to  {\sl Corollaries 2.2.3} and {\sl 3.4.2} there exists
a neighborhood $U$ of $x$ which contains no other cusp-points and no other
self-intersection points.  {\sl Theorem 3.5.1 i)}
implies that by taking the neighborhood $U$ small
enough, we may assume that any component of $M\cap U$ goes through $p\in U$.
Without losing generality we may also assume that $U$ is the unit ball $B
\subset \cc^2$, that $x$ corresponds to the center $0$ of $B$, and $J(0)=J\st$.

Denote by $\Gamma_j$ the irreducible components of $M \cap B$ and let
$u_j :\Delta \to B$ be a parameterization of $\Gamma_j$ such that $u_j(0)=0
\in B$. Denote $\mu_j \deff \ord_0 du_j +1$; thus $0\in \Delta$ is 
the cusp-point
for $u_j$ \iff $\mu_j\ge2$.

\smallskip\noindent
{\sl Step 1. Rescaling procedure}.

\nobreak
Take some cusp-component $\Gamma_j$. Due to {\sl Corollary  3.1.3} the map
$u_j$ has the form $u_j(z)= z^{\mu_j} \cdot a_j + O\bigl(|z|^{\mu_j+\alpha}
\bigr)$ with the constant $a_j\not=0\in\cc^2$. Moreover, {\sl Lemma 3.2.3}
provides that
$$
\bigl\vert du_j(z) - \mu_j z^{\mu_j-1} a_j \cdot dz \bigr\vert
\le C(\Vert u_j \Vert_{L^{1,p}(\Delta)}, \alpha, p) \cdot
|z|^{\mu_j-1+\alpha},
$$
for any $0<\alpha<1$ and $2<p<\infty$.

For $0<r\le1$ we consider the maps $\pi_r:B\to B$, $\pi_r(w)\deff r^{\mu_j}
\cdot w$, the rescaled maps $u^{(r)}_j:\Delta \to B$, $\pi_r\scirc
u^{(r)}_j(z) = u_j (rz)$, and the rescaled almost complex structures $J^{(r)}
\deff \pi_r^*J$ in $B$. The map $u^{(r)}_j$ is a parameterization of
$J^{(r)}$-holomorphic curves $\pi_r^{-1}\Gamma_j$. One can see that
$\Vert u^{(r)}_j(z) - z^{\mu_j} a_j \Vert_{C^0(\Delta)} \le C\cdot
r^{\mu_j+\alpha}$ and $\Vert du^{(r)}_j(z) - \mu_j\cdot z^{\mu_j-1} a_j \cdot
dz \Vert_{C^0(\Delta)} \le C\cdot r^{\mu_j-1+\alpha}$. In particular, there
exists $r_j>0$ such that for all $0<r<r_j$ the maps $u^{(r)}$ are transversal
to all spheres $S^3_1$. On the other hand, $\Vert J^{(r)} - J\st \Vert_{C^1(B)}
= O(r^{\mu_j})$, and for any $\epsi>0$ we can choose sufficiently small
$r_j>0$ such that $\Vert J^{(r)} - J\st \Vert_{C^1(B)} \le \epsi$ for all
$0<r<r_j$. Thus, by replacing $B$, $J$ and all the maps $u_j$ by their rescaled
counterpart, we may assume that the almost complex structure $J$ satisfies
the estimate
$$
\Vert J - J\st \Vert_{C^1(B)} \le \epsi.
\eqno(A2.2.2)
$$

\medskip\noindent
{\sl Step 2. Reduction to the case of holomorphic cusp-points}.

\nobreak
Recall that in  {\sl Lemma 1.4.1} a natural diffeomorphism
between the space $\jj$ of orthogonal complex structures in $\rr^4$ and
the unit sphere $S^2\deff \{\,(c_1, c_2, s)\in \rr^3 : c_1^2 + c_2^2
+ s^2 =1\, \}$ was established. Define the function $\Phi :\jj \to \rr^2$ by setting $\Phi(J)
=(c_1, c_2)$. The map $\Phi$ defines a diffeomorphism between the upper
half-sphere in $S^2$ and the unit disk $\rr^2$ such that the north pole
$(0,0,1)\in S^2$ corresponds to the center of the disk.

Define the function $f:B \to \rr^2$ setting $f(w)\deff \Phi(J(w))$. Due to
$(A2.2.3)$, $f$ parameterizes the given almost complex structure $J$,
$\Vert f \Vert_{C^1(B)} \le C\epsi$, and $f(0)=0$. Fix a~cut-off function
$\chi$ in $B$ such that $0\le\chi\le1$, $\chi\bigm|_{B(1/2)}\equiv 1$,
$\supp \chi \comp B$, and $\Vert d\chi \Vert_{C^1(B)} \le 3$. For $0\le t \le1$
and $0<\sigma<1$ set $f_{\sigma, t}(w) \deff (1-t\chi(w/\sigma))\cdot f(w)$
and define $J_{\sigma, t}\deff \Phi^{-1}(f_{\sigma, t})$. One can easily see
that $\Vert f_{\sigma, t} \Vert_{C^1(B)}  \le 4\cdot\Vert f \Vert_{C^1(B)}$,
and consequently $\Vert J_{\sigma, t} - J\st \Vert_{C^1(B)} \le C_1\cdot
\Vert J - J\st \Vert_{C^1(B)}\le C_1\cdot \epsi$. Here the constant $C_1$,
as well as the constants $C_2,\ldots, C_5$ below in the proof, are
independent of $\epsi$, $t$ and $\sigma$. For fixed $\sigma$, the curve
$J_{\sigma, t}$, $0\le t\le1$ is a homotopy between $J\equiv J_{\sigma, 0}$
and an almost complex structure $J_\sigma \deff J_{\sigma, 1}$ such that
$J_{\sigma, t} \equiv J$ in $B\backslash B(\sigma)$ and $J_\sigma \equiv J\st$
in $B({\sigma\over2})$. Moreover, we have
$$
\Vert J_{\sigma, t} - J\st \Vert_{C^1(B)} +
\Vert \msmall{ dJ_{\sigma, t} \over dt} \Vert_{C^1(B)} \le C_2\cdot \epsi.
\eqno(A2.2.3)
$$

\medskip
Now we fix some $p>p'>2$ and set $\alpha\deff {2\over p'} -{2\over p}$.
Due to (4.2.3) we can apply {\sl Lemma 3.3.1} \sli to the map $u_j$ with
$\nu_j=\mu_j$, $v_j\equiv0$, and with the curve of $C^1$-smooth almost
complex structures $J_{\sigma, t}$. As result we obtain a curve of maps
$u_{j, \sigma, t} = u_j + z^{\mu_j}w_{j, \sigma, t}$ with $w_{j, \sigma, t}(0)
=0$ and $\Vert w_{j, \sigma, t} \Vert_{L^{1,p}(\Delta)} \le C_3\cdot \Vert u_j
\Vert_{L^{1,p}(B)}$. The condition $\supp \dot J_{\sigma, t} \subset B(\sigma)$
together with {\sl Corollary 3.1.3} provide that $\supp (\dot J_{\sigma, t}
\scirc du_{j, \sigma, t}) \subset \Delta(\rho)$ with $\rho\sim\sigma^{1/\mu_j}$.
Due to (3.3.8) and the H\"older inequality we obtain
$$
\left\Vert z^{-\mu_j}\cdot \dot J_{\sigma, t} (u_{j, \sigma, t})
\scirc du_{j, \sigma, t} \right\Vert_{L^{p'}_{\vphantom{1}}(\Delta)}
\le C_4 \cdot \sigma^{\alpha/\mu_j}\cdot \Vert u_j  \Vert_{L^{1,p}(\Delta)}^2
\cdot \Vert \dot J_{\sigma, t} \Vert_{C^1(B)},
\eqno(A2.2.4)
$$
and consequently

\smallskip\noindent
$$
\Vert w_{j, \sigma, 1} \Vert_{L^{1,p'}_{\vphantom{1}}(\Delta)}
\le C_5 \cdot \sigma^{\alpha/\mu_j}\cdot \Vert u_j  \Vert_{L^{1,p}(\Delta).}^2
\eqno(A2.2.5)
$$
{\sl Lemma 3.2.3} provides that for $\sigma$ small enough the
$J_\sigma$-holomorphic maps $u_{j, \sigma} \deff u_{j, \sigma, 1}$ are
transversal to all spheres $S^3_r$, $0<r<{1\over2}$, have no self-intersection
points in $B({1\over2}) \backslash B({1\over4})$, and the Bennequin index
of $u_{j, \sigma}(\Delta)\cap S^3_r$, ${1\over4}<r<{1\over2}$, coincides with
the one of $\Gamma_j\cap S^3_r$, $0<r<1$. Moreover, we may match
$u_{j, \sigma}$ to the rest of $M$, changing in an appropriate way the almost
complex structure $J$, see {\sl Lemma A2.3.1}.

\smallskip
As a result, we conclude that for appropriate small enough neighborhoods $U_1
\comp U$ of any cusp-point $p\in M$ there exist a perturbed almost complex
structure $\tilde J$
and $\tilde J$-holomorphic curve $\widetilde M$, which is parameterized by
$\tilde u:
\bigsqcup_{j=1}^d S_j \longrightarrow X$ and has the following properties:

\sli $\tilde J$ and $\widetilde M$ coincide on $X\backslash U$ with $J$ and $M$
correspondingly;

\slii $U_1$ is a ball centered in $p$ and $\widetilde M \cap U_1$ is obtained
from
$M\cap U_1$ by perturbing its components in the above described way;

\sliii the (possibly reducible) curve $\tilde \gamma \deff \d U_1\cap \wt M$
is isotopic to $\gamma \deff \d U_1\cap M$, in particular all the corresponding
components $\tilde \gamma_j$ of $\tilde \gamma$ and $\gamma_j$ of $\gamma$ have
the same Bennequin index, and the linking number $l(\tilde \gamma_i,
\tilde \gamma_j)$ is equal to $l(\gamma_i,\gamma_j)$;

\sliv $\tilde u$ is homotopic to $u$;

\slv the cusp-points of $\widetilde M$ coincide with the ones of $M$ and
$\tilde\delta + \tilde\varkappa= \delta +\varkappa$,

\slvi $\tilde J$ is integrable in a neighborhood.

The last equality of  \slv follows from \sliii due to {\sl Lemma 4.1.2}.
Thus the formulas (4.2.1) for $M$ and $\widetilde M$ are equivalent.

\medskip\noindent
{\sl Step 3. Final reduction to the case of an~immersed curve}.

\nobreak
This step is rather obvious and uses the following fact, shown in [Bn]. Let
$B$
be the unit ball in $\cc^2$ and let $\Gamma_0$ be an irreducible holomorphic
curve in $B$, which is transversal to $\d B$ and is defined as a zero-divisor
of a holomorphic function $f$. Then for any sufficiently small nonzero
$\epsi\in\cc$ the curve $\Gamma_\epsi$, defined as the zero-divisor
of the function $f+\epsi$, are smooth and of the same
genus $g$. Moreover, all $\Gamma_\epsi$ are transversal to $\d B$, and the
Bennequin index of $\gamma_\epsi \deff \Gamma_\epsi \cap \d B$ equals 
$2g-1$. In particular, the conductor of a single cusp-point of a holomorphic
curve in $B$ can be defined as a genus of general small perturbation  to a
smooth curve.

In general, let $M$ be a $J$ - complex  curve in $X$, such that $J$ is
integrable in a neighborhood of every cusp-point of $M$. One can now
see that we can perturb $J$ and $M$ to an almost complex structure
$\tilde J$ and a $\tilde J$-holomorphic curve $\widetilde M$, satisfying
the conditions {\sl i)}--{\sl iv)} from above, and the desired condition

{\sl v$'$)} $\widetilde M$ has no cusp-points, $\sum\tilde g_j =\sum g_j +
\varkappa$, $\tilde \delta =\delta$.
\qed

\medskip
An important corollary of the proof of {\sl Theorem 3} is the estimate
from below of the conductor number of a cusp point.

\state Corollary A2.2.2. \it Let $J$ be an almost complex structure in
the unit ball
$B\subset \cc^2$ with $J(0)=J\st$, and let $u:\Delta \to B$ be a $J$ - complex 
curve with the cusp-point $u(0)=0$. Then $\varkappa_0$ is an integer,
$\varkappa_0 \ge \ord_0du$, or equivalently, for all  sufficiently small $r>0$
the Bennequin index of $\gamma_r\deff u(\Delta) \cap S^3_r$ is odd and
satisfies the inequality
$$
b(\gamma_r)\ge 2\cdot \ord_0du - 1. \eqno(A2.2.6)
$$

\medskip\noindent
\bf Proof. \rm Rescaling $u$ as in the {\sl Step 1} of the proof of {\sl
Theorem 3} we may assume that $\Gamma\deff u(\Delta)$  has no nodes and cusps,
excepting $0\in \Delta$ , and is transversal to all spheres $S^3_r$, $0<r<1$,
so that the Bennequin index $b(\gamma_r)$ is the same and equals 
$2\varkappa_0 -1$. We may also assume that for the rescaled almost complex
structure $J$ the estimate $\Vert J - J\st \Vert_{C^1(B)} \le \epsi$ with
the appropriate $\epsi$ is fulfilled. Applying {\sl Step 2} with
a sufficiently small $\sigma$, we can deform $u$ into a $\tilde J$ holomorphic
map $\tilde u$ which has the following properties:

\smallskip
\sli $\widetilde \Gamma \deff \tilde u(\Delta)$ is transversal to all spheres
$S^3_r$, $0<r<1$;

\smallskip
\slii $\widetilde \Gamma$ has no self-intersection points in $B\backslash
B({1\over2})$, and the Bennequin index of $\tilde\gamma_r \deff \widetilde
\Gamma\cap S^3_r$, ${1\over2}<r<1$, coincides with the one of $\gamma_r$;

\smallskip
\sliii $\tilde J$ coincides with $J$ in $B\backslash B({1\over2})$,
is integrable in the neighborhood of $0\in B$, and $\ord_0 d\tilde u =
\ord_0 du $.

\smallskip
Now the corollary follows from {\sl Lemmas A2.1.2} and the fact that for
{\sl integrable} complex structures the same statement is true, see [Bn].
\qed

\medskip
Another modification of the proof of {\sl Theorem 3} leads to the following
consequence.

\state Lemma A2.2.3 \it Let $J$ be a $C^1(X)$-almost complex
structure on $X$ and $M\subset X$ a compact $J$ holomorphic curve
parameterized by $u: S \to X$. Then

\sli $u$ can be $L^{2,p}$-approximate by $J_n$-holomorphic immersions
$u_n: S \to X$ with $J_n \longrightarrow J$ in $C^1(X)$.

\slii there exists a $C^1(X)$-approximation $J_n$ of $J$ and a sequence of
$J_n$-holomorphic imbedded curves $M_n$ which converge to $M$ in the Gromov
topology.
\rm

\state Proof. We do not need this result for the purpose of this paper; 
therefore, we give only a sketch of the proof.

In the first step, one applies the rescaling procedures in order to find
appropriate small neighborhoods of the cusp-points of $M$.

In the second step, one applies {\sl Lemma 3.1.3} to the chosen
neighborhoods, taking $\nu=1$, $v$ sufficiently small,
and $J$ unchanged. After matching procedure we obtain the statement {\sl i)}.
%%%%%%%%%%%%%%%%%%%%%%%%%page 39 question%%%%%%%%%%%%%%%%%
To obtain the statement \slii, one must first deform all the nodes
of $M$ into simple transversal ones and then find an appropriate small
neighborhood $U$ of every node. In some complex coordinates $(w_1, w_2)$
in $U$, the curve $M$ is defined by the equation $w_1\cdot w_2=0$. It remains
to replace a node $M\cap U$ by a ``small handle'' $M_\epsi\deff \bigl\{\,
(w_1, w_2) \in U  :  w_1\cdot w_2=\epsi\,\}$ with $\epsi$ sufficiently
small and to use the matching procedure once more.
\qed

%\end

\bigskip\noindent
\sl A2.3. A matching Structures Lemma.

\smallskip\nobreak
\rm Let $B(r)$ be a ball of radius $r$ in $\rr^4$ centered at zero, and $J$
a $C^1$-smooth almost-complex structure on $B(2)$, $J(0)=J\st$.
Further, let $M = u(\Delta )$ be a closed primitive $J$ - complex  disk
in $B(2)$ such that $u(0)=0$ and $M$ transversally meet $S^3_r$ for
$r\ge 1/2$. Here $S^3_r = \d B(r) $ and transversality are understood with
respect to both $TS^3_r$ and $F_r$, see {\sl paragraph 4.1.}

By $B(r_1,r_2)$ we shall denote the spherical shell $\{ x\in \rr^4:r_1<\Vert
x\Vert <r_2 \}$. In the lemma below denote by $D_{\delta }$ the pre-image of
$B(1+\delta )$ by $u$.

\state Lemma A2.3.1. \it For any positive $\delta >0$ there exists an $\epsi >0$ such
that if an almost complex structure $\tilde J$ in $B(1+\delta )$ and a closed
$\tilde J$-holomorphic curve $\tilde u : D_{\delta }\to B(1+\delta )$ satisfy
$\Vert \tilde J - J\Vert_{C^1(\bar B(1+\delta ))} < \epsi $ and
$\Vert \tilde u - u\Vert_{L^{1,p}(\bar D_{\delta })} < \epsi $, then
there exists an almost-complex structure $J_1$ in $B(2)$ and $J_1$-holomorphic
disk $M_1$ in $B(2)$ such that:

a) $J_1\ogran_{B(1-\delta )} = \tilde J\ogran_{B(1-\delta )}$ and
$J_1\ogran_{B(1+\delta ,2)} = J\ogran_{B(1+\delta ,2)}$.

b) $M_1\ogran_{B(1-\delta )} = \tilde u(D_{\delta })\cap B(1-\delta )$ and
$M_1\cap B(1+\delta ,2) = M\cap B(1+\delta ,2)$.

\state Proof. {\rm We have chosen the parameterization of $M$ to be primitive.
Thus, $u$
is an imbedding on $D_{-\delta ,\delta } = u^{-1}(B_{1-\delta ,1+\delta })$.
Let us identify a neighborhood $V$ of $u(D_{-\delta ,\delta })$ in $B_{1-
\delta ,1+\delta }$ with the neighborhood of the zero-section in the 
normal bundle
$N$ to $u(D_{-\delta ,\delta })$. Now $\tilde u\ogran_{D_{-\delta ,\delta
}}$
can be viewed as a section of $N$ over $u(D_{-\delta ,\delta })$ which is small
{\sl \i.e.,}
contained in $V$. Using an appropriate smooth function $\phi $ on
$D_{-\delta ,\delta }$ (or equivalently on $u(D_{-\delta ,\delta} )$), $\phi
\ogran_{B(1-\delta )\cap D_{\delta }}\equiv 1$, $\phi\ogran_{\partial D_{\delta
}
}\equiv 0$, $0\le \phi \le 1$ we can glue $\tilde u$ and $u$ to obtain a \sl
symplectic \rm surface $M_1$ which satisfies (b).

Patching $J$ and $\tilde J$ and simultaneously making $M_1$  complex 
can be done in an obvious way.
}
\qed

\newpage

%%%%%%%%%%%%%%%%%%%%%%%%%%%%%%%%%%%%%%%%%%%%%%%%%%%%%%%%%%%%%%%%%%%%

%%ch.2
\long\def\comment#1\endcomment{}

\noindent
{\bigbf Chapter II. Compactness Theorem.}

\smallskip\rm
The goal of this chapter is to give a proof of the Gromov Compactness 
Theorem for
continuously varying almost complex structures and for the sequences of
complex curves parameterized by some fixed real surface.

More precisely, we consider a sequence $J_n$ of continuous (\ie of class
$C^0$) almost complex structures on a manifold $X$ which converge uniformly
to some $J_\infty$, again of class $C^0$. Further, let $(C_n,j_n)$ be a
sequence of Riemann surfaces with boundaries of fixed topological type. This
means that each $(C_n,\d C_n)$ can be parameterized by the same real surface
$(\Sigma, \d\Sigma ) $ (see \S\.2 for details). Denote by $\delta_n:\Sigma
\to C_n$ some parameterizations. However, the complex structures $j_n$ on $C_n$
may vary in an arbitrary way. Finally, let some sequence of
$(j_n,J_n)$-holomorphic maps $u_n:C_n\to X$ be given.

\medskip
\state Theorem 2.1. {\it If the areas of $u_n(C_n)$ are uniformly bounded (with
respect to some fixed Riemannian metric on $X$) and the structures $j_n$ do
not degenerate at the boundary (see {\sl Definition 1.7}), then there exists
a subsequence, denoted $(C_n, u_n)$, such that

\smallskip
\sl1) \it $(C_n,j_n)$ converge to some nodal curve $(C_\infty,j_\infty)$ in
an appropriate completion of the moduli space of Riemann surfaces of given
topological type, \i.e., there exists a parameterization map $\sigma_\infty:
\Sigma \to C_\infty$ by the same real surface $\Sigma$;

\smallskip
\sl 2) \it one can choose a new parameterizations $\sigma_n$ of $C_n$ in such a
way that $\sigma_n^*j_{C_n}$ will converge to $\sigma_\infty^*j_{C_\infty}$
in the $C^\infty$-topology on compact subsets outside of the finite set of
circles on $\Sigma$, which are pre-images of the nodal points of $C_\infty$
by $\sigma_\infty$, where $j_{C_n}$ denote the complex structures on $C_n$;

{\sl 3)} the maps $u_n\scirc \sigma_n$ converge, in the $C^0$-topology on 
$\Sigma$ and in the $L^{1,p}_\loc$-topology (for all $p<\infty$)
outside of the pre-images of the nodes of $C_\infty$ to a map $u_\infty \scirc
\sigma_\infty$ such that $u_\infty$ is a $(j_\infty, J_\infty) $-holomorphic
map $C_\infty\to X$.
}

\smallskip\rm
This description of convergence is precisely the one given by Gromov in [G].
Our statement is slightly more general in two directions. First, we consider
not only the case of closed curves, but also the case where $C_n$ is open and
of a ``fixed topological type''. Second, we note that the Gromov compactness
theorem is still valid for continuous and continuously varying almost complex
structures. This could have an interesting application, because now one can
consider $C^0$-perturbations of complex structures being assured that at least
the compactness theorem still holds true. For the definitions 
involved and the formal
statement see \S4.1 and {\sl Theorem 4.1.1}.

\smallskip
We also prove the compactness theorem for curves with boundaries on a totally
real submanifold. This ``boundary'' result needs appropriate generalizations
of all ``inner'' constructions and estimates. The related considerations are
shown in Appendix III.

\medskip
Another result of this paper, which we would like to mention in the
introduction, is the improvement of the removable singularity theorem in two
directions.

\smallskip
First we prove (see {\sl Corollary 5.2.1}) the following generalization of the
removability theorem for the point singularity.

\smallskip
\state Theorem 2.2. {\it If the area of the image of $J$-holomorphic map
$u:(\check\Delta, J\st)\to (X,J)$ from the punctured disk into a compact almost
complex manifold ``is not growing too fast", \i.e., if $\area(u(R_k))\le \eps $
for all annuli $R_k \deff \{ z\in \cc :{1\over e^{k+1}}\le | z| \le {1 \over
e^k}\}$ with $k>\!>1$, then $u$ extends to the origin.
}

\smallskip
The positive constant $\eps $ here depends on the Hermitian
structure $(J,h)$ of $X$. This theorem, under the stronger assumption $\sum_k
\area(u(R_k)) \equiv \norm{du}^2_{L^2(\check\Delta)} < \infty$, was proved by
Sacks and Uhlen\-beck [S-U] for harmonic maps, and by Gromov [G] for
$J$-holomorphic maps.

\smallskip
This fact (which is proved here for continuous $J$'s) is new even in the
integrable case. In fact, it measures the ``degree of non-hyperbolicity" (in
the sense of Kobayashi) of $(X,J,h)$.

\smallskip
Another (see {\sl Corollary A3.3.5}) is a generalization of Gromov's result
about removability of the boundary point singularity, see [G]. 
An improvement is
that the statement remains valid also when one has {\sl different} boundary
conditions to the left and to the right of a singular point. Let us explain
this in more detail.

\smallskip
Define the (punctured) half-disk by setting $\Delta^+ \deff \{ z\in \Delta :
\im(z)>0\}$ and $\check\Delta^+ \deff \Delta^+\bs \{0\}$. Define $I_- \deff
]{-1},0[ \subset \d\check\Delta^+$ and $I_+ \deff ]0,+1[ \subset \d\check
\Delta^+$. Let a $J$-holomorphic map $u:(\check\Delta^+, J\st)\to (X,J)$ be
given, where $J$ is again continuous. Suppose further that $u(I_+)\subset
W_+$ and $u(I_-) \subset W_-$, where $W_+,W_-$ are totally real submanifolds
of dimension $n= \half \dim_\rr X$ and intersect transversally.

\smallskip
\state Theorem 2.3. {\it There is an $\eps^b>0$ such that if for all
half-annuli $R^+_k\deff\{ z\in \Delta^+ :e^{-(k+1)} \le | z| \le e^{-k} \}$
one has $\area(u(R^+_k))\le \eps^b$, then $u$ extends to the origin $0\in
\Delta^+$ as an $L^{1,p}$-map for some $p>2$.
}

\smallskip
As in the ``inner'' case, the necessary condition is weaker than the
finiteness of energy. But unlike ``inner" and smooth boundary cases, it is
possible that the map $u$ in the last statement is $L^{1,p}$-regular in the
neighborhood of the ``corner point'' $0\in \Delta^+$ only for some $p>2$. For
example, the map $u(z) = z^\alpha$ with $0<\alpha <1$ satisfies the totally  
real boundary conditions $u(I_+) \subset \rr$, 
$u(I_-) \subset e^{\alpha\pi\isl}\rr$ and is $L^{1,p}$-regular 
only for $p< p^* \deff {2\over1-\alpha}
\cdot$ Note also that by the Sobolev imbedding $L^{1.p}\subset C^{1,\alpha }$
with $\alpha = 1-{2\over p}$, $u$ extends to zero at least continuously. Thus
$u(0)\in W_+ \cap W_-$.

One can see such a point $x$ as a {\sl corner point} for a corresponding
complex curve. A typical example appears in symplectic geometry when
one takes Lagrangian submanifolds as boundary conditions.

\smallskip
The organization of the chapter is the following. In Lecture 4 we present,
for the convenience of the reader, the basic notions concerning the topology
on the space of stable curves and complex structure on the Teichm\"uller
space of Riemann surfaces with boundary. In Lecture 5 we give the necessary
a priori estimates for the inner case, and the proof of {\sl Theorem
2.1}, related to curves with free boundary. This includes the case of closed
curves. In Appendix III we consider curves with totally real boundary conditions,
obtain necessary a priori estimates at a ``totally real boundary'', and prove
the compactness theorem for such curves. In particular, we prove {\sl Theorem
2.3} there.

\newpage
%\bigskip
%\bigskip
%\centerline{\bigbf Table of contents}

%{\parindent=0pt
%\bigskip
%\sl 0. Introduction.

%\smallskip
%1. Stable curves and Gromov topology.

%\smallskip
%2. Complex structure on the space $\ttt_\Gamma$.

%\smallskip
%3. Apriori estimates.

%\smallskip
%4. Compactness  for the curves with free boundary.

%\smallskip
%5. Curves with boundary on totally real submanifolds.
%}

\noindent
{\bigbf Lecture 4}

\smallskip\noindent
{\bigbf Space of Stable Curves}

\medskip\noindent
{\bigsl 4.1. Stable Curves and Gromov Topology.}

\smallskip
Before stating the Gromov compactness theorem, we need to introduce an
appropriate category of  complex  curves. Since the limit of a
sequence smooth curves can be singular,i.e., a cusp-curve in Gromov's
terminology, we need to allow certain types of singularities of curves. On
the other hand, it is desirable to have singularities as simple as possible.

A similar problem appears in looking for a ``good'' compactification of
moduli spaces $\calm_{g,m}$ of abstract complex smooth closed curves of genus
$g$ with $m$ marked points. The Deligne-Mumford compactification $\barr \calm
_{g, m}$, obtained by adding {\sl stable curves}, gives a satisfactory
solution to this problem and suggests a possible way of generalizing to
other situations. In fact, the only singularity type one should allow are
nodes, or nodal points. An appropriate notion for curves in a complex
algebraic manifold $X$ was introduced by Kontsevich in [K]. Our definition of
stable curves over $(X,J)$ is simply a translation of this notion to almost
complex manifolds. The change of terminology from {\sl stable maps} to 
{\sl stable curves over $(X,J)$} is motivated by the fact that we want to
consider our objects as curves rather than maps.

Recall that a {\sl standard node} is the complex analytic set $\cala_0 \deff
\{ (z_1,z_2)\in \Delta^2 : z_1\cdot z_2 =0\}$. A point on a complex curve is
called a {\sl nodal point} if it has a neighborhood biholomorphic to the
standard node.

\state Definition 4.1.1. {\sl A {\it nodal curve} $C$ is a complex analytic
space of pure dimension 1 with only nodal points as singularities}.

In other terminology, nodal curves are called {\sl prestable}. We shall
always suppose that $C$ is connected and has a "finite topology", \i.e., 
$C$ has
finitely many irreducible components, finitely many nodal points, and that
$C$ has a smooth boundary $\d C$ consisting of finitely many smooth circles
$\gamma_i$, such that $\barr C \deff C \cup \d C$ is compact.

\state Definition 4.1.2. 
{\sl We say that a real connected oriented surface with boundary
$(\Sigma, \d\Sigma)$ {\it parameterizes} a complex nodal curve $C$ if there
is a continuous map $\sigma :\barr\Sigma \to \barr C$ such that

\sli if $a\in C$ is a nodal point, then $\gamma_a = \sigma\inv(a)$ is a
smooth imbedded circle in $\Sigma \bs \d \Sigma $, and if $a\not= b$ then
$\gamma_a \cap \gamma_b= \emptyset$;

\slii $\sigma :\barr\Sigma \bs \bigcup_{i=1}^N\gamma_{a_i}\to \barr C \bs \{
a_1,\ldots ,a_N\} $ is a diffeomorphism, where $a_1,\ldots ,a_N$ are the
nodes of $C$.
}

\bigskip
\vbox{\xsize=.54\hsize\nolineskip\rm
\putm[.13][.01]{\gamma_1}%
\putm[.236][.19]{\gamma_2}%
\putm[.487][0.195]{\gamma_3}%
\putm[.62][.20]{\gamma_4}%
\putm[.84][.15]{\gamma_5}%
\putt[1.1][0]{\advance\hsize-1.1\xsize%
\centerline{Fig.~1}\smallskip%\parindent=0pt
Circles $\gamma_1,..., \gamma_5$ are contracted by the parameterization
map $\sigma$ to nodal points  $a_1, \ldots a_5$.
}%
\putm[.56][.31]{\bigg\downarrow\sigma}%
\noindent
\epsfxsize=\xsize\epsfbox{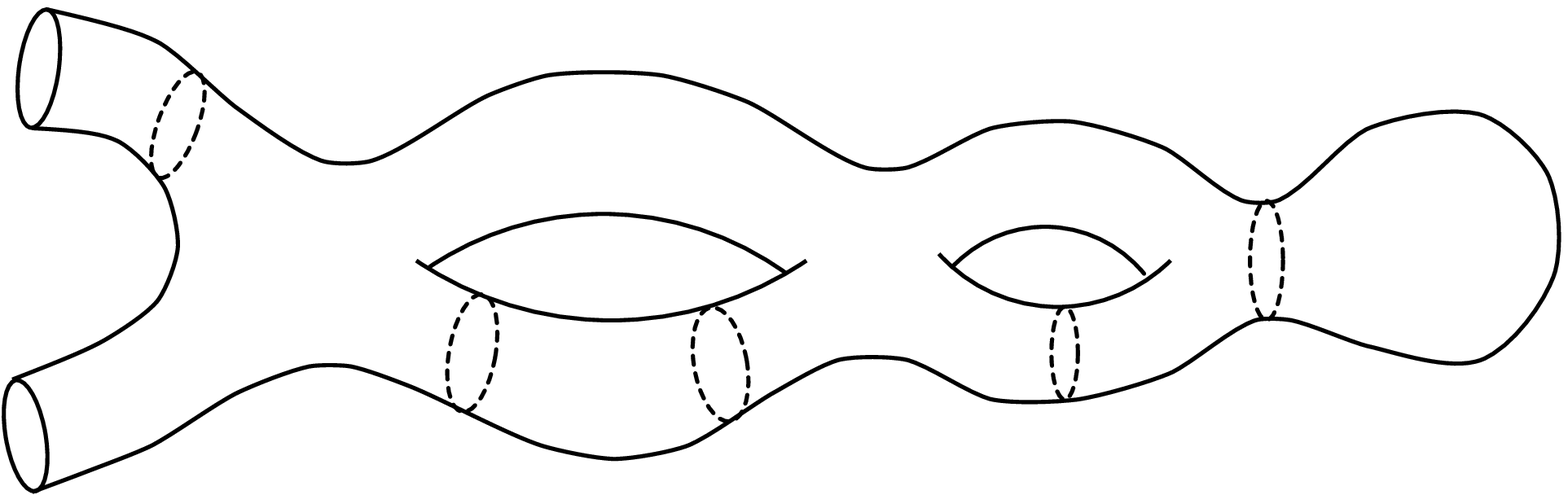}%%
\vskip15pt
%%  second picture
\putm[.14][.045]{a_1}%
\putm[.31][.25]{a_2}%
\putm[.445][0.27]{a_3}%
\putm[.67][.23]{a_4}%
\putm[.805][.20]{a_5}%
\noindent
\epsfxsize=\xsize\epsfbox{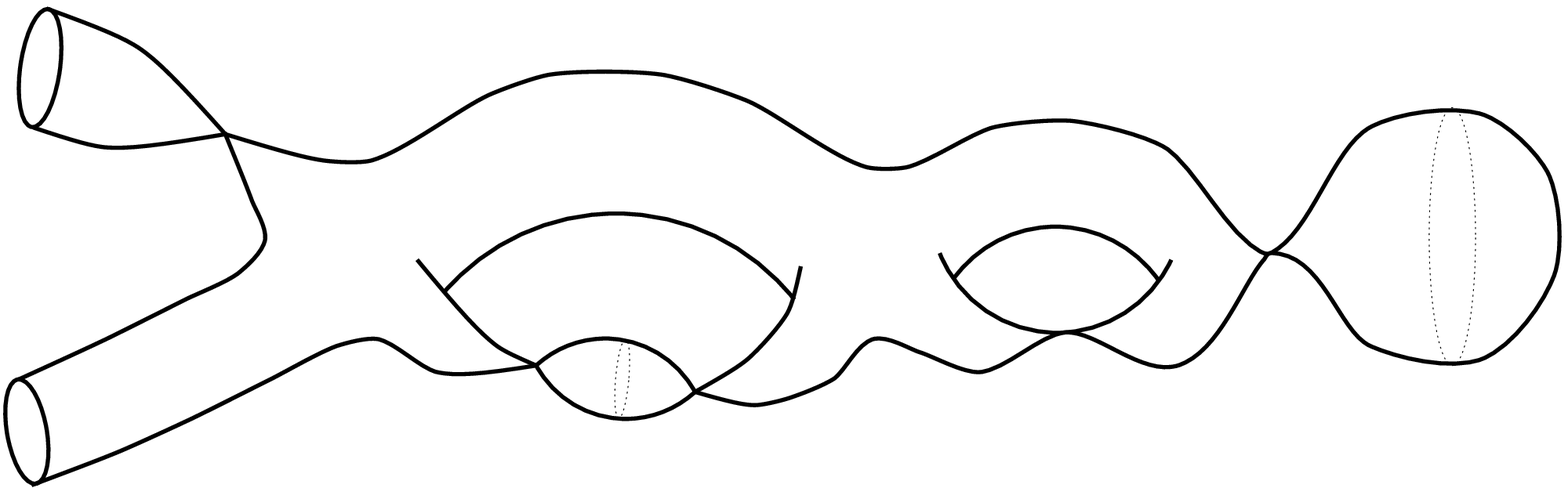}
}
\baselineskip=12.5pt plus 1.5pt

\smallskip
Note that such a parameterization is not unique: if $g:\barr\Sigma \to
\barr\Sigma$ is any orientation preserving diffeomorphism, then $\sigma \scirc
g: \barr\Sigma \to \barr C$ is again a parameterization.

\smallskip
A parameterization of a nodal curve $C$ by a real surface can be considered as
a method of ``smoothing'' of $C$. An alternative method of ``smoothing'' ---
the normalization --- is also useful for our purposes.

Consider the normalization $\hat C$ of $C$. Mark on each component of this
normalization the pre-images (under the normalization map $\pi_C: \hat C \to
C$) of nodal points of $C$. Let $\hat C_i$ be a component of $\hat C$. We can
also obtain $\hat C_i$ by taking an appropriate irreducible component $C_i$,
replacing nodes contained in $C_i$ by pairs of disks with marked points, and
marking remaining nodal points. Since it is convenient to consider components
in this form, we make the following

\state Definition 4.1.3. \it A component $C'$ \sl of a nodal curve $C$ is a
normalization of an irreducible component of $C$ with marked points selected
as above. \rm

\smallskip
This definition allows us to introduce the Sobolev and H\"older spaces of
functions and (continuous) maps of nodal curves. For example, a continuous map
$u: C \to X$ is Sobolev $L_\loc ^{1,p}$-smooth if its restrictions
on every component of $C$ is $L_\loc^{l,2}$-smooth as well. The most 
interesting case is, of course, the one of
continuous $L_\loc ^{1,2}$-smooth maps. In this case the energy functional
$\norm{du}^2_{L^2(C)}$ is defined. The definition of the energy $\norm{du}^2
_{L^2(C)}$ involves Riemannian metrics on $X$ and $C$ which are supposed to
be fixed.

\smallskip
Let $C$ be a nodal curve and $(X,J)$ an almost complex manifold with
continuous almost complex structure $J$.

\nobreak
\state Definition 4.1.4. {\sl A continuous map $u:C \to X$ is $J$-holomorphic
if $u\in L^{1,2}_\loc(C,X)$ and
$$
du_x + J\scirc du_x\scirc j_C=0
\eqno(4.1.1)
$$
for almost all $x\in C$. Here $j_C$ denotes the complex structure on $C$.
}

Recall that the area of a $J$-holomorphic map is defined as
$$
\area(u(C)) \deff \norm{du}^2_{L^2(C)}.
$$

\noindent See the end of \S 1.3.
\smallskip
We shall show later that every $J$-holomorphic $u$ is, in fact,
$L^{1,p}_\loc(C,X)$-smooth  for all $p< \infty$, see {\sl Corollary 2.4.2}.
The following notion of stability was introduced by Kontsevich in [K].

\smallskip
\state Definition 4.1.5. {\sl A {\it stable curve over $(X,J)$} is a pair
$(C,u)$, where $C$ is a nodal curve and $u:C\to X$ is a $J$-holomorphic map,
satisfying the following condition. If $C'$ is a compact component of $C$, 
such that $u$ is constant on $C'$, then there exist finitely many 
biholomorphisms of $C'$ which preserve the marked points of $C$.
}

\state Remark. One can see that the stability condition is 
nontrivial only in the
following cases:

\item{\sl 1)} some component $C'$ is biholomorphic to $\cc\pp^1$ with 1 or 2
marked points; in this case $u$ should be non-constant on any such component
$C'$;

\item{\sl 2)} some irreducible component $C'$ is $\cc\pp^1$ or a torus
without nodal points.

\noindent
Since we consider only connected nodal curves, case {\sl 2)} can happen only
if $C$ is irreducible, so that $C'=C$. In this case $u$ must be non-constant on
$C$.

\smallskip\rm
Now we are going to describe the Gromov topology on the space of stable
curves over $X$ introduced in [G]. Let a sequence $J_n$ of continuous almost
complex structures on $X$ be given. Suppose that $\{J_n\}$ converges to 
$J_\infty$ in the $C^0$-topology. Furthermore, let $(C_n,
u_n)$ be a sequence of stable curves over $(X, J_n)$ such that all $C_n$ are
parameterized by the same real surface $\Sigma$.

\state Definition 4.1.6. {\sl We say that $(C_n,u_n)$ {\it converges to a
stable $J_\infty$-holomorphic curve $(C_\infty,u_\infty)$ over $X$} if the
parameterizations $\sigma_n: \barr\Sigma \to \barr C_n$ and $\sigma_\infty:
\barr \Sigma \to \barr C_\infty$ can be chosen in such a way that the
following hold:

\sli $u_n\scirc \sigma_n$ converges to $u_\infty\scirc \sigma_\infty$ in the
$C^0( \Sigma, X)$-topology;

\slii if $\{ a_k \}$ is the set of nodes of $C_\infty$ and $\{\gamma_k\}$ are
the corresponding circles in $\Sigma$, then on any compact subset $K\comp
\Sigma \bs \cup_k \gamma_k$ the convergence $u_n\scirc \sigma_n \to u_\infty
\scirc \sigma_\infty$ is $L^{1,p}(K, X)$ for all $p< \infty$;

\sliii for any compact subset $K\comp \barr\Sigma \bs \cup_k\gamma_k$ there
exists $n_0=n_0(K)$ such that $ \sigma_n(K) \subset C_n\setminus \{ nodes\} $
for all $n\ge n_0$ and the complex structures $\sigma_n^*j_{C_n}$ converge
smoothly to $\sigma_0^*j_{C_0}$ on $K$;

\sliv the structures $\sigma_n^*j_{C_n}$ are constant in $n$ near the
boundary $\d\Sigma$.
}

\medskip
The reason for introducing the notion of a curve stable over $X$ is similar
to the one for the Gromov topology. We are looking for a completion of the
space of smooth imbedded  complex  curves which has ``nice''
properties, namely: \.1) such a completion should contain the limit of some
subsequence of every sequence of smooth curves, bounded in an appropriate
sense; \.2) the same should also hold for every sequence in the
completed space; \.3) such a limit should be unique. Gromov's compactness
theorem insures us that the space of curves stable over $X$ enjoys these nice
properties.

Condition \sliv is trivial if $\Sigma $ is closed, but it is important when
one considers the ``free boundary case'', \i.e., when $\Sigma$ (and thus all
$C_n$) are not closed and no boundary condition is imposed. However, we
would like to point out that in our approach the ``free boundary case'' is
essentially involved in the proof of the compactness theorem also in the case
of closed curves. On the other hand, in the case of curves with boundary on
totally real submanifolds  such a condition is unnecessary.

\smallskip
Recall that a complex annulus $A$ has a conformal radius $R>1$ if $A$ is
biholomorphic to $A(1,R) \deff \{ z\in \cc \,:\, 1 <|z| < R \}$. An annulus
$A$ is said to be {\sl adjacent} to a circle $\gamma$, if $\gamma$ is one of
its boundary components.

\state Definition 4.1.7. {\sl Let $C_n$ be a sequence of nodal curves,
parameterized by the same real surface $\Sigma$. We say that  complex
structures on $C_n$ {\it do not degenerate near the boundary}, if there exists
$R>1$ such that for any $n$ and any boundary circle $\gamma_{n, i}$ of $C_n$
there exists an annulus $A_{n,i} \subset C_n$ adjacent to $\gamma_{n, i}$
such that all $A_{n,i}$ are mutually disjoint, do not contain nodal points of
$C_n$, and have the same conformal radius $R$.}

Since the  conformal radii of all the $A_{n, i}$ are the same, 
we can identify them
with $A(1,R)$. This means that all changes of complex structures of $C_n$
take place away from the boundary. The condition is trivial 
if $C_n$ and $\Sigma$
are closed, $\d\Sigma = \d C_n = \emptyset$.

\state Remark. Changing our parameterizations $\sigma_n: \Sigma \to C_n$, we
can suppose that for any $i$ the pre-image $\sigma_n\inv (A_{n,i} )$ is the
same annulus $A_i$ independent of $n$.

\medskip
Now we state our main result. Fix some Riemannian metric $h$ on $X$ and some
$h$-complete set $A\subset X$.

\state Theorem 4.1.1. {\it Let $\{(C_n,u_n)\}$ be a sequence of stable
$J_n$-holomorphic curves over $X$ with parameterizations $\delta_n: \Sigma \to
C_n$. Suppose that

\item{\sl a)} $J_n$ are continuous almost complex structures on $X$,
$h$-uniformly converging to $J_\infty$ on $A$ and $u_n(C_n)\subset A$ for all
$n$;

\item{\sl b)} there is a  constant $M$ such that $\area [u_n (C_n)]\le M$
for all $n$;

\item{\sl c)} complex structures on $C_n$ do not degenerate near the boundary.

Then there is a subsequence $(C_{n_k},u_{n_k})$ and parameterizations $\sigma
_{n_k}: \Sigma \to C_{n_k}$ such that $(C_{n_k}, u_{n_k}, \sigma_{n_k})$
converges to a stable $J_\infty$-ho\-lo\-mor\-phic curve $(C_\infty, u_\infty,
\sigma_\infty)$ over $X$.

Moreover, if the structures $\delta_n^*j_{C_n}$ are constant on the fixed
annuli $A_i$, each adjacent to a boundary circle $\gamma_i$ of $\Sigma$, then
the new parameterizations $\sigma_{n_k}$ can be taken equal to $\delta_{n_k}$
on some subannuli $A'_i \subset A_i$, also adjacent to $\gamma_i$.
}

\state Remarks.~1. In the proof, we shall give a precise description of
convergence with estimates in neighborhoods of the contracted circles
$\gamma_i$. The convergence of curves with boundary on totally real
submanifolds will be studied in \S\.5.5.
%%%%%%%%%%%%%%%%%%Question%%%%%%%%%%%%%%%%%%%%%%%%
\state 2. In applications, one uses a generalized version of the Gromov
compactness theorem for nodal curves with a marked point. This version is an
immediate consequence of {\sl Theorem 1.1} due to the following construction.
Consider a nodal curve $C$ and let a $J$-holomorphic map 
$u: C \to X$. Let $\bold
x \deff \{x_1,\ldots,x_m\}$ be the set of marked points on $C$ which are
supposed to be distinct from the nodal points of $C$. Define a new curve
$C^+$ as the union of $C$ with disks $\Delta_1,\ldots, \Delta_m$ such that $C
\cap \Delta_i = \{x_i\}$ and  any $x_i$ becomes a nodal point of
$C^+$. Extend $f$ to a map $f^+: C^+ \to X$ by setting $f^+\ogran_{\Delta_i}$
to be constant and equal to $f(x_i)$. An appropriate definition of
stability, used for triples $(C, \bold x, f)$, is equivalent to stability of
$(C^+, f^+)$. Similarly, the Gromov convergence
$(C_n, \bold x_n, f_n) \to (C_\infty, \bold x_\infty, f_\infty)$
is equivalent to the Gromov convergence
$(C^+_n, f^+_n) \to (C^+ _\infty, f^+_\infty)$. Thus the Gromov compactness for
curves with marked points reduces to the case considered in our paper.
However, we shall consider curves with marked points as well.

\vskip 0pt plus 30pt\ 

\bigskip\noindent
{\bigsl 4.2. Fenchel-Nielsen Coordinates on the Space of Nodal Curves.}

\nobreak
\smallskip\nobreak
In the rest of this section we shall describe topology and conformal geometry
of nodal curves and compute the set of moduli parameterizing deformations of
a complex structure. As a basic reference we use the book of Abikoff [Ab].

\smallskip
Let $C$ be a complex nodal curve parameterized by $\Sigma$.

\state Definition 4.2.1. {\sl A component $C'$ of $C$ is called 
{\it nonstable} if one of the following two cases occurs:

\item{ 1)} $C'$ is $\cc\pp^1$ and has one or two marked points;

\item{ 2)} $C'$ is $\cc\pp^1$ or a torus and has no marked points.
}

\smallskip
This notion of stability of abstract closed curves is due to Deligne-Mumford,
see [D-M]. It was generalized by Kontsevich [K] for the case of maps $f:C
%%%%%%%%%%%%%%%%%%Question%%%%%%%%%%%%%%%%%%%
\to X$, \i.e., for curves over $X$ in our terminology. As was already noted,
the last case can happen only if $C=C'$. Strictly speaking, this case should
be considered separately. However since such considerations require only
obvious changes, we just skip them and suppose that case {\sl 2)} does not
occur.

Our first aim is to analyze the behavior of complex structures in the sequence
$(C_n,u_n)$ of $J_n$-holomorphic curves stable over $X$ with  
uniformly bounded
areas, which are parameterized by the same real surface $\Sigma$. 
At the moment,
the uniform bound on the area of $u_n(C_n)$ is needed only 
to show that the number of
components of $C_n$ is bounded. Passing to a subsequence, we can assume that
all $C_n$ are homeomorphic. This reduces the problem to a description of
complex structures on a fixed nodal curve $C$.

To obtain such a description, it is useful to cut the curve into pieces where
the behavior of a complex structure is easy to understand. Such a procedure is
a {\sl partition into pants} which is well- known in the theory of 
moduli spaces
of complex structures on curves, see, e.g., [Ab], p.93. Here we
shall make use of a related but slightly different procedure. Namely, 
we shall choose a
special covering of $\Sigma$ instead of its partition. Further, as blocks for
our construction we shall use not only pants, but also disks and annuli.
The reasons are that, first, the considered curves can have unstable
components and, second, it is convenient to use annuli for a description of
the deformation of the complex structure on curves. We start with

\state Definition 4.2.2. \it An annulus $A$ \sl on a real surface or on a
complex curve is a domain which is diffeomorphic (resp.\ biholomorphic) to
the standard annulus $A(r,R) \deff \{ z\in \cc \;:\; r<|z|<R \}$ such that
its boundary consists of smoothly imbedded circles. \it Pants \sl(also
called a \it pair of pants\sl) on a real surface or on a complex curve is a
domain which is diffeomorphic to a disk with 2 holes. \rm

The boundary of pants consists of three components, each of them being
either a smoothly imbedded circle or a point. This point can be considered as
a puncture of pants or as a marked point. An annulus or pants is {\sl
adjacent to a circle $\gamma$} if $\gamma$ is one of its boundary components.

\medskip
\vbox{%
\hbox{%
\vtop{%\putm[][]{}
\noindent
\hsize=.46\hsize\epsfxsize=\hsize\epsfbox{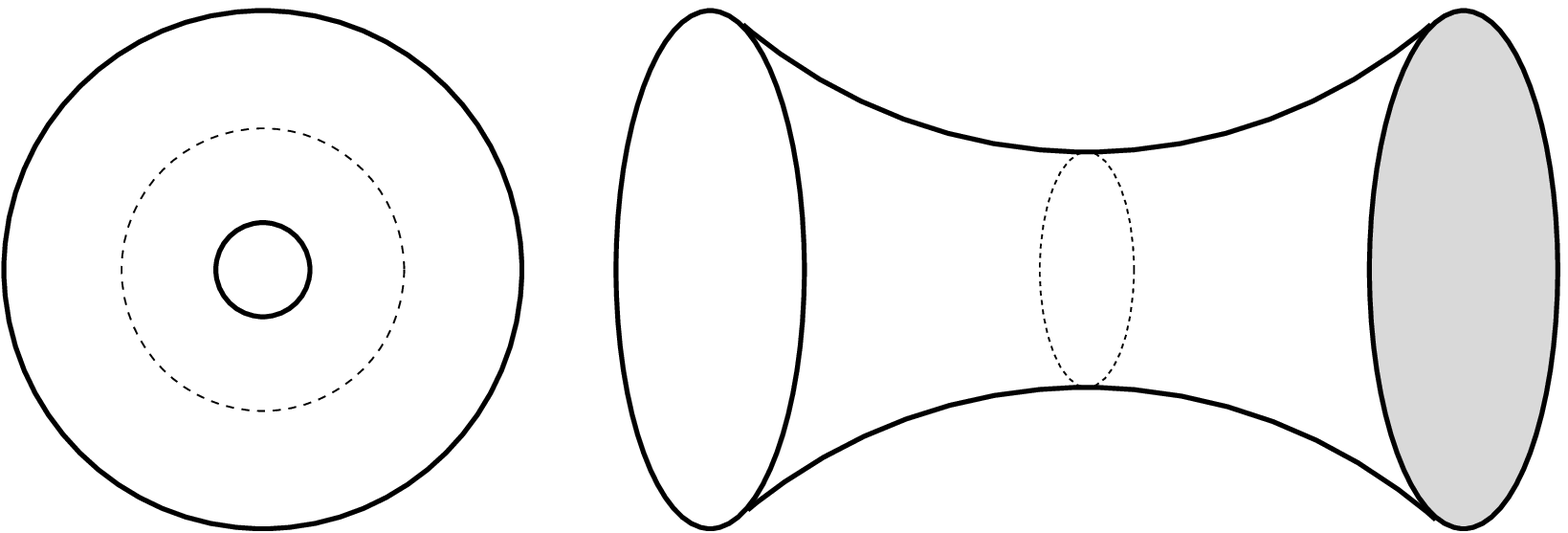}%
\smallskip
\centerline{Fig.~2. An annulus}
\smallskip
It is useful to imagine an annulus as a cylinder. After contracting the
middle circle of the annulus we get a node.
}
\hskip.05\hsize
\vtop{%\putm[][]{}
\noindent
\hsize=.49\hsize\epsfxsize=\hsize\epsfbox{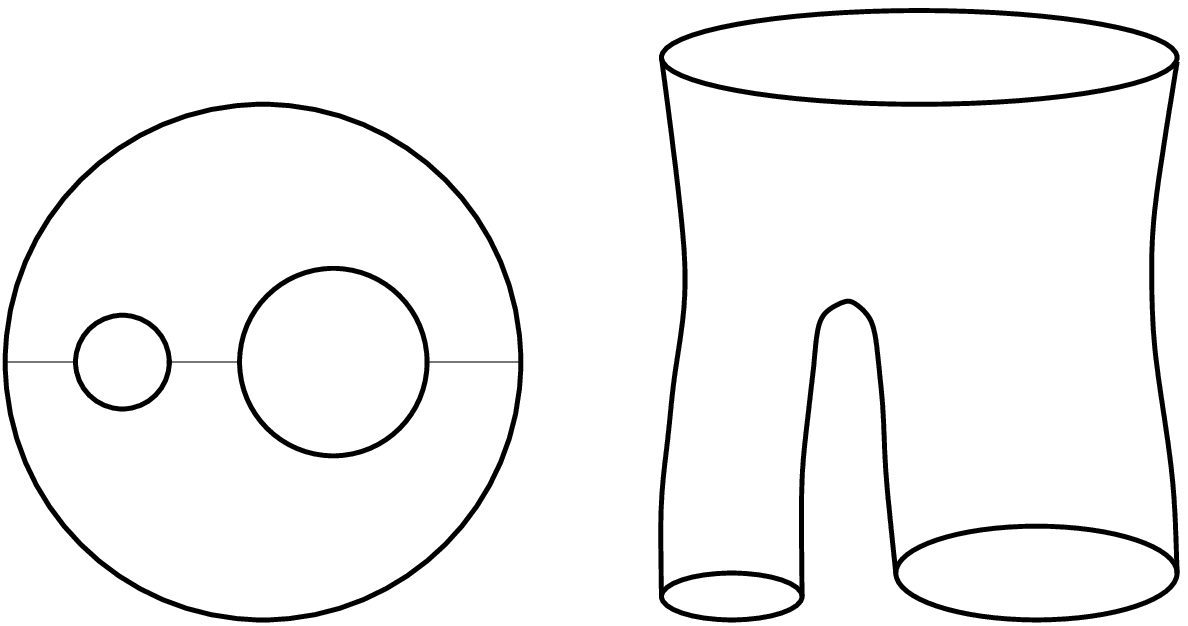}
\smallskip
\centerline{Fig.~3. Pants.}
\smallskip
One can also consider pants as a
sphe\-re with three holes. }
}}

\medskip
Let $C$ be a nodal curve parameterized by a real surface $\Sigma$. We shall
associate with every such curve $C$ a certain graph $\Gamma_C$ which
determines $C$ topologically in a unique way. In fact, $\Gamma_C$ will also
determine a decomposition of some components of $C$ into pants.

\smallskip
By definition, a compact component $C'$ is stable if it contains only a
finite number of automorphisms preserving marking points. In this case $C'
\bs\mapo$ possesses a unique so-called intrinsic metric.

\state Definition 4.2.3. {\sl The {\it intrinsic metric} for a smooth curve 
$C$ with marked points $\{ x_i \}$ and with boundary $\d C$ is a metric $g$ on
$C\bs \mapo$ satisfying the following properties:

\sli $g$ induces the given complex structure $j_C$;

\slii the Gauss curvature of $g$ is constantly -1;

\sliii $g$ is complete in a neighborhood of every marked point $x_i$;

\sliv every boundary circle $\gamma$ of $C$ is geodesic \wrt $g$.
}

\smallskip
Note that such a metric, if it exists, is unique, see \eg [Ab].

\smallskip
Now consider a component $C'$ of $C$ adjacent to some boundary circle of $C$.
Then $C' \bs \mapo$ is one of the following:

a)~a disk $\Delta$, or

b)~an annulus $A$, or

c)~a punctured disk $\check\Delta$,
\smallskip 
or else
d) $C'\bs \mapo$ admits the intrinsic metric.

\noindent
Note that if a component $C'$ is a disk or an annulus (both without marked
points), then $C'$ is the whole curve $C$. We shall consider cases a) and
b) later. Now we assume, for simplicity, that cases a) and b) do not occur.

\state Definition 4.2.4 \sl A component $C'$ of a nodal curve $C$ is called
{\it non-exceptional} \iff $C' \bs \mapo$ admits the intrinsic metric. \rm

\smallskip
In particular, nonstable components are exceptional compact ones, and
exceptional non-compact components are those of types a)--c) above.

\smallskip
Take some non-exceptional component $C'$ of $C$. There is a so- called maximal
partition of $C'\bs \mapo$ into pants $\{ C_1,\ldots C_n\}$ such that all
boundary components of these pants are either simple geodesic circles in
intrinsic metric or marked points, see [Ab]. Let us fix such a partition and
mark the obtained geodesic circles on $C'$.

\medskip
Now let $\sigma: \Sigma \to C$ be some parameterization of $C$. This defines
the set $\bfgamma'$ of the circles on $\Sigma$ which correspond to the nodes
of $C$. Let $\bfgamma''$ be the set of $\sigma $-pre-images of the geodesics
chosen above. Then $\bfgamma\deff \bfgamma' \sqcup \bfgamma''$ forms a system
of disjoint ``marked'' circles on $\Sigma$, which encodes the topological
structure of $C$. Now the graph $\Gamma_C$ in question can be constructed as
follows.

Define the set $V_C$ of vertices of $\Gamma_C$ to be the set $\{ S_j \}$ of
connected components of $\Sigma \bs \cup_{\gamma \in \bfgamma} \gamma =
\sqcup_j S_j$. Any $\gamma \in \bfgamma$ lies between 2 components, say $S_j$
and $S_k$, and we draw an edge connecting the corresponding 2 vertices.
Further, any boundary circle $\gamma$ of $\Sigma$ has the uniquely defined
component $S_j$ adjacent to $\gamma$. For any such $\gamma$ we draw a~{\sl
tail}, \i.e., an edge with one end free, attached to vertex $S_j$. Finally, we
mark all edges which correspond to the circles $\bfgamma'$, \i.e., those coming
from the nodes.

\medskip\medskip
\vbox{\xsize=.5\hsize\nolineskip\rm
\putm[.20][-.01]{\gamma_1^*}%
\putm[.46][.13]{\gamma_2^*}%
\putm[.67][.12]{\gamma_3^*}%
\putm[.21][.58]{a_1^*}%
\putm[.38][.815]{a_2^*}%
\putm[.53][.83]{a_3^*}%
\putt[1.05][0]{\advance \hsize -1.05\xsize
\centerline{Fig.~4. Graph of a curve $C$.}
\smallskip
Graph $\Gamma_C$ determines the topology of the curve $C$ in a unique way.
Take as many oriented spheres with as many vertices as 
$\Gamma_C$ has. For each edge
take a handle and join the corresponding sphe\-res by this handle. For each
tail make a hole (\ie remove a disk) in the corresponding sphere. Finally,
contract into points the circles on the handles corresponding to the marked
edges to get nodes. We obtain a topological space homeomorphic to $C$.
\smallskip
}
\noindent
\epsfxsize=\xsize\epsfbox{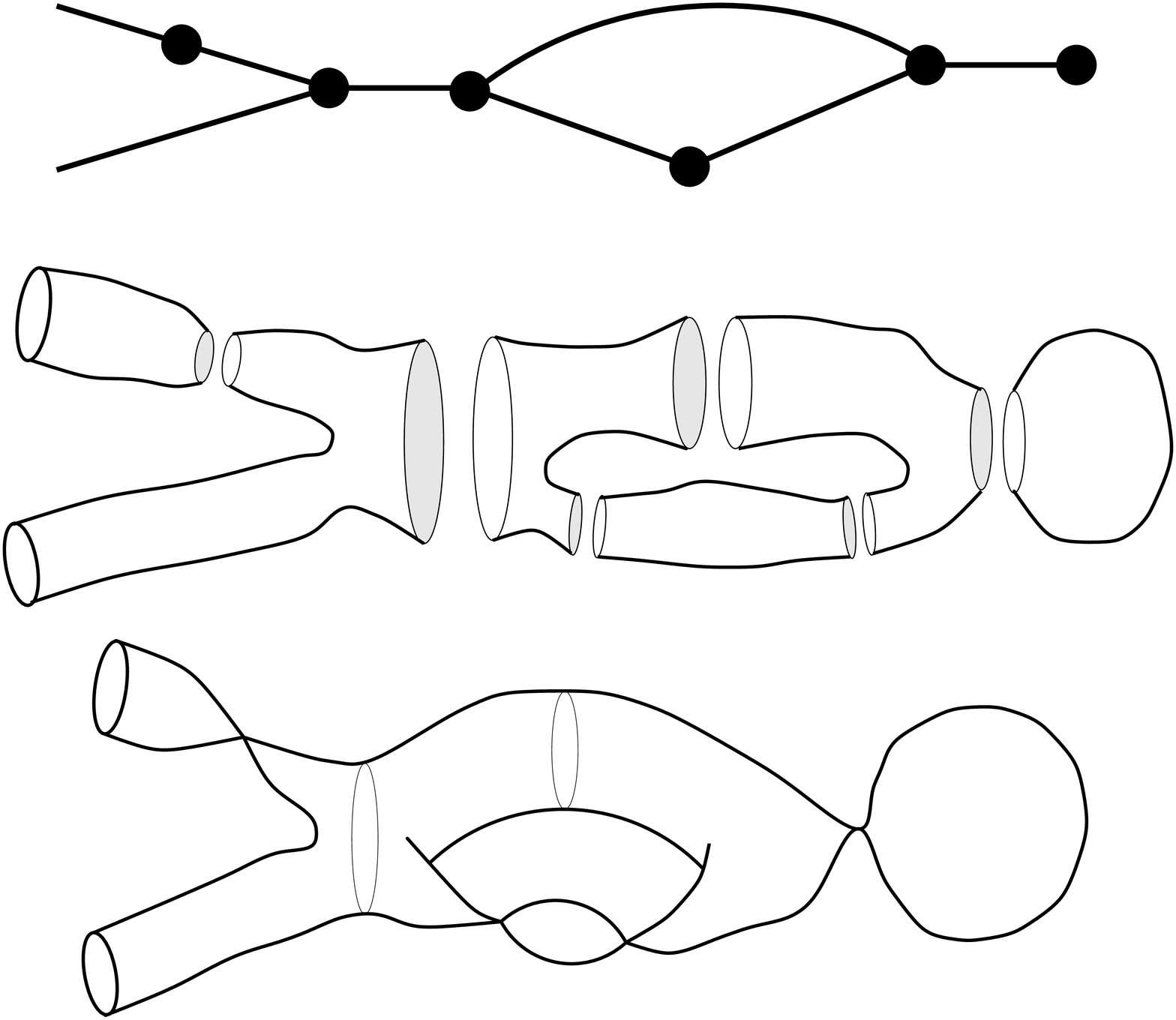}
}

\medskip
Having the graph $\Gamma$, which characterizes uniquely the topological
structure of $C$, we are now going to describe the set of parameters
defining (uniquely) the complex structure of the curves $C$. This is
equivalent to determining the complex structure and marked points on all
components of $C$. If such a component $C'$ is a sphere with 1 or 2 marked
points or a disk with 1 marked point, then its structure is defined by its
topology uniquely up to diffeomorphism. Otherwise, the component $C'$ is
non-exceptional. In this case the complex structure and the marked points can
be restored by the so- called {\sl Fenchel-Nielsen coordinates} on the
Teichm\"uller space $\ttt_{g,m,b}$. Recall that the space $\ttt_{g,m,b}$
parameterizes the complex structures on a Riemann surface $\Sigma$ of genus
$g$ with $m$ punctures (\ie marked points) and with a boundary consisting of
$b$ circles, see [Ab].

\smallskip
Let $C$ be a smooth complex curve with marked points of non-exceptional type,
so that $C$ admits the intrinsic metric. Fix some parameterization $\sigma:
\Sigma \to C$. Consider the pre-images of the marked points on $C$ as marked
points on $\Sigma$ or, equivalently, as punctures of $\Sigma$. Let $C \bs
\mapo = \cup_j C_j$ be a decomposition of $C$ into pants and $\Sigma\bs
\mapo = \cup_j S_j$ the induced decomposition of $\Sigma$.

Let $\{\gamma_i\}$ be the set of boundary circles of $\Sigma$. The boundary
of every pants $S_j$ has three components, each of them being either a marked
point of $\Sigma$ or a circle. In the last case this circle is either a
boundary component of $\Sigma$ or a boundary component of another pants, say
$S_k$. In this situation we denote by $\gamma _{jk}$ the circle lying between
the pants $S_j$ and $S_k$. Fix the orientation on $\gamma _{jk}$, induced
from $S_j$ if $j<k$ and from $S_k$ if $k<j$. For any such circle $\gamma
_{jk}$, fix a boundary component of $S_j$ different from $\gamma_{jk}$ and
denote it by $\d_k S_j$. In the same way fix a boundary component $\d_j S_k$.
Make similar notations on $C$ using primes to distinguish the circles on $C$
from those on $\Sigma$, \i.e., set $\gamma'_i \deff \sigma(\gamma_i)$ and
$\gamma'_{jk} \deff \sigma(\gamma_{jk})$.

By our construction, $\gamma'_{jk} = \sigma( \gamma_{jk})$ is a geodesic
\wrt the intrinsic metric in $C$.

\smallskip
\vbox{\xsize=.43\hsize\nolineskip\rm
\putm[.42][.12]{\gamma'_{jk}}%
\putm[.51][.21]{x^*_{j,k}}%
\putm[.33][.31]{x^*_{k,j}}%
\putm[.26][.24]{S_j}%
\putm[.70][.36]{S_k}%
\putm[-.05][.12]{\d_k S_j}%
\putm[.86][.29]{\d_j S_k}%
\putt[1.05][-.03]{\advance\hsize-1.05\xsize
\noindent
If the component $\d_k C_j$ is a marked point, we find on $C_j$ the (uniquely
defined) geodesic ray $\alpha_{j,k}$ starting at some point $x^*_{j,k} \in
\gamma'_{jk}$ and approaching $\d_k C_j$ at infinity such that $\alpha
_{j,k}$ has no self-intersections and is orthogonal to $\gamma'_{jk}$ at $x^*
_{j,k}$. Otherwise, we find on $C_j$ the shortest geodesic $\alpha _{j,k}$
which connects $\d_k C_j$ with $\gamma'_{jk}$ and we denote the point $\alpha
_{j,k} \cap \gamma'_{jk}$ by $x^*_{j,k}$. In both cases, this construction
determines a distinguished point $x^* _{j,k} \in \gamma' _{jk}$.
}
\noindent
\epsfxsize=\xsize\epsfbox{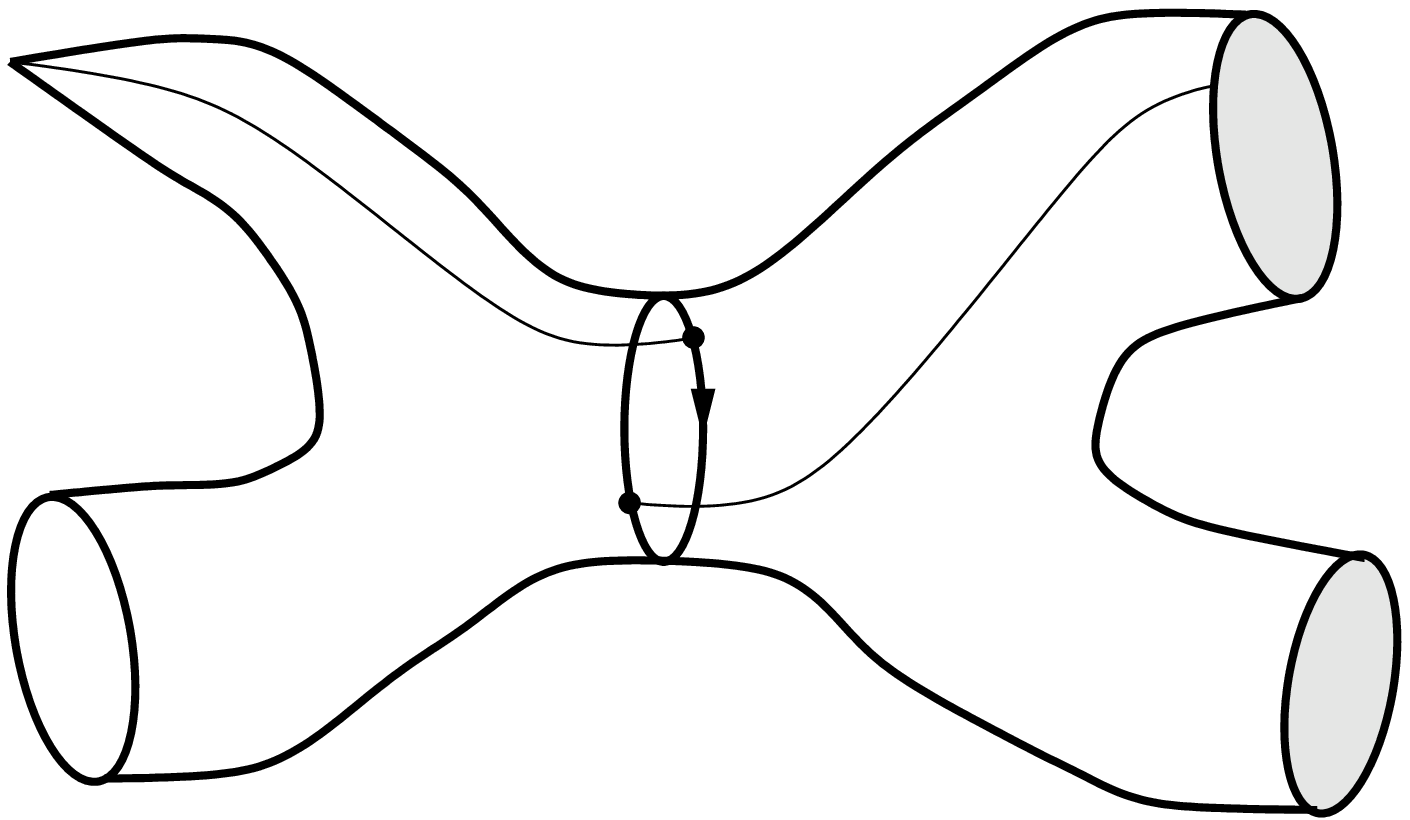}
\medskip
\vbox{\hsize=\xsize\centerline{Fig.~5. Marked points}
\centerline{on the circle $\gamma'_{jk}$.}
}
\smallskip
}

\bigskip
Using the same procedure in $C_k$, we obtain another point $x^*_{k,j} \in
\gamma' _{jk}$. Denote by $\ell_{jk}$ (resp.\ by $\ell_i$) the intrinsic
length of $\gamma'_{jk}$ (resp.\ of $\gamma'_i \deff \sigma( \gamma_i)$) in
$C$. For $j<k$ define $\lambda_{jk}$ as the intrinsic length of the arc on
$\gamma' _{jk}$, which starts at $x^* _{j,k}$ and goes to $x^* _{k,j}$ in the
direction determined by the orientation of $\gamma _{jk}$. Set $\vartheta
_{jk} \deff {2\pi \lambda _{jk} \over \ell_{jk} }$. We shall consider
$\vartheta _{jk}$ as a function of the complex structure $j_C$ on $C$ with
values in $S^1 \cong \rr/ 2\pi \zz$.

The parameters $\bfell \deff (\ell_i, \ell_{jk} )$ and
$\bfvartheta \deff(\vartheta_{jk})$ are called
{\sl Fenchel-Nielsen coordinates} of the complex
structure $j_C$. The reason is that these parameters determine up to
isomorphism the complex structure $j_C$ on the smooth complex curve with
marked points parameterized by a real surface $\Sigma$. In other words,
$(\bfell, \bfvartheta)$ can be considered as coordinates on $\ttt_{g,m,b}$.
More precisely, one has the following

\state Proposition 4.2.2. \it Let $\Sigma$ be a real surface of genus $g$ with
$m$ marked points and with the boundary consisting of $b$ circles, so that
$2g+m +b\ge3$. Let $\Sigma \bs \mapo = \cup_j S_j$ be its decomposition into
pants. Then

\item{\sl i)} for any given tuples $\bfell= (\ell_i, \ell_{jk})$ and
$\bfvartheta = (\vartheta_{jk})$ with $\ell_i, \ell_{jk} >0$ and $\vartheta_{jk} \in S^1$
there exists a complex structure $j_C$ on $\Sigma$ such that boundary circles
of all $S_j$ are geodesic \wrt the intrinsic metric on $\Sigma\bs \mapo$
defined by $j_C$ and the given $(\bfell, \bfvartheta)$
are Fenchel-Nielsen coordinates of $j$; moreover, such a structure $j_C$
is unique up to a diffeomorphism preserving the pants $S_j$ and the marked
points;

\item{\sl ii)} let $C$ be a smooth complex curve with parameterization $\sigma:
\Sigma \to C$ which has $m$ marked points; then there exists a parameterization
$\sigma_1: \Sigma \to C$ isotopic to $\sigma$, which maps boundary components
and marked points of $\Sigma$ onto the ones of $C$ in prescribed order
such that the boundary circles of $\sigma(S_j)$ are geodesic \wrt the
intrinsic metric on $C\bs \mapo$. \rm

\nobreak
\state Proof. See [Ab]. \qed

\bigskip\noindent
{\bigsl 4.3. Complex Structure on the Space $\ttt_\Gamma$.}

\smallskip
Let $\Sigma$ be a real surface of genus $g$ with $m$ marked points and with
the boundary consisting of $b$ circles. Assume that $2g+m +b\ge3$. Then there
exists a decomposition of $\Sigma \bs \mapo$ into pants, which is in general
not unique. The topological type of such a decomposition can be encoded in
graph $\Gamma$, associated with the decomposition. It is constructed in a
similar way to that above, but this time we must draw a tail  
for every marked
point, and then mark all those tails on the graph.

Let such a graph $\Gamma$ be fixed. We call two complex structures $J_1$ and
$J_2$ on $\Sigma$ isomorphic if there exists a biholomorphism $\phi: (\Sigma,
J_1) \cong (\Sigma, J_1)$ preserving the marked points of $\Sigma$ and the
decomposition of $\Sigma$ into pants given by graph $\Gamma$. Denote by
$\ttt' _\Gamma$ the space of isomorphism classes of complex structures on
$\Sigma$. By {\sl Proposition 1.2}, Fenchel-Nielsen coordinates identify
$\ttt'_\Gamma$ with the real manifold
$\rr_+^{3g-3+m+2b} \times (S^1)^{3g-3+m+b}$.

It is desirable to equip $\ttt'_\Gamma$ with some natural complex structure.
In so doing, the main difficulty is that the real dimension of $\ttt'_\Gamma$
can be odd. A possible explanation of this fact is that not all relevant
information (\ie parameters) about a complex structure has been taken into
consideration. Note that for any ``inner circle'' $\gamma_{jk}$ which appears
after the decomposition into pants, we have obtained a pair of coordinates,
mainly the length $\ell_{jk}$ and the angle $\vartheta_{jk}$. On the other
hand, for any boundary circle $\gamma_i$ of $\Sigma$ we have only the
length $\ell_i$. An obvious way to produce additional angle coordinates is to
introduce an additional marking of every boundary circle.

\state Definition 4.3.1. \sl A real surface $\Sigma$ or a nodal complex curve
$C$ is said to have a {\it marked boundary} if on every boundary circle of
$\Sigma$ (resp.\ $C$) a point is fixed. \rm

\state Remark. Later in \S\.5 we shall consider complex curves with several
marked points on boundary circles. But now we shall assume that on every
boundary circle exactly one point is marked.

\smallskip
``Missed'' angle coordinates $\vartheta_i$ can be now introduced similarly to
$\vartheta_{jk}$. For a boundary circle $\gamma_i$ we consider the adjacent
pants $S_j$. Fix a boundary component $\d_i S_j$ different from $\gamma_i$.
Let $J$ be a complex structure on $\Sigma$ such that the boundary circles 
of all
pants $S_k$ are geodesic \wrt the intrinsic metric defined by $J$. Using
constructions from above, find a geodesic (resp.\ a ray) $\alpha_i$ starting
at point $x^*_i \in \gamma_i$ and ending at the boundary circle $\d_i S_j$
(resp.\ approaching marked point $\d_i S_j$ of $\Sigma$). Take the marked
boundary point $\zeta_i$ on $\gamma_i$ and consider the length $\lambda_i$ of
the geodesic ark on $\gamma_i$, starting at $x^*_i$ and going to $\zeta_i$ in
the direction defined by the orientation of $\gamma_i$. Define $\vartheta_i
\deff {2\pi \lambda_i \over \ell_i }$, $\vartheta_i \in S^1\cong \rr/ 2\pi\zz$.
We include the coordinates $\vartheta_i$ in the system of angle coordinates
$\bfvartheta$. Denote by $\ttt_\Gamma$ the set of isomorphism classes of
complex structures on $\Sigma$ with marked boundary and with a given
decomposition into pants.

\smallskip
Let $C$ be a smooth complex curve with marked points and a marked boundary, $C
\bs \mapo = \cup_j C_j$ its decomposition into pants, $\sigma: \Sigma \to C$
a parameterization, and $\Sigma \bs \mapo = \cup_j S_j$ the induced
decomposition of $\Sigma$. To define a complex structure on $\ttt_\Gamma$, we
introduce special local holomorphic coordinates in a neighborhood of
boundary of pants on $C$. Consider some pants $S_j$ and its boundary circle
$\gamma^*$. It can be a boundary circle of $\Sigma$, $\gamma_i$ in our
previous notation, or a circle $\gamma_{jk}$ separating $S_j$ from another
pants $S_k$. Let $\ell^*$ be the intrinsic length of $\gamma^*$. Fix some
small $a>0$ and consider the annulus $A$ consisting of those $x\in S_j$ for
which the intrinsic distance $\dist(x,\gamma^*) <a$. The universal cover
$\ti A$ can be imbedded into the hyperbolic plane $\hh$ as an infinite strip
$\Theta$ of constant width $a$ such that one of its borders is a geodesic line
$L$. The action of a generator of $\pi_1(A) \cong \zz$ on $\ti A$ is defined
by the shift of $\Theta$ along $L$ by distance $\ell^*$.

Now consider the annulus $A' \deff [0, {\pi^2\over \ell^*}] \times S^1$ with
coordinates $\rho, \theta$, $0\le\rho < {\pi^2\over \ell^*}$, $0\le \theta \le
2\pi$ and with the metric $({\ell^* \over 2\pi} / \cos {\ell^* \rho\over 2\pi}
)^2 (d\rho^2 + d\theta^2)$. A direct computation shows that this metric is of
constant curvature $-1$ and that the boundary circle $\d_0A_1 \deff S^1 \times 
\{0\}$
is geodesic of length $\ell^*$, whereas $A'$ is complete in a neighborhood of
the other boundary circle. Consequently, the universal cover $\ti A{}'$ of
$A'$ can be imbedded in the hyperbolic plane $\hh$ as a hyperbolic half-plane
$\hh^+_L$ with a boundary line $L$ such that $\Theta \subset \hh^+_L$.
Moreover, the action of $\pi_1(A') \cong \zz$ on $\ti A{}' \cong \hh^+_L$ is
the same as for $\ti A \cong \Theta$. This shows that there exists an
{\sl isometric} imbedding of $A$ into $A'$ which maps $\gamma^*$ onto
$\d_0 A'$. Moreover, such an imbedding is unique up to rotations in the
coordinate $\theta$. This leads us to the following

\state Proposition 4.3.1. \it Let $C_j$ be pants with a complex structure
and $\gamma^*$ its boundary circle of the intrinsic length $\ell^*$. Let
$x^*$ be a point on $\gamma$. Then some collar annulus $A$ of $\gamma^*$
possesses the uniquely defined conformal coordinates $\theta\in S^1\cong \rr/
2\pi\zz$ and $\rho$ such that the intrinsic metric has the form $({\ell^*
\over 2\pi} / \cos {\ell^* \rho\over 2\pi})^2 (d\rho^2+ d\theta^2)$, $\rho|_{
\gamma^*} \equiv0$, $\theta(x^*) =0$ and the orientation on $S_j$
is given by $d\theta\wedge d\rho$. \rm

\medskip
We shall represent $\rho$ and $\theta$ also in the complex form $\zeta \deff
e^{-\rho + \isl \theta}$ and call $\zeta$ the {\sl intrinsic coordinate} of
the pants $C_j$ at $\gamma^*$. An important corollary of the description of
the intrinsic metric in a neighborhood of the boundary circle is the following
statement about non-degenerating complex structures in pants, see {\sl
Definition 1.7}.

\state Lemma 4.3.2. \it Let $C$ be a smooth complex curve with marked points
admitting the intrinsic metric and let $\gamma^*$ be a boundary circle of $C$
of length $\ell^*$.

\sli If there exists an annulus $A \subset C$ of conformal radius $R$ (\ie $A
\cong \{\,z \in \cc \;:\; 1<|z| < R \;\})$, adjacent to $\gamma^*$ and
containing no marked points, then $\log R \le {\pi^2 \over \ell^*}$.

\slii There exists a universal constant $a^*$ such that the condition $\ell^*
\le 1$ implies that there exists an annulus $A \subset C$ of conformal radius
$R$ with $\log R \ge { \pi^2\over \ell^*} -{2\pi\over a^*}$, which is
adjacent to $\gamma^*$, has area $a^*$ and contains no marked points of $C$.

\sliii Let $\gamma\subset C$ be a simple geodesic circle of the length $\ell$
and $A \subset C\bs\mapo$ annulus of conformal radius $R$ homotopy equivalent
to $\gamma$. Then $\log R \le {2\pi^2 \over \ell}$.
\rm

\state Proof. Let $\Omega$ be the universal cover of $C\bs \mapo$ equipped
with the intrinsic metric lifted from $C$. Then $\Omega$ can be isometrically
imbedded into the hyperbolic plane $\hh$ as a domain bounded by geodesic lines
such that each of these lines  covers some boundary circle $\gamma_i$. Take
some (not unique!) line $L$ covering the circle $\gamma^*$ and fix a hyperbolic
half-plane $\hh^+_L$ with a boundary line $L$, so that $\Omega \subset
\hh^+_L$.

Now consider the universal cover $\ti A$ of the annulus $A$ and provide it
with the metric induced from $C$. Then we can isometrically imbed $\ti A$
in $\hh^+_L$ in such a way that the line covering $\gamma^* \subset \d A$
will be mapped onto $L$. The action of a generator of $\pi_1(A) \cong \zz$ on
$\ti A$ is defined by the shift of $\hh$ along $L$ onto a distance $\ell^*$.
Consequently, $A$ can be isometrically imbedded into $\hh^+_L/ \pi_1(A)$, which
is the annulus $A'= [0, {\pi^2\over \ell^*}] \times S^1$ with coordinates
$\rho, \theta$, $0\le\rho < {\pi^2\over \ell^*}$, $0\le \theta \le 2\pi$ and
with a metric $({\ell^* \over 2\pi} / \cos {\ell^* \rho\over 2\pi})^2 
(d\rho^2 +d\theta^2)$. Note that the conformal radius of 
$A'$ is $e^{\pi^2/\ell^*}$. The
monotonicity of the conformal radius of annuli (see \eg [Ab], Ch.II, \S1.3)
yields the inequality $R\le e^{\pi^2/\ell^*}$ which is equivalent to first
assertion of the lemma.

\smallskip
Part \sliii of the lemma can be proved by same argument. More precisely,
under the hypothesis of part \sliii we imbed the annulus $A$ into the
annulus $A''= ]-{\pi^2\over \ell}, {\pi^2\over \ell}[ \times S^1$ with 
coordinates $\rho, \theta$, $- {\pi^2\over \ell} < \rho < {\pi^2\over \ell}$, 
$0\le \theta \le 2\pi$ and with a metric 
$({\ell \over 2\pi} / \cos {\ell \rho \over 2\pi})^2 (d\rho^2 +d\theta^2)$.
The conformal radius of $A$ is now estimated by the conformal radius of $A''$,
which is equal to $e^{2\pi^2\over \ell}$. 

\smallskip
The second part of our lemma follows from results of Ch.II, \S\.3.3 of [Ab].
{\sl Lemma 2}  says that there exists a universal constant $a^*$ with the
following property: If $\ell^* \le 1$, then there exists a collar neighborhood
$A$ of constant width $\rho^*$ and of area $a^*$, which is an annulus imbedded
in $C$ and contains no marked points of $C$. In particular, we can extend the
intrinsic coordinates $\rho$ and $\theta$ in $A$. Using these coordinates, we
present $A$ in the form $\{ (\rho,\theta): 0\le \rho \le \rho^*\}$ and compute
the area,
$$
a^*=\area A = 2\pi \int_{\rho=0}^{\rho^*} \left(\msmall{\ell^*/2\pi \over
\cos(\ell^*\rho/2\pi) }\right)^2 d\rho =
\ell^* \tan\left(\msmall{\ell^*\rho^* \over 2\pi}\right).
$$
Consequently, $\tan\left({\pi\over2}-{\ell^*\rho^* \over 2\pi}\right) =
\cotan\left({\ell^*\rho^* \over 2\pi}\right) = {\ell^*\over a^*}$. This
implies ${\pi\over2}-{\ell^*\rho^* \over 2\pi} \le {\ell^*\over a^*}$, which
is equivalent to $\rho^* \ge { \pi^2\over \ell^*} -{2\pi\over a^*}$. To finish
the proof we note that the conformal radius $R$ of $A$ is equal to 
$e^{\rho^*}$.
\qed

\medskip
Let $C$ be a smooth complex curve with marked points, $C_j$ a piece of a
decomposition of $C\bs \mapo$ into pants and $\gamma^*$ its boundary circle.
Then as a ``base point'' $x^*=\{\theta=0=\rho\}$ for the definition of the
intrinsic coordinate we shall use the point $x^*_{j,k}$ if $\sigma( \gamma^*
)$ is the geodesic separating $C_j$ from another pants $C_k$, or respectively,
the point $x^*_i$ if $\gamma^*$ is a boundary circle of $C$. We denote these
coordinates $\zeta_{j,k}= e^{-\rho_{j,k} + \isl \theta_{j,k}}$ and $\zeta_i=
e^{-\rho_i + \isl \theta_i}$. Note that $\vartheta_i$ is exactly the
$\theta$-coordinate of the marked boundary point $x_i\in \gamma_i$ with
respect to $x^*_i$, and $\vartheta_{jk}$ is the $\theta$-coordinate of
$x^*_{k,j}$ with respect to $x^*_{j,k}$.

Note also that any intrinsic coordinate of a pair $(\zeta_{j,k}, \zeta_{k,j}
)$ extends canonically from one collar neighborhood of $\gamma_{jk}$ to
another side in such a way that the formula for the intrinsic metric
remains valid. This extension possesses the property $\zeta_{j,k}\cdot \zeta
_{k,j} \equiv e^{\isl \vartheta_{jk}}$, where $\vartheta_{jk}$ is a constant
function. We can view this relation as the transition function from
$\zeta_{j,k}$ to $\zeta_{k,j}$.

A similar construction is possible in the case of a boundary circle
$\gamma_i$.
Namely, allowing $\rho_i$ to change also in the interval $(-{\pi^2 \over
\ell_i}, 0]$ and maintaining the formula $({\ell_i \over 2\pi} / \cos {\ell_i
\rho_i\over 2\pi})^2 (d\rho_i^2+ d\theta_i^2)$ for the metric, we can glue to
$\Sigma$ an annulus $(-{\pi^2\over \ell_i}, 0] \times S^1$ and extend 
the coordinate $\zeta_i = e^{-\rho_i + \isl \theta_i}$ there.

Making such a construction with every boundary circle $\gamma_i$, we obtain a
complex curve $C^{(N)}$ with the following properties. $C$ is relatively
compact in $C^{(N)}$ and the intrinsic metric of $C$ extends to a complete
Riemannian metric on $C^{(N)}$ with constant curvature $-1$. Such an
extension and the metric are unique. $C^{(N)}$ is called the {\sl Nielsen
extension of $C$}, see [Ab]. Note that the complex coordinate $\zeta_i$ can
be extended further to the unit disk $\{|\zeta_i|<1\}$.

\smallskip
Using the introduced complex coordinates $\zeta_i$ and $\zeta_{j,k}$, we
define a deformation family of complex structures on the curve $C$ with
marked boundary. Let $\lambda_i$ and $\lambda_{jk}$ be complex parameters
changing in small neighborhoods of $e^{\isl\theta_i}$ and $e^{\isl \theta
_{jk}}$, respectively. Having these data $\bflambda =(\lambda_i, \lambda_{jk}
)$, construct a complex curve $C_\bflambda$ in the following way. Take the
pants $\{C_j\}$ of the given decomposition of $C$ and extend all the
complex coordinates $\zeta_i$ and $\zeta_{jk}$ outside the pants.
Glue the pairs of coordinates $(\zeta_{j,k},\zeta_{k,j})$ with new transition
relations $\zeta_{j,k}\cdot \zeta_{k,j}= \lambda_{jk}$ (constant functions).
Move original boundary circles $\gamma_i= \{|\zeta_i|=1\}$ of $C$ to new
positions defined by the equations $|\zeta_i|= |\lambda_i|$ and mark the
points $\zeta_i = \lambda_i$ on them.

\state Theorem 4.3.3. \it The natural map 
$F:\bflambda \to(\bfell,\bfvartheta)$
is non-degenerated. In particular, $\bflambda$ can be considered as the set
of local complex coordinates on $\ttt_\Gamma$ and $\calc\deff \{ C_\bflambda
\}$ as a (local) universal holomorphic family of curves over $\ttt_\Gamma$.
\rm

\state Proof. Write the functions $\bflambda = (\lambda_i, \lambda_{jk})$
in the form $\lambda_i= e^{-r_i + \isl\phi_i}$, $\lambda_{jk}= e^{-r_{jk} +
\isl\phi_{jk}}$. From the definition of the map
$$
\textstyle
F: (e^{-r_i + \isl\phi_i}, e^{-r_{jk} + \isl\phi_{jk}})
\mapsto (\ell_i, \ell_{jk}; \vartheta_i, \vartheta_{jk})
$$
it is easy to see
that ${\d (\vartheta_i, \vartheta_{jk}) \over \d(\phi_i, \phi_{jk})}$ is the
identity matrix, whereas ${\d (\ell_i, \ell_{jk}) \over \d(\phi_i, \phi_{jk}
)}$ is equal to 0. So it remains to show that the matrix ${\d (\ell_i,
\ell_{jk}) \over \d(r_i, r_{jk})}$ is non-degenerate.

\smallskip
Consider a special case where $C$ is pants with the boundary circles
$\gamma_i$ (at least one) and, possibly, with marked points $x_j$. We shall
consider such points as punctures of $C$. Let $J$ denote the complex structure
on $C$ and let $\mu_0$ be the intrinsic metric. Extend the coordinates
$\zeta_i$ and the metric $\mu_0$ outside of $\gamma_i$ to some bigger complex
curve $\wt C$ with $C \comp \wt C$.

Fix real numbers $v_i$ and consider the domains $C_t$ in $\wt C$ defined in
local coordinates $\zeta_i= e^{-\rho_i + \isl\theta_i}$ by the inequalities 
$\rho_i \ge v_it$. This defines a family of deformations of $C= C_0$ 
parameterized
by a real parameter $t$, corresponding to a real curve in the parameter space
$\{ \bflambda \}$ given by $\lambda_i(t) = e^{v_it}$. Note that the
deformation is made in such a way that the original complex structure $J$ and
local holomorphic coordinates are preserved. Thus we can use them as an
``invariable basis'' in our calculations.

Let $\mu_t$ be the intrinsic metric of $C_t$. Without loss of generality we
may assume that $\mu_t$ extends to $\wt C$ as a metric with constant curvature
$-1$, which induces the original complex structure $J$ on $\wt C$. Since
the intrinsic metric depends smoothly on the operator $J$ of a complex
structure, $\mu_t$ are smooth in $t$. In a local holomorphic coordinate $z=x
+\isl y$ we can present $\mu_t$ in the form $e^{2\psi(t,z)}(dx^2 + dy^2)$.
The condition ${\sf Curv}(\mu_t) \equiv -1$ is equivalent to the differential
equation
$$
\d^2_{xx}\psi(t,\cdot) + \d^2_{yy} \psi(t, \cdot)= e^{2\psi(t, \cdot)},
$$
where $\d_x$ denotes the partial derivation ${\d \over \d x}$ and so on.
Differentiating it in $t$, we obtain $e^{-2\psi(t,\cdot)}(\d^2_{xx}+ \d^2_{yy})
\dot \psi(t,\cdot) = \dot\psi(t,\cdot)$, where $\dot\psi(t,\cdot)$
denotes the derivative of $\psi(t,\cdot)$ in $t$.

Note that $\d_t\mu_t = 2 \dot \psi(t,\cdot)\cdot \mu_t $, so $\dot \psi_t
(\cdot)= \dot\psi(t, \cdot)$ is independent of the choice of a local
holomorphic coordinate $z= x+ \isl y$ and is defined globally. The equation
$e^{-2\psi_z(t,\cdot)} (\d^2_{xx}+ \d^2_{yy})\dot\psi_z(t,\cdot) = \dot\psi_z
(t,\cdot)$ can be rewritten in the form $\Delta_t \dot\psi_t = 2 \dot \psi_t
$, with $\Delta_t$ denoting the Laplace operator for the metric $\mu_t$.

The condition that the circle $\gamma_i(t) \deff \{ \rho_i = v_it\}$ is
$\mu_t $-geodesic means that the covariant derivative $\nabla_{\theta_i}
(\d_{\theta_i})$ of the vector field $\d_{\theta_i}$, the tangent vector field
to $\gamma_i(t)$, must be parallel to $\d_{\theta_i}$ along $\gamma_i(t)$.
Expressing this relation in local coordinates $\rho_i$ and $\theta_i$,
we obtain
$\d_{\rho_i}\psi_i(t; v_it, \theta_i)=0$, where $\mu_t =e^{2 \psi_i (t;
\rho_i, \theta_i)}(\d\rho_i^2 + \d\theta_i^2)$ is a local representation of
the metric $\mu_t$. Deriving in $t$, we obtain $\d_{\rho_i}\dot \psi_i(t; v_it,
\theta_i) + v_i \d^2_{\rho_i\rho_i}\psi_i(t; v_it, \theta_i) =0$.

In the case where $t=0$ we have 
$\psi_i(0; 0, \theta_i) \equiv \log{\ell_i \over
2\pi}$, a constant. Hence $\d^2_{\theta_i\theta_i}\psi_i(0; 0, \theta_i) =
e^{2 \psi_i(0; 0, \theta_i)} $ and $\d_{\rho_i}\dot \psi_i(0; 0, \theta_i) =
- v_i e^{2 \psi_i(0; 0, \theta_i)}= - v_i \left({\ell_i \over2\pi} \right)^2
$. On the other hand, $\d_{\rho_i} = - {\ell_i \over 2\pi} \d_\nu$ on
$\gamma_i(0) = \gamma_i$, where $\nu$ denotes the unit outer normal field to
$C_0=C$. Consider the integral $\int_C |d\dot\psi {}_0|^2 + 2 \dot\psi{}_0^2
d\mu_0$. Integrating by parts, we get
$$
\eqalignno{
&\int_C |d\dot\psi{}_0|^2 + 2 \dot\psi{}_0^2 d\mu_0 =
\int_C \dot\psi_0(2 \dot\psi_0 - \Delta_0 \dot\psi_0) d\mu_0
+ \int_{\d C} \dot\psi_0 \d_\nu\dot\psi_0 dl=
\cr
\noalign{\vskip0pt\allowbreak}
=&
\sum_i \int_{\gamma_i} \dot\psi_0 \d_\nu\dot\psi_0
\msmall{\ell_i \over 2\pi}
d\theta_i =
\sum_i - \int_{\gamma_i}  \dot\psi_0 \d_{\rho_i}\dot\psi_0
d\theta_i =
\sum_i \int_{\gamma_i}  \dot\psi_0 v_i \left(\msmall
{\ell_i \over 2\pi}\right)^2 d\theta_i =
\cr
\noalign{\vskip0pt\allowbreak}
=&
\sum_i \msmall{v_i \ell_i \over 2\pi} \int_{\gamma_i}  \dot\psi_0
e^{\psi_i(0; 0, \theta_i)} d\theta_i =
\sum_i \msmall{v_i \ell_i \over 2\pi} \int_{\gamma_i}  (\dot\psi_0 +
v_i \d_{\rho_i} \psi_i(0; 0, \theta_i))e^{\psi_i(0; 0, \theta_i)} d\theta_i =
\cr
\noalign{\vskip0pt\allowbreak}
=&
\sum_i \msmall{v_i \ell_i \over 2\pi}
\left.\msmall{\d \over \d t}\right|_{t=0} \int_{\gamma_i}
e^{\psi_i(t; v_it, \theta_i)} d\theta_i =
\sum_i \msmall{v_i \ell_i \over 2\pi}
\left.\msmall{\d \over \d t}\right|_{t=0}\ell_i(t) =
\sum_i \msmall{ \ell_i \over 2\pi} v_i \dot\ell_i.
}
$$
Here $\dot\ell_i$ denotes the derivative of the length parameter $\ell_i$ for
the curve $C_t$ at $t=0$, so that
$(\dot\ell_1,\dot\ell_2,\dot\ell_3)= dF (v_1, v_2, v_3)$.
The obtained relation shows that the Jacobi matrix
$dF= {\d(\ell_i, \ell_2, \ell_3) \over \d(r_1, r_2, r_3)}$ is non-degenerate.
Otherwise there would exist a nonzero vector $(v_1, v_2, v_3)$ such that for
the deformation constructed above we get $\dot\ell_i=0$. But then $\dot\psi_0
\equiv 0$, which is a contradiction.

\smallskip
Now consider a general situation. Let $\Sigma$ be a real surface with a marked
boundary, $C$ a smooth curve with marked points, $\sigma:\Sigma \to C$ a
parameterization, and $C \bs \mapo = \cup_j C_j$ a decomposition into pants
with a given graph $\Gamma$. Let $\{\gamma_i\}$ be the set of boundary circles
and $\{ \gamma_{jk} \}$ the set of circles lying between the pants $C_j$ and
$C_k$, respectively. Consider these pants separately. Then for any circle $\gamma
_{jk} = \gamma_{kj}$ we obtain 2 distinguished ones, $\gamma_{j,k}$ considered
as a boundary circle of $C_j$, and $\gamma_{k,j}$ considered as a boundary circle
of $C_k$. Take real numbers $\bfv\deff (v_i, v_{j,k}, v_{k,j})$, where $v_i$ is
associated with the circle $\gamma_i$, $v_{j,k}$ with $\gamma_{j,k}$, and
$v_{k,j}$ with $\gamma_{k,j}$, respectively. Let $C_j(t\bfv)$ denote the pants
obtained from $C_j$ by the above construction using the corresponding parameters
$v_i$ and $v_{j,k}$. For $\bfv$ lying in a small ball $B = \{ |\bfv| < \eps \}$
all such families $C_j(t\bfv)$ can be extended for all $t \in [-1,1]$. Thus over
$B$ we obtain a collection of deformation families $C_j(\bfv)$ of complex structure
on pants $C_j$.

Let $\ell_i(\bfv)$, $\ell_{j,k}(\bfv)$, and $\ell_{k,j}(\bfv)$ denote the
lengths of circles $\gamma_i$, $\gamma_{j,k}$, and $\gamma_{k,j}$ \wrt the 
obtained intrinsic metrics $\mu_j(\bfv)$ on $C_j(\bfv)$. Denote by $\dot\ell_i$
a linear functional $\d_t|_{t=0} \ell_i(t\bfv)$, and define $\dot\ell_{j, k}$
similarly. The explicit formula for an intrinsic metric near a
boundary circle shows that $C_j(\bfv)$ can be glued to $C_k(\bfv)$ along
$\gamma_{jk}$ exactly when $\ell_{j,k}(\bfv) = \ell_{k,j}(\bfv)$. Since the
Jacobian $\d\bfell(\bfv) \over \d \bfv$ is non-degenerate, the conditions
$\ell_{j,k}(\bfv) = \ell_{k,j}(\bfv)$ define a submanifold $V \subset B$
whose tangent space $T_0V$ is given by the relations $\dot\ell_{j, k}= \dot\ell
_{k,j}$. Note that this defines a deformation family of complex structures on
$C$ over the base $V$ such that the map $\bfv \in V \mapsto \bfell(\bfv)$ is a
diffeomorphism.

We state that the set $(v_i, v_{j,k} + v_{k,j})$ is a system of coordinates
on $V$ in the neighborhood of $0\in V$. To prove this it is sufficient to
show that the linear map $\bfv=(v_i, v_{j,k}, v_{k,j}) \in T_0V \mapsto (v_i,
v_{j,k} + v_{k,j})$ is non-degenerate. If it is not true, then there
would exist a nontrivial $\bfv=(v_i, v_{j,k}, v_{k,j})\in T_0V$ with $v_i=0$
and $v_{j,k} + v_{k,j}=0$. Let $\dot\ell_i=\dot\ell_i(\bfv)$, $\dot\ell
_{j,k}= \dot\ell_{j,k}(\bfv)$ and $\dot \ell_{k,j} =\dot \ell_{k,j} (\bfv)$
be the corresponding derivatives of length. Then $\dot\ell_{j,k} = \dot
\ell_{k,j}$ and
$$
0< \sum_i \msmall{ \ell_i \over 2\pi} v_i \dot\ell_i +
\sum_{j<k} \msmall{ \ell_{jk} \over 2\pi} v_{j,k} \dot\ell_{j,k} +
\sum_{j<k} \msmall{ \ell_{jk} \over 2\pi} v_{k,j} \dot\ell_{k,j}=
\sum_{j<k} \msmall{ \ell_{jk} \over 2\pi} (v_{j,k}+ v_{k,j}) \dot\ell_{j,k}
=0.
$$
The obtained contradiction leads us to the following conclusion. The
functions $v_i$ and $v_{j,k}+ v_{k,j}$ define a coordinate system
on $V$ equivalent to $\bfell=(\ell_i, \ell_{jk})$.

\smallskip
Let us return to the holomorphic deformation family of complex structures on
$C$, defined by complex parameters $\lambda_i= e^{-r_i + \isl\phi_i}$ and
$\lambda_{jk}= e^{-r_{jk} + \isl\phi_{jk}}$. It is easy to see that the
Jacobian $\d(v_i, v_{j,k}+ v_{k,j}) \over \d( r_i , r_{jk})$  at the point
$(r_i , r_{jk})=0$ is the identity matrix. This fact proves the statement of
the lemma. \qed

\state Remark. At this point we give a possible reason why the complex
(\ie holomorphic) structure introduced by the complex coordinates
$\bflambda$ can be regarded as natural. Let $C$ be a complex curve with
marked points and a nonempty marked boundary. In the case where  
$C$ is a disk or
an annulus, assume additionally that at least one inner point of $C$ is
marked. Then in a neighborhood of every boundary circle $\gamma_i$ of $C$ we
can construct the intrinsic coordinate $\zeta_i$. Take 2 copies $C^+$ and
$C^-$ of $C$ and denote by $\tau$ the natural holomorphic map $\tau: C^\pm
\to C^\mp$ interchanging the copies. Denote by $\zeta_i^\pm$ the local
intrinsic coordinate on $C^\pm$ at boundary circles $\gamma_i^\pm$, both
corresponding to $\gamma_i$. Now we can glue $C^+$ and $C^-$ together along
every pair of circles $(\gamma_i^+, \gamma_i^-)$ by setting $\zeta_i^+ \cdot
\zeta_i^- =1$ as transition relations. We obtain a closed complex curve $C^d$
which admits a natural holomorphic involution $\tau: C^d \to C^d$. For the
constructed family $\{C_\bflambda \}$, the corresponding family $\{C^d
_\bflambda \}$ will be holomorphic. In fact, the statement of {\sl Theorem
2.3} means that $\{ C^d _\bflambda \}$ is a minimal complete family of
deformation of $C^d$ in the class of curves with holomorphic involution. This
construction of doubling should not be confused with another construction of
the {\sl Schottky double} $C^{Sch}$ of $C$ which provides an {\sl
antiholomorphic} involution $\tau^{Sch}: C^{Sch} \to C^{Sch}$. We shall use
the Schottky double $C^{Sch}$ in {\sl Section 5}, we will consider curves with
totally real boundary conditions.

\bigskip\noindent
{\bigsl 4.4. Invariant Description of the Holomorphic Structure on 
$\ttt_\Gamma$.}

\smallskip
The construction of (holomorphic) double $C^d$ shows how to give an invariant
description of a holomorphic structure on $\ttt_\Gamma$. Let $C$ be a smooth
complex curve with marked points and marked boundary, and $x\in \ttt_\Gamma$
the corresponding point on the moduli space. Denote by $D$ the 
divisor of marked
points. If the curve $C$ is not closed and $C^d$ is its double with holomorphic
involution $\tau$, we denote by $D^d \deff D + \tau(D)$ the double of $D$.

\state Lemma 4.4.1. \it If $C$ is closed, then the tangent space $T_x \ttt
_\Gamma$ is naturally isomorphic to $\sfh^1(C, \calo(TC)\otimes \calo(-D) )$.

If $C$ is not closed, then the space $T_x\ttt_\Gamma$ is naturally isomorphic
to the space $\sfh^1(C^d, \calo(TC^d)\otimes \calo(-D^d))^\tau$ of
$\tau$-invariant elements in $\sfh^1(C^d, \calo(TC^d)\otimes \calo(-D^d))$.

In both cases the complex structure on $T_x\ttt_\Gamma$ induced by local
complex coordinates $\bflambda$ coincides with those from $\sfh^1(C,
\calo(TC) \otimes \calo(-D))$ (resp.\ $\sfh^1(C^d, \calo(TC^d) \otimes
\calo(-D^d) )^{(\tau)}$). In particular, this defines a global complex
structure on space $\ttt_\Gamma$. \rm

\state Proof. The part concerning closed curves is well-known. In fact, the
natural isomorphism $\psi: T_x\ttt_\Gamma \to \sfh^1(C, \calo(TC)\otimes
\calo(-D))$ is a Kodaira-Spencer map. Its description is very simple in the
introduced local coordinates $\zeta_{j,k}$ on $C$ and $\bflambda =( \lambda
_{jk})$ on $\ttt_\Gamma$. Let $C\bs \mapo = \cup C_j$ be the decomposition of
$C$ into pants with the graph $\Gamma$. For every pants $C_j$ choose an open
set $\ti C_j$, containing a closure $\barr C_j= C_j \cup \d C_j$. Without
loss of generality we may assume that $\ti C_j$ are chosen not too big, so
that the covering $\calu\deff \{ \ti C_j \}$ is acyclic for the sheaf
$\calo(TC)$ and the local coordinates $\zeta_{j,k}$ are well-defined in
the intersections $\ti C_j \cap \ti C_k$. Then vector $v \in T_x\ttt_\Gamma$
with local representation $v= \sum_{j<k} v_{jk} {\d \over \d \lambda_{jk}}$
is mapped by the Kodaira-Spencer map $\psi$ to the \v{C}ech 1-cohomology class
$$
\psi(v) \in \sfh^1(C, \calo(TC)\otimes \calo(-D)) \cong
\check \sfh^1(\calu, \calo(TC) \otimes \calo(-D)),
$$
represented by the 1-cocycle
$$
\left(v_{jk} \zeta_{j,k} \msmall{\d \over \d \zeta_{j,k}} \right)
\in \prod_{j<k} \Gamma(\ti C_j \cap \ti C_k,  \calo(TC)\otimes \calo(-D)).
$$
For more details see [D-G].

\smallskip
Using this description of the Kodaira-Spencer map for closed curves with marked
points, it is easy to handle the case of curves with boundary. Let $C$ be
a non-compact curve with marked points and with decomposition $C \bs \mapo
= \cup_j C_j$. Take its double $C^d$ with the involution $\tau$. Then the
decomposition of $C$ induces a $\tau$-invariant decomposition
$C^d =\bigcap_j(C_j  \cap \tau C_j)$. The corresponding covering $\calu^d$
of $C^d$ can also be chosen to be $\tau$-invariant.

The local coordinates $\zeta_i$, corresponding to boundary circles $\gamma_i$
of $C$, can now be extended to a two-sided neighborhood of $\gamma_i$ in
$C^d$. The coordinates $\zeta_{j,k}$, corresponding to inner circles
$\gamma_{jk}$, induce local complex coordinates $\zeta^\tau_{j,k}\deff
\zeta_{j,k} \scirc \tau$ in $\tau( \ti C_j \cap \ti C_k)$.

Any deformation of the complex structure on $C$ induces a deformation of
the complex structure on $C^d$. This defines a map $\phi: \ttt_\Gamma \to
\ttt_{\Gamma^d}$, with $\Gamma^d$ denoting the graph corresponding to the
$\tau$-invariant decomposition of $C^d$ into pants. Using introduced
coordinates, we present a tangent vector $v \in T_x \ttt_\Gamma$ in the form
$$
v= \sum_i v_i \msmall{\d \over \d \lambda_i} +
\sum_{j<k} v_{jk} \msmall{\d \over \d \lambda_{jk}}.
$$
Then the composition of the Kodaira-Spencer map $\psi^d$ of $C^d$ with
the differential of $\phi: \ttt_\Gamma \to \ttt_{\Gamma^d}$ maps $v$ to
$$
\psi^d \scirc d\phi(v) \in \sfh^1(C^d, \calo(TC^d)\otimes \calo(-D^d)) \cong
\check \sfh^1(\calu^d, \calo(TC^d) \otimes \calo(-D^d)),
$$
represented by the \v{C}ech 1-cocycle
$$
\check v \deff \left(v_i \zeta_i \msmall{\d \over \d \zeta_i},
v_{jk} \zeta_{j,k} \msmall{\d \over \d \zeta_{j,k}},
v_{jk} \zeta^\tau_{j,k} \msmall{\d \over \d \zeta^\tau_{j,k}} \right)
\in \prod_i \Gamma(\ti C(i) \cap \tau C(i),  \calo(TC^d)\otimes \calo(-D^d))
\times
$$
$$
\prod_{j<k} \Gamma(\ti C_j \cap \ti C_k,  \calo(TC^d)\otimes \calo(-D^d))
\times \prod_{j<k} \Gamma(\tau(\ti C_j \cap \ti C_k),
\calo(TC^d)\otimes \calo(-D^d)),
$$
where $C(i)$ denotes the pants of $C$ adjacent to circle $\gamma_i$. It is
obvious that if all $v_i$ vanish, then this \v{C}ech 1-cocycle is $\tau
$-invariant. On the other hand, the relation $\zeta_i \cdot (\zeta_i \scirc
\tau) \equiv \lambda_i= \const$ implies that $\tau_*(\zeta_i {\d\over \d
\zeta_i}) = -\zeta_i {\d \over \d \zeta_i}$. The additional change of the 
sign of
the corresponding part of cocycle $\check v$ comes from the fact that $\tau$
interchange $C(i)$ with $\tau(C(i))$. This shows that $\check v$ is $\tau
$-invariant and the statement of the lemma follows. \qed

\medskip
Now we study the connection between the geometry of $\ttt_\Gamma$ and the
degeneration of complex structures on a real surface $\Sigma$ with marked
points and marked boundary. Let $\Sigma \bs \mapo= \cup_j S_j$ be a
decomposition into pants with graph $\Gamma$. The Fenchel-Nielsen coordinates
on $\ttt_\Gamma$ define a map $(\bflambda, \bfvartheta): \ttt_\Gamma \to (\rr
\times S^1)^{3g-3+ m+2b}$, which is a diffeomorphism by {\sl Proposition
1.2}. So, if $\{ j_n\}$ is a sequence of complex structures on $\Sigma$, its
degeneration means that the sequence of Fenchel-Nielsen coordinates of $\{
j_n\}$ is not bounded in $(\rr \times S^1)^{3g-3+ m+2b}$.

One can see that, in fact, we have two types of the degeneration. The first
one occurs when the maximum of the length coordinates $\ell_i$ and $\ell_{jk}$
of $j_n$ increases infinitely, and the second one is present
when a minimum of the length coordinates of $j_n$ vanishes. 
It should be pointed
out that for an appropriate sequence one can have both types of degeneration.

Note that by {\sl Proposition 1.2} the Fenchel-Nielsen coordinates of a
complex structure $j$ on $\Sigma$ are defined by a choice of a topological
type of decomposition of $\Sigma$ into pants, encoded in the graph $\Gamma$.
Thus the introduced notion of degeneration also depends on the choice of
$\Gamma$. Possibly, the best choice of such decomposition is
established by the following statement, proved in [Ab], Ch.II, \S\.3.3.

\state Proposition 4.4.2. \it Let $C$ be a complex curve with parameterization
$\sigma: \Sigma \to C$. Then

{\sl a)} there exists a universal constant $l^*>0$ such that any two geodesic
circles $\gamma'$ and $\gamma''$ on $C$ satisfying $\ell( \gamma') < l^*$ and
$\ell(\gamma'') < l^*$ are either disjoint or they coincide;

{\sl b)} there exists pants decomposition $C \bs \mapo= \cup_j S_j$ such
that the lengths of inner boundary circles $\gamma_{jk} = \barr S_j \cap \barr
S_k$ are bounded from above by a constant $L$ which depends only on 
the topology of\/ $\Sigma$ and the maximum $M$ of the lengths of 
the boundary circles of $C$;
moreover, any simple geodesic circle $\gamma$ on $C$ with $\ell(\gamma) <l^*$
occurs as a boundary circle of some $C_j$.

\smallskip
\state Corollary 4.4.3. \it 
Let $C_n$ be a sequence of nodal curves parameterized
by a real surface $\Sigma$ with uniformly bounded number of components.
Suppose that the complex structures of $C_n$ do not degenerate near the 
boundary.
Then, passing to a subsequence, one can find a decomposition $\Sigma =\cup_j
S_j$ and new parameterizations $\sigma'_n : \Sigma \to C_n$, such that

\sli the decomposition $\Sigma =\cup_j S_j$ induces a decomposition of every
non-exceptional component of $C_n$ into pants whose boundary circles are
geodesics;

\slii the intrinsic length of these geodesics are bounded uniformly in $n$.
\rm

\state Proof. By {\sl Lemma 4.3.2}, the intrinsic lengths of boundary circles
of non-exceptional components of $C_n$ are bounded uniformly in $n$. Find a
decomposition into pants of every non-exceptional component of $C_n$
satisfying the conditions of part {\sl b)} of {\sl Proposition 4.4.2}. Let 
$\Gamma_n$
denote the obtained graph of the decomposition of $C_n$. Since the number of
components of $C_n$ is uniformly bounded, we obtain a subsequence 
$C_{n_k}$ with
the same graph $\Gamma$ for all $n_k$.

It follows from the proof in [Ab] that the constant $L$ from part {\sl b)} of
{\sl Proposition 4.4.2} depends continuously on the maximum $M$ of the lengths
of the boundary circles of $C_n$. This implies condition {\sl i$\!$i)}.
Applying {\sl Proposition 4.2.2}, we complete the proof. \qed

%%%%%%%  %%%%%%%%%  %%%%%%%%  %%%%%%
\bigskip \noindent
{\bigsl 4.5. Example: Degeneration to a Half-cubic Parabola on the Language
of Stable Curves.}

\smallskip\noindent
{\sl Example.} Consider a sequence of curves $T_n$ in $\cc\pp^2$ given by
$$
T_n= \{\; y^2z-x^3=n^{-6}z^3\;\} \subset \cc\pp^2
$$
in homogeneous coordinates. In the affine chart $U_0=\{ z\not= 0\}$ the curves
are given by the equation $T_n=\{ y^2-x^3=n^{-6}\} $. In
any reasonable sense these curves should converge to a half-cubic parabola
$T_\infty:= \{ y^2=x^3\} $. We shall now explain in great detail how the
$T_n$ converge to $T_\infty$ in Gromov topology.

\smallskip
Each curve $T_n$ is the image of the curve
$$
T =  \{ y_1^2z_1-x_1^3=z_1^3\}\subset \cc\pp^2
$$
w.r.t.\ an algebraic map $f_n: [x_1, y_1,z_1] \in C
\mapsto [x:y:z]:= [n^{-2} x_1: n^{-3} y_1: z_1]$.

\smallskip We shall denote by $[x_1:y_1:z_1]$ the coordinates in the pre-image
and by $[x:y:z]$ in the range.

\medskip\noindent
{\sl Convergence of graphs.}
Consider the graphs $\Gamma_n$ of these mappings as subsets of 
$\cc\pp^2\times \cc\pp^2$.
Then the $\Gamma_n$ converge to a {\sl reducible} curve $\Gamma_\infty\in
\cc\pp^2\times \cc\pp^2$,
which consists of the graph of a {\sl constant} mapping
$[x_1: y_1:z_1] \in T \mapsto [0:0:1]$, which is a ``horizontal component'',
plus a ``vertical component'' over $[x_1:y_1:z_1] =[0:1:0]$,
which is a limit curve $\{y^2 = x^3\}$ of the sequence $T_n$.

\smallskip
The picture of the convergence is the following. Starting from $\Gamma_1$,
the curves $\Gamma_n$ transform continuously (we can define 
$\Gamma_t$ for any 
real $t$), so that for $n>\!>1$
curve $\Gamma_n$ consists of an ``almost horizontal'' torus and an ``almost
vertical'' sphere connected by a thin ``neck''. As $n\longrightarrow \infty$
this neck is shrinking to a point $P:=([0:1:0], [0:0:1]) \in \pp^2 \times
\pp^2$ which is the only singular point of $\Gamma_\infty$ where its
components intersect.

%%%%
\bigskip
\vbox{\xsize=1.1\hsize\nolineskip
\putm[.5][.053]{\Sigma}%
\putm[.57][-.02]{\gamma}%
\putm[.62][0.065]{p}%
\putm[.50][0.115]{\xbig[.04]\downarrow \sigma_n}%
\putm[.19][0.05]{\putm[-.0078][-.006]{ \xbig[.08]\downarrow \sigma_1}
        \xline[.16]}%
\putm[.63][0.05]{\xline[.30]\putm[-.0087][-.006]{%
\llap{$\sigma_\infty$}\xbig[.08]\downarrow}}%
\putm[.18][0.245]{C_1}%
\putm[.22][.20]{\scriptstyle\sigma_1(\gamma)}%
\putm[.19][0.306]{\xbig[.065]\downarrow i_1}%
\putm[.23][.51]{\Gamma_1}%
\putm[.06][.47]{\xlar[.09]\pr}%
\putm[.04][0.505]{T_1}%
\putm[.045][0.39]{\pp^2}%
\putm[.51][0.22]{C_n}%
\putm[.546][.207]{\scriptstyle\sigma_n(\gamma)}%
\putm[.57][0.306]{\xbig[.065]\downarrow i_n}%
\putm[.59][.465]{\Gamma_n}%
\putm[.473][.447]{\xlar[.09]\pr}%
\putm[.42][0.49]{T_n}%
\putm[.471][0.38]{\pp^2}%
\putm[.91][0.255]{C_\infty}%
\putm[.82][.19]{\scriptstyle\sigma_\infty(\gamma)}%
\putm[.89][0.306]{\xbig[.065]\downarrow i_\infty}%
\putm[.88][.51]{\Gamma_\infty}%
\putm[.79][.47]{\xlar[.09]\pr}%
\putm[.75][0.50]{T_\infty}%
\putm[.74][0.41]{\pp^2}%
\noindent
\centerline{\epsfxsize=\xsize\epsfbox{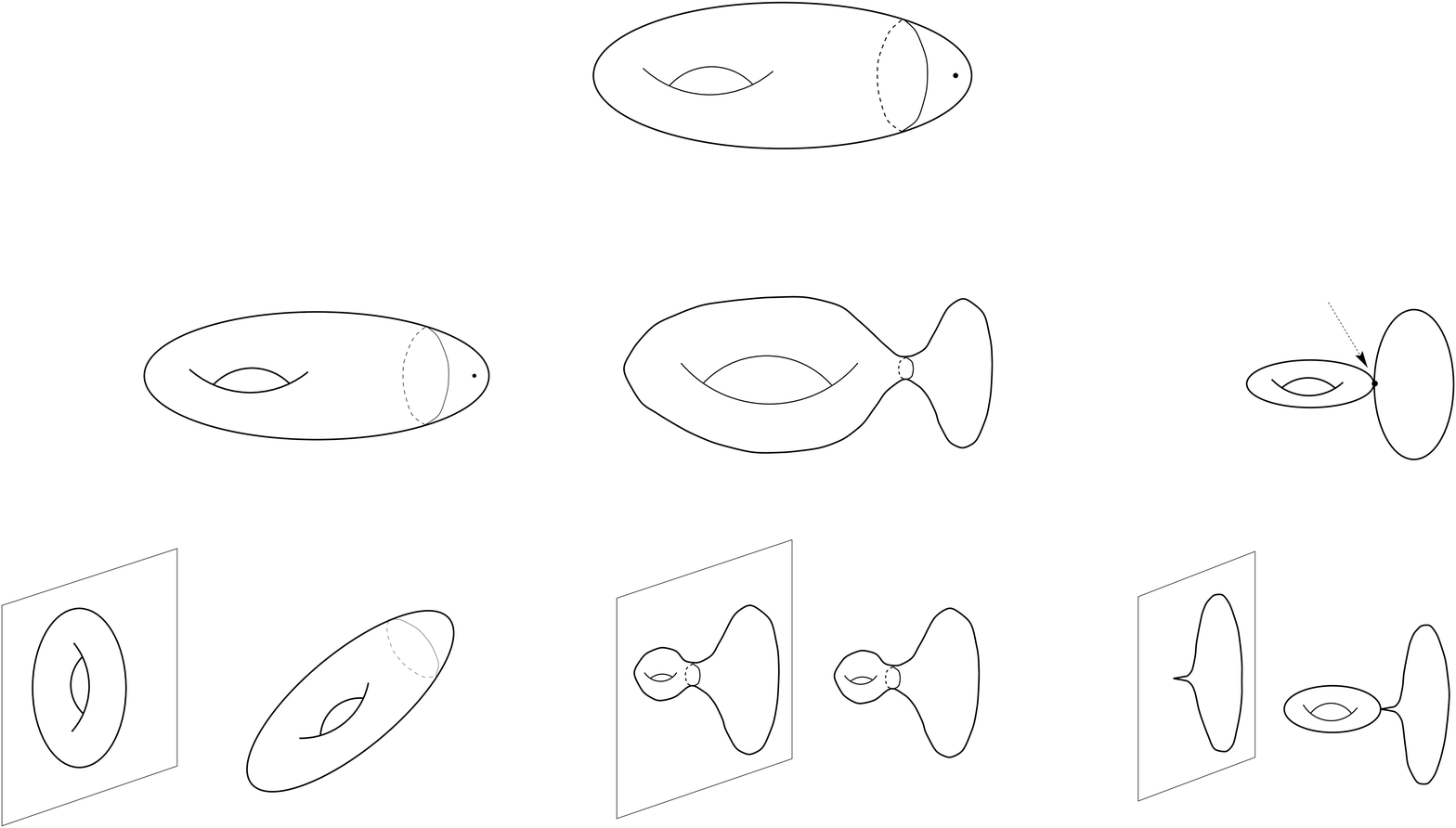}}%
\smallskip
\centerline{Fig.~6. Convergence picture}%
}

\medskip
\noindent {\sl Explanation of the picture:} $i_1,i_2,i_\infty$ are
imbeddings of $C_n$'s into $\pp^2 \times \pp^2$  as the graphs $\Gamma_n$'s.
The map $\pr:\pp^2 \times \pp^2 \to \pp^2$ is the natural projection onto the
range $\pp^2$.

\smallskip
Thus, the picture of convergency on the language of parameterized curves should
be the following: the surface parameterizing all $T_n$ (and also $T_\infty
=\{ y^2=x^3\} $) should be a torus $\Sigma $. Some fixed circle $\gamma $ on
$\Sigma $ around some fixed point $p\in \Sigma $ should be mapped by
parameterizations $\sigma_n :\Sigma \to C_n$ to ``smaller circles". Finally
the  parameterization $\sigma_\infty :\Sigma \to C_\infty$ contracts
$\gamma $ to a point. Thus $C_\infty$ should be a nodal curve which
consists of a smooth torus and a smooth sphere intersecting at one point.
The limit map $u_\infty$ should map the torus to the point $[0:0:1]$
and the sphere onto the limit curve $T_\infty$.

Now we must determine curves $C_n$, parameterizations $\sigma_n$ and
mappings $u_n$ explicitly. Note that one cannot take as $C_n$ simply the
graph $\Gamma_n$, because $\Gamma_\infty$ has a point
$P=([0:1:0], [0:0:1]) \in \pp^2 \times \pp^2$ as a cusp and a node at the
same time. However, by our definition of a nodal curve, $C_\infty$ should
have at most nodes. This is not a big problem, because one should 
``holomorphically''
parameterize $\Gamma_\infty$ by a nodal curve and at the same time
one also needs to parameterize ``coherently" $C_n$'s.

Let us do just that.

\medskip\noindent{\sl Explicit Parameterization.}

\smallskip
Now we give a description of the convergence in terms of stable curves
(terminology of [KM]; an alternative terminology speaks of ``stable maps''). 
Consider homogeneous equations
$$
y_1^2z_1 -x_1^3=z_1^3
\hbox{\qquad and \qquad}
{x\over n^{-2}x_1}= {y \over n^{-3}y_1}= {z\over z_1},
\eqno(4.5.1)
$$
defining $\Gamma_n$ in $\pp^2 \times \pp^2$. Let $Q=[x_1:y_1:z_1]$ be a point 
on $T$ with $z_1\not =0$. Then obviously $f_n(Q)= [n^{-2} x_1: n^{-3} y_1: 
z_1]$ converge to the point $[0:0:1]$ in the range $\pp^2$. Moreover, for 
a sufficiently small neighborhood $U$ of $Q$ in $T$ we obviously have 
a uniform convergence of the restriction $f_n|_U$ to a constant map $f_\infty:
 U \to \pp^2$, $f_\infty|_U  \equiv[0:0:1]$.

The set of points $[x_1:y_1:z_1]$ on $C$ with $z_1=0$ consists of one
point $[0:1:0]$. This means that for any compact $K \Subset C \bs [0:1:0]$
the sequence of restricted maps $f_n|_K : \to \pp^2$ converge uniformly
to the constant map $f_\infty: K \to \pp^2$, $u_\infty|_K \equiv[0:0:1]$. 
In terms of graphs $\Gamma_n$ this means that $\Gamma_n \cap(K\times \pp^2)$ 
converge to $\Gamma_\infty \cap(K \times \pp^2)$, which lies
in the ``horizontal'' part $T \times [0:0:1]$ of $\Gamma_\infty$.

Consider the behavior of $f_n$ in a neighborhood $V\subset T$ of the point
$[0:1:0]$. Setting $y_1=1$ and considering $x_1, z_1$ as an affine coordinate
on the first $\pp^2$ we obtain the relation
$$
z_1= x_1^3 + z_1^3.
$$
One can easily see
that $x_1$ can be chosen as a local holomorphic coordinate in a neighborhood
of $[0:1:0]$, and that there exists a holomorphic function
$\phi(x_1)$ in the disk $\Delta(a)=\{ |x_1| <a\}$ such that $\phi(0)=1$ and
$V=\{ z_1=x_1^3\phi(x_1)\} $.

\medskip
%% Fig.~cusp2.
\bigskip
\vbox{\xsize=\hsize\nolineskip
\putm[.00][.00]{z_1\xbig[.10]\uparrow}%
\putm[.32][.16]{T}%
\putm[.175][-.01]{V}%
\putm[.08][.06]{V_n}%
\putm[.025][.492]{\text{$x_1$-plane}}%
\putm[0][0]{\epsfxsize=\xsize\epsfbox{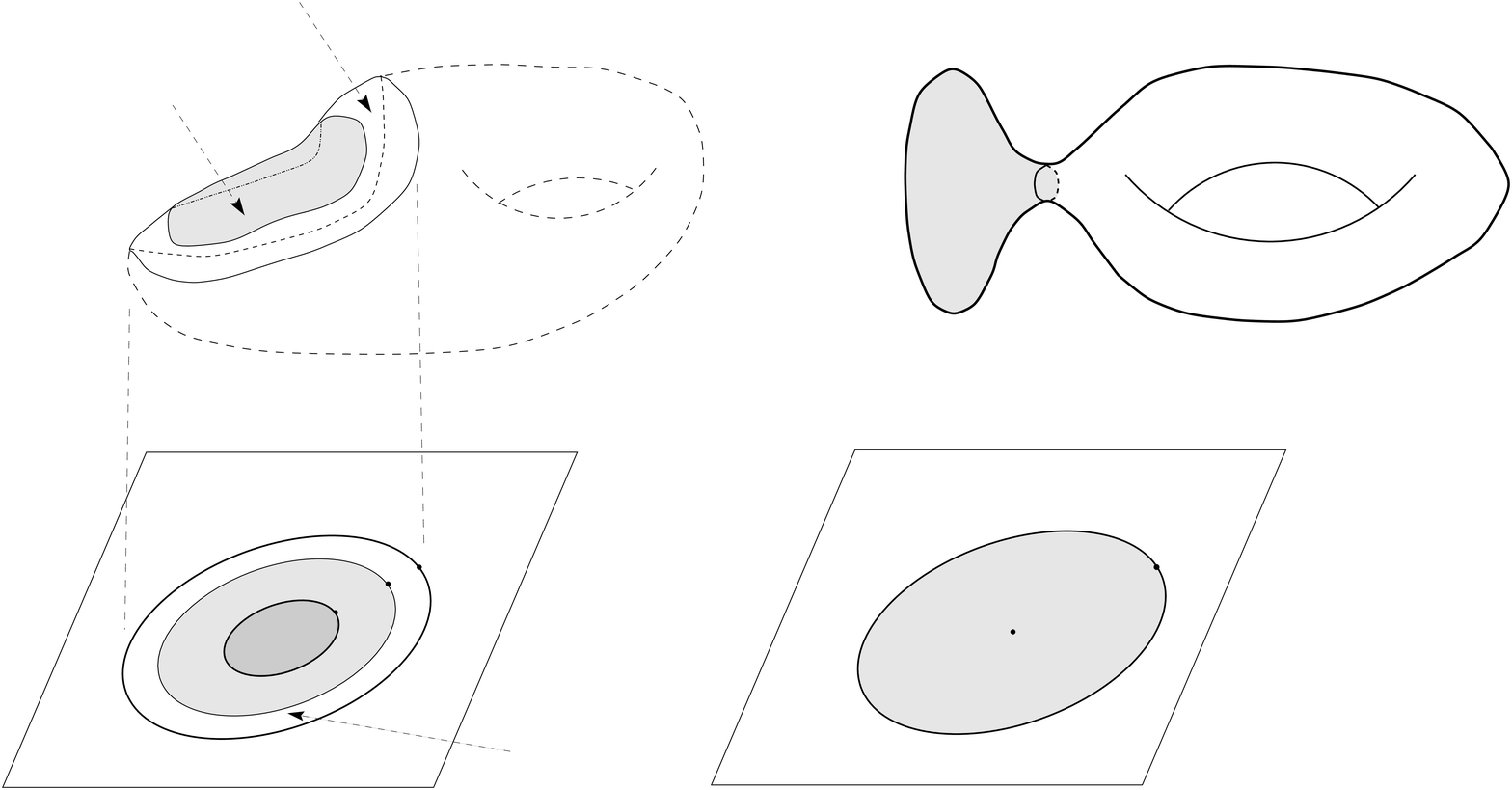}}%
\putm[.77][.36]{\scriptstyle R}%
\putm[.60][.43]{\scriptstyle \Delta(R)}%
\putm[.282][.367]{\scriptstyle a}%
\putm[.335][.49]{\scriptstyle \Delta(a)}%
\putm[.258][.378]{\scriptstyle {R\over  n}}%
\putm[.225][.393]{\scriptstyle {1\over  bn}}%
\putm[.49][.49]{\text{$\xi$-plane}}%
\putm[.33][.45]{\xlar[.16]{x_1= \xi/n \getsfrom\xi}}%
\putm[.63][.225]{\xbig[.06]\uparrow v_n}%
\vskip .53\xsize%
\smallskip
\centerline{Fig.~7. Parameterization of the Bubbled Sphere}%
}

\medskip
Fix $a>{1\over b}>0$ and take $R>{1\over b}$. Denote by $\pr:V\to \Delta (a)$
the natural projection, see {\sl Fig.~7}. For every $n$ sufficiently big,
cover $\pr(V)=\Delta(a)$ by sets $V_n$, where $V_n:=\pr^{-1}\Delta({R\over n})$
 and $V_n':=\pr^{-1}(A_n)$ with $A_n = \{ {1\over bn} <|x_1| <a\} $.
First we consider the behavior of appropriately 
rescaled maps $f_n\circ \pr^{-1}$
in $V_n$. For a complex coordinate $\xi \in \cc$
we consider maps $v_n(\xi): \Delta(R) \to \pp^2$ with $v_n(\xi):=
f_n\circ \pr^{-1}({\xi\over n})$ and the domain of definition
$\Delta(R):=\{|\xi|<R\}$.
This means that we rescale $f_n\circ \pr^{-1}$ in $\Delta({R\over n})=$ by
setting $nx_1 =\xi$. Then
$$\textstyle
v_n(\xi)=f_n\circ \pr^{-1}\left(\left[{\xi\over n}:1: {\xi^3\over n^3}
\phi\left({\xi\over n}\right) \right]\right)=
\left[{\xi\over n^3}:{1\over n^3}: {\xi^3\over n^3}
\phi\left({\xi\over n}\right) \right]=
\left[\xi :1 :\xi^3\phi\left({\xi\over n}\right) \right].
$$
Thus, images $v_n(\Delta(R))=f_n(\Delta ({R\over n}))$ lie in the affine chart
$y\not=0$ in the range
$\pp^2$, and in the affine coordinates $x,z$ ($y=1$) we obtain a representation
$v_n(\xi)= (\xi, \xi^3 \phi({\xi\over n}))$. We easily see that
$\phi({\xi\over n})$ converge on compacts $\Delta(R)$ to the constant function
1. This implies that maps $v_n$ converge on $\Delta(R)$ to the
map $v_\infty: \cc \to \pp^2$, $v_\infty(\xi)= [\xi:1:\xi^3]$. This map
is a parameterization of the ``vertical part'' of $\Gamma_\infty$.

\smallskip
Finally, we  want to describe the behavior of the ``neck'' of $\Gamma_n$
when $n\lrar \infty$. Again, take $x_1$ to be a local holomorphic coordinate
on $T$ in a neighborhood of $[0:1:0]$ setting $y_1=1$. Consider imbeddings
$h_n: A_n \to \Delta^2(a,b)= \{ |x_1| < a\} \times \{|t|<b\}$ with $h_n(x_1)
=(x_1, {1\over n x_1})$. Denote by $\cala_n$ the images $h_n(A_n)$.
Then  $\cala_n = \{ (x_1,t) \in \Delta^2(a,b):
x_1{\cdot}t= {1\over n}\}$. As $n\lrar \infty$, the annuli $\cala_n$
converge to a subset $\cala_\infty \subset \Delta^2(a,b)$,
which is the union of two disks $\Delta(a) \times \{0\}\cup \{0\} \times
\Delta(b)$.

Note that, along $\cala_n \subset \Delta^2(a,b)$
with coordinates $x_1,t$, we can represent $f_n\circ \pr^{-1}\circ h_n^{-1}$
in the form
$$
f_n\circ \pr^{-1}\circ h_n^{-1}:(x_1, t) \in A_n \mapsto \left[{x_1 \over n^2}
: {1\over n^3}: x_1^3
\phi(x_1) \right] =[ nx_1 t^3: t^3: n^3x_1^3 t^3\phi(x_1)]=
$$
$$
=[t^2:t^3: \phi(x_1)],
$$
where we use the relation $nx_1 t=1$ along $\cala_n$. Thus the restriction
of $f_n\circ \pr^{-1}$
onto $A_n$ coincides with the composition of imbeddings $h_n$ of $A_n$ into
$\Delta^2(a,b)$ with the map $F: \Delta^2(a,b) \to \pp^2$, $F(x_1,t) 
= [t^2:t^3: \phi(x_1)]$. In the affine chart with $z\not=0$ we have $F(x_1,t) =
\left({t^2 \over \phi(x_1)}, {t^3 \over\phi(x_1)}\right)$, see {\sl Fig.~8}.

\bigskip
\vbox{\xsize=\hsize\nolineskip
\putm[.00][.00]{z_1\xbig[.10]\uparrow}%
\putm[.32][.13]{T}%
\putm[.48][.11]{\xrar[.12]{f_n}}%
\putm[.45][.37]{\xrar[.12]{h_n}}%
\putm[.32][.19]{\xbig[.0348]\downarrow\pr}%
\putm[.08][.06]{V'_n}%
\putm[.025][.453]{\text{$x_1$-plane}}%
\putm[0][0]{\epsfxsize=\xsize\epsfbox{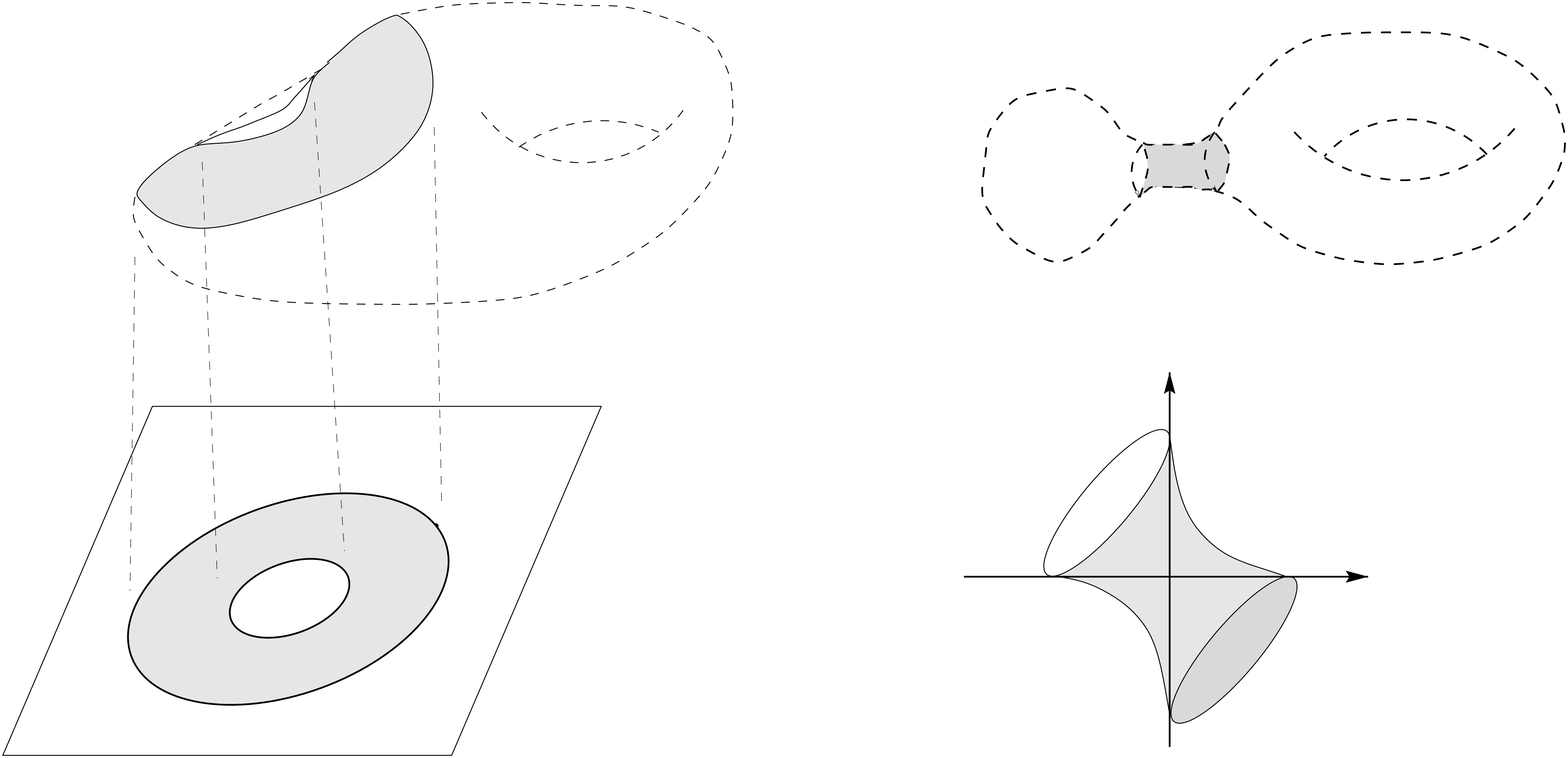}}%
\putm[.282][.33]{\scriptstyle a}%
\putm[.11][.40]{\scriptstyle A_n}%
\putm[.225][.35]{\scriptstyle {1\over  bn}}%
\putm[.77][.17]{\xbig[.07]\uparrow F}%
\vskip .51\xsize%
\centerline{Fig.~8. Parameterization of the ``Neck''.}
\smallskip
}

\medskip
Note that
$$
F\mid_{\cala_n} = f_n\circ \pr^{-1}\circ h_n^{-1}:\cala_n\to \cc\pp^2.
\eqno(4.5.2)
$$

\state Remark. Denote by $P_n$ the compact curve $\{ x_1\cdot t= {1\over n}
\}$ in $\cc\pp^2$. For us it is important that $F$ is well- defined on
$P_n\setminus \{ \vert x_1\vert >a\} $.

We interpret this picture in the following way. Annuli $\cala_n$
degenerate
to a normal crossing of two disks, a ``node'' in the terminology of [OM],
whereas maps $f_n\circ \pr^{-1}\circ h_n^{-1}:\cala_n \to \pp^2$ converge
to a (holomorphic) map from $\cala_\infty$ to $\pp^2$ which is, in this
case, a restriction of $F$ onto $\cala_\infty$, \i.e., constantly $[0:0:1]$
on the one component $\{ t = 0\} $ of $A_\infty$,
and, on the component $\{ x_1 = 0\} $, mapping $F$ is a parameterization
$t\in \Delta(b) \mapsto (t^2, t^3)$ of curve $y^2=x^3$.

Note that on $A_n$ the contractiong circle is $\gamma_n = \{ \vert x_1\vert
 = \sqrt{{1\over n}}$, because $h_n(\gamma_n) = \{ \vert x_1\vert = \vert t
\vert = \sqrt{{1\over n}}$.

Let us finally construct stable curves $(C_n,u_n)$ and $(C_\infty,
u_\infty)$. Put $W_1 = T_1\setminus V$ and $W_2 = P_n\setminus \{ \vert
x_1 \vert > a \} $. Identify the boundaries of $W_1$ and $W_2$ by
$h_n\circ \pr : \partial W_1 \to \partial W_2 $ to get the curve $C_n :=
W_1\sqcup W_2/\partial \sim \partial W_2$. Note that $C_\infty$ is a
nodal curve with one nodal point and two components: torus and sphere.

The corresponding maps are defined as follows:
$$
u(q) = \cases f_n (q) & \text{ if }q\in W_1 \\
F(q) & \text{ if }q\in W_2 .
\endcases 
$$
Note that $u_n$ is well- defined due to relation (2).

\bigskip
\vbox{\xsize=\hsize\nolineskip
\def\scst{\scriptstyle}
\putm[-.02][.27]{\Sigma}%
\putm[-.01][.07]{\scst\gamma}%
\putm[.03][.14]{\scst \wt V_1}%
\putm[.08][.03]{\scst\wt V_2}%
\putm[.13][.18]{\xrar[.09]{\sigma_n}}%
\putm[-.02][.58]{\Sigma}%
\putm[-.01][.38]{\scst\gamma}%
\putm[.03][.453]{\scst\wt V_1}%
\putm[.08][.34]{\scst\wt V_2}%
\putm[.15][.45]{\xrar[.09]{\sigma_\infty}}%
\putm[.24][.25]{C_n}%
\putm[.39][.19]{\scst W_1}%
\putm[.23][.11]{\scst\sigma_n(\gamma)}%
\putm[.31][.08]{\scst\cala_n}%
\putm[.24][.55]{C_\infty}%
\putm[.39][.51]{\scst W_1}%
\putm[.315][.403]{\scst\cala'_\infty}%
\putm[.31][.46]{\scst\cala''_\infty}%
\putm[.35][.44]{\scst 0}%
\putm[.85][.17]{T_n}%
\putm[.81][.55]{T_\infty}%
\putm[.57][.10]{\xrar[.1]{f_n}}%
\putm[.43][.40]{\xrar[.37]{F:\,t \;\mapsto\; (t^2\!\!,\, t^3)}}%
\putm[.62][.477]{\xrar[.20]{f_\infty\;\equiv\; [0,0,1]}}%
\putm[0][0]{\epsfxsize=\xsize\epsfbox{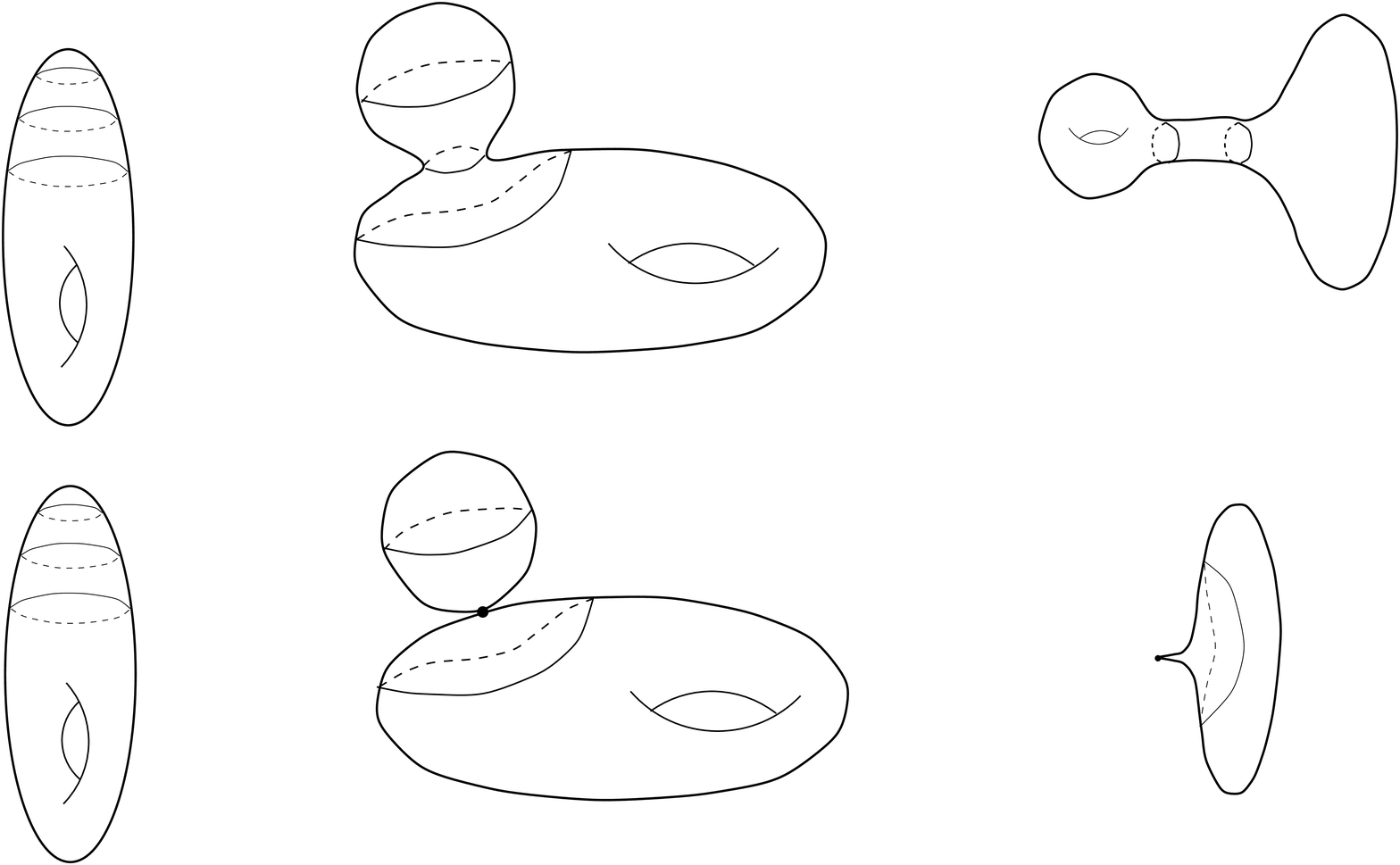}}%
\vskip .62\xsize%
\centerline{Fig.~10. }
}

\smallskip The rest is obvious, \i.e., one should take the parameterizations 
$\sigma_n:\Sigma
\to C_n$  satisfying the three following conditions:

(a) Some fixed circle $\gamma $ should be mapped by $\sigma_n$  onto $h_n(
\gamma_n)\subset \cala_n$ (also $\sigma_\infty(\gamma)= 0\in \cala_\infty$).

(b) Some fixed disk $\tilde V_2\subset \Sigma $, see {\sl Fig.~10}, should be
mapped by $\sigma_n$ onto the domain $P_n\cap \{ \vert t\vert >1\} \subset
W_2$ (also  $\sigma_\infty(\tilde V_2 = P_\infty\cap
\{ \vert t\vert >1\} $.

(c) Finally, some other domain $\tilde V_1\subset \Sigma$, which contains
a handle, should be mapped onto $W_1$.

One should also take care that the structures $\sigma_n^*j_n$ should
converge on compact subsets in $\Sigma \setminus \gamma$.

\newpage \noindent
{\bigbf Lecture 5}

\smallskip\noindent
{\bigbf Gromov Compactness Theorem}

\medskip\noindent{\bigsl 5.1. Second A priori Estimate.}

\smallskip
Let $(X,J)$ be an almost complex manifold. In what follows the  tensor $J$ is
supposed to be only continuous, \i.e., of class $C^0$. Fix some 
Riemannian metric
$h$ on $X$. All norms and distances will be taken with respect to $h$.
In particular, we have the following

\state Definition 5.1.1. {\sl A continuous almost complex structure $J$ is
called {\it uniformly continuous on $A \subset X$ with respect to $h$}, if
$\norm{J}_{L^\infty(A)} < \infty$ and for any $\epsi >0$ there exists
$\delta=\delta(J, A, h) >0$ such that for any $x\in A$ one can find a
$C^1$-diffeomorphism $\phi: B(x, \delta) \to B(0, \delta)$ from the ball
$B(x, \delta)\deff \{y\in X \;:\; \dist_h(x,y) <\delta \}$ onto the standard
ball in $\cc^n$ with the standard metric $h\st$ such that
$$
\norm{J -\phi^*J\st}_{L^\infty(B(x,\delta)\cap A)} +
\norm{h -\phi^*h\st}_{L^\infty(B
(x,\delta)\cap A)}
\le\epsi.
\eqno(5.1.1)
$$
}

Roughly speaking, this means that on the set $A$ we can $\norm\cdot _{L^\infty}
$-approximate $J$ by an integrable structure in $h$-metric balls of a radius
independent of $x\in A$. The function $\mu(J,A,h)$ whose value at $\epsi>0$ is
the biggest possible $\delta\le1$ with the above property is called the {\it
modulus of uniform continuity of $J$ on $A$}. Note that every continuous almost
complex structure $J$ is always uniformly continuous on {\it relatively
compact} subsets $K \Subset X$.

\smallskip
Let $J^*$ be a continuous almost complex structure on $X$ and $A\subset X$
a subset. Assume that $J^*$ is uniformly continuous on $A$ and denote
by $\mu_{J^*}= \mu( J^*,A, h)$ the modulus of uniform continuity of
$J^*$ on $A$. Recall that in {\sl Lemma 2.4.1} we proved the following

\state Proposition 5.1.1. {\sl (First A priori Estimate). \it For every $p$ 
with
$2< p<\infty$ there exists an $\eps_1 =\eps_1(\mu_{J^*}, A, h)$ (independent
of $p$) and $C_p=C(p,\mu_{J^*}, A, h)$ such that for any continuous almost
complex structure $J$ with $\norm{ J- J^*} _{L^\infty(A)}<\eps_1 $ and
for every $J$-holomorphic map $u\in C^0\cap L^{1,2}(\Delta ,X)$, satisfying
$u(\Delta )\subset A$ and $\norm{du} _{L^2( \Delta )} <\eps_1 $, one has
the estimate
$$
\norm{du}_{L^p({1\over 2}\Delta )}\le C_p\cdot
\norm{du}_{L^2(\Delta )}.
\eqno(5.1.2)
$$
}

\smallskip
\state Definition 5.1.2. {\sl Define a {\it cylinder} $Z(a,b)$ by
$Z (a,b) \deff S^1 \times [a,b]$, equipping it with coordinates
$\theta \in [0,2\pi]$, $t\in[a,b]$, with the metric 
$ds^2= d\theta^2 + dt^2$ and
the complex structure $J\st({\d \over \d\theta}) ={\d \over \d t}$. Denote
$Z_i\deff Z(i-1,1)= S^1 \times [i-1,i]$.}

\smallskip
Let $J^*$ be some continuous almost complex structure on $X$ and $A$ a
subset of $X$ such that $J^*$ is uniformly continuous on $A$. Let
$\mu_{J^*}$ denote the modulus of uniform continuity of $J^*$ on $A$.

\smallskip
\state Lemma 5.1.2. {\sl (Second A priori Estimate). \it There exist constants
$\gamma \in \,]\,0,1\,[$ and $\eps_2=\eps_2(\mu_{J^*}, A,h)>0$ such that for 
any $J$ with $\norm{ J- J^*} <\eps_2$ and every $J$-holomorphic map $u:Z(0,5)
\to X$ with $u(Z(0,5)) \subset A$ the condition $\norm{du}_{ L^2(Z_i)}<\eps_2$
for $i=1,\ldots,5$ implies
$$
\norm{du}^2_{L^2(Z_3)}\le {\gamma \over 2}\bigl(\norm{du}^2_{L^2(Z_2)} +
\norm{du}^2_{L^2(Z_4)}\bigr).
\eqno(5.1.3)
$$
}

\state Proof. Take $\eps_2>0$ small enough such that $\mu_{J^*}(\epsi_2) <
\epsi_1$, where $\eps_1$ is the constant from {\sl Lemma 5.1.1}. Then for any
$A' \subset A$ the condition $\diam(A') \le \epsi_2$ implies that $\osc(J^*,
A') \le \eps_1$. Due to {\sl Lemma 5.1.1}, we may assume that $u(Z_i)\subset 
B$ for $i=2,3,4$, where $B$ is a small ball in $\rr^{2n} =\cc^n$ with the
structure $J\st$. Moreover, we may assume that $\norm{ J^* - J\st}_{
L^\infty(B)} \le \eps_1$.

Find $v\in C^0\cap L^{1,2}(Z(1,4),\cc^n)$ such that $\dbar_{J\st}v=0$ and
$\norm{du - dv}_{L^2(Z(1,4))}$ is minimal. We have
$$
\norm{  \dbar_{J\st}(u-v)}_{L^2(Z_i)} = \norm{(J\st-J(u))\d_y u
}_{L^2(Z_i)}
\le \norm{J\st - J}_{L^\infty(B)}\norm{du}_{L^2(Z_i)}.
$$
So for $i=2,3,4$ we get
$$
\norm{du - dv}_{L^2(Z_i)}\le C\norm{J\st -J}_{L^\infty(B)}
\norm{ du
}_{L^2(Z(1,4))}.\eqno(5.1.4)
$$
Now let us check the inequality (4.7) for $v$. Write $v(z)= \Sigma_{k=
-\infty} ^\infty v_k e^{k(t+i\theta)}$. Then $\norm{dv}^2_{L^2 (S\times \{
t\} )}=4\pi \Sigma_{k=-\infty}^\infty k^2| v_k|^2 e^{2kt}$. Since obviously
$$
\int_2^3e^{2kt}\le {\gamma_1\over 2}
\left( \int_1^2e^{2kt}dt+\int_3^4e^{2kt}dt \right)
%\eqno(3.8A)
$$
for all $k \not=0$ with $\gamma_1={2\over e^2}$, one gets the required
estimate for all holomorphic $v$.

Using (5.1.4) with $\norm{J\st - J}_{L^\infty}$ sufficiently small, we conclude
that the estimate (5.1.3) holds for $u$ with appropriate $\gamma > \gamma_1$.
\qed

\smallskip
\state Corollary 5.1.3. {\it Let $X$, $h$, $J^*$, $A$, and the constants
$\eps_2$ and $\gamma$ be as in {\sl Lemma 5.1.2}. Suppose that $J$ is a
continuous almost complex structure on $X$ with $\norm{J-J^*}_{L^\infty(A)}
<\eps_2$ and $u \in C^0\cap L^{1,2}(Z(0,l),X)$ a $J$-holomorphic map such
that $u(Z)\subset A$ and $\norm{du}_{ L^2(Z_i)}<\eps_2$ for any $i=1,\ldots
,l$. Let $\lambda>1$ be (the uniquely defined) real number with
$\lambda = {\gamma \over 2} (\lambda^2+ 1)$.

Then for $2\le k\le l-1$ one has
$$
\norm{du}^2_{L^2(Z_k)} \le \lambda^{-(k-2)} \cdot \norm{du}^2_{L^2(Z_2)}
+ \lambda^{-(l-1-k)} \cdot \norm{du}^2_{L^2(Z_{l-1})}.
\eqno(5.1.5)
$$
}

\state Proof. The definition of $\lambda$ implies that for any $a_+$ and
$a_-$ the sequence $y_k \deff a_+ \lambda^k + a_- \lambda^{-k}$ satisfies the
recurrent relation $y_k = {\gamma\over2}(y_{k-1} + y_{k+1})$. In particular,
so does the sequence
$$
A_k \deff
\msmall{ \lambda^{-(k-2)} - \lambda^{6-2l+ k-2}
\over 1- \lambda^{6-2l} } \norm{du}^2_{L^2(Z_2)}
+ \msmall{ \lambda^{-(l-1-k)} - \lambda^{6-2l+ l-1-k}
\over 1- \lambda^{6-2l} } \norm{du}^2_{L^2(Z_{l-1})},
$$
which is determined by the values $A_2 = \norm{du}^2_{L^2(Z_2)}$ and $A_{l-1}
= \norm{du}^2_{L^2(Z_{l-1})}$.

We claim that for $2\le k\le l-1$ one has the estimate $\norm{du}^2 _{L^2(
Z_k)} \le A_k$, which is obviously stronger than (5.1.5). Suppose that there
exists a $k_0$ such that $2\le k_0\le l-1$ and
$\norm{du}^2 _{L^2(Z_{k_0})}>A_{k_0}$.
Choose $k_0$ so that the difference $\norm{du }^2_{L^2(Z_{k_0})}-A_{k_0}$
is maximal. By {\sl Corollary 2.5.1} and by our recurrent definition of $A_k$, 
we have $2< k_0 < l-1$ and
$$
\norm{du }^2_{L^2(Z_{k_0})}-A_{k_0}\le {\gamma \over 2}(
\norm{du }^2_{L^2(Z_{k_0+1})}-A_{k_0+1}+
\norm{du }^2_{L^2(Z_{k_0-1})}-A_{k_0-1})\le
$$
$$
\le {\gamma \over 2}2(\norm{du }^2_{L^2(Z_{k_0})}-A_{k_0}).
$$
The second inequality follows from the fact that $\norm{du}^2_{L^2(
Z_{k_0})} - A_{k_0}$ is maximal. This gives a contradiction. \qed

\bigskip\noindent
{\bigsl 5.2. Removal of Point Singularities.}

\smallskip
An immediate corollary of this estimate is the following improvement of the
Sacks-Uhlenbeck theorem about removability of a point singularity, see [S-U]
and [G].

\smallskip
\state Corollary 5.2.1. {\sl (Removal of Point Singularities). 
\it Let $X$ be a
manifold with a Riemannian metric $h$, $J$ a continuous almost complex
structure, and $u:(\check\Delta, J\st)\to (X,J)$ a $J$ - holomorphic map
from the punctured disk. Suppose that

\item\sli $J$ is uniformly continuous on $A \deff u(\check\Delta)$ \wrt
$h$ and the closure of $A$ is $h$-complete;

\item \slii there exists $i_0$ such that, for all annuli $R_i\deff\{ z\in \cc
:{1\over e^{i+1}}\le | z| \le {1\over e^i}\} $ with $i\ge i_0$, one has
$\norm{du}^2_{L^2(R_i)}\le \eps_2$, where $\eps_2$ is defined in 
{\sl Lemma 5.1.2}.

\smallskip\noindent
Then $u$ extends to the origin.
}

\smallskip
Condition \sli is automatically satisfied if $A= u(\check\Delta)$ is
relatively compact in $X$. Condition \slii of ``slow growth'' is clearly weaker
than just the boundedness of the area, see, e.g.,\ [S-U], [G]. It is sufficient
to have $\lim_{i\lrar\infty} \norm{du}^2_{L^2(R_i)} =0$, whereas boundedness
of the area means $\sum_{i=1}^\infty \norm{du}^2_{L^2(R_i)} <\infty$.

\state Proof. The exponential map $\exp(t,\theta) \deff e^{-t+i\theta}$
defines a biholomorphism between the infinite cylinder $Z(0,\infty)$ and
the punctured disk $\check\Delta$, identifying every annulus $R_i$ with
the cylinder $Z(i,i+1)$. Applying {\sl Corollary 5.1.3} to the map
$u \scirc \exp$ on cylinders $Z(i_0,l)$ and setting $l\lrar \infty$, we obtain
the estimate
$$
\norm{du}^2_{L^2(R_i)} \le \lambda^{-(i-i_0)} \cdot
\norm{du}^2_{L^2(R_{i_0})}, \qquad i>i_0.
$$
Using this and {\sl Lemma 5.1.1} we conclude that $\diam(u(R_i)) \le C\cdot
\lambda^{-i/2}$ for $i>i_0$. Since $\sum \lambda^{-i/2} <\infty$, $u$ extends
continuously into $0\in \Delta$. \qed

\medskip
In the proof of the compactness theorem we shall use the following corollary
of {\sl Corollary 2.5.1}. Let $X$ be a manifold with a Riemannian metric 
$h$, $J$ a continuous almost complex structure on $X$, $A\subset X$ a closed 
$h$-complete
subset such that $J$ is $h$-uniformly continuous on $A$. Furthermore, let
$\{J_n\}$ be a sequence of almost complex structures uniformly converging to
$J$, $\{l_n\}$ a sequence of integers with $l_n\to \infty$ and $u_n:Z(0,l_n)
\to X$ a sequence of $J_n$-holomorphic maps.

\medskip
\state Lemma 5.2.2. 
{\it Suppose that $u_n(Z(0,l_n))\subset A$ and $\norm{du_n}_{L^2 (Z_i)}\le 
\eps_2$ for all $n$ and $i\le l_n$. Take a sequence $k_n\to \infty$ such that 
$k_n<l_n-k_n\to \infty$. Then

\noindent
{\sl1)} $\norm{du_n}_{L^2(Z(k_n,l_n-k_n))}\to 0$ and $\diam\bigl(u_n
(Z(k_n,l_n-k_n))\bigr)\to 0$;

\noindent
{\sl2)} if, in addition, all images $u_n(Z(0,l_n))$ are contained in some 
bounded subset of $X$, then there is a subsequence $\{ u_n \}$, still denoted 
$\{u_n \}$ such that both $u_n |_{Z(0,k_n)}$ and $u_n |_{Z(k_n,l_n)}$ converge
in $L^{1,p}$-topology on compact subsets in $ \check \Delta \cong Z(0, +\infty )$
to $J$-holomorphic maps $u^+_\infty: \check \Delta \to X$ and $u^-_\infty:
\check \Delta \to X$. Moreover, both $u^+_\infty$ and $u^-_\infty$ extend to
the origin and $u^+_\infty(0)= u^-_\infty(0)$.
}

\state Remarks.~1. The punctured disk $\check\Delta$ with the standard
structure $J_{\Delta } {\d \over \d r}={1\over r}{\d \over \d \theta }$ is
isomorphic to $Z(0,\infty )$ with the structure $J_Z{\d \over \d t}=-{\d \over
\d \theta }$ under a biholomorphism $(\theta ,t)\mapsto e^{-t+ \isl\theta}$.
Thus, statement (2) of this corollary is meaningful.

\state 2. {\sl Lemma 5.2.2} describes explicitly how the sequence of
$J_n$-holomorphic maps of the cylinders of growing conformal radii
converges to a $J$-holomorphic map of the standard node.

\smallskip
\state Lemma 5.2.3. {\it There is an $\eps_3 =\eps_3(\mu_{J_\infty},A,h)$ 
such that for any continuous almost-complex structure $J$ on $X$ with $\norm{J 
- J_\infty}_{L^\infty}\le \eps_3 $ and any non-constant $J$-comp\-lex
sphere $u: \cc\pp^1 \to X$, $u(\cc\pp^1)\subset A$ one has the inequalities
$$
\area(u(\cc\pp^1 ))\ge \eps_3
\qquad\text{and}\qquad
\diam(u(\cc\pp^1))\ge \eps_3.
$$
}

\state Proof. Let $\eps_1$ be the constant from {\sl Lemma 5.1.1}. 
Suppose that
$\area u(\cc\pp^1)=\norm{du}^2_{L^2(\cc\pp^2)}\le \eps_1^2$. Cover $\cc\pp^1$
by two disks $\Delta_1$ and $\Delta_2$. By (3.1) and the Sobolev imbedding
$L^{1,p} \subset C^{0,1-{2\over p}}$, we obtain that $\diam(u(\Delta_1))$ and
$\diam(u (\Delta_2) )$ are smaller than $const\cdot \eps_1$. Thus the diameter
of the image of the sphere is smaller than $const\cdot \eps_1$.

Therefore, we can suppose that the image $u(S^2)$ is 
contained in the coordinate
chart, \i.e., in a subdomain in $\cc^n$, and the structures $J$ and $J_\infty$ are
$L^\infty$-close to a standard one. Consider now $u:S^2\to U\subset \cc^n$
as a solution of the linear equation
$$
\d_xv(z) + J(u(z))\cdot\d_yv(z) = 0\eqno(5.2.1)
$$
on the sphere. The operator $\dbar _J(v) = \d_xv(z) + J(u(z))\cdot
\d_yv(z)$ acts from $L^{1,p}(S^2, \cc^n)$ to $L^p(S^2, \cc^n)$ and is a small
perturbation of the standard $\dbar$-operator. Note that the standard $\dbar$
is surjective and Fredholm. Thus small perturbations are also surjective
and Fredholm, having the kernel of the same dimension. But the kernel of
$\dbar$ consists of constant functions. Since all constants are in the
kernel of (5.2.1), our $u$ should be a constant map.

We have proved that if the area or a diameter of a $J$-holomorphic map is
sufficiently small then this map is constant.\qed

\smallskip
\state Remark. The same statement is true for the curves of arbitrary genus
$g$. In that case, in addition to the estimate (5.1.2), one should also use the
estimate (5.1.3). This yields the existence of an $\eps $ which depends on $g$
(and, of course, on $X$, $J$, and $K$), but not on the complex structure on the
parameterizing surface.

\bigskip\noindent
{\bigsl 5.3. Compactness for Curves with Free Boundary.}

\smallskip
In this section we give a proof of the Gromov compactness theorem 
for the curves
with boundaries of fixed finite topological type and without boundary
conditions on maps. The case of closed curves is obviously included in this.

Throughout this section we assume that the following setting holds.

\sl Let $X$ be a manifold with a Riemannian metric $h$, $J_\infty$ a continuous
almost complex structure on $X$, $A\subset X$ an $h$-complete subset, $\{ C_n
\}$ a sequence of nodal curves parameterized by a real surface $\Sigma$ with
parameterizations $\delta_n:\Sigma \to C_n$, and $u_n:(C_n,j_n) \to (X, J_n)$
a sequence of  holomorphic  maps. Further, $J_\infty$ is $h$-uniformly
continuous on $A$, $J_n$ are also continuous and converge to $J_\infty$,
$h$-uniformly on $A$, $u_n(C_n)\subset A$ for all $n$. \rm

\smallskip
Let us explain the main idea of the proof of {\sl Theorem 1.1}. The Gromov
topology on the space of stable curves over $X$ is introduced in order to
recover natural convergence of sequences $(C_n,u_n)$ of  complex 
curves of bounded area which do not converge in the ``strong" 
(\i.e., $L^{1,p}$-type) sense. The are two reasons
for this. The first is that a sequence of (say, smooth) curves $C_n$ could
diverge in an appropriate moduli space and the second is a phenomenon of
``bubbling". In both cases one has to deal with the appearance of 
new nodes, \i.e.,
with a certain degeneration of the complex structure on curves. The ``model''
situation of {\sl Lemma 5.2.2} describes a convergence of ``long cylinders''
$u_n :Z(0, l_n) \to X$, $l_n \lrar \infty$, to a node $u_\infty: \cala_0 \to
X$. In our proof we cover curves $C_n$ by pieces which are either ``long
cylinders'' converging to nodes or have the property that complex structures
and maps ``strongly'' converge. Here the ``strong'' convergence means the
usual one, \i.e., \wrt the $C^\infty$-topology for complex structures, and \wrt
the $L^{1,p}$-topology with some $p>2$ for maps. In fact, the strong
convergence of maps is equivalent to the uniform one, \i.e., \wrt the
$C^0$-topology, and implies further regularity in the case when $J_n$
and $J_\infty$ have more smoothness. One consequence of this is that we remain
in the category of nodal curves. Another is that we treat degeneration
of a complex structure on $C_n$ and the
``bubbling'' phenomenon  in a uniform framework of ``long cylinders''.

\smallskip
For the proof we need some additional results.

\state Lemma 5.3.1. \it For any $R>1$ there exists an $a^+= a^+(R)>0$ with 
the following property. For any cylinder $Z = Z(0,l)$ with $0 <l\le +\infty$
and any annulus $A \subset Z(0,l)$, which is adjacent to $\d_0 Z = S^1\times
\{0\}$ and has a conformal radius $R$, one has $Z(0,a^+) \subset A$.
\rm

\state Proof. Without loss of generality we may assume that $l=+\infty$ and
identify $Z$ with the punctured disk $\check\Delta$ via the exponential map
$(t+ \isl \theta) \mapsto e^{ -t+ \isl \theta}$ such that $\d_0 Z$ is
mapped onto $S^1= \d\Delta$.

Suppose that the statement is false. Then there would exist holomorphic
imbeddings $f_n : A(1, R) \to \check\Delta$ and points $a_n \in \Delta \bs
f_n(A(1, R))$ such that $f_n(A(1, R))$ are adjacent to $\d\Delta$ and $a_n
\lrar a\in \d\Delta$. Passing to a subsequence, we may assume that $\{f_n\}$
converges uniformly on compact subsets in $A(1, R)$ to a holomorphic map
$f: A(1, R) \to\Delta$.

If $f$ is not constant, then $f(A(1, R))$ must contain some annulus $\{ b
<|z|<1 \}$ with $b<1$. But then $\{ \sqrt{b} <|z|<1\} \subset f_n(A(1, R))$
for $n>\!>1$, which is a contradiction.

If $f$ is constant, then the diameter of images of the middle circle
$\gamma\deff \{ |z|= \sqrt{R} \} \subset A(1, R)$ must converge to $0$. But
$\diam(f_n(\gamma)) \ge \dist(0, a_n) \lrar1$. The obtained contradiction
finishes the proof. \qed

\medskip
For the proof of {\sl Theorem 1.1} we need a special covering of\/ $\Sigma$
which will be constructed in the following theorem.

\state Theorem 5.3.2. {\it Under the conditions of {\sl Theorem 1.1}, after
passing to a subsequence, there exist a finite covering $\calv$ of\/ $\Sigma$
by open sets $V_\alpha$ and parameterizations $\sigma_n:\Sigma \to C_n$
such that

{\sl(a)} all $V_\alpha$ are either disks, or annuli or pants;

{\sl(b)} for any boundary circle $\gamma_i$ of\/ $\Sigma$ there is some
annulus $V_\alpha$ adjacent to  $\gamma_i$;

{\sl(c)} $\sigma_n^*j_n\ogran_{V_\alpha}$ does not depend on $n$ if
$V_\alpha$ is a disk, pants or an annulus adjacent to a boundary circle
of\/ $\Sigma$;

{\sl(d)} all non-empty intersections $V_\alpha \cap V_\beta$ are annuli,
where the structures $\sigma^*j_n$ are independent of $n$;

{\sl(e)} if $a$ is a node of $C_n$ and $\gamma^n_a=\sigma_n\inv (a)$ the
corresponding circle, then $\gamma_a^n=\gamma_a$ does not depend on $n$ and is
contained in some annulus $V_\alpha$, containing only
one such ``contracting" circle for any $n$; moreover, the structures
$\sigma_n^*j_n\ogran_{V_\alpha\bs \gamma_a}$ are independent of $n$;

{\sl(f)} if $V_\alpha$ is an annulus and $\sigma_n( V_\alpha)$ doesn't 
contains nodes, then the conformal radii of $\sigma_n( V_\alpha)$ converge to some
positive $R_\alpha^\infty>1$ or to $+\infty$;

{\sl(g)} if for initial parameterizations $\delta_n$ and fixed annuli
$A_i$, each adjacent to the boundary circle $\gamma_i$ of\/ $\Sigma$,
the structures $\delta_n^* j_n\ogran_{A_i}$ do not depend on $n$, then
the new parameterizations $\sigma_{n}$ can be taken equal to $\delta_{n}$
on some subannuli $A'_i \subset A_i$ also adjacent to $\gamma_i$.  }

\state Proof. We shall prove the properties {\sl(a)--(f)}. The property
{\sl(g)} will follow from {\sl Lemma 5.3.3} below.

There are four cases where the existence of such a covering is obvious. If
all $C_n$ are disks or annuli without nodal points, there is nothing to
prove. In the third case each $C_n$ is a sphere, and we cover it by two
disks.

In the forth case each $C_n$ is a torus without marked points. Then
any complex torus can be represented by the form $\cc{\bigm/}(\zz+
\tau \zz)$ with $|\re\tau| \le {1\over 2}$ and $\im\tau > {1\over 2}$.
Considering the map $z\in \cc \mapsto e^{2\pi\isl z} \in \check\cc \deff
\cc\bs\{0\}$, we represent $(T^2,j)$ as the quotient $\cc\bigm/ \{ z \sim
\lambda^2 z\}$ with $\lambda = e^{\pi\isl\tau}$, so that $|\lambda|<
e^{-\pi/2}< {1\over 3}$. The annuli $\{ {|\lambda| \over2} <|z| <1\}$ and
$\{{|\lambda|^2 \over2} <|z| <|\lambda| \}$ form the needed covering.

\smallskip
In all remaining cases we start with the construction of appropriate graphs
$\Gamma_n$ associated with some decomposition of $C_n$ into pants. {\sl Lemma
5.1.2} and a non-degeneration of the complex structure 
$j_n$ on $C_n$ shows that
lengths of all boundary circles of all non-exceptional components $C_{n,i}$ of
$C_n$ are uniformly bounded from above. At this point we make the following

\state Remark. The Collar Lemma from [Ab], Ch.II, \S\.3.3 yields the existence
of the universal constant $l^*$ such that for any simple geodesic circles
$\gamma'$ and $\gamma''$ on $C_{n,i}$ the conditions $\ell(\gamma') <l^*$ and
$\ell( \gamma'') <l^*$ imply $\gamma' \cap \gamma'' = \emptyset$. We shall
call geodesic circles $\gamma$ with $\ell(\gamma) <l^*$ {\sl short geodesics}.

\smallskip
The fact that $(C_n,u_n)$ are $J_n$-complex and of bounded area shows
that $C_n$ have a uniformly bounded number of components. Indeed, the number of
exceptional components, which  are spheres and disks, is bounded by the energy
(see {\sl Lemma 5.2.3}), and  the  number of boundary circles of $C_n$ is 
equal to that of $\Sigma$.
%%%%%%%%%%%%%%%%%%%%%Question%%%%%%%%%%%%%%%%%%%

Furthermore, the operation of contracting a circle on $\Sigma$ to a nodal point
either diminishes the genus of some component of $C_n$ or increases
 the number of
components. Thus, the number of nodal points on $C_n$ and the total number of
marked points on its components are also uniformly bounded. This implies that
the number of possible topological types of components is  finite.

In this situation the Teichm\"uller theory (see [Ab], Ch.II, \S\.3.3) states
the existence of decomposition of every non-exceptional component $C_{n,i}\bs
\mapo$ into pants with the following properties:

\sli Every short geodesic is a boundary circle of some pants of the
decomposition.

\slii The intrinsic length of every boundary circle is bounded from above by
a (uniform) constant depending only on an upper bound of lengths of boundary
circles and possible topological types of $C_{n,i} \bs \mapo$.

Having decomposed all $C_{n,i} \bs \mapo$ into pants, we associate with every
curve $C_n$ its graph $\Gamma_n$. As was noted above, the number of
vertices and edges of $\Gamma_n$ is uniformly bounded. Thus, after passing to
a subsequence, we can assume that all $\Gamma_n$ are isomorphic to each other
(as marked graphs). Denote this graph by $\Gamma$. Now, the parameterizations
$\sigma_n: \Sigma \to C_n$ can be found in such a way that the decompositions
of $C_{n,i} \bs\mapo$ into pants define the same set $\bfgamma =\{ \gamma
_\alpha \}$ of circles on $\Sigma$ and induce the same decomposition $S\bs
\cup_\alpha \gamma_\alpha = \cup_j S_j$ with the graph $\Gamma$.

\smallskip
By our construction of the graph $\Gamma$, each edge of $\Gamma$ corresponds
either to a circle in $\Sigma$ contracted by every parameterization $\sigma_n$
to a nodal point, or to a circle mapped by every $\sigma_n$ onto a geodesic
circle separating two pants. Furthermore, each tail of $\Gamma$ corresponds
to a boundary circle of $\Sigma$. Thus, we shall use the same notation $\gamma
_\alpha$ for an edge or a tail of $\Gamma$ and for the corresponding circle
on $\Sigma$. If $\sigma_n(\gamma_\alpha)$ is a boundary circle of some pants
$C_j$, then the intrinsic length $\ell_{n,\alpha}= \ell_n (\gamma_\alpha)$ of
$\sigma_n(\gamma_\alpha)$ is well- defined. This happens in the following two
cases:

{\sl \.a)}~$\sigma_n(\gamma_\alpha)$ separates two pants, or else

{\sl \.b)}~$\gamma _\alpha$ is a boundary circle $\Sigma$ and the irreducible
component of $C_n$ attached to $\sigma_n(\gamma_\alpha)$ is not a disk with a
single nodal point. Note that the appearance of these two cases is 
independent of $n$.

By our choice of $\gamma_\alpha$, the lengths $\ell_n (\gamma_\alpha)$ are
uniformly bounded from above. Passing to a subsequence, we may assume that
for any fixed $\alpha$ the sequence $\{ \ell_{n,\alpha} \}$ converges to
$\ell_{\infty,\alpha}$.

\smallskip
As one can expect, the condition $\ell_{n, \alpha} \lrar 0$ means that the
circle $\gamma_\alpha$ is shrunk to a nodal point on the limit curve. We
shall prove the statement of the theorem by induction in the number $N$ of
those circles $\gamma_\alpha$ for which $\ell_{\infty,\alpha} =0$.

\smallskip
The case $N=0$, where there are no such circles, is easy. Passing to a
subsequence, we may assume that the Fenchel-Nielsen coordinates $(\bfell_n,
\bfvartheta_n)$ of any non-exceptional component $C_{n,i}$ of $C$ converge to
the Fenchel-Nielsen coordinates $(\bfell_\infty, \bfvartheta_\infty)$ of some
smooth curve $C_{\infty, i}$ with marked points. Gluing together appropriate
pairs of marked points, we obtain a nodal curve $C_\infty$, which admits a
suitable parameterization $\sigma_\infty: \Sigma \to C_\infty$ and has the
same graph $\Gamma$. {\sl Lemma 5.1.3} shows that for $n>\! >1$ the curves
$C_n$ can be obtained from $C_\infty$ by deformation of the transition
functions for the intrinsic local coordinates on non-exceptional components
of $C_\infty$. Note that such a deformation can be realized as a deformation
of the operator $j_\infty$ of a complex structure on $C_\infty$, 
localized in small
neighborhoods of circles $\sigma_\infty(\gamma_\alpha)$, see Fig.~10. In the
case where $\gamma_\alpha$ is a boundary circle we may additionally assume
that the annulus, where $j_n$ changes, lies away from $\gamma_\alpha$. Now
the existence of the covering with desired properties is obvious.

\medskip
\vbox{\nolineskip\xsize.54\hsize%
\putm[.01][.38]{\underbrace{\hskip.35\xsize}_{V_\beta}}%
\putm[.54][.44]{\underbrace{\hskip.46\xsize}_{V_\gamma}}%
\putm[.27][.06]{\overbrace{\hskip.36\xsize}^{V_\alpha}}%
\putm[.27][.3]{\underbrace{\hskip.09\xsize}_{W_\beta}}%
\putm[.53][.3]{\underbrace{\hskip.09\xsize}_{W_\gamma}}%
\putt[1.05][0]{\advance\hsize-1.05\xsize\parindent=0pt%
\centerline{Fig.~10. }
\smallskip
$V_\alpha$ and $V_\beta$ represent elements of the covering where the complex
structure is constant. The change of complex structure is done in the shaded
part of  $V_{\gamma}$. $W_\beta\deff V_\alpha \cap V_\beta$ and $W_\gamma
\allowbreak \deff V_\alpha \cap V_\gamma$ represent the annuli with the
constant complex structure.
}
\noindent
\epsfxsize=\xsize\epsfbox{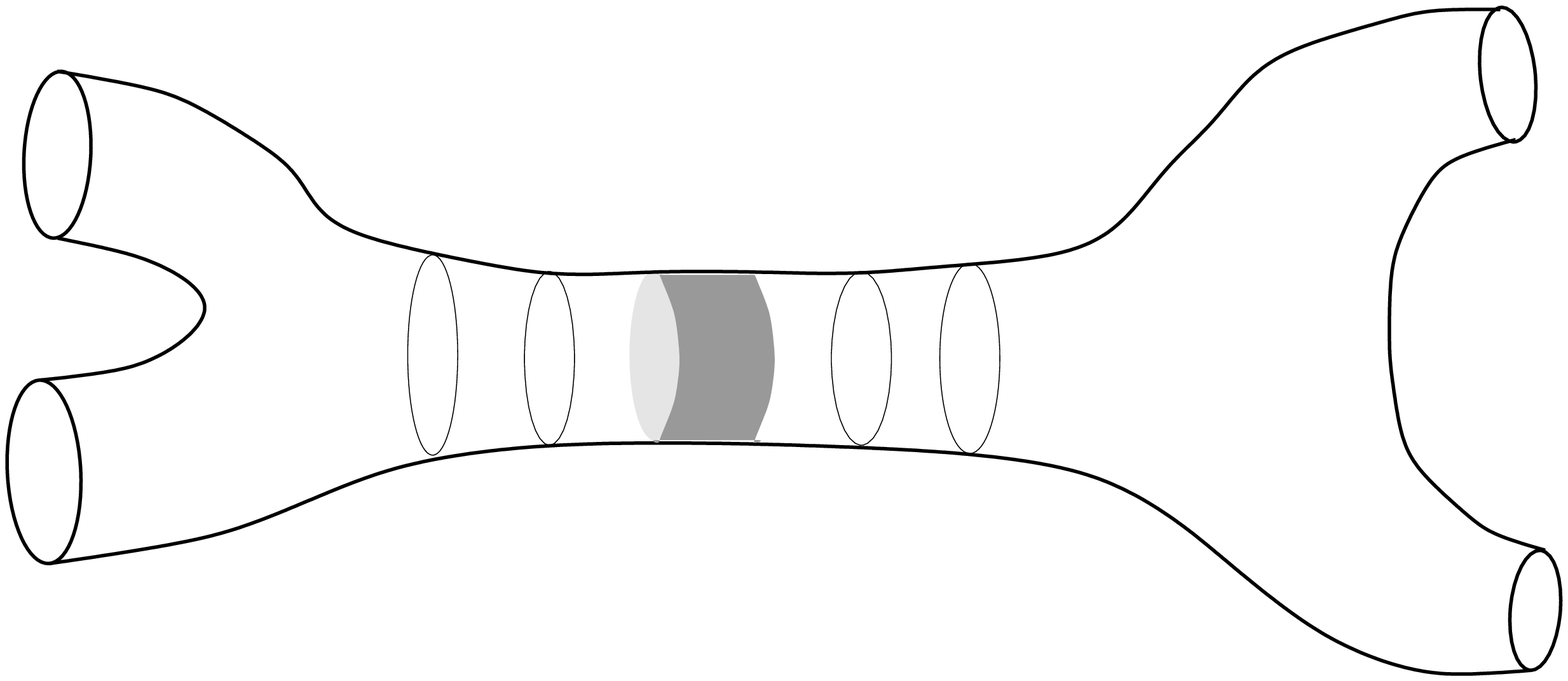}
}
\bigskip\smallskip

\medskip
Let us now consider the general case where the number $N$ of ``shrinking
circles'' is not zero. Take a circle $\gamma_\alpha$ with $\ell_\infty
(\gamma _\alpha) =0$. Let $S_j$ be pants adjacent to $\gamma_\alpha$.
Consider the intrinsic coordinates $\rho_\alpha$ and $\theta_\alpha$ at
$\sigma_n(\gamma_\alpha)$ and the annuli
$$
\eqalign{
A_{n,\alpha,j} &\deff \{\, (\rho
_\alpha, \theta_\alpha) \in \sigma_n(S_j) \;:\; 0\le \rho_\alpha \le {\pi^2
\over \ell_{n,\alpha} } - {2\pi \over a^*} \;\}
\cr
A^-_{n,\alpha,j} &\deff
\{\, (\rho_\alpha, \theta_\alpha) \in \sigma_n(S_j) \;:\;
0\le \rho_\alpha \le {\pi^2 \over \ell_{n,\alpha} } - {2\pi \over a^*}-1
\;\},
}
$$
adjacent to $\sigma_n( \gamma_\alpha)$. Note that ${\pi^2 \over \ell_{n,
\alpha} }- {2\pi \over a^*}$ (resp.\ ${\pi^2 \over \ell_{n,\alpha} }- {2\pi
\over a^*}-1$) is the logarithm of the conformal radius of $A_{n, \alpha, j}$
(resp.\ of $A^-_{n, \alpha,j}$). Consequently, we can use {\sl Lemma 4.3.2} 
to show that these annuli are well-defined.

If $\gamma_\alpha$ is a boundary circle, we set $C^-_n \deff C_n \bs A^-_{n,
\alpha,j}$. Otherwise $\gamma_\alpha$ separates two pants, say $S_j$ and
$S_k$. Then we define in a similar way the annuli $A_{n,\alpha, k}$ and $A^-
_{n,\alpha,k}$, set $A_{n,\alpha} \deff A_{n, \alpha,j} \cup A_{n, \alpha,
k}$ and $A^-_{n,\alpha} \deff A^-_{n,\alpha,j} \cup A^-_{n, \alpha, k}$ and
put $C^-_n \deff C_n \bs A^-_{n,\alpha}$.

The parameterizations $\sigma_n: \Sigma \to C_n$ can be chosen in such a way
that the annuli $\sigma_n\inv(A^-_{n,\alpha,j})$ (resp.\ $\sigma_n\inv(A^-
_{n, \alpha,k})$) define the same annulus $A^-_{ \alpha,j}$ (resp.\ $A^-
_{\alpha, k}$) on $\Sigma$. Let $\gamma^-_{ \alpha,j}$ (resp.\ $\gamma^-
_{\alpha,k}$) denote its boundary circles different from $\gamma_\alpha$.
Thus, the curves $C^-_n$ are parameterized by a real surface $\Sigma^- \deff
\Sigma \bs A^-_{ \alpha,j}$ (resp.\ $\Sigma^- \deff \Sigma \bs (A^-_{\alpha,
j} \cup A^- _{\alpha, k}$), and the restrictions of $\sigma_n$ can be chosen
as parameterization maps. Thus, the decompositions of components of $C_n$ into
pants define the combinatorial type of decompositions of components of
$C^-_n$ into pants. Moreover, the corresponding graph $\Gamma^-$ will be the
same for all $C^-_n$. It coincides with $\Gamma$ if $\gamma_\alpha$ is a
boundary circle. Otherwise $\Gamma^-$ can be obtained from $\Gamma$ by
replacing the edge corresponding to $\gamma_\alpha$ by two new tails for 
two new boundary components.

Take some non-exceptional component $C^-_{n,i}$ of $C^-_n$ and decompose it
into pants according to graph $\Gamma^-$ in the canonical way, so that the
boundary circles of the obtained pants are geodesic. Note that even if the
constructed pants are in combinatorial one-to-one correspondence with the
pants of $C_n$, the intrinsic metric on $C^-_{n,i}$ and the obtained geodesics
circle $\gamma^-_\beta$ differ from the corresponding objects on $C_{n,i}$

Nevertheless, we claim that for the obtained decomposition of $C^-_n$ the
intrinsic lengths are uniformly bounded from above (possibly by a new
constant) and that the sequence $\{ C^-_n \}$ has fewer ``shrinking circles''
than $\{ C_n \}$. The meaning of the above construction is the following.
The curves $C^-_n$ are obtained from $C^-_n$ by cutting off the annuli
$A^-_{n,\alpha,j}$ (and resp.\ the annuli $A^-_{n,\alpha,j}$). These annuli
are sufficiently long so that one ``shrinking circle'' disappears, but not
too long so that the complex structures of the curves $C^-_n$ remain
non-degenerating near the boundary.

\smallskip
Indeed, the complex structures on $C^-_n$ do not degenerate at the boundary
circle $\gamma^-_{ \alpha,j}$ (resp.\ at $\gamma^-_{\alpha,k}$), because
$C^-_n$ contain annuli $A_{n,\alpha,j} \bs A^-_{n,\alpha,j}$ (resp.\
$A_{n,\alpha,k} \bs A^- _{n,\alpha,k}$) of the constant conformal radius
$R=e>1$. This means that the lengths $\ell^-_n(\gamma^-_{\alpha,j})$ of
$\sigma_n( \gamma^- _{\alpha,j})$ (resp.\ $\ell^-_n(\gamma^-_{\alpha,k})$ of
$\sigma_n(\gamma^- _{\alpha,k})$) with respect to the intrinsic metrics on
$C^-_n$ are uniformly bounded.

On the other hand, the lengths $\ell^-_n(\gamma^-_{\alpha,j})$ (resp.\
$\ell^-_n(\gamma^-_{\alpha,k})$) are also uniformly bounded from below by a
positive constant. Otherwise, by {\sl Lemma 3.4}, after passing to a
subsequence, there would exist annuli $A_n\subset C^-_n$ of infinitely
increasing radii $R_n$, adjacent to $\sigma_n(\gamma^-_{\alpha,j})$ (resp.\
to $\sigma_n(\gamma^-_{ \alpha,k})$). The superadditivity of the logarithm of
the conformal radius of annuli, see [Ab], Ch.II, \S\.1.3, shows that the
conformal radius $R^+_n$ of the annulus $A^-_{n,\alpha,j} \cup A_n$ satisfies
the inequality $\log R^+_n \ge {\pi^2 \over \ell_{n,\alpha} }- {2\pi \over
a^*}-1 + \log R_n$, which contradicts {\sl Lemma 4.3.2}, part {\sl i)}.

\smallskip
Now we estimate the intrinsic lengths of boundary circles and the number of
``shrinking circles'' on $C^-_n$. Denote $\ell_n \deff \ell_n(\gamma_\alpha)$,
where $\gamma_\alpha$ is the circle on $\Sigma$ used in the above
constructions. Compute the width $L_n$ of $A_{n,\alpha,j} \bs A^- _{n, \alpha,
j}$ \wrt the intrinsic metric on $C_n$. Using $\ell_n \lrar0$, we obtain
$$
\eqalign{
L_n &= \int_{\rho={\pi^2 \over \ell_n }- {2\pi \over a^*}-1}
^{ {\pi^2 \over \ell_n }- {2\pi \over a^*}}
\left(\msmall{ {\ell_n \over 2\pi} \over \cos {\ell_n \rho \over 2\pi}}
\right) d\rho =
\left[ \log\;\cotan \left( {\pi \over 4} -{\ell_n \rho \over 4\pi} \right)
\right]_{\rho={\pi^2 \over \ell_n }- {2\pi \over a^*}-1}
^{{\pi^2 \over \ell_n }- {2\pi \over a^*}}
\cr
&= \log \msmall{ \cotan{\ell_n  \over2a^*} \over
\cotan({\ell_n \over 2a^*} + {\ell_n \over4\pi}) }
\approx \log \msmall{ {\ell_n  \over2a^*} + {\ell_n \over4\pi} \over
{\ell_n \over2a^* } } =
\log \left( 1 + \msmall{ a^* \over 2\pi } \right) >0.
}
$$

Thus we can find annular neighborhoods $A_n\subset C^-_n$ of
$\sigma_n(\gamma_\alpha)$ with constant conformal radius $R>1$.
%%%%%%%%%%%%%%%%%%%%%Question%%%%%%%%%%%%%%%%%%

Let $\gamma^-_{n, \alpha} \subset C^-_n$ be the geodesic circle corresponding
to $\gamma_\alpha$. It is the unique circle, which is homotopic to $\gamma_{n,
\alpha} = \sigma_n( \gamma_\alpha)$ and is geodesic \wrt the
intrinsic metric on $C^-_n$. Let 
$\ell^-_{n,\alpha}$ denote its $C^-_n$-intrinsic length. Using the monodromy
argument as in the proof of {\sl Lemma 5.1.2}, we see that every $A_n$ can be
isometrically imbedded into the annulus $A^+_n \deff (-{\pi^2\over \ell^-
_{n,\alpha}} , {\pi^2\over \ell^-_{n,\alpha}}) \times S^1$ with coordinates $
-{\pi^2\over \ell^-_{n,\alpha}} <\rho < {\pi^2\over \ell^-_{n,\alpha}}$, $0\le
\theta \le 2\pi$ and with metric $\Bigl({\ell^-_{n,\alpha} \over 2\pi} / 
\cos{\ell^-_{n,\alpha} \rho\over 2\pi}\Bigr)^2 (d\rho^2 + d\theta^2)$ of 
conformal radius $e^{2\pi^2/\ell^-_{n,\alpha}}$. The monotonicity of 
conformal radius of
annuli yields the upper bound $\ell^-_{n,\alpha} \le {2\pi^2 \over \log\,R}$.

The same argumentation shows that if $\ell_{\infty,\alpha} \not=0$, then
$\ell^- _{n,\alpha}$ also do not vanish. Indeed, if $\ell^-_{n,\alpha}\lrar 0$,
then we could find annular neighborhoods $A'_n\subset C^-_n \subset C_n$ of
$\gamma_{n,\alpha}$ with infinitely increasing radii $R'_n$. Using the 
homotopic
equivalence of $\gamma^-_{n,\alpha}$ with $\sigma_n(\gamma_\alpha)$ as above,
we would obtain the estimate $\ell_{n,\alpha} \le {2\pi^2 \over \log\,R'_n}
\lrar 0$, which is a contradiction.

\smallskip
Thus, we have shown that $C^-_n$, with the intrinsic metric and defined by
the $\Gamma^-$ decomposition, have uniformly bounded lengths of 
marked circles and
less ``shrinking circles'' than $C_n$. Using induction, we may assume that
the postulated covering of $\Sigma^-$ and parameterizations of $C^-_n$ by
$\Sigma^-$ exist. Since $C_n \bs C^-_n$ is an annulus of increasing conformal
radius, the statement of the theorem is valid for $C_n$. \qed

\state Lemma 5.3.3. \it Let $C_n$ be a sequence of complex annuli with
structures $j_n$, $\Sigma$ some fixed annulus and $\delta_n: \Sigma \to C_n$
parameterizations. Suppose that for some fixed annuli $A_1, A_2 \subset
\Sigma$ adjacent to the boundary circles of $\Sigma$ restrictions
$\delta_n^*j_n|_{A_i}$ do not depend on $n$.

Then one can find parameterizations $\sigma_n: \Sigma \to C_n$ such that
$\sigma_n$ coincide with $\delta_n$ on some (possibly smaller) annuli $A'_i$,
are also adjacent to the boundary circles of $\Sigma$ restrictions, and 
that

\smallskip
\sli if conformal radii $R_n$ of $C_n$ converge to $R_\infty <\infty$, then
$\sigma_n^*j_n$ converge to some complex structure;

\smallskip
\slii if conformal radii $R_n$ of $C_n$ converge to $\infty$, then for some
circle $\gamma \subset \Sigma$ structures $\sigma_n^*j_n$ converge on
compact subsets $K \Subset \Sigma \bs \gamma$ to a complex structure of 
disjoint
union of two punctured disks. Moreover, such a $\gamma$, an arbitrary
imbedded circle generating $\pi_1(\Sigma)$, can be chosen. \rm

\state Proof. Without loss of generality we may assume that $\Sigma =
A(1,10)$, $A_1=A(7,10)$, $A_2=A(1,4)$ and that a given circle $\gamma$ lies in
$A(3,7)$. Let $\delta_n : \Sigma \to C_n$ be the given parameterizations.
There exist biholomorphisms $\phi_n: C_n \to A(r_n,1)$ with $r_n\inv=R_n$
being the conformal radii of $C_n$ such that $\phi_n (\delta_n(A_1))$ is
adjacent to $\{ |z|=1\} = \d\Delta$ and $\phi_n (\delta_n(A_2))$ is adjacent
to $\{ |z|=r_n\}$. Define $\phi'_n(z) \deff {r_n \over \phi_n(z)}$.

Recall that structures $\delta_n^*j_n|_{A_i}$ are independent of $n$. We call
this structure $j$. Consider maps $\phi_n \scirc \delta_n : (A_1, j) \to
A(r_n,1) \subset \Delta$ and $\phi'_n \scirc \delta_n : (A_2, j) \to A(r_n,1)
\subset \Delta$. Passing to a subsequence we can suppose that $r_n \lrar
r_\infty <1$ and that maps $\phi_n \scirc \delta_n$, $\phi'_n \scirc \delta_n$
converge on $A_1\cup A_2$ to holomorphic maps $\psi: (A_1\cup A_2, j) \to
\Delta$ and $\psi': (A_1\cup A_2, j) \to \Delta$, respectively. This means
that maps $(\phi_n \scirc \delta_n, \phi'_n \scirc \delta_n): (A_1 \cup A_2)
\to \Delta^2$ take values in $\{(z,z') \in \Delta^2 \;:\; z{\cdot} z'= r_n
\}$ and converge to the map $(\psi, \psi'): (A_1 \cup A_2) \to \Delta^2$ with
values in $\{(z,z') \in \Delta^2 \;:\; z{\cdot}z' = r_\infty \}$.

Arguments from the proof of {\sl Lemma 5.3.1} show that $\psi(A_1)$ and
$\psi'(A_2)$ are annuli adjacent to $\d\Delta$. This implies that for $n>\!>
1$ there exist diffeomorphisms $(\psi_n, \psi'_n): \Sigma \to \{(z,z') \in
\Delta^2 \;:\; z{\cdot}z' = r_n \}$ such that $\psi_n \equiv \phi_n \scirc
\delta_n$ on $A(9,10)$, $\psi'_n(z) \equiv \phi'_n \scirc \delta_n(z)$ for
$z\in A(1,2)$, and $(\psi_n, \psi'_n)$ converge to $(\psi, \psi')$ on
$\Sigma$. Moreover, we may assume that $|\psi_n(t)|=|\psi'_n(t)|= \sqrt{r_n}$
for any $t\in \gamma$. This means that $(\psi_n, \psi'_n)(\gamma)$ lies on
the middle circle $\{ (z,z'): |z|= \sqrt{r_n},\; z'= {r_n\over z} \}$ of
$\{(z,z') \in \Delta^2 \;:\; z{\cdot}z' = r_n \}$.

Set $\sigma_n \deff \phi_n\inv \scirc \psi_n: \Sigma \to C_n$. Then, obviously,
$\sigma_n \equiv \delta_n$ on $A(1,2)$ and on $A(9,10)$, and $\sigma_n^*j_n
= (\psi_n,\psi')_n^* J\st \lrar (\psi,\psi')^* J\st$, where $J\st$
denotes the standard complex structure on $\Delta^2$.
\qed

\medskip
\state Proof of Theorem 4.1.1. Let $\{(C_n, u_n)\}$ be the sequence from the
hypothesis of the theorem. Then the condition {\sl c)\/} of the theorem and
{\sl Lemma 5.3.1} provide the existence of parameterizations $\delta_n: \Sigma \to
C_n$ and annuli $A_i$, adjacent to each boundary circle $\gamma_i$, such that 
$\delta_n^* j_{C_n}$ are constant in every $A_i$. Thus we may assume that
such $\delta_n$ are given.
%%%%%%%%%%%%%%%%%%%%%%%Question%%%%%%%%%%%%%%%%%%%55

Take a covering $\calv= \{V_\alpha\}$ and parameterizations $\sigma_n$ as in
{\sl Theorem 4.3}. With every such covering we can associate the curves
$C_{\alpha,n} \deff \sigma_n(V_\alpha)$, the parameterizations $\sigma_{\alpha,
n} \deff \sigma_n\ogran_{V_\alpha}: V_\alpha \to C_{\alpha,n}$ and the maps
$u_{\alpha,n} \deff u_n \ogran_{C_{\alpha, n}} : C_{\alpha, n} \to X$. Consider
the following type of convergence of sequences $\{ (C_{\alpha, n}, u_{\alpha,
n}, \sigma_{\alpha, n} )\}$ with the fixed $\alpha$:

\sl
\item{A)} $C_{\alpha, n}$ are annuli of infinitely growing conformal radii
$l_n$ and the conclusions of {\sl Lemma 5.2.2} hold;

\item{B)} every $C_{\alpha, n}$ is isomorphic to the standard node $\cala_0=
\Delta \cup_{ \{0\} } \Delta$ such that the compositions $V_\alpha \buildrel
\sigma _{\alpha, n} \over \lrar C_{\alpha, n} \buildrel \cong \over \lrar
\cala_0$ define the same parameterizations of $\cala_0$ for all $n$;
furthermore, the induced maps $\ti u_{\alpha, n}: \cala_0 \to X$ strongly
converge;

\item{C)} the structures $\sigma_n^*j\vph_n \ogran_{V_\alpha}$ and the maps
$u_{\alpha, n}\scirc \sigma _{\alpha, n}: V_\alpha \to X$ strongly converge.

\rm\noindent
Here the strong convergence of maps is the one in the $L^{1,p}$-topology on
compact subsets for some $p>2$ (and hence for all $p<\infty$), and the
convergence of structures means the usual $C^\infty $-convergence.

\smallskip
Suppose that there is a subsequence, still indexed by $n\to \infty$ such
that for any $V_\alpha$ we have one of the convergence types {\sl A)--C)}.
Then the sequence of global maps $\{ (C_n , u_n, \sigma _n )\}$ converges in
the Gromov topology which gives us the proof, as well as a 
precise description of the convergence picture.

Otherwise, we want to find a refinement of our covering $\calv$ and
parameterizations $\sigma_n$ which have the needed properties. We shall proceed
by induction, estimating an area of pieces of coverings of $\Sigma$. 
To do this, we
fix $\eps>0$ satisfying $\eps \le {\eps_1^2 \over2 }$ with $\eps_1$ from {\sl
Lemma 5.1.1}, $\eps \le {\eps_2^2 \over2 }$ with $\eps_2$ from 
{\sl Lemma 5.1.2}
and $\eps \le {\eps_3\over 3}$ with $\eps_3$ from {\sl Lemma 5.2.3}. First
consider the

\smallskip\noindent{\sl
Special case: $\area(u_n(\sigma_n(V_\alpha))) \le \eps$ for any $n$ and any
$V_\alpha \in \calv$}. We can consider every $V_\alpha$ separately. If the
structures $\sigma_n^* j_n |_{V_\alpha}$ are constant, then some subsequence
of $u_n \scirc \sigma_n$ strongly converges due to {\sl Corollary 2.5.1}.

If structures $\sigma_n^* j_n |_{V_\alpha}$ are not constant, then $V_\alpha$
must be an annulus. Fix biholomorphisms $\phi_n: Z(0,l_n) \buildrel \cong
\over \lrar \sigma_n(V_\alpha)$. If $l_n \lrar \infty$, {\sl Lemma 5.2.2} 
shows
that (and describes how!) an appropriate subsequence of $u_n \scirc \phi_n$
converges to a $J_\infty$-holomorphic map of a standard node. Otherwise, we can
find a subsequence, still denoted $(C_n, u_n)$, for which $l_n \lrar l_\infty
<\infty$ and $u_n \scirc \phi_n$ converge to a $J_\infty$-holomorphic map of
$Z(0,l_\infty)$ in $L^{1,p}$-topology on compact subsets $K \Subset Z(0,l_\infty)$
for any $p<\infty$. To construct refined parameterizations $\ti\sigma _{\alpha,
n} : V_\alpha \to C_{\alpha, n}= \sigma_n(V_\alpha)$, we use property
{\sl(d)} from {\sl Theorem 5.3.2} and apply {\sl Lemma 5.3.1}.

Thus, we reach one of the convergence types {\sl A)--C)} which gives the proof
in this {\sl Special case}.

\smallskip\noindent
{\sl General case.}
Suppose that the theorem is proved for all sequences of $J_n$-holo\-mor\-phic
curves $\{(C_n, u_n)\}$ with parameterizations $\delta_n: \Sigma \to C_n$
which satisfy the additional condition $\area(u_n(C_n)) \le (N-1) \eps$ for
all $n$. We see this as the hypo\-thesis of the induction in $N$, so that our
{\sl Special case} is the base of the induction.

\smallskip
Assume that there exists a subsequence, still indexed by $n\to \infty$ such
that for every $V_\alpha$ and for the curves $C_{\alpha, n} = \sigma_n(
V_\alpha )$ the statement of the theorem holds. This means that there exist
refined coverings $V_\alpha= \cup_i V _{\alpha, i}$ and new parameterizations
$\ti\sigma _{\alpha, n}: V_\alpha \to C_{\alpha, n}$ such that the 
$\ti\sigma_n$
coincide with the $\sigma_n$ near the boundary of every $V_\alpha$ and 
such that
for the curves $C_{\alpha, i,n} \deff \ti\sigma_n (V _{\alpha ,i})$ we have 
convergence of one of the types {\sl A)--C)}. Then we can glue $\ti\sigma
_{\alpha, n}$ together to global parameterizations $\ti\sigma_n: \Sigma \to
C_n$ and set $\wt\calv \deff \{ V _{\alpha, i} \}$, obtaining the proof.

In particular, due to the inductive hypothesis, this is also true for 
any $V_\alpha$ such that 
$\area(u_n(\sigma_n(V_\alpha))) \le (N-1) \eps$ for all $n$.

This implies that it is sufficient to consider only those $V_\alpha$
for which $(N-1) \eps\le \area(u_n(\sigma_n( V_\alpha ))) \le N \eps$
for all $n$. Obviously, it is sufficient to show the desired property only
for such a piece of covering, say for $V_1$. To construct the refined
parameterizations $\ti\sigma_{1,n}$ and the covering $V_1 = \cup_i V _{1, i}$,
we consider four cases.

\smallskip
{\sl Case 1):  The structures $\sigma_n^*j_n|_{V_1}$ do not change and
$C_{1,n}$ are not isomorphic to the standard node $\cala_0$}.
Then we can realize $(V_1, \sigma_n^*j_n)$ as a constant bounded
domain $D$ in $\cc$. Hence we can consider $u_n \scirc \sigma_n: V_\alpha \to
X$ as holomorphic maps $u_n :D \to (X,J_n)$. Now we use the ``patching
construction'' of Sacks-Uhlenbeck [S-U].

Fix some $a>0$. Denote $D_{-a} \deff \{ z\in D\;:\; \Delta(z,a) \subset D\}$.
Find a covering of $D_{-a}$ by open sets $U_i \subset D$ with $\diam(U_i) <a$
such that any $z\in D$ lies in at most 3 pieces $U_i$. Then for any $n$ there
exists at most $3N$ pieces $U_i$ with $\area( u_n(U_i)) > \eps$. Taking a
subsequence, we may assume that the set of such ``bad'' pieces $U_i$ is the
same for all $n$. Repeat successively the same procedure for ${a\over2}$,
${a\over4}$, and so on, and then take a diagonal subsequence. We obtain at
most $3N$ ``bad'' points $y^*_1, \ldots, y^*_l$ such that a subsequence of
$u_n$ converges in $D\bs \{y^*_1, \ldots,y^*_l\}$ strongly, \i.e., in 
the $L^{1,p}$-topology on compact subsets 
$K \Subset D\bs \{y^*_1, \ldots,y^*_l\}$.
These ``bad'' points $y^*_1, \ldots, y^*_l$ are characterized by the property
$$
\text{\ \ for any $r>0$\ \ }
\area(u_n(\Delta(y^*_i, r)) >\eps \text{\ \ for $n$ all sufficiently big}.
\eqno(5.3.1)
$$

%%%%%%%%%%%%%%%%%%%%page 74%%%%%%%%%%%%%%%%%%%%%%%

\state Remark. As we shall now see, every such point is the place where the
``bubbling phenomenon'' occurs. Therefore we shall call $y^*_i$ {\sl bubbling
points}. The characterization property of a bubbling point is (4.1).

\smallskip
If there are no bubbling points, \i.e., $l=0$, then the chosen subsequence
$u_n$ converges strongly and we can finish the proof by induction.

Otherwise, we consider the first point $y^*_1 \in D$. Take a disk $\Delta
(y^*_1,\varrho)$ which doesn't contain any other bubbling points $y^*_i$,
$i>1$.

Then for any $n$ we can find the unique $r_n$ such that

(1) $r_n \le {\varrho \over 2}$ and $\area(u_n(\Delta(x, r_n)))\le\eps$ for
any $x\in \barr\Delta (y^*_1,{\varrho \over 2})$;

(2) $r_n$ is maximal \wrt (1).

\smallskip
Then $r_n\lrar 0$, because otherwise for $r^+ \deff {\sf lim\,sup}\; r_n >0$
and for some subsequence $n_k \lrar \infty$ with $r_{n_k} \lrar r^+$
we would have
$$
\area\bigl(u_{n_k}(\Delta(y^*_1, r^+))\bigr) \le \eps,
$$
which would contradict (5.3.1).

\state Lemma 5.3.4. {\it For every $n>\!>1$ there exists $x_n\in \barr 
\Delta(y^*_1, {\varrho \over2})$, such that $x_n\to y^*_1$ and $\area( u_n( 
\Delta(x_n, r_n))) =\eps$.}

\nobreak
\state Proof. If not, then for some subsequence $n_k\lrar \infty$ and every
$x\in \barr\Delta (y^*_1, {\varrho \over2})$ we would have $\area(u_{n_k}(
\Delta (x, r_{n_k})))< \eps $. Since $r_n \lrar0$, this would contradict
the maximality of $r_n$. In particular, there exists the postulated
sequence $\{ x_n \}$.

If $x_n$ do not converge to $y^*_1$, then after going to a subsequence we
would find $y'=\lim_{n\to \infty } x_n\not= y^*_1$. By our construction, $y'$
does not coincide with any other bubbling point $y^*_i$. Take $a>0$ such that
$\Delta(y', a)$ contains no bubbling point. Then by {\sl Corollary 2.5.1} some
subsequence $u_{n_k}$ would converge to some $u' \in L^{1,p}_\loc (\Delta(y',
a), X)$ in strong $L^{1,p}(K)$-topology for any compact subset $K \Subset
\Delta(y', a)$ and any $p<\infty$. In particular, for sufficiently small
$b<a$ we would have $\area(u_{n_k}(\Delta(y', b))) \lrar \area(u'( \Delta(y',
b))) < \eps$, which would contradict the choice of $r_n$ and $x_n$.
\qed

\smallskip
Using constructed $r_n$ and $x_n$, define maps $v_n:\Delta (0,{\varrho \over
2r_n}) \to (X,J_n)$ by $v_n(z)\deff u_n(x_n +r_n z)$. By the definition of
$r_n$ we have
$$
\area(v_n(\Delta (x ,1))\le \eps \ \ \text{for all $x\in \Delta
(0,{\textstyle{\varrho \over 2r_n}}-1)$}.
\eqno(5.3.2)
$$
On the other hand, $\area(v_n(\Delta(0,1)) = \area(u_n( \Delta (x_n, r_n)) =
\eps$ by {\sl Lemma 5.3.4}. Thus $v_n$ converge (after going to a subsequence)
on compact subsets in $\cc$ to a nonconstant $J_\infty $-holomorphic map
$v_\infty$ with finite energy. Consequently, $v_\infty$ extends onto $S^2$ by
the removable singularity theorem of {\sl Corollary 5.2.1}.

Since $v_\infty$ is nonconstant, $\norm{ dv_\infty }^2_ {L^2( S^2)}= \area(
v_\infty(S^2))\ge 3\eps$ by {\sl Lemma 5.2.3} and the choice of $\eps$.
Choose $b>0$ in such a way that
$$
\area(v_\infty (\Delta (0,b)) =
\norm{ dv_\infty}_{L^2(\Delta (0,b))}^2 \ge 2\eps.
\eqno(5.3.3)
$$
By {\sl Corollary 2.5.1} this provides that
$$
\norm{du_n}_{L^2(\Delta (x_n,br_n))}^2=
\norm{ dv_n}_{L^2(\Delta (0,b))}^2 \ge\eps.
\eqno(5.3.4)
$$

For $n>\!>$ we consider the coverings of $V_1$ by 3 sets
$$
\textstyle\mathsurround=0pt
V^{(n)}_{1,1} \deff V_1 \bs \barr\Delta (y^*_1, {\varrho\over2}),
\qquad
V^{(n)}_{1,2}\deff \Delta (y^*_1, \varrho) \bs \barr \Delta (x_n, br_n),
\qquad
V^{(n)}_{1,3} \deff \Delta (x_n, 2br_n).
$$
Fix $n_0$ sufficiently big. Denote $V_{1,1} \deff V^{(n_0)}_{1,1}$, $V_{1,2}
\deff V^{(n_0)}_{1,2}$, and $V_{1,3} \deff V^{(n_0)}_{1,3}$. There exist
diffeomorphisms $\psi_n: V_1 \to V_1$ such that $\psi_n: V_{1,1} \to V^{(n)}
_{1,1}$ is identity, $\psi_n: V_{1,2} \to V^{(n)}_{1,2}$ is a diffeomorphism
and $\psi_n: V_{1,3} \to V^{(n)}_{1,3}$ is biholomorphic \wrt the complex
structures, induced from $C_{1,n}$.

Thus, we have constructed the covering $\{ V_{1,1}, V_{1,2}, V_{1,3} \}$ of
$V_1$ and parameterizations $\sigma'_n \deff \sigma_{1,n} \scirc \psi_n: V_1
\to C_{1, n}$ such that the conditions of {\sl Theorem 5.3.2} are satisfied.
Moreover, $\area(u_n(\sigma'_n(V_{1,i}))) \le (N-1)\eps$ due to inequality
(5.3.4). Consequently, we can apply the inductive assumptions for the sequence
of curves $\sigma'_n(V_{1,i})$ and finish the proof by induction.

\smallskip\noindent
{\sl Case 2):} $V_1$ is a cylinder, structures $\sigma_n^*j_n|_{V_1}$ vary
with $n$, but conformal radii of $(V_1, \sigma_n^*j_n)$ are bounded uniformly
in $n$. Applying {\sl Lemma 5.3.3}, we can assume that structures $\sigma_n^*
j_n$ converge to a structure of an annulus with finite conformal radius. The
constructions of {\sl Case 1)} are used here with the following minor
modifications. First, we find the set of the bubble points $ y^*_i\in V_1$,
using the same patching construction and the characterization (5.3.1). Then
we find diffeomorphisms $\phi_n : V_1 \to V_1$ such that {\sl\.a)} $\phi_n$
converge to the identity map $\id: V_1 \to V_1$; {\sl\.b)} $\phi_n$ are
identical in fixed (\ie independent of $n$) annuli adjacent to the boundary
circles of $V_1$; {\sl\.c)} $\phi_n$ preserve every bubble point, $\phi_n(
y^*_i) = y^*_i$; and finally {\sl\.d)} for the ``corrected'' parameterizations
$\ti\sigma_n \deff \sigma_n \scirc \phi_n$ the structures $\ti\sigma_n^*j_n|
_{V_1}$ are constant in a neighborhood of every bubble point $y^*_i$. Then
we repeat the rest of the constructions of {\sl Case 1)} using the new
parameterizations $\ti\sigma_n$.

\smallskip\noindent
{\sl Case 3): Every $C_{1,n} = \sigma_n(V_1)$ is isomorphic to the standard
node $\cala_0$}. Fix identifications $C_{1,n} \cong \cala_0$ such that the
induced  parameterization maps $\sigma_{1, n} : V_1 \to \cala_0$ are the same
for all $n$. Represent $\cala_0$, and hence every $C_{1,n}$, as the union of
two discs $\Delta'$ and $\Delta''$ with identification of the centers $0\in
\Delta'$ and $0\in \Delta''$ into the nodal point of $\cala_0$, still denoted
by $0$. Let $u'_n :\Delta' \to X$ and $u''_n :\Delta'' \to X$ be the
corresponding ``components'' of the maps $u_{1,n} : C_{1,n} \to X$. Find the
common collection of bubbling points $y^*_i$ for both maps
$u'_n :\Delta' \to X$ and $u''_n :\Delta'' \to X$.
If there are no bubble points, then we obtain the convergence type
{\sl B)} and the proof can be finished by induction.
Otherwise consider the first such point $y^*_1$ which lies, say, on
$\Delta'$. If $y^*_1$ is distinct from the nodal point $0 \in
\Delta'$, then we simply repeat all the constructions from {\sl Case 1)}.

It remains to consider the case $y^*_1=0 \in \Delta'$. Now one should
modify the argumentations of {\sl Case 1} in the following way. Repeat
the construction of the radii $r_n\lrar 0$, the points $x_n \lrar y^*_1=0$,
and the maps $v_n: \Delta(0, {\varrho \over 2r_n}) \to X$, $v_\infty: S^2 \to
X$ from {\sl Case 1}. Set $R_n \deff |x_n|$, so that $R_n$ is the distance
from $x_n$ to point $0=y^*_1 \in \Delta'$. After rescaling $u_n$ to the maps
$v_n$, the point $0\in \Delta'$ will correspond to the point $z^*_n \deff
-{x_n \over r_n}$ in the definition domain $\Delta(0, {\varrho \over 2r_n})$
of the map $v_n$. We will now consider two subcases.

\smallskip\noindent
{\sl Subcase 3\/$'$): The sequence ${R_n\over r_n}$ is bounded}. This is
equivalent to boundedness of the sequence $ z^*_n$. Going to a
subsequence we may assume that the sequence $z^*_n$ converges to a point
$z^*\in \cc$. This point will be a nodal one for $(S^2, v_\infty)$.
As above, $v_\infty$ is nonconstant and $\norm{ dv _\infty} ^2 _{L^2( S^2)}=
\area( v_\infty( S^2) )\ge 3\eps$. Choose $b>0$ in such a way that
$$
\norm{ dv_\infty}_{L^2(\Delta (0,b))}^2 \ge 2\eps
\eqno(5.3.5)
$$
and $b \ge 2|z^*|+2$. Due to {\sl Corollary 2.5.1} for $n>\!>1$ we obtain
the estimate
$$
\norm{du'_n}_{L^2(\Delta' (x_n, br_n))}^2 =
\norm{ dv_n}_{L^2(\Delta (0,b))}^2 \ge \eps.
\eqno(5.3.6)
$$
Here $\Delta' (x, r)$ denotes the subdisc of $\Delta'$ with the
center $x$ and the radius $r$. Furthemore, for $n>\!>1$ we have the relation
$z^*_n \in \Delta( 0, b-1)$, or equivalently, $0 \in  \Delta' (x_n, (b-1)
r_n))$.

Define the coverings of $\cala_0$ by four sets
$$\mathsurround=0pt
\matrix\format\l\ \ &\l\\
W^{(n)}_1 \deff  \Delta'
\bs \barr\Delta' (0, {\textstyle{\varrho\over2}}),
&
W^{(n)}_2 \deff \Delta' (0, \varrho)
\bs \barr \Delta' (x_n, br_n),
\cr
\noalign{\vskip5pt}
W^{(n)}_3 \deff \Delta' (x_n, 2br_n)
\bs \barr\Delta'(0,{\textstyle {r_n \over 2}}),
&
W^{(n)}_4 \deff \Delta' (0, r_n) \cup
\Delta'',
\endmatrix
$$
and lift them to $V_1$ by putting $V^{(n)}_{1,i} \deff \sigma\inv_{1,n}(
W^{(n)}_i)$. Choose $n_0 >\!> 0$ such that $z^*_{n_0} \in \Delta( 0, b-1)$
and the relation (5.3.6) holds. Set $V_{1,i} \deff V^{(n_0)} _{1, i}$. Choose
diffeomorphisms $\psi_n: V_1 \to V_1$ such that $\psi_n: V_{1,1} \to V^{(n)}
_{1,1}$ is the identity map, $\psi_n: V_{1,2} \to V^{(n)}_{1,2}$ and $\psi_n:
V_{1,3} \to V^{(n)}_{1,3}$ are diffeomorphisms, and $\psi_n: V_{1,4} \to
V^{(n)}_{1,4}$ corresponds to isomorphisms of nodes $W^{(n)}_4 \cong
\cala_0$. Set $\sigma'_n \deff \sigma_n \scirc \psi_n$. The choice above can
be done in such a way that the refined covering $\{ V_{1,i} \}$ of $V_1$ and
parameterization maps $\sigma'_n: V_1 \to C_{1,n}$ have the properties of {\sl
Theorem 5.3.2}. Moreover, relations (5.3.2) and (5.3.6) imply the estimate $\area(
u_n (\sigma'_n(V_{1,i})) \le (N-1)\,\eps$. This provides the inductive
conclusion for {\sl Subcase 3\/$'$)}.

\smallskip\noindent
{\sl Subcase 3\/$''$): The sequence $R_n \over r_n$ increases infinitely}.
This means that the sequence $ z^*_n$ is not bounded. Nevertheless $R_n \lrar
0$ since $x_n \lrar 0$. We proceed as follows. Repeat the construction of the
raduis $b$ from {\sl Case 1)}. For $n>\!>0$ define the coverings of $\cala_0$
by six sets
$$\mathsurround=0pt
\matrix\format\l\ \ &\l\\
W^{(n)}_1 \deff  \Delta'
\bs \barr\Delta' (0, {\varrho\over2}),
&
W^{(n)}_2 \deff \Delta' (0, \varrho)
\bs \barr \Delta' (x_n, 2R_n),
\cr\noalign{\vskip5pt}
W^{(n)}_3 \deff  \Delta' (x_n, 4R_n)
\bs \bigl(\barr \Delta' (x_n, {R_n \over 6})
\cup \barr \Delta' (0,{R_n\over 6}) \bigr)
&
W^{(n)}_4 \deff \Delta' (0, {R_n\over 3}) \cup \Delta'',
\cr\noalign{\vskip5pt}
W^{(n)}_5 \deff \Delta' (x_n, {R_n \over 3})
\bs \barr\Delta'(x_n, br_n),
&
W^{(n)}_6 \deff \Delta' (0, 2br_n ),
\endmatrix
$$
and lift them to $V_1$ by putting $V^{(n)}_{1,i} \deff \sigma\inv_{1,n}(
W^{(n)}_i)$. Choose $n_0 >\!> 0$ such that $R_{n_0} >\!> br_{n_0}$ , and set
$V_{1,i} \deff V^{(n_0)} _{1, i}$. Choose diffeomorphisms $\psi_n: V_1 \to
V_1$ such that $\psi_n: V_{1,1} \to V^{(n)} _{1,1}$ is the identity map,
$\psi_n: V_{1,2} \to V^{(n)}_{1,2}$, $\psi_n: V_{1,4} \to V^{(n)}_{1,4}$ and
$\psi_n: V_{1,5} \to V^{(n)}_{1,5}$ are diffeomorphisms, and finally,
$\psi_n: V_{1,6} \to V^{(n)}_{1,6}$ corresponds to isomorphisms of nodes
$W^{(n)}_6 \cong \cala_0$. Set $\sigma'_n \deff \sigma_n \scirc \psi_n$.
Again, this choice can be done in such a way that $\{ V_{1,i} \}$ and
parameterization maps $\sigma'_n: V_1 \to C_{1,n}$ have the properties of
{\sl Theorem 5.3.2}. As above, we get the estimate $\area(u_n (\sigma'_n(
V_{1,i} )) \le (N-1)\,\eps$ due to (5.3.2). Thus we get the inductive
conclusion for {\sl Subcase 3\/$''$)} and can proceed further.

\smallskip\noindent
{\sl Case 4):} $V_1$ is a cylinder, structures $\sigma_n^*j_n|_{V_1}$ vary
with $n$, and conformal radii of $(V_1, \sigma_n^*j_n)$ converge to
$+\infty$. Using {\sl Lemma 5.3.3}, we can assume that structures $\sigma_n^*
j_n$ satisfy property \slii of this lemma.

Fix biholomorphisms $\sigma_n(V_1) \cong Z(0, l_n)$. If $\area(u_n(Z(a-1,
a))) \le \eps$ for any $n$ and any $a \in [1, l_n]$, then {\sl Lemma 5.2.2}
shows that $u_n: \sigma_n(V_1) \to X$ converge to a $J_\infty$-holomorphic map
from a node.

If not, then, after passing to a subsequence, we can find a sequence $\{a_n
\}$ with $a_n \in [1, l_n]$ such that $\area(u_n(Z(a_n \allowbreak -1, a_n)))
\ge \eps$. If $a_n$ is bounded, say $a_n \le a^+$, then we cover $Z(0,
l_n)$ by the sets $V_{1,1} \deff Z(0, a^+ +2)$ and  $V_{1,2} \deff Z(a^+ +1,
l_n)$. If $l_n - a_n$ is bounded, say $l_n - a_n \le a^+$, then we cover $Z(0,
l_n)$ by the sets $V_{1,1} \deff Z(0, l_n - a^+ +2)$ and  $V_{1,2} \deff
Z(l_n -a^+ +1, l_n)$. In the remaning case, when both $a_n$ and $l_n - a_n$
increase infinitely, we cover $Z(0, l_n)$ by 3 sets $V_{1,1}^{(n)} \deff
Z(0,a_n-1)$, $V_{1,2}^{(n)} \deff Z(a_n-2, a_n+1)$, and $V_{1,3}^{(n)} \deff
Z(a_n, l_n)$. Cover $V_1$ by 2 or, respectively, 3 successive cylinders $V_{1,
i}$ in an obvious way. Find diffeomorphisms $\psi_n : V_1 \to V_1$ identical
in the neighborhood of the boundary of $V_1$ and such that $\psi_n(V_{1,i}) =
\sigma_n ^{-1} V^{(n)} _{1,i}$. Define the new parameterizations $\sigma'_n
\deff \sigma_n \scirc \psi_n$. Note that we may additionally assume that if
the conformal radius of $\sigma'_n(V_{1,i})$ is independent of $n$ then the
structure $\sigma'_n{}^*j_n \ogran_{V_{1,i}}$ is also independent of $n$.

By this construction we obtain the following property of the 
covering $\{ V_{1,i}\}$ and new parameterizations 
$\sigma'_n$. For any $V_i$ we have either
the estimate
$$
\area(u_n(\sigma'_n(V_{1,i}))) \le (N-1)\,\eps.
$$
or the structures $\sigma'_n{}^*j_n|_{V_{1,i}}$ do not depend on $n$. Thus we
reduce our case to the situation which is covered either by the inductive
assumption or by {\sl Case 1)}.

\medskip
The proof of the theorem can be finished by induction. The fact that the
limit curve $(C_\infty, u_\infty)$ remains stable over $X$ is proved in {\sl
Lemma 5.3.5} below. \qed

\medskip
\state Remark. Here we explain the meaning of the constructions used in the
proof. We start with {\sl Case 1)}, where $J_n$-holomorphic maps from a fixed
domain $D\subset \cc$ are treated. Bubbling points $y^*_i$ appear in this case
as those where the strong convergence of maps $u_n: D \to (X, J_n)$
fails. The patching construction of Sacks and Uhlenbeck insures that the
``convergence failure'' set is finite and insures an effective estimate on
the number of bubbling points by the upper bound of the area, $l \le 3N$ in
our situation. The characterization property (5.3.2) of bubbling points is
essentially due to Sacks and Uhlenbeck; the only difference is that we use
the area of the map $u$ (which is equivalent to the energy of $u$, \i.e.,
$L^2$-norm of $du$) , whereas in [S-U] the $\norm{du}_{L^\infty}$ is used. The
next step, the construction of maps $v_n$ as the rescaling of the $u_n$ and the
existence of the limit $v_\infty$ is also due to Sacks and Uhlenbeck.

Due to the explicit construction of the map $v_\infty$, it is useful to
imagine the curve $(S^2, v_\infty)$ as a ``bubbled sphere'' and $y^*_1$ as
the point where the ``bubbling'' occurs. Moreover, one obtains natural
partitions (one for each $n>\!>0$) of $D$ into three pieces: $D$ minus a fixed
small neighborhood of $y^*_1$; disks $(\Delta(x_n, br_n), u_n)$ representing
pieces $(\Delta(0, b), v_\infty)$ and approximating a sufficiently big part
$(\Delta(0, b), v_\infty)$ of the bubbled sphere and the ``part between''.

These latter ``parts between'' appear to be the annuli of infinitely growing
conformal radii, the situation considered in {\sl Case 4)}. Since neither
outer nor inner boundary circle should be preferred in some way, we consider
them as long cylinders $C_n = Z(0,l_n)$ with $l_n \lrar \infty$, in the 
spirit of
{\sl Definition 3.2}. {\sl Lemma 5.2.2} provides a ``good'' convergence model
for long cylinders, stated above as convergence type {\sl A)}. If such
a convergence for a sequence $(C_n, u_n)$ fails, then there must exist
subannuli $A_n \subset C_n$ of a constant conformal radius for which
$\area(u_n(A_n)) \ge \eps$.

In both cases --- a constant domain $D$ or long cylinders --- we proceed by
cutting the curves into smaller pieces. We arrive at a situation that is 
simpler
in the following way. The obtained curves either converge or have the
upper bound for the area smaller in the fixed constant $\eps$. Thus, the use
of induction leads finally to a decomposition of the curves into pieces for
which one of the convergence types {\sl A)--C)} holds. The possibility of
gluing these final pieces together is insured by the fact that the partitions
above are represented by appropriate coverings satisfying the conditions of
{\sl Theorem 5.3.2}.

Considering curves with nodes, additional attention should be paid to the
case where bubbling appears at a nodal point. This situation is considered in
{\sl Case 3)}. The constructed points $x_n$ and radii $r_n$ describe the
``center'' and the ``size'' of energy localization of the bubbling,
represented by the sequence of the maps $v_n \lrar v_\infty$. So the
convergence picture depends on whether the energy localization occurs
near the nodal point ({\sl Subcase 3\/$'$)}) or away from it ({\sl Subcase
3\/$''$)}). As a result, the nodal point can either remain on the ``bubbled''
sphere $(S^2, v_\infty)$ or move into the ``part between'', which is
represented by long cylinders.

In {\sl Subcase 3\/$'$)} we remove neighborhoods of the nodal point from the
disks $(\Delta(0, b), v_n)$ and thus get four pieces of converging instead 
of three as
in {\sl Case 1)}. In {\sl Subcase 3\/$''$)} the situation is more complicated,
because we must take into account the position of the nodal point in the long
cylinders --- the ``parts between''. Thus, we must consider now the sequence
of cylinders with one marked point, \i.e., the sequence of pants. 
The vanishing
$R_n\to 0$ and $r_n / R_n \to$ mean that conformal structure of those pants
is not constant and converges to one of the spheres with three punctures.

In order to have the covering pieces with the convergence types {\sl A)--C)},
we choose an appropriate refinement of the covering. After that, we obtain two
sequences of long cylinders, describing the appearance of two new nodal points.
The first one corresponds to the part between the ``original'' nodal point
and bubbled sphere and is represented by $V_5$, whereas the other one,
represented by $V_2$, lies on the other side of the ``original'' nodal
point. In addition, we fix a neighborhood of the ``original'' nodal point which
has a constant complex structure and is topologically an annulus with the
disc $\Delta''$ attached to the nodal point. To satisfy the requirements 
of {\sl Theorem 5.3.2}, we cover the neighborhood by two pieces, 
the pants $V_3$ and the
piece $V_4$ parameterizing the nodes $W^{(n)}_4$. This explains the 
appearance of
six pieces of covering in {\sl Subcase 3\/$''$)}.

\smallskip
\state Lemma 5.3.5. \it The limit curve $(C_\infty, u_\infty)$ constructed
in the proof of {\sl Theorem 1.1} is stable over $X$. \rm

\state Proof. If $(C_\infty, u_\infty)$ is unstable over $X$, then either
$C_\infty$ is a torus with $u_\infty$ constant, or $C_\infty$ should
have a component $C' \subset C_\infty$ such that $\ti C'$ is a sphere with
at most two marked points, and $u_\infty(C')$ is a point.

The case of a constant map from a torus is easy to exclude. In fact, in this
case all $C_n$ must also be tori with $\area(u_n(C_n))$ sufficiently small
for $n>\!>1$. Cover every $C_n$ by an infinite cylinder $Z(-\infty, +\infty)$
and consider compositions $\ti u_n: Z(-\infty, +\infty) \to X$ of $u_n$ with
the covering maps. Since $\area(u_n(C_n))\approx 0$, {\sl Corollary 5.2.1} can
be applied to show that every $\ti u_n$ extends to a $J_n$-holomorphic
map from $S^2$ to $X$. Consequently, $\area(\ti u_n(Z(-\infty, +\infty)))$
must be finite. On the other hand, $\area(u_n(C_n))>0$ due to the stability
condition; hence $\area(\ti u_n(Z(-\infty, +\infty)))$ must be infinite.
This contradiction excludes the case of a torus.

The same argumentations are valid in the case, where $C_\infty$ is the sphere
with no marked points. Then the curves $C_n$ are also parameterized by
the sphere $S^2$. The condition of instability means that
$\area(u_\infty(C_\infty))=0$. Due to {\sl Corollary 2.5.1}, $\area(u_n(C_n))$
must be sufficiently small for $n>\!>1$. Now {\sl Lemma 5.2.3} and  the
stability of $(C_n,u_n)$ exclude this possibility.

Now consider the cases where the limit curve $C_\infty$ has a ``bubbled''
component $C'$, which  is the sphere with one or two marked points. 
If $C'$ has one
marked point, then $C'$ must appear as a ``bubbled'' sphere $(S^2, v_\infty)$
in the constructions of {\sl Cases 1)--3)} in the proof of {\sl Theorem 1.1}.
But these constructions yield only non-trivial ``bubbled'' spheres, for which
$v_\infty \not= const$. Thus, such a component $C'$ must be stable.

\smallskip
In the remaining case, a component $C'$ with two marked points, we consider a
domain $U\subset C_\infty$, which is the union of the component $C'$ and
neighborhoods of the marked point on $C'$. If $C'$ is an unstable component,
then $\area(u_\infty(C'))=0$, and we can achieve the estimate $\area( u_\infty(
U))< \eps$ taking $U$ sufficiently small. Set $\Omega \deff \sigma_\infty
\inv(U)$, where $\sigma_\infty: \Sigma \to C_\infty$ is the parameterization of
$C_\infty$. Let $\gamma_1$ and $\gamma_2$ be the pre-images of marked points
on $C'$. Then $\Omega$ must be a topological annulus, and $\gamma_i$, $i=1,2$,
disjoint circles generating the group $\pi_1 (\Omega) = \zz$. Further, $C'$
must be a ``bubbling'' component of $C_\infty$, \i.e., at least for one of the
circles $\gamma_i$, $i=1,2$, the images $\sigma_n(\gamma_i)$ are not nodal
points of $C_n$ but smooth circles.

If these are both circles $\gamma_1$ and $\gamma_2$, then $U_n \deff 
\sigma_n (\Omega)$ would satisfy the conditions of {\sl Lemma 5.2.2}. In this
case we should have the convergence type {\sl A)}, and hence the limit piece
$\sigma_\infty(\Omega)$ should be isomorphic to the node $\cala_0$.

In the case where only one circle, say $\gamma_1$, corresponds to nodal points
on $C_n$, and for the other one the images $\sigma_n(\gamma_2)$ are smooth
circles, then the domains $\sigma_n(\Omega)$ must be isomorphic to the node
$\cala_0$. Furthermore, due to the condition $\area(u_\infty(\sigma_\infty(
\Omega))) < \eps$, we would have $\area(u_n(\sigma_n (\Omega))) < \eps$.
Consequently, we would have the convergence type {\sl B)}, and the
unstable component $C'$ could not appear. \qed

%%%%%%%%%%%%%%
\newpage 
\noindent
{\bigbf Appendix 3}

\smallskip\noindent
{\bigbf Compactness with Totally Real Boundary Conditions.}

\medskip\noindent
{\bigsl A3.1. Curves with Boundary on Totally Real Submanifolds.}

\smallskip
In this section we consider the behavior of  complex  curves over an
almost complex manifold $(X,J)$ with a boundary on totally real submanifold(s).
As in the ``interior" case, we need to allow some type of boundary
singularity.

\state Definition A3.1.1. {\sl The set $\cala^+ \deff \{ (z_1,z_2) \in 
\Delta^2\;:\; z_1 \cdot z_2=0, \im z_1\ge 0, \im z_2\ge 0\}$ is called the 
{\it standard boundary node}. A curve $\barr C$ with boundary $\d C$ is called
a {\it nodal curve with boundary} if

\sli $C$ is a nodal curve, possibly not connected;

\slii $\barr C=C\cup \d C$ is connected and compact;

\sliii every boundary point $a\in \d C$ has a neighborhood homeomorphic
either to the half-disk $\Delta^+ \deff \{ z \in \Delta
\;:\; \im z\ge0 \}$, or to the standard boundary node $\cala^+$.

\smallskip\rm
In the last case $a\in \d C$ is called a {\it boundary nodal point}, whereas
nodal points of $C$ are called {\it interior nodal points}. Both boundary and
interior nodal points are simply called nodal points.
}

\state Definition A3.1.2. 
{\sl Let $(X,J)$ be an almost complex manifold. A pair
$(\barr C, u)$ is called a {\it curve with boundary over $(X,J)$} if $\barr
C=C\cup \d C$ is a nodal curve with boundary, and $u: \barr C \to (X,J)$ is a
continuous $L^{1,2}$-smooth map, which is holomorphic on $C$.
}

\smallskip
A curve $(C,u)$ with boundary is stable if the same condition as in {\sl
Definition 1.5} on the automorphism groups of compact irreducible components
is satisfied.

\smallskip
\state Remark. One can see that $\barr C$ has a uniquely defined real
analytic structure such that the normalization $\barr C{}^{\sf nr}$ is a
real analytic manifold with a boundary. More
precisely, the pre-image of every boundary nodal point $a_i$ consists of two
points $a'_i$ and $a''_i$. The normalization
map $s: \barr C{}^{\sf nr} \to \barr C$ glues pairs $(a'_i, a''_i)$ of points
on $\barr C{}^{\sf nr}$ into nodal points $a_i= s(a'_i)= s(a''_i)$ on $\barr
C$. This implies that the notion of a $L^{1,p}$-smooth map, $p>2$, as well as
a continuous $L^{1,p}$-smooth map $u: \barr C \to X$ is well- defined.

\state Definition A3.1.3. 
{\sl We say that a real oriented surface with boundary
$(\Sigma, \d\Sigma)$ {\it parameterizes} a nodal curve with boundary $C$ if
there is a continuous map $\sigma :\barr\Sigma \to \barr C$ such that

\sli if $a\in C$ is an interior nodal point, then $\gamma_a \deff \sigma\inv
(a)$ is a smooth imbedded circle in $\Sigma$;

\slii if $a\in\d C$ is a boundary nodal point, then $\gamma_a\deff \sigma
\inv (a)$ is a smooth imbedded arc in $\Sigma$ with end points on $\d\Sigma$,
transversal to $\d\Sigma$ at this point;

\sliii if $a,b \in \barr C$ are distinct (interior or boundary) nodal points,
then $\gamma_a\cap \gamma_b= \emptyset$;

\sliv $\sigma :\barr\Sigma \bs \bigcup_{i=1}^N\gamma_{a_i}\to \barr C \bs \{
a_1,\ldots ,a_N\} $ is a diffeomorphism, where $a_1,\ldots ,a_N$ are all
(interior and boundary) nodal points of $\barr C$.
}

\smallskip\rm
Recall that a real subspace $W$ of a complex vector space is called {\it
totally real} if $W \cap \isl W =0$. Similarly, a $C^1$-immersion $f:W \to X$
is called {\it totally real} if for any $w\in W$ the image $d f(T_wW)$ is
a totally real subspace of $T_{f(w)}X$.

Let $(\barr C, u)$ be a stable curve with boundary over an almost complex
manifold $(X,J)$ of a complex dimension $n$.

\state Definition A3.1.4. 
{\sl We say that $(\barr C, u)$ satisfies the {\it totally
real boundary condition $\mib W$ of type $\bfbeta$} if

\sli $\bfbeta= \{\beta_i\}$ is a collection of arcs with disjoint interiors,
which defines a decomposition of the boundary $\d C= \cup_i \beta_i$;
moreover, we assume that every boundary nodal point is an endpoint for four 
arcs $\beta_i$;

\slii ${\mib W} = \{(W_i, f_i)\}$ is a collection of totally real immersions
$f_i: W_i \to X$, one for every $\beta_i$;

\sliii there are given continuous maps $u^{(b)}_i :\beta_i \to W_i$ {\it
realizing conditions $\mib W$}, \i.e., $f_i \scirc u^{(b)}_i = u|_{\beta_i}$.
}

\state Remarks.~1. We shall consider (immersed) totally real submanifolds
only of {\sl maximal real dimension} $n= \dimc X$.

\state 2. If $\bfbeta$ is a collection of arcs as above, a parameterization
$\sigma: \barr \Sigma \to \barr C$ induces a collection of arcs $\sigma\inv(
\bfbeta) \deff \{ \sigma\inv(\beta_i)\;:\; \beta_i \in \bfbeta \} $ with the
properties similar to \sli of {\sl Definition A3.1.4}. Thus, $\sigma\inv(
\beta_i)$ have disjoint interiors, $\cup_i\sigma\inv(\beta_i) = \d\Sigma$,
and for any boundary node $a\in \barr C$ every endpoint of the arc $\beta_a=
\sigma\inv(a)$ is an endpoint of two arcs $\sigma\inv(\beta_i)$. Since
$\bfbeta$ is completely determined by $\sigma\inv(\bfbeta)$, we shall
denote both collections simply by $\bfbeta$ and shall not distinguish them
when considering boundary conditions.

\bigskip\noindent
{\bigsl A3.2. A priori Estimates near a Totally Real Boundary.}

\smallskip
A totally real boundary condition is a suitable elliptic boundary condition for
an elliptic differential operator $\dbar$ of the Cauchy-Riemann type. In
particular, all statements about ``inner" regularity and convergence for
 complex  curves remain valid near ``totally real" boundary. As in
the ``inner'' case, to get some ``uniform'' estimate at boundary 
one needs $W$ to be ``uniformly totally real''.

\state Definition A3.2.1. {\sl Let $X$ be a manifold with a Riemannian metric
$h$, $J$ a continuous almost complex structure, $W$ a manifold, and $A_W
\subset W$ a subset. We say that an immersion $f: W \to X$ is {\it
$h$-uniformly totally real along $A_W$ with a lower angle $\alpha = \alpha(W,
A_W, f) > 0$}, \iff

\sli $df: TW \to TX$ is $h$-uniformly continuous along $A_W$;

\slii for any  $w \in A_W$ and any $\xi \not =0 \in T_wW$ the angle
$\angle_h\bigl(Jdf(\xi), df(T_wW) \bigr) \ge \alpha$.
}

\medskip
We start with an analog for the First A priori Estimate. Define the
half-disks $\Delta^+(r) \deff \{ z\in \Delta(r) \;:\; \im z \ge 0 \}$ with
$\Delta^+ = \Delta^+(1)$ and $\Delta^- \deff \{ z \in \Delta \; :\; \im z \le
0 \}$. Set $\beta_0 \deff (-1, 1) \subset \d\Delta^+$. Let $X$ be a manifold
with a Riemannian metric $h$, $A\subset X$ a subset, $J^*$ a continuous
almost complex structure, $f: W \to X$ a totally real immersion and $A_W
\subset W$ a subset.

\state Lemma A3.2.1. 
{\it Suppose that $J^*$ is $h$-uniformly continuous on $A$
with the uniform continuity modulus $\mu_{J^*}$ and that $f : W \to X$ is
$h$-uniformly totally real along $A_W$ with a lower angle $\alpha_f>0$ and
the uniform continuity modulus $\mu_f$. Then for every $2< p<\infty $ there
exists an $\eps^b_1 =\eps^b_1(\mu_{J^*}, \alpha_f, \mu_f )$ (independent of
$p$) and $C_p= C(p, \mu_{J^*}, \alpha_f, \mu_f )$ such that the following
holds.

If $J$ is a continuous almost complex structure on $X$ with $\norm{ J- J^*}
_{L^\infty(A)} <\eps^b_1$, and if $u\in C^0 \cap L^{1,2}(\barr\Delta^+ ,X)$ is
a $J$-holomorphic map with $u(\Delta) \subset A$ and with the boundary
condition $u|_{\beta_0} =f \scirc u^b$ for some continuous $u^b: \beta_0
\to A_W \subset W$, then the condition $\norm{du}_{ L^2( \Delta^+ )}
<\eps^b_1$ implies the estimate
$$
\norm{du}_{L^p(\Delta^+({1\over 2}) )}\le C_p\cdot
\norm{du}_{L^2(\Delta^+ )}.\eqno(A3.2.1)
$$
}

\state Proof. Suppose additionally that $\diam(u(\Delta^+))$ is sufficiently
small. Then we may assume that $u(\Delta^+)$ is contained in some chart $U
\subset \cc^n$ such that $\norm{ J- J\st} _{L^\infty(U)}$ is also small
enough. Let $z=(z_1, \ldots, z_n)$ be $J\st$-holomorphic coordinates in $U$
such that $u(0)= \{z_i= 0\}$. Making an appropriate diffeomorphism of $U$, we
may additionally assume that $W_0 \deff f(W) \cap U$ lies in $\rr^n$ and that
$J = J\st$ along $W_0$.

Consider a trivial bundle $E \deff \Delta \times \cc^n$ over $\Delta$ with
complex structures $J\st$ and $J_u \deff J \scirc u$. We can consider $u$ as
a section of $E$ over $\Delta^+$ satisfying the equation 
$\dbar_{J_u} u\deff \d_x u+ J_u \d_y u=0$. Over $\beta_0$ we obtain 
a $J_u$-totally real subbundle $F \deff \beta_0 \times \rr^n$ 
such that $u(\beta_0) \subset F$. Let $\tau$
denote a complex conjugation in $\Delta$ as well as a 
complex conjugation in $E$
with respect to $J\st$. Extend $J_u$ on $E|_{\Delta^-}$ as the composition
$- \tau\scirc J_u \scirc \tau$. This means that for $z\in Delta^-$ we obtain
$$
J_u(z): v \mapsto  \tau v \in E_{\tau z}
\mapsto J_u(\tau z)(\tau v) \in E_{\tau z}
\mapsto  -\tau J_u(\tau z)(\tau v) \in E_z .
\eqno(A3.2.2)
$$

Since $J_u=J \scirc u$ coincides with $J\st$ along $\beta_0$, this extension
is also continuous. Further, for $v\in L^1( \Delta^+, E)$ define the
extension $\ex(v)$ by setting $v(z)\deff \tau v(\tau z)$. This gives a
continuous linear operator $\ex: L^p( \Delta^+, E) \to L^p( \Delta^+, E)$
for any $1\le p\le \infty$. An important property of operator $\ex$ is that
if $v \in L^{1,p} (\Delta^+, E)$ with $1\le p\le \infty$ (resp.\ $v\in C^0(
\barr \Delta^+ )$) satisfies the boundary condition $v|_{\beta_0} \subset
F$, then $\ex v \in L^{1,p}(\Delta, \cc^n)$ (resp.\ $\ex v \in C^0 (\barr
\Delta)$). Let us denote by $L^{1,p} (\Delta^+, E,F)$ (resp.\ by $C^0
(\Delta^+, E,F) $) the corresponding spaces of $v$ with the boundary
condition $v|_{\beta_0} \subset F$.

Since for $v\in L^{1,1}(\Delta^+, E)$ holds $\d_x(\tau v) =\tau (\d_x v)$ and
$\d_y(\tau v) = - \tau (\d_y v)$, we obtain 
$\dbar_{J_u} (\ex v) =\ex (\dbar_{J_u} v)$ 
for any $v \in L^{1,p} (\Delta^+, E,F)$. In particular, for 
$\ti u \deff \ex u \in C^0 \cap L^{1,2} (\Delta, E)$ 
we have $\dbar_{J_u} \ti u=0$.

From this point we can repeat the steps of the proof of {\sl Lemma 2.4.1}.
\qed

\state Remark. We shall refer to the construction of a complex structure
$J_u$ in $E$ over $\Delta^-$ and (resp.\ of a section $\ti u$ of $E$ over
$\Delta^-$) as an {\sl extension of $J \scirc u$ (resp.\ of $u$) by the
reflection principle}.

\smallskip\rm
Let $X$ be a manifold with a Riemannian metric $h$, $J^*$ a continuous almost
complex structure on $X$, $A\subset X$ a closed $h$-complete subset such
that $J^*$ is $h$-uniformly continuous on $A$, and $f_0 :W \to X$ an
immersion, which is $h$-uniformly totally real on some closed $f^*h$-complete
subset $A_W\subset W$.

\state Corollary A3.2.2. 
{\it Let $\{ J_n \}$ be a sequence of continuous almost
complex structures on $X$ such that $J_n$ converge $h$-uniformly on $A$ to
$J$, $f_n: W \to X$ a sequence of totally real immersion such that $df_n$
converge $h$-uniformly on $A_W$ to $df_0$, and $u_n\in C^0\cap L^{1,2}(
\Delta^+ ,X)$ a sequence of $J_n$-holomorphic maps such that $u_n(\Delta^+)
\subset A$, $\norm{du_n}_{L^2 (\Delta^+ )} \le \eps^b_1$, $u_n(0)$ is bounded
in $X$, and $u_n| _{\beta_0} = f_n \scirc u^b_n$ for some continuous $u^b_n
:\beta_0 \to A_W$.

Then there exists a subsequence $u_{n_k}$ where $L^{1,p}_\loc (\Delta^+
)$-converge to a $J$-holomorphic map $u_\infty$ for all $p < \infty$.
}

Here $\beta_0 = (-1,1) \subset \d \Delta^+$ and $L^{1,p}_\loc (\Delta^+
)$-convergence means $L^{1,p} (\Delta^+(r) )$-con\-ver\-gence for all $r<1$,
\i.e., convergence up to boundary component $\beta_0$.

\state Proof. As a statement itself, the proof of {\sl Corollary A3.2.2} 
copies the one of {\sl Corollary 2.5.1} with appropriate modifications and 
using the reflection principle.\qed

\bigskip\noindent
{\bigsl A3.3. Long Strips and the Second A priori Estimate.}

\smallskip
Consider now an analog of the Second A priori Estimate. An analog of ``long 
cylinders'' is now ``long strips" satisfying appropriate boundary conditions.

\state Definition A3.3.1. {\sl Define a {\it strip $\Theta(a,b) \deff (a,b)
\times [0,1]$} with the complex coordinate $\zeta \deff t - i\theta$, 
$t \in (a,b)$, $ \theta \in [0,1]$. Define also $\Theta_n \deff \Theta(n-1,n)$,
$\d_0 \Theta(a,b) \deff (a,b) \times \{0\}$, and $\d_1 \Theta(a,b) \deff (a,b)
\times \{1\}$.
}

We are interested in maps $u: \Theta(a,b) \to X$, which are holomorphic with
respect to the complex coordinate $\zeta$ on $\Theta(a,b)$ and a continuous
almost complex structure $J$ on $X$, which satisfy the boundary conditions
$$
u\ogran_{\d_0\Theta(a,b)} = f_0 \scirc u^b_0,
\hskip6em
u\ogran_{\d_1\Theta(a,b)} = f_1 \scirc u^b_1,
$$
with some $J$-totally real immersions $f_{0,1} :W_{0,1} \to X$ and continuous
maps $u_{0,1} :\d_{0,1}\Theta(a,b) \to W_{0,1}$. First we consider the special
linear case.

\state Lemma A3.3.1. {\it Let $W_0$ and $W_1$ be $n$-dimensional totally real
subspaces in $\cc^n=(\rr^n,J\st )$. Then there exist a constant $\gamma_W=
\gamma(n, W_0, W_1)$ with $0< \gamma_W <1$ such that for any holomorphic map
$u: \Theta(0,3) \to \cc^n$ with the boundary conditions
$$
u(\d_0 \Theta(0,3)) \subset W_0 \qquad
u(\d_1 \Theta(0,3)) \subset W_1
\eqno(A3.3.1)
$$
we have the following estimate:
$$
\int_{\Theta_2} |du|^2 dt\,d\theta \le \msmall{\gamma_W\over2}
\left( \int_{\Theta_1} |du|^2 dt\,d\theta +
\int_{\Theta_3} |du|^2 dt\,d\theta \right).
\eqno(A3.3.2)
$$
}

\state Proof. Let $L^{1,2}_W([0,1], \cc^n)$ be a Banach manifold of those
$v(\theta) \in L^{1,2}([0,1], \cc^n)$, $v(0) \in W_0$ and $v(1) \in W_1$.
Consider a nonnegative quadratic form $Q(v) \deff \int_0^1 |\d_\theta v(
\theta)|^2 d\theta$. Since $Q(v) + \norm{v} ^2_{L^2} =\norm{v}^2 _{L^{1,2}}$
and the imbedding $L^{1,2}_W([0,1], \cc^n) \hook L^2([0,1],\cc^n)$ is compact,
we can decompose $L^{1,2}_W([0,1], \cc^n)$ into a direct Hilbert sum of
eigenspaces $\ee_\lambda$ of $Q$ \wrt $\norm{v}^2_{L^2}$. This means that
$v_\lambda$ belongs to $\ee_\lambda$ iff for any $w\in L^{1,2}_W([0,1],
\cc^n)$ 
$$
\int_0^1 \<
\d_\theta v_\lambda(\theta), \d_\theta w(\theta) \> d\theta = \int_0^1
\lambda\< v_\lambda(\theta), w(\theta) \> d\theta ,\eqno(A3.3.3)
$$
where $\<\cdot,\cdot \>$ denotes a standard {\sl $\rr$-valued} scalar product
in $\cc^n$. Integrating by parts yields
$$
\int_0^1 \< \d^2_{\theta\theta} v_\lambda(\theta)
+ \lambda v_\lambda(\theta), w(\theta) \> d\theta + \<\d_\theta v_\lambda(
\theta), w(\theta) \>|_{\theta=1} - \<\d_\theta v_\lambda(
\theta), w(\theta) \>|_{\theta=0} =0.
$$
This implies that $v_\lambda$ belongs to $\ee_\lambda$ \iff
$\d^2_{\theta\theta} v_\lambda(\theta) + \lambda v_\lambda( \theta)=0$,
$\d_\theta v_\lambda(1) \perp W_1$, and $\d_\theta v_\lambda(0) \perp W_0$.
Since $J\st$ is $\<\cdot, \cdot\>$-orthogonal, we can conclude that
$J\d_\theta v_\lambda(\theta) \in \ee_\lambda$.

Positivity and compactness of $Q$ \wrt $\norm\cdot_{L^2}$ imply that
all $\ee_\lambda$ are finite dimensional and empty for $\lambda<0$. Further,
since $\d_\theta v=0$ for any $v\in \ee_0$, $\ee_0$ consists of constant
functions with values in $W_0 \cap W_1$.

Now let $u: \Theta(0,3) \to \cc^n$ be a holomorphic map with the boundary
condition (5.2). We can represent $u$ in the form $u(t, \theta)= \sum_\lambda
u_\lambda(t,\theta)$ with $u_\lambda(t, \cdot) \deff \pr_\lambda (u(t,\cdot))
\in L^{1,2}([0,3],\ee_\lambda)$. Since $J\d_\theta$ is an endomorphism of
every $\ee_\lambda$, every $u_\lambda(t,\theta)$ is also holomorphic. In
particular, $u_0$ is also holomorphic and constant in $\theta$. Thus $u_0$ is
constant.

Since $u$ is harmonic, \i.e., $(\d^2_{tt} + \d^2_{\theta \theta})u =0$, it
follows that 
$\d^2_{tt} u_\lambda(t,\theta)\allowbreak = \lambda u_\lambda(t,\theta)$. For
$\lambda > 0$ this yields $u_\lambda(t,\theta) = e^{+\sqrt \lambda t}
v^+_\lambda(\theta) + e^{-\sqrt\lambda t} v^-_\lambda(\theta)$ with $v^\pm
_\lambda(\theta) \in \ee_\lambda$. Fixing an orthogonal $\rr$-basis of
$\ee_\lambda $ $v^i_\lambda$, we write every $u_\lambda$ in the form
$$
u_\lambda(t,\theta)= \sum_i (a^i_\lambda e^{+\sqrt \lambda t} +
b^i_\lambda e^{-\sqrt \lambda t}) v^i_\lambda(\theta)
$$
with real constants $a^i_\lambda$, $b^i_\lambda$. Since $u_0(t,\theta)$
is constant, 
$\norm{du}^2_{L^2(\Theta_k)} = 2\norm{\d_\theta u}^2_{L^2(\Theta_k} = 
\sum_{\lambda,i} 2\lambda \int_{k-1}^k (a^i_\lambda e^{+\sqrt\lambda t}
+b^i_\lambda e^{-\sqrt\lambda t})^2 d\theta $. Here we used $(5.4)$ and
$L^2$-orthonormality of $v^i_\lambda $. This leads us to the problem of 
determining the smallest
possible constant $\gamma$ for the inequality
$$
\int_1^2 (ae^{\alpha t} + be^{-\alpha t})^2 dt
\le \msmall{\gamma\over2} \left(
\int_0^1 (ae^{\alpha t} + be^{-\alpha t})^2 dt +
\int_2^3 (ae^{\alpha t} + be^{-\alpha t})^2 dt \right)
\eqno(A3.3.4)
$$
with $a,\,b\in\rr$ for given $\alpha>0$. Integration gives
$$
a^2 e^{3\alpha}\msmall{ e^\alpha - e^{-\alpha} \over 2\alpha}
+ b^2 e^{-3\alpha}\msmall{ e^\alpha - e^{-\alpha} \over 2\alpha}
+2ab \le
$$
$$
\le \msmall{\gamma\over2} \left(
a^2 e^{3\alpha}\msmall{ (e^\alpha - e^{-\alpha}) (e^{2\alpha} + e^{-2\alpha})
\over 2\alpha}+
b^2 e^{-3\alpha}\msmall{ (e^\alpha
- e^{-\alpha}) (e^{2\alpha} + e^{-2\alpha})
\over 2\alpha}
+ 4ab \right)
$$
or equivalently
$$
a^2 e^{3\alpha}\msmall{ (e^\alpha - e^{-\alpha})
(e^{2\alpha} + e^{-2\alpha} -2/\gamma ) \over 2\alpha}
+
b^2 e^{-3\alpha}\msmall{ (e^\alpha - e^{-\alpha})
(e^{2\alpha} + e^{-2\alpha} -2/\gamma) \over 2\alpha}
$$
$$
+ 4ab(1-1/\gamma)
\ge 0
$$
The determinant of the last quadratic form in $a,b$ is
$$
\left(\!\!\msmall{ (e^\alpha - e^{-\alpha})
(e^{2\alpha} + e^{-2\alpha} -2/\gamma )
\over 2\alpha}\!\!\right)^{\!\!2} - 4(1-1/\gamma)^2
=
4\left(\!\!\msmall{ \sh \alpha (\ch\,2\alpha -1/\gamma )
\over \alpha} \!\!\right)^{\!\!2} - 4(1-1/\gamma)^2
$$
$$
\ge 4(\ch\,2\alpha -1/\gamma )^2 -4(1-1/\gamma)^2 =
4(\ch\,2\alpha -1)(\ch\,2\alpha +1 -2/\gamma).
$$
Thus, the  inequality (5.5) holds for every $a,b \in\rr$ provided $\gamma \ge
{2 \over 1 + \ch\,2\alpha}<1$.

Note that there exists a minimal positive eigenvalue $\lambda_1>0$ of $Q$.
Thus, the estimate (5.3) holds for $\gamma_W \deff {2 \over 1 + \ch(2
\sqrt{\lambda_1} )}<1$. \qed

\smallskip
\state Remark. A dependence of $\gamma_W$ as a function of $\lambda_1=
\lambda_1(W_0, W_1)$ shows that $\gamma_W<1$ can be chosen the same for all
pairs $(\wt W_0,\wt W_1)$ sufficiently close to $(W_0, W_1)$, provided $\dim
(\wt W_0 \cap \wt W_1) = \dim (W_0 \cap W_1)$. Vice versa, if we perform a
sufficiently small deformation of $(W_0, W_1)$ into $(\wt W_0, \wt W_1)$ with
$\dim (\wt W_0 \cap \wt W_1) < \dim (W_0 \cap W_1)$, then some $v \in
\ee_0(W_0, W_1)$ will wander into an eigenvector $\ti v \not\in \ee_0 (\wt
W_0, \wt W_1) = \wt W_0 \cap \wt W_1$, but with a sufficiently small eigenvalue
$\lambda_1(\wt W_0,\wt W_1) >0$, so that the best possible $\gamma_{\wt W}$
will be arbitrarily close to $1$. Thus the uniform separation of $\gamma_W$
from 1 under a small perturbation of $(W_0, W_1)$ is equivalent to a uniform
separation of $\lambda_1(W_0, W_1)$ from 0, and is equivalent to stability of
$\dim(W_0 \cap W_1)$.

Another phenomenon, also connected with spectral behavior, is that
for {\sl harmonic} $u(t, \theta)$ from $\ee_0= \{\,const\, \}$ does not
follow $u_0 = const$, but merely $u_0(t, \theta) = v_0 + v_1 t$ with $v_0,
v_1 \in \ee_0$, so that the inequality (5.3) is not true. This leads 
to a much more
complicated bubbling phenomenon with energy loss for harmonic and
harmonic-type maps, compare [S-U], [P-W], [Pa].

Note also that if $W_0$ and $W_1$ are {\sl affine} totally real subspaces of
$\cc^n$, then the inequality (5.3) for holomorphic 
$u: \Theta(0,3)\to \cc^n$ with
boundary condition (5.2) is in general not true. An easy example is a natural
imbedding $u:\Theta(0,3) \hook \cc$, with the nonconstant component $u_0 \equiv
u$ and with $\int_{\Theta(0,1)} |du|^2 = \int_{\Theta (1,2)} |du|^2 = \int_{
\Theta(2,3)} |du|^2 \not=0$. In general, those are $n$-dimensional totally
real affine planes $W_0$ and $W_1$ in $\cc^n$ with an empty intersection. This
can happen if $W_0$ and $W_1$ are parallel or {\sl skew}. The latter means
that the corresponding vector spaces $V_0$ and $V_1$ ($W_i = V_i + w_i$ for
some $w_i \in \cc^n$) are different. In both cases the intersection $V_0 \cap
V_1$ is not zero, because otherwise $W_0 \cap W_1 \not= \emptyset$ by
dimensional argumentation.

\smallskip
The considerations above show which properties should be controlled in order   to obtain a reasonable statement in the nonlinear case.

\state Definition A3.3.2. {\sl Let $X$ be a manifold with a Riemannian metric
$h$, $f_0 : W_0 \to X$ and $f_1 : W_1 \to X$ immersions and $A_0 \subset
W_0$, $A_1 \subset W_1$ subsets. We say that {\it $f_0 : W_0 \to X$ and $f_1
:W_1 \to X$ are $h$-uniformly transversal along $A_0$ and $A_1$ with
parameters $\delta >0$ and $M$} if for any $x_0 \in A_0$ and $x_1 \in A_1$
the following holds:

\smallskip
\sli either $\dist_h(f_0(x_0), f_1(x_1)) > \delta$;

\slii or there exists $x'_i \in A_i$ with $f_0( x'_0) =f_1(x'_1)$ such
that
$$
\dist_h(x_0,x'_0) + \dist_h(x_1,x'_1) \le M\,\dist_h(f_0(x_0), f_1(x_1)).
$$
}

\state Remark. Roughly speaking, the condition excludes the appearance 
of points,
where $W_0$ and $W_1$ are ``asymptotically parallel or skew" and ensures a
uniform lower bound for the angle between $W_0$ and $W_1$.

\medskip
Now let $X$ be a manifold with a Riemannian metric $h$, $J^*$ a continuous
almost complex structure on $X$, $f_0 : W_0 \to X$ and $f_1 : W_1 \to X$
immersions, and $A\subset X$, $A_0 \subset W_0$, $A_1 \subset W_1$ subsets.
Suppose that $J^*$ is $h$-uniformly continuous on $A$, $df_i: TW_i \to TX$
are $h$-uniformly continuous on $A_i$ and that $f_i$ are $h$-uniformly
transversal along $A_i$ with parameters $\delta= \delta(f_0,f_1)>0$ and $M=
M(f_0,f_1)$.

\medskip
\state Lemma A3.3.2. 
{\it There exist constants $\eps^b_2 =\eps^b_2 (\mu_{J^*},
\, f_0,\, f_1,\, \delta,\, M) >0$ and $\gamma^b = \gamma^b(\mu_{J^*},\, f_0,\,
f_1,\, \delta,\, M) <1$ such that for any continuous almost complex $\ti J$
with $\norm{\ti J-J^*}_{L^\infty(A)}<\eps_2^b$, any immersions $\ti f_i :W_i
\to X$ with $\dist(\ti f_i, f_i)_{C^1(A_i)} <\eps_2^b$, and any $\ti J
$-holomorphic map $u \in C^0 \cap L^{1,2}( \Theta(0,5), X)$ with $u(
\Theta(0,5)) \subset A$, $u|_{\d_i\Theta(0,5)} = f_i \scirc u^b_i$ for some
continuous $u^b_i : \d_i\Theta(0,5) \to A_i \subset W_i$ the conditions

\itemitem\sli $\norm{ du}_{ L^2(\Theta_i)} <\eps_2^b$;

\itemitem\slii $\ti f_i: W_i \to X$ are $h$-uniformly transversal along
$A_i$ with the same parameters $\delta$ and $M$

\noindent
imply the estimate
$$
\norm{du}^2_{L^2(\Theta_3)}\le \msmall{\gamma^b \over 2} \cdot
\left( \norm{du}^2_{L^2(\Theta_2)} +  \norm{du}^2_{L^2(\Theta_4)} \right)
\eqno(A3.3.5).
$$
}

\state Proof. Suppose the statement of the lemma is false. Then there should
exist a sequence of continuous almost complex structures $J_k$ with $\norm{
J_k- J^*} _{L^\infty (A)} \lrar 0$, a sequence of immersions $f_{k,i}: W_i
\to X$ with $f_{k,i} \lrar f_i$ in $C^1(A_i)$ such that $f_{k,i} :W_i \to X$
are $h$-uniformly transversal with the same parameters $\delta$ and $M$, and
a sequence of $J_k$-holomorphic maps $u_k \in C^0\cap L^{1,2}(\Theta(0,5),
X)$ with $u_k(\Theta(0,5)) \subset A$ and $u_k |_{\d_i\Theta(0,5)} = f_{k,i}
\scirc u^b_{n,i}$ for some continuous $u^b_{n,i} : \d_i\Theta(0,5) \to A_i
\subset W_i$ such that $\norm{du_k}^2_{L^2(\Theta(0,5))} \lrar 0$ and
$$
\norm{du_k}^2_{L^2(\Theta_3)}\ge \msmall{\gamma_k \over 2} \cdot
\left( \norm{du_k}^2_{L^2(\Theta_2)} +  \norm{du_k}^2_{L^2(\Theta_4)} \right)
\eqno(A3.3.6)
$$
with $\gamma_k = 1- 1/k$. {\sl Lemmas 5.3.1} and {\sl A3.2.1} provide that 
in this case 

\noindent $\diam_h(u_k( \Theta(1,4))) \allowbreak \lrar 0$.

Since $f_{k,i} :W_i \to X$ are $h$-uniformly transversal with the same
parameters $\delta$ and $M$, there should exist sequences $x_k \in A$, $x_{k,
0} \in A_0$, and $x_{k,1} \in A_1$ such that $x_k = u_0(x_{k,0}) =u_1(x_{k,
1})$ and $u_k(\Theta(1,4)) \subset B(x_k, r_k)$ with $r_k \lrar 0$. The $h
$-uniform continuity of $J^*$ implies that there exist $C^1$-diffeomorphisms
$\phi_k: B(x_k, r_k) \to B(0, r_k) \subset \cc^n$ with $\norm{J_k - \phi^*_k
J\st}_{L^\infty(B(x_k, r_k))} + \norm{h -\phi^*_k h\st} _{L^\infty(B(x_k,
r_k))} \lrar 0$.

Using $\phi_k$, we transfer our situation into $B(0, r_k) \subset \cc^n$ and
rescale it. Namely, we set $\alpha_k \deff \norm{du_k}_{L^2(\Theta_3)}$ and
define diffeomorphisms $\psi_k \deff {1\over \alpha_k} \scirc \phi_k: B(x_k,
r_k) \to B(0, R_k) \subset \cc^n$ with $R_k \deff \alpha_k\inv \cdot r_k$.
Note that by {\sl Lemmas 5.3.1} and {\sl A3.2.1} we have $\alpha_k = 
\norm{du_k} _{L^2(\Theta(2,3))} \le C \diam_h(u_k( \Theta(1,4))) \le C' r_k$, 
so that $R_k$ are uniformly bounded from below.

In $B(0, R_k)$ we consider Riemannian metrics $h_k \deff \alpha^{-2}_k \cdot
\psi_{k\,*}h_k$ (\ie pushed forward and $\alpha^{-2}_k$-rescaled $h_k$),
almost complex structures $J^*_k \deff \psi_{k\,*} J_k$ and $J^*_k
$-holomorphic maps $u^*_k \deff \psi_k \scirc u_k : \Theta(1,4) \to
B(0,R_k)$. Note that here we consider $h$ as a metric tensor; thus, multiplying
$h$ by $\alpha^{-2}$, we increase $h$-norms and $h$-distances in $\alpha^{-1}$
and not in $\alpha^{-2}$ times.

Then $\norm{du^*_k}_{L^2(\Theta_3, h_k)}=1$, $\norm{du^*_k}^2_{L^2(\Theta_2,
h_k)} + \norm{du^*_k}^2_{L^2(\Theta_4, h_k)} \le {2k\over k-1}$, and
$$
\norm{J^*_k - J\st}_{L^\infty(B(0, R_k), h_k)} =
\norm{J_k - \phi^*_kJ\st}_{L^\infty(B(x_k, r_k), h)} \lrar 0.
$$
The last equality uses the obvious relation
$$
\msmall{|F(\xi)|_{\alpha^{-2} \cdot h} \over |\xi|_{\alpha^{-2} \cdot h} }
= \msmall{\alpha\inv \cdot |F(\xi)|_h \over \alpha\inv \cdot |\xi|_\cdot h }
= \msmall{|F(\xi)|_h \over |\xi|_h }
$$
for any linear $F: T_xX \to T_xX$ and $\xi\not = 0\in T_xX$. In a similar way
we also obtain $\norm{h_k - h\st}_{L^\infty(B(0, R_k), h_k)}\lrar 0$.

Going to a subsequence, we may additionally assume that the tangent spaces
$d\psi_k \scirc df_i (T_{x_{k,i}} W_i)$, $i=1,2$, converge at some spaces
$W_i^* \subset \cc^n$. Since $W_i$ are uniformly totally real, $W_i^*$ are
also totally real linear subspaces in $\cc^n$. Since the maps $df_i : TW_0
\to TX$ are uniformly continuous on $A_i \subset W_i$, $f_{k,i} \lrar f_i$ in
$C^1(A_i)$, and since $r_k \lrar 0$, the images $W^*_{k,i} \deff \psi_k
\scirc f_{k,i}(B_{W_i}(x_{k,i}, r_k))$ of the balls $B_{W_i}(x_{k,i}, r_k)
\subset W_i$ are imbedded submanifolds of $\cc^n$ with $0 \in W^*_{k,i}$
which converge to $W^*_i$ in Hausdorff topology. Moreover, we can consider
$W^*_{k,i}$ as graphs of maps $g_{k,i}$ from subdomains $U_{k,i} \subset
W^*_i \cap B(0,R_k)$ to $W^*_i {}^\perp$ and for any fixed $R \le \inf
\{R_k\}$ the restrictions $g_{k,i} |_{W^*_i \cap B(0,R)}$ converge to zero
map from $W^*_i \cap B(0,R)$ to $W^*_i {}^\perp$.

The a priori estimates for the maps $u^*_k : \Theta(1,4) \to \cc^n$ provide
that for any $p<\infty$ the maps $u^*_k$ converge in weak- $L^{1,p}$-topology
to some $J\st$-holomorphic map $u^*: \Theta(1,4) \to \cc^n$. Furthermore, since
$u^*_k$ satisfy totally real boundary conditions $u^*_k|_{\d_i\Theta(1,4)}
\subset W^*_{k,i}$, the same is true for $u^*$, \i.e., $u^*|_{\d_i\Theta(1,4)}
\subset W^*_i$. Nice behavior of $W^*_{k,i}$ provides that on $\Theta_3$ we
also have a strong convergence; hence $\norm{du^*}_{L^2(\Theta_3)} = \lim
\norm{du^*_k}_{L^2(\Theta_3)} =1$. In particular, $u^*$ is not constant. On
the other hand, $\norm{du^*}^2_{L^2(\Theta_2)} + \norm{du^*}^2_{L^2(\Theta_4)}
\le \lim \norm{du^*_k}^2_{L^2(\Theta_2)} + \norm{du^*_k}^2_{L^2(\Theta_4)} \le
2$. The obtained contradiction to {\sl Lemma A3.3.1} shows that 
{\sl Lemma A3.3.2} is true. \qed

\medskip
Let $X$, $h$, $J$, $A$, $f_i: W_i \to X$, $A_i$, and the constant $\eps^b_2$
and $\gamma^b$ be as in {\sl Lemma A3.3.2}. Suppose that $\ti J$ is a 
continuous
almost complex structure on $X$ with $\norm{\ti J-J}_{L^\infty(A)} <\eps^b_2$,
and $\ti f_i: W_i \to X$ are totally real immersions with $\dist(\ti f_i, f_i)
_{C^1(A_i)} \le \eps^b_2$ such that $f_i$ are $h$-uniformly transversal
along $A_i$ with the same parameters $\delta$ and $M$ as $f_i: W_i \to X$.

\medskip
\state Corollary A3.3.3. {\it Let $u \in C^0\cap L^{1,2}(\Theta(0,l),X)$ be a
$\ti J$-holomorphic map such that $u(\Theta(0,l))\subset A$, $u|_{\d_i
\Theta(0,l)} = \ti f_i \scirc u^b_i$ for some continuous $u^b_i: \d_i\Theta(
0,l) \to A_i \subset W_i$ and $\norm{du}_{ L^2(Z_k)}<\eps_2$ for
any $k=1,\ldots,l$.

Let $\lambda_b>1$ be the (uniquely defined) real number with $\lambda_b 
= {\gamma^b \over 2}(\lambda_b^2+ 1)$. Then for $2\le k\le l-1$ the following  holds:
$$
\norm{du}^2_{L^2(\Theta_k)}
\le \lambda_b^{-(k-2)} \cdot \norm{du}^2_{L^2(\Theta_2)} +
\lambda_b^{-(l-1-k)} \cdot \norm{du}^2_{L^2(\Theta_{n-1})}.
\eqno(A3.3.7)
$$
}
%%%%%%%%%%%%%%%%%%%%%question%%%%%%%%%%%%%%%%%%%%

\state Proof. It is the same as in {\sl Lemma 5.3.5}. \qed

\medskip
The immediate corollary of this estimate is a lower bound of energy on
a nonconstant ``infinite strip''.

\smallskip
\state Lemma A3.3.4. {\it Let $X$, $h$, $J$, $A$, $f_i: W_i \to X$, $A_i$, the
constant $\eps^b_2$ and $\gamma^b$, also a structure $\ti J$, and
immersions $\ti f_i : W_i \to X$ be the same as in {\sl Corollary A3.3.3}. Let
$u \in C^0\cap L^{1,2}(\Theta(-\infty, +\infty),X)$ be a nonconstant $\ti
J$-holomorphic map such that $u(\Theta(-\infty, +\infty))\subset A$ and
$u|_{\d_i\Theta(0,l)} = \ti f_i \scirc u^b_i$ for some continuous $u^b_i :
\d_i\Theta(-\infty, +\infty) \to A_i \subset W_i$. Then $\norm{du}_{ L^2(
\Theta_k)} > \eps^b_2$ for some $k$. In particular, $\norm{du}_{ L^2(
\Theta(-\infty, +\infty))} > \eps^b_2$.
}

\state Proof. {\sl Corollary A3.3.3\/} provides that if $\norm{du}_{L^2(
\Theta_k)} \le \eps^b_2$ for all $k$, then $\norm{du}_{ L^2(\Theta_k)} 
\allowbreak =0$, \i.e., $u$ is constant. \qed

\medskip
Another consequence of {\sl Corollary A3.3.3} is a generalization of Gromov's
result about removability of boundary point singularity, see [G]. An
important improvement is the fact that the statement remains valid also when
one has {\sl different} boundary conditions on the left and on the
right of a singular point. One can see such a point $x$ as a {\sl corner
point} for the corresponding  complex  curve. A typical example appears
in symplectic geometry where one takes Lagrangian submanifolds as boundary
conditions.

\smallskip
Define the punctured half-disk by setting $\check\Delta^+ \deff \Delta^+\bs
\{0\}$. Define $I_- \deff (-1,0) \subset \d\check\Delta^+$ and
$I_+ \deff (0,+1) \subset \d\check\Delta^+$.

\smallskip
\state Corollary A3.3.5. {\sl (Removal of boundary point singularities).
\it Let $X$ be a manifold with a Riemannian metric $h$, $J$ a continuous
almost complex structure, $f_i : W_i \to X$, $i=1,2$, totally real
immersions and $A_i \subset W_i$ subsets. Let $u:(\check\Delta^+, J\st)
\to (X,J)$ be a holomorphic map. Suppose that

\item\sli $J$ is uniformly continuous on $A \deff u(\check\Delta)$ \wrt
$h$, and closure of $A$ is $h$-complete;

\item\slii $u$ satisfies boundary conditions of the form $u|_{I_+} = f_0
\scirc u^b_+$ and $u|_{I_-} = f_1 \scirc u^b_-$ with some continuous $u^b_+ :
I_+ \to A_0 \subset W_0$ and $u^b_- : I_- \to A_1 \subset W_1$;

\item\sliii $f_i$ are $h$-uniformly totally real on $A_i$ and $h$-uniformly
transversal along $A_i$;

\item\sliv there exists $k_0$ such that for all half-annuli $R^+_k\deff \{
z\in \Delta^+ :{1\over e^{\pi(k+1)}}\le | z| \le {1\over e^{\pi k}}\}$ with
$k\ge k_0$ one has $\norm{du}^2_{L^2(R^+_k)}\le \eps^b_2$, $\eps^b_2$ as in
{\sl Lemma A3.3.2}.

\smallskip\noindent
Then $u$ extends to origin $0\in \Delta^+$ as an $L^{1,p}$-map for some
$p>2$.
}

\state Proof. Using the holomorphic map $\exp : \Theta(0,\infty) \to
\check\Delta^+$, $\exp(\theta + \isl t) \deff e^{\pi(-t + \isl\theta)}$, we
can reduce our situation to the case of the holomorphic map $u^* \deff u
\scirc \exp$ form ``infinite strip'' $\Theta(0,\infty)$. By {\sl Corollary
A3.3.3}, for $k\ge k_0$ we obtain the estimate 
$\norm{ du^*} _{L^2(\Theta_k)} \le
\lambda_b^{-(k-k_0)/2} \norm{ du^*}_{L^2(\Theta_{k_0})}$ with some $\lambda_b
>1$. This is equivalent to the estimate $\norm{ du} _{L^2(R^+_k)} \le
\lambda_b^{-(k-k_0)/2} \norm{ du}_{L^2(R^+_{k_0})}$. {\sl Lemmas 5.3.1} and
{\sl A3.2.1} and the scaling property of $L^p$-norms provide the estimate
$$
\norm{ du}_{L^p(R^+_k)} \le C e^{-k(\log\lambda_b/2 + \pi(2/p-1))}.\eqno(A3.3.8)
$$
Thus, $du \in L^p(\Delta^+)$ for any $p$ with $\log\lambda_b/2 > \pi(1-2/p)$,
which means $p < {4\pi \over 2\pi - \log \lambda_b}\cdot$ \qed

\smallskip
\state Remark. Unlike the ``inner" and smooth boundary cases, it is possible
that the map $u$, as in {\sl Corollary A3.3.5\/}, is not $L^{1,p}$-regular in 
the
neighborhood of a ``corner point'' $0\in \Delta^+$ for some $p>2$. For example,
the map $u(z) = z^\alpha$ with $0<\alpha <1$ satisfies totally real boundary
conditions $u(I_+) \subset \rr$, $u(I_-) \subset e^{\alpha\pi \isl}\rr$ and
is $L^{1,p}$-regular only for $p< p^* \deff {2\over1-\alpha}\cdot$

\medskip
As in the ``inner'' case, for the proof of the boundary compactness theorem we
need a description of a convergence of a ``infinitely long strip''. Let
$X$ be a manifold with a Riemannian metric $h$, $J$ a continuous almost
complex structure, $A\subset X$ a closed $h$-complete subset such that $J$
is $h$-uniformly continuous on $A$, and let $\{J_n\}$ be a sequence of almost
complex structures converging $h$-uniformly on $A$ to $J$. Also let $f_0: W_0
\to X$ and $f_1: W_1\to X$ be immersions, $A_i \subset W_i$ subsets such
that $df_i$ are uniformly $h$-uniformly totally real on $A_i$ and $f_i$ are
$h$-uniformly transversal along $A_i$. Let $f_{n,i} : W_i \to X$ be totally
real immersions, which $C^1$-converge to $f_i$ on $A_i$ such that $f_{n,0}$
and $f_{n,1}$ are $h$-uniformly transversal along $A_i$ with uniform in $n$
parameters $\delta$ and $C^*$. Finally, let $\{l_n\}$ be a sequence of
integers with $l_n\to \infty$, and $u_n: \Theta(0,l_n) \to X$ a sequence of
$J_n$-holomorphic maps, satisfying boundary conditions $u_n|_{\d_i\Theta(0,
l_n)} = f_{n,i}\scirc u^b_{n,i}$ with some continuous $u^b_{n,i}:
\d_i\Theta(0,l_n) \to A_i \subset W_i$.

\state Lemma  A3.3.6. {\it In the described situation, suppose additionally 
that $u_n(\Theta(0,l_n)) \allowbreak \subset A$ and $\norm{du_n} _{ L^2( 
\Theta_k) }\le \eps^b_2$ for all $n$ and $k\le l_n$. Take a sequence $k_n\to
\infty$ such that $k_n<l_n-k_n\to \infty$. Then

\item{\sl1)} $\norm{du_n}_{L^2(\Theta(k_n,l_n-k_n))}\to 0$ and
$\diam\bigl(u_n (\Theta(k_n,l_n-k_n))\bigr)\to 0$.

\item{\sl2)}
There is a subsequence $\{ u_n \}$, still denoted $\{ u_n \}$ such that both
$u_n |_{\Theta(0,k_n)}$ and $u_n |_{\Theta(k_n,l_n)}$ converge in $L^{1,p}
$-topology on compact subsets in $\check\Delta^+ \cong \Theta(0, +\infty)$ to 
$J^*$-holomorphic maps $u^-_\infty$ and $u^+_\infty$. Moreover, both
$u^+_\infty$ and $u^-_\infty$ extend to the origin and $u^+_\infty(0)=
u^-_\infty(0)$.
}

\state Proof. It follows from the above considerations.

\bigskip\noindent
{\bigsl A3.4. Gromov Compactness for Curves with Totally Real Boundary 
Conditions.}

\smallskip
Let us turn to the Gromov Compactness Theorem for curves with boundary on
totally real submanifolds. To give a precise statement we need to modify the
definition of the Gromov convergence ({\sl Definition 4.1.6}). The reason to do
this is the following. Considering {\sl open} curves $C_n$ 
with changing complex
structures, we want to fix some kind of a common ``neighborhood of infinity"
$i_n: C^*\hookrightarrow C_n$ of every $C_n$. Thus, we can imagine that all
changes of complex structure take place ``outside of infinity'', \i.e., in a
relatively compact part $C_n \bs i_n(C^*) \Subset C_n$. This is done to
insure that $C_n$ do not approach infinity in an appropriate moduli space.

On the other hand, it is more natural to consider curves $(\barr C_n, u_n)$
with totally real boundary conditions as compact objects without ``infinity''.
In fact, in this case the behavior of $u_n$ near the boundary $\d C_n$ can be
controlled. The obtained a priori estimates near a ``totally real boundary''
can be viewed as a part of such a ``control''. So for curves with totally real
boundary conditions we can hope to extend the Gromov convergence up to the 
boundary.

Further, as in the ``inner case'', an appropriate modification of the Gromov
convergence in this case should allow boundary bubbling and the appearance of
boundary nodes. This means, however, that the structure of the boundary can
change during the approach to the limit curve and 
cannot be considered as fixed.
Instead, one should fix a type of boundary conditions. We shall
consider the following general situation.

Let $u_n: \barr C_n \to X$ be a sequence of stable $J_n$-complex curves
 over $X$
with parameteri\-zations $\delta_n: \barr\Sigma \to \barr C_n$. Also let
$\bfbeta =\{ \beta_i \}_{i=1}^m$ be a collection of arcs $\beta_i$ in $\d
\Sigma$ such that $\cup_{i=1}^m \beta_i= \d\Sigma$ and that
the interiors of $\beta_i$ are mutually disjoint and do not intersect
pre-images of boundary nodal points of $C_n$. Let further $\{ W_i\} _{i=1}^m$ 
be a collection of real $n$-dimensional manifolds,
$f_{n,i} : W_i \to X$ a sequence of totally real immersions and $u^b_{n,i}:
\beta_{n,i} \to W_i$ a sequence of continuous maps from $\beta_{n,i} \deff
\delta_n( \beta_i)$.  Then ${\mib W}_n \deff \{ (W_i,
f_{n,i} )\} _{i=1}^m$ are totally real boundary conditions on $(\barr C_n,
u_n)$ of the same type $\bfbeta$.

\state Definition A3.4.1. {\sl In the situation above we say that the sequence
of boundary conditions ${\mib W}_n$ of the same type $\bfbeta$ {\it
converges $h$-uniformly transversally} to $J^*$-totally real boundary
conditions ${\mib W}$ on subsets $A_i \subset W_i$ if

\sli ${\mib W}= \{ (W_i, f_i) \}_{i=1}^m$, where $f_i: W_i \to X$ are
$J^*$-totally real immersions;

\slii $f_{n,i}$ converge to $f_i$ in $C^1$-topology and this
convergence is $h$-uniform on $A_i$;

\sliii for any $n$ immersions $\{f_{n,i}\}_{i=1}^m$ are mutually
$h$-uniformly transversal along $A_i$ with parameters $\delta>0$ and $M$, and
these parameters are independent of $n$.
}

\smallskip
Note that condition \sliii implies that the limit immersions $f_i$ are
also mutually $h$-uniformly transversal along $A_i$ with the same parameters
$\delta>0$ and $M$.

\state Definition A3.4.2. 
{\sl We say that the sequence $(\barr C_n, u_n)$ {\it
converges up to the boundary} to a stable $J^*$-holomorphic curve $(\barr
C_\infty, u_\infty)$ over $X$ if the parameterizations $\sigma_n: \barr\Sigma
\to \barr C_n$ and $\sigma_\infty: \barr\Sigma \to \barr C_\infty$ can be
chosen in such a way that the following holds:

\sli $u_n\scirc \sigma_n$ converges to $u_\infty\scirc \sigma_\infty$ in
$C^0( \barr\Sigma, X)$-topology;

\slii if $\{ a_k \}$ is the set of the nodes of $C_\infty$ and $\{ \gamma_k
\}$, $\gamma_k \deff \sigma_\infty\inv(a_k)$ are the corresponding circles
and arcs in $\barr\Sigma$, then on any compact subset $K\comp \barr\Sigma \bs
\cup_k\gamma_k$ the convergence $u_n\scirc \sigma_n\to u_\infty\scirc
\sigma_\infty$ is $L^{1,p}(K, X)$ for all $p< \infty$;

\sliii for any compact subset $K\comp \barr\Sigma \bs \cup_k\gamma_k$ there
exists $n_0=n_0(K)$ such that $ \sigma_n^{-1}(\{ a_k \}) \cap K= \emptyset$
for all $n\ge n_0$ and complex structures $\sigma_n^*j_{C_n}$ smoothly
converge to $\sigma_\infty^*j_{C_\infty}$ on $K$.
}

\state Theorem A3.4.1. 
{\it Fix a metric $h$ on $X$, and an $h$-complete subset
$A\subset X$, and subsets $A_i \subset W_i$. Suppose that

\item{\sl a)} $J_n$ are continuous almost complex structures on $X$,
converging $h$-uniformly on $A$ to a continuous almost complex structure
$J^*$;

\item{\sl b)} $u_n(C_n) \subset A$ and $\area [u_n (C_n)]\le M$ with a
constant $M$ independent of $n$;

\item{\sl c)} ${\mib W}_n \deff \{(W_i, f_{n,i}) \}_{i=1}^m$ are totally real
boundary conditions of the same type $\bfbeta = \{ \beta_i \}_{i=1}^m$ such
that ${\mib W}_n$ converge $h$-uniformly transversally to a boundary condition
${\mib W}= \{(W_i, f_i) \}_{i=1}^m$ on subsets $A_i \subset W_i$;

\item{\sl d)} immersions $f_i: W_i \to (X,J^*)$ are $h$-uniformly totally
real along $A_i$;

\item{\sl e)} there exist maps $u^b_{i,n}: \beta_i \to  A_i\subset W_i$,
realizing boundary conditions ${\mib W}_n$.

\smallskip
Then there exits a subsequence of $\{( \barr C_n, u_n )\}$, still denoted
$\{( \barr C_n, u_n )\}$, and para\-meteri\-zations $\sigma_n: \barr \Sigma
\to \barr C_n$ such that $(C_n, u_n, \sigma_n)$ converges up to the boundary to
a stable $J^*$-holomorphic curve $(\barr C_\infty, u_\infty, \sigma_\infty)$
over $X$.

If, in addition, $A_i \subset W_i$ are $f_i^*h$-complete, then the limit
curve $(\barr C_\infty, u_\infty)$ satisfies real boundary conditions $\mib
W$ with maps $u^b_i: \beta_i \to A_i \subset W_i$.
}

\medskip
Our main idea of the proof is to apply arguments used in the demonstration of
{\sl Theorem 1.1}. To realize this, we use the following trick. We replace 
every pair $(C_n, u_n)$ by a triple $(C^d_n, \tau_n, u^d_n)$, where $C^d_n$ 
is the {\sl Schottky double} of $C_n$ with an antiholomorphic involution 
$\tau_n$ and $u^d_n: C_n^d \to X$ a $\tau_n$-invariant map. Then we shall show
the changes of all the constructions in the proof to make them 
$\tau_n$-invariant in an appropriate sense. In particular, the convergence 
$(C^d_n, \tau_n, u^d_n) \lrar (C^d_\infty, \tau_\infty, u^d_\infty)$ will be 
equivalent to the convergence $(C_n, u_n) \lrar (C_\infty, u_\infty)$.

\smallskip
We start with a construction of the Schottky double of a nodal curve $\barr C$
with boundary. Take two copies $\barr C^+ \equiv \barr C$ and $\barr C^-$ of
$\barr C$. Equip $\barr C^-$ with the opposite complex structure, so that the
identity map $\tau: \barr C^+ \to \barr C^-$ now becomes antiholomorphic.
Glue $\barr C^+$ and $\barr C^-$ together along their boundaries, identifying
$\d C^+$ and $\d C^-$ by means of the identity map $\tau: \d C^+ \buildrel
\cong \over \lrar \d C^-$. The union $C^d \deff \barr C^+ \cup_{\d C} \barr
C^-$ obeys the unique structure of a closed nodal curve compatible with
imbeddings $\barr C^\pm \hook C^d$. The boundary $\d C$ becomes the fixed
point set of $\tau$.

The map $\tau$ induces an antiholomorphic involution of $C^d$ which we also
denote by $\tau$. We call the obtained curve $C^d$ the {\sl Schottky double} of
$\barr C$. Note that every boundary nodal point $a_i\in \d C$ defines a $\tau
$-invariant nodal point $a_i$ on $C^d$, whereas an inner nodal point $b_i \in
C$ defines a pair of nodal points $b_i^\pm$ on $C^d$ interchanged by $\tau$.
If $\sigma: \barr\Sigma \to \barr C$ is a parameterization of $\barr C$, then
we obtain in an obvious way the double $\Sigma^d$ with the involution $\tau:
\Sigma^d \to \Sigma^d$ and the parameterization $\sigma^d: \Sigma^d \to C^d$
compatible with the involutions.

\state Remark.
The introduced notation $C^d$ for the {\sl Schottky double} of a nodal curve
$\barr C$ with boundary coincides with that for the {\sl holomorphic
double}, used in {\sl Section 2}. Since in the present section 
only the Schottky
double is considered, this should not lead to confusion.

\smallskip
Suppose, additionally, that an almost complex structure $J$ on $X$ and a
$J$-holo\-mor\-phic map $u: \barr C \to X$ are given. Suppose, also, that the
curve $(\barr C,u)$ satisfies the totally real boundary conditions $\mib W$ of
type $\bfbeta$. In particular, $\bfbeta$ defines a certain system of arcs
$\{\beta_i\}$ on $\d C$. In order to take into account the type of boundary
conditions, we fix the ends of $\beta_i$ which are not boundary nodal points
of $\barr C$ and declare these points as marked points of $C^d$. Note that
these and the nodal points are the only ``corner'' points of $(\barr C,
u)$. The latter means that in a neighborhood of these points the map $u$ cannot
be $L^{1,p}$-smooth for all $p< \infty$. The example in the {\sl Remark} 
following {\sl Corollary A3.3.5} explains the notion of a ``corner point''. 
Considering the
Schottky double, we shall always equip $C^d$ with this set of marking points.
Note also that every boundary circle of $\barr C$ contains at least one nodal
or marked point as above.

For $(\barr C,u)$ as above, we extend the $J$-holomorphic map $u: \barr C \to
X$ to a map $u^d: C^d \to X$ by setting $u^d(x) \deff u(\tau(x))$ for $x\in
C^-$. By the construction, $u^d$ is $\tau$-invariant, $u^d \scirc \tau= u^d$,
but $u^d$ {\sl is not $J$-holomorphic} (with the only trivial exception $u
\equiv const$). However, the analysis already described in this 
section provides
necessary $L^{1,p} $-estimates for $u^d$, at least for some $p^*>2$.

\smallskip
In the situation of {\sl Theorem A3.4.1}, such an exponent $p^*>2$ can be chosen
to be the same for all curves $(\barr C_n, u_n)$. 
This depends only on the topology
of $\barr C_n$ and the geometry of immersions $f_n: W_n \to X$. In
particular, every $u_n^d$ is continuous.

\smallskip
The next step of the proof is to find a $\tau_n$-invariant decomposition of
$C^d_n$ into pants. This implies that the corresponding graph $\Gamma_n$
becomes $\tau_n$-invariant. In the construction which follows we shall use
the fact that $\tau_n$ is an isometry on the union of the non-exceptional
components of $C^d_n$. This is provided by the uniqueness of the intrinsic
metric.

\state Lemma A3.4.2. 
\it Let $C$ be a nodal curve with boundary, $\sigma: \barr
\Sigma \to \barr C$ a parameterization, and $\{ x_i \}_{i=1}^m$ a set of
marked points on the boundary $\d C$. Let $C^d$ be the Schottky double of $C$
with the antiholomorphic involution $\tau$.

Then there exists a $\tau$-invariant decomposition of $C^d\bs \mapo$ into
pants such that the intrinsic length of corresponding boundary circles is
bounded by a constant $l^+$ depending only on genus $g$ of $\Sigma^d$ and
the number of marked points $m$.

Moreover, every {\sl short} geodesic appears as a boundary circle of some
pants of the decomposition.\rm

\state Remark. Recall (see {\sl Remark} on page 25) that a closed geodesic
$\gamma$ is called {\sl short} if $\ell(\gamma) < l^*$, where $l^*$ is the
universal constant $l^*$ with the following property. For any simple closed
geodesics $\gamma'$ and $\gamma''$ on the conditions $\ell(\gamma') <l^*$ and
$\ell( \gamma'') <l^*$ imply $\gamma' \cap \gamma'' = \emptyset$.

\state Proof. Since the genus of the parameterizing real surface $\Sigma^d$ and
the number of marked points is fixed, we obtain a uniform upper bound on the
possible genera and the number of marked points of non-exceptional components
of $C^d$, as well as on the number of exceptional components. This implies that
there exists a decomposition of every non-exceptional component $C_i$ of
$C^d$ into pants $S_\alpha$ such that the intrinsic length of boundary
circles of $S_\alpha$ is bounded by the constant $l^+$ depending only on $g$
and $m$. The idea of the proof of our lemma is to show that the construction
of such a decomposition, given in [Ab], Ch.II, \S\.3.3, can be modified to
produce a $\tau$-invariant decomposition.

Let us first describe the construction itself, say, for a given smooth curve
$C^*$ with marked points $\{x_i\}$ of non-exceptional type. The procedure is
done inductively by choosing at every step a non-trivial simple closed
geodesic $\gamma_{J^*} \subset C^* \bs\mapo$, disjoint from an already chosen
geodesic $\gamma_j$, $j<J^*$. Moreover, at every step there exists a geodesic
$\gamma_{J^*}$ as above whose intrinsic length is bounded by a constant
$l^+_{J^*}$ depending only on the genus of $C^*$, the number of marked
points, and the maximum of the lengths of the already chosen geodesics
$\gamma_j$, $j<J^*$.

\smallskip
Take any non-exceptional component $C^d_i$ of $C^d$. Two cases can occur:
either $C^d_i$ is $\tau$-invariant, or $\tau \bigl( C^d_i \bigr)$ is another
component $C^d_{i'}$. Two separate cases are distinguished:
$C^d_i$ intersects the boundary $\d C$ (first case) or not (second one).

The existence of $\tau$-invariant decomposition into pants for every pair of
non-exceptional components $C^d_i$ and $\tau(C^d_i) \not= C^d_i$ is obvious.
We choose an appropriate decomposition of $C^d_i$ and transfer it on $\tau(
C^d_i)$ by means of $\tau$.

\smallskip
It remains to consider the case of a $\tau$-invariant non-exceptional
component $C^d_i$.

Suppose that at some step we have already chosen a $\tau$-invariant set
$\{\gamma_1, \ldots, \allowbreak \gamma_{J^*-1}\}$ of simple disjoint
geodesics on $C^d_i \bs \mapo$. Take a simple geodesic $\gamma$ of the length
$\ell(\gamma) \le l^+_{J^*}$, where $l^+_{J^*}$ is the upper bound introduced
above. By the construction of the double $C^d$, the fixed point set of $\tau$
on $C^d_i$ is $C^d_i \cap \d C$ and is non-empty. Denote $C_i \deff C \cap
C^d_i$, so that $C^d_i \cap \d C= \d C_i$. Note that any boundary circle of
$C_i$ contains at least one marked point of $C^d_i$. Consequently, it has an
infinite length \wrt the intrinsic metric on $C^d_i \bs\mapo$. Thus the
chosen geodesic $\gamma$ cannot lie on $\d C_i$. Only three cases can happen.

\smallskip\noindent
{\sl Case 1.} $\gamma$ is disjoint from $\d C_i$. Then $\gamma$ lies either in
$C_i$ or in $\tau (C_i)$. In any case, $\gamma\cap \tau(\gamma) = \emptyset$.
Thus, we can set $\gamma_{J^*} = \gamma$ and $\gamma_{J^*+1}= \tau(\gamma)$,
obtaining the $\tau$-invariant set $\{\gamma_1, \ldots, \gamma _{J^*+1}\}$ of
simple disjoint geodesics. This will be discussed in the next two 
steps of our construction.

\smallskip\noindent
{\sl Case 2.} $\gamma \cap \d C_i \not=0$ and $\gamma$ is $\tau$-invariant.
We set $\gamma_{J^*}=\gamma$ and proceed inductively. Note that in
this case $\gamma \cap \d C_i$ consists of 2 points, in which $\gamma$ is
orthogonal to $\d C_i$.

\smallskip\noindent
\line{\vtop{\hsize=.5\hsize\noindent
{\sl Case 3.} This time $\gamma \cap \d C_i \not=0$, but $\gamma \not=
\tau(\gamma)$. Define arcs $\gamma^+ \deff \gamma \cap \barr C_i$ and
$\gamma^- \deff \gamma \cap \tau(\barr C_i)$, the parts of $\gamma$ inside
and outside of $C_i$ (see Fig.~11). Consider the following free homotopy
classes of closed circles on $C^d_i$:

1) $[\ti\gamma_1] \deff [\gamma^+ \cup \tau(\gamma^-)]$;

2) $[\ti\gamma_2] \deff [\gamma^+ \cup \tau(\gamma^-)]$;

3) $[\ti\gamma_3] \deff [\gamma^- \cup \tau(\gamma^-)]$;

4) $[\ti\gamma_4] \deff [\gamma^+ \cup \tau(\gamma^-) \cup \gamma^- \cup
\tau(\gamma^+)]$.
%\par
%\parshape=2 0pt \hsize 0pt 2\hsize \noindent
}
\hss
\vtop{\xsize=.46\hsize\hsize=\xsize\nolineskip\rm
%\vskip15pt
%\parshape=3 0pt .5\hsize 0pt .5\hsize 0pt \hsize
\putm[.47][-.01]{\ti\gamma_1}%
\putm[.02][.32]{\ti\gamma^+_2}%
\putm[.95][.42]{\ti\gamma^+_3}%
\putm[.15][.25]{\gamma^+}%
\putm[.61][.23]{\tau(\gamma^-)}%
\putm[.40][.17]{\ti\gamma^+_4}%
\noindent
\epsfxsize=\xsize\epsfbox{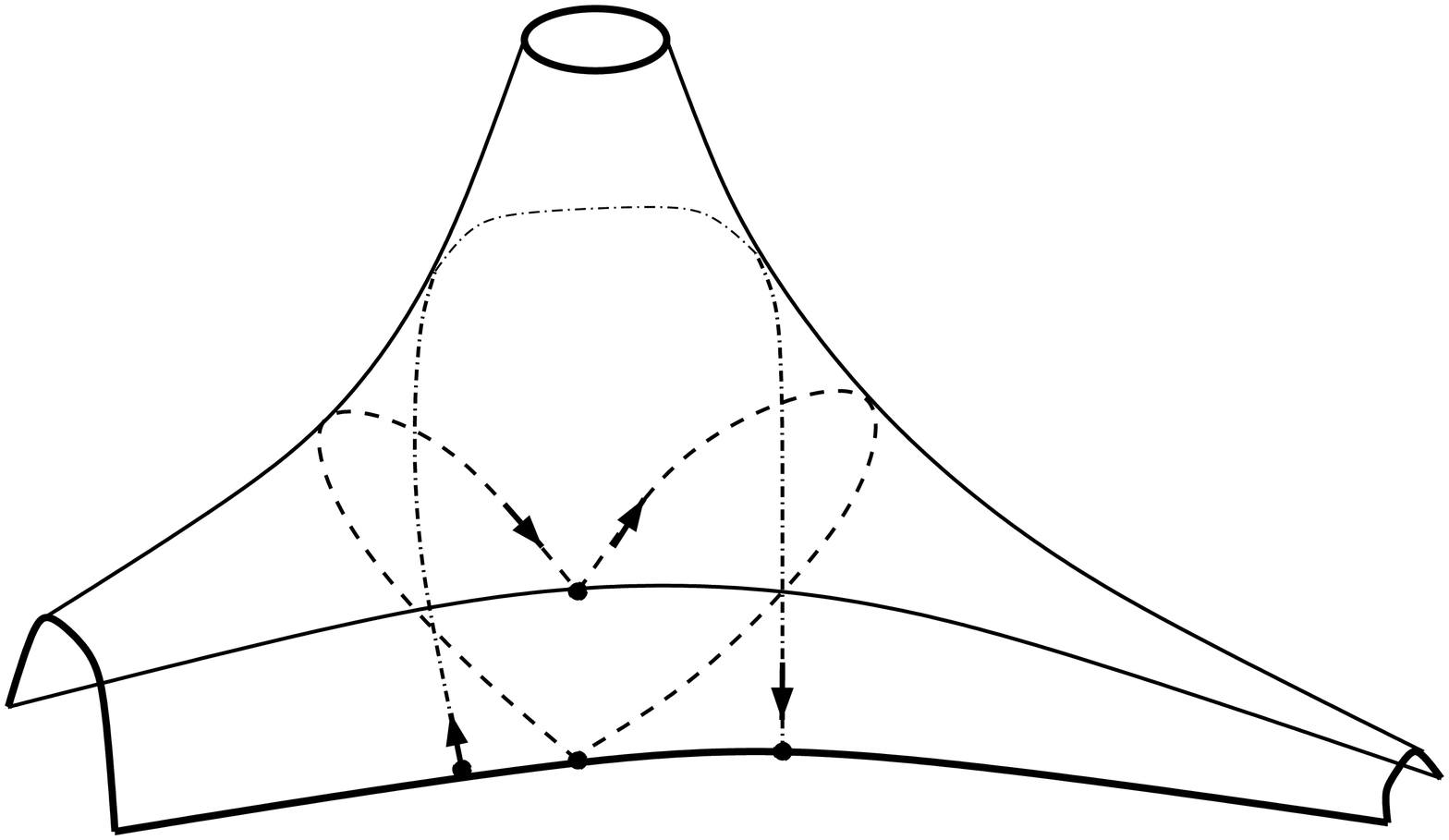}
\smallskip
\centerline{Fig.~11. Geodesics on $C_i$.}
}}

\vskip2pt\noindent
The last expression means that we move along corresponding arcs in the
prescribed order, as it is shown on Fig.~11. Note that only one part of
$C^d_i$ is drawn, namely $C_i$. The rest of the picture is symmetric \wrt
the involution $\tau$. Thus we can see only half of the geodesics in classes
$[\ti\gamma_i]$, $i=2,3,4$.

Each of the classes $[\ti\gamma_k]$ is either represented 
by a closed geodesic or
corresponds to a curve which winds around some marked point of 
$C^d_i$. To shorten
notation, we say in the last case that the class $[\ti\gamma_i]$ corresponds
to a marked point of $C^d_i$.

If one of the classes $[\ti\gamma_k]$, $k=1,2,3$, is represented by the
geodesic $\ti\gamma_k$, which is different and disjoint from the already
chosen geodesics $\gamma_j$, $j<J^*$, then we can set $\gamma_{J^*} =
\ti\gamma_k$. If $k=1$ we also set $\gamma_{J^*} = \ti\gamma_1$ and
$\gamma_{J^*+1} = \tau(\ti\gamma_1)$. Then we proceed inductively.

\smallskip
To finish the proof it remains to consider the following situation. Under
the conditions of {\sl Case 3}, each of the classes $[\ti\gamma_k]$, $k=1,2,3$,
either corresponds to a marked point or is represented by a closed geodesic
$\ti\gamma_k$, which intersects or coincides with one from the already chosen
geodesics $\gamma_j$, $j<J^*$.

We state that, in fact, a proper intersection cannot happen, \i.e., each class
$[\ti\gamma_k]$, $k=1,2,3$, either corresponds to a marked point or is
represented by an already chosen geodesic $\gamma_j$, $j<J^*$. To show this
we note that $\gamma_j \cap \tau(\gamma) =\emptyset$ for all $j<J^*$.
Otherwise, we could have a contradiction with the conditions $\gamma_j \cap
\gamma =\emptyset$ and $\tau$-invariance of the set of the geodesics
$\gamma_j$, $j<J^*$. Consequently, each class $[\ti\gamma_k]$ is represented
by a circle $\alpha_k \subset C^d_i \bs\mapo$, $k=1,2,3,4,$ with $\alpha_k
\cap \gamma_j = \emptyset$.

Now assume that the proper intersection of $\ti\gamma_k$ and some $\gamma_j$,
$j<J^*$ occurs. Let $\ell_k\deff \ell(\ti\gamma_k)$ be the intrinsic
metric of $\ti\gamma_k$. As in the proof of {\sl Lemma 4.3.2} construct the
annulus $A= \{ (\rho, \theta) : |\rho| < {\pi^2 \over \ell} \} \times \{0\le
\theta \le 2\pi\}$ with the metric $({\ell_k \over 2\pi} / \cos {\ell_k \rho
\over 2\pi})^2 (d\rho^2 + d\theta^2)$ and an isometric covering of $C^d_i
\bs\mapo$ by $A$, which sends the geodesic $\beta_k \deff \{ \rho=0\} \subset
A$ onto $\ti\gamma_k \subset C^d_i$. Find a lift of $\gamma_j$ to a geodesic
line $L_j \subset A$ with $L_j\cap \beta_k \not= \emptyset$ and a lift of a
circle $\alpha_k$ to a circle $\ti\alpha_k \subset A$ homotopic to $\beta_k$.
Then the intersection $L_j \cap \beta_k$ must consist of exactly one point,
and, consequently, the homology intersection index $[L_j] \cdot [\beta_k]$ 
is equal to $\pm1$. This would imply that $[L_j]\cdot [\ti\alpha_k] =[L_j] 
\cdot [\beta_k] \not =0$ and consequently $L_j \cap \ti\alpha_k \not= 
\emptyset$. But this would contradict $\gamma_j \cap \alpha_k =\emptyset$.

Summing up, we see that in the situation we are considering we must have 
the picture
of {\sl Fig.~11}. Namely, both geodesics $\gamma$ and $\tau(\gamma)$ lie in a
$\tau $-invariant domain $\Omega$ on $C^d_i$ with four components of the
boundary; these components of $\d\Omega$ are either marked points or geodesics
corresponding to the classes $[\ti\gamma_1]$, $[\tau(\ti\gamma_1)]$, $[\ti
\gamma_2]$, $[\ti\gamma_3]$; finally, every boundary circle of $\Omega$
is one of the geodesic $\gamma_j$. We conclude that the class 
$[\ti\gamma_4]$ is
represented by a $\tau$-invariant geodesic $\ti\gamma_4$, which can be chosen
at this step of the construction of $\tau$-invariant decomposition of $C^d_i$
into pants.

Note that by construction for the intrinsic length of $\gamma_{J^*}$ we obtain
$\ell(\gamma_{J^*}) \le 2\ell(\gamma) \le 2l^+_{J^*}$. This means that in our
construction we do not lose control of the intrinsic length of the chosen
geodesics. This shows the existence of a constant $l^+$ stated in the
lemma.

Finally, the definition of a {\sl short} geodesic provides that the geodesic
$\gamma$ in {\sl Case 3} above cannot be short. This implies that the set of
short geodesics on $C^d$ is disjoint. Since the involution $\tau$ is an
isometry, the set of short geodesics on $C^d$ is also $\tau$-invariant. Thus,
in our construction of decomposition into pants we can start with this set of
geodesics. This shows the last statement of the lemma. \qed

\smallskip
\state Remark. To explain the meaning of {\sl Lemma A3.4.2}, let us consider
pants $S$ with a complex structure $J_S$ and an antiholomorphic involution
$\tau$ acting on $S$. It is easy to see that only two types of such an
action, illustrated by Figs.\ 8\.a) and 8\.b), are possible.

\medskip\smallskip
\line{%
\vbox{\xsize=.47\hsize
\hsize=\xsize\nolineskip
\putm[.14][.11]{\gamma_1}%
\putm[.14][.60]{\gamma_2}%
\putm[.92][.37]{\gamma_3}%
\putm[.52][.36]{\beta}%
\putm[.52][.10]{S^+}%
\putm[.52][.57]{S^-}%
\noindent
\epsfxsize=\hsize\epsfbox{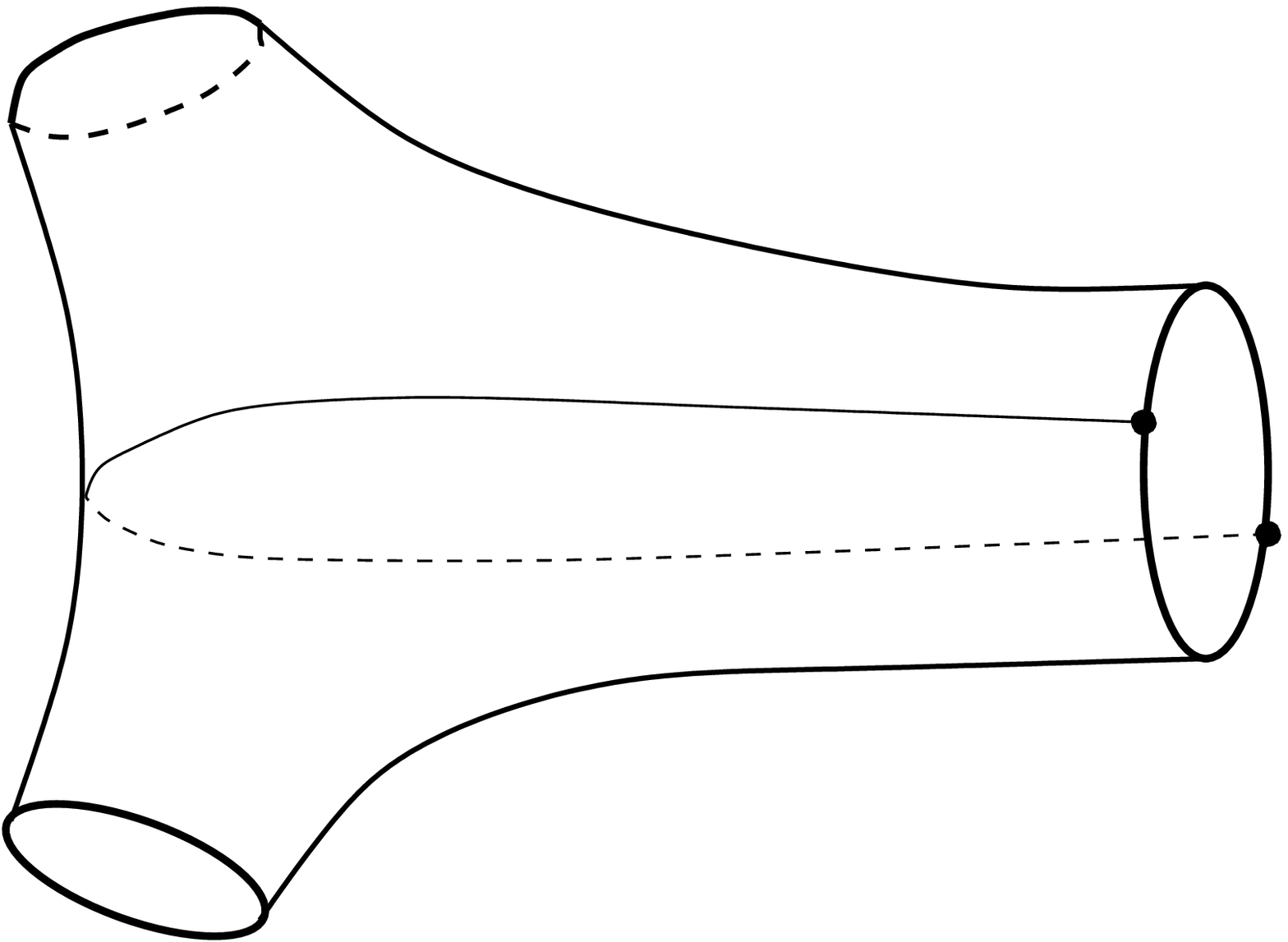}}%
\hfil
\vbox{\xsize=.47\hsize
\hsize=\xsize\nolineskip
\putm[.0][-.035]{\gamma_1}%
\putm[.28][.45]{\gamma_2}%
\putm[.90][-.045]{\gamma_3}%
\putm[.60][.22]{\beta_1}%
\putm[.52][.14]{\beta_2}%
\putm[.28][.12]{\beta_3}%
\noindent
\epsfxsize=\hsize\epsfbox{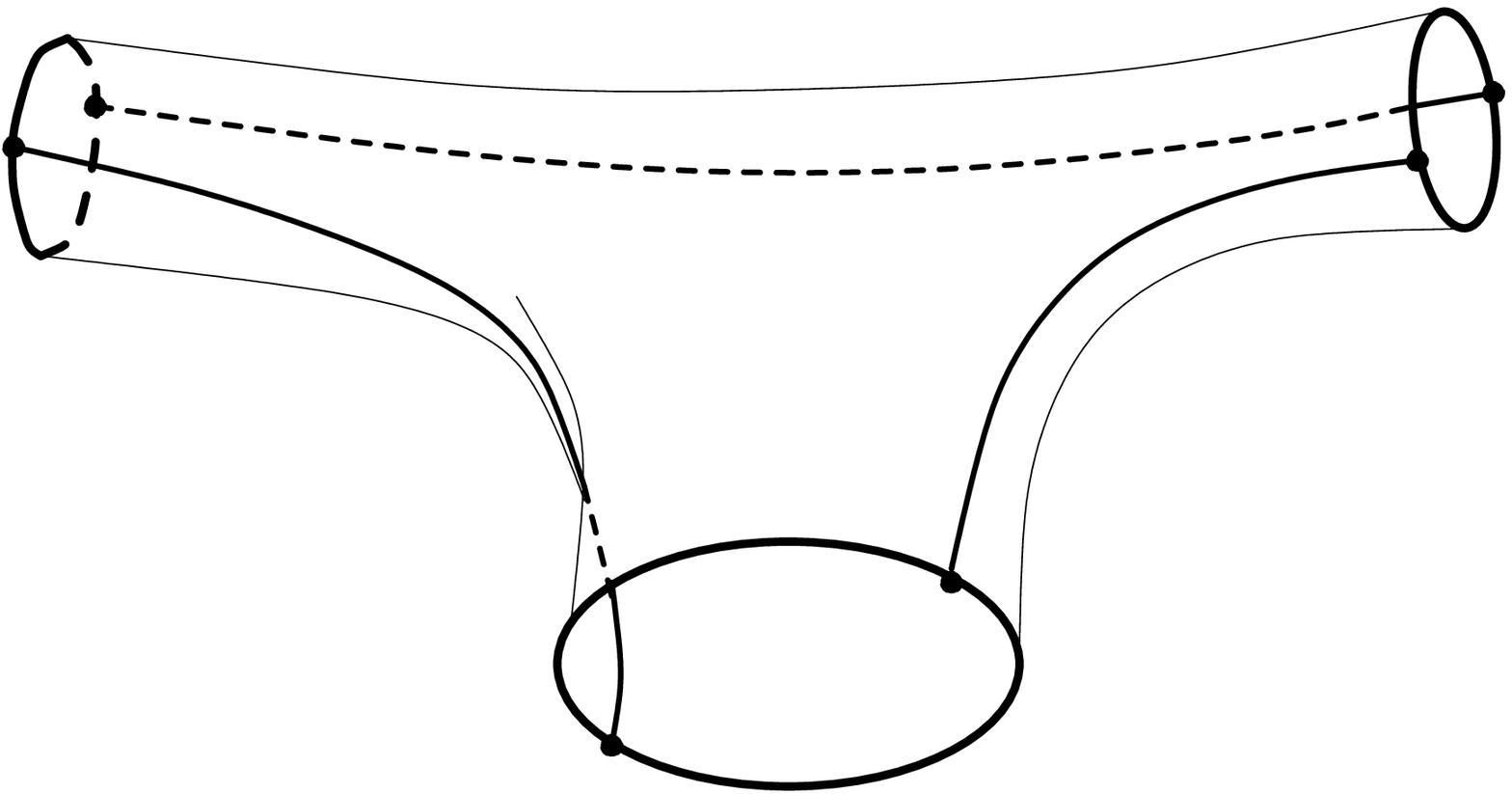}}}
\smallskip
\line{%
\vtop{\hsize=.47\hsize
\centerline{Fig.~12\.a)}
}%
\hfil
\vtop{\hsize=.47\hsize
\centerline{Fig.~12\.b)}
}%
}

\smallskip
\medskip
In the first case, Fig.~12\.a), the involution $\tau$ interchanges two
boundary components $\gamma_1$ and $\gamma_2$ of $S$ and leaves the third one
$\gamma_3$ invariant. The fixed point set $\beta$ of $\tau$ is a geodesic arc
with both ends on the $\tau$-invariant boundary component $\gamma_3$. This
case includes subcases where some boundary components of $S$ are not geodesics
but marked points (\ie punctures). In particular, if $\gamma_3$ is a marked
point, then the set $\beta$ is an (infinite) geodesic line with both ends
approaching  $\gamma_3$. The set $\beta$ divides $S$ into two parts, $S^+$
and $S^-$ (see Fig.~12\.a)), which are interchanged by $\tau$. Topologically,
each part $S^\pm$ is an annulus.

\smallskip
In the second case, Fig.~12\.b), all three boundary components $\gamma_1$,
$\gamma_2$, and $\gamma_3$ are invariant. The fixed point set of $\tau$
consists of geodesic arcs $\beta_1$, $\beta_2$, and $\beta_3$. These are the
shortest simple geodesics between $\gamma_2$ and $\gamma_3$, resp.,\ $\gamma_3$
and $\gamma_1$ and resp.,\ $\gamma_2$ and $\gamma_3$. If some boundary
component of $S$ is not a geodesic but a marked point, then corresponding
arcs have ends of infinite length approaching this boundary component. The
arcs $\beta_k$, $k=1,2,3$, divide $S$ into two parts, $S^+$ and $S^-$ (see
Fig.~12\.b)), which are interchanged by $\tau$. In this case, each part
$S^\pm$ is topologically a disc.

\smallskip
We call pieces $S^\pm$ {\sl half-pants of the first or second type},
respectively.
Note that in both cases $\tau$-invariant arcs $\beta$ or $\beta_i$ are
orthogonal to corresponding boundary circles $\gamma_j$.

\medskip
Return to the situation of a nodal curve $\barr C$ with boundary and marked
points. Let $C^d$ be the Schottky double and $\tau$ the antiholomorphic
involution. Suppose that $C^d \bs\mapo$ is non-exceptional. Use {\sl Lemma
A3.4.2} and find a $\tau$-invariant decomposition into pants $C^d = \cup_j
S_j$. Set $S_j^+ \deff S_j \cap \barr C$. Then we obtain a decomposition
$\barr C= \cup_j S^+_j$ such that the pieces $S^+_j$ are either pants (which
means $S^+_j = S_j$), or half-pants of the first or second type. This
decomposition is a suitable one for the situation of Gromov convergence up to
the boundary of curves with totally real boundary conditions. In particular, we
obtain arcs $\beta_{j,k}$ as $\tau$-fixed point sets of $S^+_j$, which define
a decomposition $\d C= \cup_{j,k} \beta_{j,k}$ of the boundary of $\barr C$.
The collection $\bfbeta' \deff \{ \beta_{j,k} \}$ of these arcs satisfies 
condition \sli of {\sl Definition 5.4}, but it can be different from the
collection $\bfbeta =\{ \beta_i \}$ which was given. The reason is that in the
construction of the pants-decomposition $C^d = \cup_j S_j$ we can subdivide
original arcs $\beta_i \in \bfbeta$ into smaller pieces, so that every arc
$\beta_i \in \bfbeta$ is a union of arcs $\beta_{j,k}$ from $\bfbeta'$. This
means compatibility of $\bfbeta$ and $\bfbeta'$.

\medskip\smallskip
The next step is to establish an analog of {\sl Theorem 5.3.2}. Assume that the
hypothesis of {\sl Theorem A3.4.1} is fulfilled. For each curve $C_n$ denote by
$C^d_n$ its Schottky double and by $\tau_n$ the corresponding involution.

\state Lemma A3.4.3. \it In the situation above, after passing to a
subsequence, there exist parameterizations $\sigma^d_n: \Sigma^d \to C^d_n$,
a finite covering $\calv$ of\/ $\Sigma^d$ by open sets $\{ V_\alpha \}$, and
a set $\{x^*_1, \ldots, x^*_m\}$ of marked points on $\Sigma$ such that the
conditions {\sl(a), (c)--(f)} of {\sl Theorem 4.2} and the following
additional conditions {\sl(b')} and {\sl(h\/)} are satisfied:

{\sl(b')} $\sigma_n\{x^*_1, \ldots, x^*_m\}$ is the set of marked points
on $C^d_n$ corresponding to the decomposition of the boundary 
$\d C_n$ into arcs
$\beta_{n,i}$; moreover, each such point $x^*_j$ lies in a single piece of
covering $V_\alpha$ which is a disc;

{\sl(h\/)} there exists an involution $\tau : \Sigma^d \to \Sigma^d$ which is
compatible with the covering $\calv$ and with parameterizations $\sigma^d_n$,
\i.e.,$\calv$ is $\tau$-invariant and $\tau_n \scirc \sigma^d_n = \sigma^d_n
\scirc \tau$. In particular, each marked point $x^*_i$ of $\Sigma^d$ is fixed
by $\tau$. \rm

\state Remark. The condition {\sl(g)} of {\sl Theorem 5.3.2} is trivial in this
case, because $C^d_n$ and $\Sigma^d$ are closed.

\state Proof. One can use the proof of {\sl Theorem 5.3.2} with minor
modifications. Note that the starting points of that proof were the intrinsic
metric on non-exceptional components of nodal curves $C_n$  and the
decomposition of $C_n$ into pants. Now the existence of a $\tau$-invariant
decomposition of the curves $C^d_n$ into pants is provided by {\sl Lemma
A3.4.2}, whereas the $\tau$-invariance of the intrinsic metrics follows from
the fact that any (anti)holomorphic isomorphism of curves with marked points
is an isometry \wrt the intrinsic metric. Thus the constructions of the proof
of {\sl Theorem 5.3.2} yield $\tau$-invariant objects. Condition {\sl(b')}
does not cause much difficulty. \qed

\medskip
Now we are ready to finish the

\state Proof of Theorem A3.4.1. 
As was mentioned, our main idea is to modify
the construction used in the proof of {\sl Theorem 4.1.1} to make it
$\tau$-invariant. The main work has already been done. We have the 
necessary a priori
estimates, the construction of a $\tau$-invariant pants-decomposition of the
double $C^d_n$ of the curve $\barr C_n$ and the appropriate covering $\calv$
of the real surface $\Sigma^d$ parameterizing the doubles $C^d_n$.

As in the proof of {\sl Theorem 4.1.1}, we consider the curves $C_{ \alpha ,n}
\deff \sigma^d_n(V_\alpha)$. Due to the presence of the involutions 
$\tau_n$, the
geometric situation in now different. This involves new phenomena and needs
additional considerations and constructions. In particular, the pieces $C_{
\alpha ,n}$ are divided into two groups depending on whether they are
disjoint from the boundary $\d C_n$ or intersect it. In the last case $C_{
\alpha ,n}$ is $\tau_n$-invariant. In this case we shall use the notation
$C^+_{ \alpha ,n} \deff C_{ \alpha,n} \cap C_n$ for the part of $C_{ \alpha
,n}$ lying in $C_n$. In addition, we denote $V^+_\alpha \deff V_\alpha \cap
\Sigma$. Then $V_\alpha$ appear as the union of domains $V^+_\alpha$ and
$\tau(V^+_\alpha)$, interchanged by $\tau$. Similarly, it is true 
for $C_{ \alpha,n}$.

\smallskip
To prove the theorem, we want to construct a refined covering $\wt\calv$ of
$\Sigma$ and refined parameterizations $\ti\sigma_n: \Sigma \to C_n$ such
that for every $V_\alpha \in \wt\calv$ the sequence $(C_{ \alpha ,n},
u_{ \alpha ,n})$ with $C_{ \alpha ,n} \deff \ti\sigma_n(V_\alpha)$
one the following convergence types holds:

\smallskip{\parindent=1.5\parindent\sl
\item{A\/$'$)} $C_{\alpha, n}$ are annuli of infinitely growing conformal
radii $l_n$ disjoint from $\d C_n$, and the conclusions of {\sl Lemma 5.2.2}
hold.

\item{A\/$''$)} $C_{\alpha, n}$ are $\tau_n$-invariant annuli of infinitely
growing conformal radii $l_n$, and the conclusions of {\sl Lemma A3.3.6} are
valid for $\Theta(0,l_n) \cong C^+_{\alpha, n} \deff C_{\alpha, n} \cap C_n$.

\item{B\/$'$)} Every $C_{\alpha, n}$ is disjoint from $\d C_n$ and isomorphic
to the standard node $\cala_0 =\Delta \cup_{ \{0\} } \Delta$ such that the
compositions $V_\alpha \buildrel \sigma _{\alpha, n} \over \lrar C_{\alpha,
n} \buildrel \cong \over \lrar \cala_0$ define the same parameterizations of
$\cala_0$ for all $n$; furthermore, the induced maps $\ti u_{\alpha, n}:
\cala_0 \to X$ strongly converge;

\item{B\/$''$)} Every $C_{\alpha, n}$ is $\tau_n$-invariant and $C^+_{\alpha,
n} \deff C_{\alpha, n} \cap C_n$ is isomorphic to the standard boundary node
$\cala^+ _0 =\Delta^+ \cup_{ \{0\} } \Delta^+$ such that for $V^+_\alpha
\deff V_\alpha \cap \Sigma$ the compositions $V^+_\alpha \buildrel \sigma^+
_{\alpha, n} \over \lrar C^+_{\alpha, n} \buildrel \cong \over \lrar
\cala^+_0$ define the same parameterizations of $\cala^+_0$ for all $n$;
furthermore, the induced maps $\ti u^+_{\alpha, n}: \cala^+_0 \to X$ strongly
converge.

\item{C)} The structures $\sigma_n^*j\vph_n \ogran_{V_\alpha}$ and the maps
$u_{\alpha, n}\scirc \sigma _{\alpha, n}: V_\alpha \to X$ strongly
converge.
}

\noindent
In the case {\sl B\/$''$)} the strong convergence of maps $\ti u^+_n: \cala^+
_0 \to X$ is the one in the $L^{1,p^*}$-topology for some $p^*>2$ {\sl up to
the boundary intervals containing the nodal point}. An equivalent requirement
is the usual $L^{1,p^*}$-convergence of the doubles $\ti u^d_n: C_{\alpha, n}
\to X$ on compact subsets of $C_{\alpha, n} \cong \cala_0$.

\smallskip
To obtain a desired refinement, we use the same inductive procedure as in the
proof of {\sl Theorem 1.1}. To insure convergence near the boundary 
$\d C_n$, we
take a new value for the constant determining the inductive step. We choose a
positive $\eps^b$ such that $\eps^b \le \eps$ and all a priori
estimates of this section are valid for maps with area\.\.$\le 3\eps^b$. This
will yield the convergence of type {\sl A)--C)} for sequences of curves with
totally real boundary conditions and with the upper bound $\eps^b$ on the area.

In fact, essential modifications of the constructions of {\sl Theorem 1.1} are
needed only if the covering piece $V_\alpha$ is $\tau$-invariant. Indeed, if
$V_\alpha$ is not $\tau$-invariant, then we can apply all the argumentations
and constructions used in {\sl Cases 1)--4)} in the proof of {\sl Theorem
1.1}, and then ``transfer'' them onto $\tau(V_\alpha)$ by means of $\tau$.
This gives the inductive step preserving $\tau$-invariance.

Hence, it remains to consider the situation when the covering piece $V_\alpha$
is $\tau$-invariant. As in {\sl Theorem 1.1}, we must consider four cases.

{\sl\smallskip
Case 1$_b$): $C_{\alpha, n}$ have constant complex structure different from
the one of the standard node.

\smallskip
Case 2$_b$): $C_{\alpha, n}$ are annuli of changing conformal radii $R_n$
such that $R_n \to R <\infty$.

\smallskip
Case 3$_b$): $C_{\alpha, n}$ are isomorphic to the standard node, so that
$C^+_{\alpha, n}$ are isomorphic to the standard boundary node $\cala_0^+$.

\smallskip
Case 4$_b$): $C_{\alpha, n}$ are annuli of infinitely growing conformal radii
$R_n$.
}

The subindex $(\cdot)_b$ indicates that we will consider the cases where 
$V_\alpha$
intersects the {\it b}oundary of $\Sigma$. As was mentioned, the last
property is equivalent to the fact that $V_\alpha$ is $\tau$-invariant. 
References to {\sl Cases 1)--4) \it without} the subindex will mean 
the corresponding parts of the proof of {\sl Theorem 1.1}.

\smallskip\noindent
{\sl Case 1$_b$)}. Without loss of generality we may assume that $V_\alpha$
is a domain with a fixed complex structure and a fixed antiholomorphic
involution $\tau$, and that $u_{\alpha, n}: V_\alpha \to X$ is a sequence of
$\tau$-invariant maps which are  (anti)holomorphic outside the set of $\tau
$-invariant points of $V_\alpha$. If we have the convergence of type {\sl C)},
there is nothing to do. Otherwise, we fix a $\tau$-invariant metric on
$V_\alpha$ compatible with the complex structure. Repeating the constructions
from {\sl Case 1)}, we distinguish the ``bubbling'' points
$y^*_1, \ldots, y^*_l$ where the strong convergence fails.

Take the first point $y^*_1$. Suppose $y^*_1$ is disjoint from
$\d\Sigma$. Then we may assume that $y^*_1\in V^+_\alpha$. Thus, we can repeat
the rest of the constructions from {\sl Case 1)}. The only correction needed
at this point is that the neighborhood $\Delta( y^*_1, \varrho)$ of $y^*_1$
must be small enough and lie in $V^+_\alpha$. Transferring all these
constructions into $\tau(V_\alpha)$, we realize the inductive step  
preserving the $\tau$-invariance.

It remains to consider the case, where $y^*_1 \in \d\Sigma$. This means that
$y^*_1$ is $\tau$-invariant. Let $z$ be a holomorphic coordinate in a
neigborhood of $y^*_1$ on $V_\alpha$ such that $z=0$ in $y^*_1$, the
involution $\tau$ corresponds to the conjugation $z \mapsto \bar z$ and $\im
z >0$ in $\Sigma$. Find the sequences $r_n\lrar0$ of the radii and $x_n \to
y^*_1$, using the constructions from {\sl Case 1)}. Note that the sequence
$\tau(x_n)$ has the same property. Thus, replacing some points $x_n$ by
$\tau(x_n)$, we may additionally assume that all $x_n$ lie in $\barr
V^+_\alpha$. Let $v_n:\Delta (0,{\varrho \over 2r_n}) \to (X,J_n)$ be the
rescalings of maps $u_n$ defined by $v_n(z) \deff u_n(x_n + {z\over r_n})$.
Argumentations of {\sl Case 1} show that there exists the limit $v_\infty:
\cc \to X$ of (a subsequence of) $\{ v_n \}$ which extends to a map $v_\infty
: S^2 \to X$.

Denote by $\rho_n$ the distance from $x_n$ to $\d\Sigma$ and by $\ti x_n$ the
point on $\d\Sigma$ closest to $x_n$. Then $x_n = \ti x_n + \isl \rho_n$ in the
coordinate $z$ introduced above. In addition, $\lim \ti x_n = y^*_1$. 
We consider two subcases according to the possible behavior of 
$\rho_n$ and $r_n$.

\smallskip\noindent
{\sl Subcase 1$'_b$): $\bigl\{{ \rho_n \over r_n }\bigr\}$ is bounded}.
Passing to a subsequence, we may assume that ${\rho_n \over r_n}$ converges.
Fix an upper bound $b$ for the sequence ${\rho_n \over r_n}$. In particular,
$b \ge \lim {\rho_n \over r_n}$.

For $n>\!>1$ define maps $v_n:\Delta (0,{\varrho \over 2r_n}-b) \to (X,J_n)$
and $\ti v_n:\Delta (0,{\varrho \over 2r_n}-b) \to (X,J_n)$ setting $v_n(z)
\deff u_n (x_n +r_n z)$ and $\ti v_n(z)\deff u_n(\ti x_n +r_n z)$,
respectively. Then every $\ti v_n$ is the shift of the map $v_n$ by
$\isl{\rho_n\over r_n}$, \i.e., $\ti v_n(z)= v_n\bigl(z + \isl {\rho_n\over
r_n}\bigr)$. The arguments of {\sl Case 1)} show that $v_n$ converge on
compact subsets of $\cc$ to a non-constant map. Consequently, $\ti v_n$ also
converge on compact subsets of $\cc$ to a non-constant map $\ti v_\infty: \cc
\to X$. Moreover, since $\area( \ti v_\infty( \cc))$ is finite, $\ti v_\infty$
extends to a map $\ti v_\infty: S^2 \to X$. By the choice of $\eps^b$,
$\area(\ti v_\infty (S^2)) \ge 3\eps^b$. Changing the choice of the constant
$b$, we can additionally assume that $\area(\ti v_\infty( \Delta(0,b)) \ge
2\eps^b$. Then for all sufficiently big $n$ we obtain
$$
\area(\ti v_n(\Delta(0,b))  \ge \eps^b.
\eqno(A3.4.1)
$$
%%%%%%%%%%%%%%%%%%%%%question%%%%%%%%%%%%%%%%%

For $n>\!>1$ we define the coverings of $V_\alpha$ by three sets
$$
V^{(n)}_{\alpha,1} \deff  V_\alpha
\bs \barr\Delta (0, {\textstyle{\varrho\over2}}),
\qquad
V^{(n)}_{\alpha,2} \deff \Delta (0, \varrho)
\bs \barr \Delta (\ti x_n, br_n),
\qquad
V^{(n)}_{\alpha,3} \deff \Delta (\ti x_n, 2br_n).
$$
Fix $n_0$ sufficiently big. Denote $V_{\alpha,1} \deff V^{(n_0)}_{\alpha,1}$,
$V_{\alpha,2} \deff V^{(n_0)}_{\alpha,2}$ and $V_{\alpha,3} \deff V^{(n_0)}
_{\alpha,3}$. There exist diffeomorphisms $\psi_n: V_1 \to V_1$ such that
$\psi_n: V_{\alpha,1} \to V^{(n)}_{\alpha,1}$ is an identity, $\psi_n:
V_{\alpha, 2} \to V^{(n)}_{\alpha,2}$ is a diffeomorphism, and $\psi_n:
V_{\alpha,3} \to V^{(n)}_{\alpha,3}$ is biholomorphic \wrt the complex
structures, induced from $C_n$ by means of $\sigma^d_n$. Note that the sets
$V^{(n)}_{\alpha,i}$ are $\tau$-invariant. Moreover, we can choose the maps
$\psi_n$ in such a way that $\psi_n$ are also $\tau$-invariant.

The covering $\{ V_{\alpha,1}, V_{\alpha,2}, V_{\alpha,3} \}$ of $V_1$ and
parameterizations $\ti\sigma_n \deff \sigma_{\alpha,n} \scirc \psi_n: V_1 \to
C_{\alpha, n}$ satisfy the conditions of {\sl Lemma A3.4.3}. Moreover,
inequality (A3.2.1) implies $\area(u_n(\ti\sigma_n(V_{\alpha,i}))) \le
%%%%%%%%%%%%%%%%%%%%%%%question%%%%%%%%%%%%%%%%%%%%%
(N-1)\eps^b$. Consequently, we can apply the inductive assumptions for the
sequence of curves $\ti\sigma_n(V_{\alpha,i})$ and finish the proof by
induction.

\smallskip\noindent
{\sl Subcase 1$''_b$): $\bigl\{{ \rho_n \over r_n}\bigr\}$ is unbounded}.
Passing to a subsequence, we may assume that ${\rho_n \over r_n}$ increases
infinitely. However, $\rho_n \lrar 0$ since $x_n \lrar y^*_1 \in \d\Sigma$.

Define maps $v_n:\Delta (0,{\varrho \over 2r_n}) \to (X,J_n)$, setting $v_n(z)
\deff u_n(x_n +r_n z)$. As in {\sl Case~1)}, $v_n$ converge on compact
subsets of $\cc$ to a non-constant map $v_\infty: \cc \to X$, which extends
to a map from the whole sphere $S^2$. Choose $b>0$ satisfying (A3.2.1).
%%%%%%%%%%%%%%%%%question%%%%%%%%%%%%%%%%%%%
For $n>\!>1$ we define the coverings of $V_\alpha$ by five sets
$$\mathsurround=0pt
\matrix\format\l\ \ &\l\\
V^{(n)}_{\alpha,1} \deff  V_\alpha
\bs \barr\Delta (0, {\textstyle{\varrho\over2}}),
&
V^{(n)}_{\alpha,2} \deff \Delta(0, \varrho) \bs
\barr \Delta(0, 2\rho_n)
\cr
\noalign{\vskip5pt}
\rlap{$
V^{(n)}_{\alpha,3} \deff \Delta(0, 4\rho_n) \bs
\bigl(\barr \Delta (x_n, br_n) \cap  \barr \Delta (\tau(x_n), br_n) \bigr)
$}&
\cr
\noalign{\vskip5pt}
V^{(n)}_{\alpha,4} \deff \Delta (x_n, 2br_n),
&
V^{(n)}_{\alpha,5} \deff \Delta(\tau(x_n), 2br_n).
\endmatrix
$$
Fix $n_0$ sufficiently big. Denote $V_{\alpha,i} \deff V^{(n_0)}_{\alpha,i}$,
$i=1,\ldots,5$. Then for every $n>\nobreak\!>1$ there exists a diffeomorphism
$\psi_n: V_1 \to V_1$ with the following properties:

\sli $\psi_n$ maps $V_{\alpha,i}$ onto $V^{(n)}_{\alpha,i}$ diffeomorphically;

\slii $\psi_n: V^{(n)}_{\alpha,1} \to V^{(n)}_{\alpha,1}$ is the identity;

\sliii $\psi_n: V_{\alpha, 2} \to V^{(n)}_{\alpha,2}$ and $\psi_n:
V_{\alpha,3} \to V^{(n)}_{\alpha,3}$ are diffeomorphisms;

\sliv $\psi_n: V_{\alpha,3} \to V^{(n)}_{\alpha,3}$ is biholomorphic \wrt
the complex structures, induced from $C_n$ by means of $\sigma^d_n$; and,
finally

\slv $\psi_n$ are $\tau$-invariant: $\tau \scirc\psi_n =\psi_n \scirc\tau$.
\newline\noindent
Note that the last property is obtained due to the fact that the sets
$V^{(n)} _{\alpha, i}$ are $\tau$-invariant. The remaining constructions 
are the same as in {\sl Subcase 1$'_b$)}.

\medskip\noindent
{\sl Case 2$_b$)}. Consider the parameterizations $\sigma_n: V_\alpha \to
C_{\alpha,n}$. Without loss of generality we may assume that the complex
structures $\sigma_n^* j_n\ogran_{V_\alpha}$ are constant near the boundary 
$dV_\alpha$ and converge to some complex structure. If we have the convergence
of type {\sl C)}, i.e., the strong convergence, there is nothing to do.
Otherwise, there exists only a finite set of points $\{y^*_1,\ldots, y^*_l\}$
where the strong convergence fails. Changing the parameterizations $\sigma_n$,
we may additionally assume that the structures $\sigma_n^* j_n \ogran
_{V_\alpha}$ are constant in the neighborhood of these points. Then we repeat
the argumentations of {\sl Case 1$_b$)}.

\medskip\noindent
{\sl Case 3$_b$)}. Fix identifications $C_{\alpha,n} \cong \cala_0$ such that
every $C_{\alpha,n}^+$ is mapped onto $\cala_0^+$ and the induced
parameterization maps $\sigma_{\alpha, n} : V_\alpha \to \cala_0$ are the same
for all $n$ and $\tau$-invariant. Fix the standard representation of $\cala_0$
as the union of two discs $\Delta'$ and $\Delta''$ with identification of the
centers $0\in \Delta'$ and $0\in \Delta''$ into the nodal point of $\cala_0$,
still denoted by $0$. Let $\Delta' (x, r)$ denote the subdisc of $\Delta'$
with the center $x$ and the radius $r$.

Denote by $u'_n :\Delta' \to X$ and $u''_n :\Delta'' \to X$ the corresponding
``components'' of the maps $u_{\alpha,n} : C_{\alpha,n} \to X$. Find the
common collection of bubbling points $y^*_i$ for both sequences of maps $u'_n
:\Delta' \to X$ and $u''_n :\Delta'' \to X$. If there are no bubbling points,
then we obtain a convergence of type {\sl B)}, and the proof can be finished by
induction. Otherwise, consider the first such point $y^*_1$, which lies, say,
on $\Delta'$. If $y^*_1$ is distinct from the nodal point $0 \in \Delta'$,
then we simply repeat all the constructions in {\sl Case 1$_b$)}.

It remains to consider the case $y^*_1=0 \in \Delta'$. The following
modifications of the argumentations are needed. Repeat the construction of
the radii $r_n\lrar 0$ and the points $x_n \lrar y^*_1=0$ from {\sl Case
1$_b$}. Then $\{x_n \}$ is a sequence in the half-disk $\delta^{\prime+}
\deff \{ z\in \Delta' : \im z \ge0 \}$. Set $\ti x_n \deff \re(x_n)$, $\rho_n
\deff \im(x_n)$ and $R_n \deff |x_n|$. Thus, $x_n= \ti x_n +\isl\rho_n$, $R_n$
is the distance from $x_n$ to the point $0=y^*_1 \in \Delta'$, whereas
$\rho_n$ is the distance from $x_n$ to the interval $]-1,1[ \subset \Delta'$,
the set $\tau$-invariant points of $\Delta'$. Thus $\rho_n \le R_n$. Fix
$\varrho>0$ such that the disc $\Delta'(0,\varrho)$ contains no bubble points
$y^*_i\not=0 \in \Delta'$.

Depending on the behavior of the sequences $r_n$, $\rho_n$ and $R_n$, we
consider four subcases.
%%%%%%%%%%%%%%%%%%%%page 100%%%%%%%%%%%%%%%%%%%%%%%%%
\medskip\noindent
{\sl Subcase 3\/$'_b$): The sequence $\bigl\{{ R_n \over r_n }\bigr\}$ is
bounded.} Then the sequences $\bigl\{{ {\rho_n\over r_n} }\bigr\}$ and $\bigl
\{ {\ti x_n\over r_n} \bigr\}$ are also bounded. Passing to a subsequence, we
may assume that the corresponding limits exist. Let $b$ be some upper bound
for the sequence $\bigl\{{ R_n \over r_n }\bigr\}$. Consider the maps $\ti
v_n: \Delta(0, {\varrho \over 2r_n}-b) \to X$ defined by $\ti v_n(z) \deff
u_n(\ti x_n+ {z\over r_n})$. Then $\ti v_n$ are $\tau$-invariant, $\ti v_n
\scirc \tau= \ti v_n$, $\ti v_n$ converge to a nonconstant map $\ti v_\infty:
\cc \to X$ on compact subsets of $\cc$, and $\ti v_\infty$ extends to a map
$\ti v_\infty: S^2 \to X$.

Since  $\ti v_\infty$ is nonconstant, $\norm{ d\ti v _\infty} ^2 _{L^2( S^2)}
= \area( \ti v_\infty( S^2) )\ge 3\eps_b$. Choose $b>0$ in such a way that
$$
\norm{d\ti v_\infty}_{L^2(\Delta (0,b))}^2 \ge 2\eps_b
\eqno(A3.4.2)
$$
and $b \ge 2\lim {R_n \over r_n}+2$. Due to {\sl Corollary A3.2.2}, for 
$n>\!>1$ we obtain the estimate
$$
\norm{du'_n}_{L^2(\Delta' (\ti x_n, br_n))}^2 =
\norm{ d\ti v_n}_{L^2(\Delta (0,b))}^2 \ge \eps_b.
\eqno(A3.4.3)
$$
Note that $0\in \Delta' (\ti x_n, (b-1)r_n))$ for $n>\!>1$ by the choice of
$b$.

Define the coverings of $\cala_0$ by four sets
$$\mathsurround=0pt
\matrix\format\l\ \ &\l\\
W^{(n)}_1 \deff  \Delta'
\bs \barr\Delta' (0, {\textstyle{\varrho\over2}}),
&
W^{(n)}_2 \deff \Delta' (0, \varrho)
\bs \barr \Delta' (\ti x_n, br_n),
\cr
\noalign{\vskip4.5pt}
W^{(n)}_3 \deff \Delta' (\ti x_n, 2br_n)
\bs \barr\Delta'(0,{\textstyle {r_n \over 2}}),
&
W^{(n)}_4 \deff \Delta' (0, r_n) \cup
\Delta'',
\endmatrix
$$
and lift them to $V_\alpha$ by putting $V^{(n)}_{\alpha,i} \deff \sigma\inv
_{\alpha,n}(W^{(n)}_i)$. Choose $n_0 >\!> 0$ such that $|x_n| < (b-1)r_n$
and the relation (A3.4.3) holds for all $n \ge n_0$. Set $V_{\alpha,i} \deff
%%%%%%%%%%%%%%%%%%%question%%%%%%%%%%%%%%%%%%%%%%%%%
V^{( n_0)} _{\alpha, i}$. Fix diffeomorphisms $\psi_n: V_\alpha \to V_\alpha$
such that $\psi_n: V_{\alpha,1} \to V^{(n)} _{\alpha,1}$ is the identity map,
$\psi_n: V_{\alpha,2} \to V^{(n)}_{\alpha,2}$ and $\psi_n: V_{\alpha,3} \to
V^{(n)} _{\alpha,3}$ are diffeomorphisms, and $\psi_n: V_{\alpha,4} \to
V^{(n)} _{\alpha,4}$ correspond to isomorphisms of nodes $W^{(n)}_4 \cong
\cala_0$. Set $\sigma'_n \deff \sigma_n \scirc \psi_n$. The choice above can
be made in such a way that the refined covering $\{ V_{\alpha,i} \}$ of
$V_\alpha$ and parameterization maps $\sigma'_n: V_\alpha \to C_{\alpha,n}$
have the properties of {\sl Lemma A3.4.3}. Relation (A3.4.3) implies the
%%%%%%%%%%%%%%%%%%%%question%%%%%%%%%%%%%%%%%%%%%%
estimate $\area(u_n (\sigma'_n(V_{\alpha,i})) \le (N-1)\,\eps$. This provides
the inductive conclusion for {\sl Subcase 3\/$'_b$)}.

\medskip\noindent
{\sl Subcase 3\/$''_b$): The sequence $\bigl\{{R_n \over r_n }\bigr\}$
increases infinitely but $\bigl\{{\rho_n \over r_n }\bigr\}$ remains
bounded}.  Note that in this subcase we still have the relation $R_n \lrar
0$, or equivalently, $x_n \lrar 0$. On the other hand, $\lim {\rho_n \over
R_n} =0$.  This implies that for $\ti R_n \deff |\ti x_n|$ we have $\lim
{\ti R_n \over R_n} =1$ since $R_n^2 = \ti R_n^2 + \rho_n^2$.

We proceed as follows. Define the maps $\ti v_n: \Delta(0, {\varrho \over
2r_n}-b) \to X$ setting $\ti v_n(z) \deff u'_n(\ti x_n + {z\over r_n})$. Then
the $\ti v_n$ have the same properties as in {\sl Subcase 3\/$'_b$)}. 
Choose $b>0$
obeying the relation (A3.4.2). Then for $n>\!>0$ the property (A3.4.3) follows.

%%%%%%%%%%%%%%%%%2 questions%%%%%%%%%%%%%%%%%%%%

For $n>\!>0$ define the coverings of $\cala_0$ by six sets
$$\mathsurround=0pt
\matrix\format\l\ \ &\l\\
W^{(n)}_1 \deff  \Delta'
\bs \barr\Delta' (0, {\textstyle{\varrho\over2}}),
&
W^{(n)}_2 \deff \Delta' (0, \varrho)
\bs \barr \Delta' (\ti x_n, 2\ti R_n),
\cr\noalign{\vskip5pt}
W^{(n)}_3 \deff  \Delta' (\ti x_n, 4\ti R_n)
\bs \bigl(\barr \Delta' (\ti x_n, {\ti R_n \over 6})
\cup \barr \Delta' (0, {\ti R_n \over 6}) \bigr)
&
W^{(n)}_4 \deff \Delta' (0, {\ti R_n \over 3}) \cup \Delta'',
\cr\noalign{\vskip5pt}
W^{(n)}_5 \deff \Delta' (\ti x_n, {\textstyle{\ti R_n \over 3}})
\bs \barr\Delta'(\ti x_n, br_n),
&
W^{(n)}_6 \deff \Delta' (0, 2br_n ),
\endmatrix
$$
and lift them to $V_\alpha$ by putting $V^{(n)}_{\alpha,i} \deff \sigma\inv
_{\alpha,n}(W^{(n)}_i)$. Choose $n_0 >\!> 0$ such that $R_{n_0} >\!> b
r_{n_0} $, and set $V_{\alpha,i} \deff V^{(n_0)} _{\alpha, i}$. Choose
diffeomorphisms $\psi_n: V_\alpha \to V_\alpha$ such that $\psi_n: V
_{\alpha, 1} \to V^{(n)} _{\alpha,1}$ is the identity map, $\psi_n: V
_{\alpha, 2} \to V^{(n)}_{\alpha,2}$, $\psi_n: V_{\alpha,4} \to V^{(n)}
_{\alpha,4}$ and $\psi_n: V_{\alpha,5} \to V^{(n)}_{\alpha,5}$ are
diffeomorphisms, and finally, $\psi_n: V_{\alpha,6} \to V^{(n)}_{\alpha,6}$
corresponds to isomorphisms of nodes $W^{(n)}_6 \cong \cala_0$. Set
$\sigma'_n \deff \sigma_n \scirc \psi_n$. Note that the choices can be made
in such a way that $\{ V_{\alpha,i} \}$ and parameterization maps $\sigma'_n:
V_\alpha \to C_{\alpha,n}$ have the properties of {\sl Lemma A3.4.3}. 
As above,
we get the estimate $\area(u_n (\sigma'_n(V_{\alpha,i} )) \le (N-1)\, \eps$
due to (A3.4.3). Thus we obtain the inductive conclusion for {\sl Subcase
3\/$''_b$)} and can proceed further.

%%%%%%%%%%%%%%%%%%question%%%%%%%%%%%%%%%%%%%

\smallskip\noindent
{\sl Subcase 3\/$'''_b$): The sequence $\bigl\{{\rho_n \over r_n} \bigr\}$
increases infinitely, but $\bigl\{{ R_n \over \rho_n} \bigr\}$ remains
bounded}. Then $\bigl\{{R_n \over r_n} \bigr\}$ also increases infinitely, but
both sequences $\{R_n \}$ and $\bigl\{{ \rho_n }\bigr\}$ converge to 0.
We may also assume
that $\bigl\{{ \rho_n \over R_n }\bigr\}$ and $\bigl\{{ \ti x_n \over R_n }
\bigr\}$ also converge. Set $a_1 \deff \lim {\ti x_n \over R_n}$, $a_2 \deff
\lim{ \rho_n \over R_n}$, $a\deff a_1 + \isl a_2$ and $\barr a \deff a_1
- \isl a_2$. Note that $0 < a_2 \le 1$ and the involutions $\tau_n$ in
$C_{\alpha, n}$ correspond to the complex conjugation $z\to \barr z$ in
$\Delta'$. In particular, $\barr x_n = \tau_n(x_n)$.

Consider maps $v_n: \Delta(0, {\varrho \over 2r_n}) \to X$ defined by $v_n(z)
\deff u'_n(x_n + {z\over r_n})$. Then the sequence $\{ v_n \}$ converges on
compact subsets to a nonconstant map which extends to the map $v_\infty :S^2
\to X$. Moreover, we can fix sufficiently big $b>0$ such that for $n>\!>0$ we
get the property (A3.4.3).

%%%%%%%%%%%%%%%%%question%%%%%%%%%%%%%%%%

For $n>\!>0$ define the coverings of $\cala_0$ by eight sets
$$\mathsurround=0pt
\matrix\format\l\ \ \ \ \ &\l\\
W^{(n)}_1 \deff  \Delta'
\bs \barr\Delta' (0, {\textstyle{\varrho\over2}}),
&
W^{(n)}_2 \deff \Delta' (0, \varrho)
\bs \barr \Delta' (0, 2R_n),
\cr\noalign{\vskip6pt}
\rlap{$W^{(n)}_3 \deff  \Delta' (0, 4R_n)
\bs \bigl(\, \barr \Delta' ( a\, R_n , {a_2 R_n \over 4})
\cup \barr \Delta' ( \barr a\, R_n , { a_2 R_n \over 4})
\cup \barr \Delta' (0, {a_2 R_n \over 4})
\,\bigr)\qquad$}
&
\cr\noalign{\vskip6pt}
W^{(n)}_4 \deff \Delta' (0, {a_2 R_n \over 3}) \cup \Delta''.
\cr\noalign{\vskip5pt}
W^{(n)}_5 \deff \Delta' (a R_n, {a_2 R_n \over 3})
\bs \barr\Delta'(x_n, br_n),
&
W^{(n)}_6 \deff \Delta' (x_n, 2br_n ),
\cr\noalign{\vskip5pt}
W^{(n)}_7 \deff \Delta' (\barr a R_n, {a_2 R_n \over 3})
\bs \barr\Delta'(\barr x_n, br_n),
&
W^{(n)}_8 \deff \Delta' (\barr x_n, 2br_n ),
\endmatrix
$$
and lift them to $V_\alpha$ by putting $V^{(n)}_{\alpha,i} \deff \sigma\inv
_{\alpha,n}(W^{(n)}_i)$. Fix sufficiently big $n_0 >\!> 0$, 

%%%%%%%%%%%%%%%%%question%%%%%%%%%%%%%%%%%%

and set
$V_{\alpha, i} \deff V^{(n_0)} _{\alpha, i}$. Choose diffeomorphisms $\psi_n:
V_\alpha \to V_\alpha$ mapping $V _{\alpha, i}$ diffeomorphically onto
$V^{(n)} _{\alpha, i}$ such that the assertions of {\sl Lemma A3.4.3} are
fulfilled. As above, we obtain the estimate 
$\area(u_n (\sigma'_n(V_{\alpha,i} )) \le (N-1)\, \eps$. 
This gives the inductive conclusion for {\sl Subcase 3\/$'''_b$)}.

\smallskip\noindent
{\sl Subcase 3\/$''''_b$): The sequences $\bigl\{{\rho_n \over r_n} \bigr\}$
and  $\bigl\{{ R_n \over \rho_n} \bigr\}$ increase infinitely}.
Thus $\lim {\ti R_n \over R_n} =1$.
We consider the sequence of maps $\{ v_n \}$. It is defined in the same way
as in the previous subcase and has the same properties. In particular,
$\{ v_n \}$ converges to the map $v_\infty :S^2 \to X$ and there exists
a sufficiently big $b>0$ such that for $n>\!>0$ we get the property (A3.4.3).

%%%%%%%%%%%%%%%question%%%%%%%%%%%%%%

For $n>\!>0$ define the coverings of $\cala_0$ by ten sets
$$\mathsurround=0pt
\matrix\format\l\ \ \ \ \ &\l\\
W^{(n)}_1 \deff  \Delta'
\bs \barr\Delta' (0, {\textstyle{\varrho\over2}}),
&
W^{(n)}_2 \deff \Delta' (0, \varrho)
\bs \barr \Delta' (0, 2R_n),
\cr\noalign{\vskip6pt}
\rlap{$W^{(n)}_3 \deff  \Delta' (0, 4R_n)
\bs \bigl(\, \barr \Delta' ( 0, {\ti R_n \over 4})
\cup \barr \Delta' ( \ti x_n, {\ti R_n \over 4})
\bigr),$}
&
\cr\noalign{\vskip6pt}
W^{(n)}_4 \deff \Delta' (0, {\ti R_n \over 3}) \cup \Delta'',
&
W^{(n)}_5 \deff \Delta' (\ti x_n, {\ti R_n \over 3})
\bs \barr\Delta'(\ti x_n, 2\rho_n),
\cr\noalign{\vskip6pt}
\rlap{$W^{(n)}_6 \deff \Delta' (\ti x_n, 4\rho_n )
\bs \bigl(\, \barr \Delta' ( x_n, {\rho_n \over 4})
\cup \barr \Delta' ( \barr x_n, {\rho_n \over 4})
\bigr), $}
\cr\noalign{\vskip5pt}
W^{(n)}_7 \deff \Delta' (x_n, {\rho_n \over 3})
\bs \barr\Delta'(x_n, br_n),
&
W^{(n)}_8 \deff \Delta' (x_n, 2br_n ),
\cr\noalign{\vskip5pt}
W^{(n)}_9 \deff \Delta' (\barr x_n, {\rho_n \over 3})
\bs \barr\Delta'(\barr x_n, br_n),
&
W^{(n)}_{10} \deff \Delta' (\barr x_n, 2br_n ).
\endmatrix
$$
The remaining ``manipulations'' with $W^{(n)}_i$ are the 
same as in the previous
subcases. As a result, we obtain the covering of $V_\alpha$ by sets $V_{\alpha,
i} \deff \sigma\inv _{\alpha, n_0}(W^{(n_0)}_i)$ with an appropriate $n_0 >\!>
0$ and refined parameterizations $\sigma'_n: V_\alpha \to C_{\alpha, n}$, for
which the assertions of {\sl Lemma A3.4.3} are fulfilled.
As above, we obtain the estimate $\area(u_n (\sigma'_n(V_{\alpha,i} ))
\le (N-1)\, \eps$. This gives the inductive conclusion for {\sl Subcase
3\/$''''_b$)}.

\smallskip\noindent
{\sl Case 4$_b$): $V_\alpha$ is a cylinder such that conformal radii
of $(V_\alpha, \sigma_n^*j_n)$ increase infinitely}. We can simply repeat the
contructions made in {\sl Case 4)} from the proof of {\sl Theorem 1.1}.
Additional attention is  needed to preserve the $\tau$-invariantness.

\smallskip
The proof of the theorem  can now be finished by induction.
\qed

\medskip
\state Remark. Here we give some explanation of the geometric meaning of
the constructions of the proof of {\sl Theorem 5.9} and describe the picture
of the bubbling. We restrict ourselves to {\sl Case 3$_b$)} as the most
complicated one; the constructions of the other cases can be treated similarly.
The reflection principle allows us to reduce  {\sl Case 3$_b$)} to the 
consideration of $\tau$-invariant maps $u^d_n: \cala_0 \to X$ from the standard
node which are $J_n$-holomorphic on $\cala^+_0$. The situation is
different from those in {\sl Theorem 1.1}, where the bubbling appears
in the nodal point. In this case we must take into consideration not only
the parameters $r_n$ describing the size of energy localization of the bubbled
sphere, but also the additional parameters $R_n$ and $\rho_n$. These describe
the position of the localization centers $x_n$ \wrt the nodal point and the
set of $\tau$-invariant points of $\cala_0$. Depending on the behavior of $r_n$
$\rho_n$ and $R_n$, we can have four different types of bubbling and
corrsponding {\sl Subcases 3\/$'_b$)--3\/$''''_b$)}.

\smallskip
\line{%
\vtop{\hsize=.52\hsize
In {\sl Subcase 3\/$'_b$)} the bubbling takes place in the nodal point, so that
the nodal point remains on the bubbled sphere (see region $W^{(n)}_3$ in
Fig.~13). Furthermore, the bubbled sphere contains another nodal point.
This one  appears in the limit of long cylinders $W^{(n)}_2$. Note $W^{(n)}_2$
can either strongly converge to a boundary node or have additional bubblings.
}
\hfil
\vtop{\hsize=.43\hsize\xsize=\hsize\nolineskip\rm
\putm[-.03][.17]{W^{(n)}_1}%
\putm[.23][.09]{W^{(n)}_2}%
\putm[.53][.17]{W^{(n)}_3}%
\putm[.85][.14]{W^{(n)}_4}%
\noindent
\epsfxsize=\xsize\epsfbox{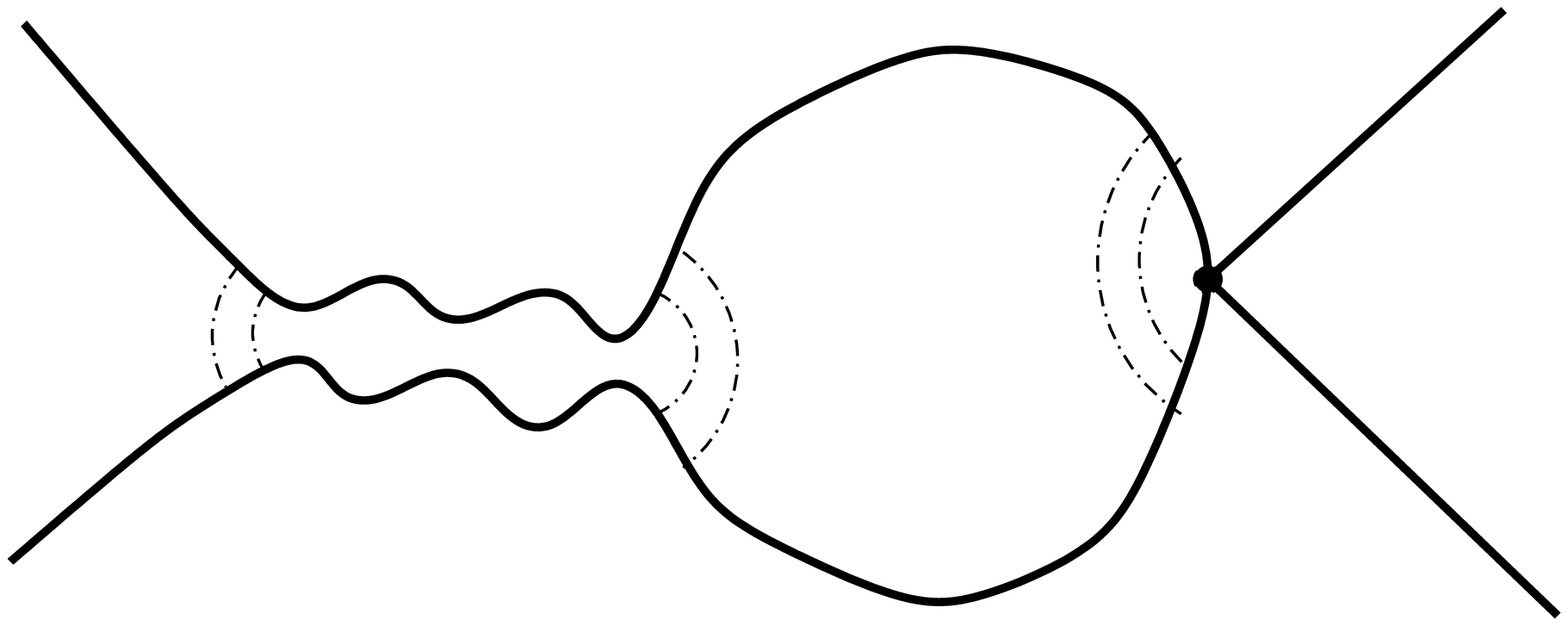}
\smallskip\smallskip
\putt[.0][.0]{
\centerline{Fig.~13. Bubbling in {\sl Subcase 3\/$'_b$)}.}
}
}}

\bgroup\baselineskip=12.3pt plus 1.5pt
Turning back from a ``doubled'' description by $\tau$-invariant objects to the
original  maps $u_n: \cala^+_0 \to X$ with a totally real boundary condition,
we
obtain the following picture. Since every covering piece $W^{(n)}_i$ is
$\tau$-invariant, for $\cala_0^+$ we obtain the covering piece $W^{(n)\,+}_i
\deff W^{(n)}_i \cap \cala_0^+$. Thus we obtain a bubbled disk represented by
$W^{(n)\,+}_3$ instead of the bubbled sphere represented by $W^{(n)}_3$,
the sequence of long strips $W^{(n)\,+}_2$ instead of the  sequence of long
cylinder $W^{(n)\,+}_2$ and so on.

\medskip
\line{%
\vtop{\hsize=.533\hsize \baselineskip=13.2pt
In {\sl Subcase 3\/$''_b$)} the bubbling happens at the boundary but away
from the nodal point. In the limit we obtain two bubbled spheres. The first one
is the limit of the re\-scaled maps $v_n$ (region $W^{(n)}_6$ in the Fig.~14).
The appearance of the second sphere can be explained as follows. The part of
the node $\cala_0$ between the first bubbled sphere $W^{(n)}_6$ and the
``constant part'' $W^{(n)}_1$ of the node is a long cylinder, represented by
pieces $W^{(n)}_2$, $W^{(n)}_3$ and $W^{(n)}_5$.
}
\hfil \vtop{\hsize=.44\hsize\xsize=\hsize\nolineskip\rm
\putm[-.02][.40]{W^{(n)}_1}%
\putm[.21][.31]{W^{(n)}_2}%
\putm[.53][.43]{W^{(n)}_3}%
\putm[.50][.06]{W^{(n)}_6}%
\putm[.63][.22]{W^{(n)}_5}%
\putm[.85][.415]{W^{(n)}_4}%
\noindent
\epsfxsize=\xsize\epsfbox{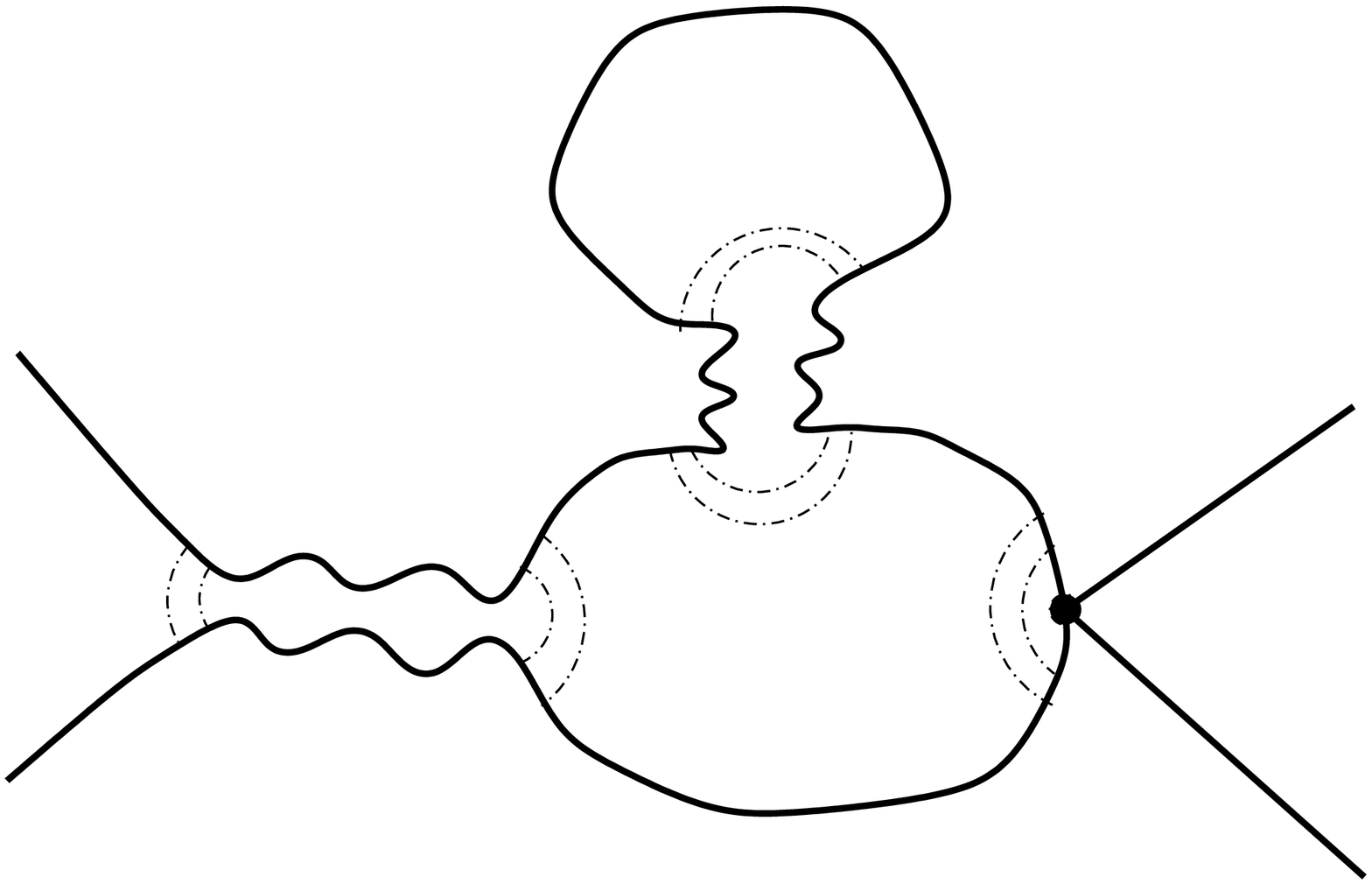} 
\smallskip\smallskip 
\putt[.0][.0]{
\centerline{Fig.~14. Bubbling in {\sl Subcase 3\/$''_b$)}.}
}
}}

However, because of the presence of the nodal point (piece $W^{(n)}_4$ on
the figure), this ``part inbetween'' is topologically not a cylinder (\ie an
annulus) but pants. Furthermore, the complex structures on the pants are not
constant. To get pants with a constant structure (piece $W^{(n)}_3$), we 
cut off
the annuli $W^{(n)}_2$ and $W^{(n)}_5$. Since $\lim R_n =0 = \lim {r_n\over
R_n}$, the conformal radii of these annuli increase infinitely. This shows
that $W^{(n)}_2$ and $W^{(n)}_5$ are sequences of long cylinders and that the
sequence $W^{(n)}_3$ defines in the limit a sphere with three nodal points.

As in {\sl Subcase 3\/$'_b$)} every covering piece $W^{(n)}_i$ is
$\tau$-invariant, whereas $W^{(n)\,+}_i \deff W^{(n)}_i \cap \cala_0^+$
is the ``half'' of $W^{(n)}_n$. Thus, for a sequence of undoubled maps $u_n:
\cala^+_0 \to X$ we obtain the following bubbling picture. The limit contains 
two bubbled disks, represented by $W^{(n)\,+}_6$ and $W^{(n)\,+}_3$, 
a boundary node
$W^{(n)\,+}_4$ and possibly further bubbled pieces which can appear in the
limit of long strips $W^{(n)\,+}_2$ and $W^{(n)\,+}_5$. Note also that the
action of the involution $\tau$ on the pants $W^{(n)}_3$ is described in
Fig.~12\.b).

\medskip
In {\sl Subcase 3\/$'''_b$)} the bubbling takes place near but not at the
boundary. Indeed, since ${\rho_n \over r_n}\lrar \infty$, the bubbled sphere
which appears as the limit of the sequence $\{v_n \}$ is not $\tau$-invariant.
To see this phenomenon, we note that for any fixed $b>0$ the covering pieces
$W^{(n)}_5= \Delta' (x_n,2br_n)$ representing a sufficient big part of this
sphere lie in $\cala^+_0$ for $n>\!>0$. This implies that the sequence $v_n
\scirc \tau$ converges to another bubbled sphere, which is $\tau$-symmetric
to the first and represented by $W^{(n)}_7$.

Another bubbled sphere, represented by $W^{(n)}_3$, appears from
pants between the first two spheres and the disk $\Delta'$. Since $\{ {R_n
\over \rho_n}\}$ remains bounded, the original nodal point remains on
this latter sphere.
\egroup

%\smallskip
\line{%
\vtop{\hsize=.53\hsize
The corresponding bubbling picture for undoubled maps $u_n: \cala^+_0 \to X$
is shown in Fig.~15. The boundary of $\cala^+_0$ is shown by a thick line.
We obtain the bubbled sphere represented by $W^{(n)}_6$, the sequence of long
cylinders $W^{(n)}_8\!$, the bubbled disk $W^{(n)+}_3\!\!$ and the sequence
of long strips $W^{(n)+}_2$. Note that both sequences of long cylinders
and long strips can yield  further bubblings in the limit.
}
\hfil
\vtop{\hsize=.44\hsize\xsize=\hsize\nolineskip\rm
\putm[-.02][.34]{W^{(n)+}_1}%
\putm[.21][.26]{W^{(n)+}_2}%
\putm[.51][.42]{W^{(n)+}_3}%
\putm[.84][.34]{W^{(n)+}_4}%
\putm[.77][.07]{W^{(n)}_5}%
\putm[.50][.065]{W^{(n)}_6}%
\noindent
\epsfxsize=\xsize\epsfbox{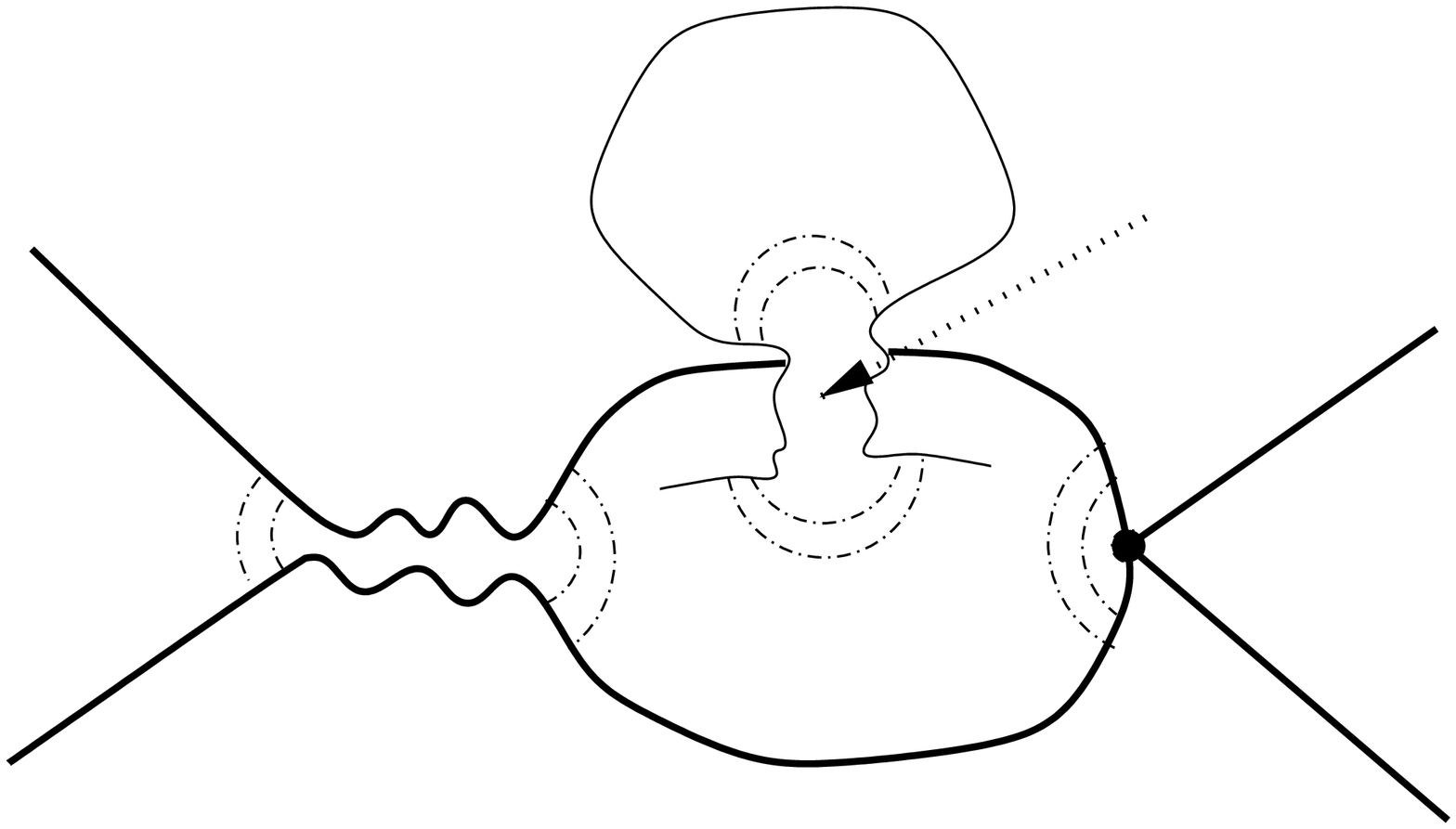}
\smallskip\smallskip
\putt[.0][.0]{
\centerline{Fig.~15. Bubbling in {\sl Subcase 3\/$'''_b$)}.}
}
\vskip\baselineskip\vskip11pt
}}

\medskip\smallskip
\line{%
\vtop{\hsize=.53\hsize
The bubbling picture in {\sl Subcase 3\/$''''_b$)} is similar to the one of
the previous subcase; therefore, we explain only the difference. It comes from
the fact that the sequence $\{ {R_n \over \rho_n}\}$ is now unbounded, \i.e.,
$\lim {\rho_n \over R_n}=0$. Informally speaking, this means that the long
cylinder from {\sl Subcase 3\/$'''_b$)} (piece $W^{(n)}_5$ in Fig.~15)
moves to the boundary of the bubbled disk (piece $W^{(n)+}_3$ in Fig.~15).
The procedure of additional re\-scaling divides such a disk into two
new disks connected by a strip (pieces $W^{(n)+}_3\!\!$, $W^{(n)+}_6\!\!$ and
$W^{(n)+}_5$ in Fig.~16, respectively). The infinite growth ${R_n \over
\rho_n} \to \infty$ means that $W^{(n)+}_5$ form a sequence of long strips.
}
\hfil
\vtop{\hsize=.44\hsize\xsize=\hsize\nolineskip\rm
\putm[-.02][.70]{W^{(n)+}_1}%
\putm[.22][.62]{W^{(n)+}_2}%
\putm[.51][.75]{W^{(n)+}_3}%
\putm[.86][.71]{W^{(n)+}_4}%
\putm[.82][.36]{W^{(n)+}_5}%
\putm[.51][.39]{W^{(n)+}_6}%
\putm[.80][.07]{W^{(n)}_7}%
\putm[.51][.065]{W^{(n)}_8}%
\noindent
\epsfxsize=\xsize\epsfbox{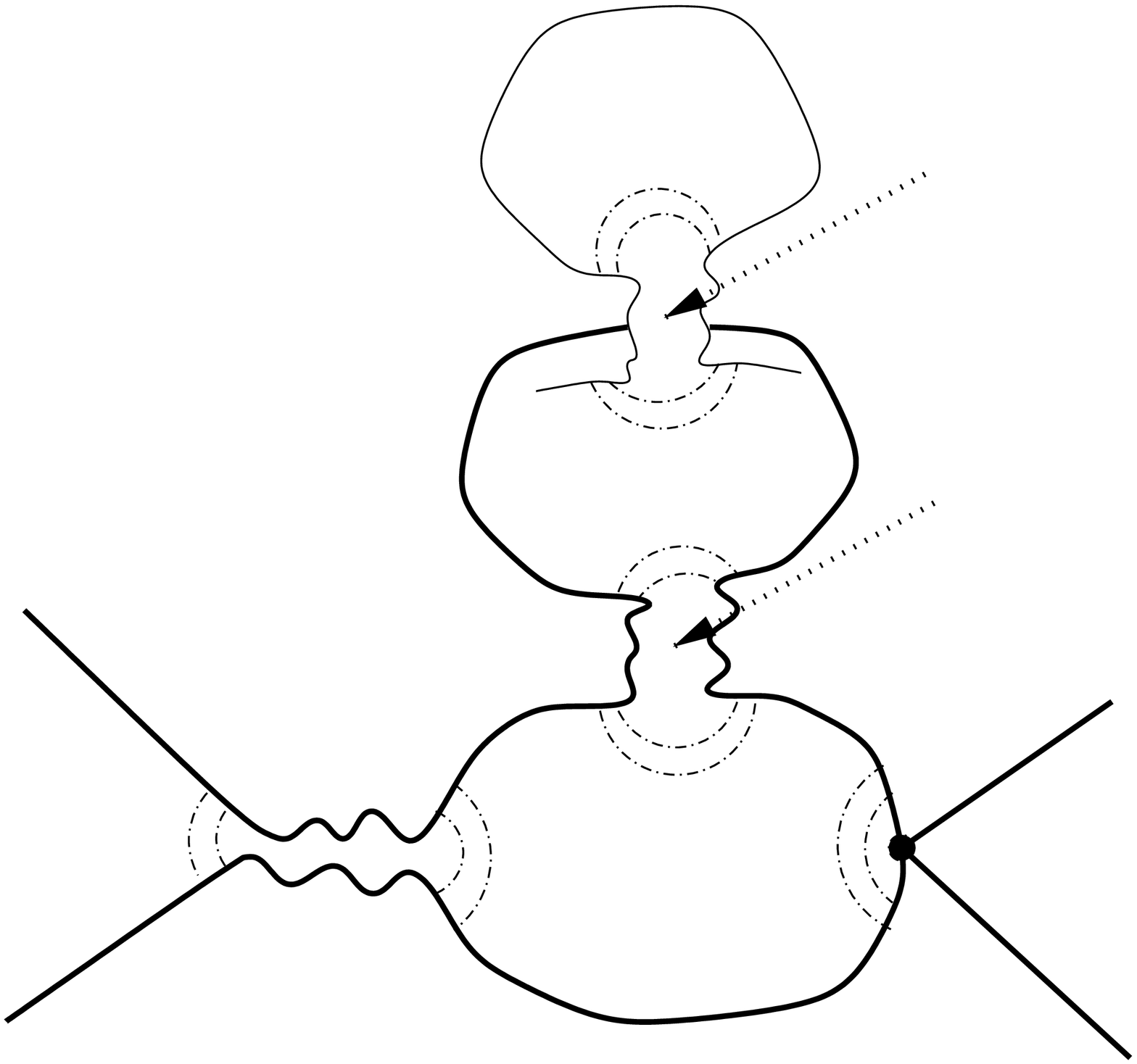}
\smallskip\smallskip
\putt[.0][.0]{
\centerline{Fig.~16. Bubbling in {\sl Subcase 3\/$''''_b$)}.}
}
\vskip\baselineskip\vskip11pt
}}

\bigskip\noindent{\bigsl
A3.5. Attaching an Analytic Disk to a Lagrangian Submanifold of $\cc^n$}.

\medskip\rm
As in the first lecture we consider $\cc^n \approx \rr^{2n}$ together with  
some symplectic form $\omega $ taming the standart complex structure 
$J\st $. An $n$-dimensional submanifold 
$W\subset \cc^n$ is called $\omega $-Lagrangian if $\omega\mid_W\equiv 0$. 

\smallskip\noindent\bf
Exercise 1. \rm Prove that every  Lagrangian submanifold of $\cc^n$ is totally 
real.

\noindent
{\bf 2.} Prove that for a Lagrangian manifold $W\subset \cc^n$ and the unit 
circle $\ss^1\subset \cc $ the manifold $\ss^1\times W$ is Lagrangian in 
$\cc^{n+1}$ with respect to $\hat\omega = {i\over 2}dz_1\wedge d\bar z_1 + 
\omega $.

\smallskip A holomorphic map $u:\Delta \to \cc^n$ "sufficiently smooth" up 
to the boundary and such that $u(\d \Delta )\subset W$ we shall call \it an 
analytic disk attached to $W$. \rm 

Our goal in this paragraph is to prove the following theorem of Gromov:

\smallskip
\state Theorem A3.5.1. {\it Let $W$ be a compact Lagrangian submanifold 
of $\cc^n$. Then there exists a non-constant analytic disk attached to 
$W$.
}

\smallskip We shall closely follow the exposition of H. Alexander, [Al]. 
Fix some point $w_0\in W$ and denote by $u_0(z)\equiv w_0$ the  
constant holomorphic map. Fix $p>2$ and consider the Banach manifold 
$L^{2,p}(\Delta , \d \Delta , 1; \cc^n , W, w_0)$ of $L^{2,p}$-maps from 
$\bar\Delta $ to $\cc^n$ which map $\d \Delta $ to $W$ and $1$ to $w_0$, 
and which are homotopic to the constant map $u_0\equiv w_0$  
as $(\Delta , \d \Delta , 1)\to (\cc^n , W, w_0)$ mappings. 
Note that due to the Sobolev imbedding $L^{2,p}\subset C^{1,1-{2\over p}}$ 
our mappings are smooth up to the boundary.

Take $u\in L^{2,p}(\Delta , \d \Delta , 1; \cc^n , W, w_0)$. Denote by 
$E$ as usually the pull-back by $u$ of the tangent bundle of $\cc^n$. In 
fact $E=\Delta\times \cc^n \to \Delta $, the trivial bundle over $\Delta $. 
Denote by $F$ the pull-back $u^*TW$ of the tangent to $W$ bundle. $F$ is a 
totally real subbundle of $E$ of real dimension $n$.   
The tangent space to $L^{2,p}(\Delta , \d \Delta , 1; \cc^n , W, w_0)$ at $u$ 
is $T_uL^{2,p}(\Delta , \d \Delta , 1; \cc^n , W, w_0) = 
\{ h\in L^{2,p}(\Delta ,\cc^n): h\mid_{\d \Delta }\in F, h(1)=0\} $.
In the cartesian product $L^{2,p}(\Delta , \d \Delta , 1; \cc^n , W, w_0)
\times L^{1,p}(\Delta ,\cc^n)$  consider the submanifold $\cale =\{ (u,v):
\dbar u=v\} $ 
with the natural projection $\pi :\cale \to L^{1,p}(\Delta ,\cc^n)$. 

\smallskip\noindent\bf
Exercise 1. \rm In $\rr^{2n+2}$ consider the operator given by the 
matrix 

$$
J_v = \left(
\matrix
0 & -1 & ......  & v_2       & -v_1 \cr
1 &  0 & ......  & -v_1      & -v_2 \cr
. &  . & ....... & ......    & ... \cr
0 &  0 & ....... & v_{2n}    & -v_{2n-1} \cr
0 &  0 & ....... & -v_{2n-1} & -v_{2n} \cr 
0 &  0 & ....... & 0         & -1 \cr
0 &  0 & ....... & 1         & 0 \cr
\endmatrix
\right) 
$$
Prove that $J_v$ defines an almost complex structure in $\rr^{2n+2}$, which 
for every $z\in \rr^2$ on the vertical slice $\{ z\} \times \rr^{2n}$ coincides 
with the standart structure $J\st$ of $\cc^n$.

\noindent
{\bf 2.} Prove that the equation $\dbar u=v$ for a $C^1$-map $u:\cc \to \cc^n$
is equivalent to the $J_v$-holomorphicity of the section $\hat u :z\to (z,u(z))$ of 
the fibration $(\rr^{2n+2},J_v)\to (\rr^2, J\st)$, i.e. to the equation 

$$
{\d (z,u)\over \d x} + J_v(z,u)[{\d (z,u)\over \d y}] = 0.
$$

\smallskip Now we shall prove the following alternative:

\state Lemma A3.5.2. {\it If there is no  nonconstant analytic disk 
$u\in  L^{2,p}(\Delta , \d \Delta , 1;$ $ \cc^n , W, w_0)$, then the projection 
$\pi :\cale \to L^{1,p}(\Delta ,\cc^n)$ is surjective.
}
\state Proof. {\rm Suppose that a nonconstant analytic disk 
$u\in L^{2,p}(\Delta , \d \Delta , 1; \cc^n , W, w_0)$ doesn't exists.  We are 
going to prove that in this case $\pi $ is surjective. 

\smallskip\noindent
{\sl Step 1. \it $\pi $ is a proper mapping, i.e. for the converging sequence 
$v_k\to v_0$ in $L^{1,p}(\Delta ,\cc^n)$   
and for the sequence $u_k$ with $(u_k,v_k)\in \cale $, there is a converging 
subsequence $u_{k_n}$.
}

\smallskip
Note that $\dbar u_k=v_k$. According to the Exercise above this means that 
the sections $\hat u_k:=(z,u_k)$ are $J_{v_k}$-holomorphic with $J_{v_k}$ 
converging to $J_{v_0}$ in $C^0$-sence. Note also that boundaries of our disks 
are 
on the Lagrangian (and thus totally real) submanifold $\hat W:=\ss^1\times W$ and they 
are homotopic to each other. From here we see that 

$$
\area(\hat u_k(\Delta )) = \int_{\hat u_k(\Delta )}\hat\omega = 
\int_{\hat u_k(\d \Delta )}\d \lambda ,
$$
where $\lambda $ is some primitive of $\hat\omega $. The second integral does
not depend on the homology class of $\hat u_k(\Delta )$ in $\sfh_1(\hat W, 
\rr)$, because $\hat\omega\mid_W = \d \lambda\mid_{\hat W} \equiv 0$ ($\hat W$
is $\hat\omega $-Lagrangian!).
 
So by the Theorem A3.4.1 either the limit of some subequence, still denoted as 
$u_k(\Delta )$, contains a nonconstant complex sphere (this is impossible in 
$\cc^n$), or the limit of $u_k(\Delta )$ contains some nonconstant analytic 
disk with boundary on $W$ (this is prohibited by our assumption), or $u_k$ 
$C^1$-converge. 

\smallskip\noindent
{\sl Step 2. \it $d\pi_{(u,v)} : T_{(u,v)}\cale \to L^{1,p}(\Delta ,\cc^n)$ is 
Fredholm of index zero for every $u\in \cale_v :=\pi^{-1}(v)$.
}
\smallskip
The fact that the boundary value problem $\dbar h=v, h\mid_{\d \Delta }\in F$ 
is Fredholm is classical, see [Ga]. Our manifold $\cale $ is connected, 
so the index of $d\pi_{(u,v)}$ doesn't depend on $(u,v)$ and can be calculated 
at $(u_0,0)\in \cale $ where $d\pi $ is a bijection. Therefore $ind(d\pi )$ is 
everywhere zero.
  
\smallskip One says that for the smooth  mapping 
$\pi :\cale \to L^{1,p}(\Delta ,\cc^n)$ the point $(u,v)$ is regular if 
$d\pi_{(u,v)}$ is surjective. $v$ is a regular value is it is not an image of 
a nonregular point.

\smallskip\noindent
{\sl Step 3. (Smale's theorem). \it Let $\pi :\cale \to \calm $ be a proper 
Fredholm map (i.e. $d\pi_x $ is Fredholm for all points $x\in \cale $. Then 
the set of regular values is dence in $\calm $. Moreover, for every regular 
value $v\in \calm $ the set $\cale_v:=\pi^{-1}(v)$ is a manifold of 
dimension equal to the $\ind(d\pi_x) $ at $x\in \cale_v$. Moreover, for 
any two regular values $v_1$ and $v_2$ the manifolds $\cale_{v_1}$ and 
$\cale_{v_2}$ are cobordant.  
} 
See [Sm].

\smallskip In our case $0$ is a regular value, so for a dence subset of 
$v$'s $\cale_{v}$ is cobordant to a point, therefore is a point intself. 
Properness of $\pi $ implies now that $\cale_v$ is allways a point.
}

\qed
 
\smallskip It is not difficult to show that $\pi $ cannot be surjective. 
This will imply the existence of non constant analytic disk attached to $W$.

\state Lemma A3.5.3. {\it The projection $\pi :\cale \to L^{1,p}(\Delta ,\cc^n)$ 
is not surjective.
}
\state Proof. {\rm Othervice, for $v^C=(C,0,....,0)$ find 
$u^C=(u_1^C,...,u_n^C)$ with 
$\dbar u^C = v^C$. Therefore $\bar \d u_1^C=C$. This implies that 
$u_1^C=C\bar z - h^C$, where $h^C$ is a holomorphic function on $\Delta $. 
Since $u^C(\d \Delta )\subset W$ the family $u^C$ is uniformly bounded on 
$\d \Delta $ by a constant $k$ independent of $C$. Therefore $\vert \bar z - 
h^C(z)/C\vert \le k/C$ for $z\in \d \Delta $. Since $\bar z - h^C(z)/C$ is 
harmonic on $\Delta $, the bound holds for all $z\in \Delta $. This implies 
that $\bar z$ can be uniformly approximated on $\Delta $ by holomorphic 
functions. Contradiction.
}

\qed

%% ch.3

\newpage

\bigskip\noindent
{\bigbf Chapter III. Global Properties and Moduli Spaces.}

\medskip
Let $(X, J)$ be an almost-complex manifold of complex dimension $n$. An
almost complex structure is always assumed to have smoothness of class $C^1$.
Further, let $(S, J_S)$ denote a Riemann surface with complex structure $J_S$.

In Lecture 7 we show that for a $J$-complex curve $u:(S,J_S) \to (X,J_0)$
the pulled-back bundle $E\deff u^*TX$ possesses a natural {\sl holomorphic}
structure (the corresponding sheaf of holomorphic sections will be denoted as
$\calo(E)$) such that the differential $du:TS\to E$ is a {\sl holomorphic}
homomorphism. This allows us to define the order of vanishing of
the differential $du$ at point $s\in S$. We denote this number by
${\sl ord}_sdu$.

This also gives the following short, exact sequence:
$$
0\longrightarrow \calo(TS) \buildrel du \over\longrightarrow
\calo(E) \buildrel \pr \over\longrightarrow
\calo(N_0)\oplus \caln_1 \longrightarrow 0,\eqno(3.1)
$$
where  $\calo(N_0)$ denotes a free part of the quotient
$\calo(E)/du(\calo(TS))$, and
$\caln_1$ is supported on a finite set of cusps of $u$ ({\sl \ie}
points of vanishing of $du$).

On the Sobolev space $L^{1,p}(S,N_0)$ of $L^{1,p}$-smooth sections
of the bundle $N_0$ the natural Gromov operator $D_N: L^{1,p}(S,N_0) \to
L^p(S,\Lambda^{0,1}S \otimes N_0)$ is defined. Put $\sfh^0_D(S, N_0):=
\ker D_N$ and $\sfh^1_D(S, N_0):=\coker D_N$.

We prove the following
\state Theorem 3.1. {\it  Let $u: (S, J_S) \to (X, J_0)$ be a nonconstant
irreducible and non- multiply covered holomorphic map such that
$\sfh^1_D(S, N_0)=0$. Then

\smallskip
\sli  in a neighborhood of $M\deff u(S)$ the Moduli space of
nonparameterized $J_0$-holomorphic curves $\calm_{[\gamma], g, J_0}$ is
a manifold whose tangent space is $T_M\calm_{[\gamma], g, J_0} =
\sfh^0_D(S, N_0) \oplus \sfh^0(S, \caln_1)$;

\smallskip
\slii further, there is a neighborhood $V\ni J$ in the Banach manifold
$\calj$ of $C^1$-smooth almost-complex structures on $X$ and a neighborhood
$W$ of $M$ in $\calm_{[\gamma], g, V}:=\bigcup_{J\in V}
\calm_{[\gamma], g, J}$ such that the natural projection $\pr_\calj:
W\to V$ is a trivial Banach bundle;

\smallskip
\sliii if $\dim_{\rr }X=4$ and $c_1(E)[M] > \sum_{s\in S}{\sl ord}_sdu$, then
$\sfh^1_D(S, N_0)=0$; thus, the conclusions \sli and \slii hold.
}

\smallskip
Let us explain our second result in this chapter, which will be strongly
relied upon in Chapter IV and also in the proof of the
Gromov non-squeezing theorem.
Denote by $\calj $ some Banach manifold of $C^k$-smooth almost-complex
structures on a manifold $X$, $k\ge 1$. In our applications $\calj $ can
be

\smallskip
1) the manifold $\calj_{\omega }$ of all a.-c. structures $J$ on $X$, which
are tamed by some fixed symplectic form $\omega $; \i.e., $\omega (u,Ju)>0$
for all $u\in TX\setminus \{ 0\} $.

\smallskip
2) the manifold $\calj_U$ of all a.-c. structures $J$ on a complex K\"ahler
manifold $(X,J\st ,\omega )$ tamed by $\omega $ and which are different
from $J\st $ only on some relatively compact subset of $U$, \i.e.,
$\{ x\in X:J_x\not= J\st \} \comp U$.

\smallskip Denote by $\calm =\calm_{\calj ,[\gamma ],\Sigma }$ the
(topological)
space of all $J$-complex curves in $X$ ($J$ runs over all structures in
$\calj $) representing some fixed homology class $[\gamma ]\in \sfh_2(X,\zz )$,
and parameterized by some fixed compact surface $\Sigma $, see Chapter II.
For $\calj $, as in 1) and 2), above the space $\calm $ is a manifold in
all its non- multiply covered points, see {\sl Lemma 8.2.2}.

Also let $h:[0,1]\to \calj $ be some smooth map (more generally one
can consider a smooth map $h:Y\to \calj $, where $Y$ is a compact real
manifold with boundary). By $\calm_h$ we denote the subset of $\{ (u,J_S,t)
\in \cals \times \ttt_{\Sigma }\times [0,1]\} $ which consists of
$h(t)$-complex curves $u:(S,J_S)\to (X,h(t))$ modulo the natural action
of the group $G$ of reparameterizations, see \S 7.4 and Definition 8.3.1.
Roughly speaking, $\calm_h$ is a moduli space of $h(t)$-complex curves in
$X$, $t\in [0,1]$. This is a closed subset of $\calm $. For a point
$(M_0,t_0)$ in $\calm_h$ (\ie for a $h(t_0)$-complex curve $M_0$) one can
define a component $\calm_h(M_0,t_0)$ of $\calm_h$ containing $(M_0,t_0)$
as was shown in {\sl Definition 8.3.2}.

\smallskip
\state Theorem 3.2. {\it If the component $\calm_h(M_0,t_0)$ is not compact,
then there exists a continuous path $\beta :[0,1]\to \calm_h(M_0,t_0)$,
$\beta (t) =(M_t,J_t)$
starting at $(M_0,t_0)$ such that $\beta (t_n)$ is not bounded in
$\calm_h(M_0,t_0)$ for some sequence $t_n\to 1$, but $J_{t_n}$
converge to some $J^*\in \calj $.
}

\smallskip This means, in fact, that in the limit the sequence $M_{t_n}$ 
breaks into several irreducible $J^*$-complex components.

\smallskip Let us consider an example where such behavior is not possible.

\smallskip
\state Example. {\rm Suppose that $X=S^2\times Y$ and $\omega =\omega_1
+ \omega_2$, where $\omega_1$ is some positive form on $S^2$  and
$\omega_2$ is some symplectic form on a compact manifold $Y$. As a homology
class $[\gamma ]$ take $S^2\times \{ pt\} $ and as a parameterizing surface
for our curves again $\Sigma :=S^2$.

Observe that for any structure $J$ on $X$ tamed by $\omega $ no $J$-complex
curve $C\in [\gamma ]$ can be decomposed as $C=C_1\cup C_2$ for some
$J_k$-complex curves $C_k$, where $J_k\in \calj_{\omega }$!

\smallskip
\state Exercise. {\rm Prove this.}

\smallskip
In this example $\calm_h(M_0,t_0)$ will be always be compact.
}

\smallskip\noindent
\state Theorem 3.3. {\it Under the above conditions suppose that $S$ is a
sphere $S^2$ and that $\calm_h(M_0,t_0)$ is compact. Then the path $h$
can be $C^1$-approximated by the paths $h_n$, all starting at $(M_0,t_0)$
such that for all points
$(M,t)\in \calm^0_{h_n}$
the associated $D_N$-operator is surjective, \i.e., 
$\sfh^1_{D_N}(S^2,N_{M_t})=0$.

Moreover, $\calm_{h_n}$ is a trivial bordism: $\calm_{h_n(0)}\times [0,1]$.
In particular for every $h_n(0)$-holomorphic sphere $M_0\in \calm_{h_n(0)}$
there exists a continuous family of $h_n(t)$-holomorphic spheres $M_{n,t}=
u_{n,t}(S^2)$ with $M_{n,0}=M_0$.
}

\newpage

\bigskip\noindent
{\bigbf Lecture 6}

\smallskip\noindent
{\bigbf  First Variation of the $\dbar_J$-equation.}

\medskip
\rm
In this paragraph we want to introduce, following Gromov, a $D_{u,J}$-operator
associated to the $J$ - complex  curve $u : S \to X$ in an almost complex
manifold $(X,J)$. For the convenience of the reader we now recall the
notion of a connection.

\smallskip\noindent\sl
6.1. Symmetric Connections.

\smallskip\rm
Let $M$ be a real manifold, $E$ a real vector bundle over $M$. By
$\Gamma (E)$ we denote the space of smooth sections of $E$ and by $\Lambda^p$,
a space of smooth $p$-forms on $M$.

\smallskip
\state Definition 6.1.1. {\it A connection on the real vector bundle $E\to M$
is a mapping $\nabla :\Gamma (E) \to \Lambda^1\otimes \gamma (E)$ such that

\smallskip\noindent
$\nabla 1$. $\nabla $ is $\rr $-linear;

\smallskip\noindent
$\nabla 2$. for any smooth function $f$ on $M$ and any section $e\in \Gamma (E)$

$$
\nabla (f\cdot e) = df\otimes e + f\cdot \nabla e.\eqno(6.1.1)
$$
}

\smallskip
If $X\in \Gamma (TM)$ is a vector field on $M$, one denotes $\nabla_Xe :=
\nabla (e)(X)$, a covariant derivative of $e$ along $X$.

\smallskip
\state Remarks.~1. 
Choose a local frame $e_1,...,e_n$ of $E$ and write $\nabla (e_i) =
\Gamma_i^j\otimes e_j$, where $\Gamma_i^j$ are some $1$-forms. Then
$\nabla $ is uniquely determined by $\{ \Gamma_i^j\}_{i,j=1}^n$.

\state 2. The dual connection $\nabla^*$ is defined on the 
dual bundle $E^*\to M$ by
$$
d\< h,e\> = \<\nabla^*h,e\> + \<h, \nabla e\>, \eqno(6.1.2)
$$
where $h\in \Gamma (E^*)$.

\state 3. 
One defines an extension of $\nabla $ onto $E^{\otimes p}\otimes
E^{\otimes q}\to M$ as
$$
\nabla (e\otimes h) = \nabla e\otimes h + e\otimes \nabla h. 
\eqno(6.1.3)
$$
In particular, since $\hom (E,E) \equiv E\otimes E^*$, $\nabla $ extends
onto $\hom (E,E)$.

We shall usually denote all these extensions by the same symbol $\nabla$. One
remarks, of course, that (6.1.1), (6.1.2) and (6.1.3) are nothing  but the 
Leibnitz formula.

\smallskip\noindent
4. Mapping $\nabla^2 = \nabla\scirc \nabla :\Gamma (E) \to \Lambda^2\otimes
\Gamma (E)$ is called the curvature of $\Gamma $ and is, in fact, linear
over $\cala (M)$: $\nabla^2(fe)=\nabla (df\otimes e+f\nabla e)=
d^2f\otimes e-df\wedge \nabla e+df\wedge \nabla e+f\nabla^2e= f\nabla^2e$.
}

\smallskip
\state Definition 6.1.2. {\it Let $E=TM$. A connection 
$\nabla $ is called symmetric if
$\nabla_XY-\nabla_YX = [X,Y]$.
}

\smallskip
\state Exercise. {\rm Take $e_i={\d \over \d x_i}$ for some local coordinates
$\{ x_i\} $. Write $\Gamma_i^j=\gamma_{ik}^jdx^k$, \i.e., $\nabla_{e_k}e_i =
\Gamma_{ik}^je_j$. What does symmetricity of $\nabla $ mean in terms of $\{
\Gamma_{ik}^j\} $?
}

\smallskip
\state Exercise. {\rm Choose a symmetric connection $\nabla $. Then the 
Nijenhuis tensor of an almost-complex structure $J$ is defined as

$$
4N_J(X,Y) = (\nabla_XJ\scirc J)Y - (\nabla_YJ\scirc J)X + (\nabla_{JX}J)Y -
(\nabla_{JY}J)X.
$$
Prove that this defintion doesn't depend on the choice of a symmetric 
connection, i.e. express $N_J(X,Y)$ in the terms of Lee brackets.
}

\medskip\noindent
\sl 6.2. Definition of the $D_{u, J}$-operator.

\smallskip\rm
Recall that a $C^1$ map $u:S\to X$ from a Riemann surface $S$ with a complex
structure $J_S$ is called holomorphic with respect
to $J_S$ and $J$ if it satisfies the equation

\smallskip
$$
du + J\scirc du\scirc J_S = 0. \eqno(6.2.1)
$$

This simply means that $du\scirc J_S = J\scirc du$, {\sl \i.e.,} $du:T_pS \to
T_{u(p)}X$ is a {\sl complex} linear map for every $p\in S$. Equation
(6.2.1)
is an elliptic quasi-linear PDE of order one. We are interested in the behavior
of the solutions of (6.2.1), in particular, when the structures $J$ and $J_S$
change. So we need to choose appropriate functional spaces both for solutions
and for the coefficient of (6.2.1). Our choice is based on the following
facts:

{\sl a)} The minimal reasonable smoothness of an almost complex structure $J$
on $X$, for which the Gromov operator $D_{u, J}$ can be defined, is $C^1$,
see the explicit formula (6.2.5).

{\sl b)} It is more convenient to operate with Banach spaces and manifolds and
thus with finite regularity like $C^k$, $C^{k,\alpha}$ or $L^{k,p}$ than
with $C^\infty$-smoothness defining only a Frechet-type topology.

{\sl c)} Equation (6.2.1) is defined also for $u$ lying in Sobolev-type
spaces $L^{k,p}(S, X)$ with $k\ge1$, $1\le p\le\infty$ and $kp>2$; such
solutions are $C^1$-smooth and the topology on the space of a solution is,
in fact, independent of the particular choice of such a Sobolev space.

{\sl d)} For $J\in C^k$ the coefficients of the Gromov operator $D_{u, J}$
are, in general, only $C^{k-1}$-continuous, see (6.2.5), and hence solutions
of the ``tangential equation'' $D_{u, J}v=0$ are only $L^{k,p}$-smooth, $1\le
p<\infty$; thus for obtaining a {\sl smooth} structure on a space
of (parameterized) $J$-complex curves, one should use Sobolev spaces
$L^{k', p}(S, X)$ with $k'\le k$.

\medskip
Fix a compact Riemann surface $S$, {\sl \i.e.,}\ a compact, connected,
oriented, smooth manifold of real dimensions 2. Recall that the Sobolev space
$L^{k, p}(S, X)$, $kp>2$, consists of those continuous maps $u: S \to X$,
which are represented by $L^{k, p}$-functions in local coordinates on $X$ and
$S$. This is a Banach manifold, and the tangent space $T_uL^{k, p}(S, X)$
to $L^{k, p}(S, X)$ in $u$ is the space $L^{k,p }(S, u^*(TX))$ of all
$L^{k,p}$-sections of the pull-back under $u$ of the tangent bundle $TX$.
One has the Sobolev imbeddings
$$
\matrix
L^{k,p}(S, X)& \hookrightarrow & L^{k-1, q}(S, X),&
\quad & \hbox{for $1\le p<2$}   &  \hbox{and} & 1\le q\le {2p\over2-p},
\cr
L^{k, p}(S, X)& \hookrightarrow & C^{k-1, \alpha}(S, X),&
\quad & \hbox{for $2<p\le\infty$}& \hbox{and} & 0\le\alpha\le 1-{2\over p}.
\cr
\endmatrix
$$

\smallskip
Let $[\gamma ]$ be some homology class in  $\sfh_2(X, \zz)$. Fix $p$ with
$2<p<\infty$ and consider the Banach manifold
$$
\cals = \{u\in L^{1,p}(S, X):u(S)\in [\gamma]\}
$$
of all $L^{1,p }$-smooth mappings from $S$ to $X$, representing the~class
$[\gamma]$. This makes sence, because $L^{1,p}\subset C^{0,1-{2\over p}}$, 
and therefore all $u$ from $\cals $ are continuous. Remark that the tangent 
space to $\cals $ at $u$ is 

$$
T_u\cals = L^{1,p}(S,u^*(TX))
$$
the space of $L^{1,p}$-sections of the pulled-back by $u$ the tangent bundle 
of $X$.

Denote by $\calj$ the Banach manifold of $C^1$-smooth almost
complex structures on $X$. In other words, $\calj = \{ J\in C^1
(X, \endo (TX)):J^2=-\id  \} $. The tangent space to $\calj$ at $J$
consists of $C^1$-smooth $J$-antilinear endomorphisms of $TX$,
$$
T_J\calj = \{ I\in C^1(X, \endo (TX)) : JI + IJ = 0\}
\equiv C^1(X, \Lambda^{0,1}X \otimes TX),
$$
where $\Lambda^{0,1}X$ denote the complex bundle of $(0,1)$-form on $X$.

Denote by $\calj_S$ the Banach manifold of $C^1$-smooth complex structures
on $S$. Thus, $\calj_S = \{ J_S\in C^1(S, \endo (TS)):J_S^2=-\id  \}$
and the tangent space to $\calj_S$ at $J_S$ is
$$
T_{J_S}\calj = \{ I\in C^1(S, \endo (TS)) : J_S I + IJ_S = 0\}
\equiv C^1(S, \Lambda^{0,1}S\otimes TS).
$$

\smallskip
Consider also the subset $\calp\subset \cals\times \calj_S \times
\calj$ consisting of all triples $(u, J_S, J)$ with $u$ being
$(J_S,J)$-holomorphic, {\sl \i.e.,}

\smallskip
$$
\calp = \{(u, J_S, J)\in \cals\times \calj_S\times \calj: du +
J\scirc du\scirc  J_S = 0 \}. \eqno(6.2.2)
$$

\state Lemma 6.2.1.
\it Let $J$ and $J_S$ be continuous almost-complex structures on $X$ and $S$,
respectively, and let $u\in L^{1,p}(S, X)$. Then
$$
du + J\scirc du\scirc J_S \in L^p(S,\Lambda^{0, 1}S \otimes u^*(TX)).
$$

\state Proof. \rm One can easily see that $du \in L^p(S,\hom_\rr(TS, u^*(TX))$.
On the other hand,
$$
(du + J\scirc du\scirc J_S ) \scirc J_S= -J \scirc
(du + J\scirc du\scirc J_S ),
$$
which means that $du + J\scirc du\scirc J_S$ is $L^p$-integrable
$u^*(TX)$-valued (0,1)-form.

\qed

\smallskip
Consider a Banach bundle $\calt\to \cals\times \calj_S\times \calj$
with a fiber
$$
\calt_{(u, J_S, J)} = L^p(S,\Lambda^{0, 1}S\otimes
u^*(TX)),
$$
where $TS$ and $TX$ are equipped with complex structures $J_S$ and $J$,
respectively. $\calt$ has two distinguished sections:

\smallskip
1) $\sigma_0 \equiv 0$, the~zero section of $\calt$;

2) $\sigma_{\dbar}(u, J_S, J) = du + J\scirc du\scirc J_S$.

\smallskip
\noindent
By definition $\calp$ is the zero-set of $\sigma_{\dbar}$.

\medskip
Let us compute the tangent space to $\calp$ at the point $(u, J_S, J)$. Let
$(u_t, J_S(t),J(t))$ be a curve in $\calp$ such that $(u_0, J_S(0), J(0)) =
(u, J_S, J)$. Let
$$
(v,\dot J_S,\dot J) \deff
\left({du_t\over dt}|_{t=0},
{dJ_S(t)\over dt}|_{t=0},
{dJ(t)\over dt}|_{t=0}\right)
$$
be the tangent vector to this curve and hence to $\calp$ at $t=0$.
The condition $(u_t, J_S(t), J(t)) \in \calp$ means that
$$
du_t + J(u_t, t)\scirc du_t\scirc J_S(t) = 0 \eqno(6.2.3)
$$
in $ L^p(S,\Lambda^{0, 1}(S)\otimes u_t^*(TX))$. Let $\nabla $ be
some symmetric connection on $TX$, {\sl \i.e.,} $\nabla_YZ - \nabla_ ZY = [Y,
Z]$. The co-variant differentiation of (6.2.3), with respect to $t$, gives
$$
\nabla_{\partial \over \partial t}(du_t) + (\nabla_vJ)(du_t\scirc J_S)
+ J(u_t, t)\scirc \nabla_{\partial \over \partial t}(du_t)\scirc J_S +
$$
$$
+ J\scirc du_t\scirc \dot J_S +
\dot J\scirc du_t\scirc J_S = 0.
$$
Let us show that $\nabla_{\partial \over \partial t}(du_t) = \nabla v$.
Indeed, for $\xi \in TS$ one has
$$
(\nabla_{\partial \over \partial t}
du_t)(\xi ) = \nabla_{\partial \over \partial t}[du_t(\xi )] =
\nabla_{\partial \over \partial t}({\partial u_t\over \partial \xi }) =
\nabla_{\xi }({du_t\over dt}) = \nabla_{\xi }v.
$$
So every vector $(v,\dot J_S,\dot J)$ which is tangent to $\calp$
satisfies the equation

\smallskip
$$
\nabla v + J\scirc \nabla v\scirc J_S + (\nabla_vJ)\scirc (du\scirc
J_S) + J\scirc du\scirc \dot J_S +
\dot J\scirc du\scirc J_S = 0. \eqno(6.2.4)
$$

\smallskip
\noindent
{\bf Definition 6.2.1.} \rm
Let $u$ be a $J$ - complex  curve in $X$. Define the operator $D_{u, J}$
on $L^{1,p}$-sections $v$ of
$u^*(TX)$ as
$$
D_{u, J}(v) = \msmall{1\over2}\bigl(\nabla v + J\scirc\nabla v\scirc J_S
+ (\nabla_vJ) \scirc (du\scirc J_S).  \bigr)
\eqno(6.2.5)
$$

\bigskip
\noindent
{\bf Remark.} This operator plays a crucial role in studying  properties of
complex curves. In its definition we use the symmetric connections
instead of those compatible with $J$, as is shown in [G]. The matter is that
one can use the same connection $\nabla$ for changing almost complex
structures $J$. The lemmas below justify our choice.

\bigskip
\noindent
{\bf Lemma 6.2.2. }\it $D_{u, J}$ does not depend on the choice of a~symmetric
connection $\nabla $ and is an $\rr$-linear operator from $L^{1,p}(S,
u^*(TX))$ to $L^p(S, \Lambda^{0, 1}S \otimes u^*(TX))$.

\smallskip
\rm
\noindent
{\bf Proof.} Let $\widetilde \nabla $ be another symmetric connection on $TX$.
Consider the bilinear tensor on $TX$, given by formula $Q(Z, Y) {:=} \nabla_ZY
- \widetilde \nabla_ZY$. It is easy to see that $Q$ is symmetric on $Z$ and
$Y$.
Also note that $\nabla_{\xi }v - \widetilde\nabla_{\xi }v =
\nabla_{du(\xi )}v
- \widetilde \nabla_{du(\xi )}v = Q(du(\xi), v)$ and, in addition, that
$(\nabla_ZJ)(Y) - (\widetilde \nabla_ZJ)(Y) = Q(Z, JY) - JQ(Z, Y)$. From here
one
obtains
$$
2(D_{u, J}v)(\xi ) - 2(\widetilde D_{u, J}v)(\xi) =
$$
$$
= \nabla_{\xi }v - \widetilde \nabla_
{\xi }v + J(\nabla_{J_S\xi }v - \widetilde \nabla_{J_S\xi }v) + (\nabla_vJ -
\widetilde \nabla_vJ)du(J_S\xi ) =
$$
$$
= Q(du(\xi ), v) + JQ(du(J_S\xi ), v) + Q(v, JduJ_S
\xi ) - JQ(v, du(J_S\xi)) =
$$
$$
=Q(du(\xi ), v) + Q(JduJ_S\xi , v) = Q((du + JduJ_S)(\xi ), v) = 0.
$$

Now let us show that $D_{u, J}(v)$ is $J_S$-antilinear:
$$
2D_{u, J}(v)[J_S\xi] = \nabla_{J_S\xi }v + J\scirc \nabla_{J^2
_S\xi }v + (\nabla_vJ)\scirc (du\scirc J^2_S)(\xi ) =
$$
$$
= \nabla_{J_S\xi }v -
J(\nabla_{\xi }v) - (\nabla_vJ)(du(\xi )) = -J[\nabla_{\xi }v + J(\nabla
_{J_S\xi }v) - J\nabla_vJdu(\xi )] =
$$
$$
= -J[\nabla_{\xi }v + J\scirc \nabla_{J_S\xi}v + (\nabla_vJ)
(du\scirc J_S(\xi )] = -2J\, D_{u, J}(v)[\xi].
$$
Here we use the~fact that $J\scirc \nabla_vJ + \nabla_vJ\scirc J = 0$ and 
$du\scirc J_S = J\scirc du $.

\qed

\smallskip
This lemma allows us to obtain expression (6.2.4) also by computation in local
coordinates $x_1, x_2$ on $S$ and $u_1,\ldots, u_{2n}$ on $X$, choosing as 
connection $\nabla$ the de Rham differential $d$. Really, wright $u_t(x) = 
u(x,t) = (u_1(x_1,x_2,t),...,u_{2n}(x_1,x_2,t))$ for the $J_t$-holomorphic map
$u:S\to X$ (note that $u_t$ and structure $J_t$ are both time dependent). 
In the local basis ${\d \over \d u_1},...,{\d \over \d u_{2n}}$ of the tangent
bundle $TX$ $J_t$ is represented by a time dependent $2n\times 2n$ matrix, 
which we will also denote as $J(\cdot ,t)$. Differential of $u_t$ will be 
denoted as $du_t(x)$ or as ${du_t\over dx}$ and is $2\times 2n$ matrix 
$$
du_t(x) =
\left( \matrix
{\d u_1\over \d x_1}(x,t)&{\d u_1\over \d x_2}(x,t)\cr
\cdot \cdot \cdot \cdot &\cdot \cdot \cdot \cdot \cr
{\d u_{2n}\over \d x_1}(x,t)&{\d u_{2n}\over \d x_2}(x,t)\cr
\endmatrix \right).
\eqno(6.2.6)
$$

Complex structure $J_S(t)$ on the surface is also time dependent and in the 
local frame ${\d \over \d x_1},{\d \over \d x_2}$ of $TS$ is represented by a 
$2\times 2$ matrix, which we shall denote also as  $J_S(t)$. Equation of 
holomorphicity (6.2.3) now reeds in matrix form as 
$$
{du_t\over dx} + J(u_t,t)\cdot {du_t\over dx}\cdot J_S(t)\equiv 0.
\eqno(6.2.7)
$$
Differentiating this with respect to $t$ we get 

$$
{d\over dx}\left[{du_t\over dt}\right] + 
{\d J\over \d t}(u_t,t)\cdot {du_t\over dx}\cdot J_S(t) + 
\left\<{\d J(u_t,t)\over \d u},{du_t\over dt}\right\>
\cdot {du_t\over dx} \cdot J_S(t) 
+
$$
$$ 
+ J(u_t,t)\cdot {d\over dx}
\left[{du_t\over dt}\right]\cdot J_S(t) + J(u_t,t)\cdot 
{du_t\over dx}\cdot {dJ_S(t)\over dt}  \equiv 0.
\eqno(6.2.8)
$$
Here we put $\<{\d J(u_t,t)\over \d u},{du_t\over dt}\>:=\sum_{j=1}^{2n}{\d 
J\over \d u},{du_j(t)\over dt}$. Therefore the vector 
$$
(v,\dot J_S,\dot J) \deff
\left({du_t\over dt}|_{t=0},
{dJ_S(t)\over dt}|_{t=0},
{dJ(t)\over dt}|_{t=0}\right)
$$ 
is tangent to $\calp $ at $(u_0,J_S(0),J_0)$ if and only if  
$$
{dv\over dx} + \left\<{\d J(u_t,t)\over \d u},v \right\>
\cdot {du_t\over dx}\cdot J_S + 
J(u_0,0)\cdot {dv\over dx}\cdot J_S(0) + 
$$
$$
+ \dot J\cdot {du\over dx}\cdot J_S + 
J\cdot {du\over dx}\cdot \dot J_S = O.\eqno(6.2.9)
$$
Note that $\<{\d J(u_t,t)\over \d u},v\>={dJ\over dv}$ is the derivative of 
$J$ along the vector $v$. Therefore (6.2.9) with $\nabla =d$ gives the same 
expression as (6.2.4).

\medskip\noindent\sl
6.3. $\dbar $-type Operators.

\smallskip\rm 
Now we need to understand the structure of the operator $D_{u,J}$ in more
detail. The problem arising here is that $D\deff D_{u,J}$ is only
$\rr$-linear. So we decompose it into $J$-linear and $J$-antilinear parts. 
For the rest of this text we shall denote by $E$ the pulled-back $u^*TX$ by 
$u$ tangent bundle of $X$. 
For $\xi \in C^1 (S, TS)$ and $v\in L^{1,p}(S, E)$ write $D_\xi v= {1\over2}[ 
D_{\xi }v - JD_{\xi }(Jv)] + {1\over2} [D_{\xi }v + JD_\xi(Jv)] = 
\dbar_{u, J}[v](\xi) + R(v,\xi )$.

\smallskip
\state Definition 6.3.1. The operator $\dbar_{u, J}$, introduced above 
as the $J$-linear part of $D_{u, J}$, we shall call the {\sl $\dbar$-operator
for a $J$-complex curve $u$}. 

\smallskip
\state Lemma 6.3.1. {\it $\dbar_{u,J}: L^{1,p}(S,E) L^p(S,\Lambda^{0,1}S
\otimes E)$ is a first order differential operator satisfying 
$$
\dbar_{u,J}(fv) = \dbar_{J_S}f\otimes v + f\otimes \dbar_{u,J}v .\eqno(6.3.1)
$$
where $\dbar_{J_S}f={1\over 2}[df + i\cdot df\circ J_S]$ is $\dbar $-operator 
on $S$. 
}

\state Proof. \rm Really, 
$$
2D_{u,J}(fv)=\nabla (fv)  + J\circ \nabla (fv)\circ 
J_S + (\nabla_{fv}J)\circ (du\circ J_S) = 
$$ 

$$
= f\nabla v + fJ\circ \nabla v\circ J_S + f\nabla_vJ\circ (du\circ J_S) + 
$$
$$
+ df\otimes v + Jdf\otimes vJ_S = f2D_{u,J}v + 2\dbar_{J_S}f\otimes v.
$$
From here we see that $2JD_{u,J}(Jfv) = 2fJD_{u,J}(Jv) + 2J\dbar f\otimes Jv$. So
$$
4\dbar_{u,J}(fv) = 2[D_{u,J}(fv) - JD_{u,J}(Jfv)] = 2f[D_{u,J}v-JD_{u,J}(Jv)] 
+ 4\dbar f\otimes v = 
$$
$$
= 4f\dbar_{u,J}v + 4\dbar_{J_S}f\otimes v .
$$
\qed

 \smallskip
 The following statement is well known in the smooth case. The line
 bundles case can be found in [Hf-L-Sk].
 
 \smallskip\noindent
 {\bf Lemma 6.3.2.} {\it Let $S$ be a Riemann surface with a complex
 structure $J_S$ and $E$ a $L^{1,p}$-smooth complex vector bundle of rank
 $r$ over $S$. Let also  $\dbar_E : L^{1,p}(S, E) \to
 L^p(S,\Lambda^{(0, 1)}S\otimes E)$ be a differential operator, satisfying
 the condition
 $$
 \dbar_E(f\xi) = \dbar_S f \otimes \xi
 + f\cdot \dbar_E\xi,  \eqno(6.3.2)
 $$
 where $\dbar_S$ is the Cauchy-Riemann operator, associated to $J_S$.
 Then the sheaf
 $$
 U\subset S \mapsto {\cal O}(E)(U) := \{\,\xi \in  L^{1,p}(U, E) \, :\,
 \dbar_E\xi=0 \,\}   \eqno(6.3.3)
 $$
 is analytic and locally free of rank $r$. This defines the holomorphic structure
 on $E$, for which $\dbar_E$ is an associated  Cauchy-Riemann operator.
 }
 
 \smallskip\noindent
 {\bf Remark.} The condition $(6.3.2)$ means that $\dbar_E$
 of order 1 and has the Cauchy-Riemann symbol.
 
 \smallskip\noindent
 {\bf Proof.} It is easy to see that the problem is essentially local. So we
 may assume that $S$ is a unit disk $\Delta$ with the standard complex
 structure. Let $\xi=(\xi_1,\ldots,\xi_r)$ be a $L^{1,p}$-frame of $E$
 over $\Delta$. Let also $\Gamma \in L^p(\Delta,\Lambda^{0, 1}\Delta
 \otimes Mat(r,\cc))$ be defined by relation
 $$
 \dbar_E \xi_i = \sum_j \Gamma_i^j \xi_j,
 \qquad \hbox{or in matrix form, } \qquad
 \dbar_E \xi =  \xi \cdot \Gamma.
 \eqno(6.3.4)
 $$
 Then for any section $\eta=\sum g^i \xi_i$ the equation $\dbar_E
 \eta =0$ is equivalent to
 $$
 \dbar g^i + \sum_j \Gamma^i_j g^j=0,
 \qquad \hbox{or in matrix form, } \qquad
 \dbar g + \Gamma \cdot g =0.
 \eqno(6.3.5)
 $$
 
 \smallskip
  Let the map $\tau_t : \Delta \to \Delta$ be defined by
 formula $\tau_t(z) = t\cdot z$, $0<t<1$. One can easily check that
 $$
 \Vert \tau_t^*\Gamma \Vert_{L^p(\Delta)}
 \le t^{1-2/p} \Vert \Gamma \Vert_{L^p(\Delta)}.
 $$
 So taking pull-backs $\tau_t^*E$ with $t$ sufficiently small we
 may assume that $\Vert \Gamma \Vert_{L^p(\Delta)}$ is small
 enough. Now consider the mapping $F$ from $L^{1,p}(\Delta,
  Gl(r,\cc))\subset $
  
  \noindent $ \subset L^{1,p}(\Delta,
 Mat(r,\cc))$ to $L^p(\Delta,\Lambda^{0, 1}\Delta
 \otimes Mat(r,\cc))$ defined by formula $F(g):=
 \dbar g \cdot g^{-1}$. It is easy to see that
 the derivation of $F$ in $g\equiv \id$ equals to $\dbar$.
 Due to the {\sl Lemma 3.2.1}  and the implicit function theorem
 for any $\Gamma$ with $\Vert \Gamma \Vert_{L^p(\Delta)}$ small
 enough there exists $g$ in $L^{1,p}(\Delta, Gl(r,\cc))$ with
 $g(0)=\id$,
 satisfying the equality $F(g)=-\Gamma$, which is equivalent to 
 (6.3.4). Consequently, in the neighbourhood of every point $p\in S$ there
 exists a frame $\eta =(\eta_1,\ldots,\eta_r)$ of $E$ consisting of sections
 of ${\cal O}(E)$. This implies that ${\cal O}(E)$ is analytic  and locally
 free of rank $r$.
 \qed
 \smallskip
 This lemma has an interesting corollary, which will be not used in this 
 notes.  
 Let $E\to S$ be real vector bundle over a Riemann 
 surface $(S,J_S)$, and let $J$ be a complex structure on $E$, \i.e. $J\in 
 End(E)$ with $J^2=-\id $. Take some connection $\nabla $ on $E$ and 
 consider a $J$-linear first order differential operator 
 $D^{\nabla }v:={1\over 2}[\nabla v - J\nabla (Jv)]:\Gamma^{1,p}(S,E)\to 
 \Gamma^p(S,\Lambda^1_S\otimes E)$. Define 
 
 $$
 \dbar^{\nabla }v = {1\over 2}(D^{\nabla }v - JD^{\nabla }_{J_S\cdot }v).
 \eqno(6.3.6)
 $$

 \smallskip
 \state Lemma 6.3.3. {\it  Operator $\dbar^{\nabla }$ acts from  
 $\Gamma^{1,p}(S,E)$ to $\Gamma^p(S,\Lambda^{0,1}_S\otimes E)$ and satisfies 
 condition (6.3.2) \i . e. 
  $$
 \dbar^{\nabla }(fv) = \dbar_S f \otimes v
 + f\cdot \dbar^{\nabla }v , \eqno(6.3.7)
 $$
 and therefore defines by  {\sl Lemma 6.3.2} a holomorphic structure on 
 $E$ compatible with a given complex structure $J$.
 }
 \state Proof. {\rm We need to check only the antilinearity of the form 
 $\dbar^{\nabla }v$and relation (6.3.7). Take a $\xi \in TS$, then 
 $$
 (\dbar^{\nabla }v)(J_S\xi ) = {1\over 2}[D^{\nabla }_{J_S\xi }v - 
 JD^{\nabla }_{\xi }v] = -J{1\over 2}[D^{\nabla }_{J\xi }v + 
 JD^{\nabla }_{J_S\xi }v] = -J\dbar^{\nabla }v(\xi ).
 $$
 This shows antilinearity of $\dbar^{\nabla }v$ with respect to $TS$-variable 
 $\xi $, and therefore  $\dbar^{\nabla }$ acts from $\Gamma^{1,p}(S,E)$ to 
 $\Gamma^p(S,\Lambda^{0,1}_S\otimes E)$. Note that $\dbar^{\nabla }$ is 
 $J$-linear on the variable $v$ as $D^{\nabla }$ is. Further,
 
 $$
 D^{\nabla }(fv) = {1\over 2}[df\otimes v + f\nabla v - Jf\nabla (Jv) - 
 Jdf\otimes Jv] = 
 $$
 
 $$
 = fD{\nabla }v + {1\over 2}(df\otimes v + df\otimes v) = fD^{\nabla }v + 
 df\otimes v.\eqno(6.3.8)
 $$
 Therefore 
$$
 \dbar^{\nabla }(fv) = {1\over 2}[D^{\nabla }(fv) - JD^{\nabla }_{J_S\cdot }
 (fv)] = {1\over }[fD^{\nabla }v + df\otimes v - fJD^{\nabla }_{J_S\cdot }v 
 - Jdf\otimes v] = 
$$
 
$$
 = f\dbar^{\nabla }v + {1\over 2}[df\otimes v -Jdf\otimes v] = 
 f\dbar^{\nabla }v + \dbar _{J_S}f\otimes v .
$$
This is because 
$$
 D^{\nabla }_{J_S\cdot }(fv) = {1\over 2}[\nabla_{J_S\cdot } (fv) - 
 J\nabla_{J_S\cdot }(fv)] = fD^{\nabla }_{J_S\cdot }v + {1\over 2}[df\circ J_S
 \otimes v - Jdf\circ J_S\otimes v] = 
$$
$$
 = fD^{\nabla }_{J_S\cdot }v + Jdf\circ J_S\otimes v =  
 fD^{\nabla }_{J_S\cdot }v  + \dbar_Sf\otimes v.
$$
We used here the definition $if\otimes v = f\otimes Jv$ of multiplication 
by $i$ in the bundle $E$. \qed

\bigskip\noindent
\sl 6.4. Holomorphic Structure on the Induced Bundle.

\smallskip\rm
Relation (6.3.1) tells us that the $J$-linear operator $\dbar_{u, J} : 
L^{1,p}(S, E) \to L^p(S, \allowbreak 
\Lambda^{0, 1}S\otimes E)$ satisfies the condition 
(6.3.2) and therefore defines in way explained in the previous 
paragraph a holomorphic structure on the bundle $E$. We shall denote by ${\cal
O}(E)$ the sheaf of holomorphic sections of $E$.

\smallskip Now let us turn to the antilinear part $R$ of Gromov's operator.

\medskip
\noindent
\bf Lemma 6.4.1. \it $R$ is a continuous $J$-antilinear operator from
$E$ to $\Lambda^{0,1}\otimes E$ of the order zero, satisfying
$$
 R(v,\xi ) = N(v, du(\xi )).\eqno(6.4.1) 
$$
and consequently 
$$
 R\scirc du \equiv 0.\eqno(6.4.2)
$$
Here the second relations means that for all $\xi , \eta \in TS$ we have 
$R(du(\eta ),\xi )=0$.

\noindent \rm
{\bf Proof.} $J$ anti-linearity of $R$ is given by its definition. Compute
$R(v,\xi )$ for $v\in L^{1,p}(S, E)$ and $\xi \in C^1(S, TS)$, setting
$w\deff du(\xi)$ and $D\deff D_{u, J}$ to simplify the notations

\smallskip
\noindent
$$
4R(v,\xi ) = 2D[v](\xi ) + 2JD[Jv](\xi ) =
$$

\smallskip
$$
=\nabla_{\xi }v + J\nabla_{J_S\xi }v + \nabla_vJ\scirc du\scirc J_S(\xi )
$$
$$
+J\nabla_{\xi }(Jv) + J^2\nabla_{J_S(\xi )}(Jv) +
J \nabla_{Jv}J\scirc du\scirc J_S(\xi ) =
$$

\smallskip
$$
= \nabla_\xi v + J\nabla_{J_S\xi }v + \nabla_vJ \scirc du\scirc J_S(\xi ) +
$$
$$
+ J^2\nabla_{\xi }v + J(\nabla_wJ)v + J^3\nabla_{J_S\xi }v +
J^2(\nabla_{Jw}J)v + J\scirc \nabla_{Jv}J\scirc Jw.
$$

\smallskip
Here we used the relations $\nabla_{du(\xi )}J = \nabla_wJ$, $du(J_S\xi) =Jw$
and $\nabla_{du(J_S\xi )}J = \nabla_{Jw}J$. Contracting terms,
we obtain 
$$
4R(v,\xi ) = \nabla_vJ(Jw) + J(\nabla_wJ)v - (\nabla_{Jw}J)v +
J(\nabla_{Jv}J(Jw)) =
$$
$$
= (\nabla_vJ\scirc J)w - (\nabla_wJ\scirc J)v - (\nabla_{Jw}J)v +
(\nabla_{Jv}J)w = 4N(v, w),
$$
where $N(v, w)$ denotes the torsion tensor of the almost-complex structure
$J$, see [Li], p.183, or [Ko-No], vol.II., p.123, where another normalization
constant for the almost complex torsion is used. Finally we obtain
$$
R(v,\xi ) = N(v, du(\xi )). 
$$

\smallskip
$N$ is antisymmetric and $J$-antilinear on both arguments, so
$$
-JR\bigl(du(\eta ),\xi \bigr) = R\bigl(du(\eta ), J_S\xi \bigr) =
N\bigl(du(\eta ), du(J_S\xi ) \bigr)
= N\bigl(du(\eta ), du(\eta )\bigr) =0
$$
if $\xi$ and $\eta$ were chosen in such a way that $J_S(\xi )=\eta $.
The relation $R(du(\xi ),\xi )=0$  obviously follows from (6.4.1).
\qed

\smallskip
\state Remarks. 1. {\rm If our structure $J$ is integrable, i.e. $N_J\equiv 0$ 
then by (6.4.2) $R\equiv 0$ and therefore $D_{u,J}=\dbar_{u,J}$ is $J$-linear.
In fact we know that $J$ is integrable iff there exists a symmetric connection 
$\nabla $ on $TX$ compatible with $J$, i.e. $\nabla J=0$, see [Ko-No]. From 
(6.2.5) we see that in this case $D_{u,J}(v)={1\over 2}(\nabla v + J\circ 
\nabla v\circ J_S)$. So $D_{u,J}$ is onviously $J$-linear.
}

\noindent\bf
2. \rm  The tangent bundle $TS$ to the Riemann surface  $S$ carries a
natural holomorphic structure. We shall denote by $\calo(TS)$
the~corresponding analytic sheaf. If the structure $J$ on $X$ is integrable 
then  the differential $du:(TS,\calo(TS))\to (S,E)$ of a holomorhpic map 
$u:S\to X$ is an analytic morphism of sheaves. We shall see know that this 
fact is stil true in nonintegrable case.

\medskip
\noindent
\bf Lemma 6.4.2. \it Let $u : (S, J_S) \to (X, J)$ be a non-constant
 complex  curve in almost complex manifold $X$. Then $du$ defines
an injective analytic morphism of analytic sheaves
$$
0\longrightarrow \calo(TS)
\buildrel du \over \longrightarrow \calo(E),  \eqno(6.4.3)
$$

\smallskip
\noindent
where $E = u^*(TX)$ is equipped with a holomorphic structure defined as above
by the operator $\dbar_{u, J}$.

\smallskip
\noindent
\rm {\bf Proof.} Injectivity of a sheaf homomorphism is equivalent to its
nondegeneracy, which is our case.

To prove holomorphicity of $du$ it is sufficient to show that for $\xi,\eta \in
C^1(S, TS)$ one has
$$
\bigl(\dbar_{u, J}(du(\xi ))\bigr)(\eta ) =
du\bigl((\dbar_S\xi ) (\eta )\bigr),  \eqno(6.4.4)
$$

\smallskip
\noindent
where $\dbar_S$ is the~usual $\dbar $-operator on $TS$. We shall use
the relation which is, in fact, the~definition for $\dbar_S$:

\smallskip
$$
(\dbar_S\xi )(\eta ) =
\msmall{1\over2}\left( \nabla_{\eta }\xi + J_S\nabla_{J_S\eta}\xi \right).
\eqno(6.4.5)
$$
\noindent Here $\nabla$ is a symmetric connection on $S$, compatible with
$J_S$. One has
$$
2\cdot (\dbar_{u, J}du(\xi ))(\eta ) = \nabla_{\eta }(du(\xi) ) +
J\nabla_{J_S\eta }(du(\xi))  + (\nabla_{du(\xi )}J)(du(J_S\eta )) =
$$
$$
= (\nabla_{\eta }du)(\xi ) + du(\nabla_{\eta }\xi ) +
J(\nabla_{J_S\eta}du)(\xi ) + J(du(\nabla_{J_S\eta }\xi ))
+ (\nabla_{du(\xi )}J)(du(J_S\eta )) =
$$
$$
= du\bigl[ \nabla_{\eta }\xi + J_S \nabla_{J_S\eta }\xi \bigr] +
\bigl[ (\nabla_{\eta }du)(\xi ) + J(\nabla_{J_S\eta}du)(\xi )
+ (\nabla_{du(\xi )}J)(du(J_S\eta )) \bigr]. \eqno(6.4.6)
$$
The first term of $(6.4.6)$ is $2\cdot du(\dbar_S\xi)(\eta)$.
To cancel the second one we use the identities $(\nabla_\xi du)[\eta] 
= (\nabla_\eta du)[\xi]$, $\nabla_wJ\scirc J = - J\scirc \nabla_w J$, and
$(\nabla_\xi du)\scirc J_S= J\scirc(\nabla_\xi du) + \nabla_{du(\xi)}J\scirc 
du$.
The last identity is obtained via co-variant differentiation of $du\scirc J_S=
J\scirc du$. Consequently, we obtain
$$
(\nabla_{\eta }du)(\xi ) + J(\nabla_{J_S\eta}du)(\xi ) +
(\nabla_{du(\xi )}J)(du(J_S\eta ))=
$$
$$
(\nabla_{\xi }du)(\eta ) + J(\nabla_\xi du)(J_S\eta ) +
(\nabla_{du(\xi )}J)(du(J_S\eta ))=
$$
$$
(\nabla_{\xi }du)(\eta ) +
J^2(\nabla_\xi du)(\eta ) +
J (\nabla_{du(\xi)}J)(du(\eta))
+ (\nabla_{du(\xi )}J)(du(J_S\eta) )=
$$
$$
=(J \scirc\nabla_{du\xi}J)(du\eta)
+ (\nabla_{du(\xi )}J\scirc J)(du(\eta ))=0.
$$
\qed

\smallskip\state Remark. 
We can give an alternative proof to both {\sl Lemmas 6.3.1} and
{\sl 6.3.2}, which does not use direct calculation. Fix a~complex structure
$J_S$ on $S$ and let $\phi_t$ be the~one
parameter group of diffeomorphisms of $S$, generated by a~vector field
$\xi$. Then ${d\over dt}|_{t=0}(\dbar\phi_t)=D_{J_S, \id }\xi=\dbar_S \xi$.
Let a~$J$-holomorphic map $u:S\to X$ also be fixed. Then $
{d\over dt}|_{t=0}(u\scirc \phi_t)=du(\xi)$ and consequently
$$
D_{J, u}(du(\xi))={d\over dt}\bigm|_{t=0}\dbar_J(u\scirc \phi_t)=
{d\over dt}\bigm|_{t=0}(du \scirc \dbar \phi_t)=
du\bigl({d\over dt}\bigm|_{t=0}(\dbar_{J_S} \phi_t)=
$$
$$
= du(\dbar_S\xi)
$$
or equivalently $D_{J, u}\scirc du=du \scirc \dbar$. Taking the $J$-antilinear
part of the last equality we obtain $R\scirc du=0$. 
Nevertheless we shall use the
explicit form of the Gromov operator:

$$
D_{u,J}(v)[\xi ] = \dbar_{u,J}(v)[\xi ] + N_J(v,du(\xi )).\eqno(6.4.6)
$$

\bigskip
The zeroes of analytic morphism $du : \calo(TS) \to \calo
(E)$ are isolated. So we have the following

\medskip
\noindent
\bf Corollary 6.4.3. \rm ([Sk]). \it The set of critical points of a
 complex  curve in almost the complex manifold $(X, J)$ is discrete,
provided $J$ is of class $C^1$.

\smallskip\rm
For $C^\infty$-structures the result is due to McDuff, see [McD-1].
%%%%%%%%%%%%%%%%%Reference??%%%%%%%%%%%%%%%%%%
\newpage

\bigskip\noindent
{\bigbf Lecture 7}

\smallskip\noindent
{\bigbf Fredholm Properties of the Gromov Operator}

\medskip\noindent\sl
7.1. Generalized Normal Bundle.

\smallskip\rm
{\sl Lemma 6.3.2} makes it possible to define the order of vanishing of the
differential of $(J_S,J)$-holomorphic map $u:(S,J_S)\to (X,J)$, provided
$J\in C^1$.

\smallskip
\noindent\bf
Definition 7.1.1. \sl By the order of zero ${\sf ord}_p du$ of the differential
$du$ at a point $p\in S$, we  understand the order of vanishing at $p$ of
the holomorphic morphism $du : \calo(TS)\to \calo(E)$.
\rm

\smallskip
From (6.3.2) we obtain the following short exact sequence:
$$
0\longrightarrow \calo(TS) \buildrel du \over\longrightarrow \calo(E)
\longrightarrow \caln\longrightarrow 0. \eqno(7.1.1)
$$

\smallskip
Here $\caln$ is a quotient-sheaf $\calo(E)/du(TS)$. We can decompose
$\caln = \calo(N_0)\oplus \caln_1$, where $N_0$ is a holomorphic
vector bundle and $\caln_1 =
\bigoplus_{i=1}^P \cc_{a_i}^{n_i}$. Here $\cc_{a_i}^{n_i}$ denotes the sheaf,
supported at the~critical points $a_i\in S$ of $du$ and having a stalk
$\cc^{n_i}$ with $n_i = \ord_{a_i}du$,
the~order of zero of $du$ at $a_i$. We shall call $\caln$ a normal
sheaf and $N_0$ a generalized normal bundle.

Denote by $[A]$ the divisor $\sum_{i=1}^Pn_i[a_i]$, and by $\calo([A])$ a
sheaf of meromorphic functions on $S$ having poles in $a_i$ of order at most
$n_i$. Then (7.1.1) gives rise to the exact sequence

\smallskip
$$
0\longrightarrow \calo(TS)\otimes \calo([A])\buildrel{du}\over
{\longrightarrow} \calo(E)
\longrightarrow \calo(N_0)\longrightarrow 0. \eqno(7.1.2)
$$

\smallskip
Denote by $L^p_{(0, 1)}(S, E)$ the space of $L^p$-integrable (0, 1)-forms
with coefficients in $E$. Then (7.1.2) together with {\sl Lemma 6.3.1}
implies that the following diagram is commutative

\smallskip
$$

\def\mapright#1{\smash{\mathop{\longrightarrow}\limits^{#1}}}
\def\mapdown#1{\Big\downarrow\rlap{$\vcenter{\hbox{$\scriptstyle#1$}}$}}
\matrix
0&\mapright{}&L^{1,p}(S, TS\otimes [A])&\mapright{du}&
L^{1,p}(S, E)&\mapright{\pr}&L^{1,p}(S, N_0)&\mapright{}&0\cr
& &\mapdown{\dbar_S}& &\mapdown{D_{u, J}}& &\mapdown{}& & \cr
0&\mapright{}&L^p_{(0,1)}(S, TS\otimes [A])&\mapright{du}&
L^p_{(0,1)}(S, E)&\mapright{}&
L^p_{(0,1)}(S, N_0)&\mapright{}&0. \cr
\endmatrix
\eqno(7.1.3)
$$

\smallskip
This defines an operator $D_{u, J}^N : L^{1,p}(S, N_0) \longrightarrow
L^p_{(0,1)}(S, N_0)$ which has the form $D_{u, J}^N =
\dbar_N + R$. Here $\dbar_N$ is a usual $\dbar $-operator on $N_0$
and $R\in C^0(S, \hom_\rr(N_0,\Lambda^{0,1} \otimes N_0))$. This
follows from the fact that $D_{u, J}$ has the same form.

\medskip
\state Definition 7.1.2. Let $E$ be a holomorphic vector bundle over
a~compact Riemann surface $S$ and let $D:L^{1,p}(S, E)\to L^p(S,
\Lambda^{0,1}S \otimes E)$ be an operator of the~form $D=\dbar + R$, where
$R\in L^p\bigl(S,\,\hom_\rr(E,\,\Lambda^{0, 1}S\otimes E) \bigr)$
with $2<p<\infty$. Define $\sfh^0_D(S, E)\deff \ker D$ and $\sfh^1_D(S, E)
\deff \coker D$.

\smallskip
\state Remark. It is shown in {\sl Lemma 7.2.2} below that given $S$, $E$
and $R\in L^p$, $2<p<\infty$, one can define $\sfh^i_D(S, E)$
as a (co)kernel of the operator $\dbar +R: L^{1,q}(S, E) \to
L^q(S,\,\Lambda^{0, 1}S\otimes E)$ for any $1<q\le p$. Thus, the definition
is independent of the choice of a functional space.

By the standard lemma of homological algebra we obtain from (7.1.2) 
the following
long exact sequence of $D$-cohomologies.

\smallskip
$$

\def\mapright#1{\smash{\mathop{\longrightarrow}\limits^{#1}}}
\def\mapdown{\Big\downarrow}
\matrix
0& \mapright{}& \sfh^0(S, TS\otimes [A]) &\mapright{} & \sfh^0_D(S, E)
& \mapright{}& \sfh^0_D(S, N_0)   &\mapright{\delta} &\cr
\vphantom{\mapdown}&&&&&&&&\cr
& \mapright{}& \sfh^1(S, TS\otimes [A]) &\mapright{} & \sfh^1_D(S, E)
& \mapright{}& \sfh^1_D(S, N_0)         &\mapright{} &0.
\endmatrix
\eqno(7.1.4)
$$

\bigskip\smallskip
\noindent \sl 7.2. Surjectivity of $D_{u, J}$.

\nobreak\smallskip\rm
We shall use a result of
Gromov ([G]) and Hofer-Lizan-Sikorav ([Hf-L-Sk]) about
surjectivity of  $D^N_{u, J}$, namely a vanishing theorem
for $D$-cohomologies. First we prove some technical statements.

\medskip
\state Lemma 7.2.1. {\it Let $X$ and $Y$ be Banach spaces and $T:X\to Y$
a~closed dense defined unbounded operator with the~graph $\Gamma=\Gamma_T$
endowed with the~graph norm $\Vert x \Vert_\Gamma= \Vert x \Vert_X + \Vert Tx
\Vert_Y$. Suppose that the~natural map $\Gamma\to X$ is compact. Then

\sli $\ker(T)$ is finite-dimensional;

\slii $\im(T)$ is closed;

\sliii the~dual space $\bigl(Y/\im(T)\bigr)^*$ is naturally isomorphic
to $\ker(T^*:Y^*\to X^*)$.
}

\state Proof. Obviously, for $x\in \ker(T)$ one has $\Vert x \Vert_X =
\Vert x \Vert_\Gamma$. Let $\{x_n\}$ be a sequence in $\ker(T)$ which is
bounded in the $\Vert \cdot \Vert_X$- norm. Then it is bounded in the 
$\Vert \cdot
\Vert_\Gamma$- norm and hence relatively compact with respect to the
$\Vert \cdot \Vert_X$- norm. Thus, the~unit ball in $\ker(T)$ is compact
which implies the~statement \sli of the~lemma.

Due to finite-dimensionality, there exists a~closed complement $X_0$ to
$\ker(T)$ in $X$. Let $\{x_n\}$ be a~sequence in $X$ such that $Tx_n \lrar y
\in Y$. Without losing generality we may assume that $x_n$ belong to $X_0$.
Suppose that $\Vert x_n \Vert_X \lrar \infty$. Denote $\tilde x_n \deff {x_n
\over \Vert x_n \Vert_X }$. Then $\Vert \tilde x_n \Vert_\Gamma$ is bounded
and hence some subsequence of $\{\tilde x_n\}$, still denoted by $\{\tilde
x_n\}$, converges in $X_0$ to some $\tilde x$. Note that $\tilde x\not=0$,
because $\Vert \tilde x \Vert_X = \lim \Vert \tilde x_n \Vert_ X =1$. On
the~other hand, one can see that $T\tilde x_n \lrar 0\in Y$. Since $\Gamma$ is
closed, $(\tilde x, 0)\in \Gamma$ and hence $\tilde x\in \ker(T) \cap X_0
=\{0\}$. The~contradiction shows that the~sequence $\{x_n\}$ must be bounded
in $X$. Since $\{Tx_n \}$ is also bounded in $Y$, some subsequence of
$\{x_n\}$, still denoted by $\{x_n\}$, converges in $X_0$ to some $x$. Due to
the~closeness of $\Gamma$, $Tx=y$. Thus $\im (T)$ is closed in $Y$.

Denote $Z\deff\ker(T^*:Y^*\to X^*)$ and let $h\in \bigl(Y/\im(T)
\bigr)^*$. Then $h$ defines a~linear functional on $Y$,\i.e., some element
$h'\in Y^*$, which is identically zero on $\im(T)$. Thus, for any $x$ from
the~domain of the definition of $T$ one has $\langle h', Tx \rangle =0$, which
implies $T^*(h')=0$. Consequently, $h'$ belongs to $Z$. Conversely, every
$h'\in Z$ is a~linear functional on $Y$ with $h'(Tx)= \langle T^*h', x\rangle
=0$ for every $x$ from the~domain of the definition of $T$. Thus, $h'$ is
identically zero on $\im(T)$ and is defined by some unique $h\in
\bigl(Y/\im(T)\bigr)^*$.\qed

\bigskip
\state Lemma 7.2.2. {\sl (Serre Duality for $D$-cohomologies.)}
{\it Let $E$ be a holomorphic vector bundle over a~compact Riemann
surface $S$,
and let $D:L^{1,p}(S, E)\to L^p(S, \Lambda^{0, 1}S\otimes E)$ be
an operator of the~form $D=\dbar + R$, where $R\in L^p\bigl(S,\,
\hom_\rr(E,\Lambda^{0, 1}S\otimes E) \bigr)$ with $2<p <\infty$. Also let
$K\deff \Lambda^{1, 0}S$ be the~canonical holomorphic line
bundle of $S$. Then there exists the~naturally defined operator
$$
D^*=\dbar- R^* : L^{1,p}(S, E^* \otimes K)
\to L^p(S,\Lambda^{0, 1} \otimes E^* \otimes K)
$$
with $R^* \in L^p\bigl(S,\,\homr(E^*\otimes K,\,
\Lambda^{0, 1}S \otimes E^*\otimes K) \bigr)$ and the~natural isomorphisms
$$
\sfh^0_D(S,\, E)^*\cong \sfh^1_{D^*}(S,\, E^*\otimes K),
$$
$$
\sfh^1_D(S,\, E)^*\cong \sfh^0_{D^*}(S,\, E^*\otimes K).
$$
If, in addition, $R$ is $\cc$-antilinear, then $R^*$ is also
$\cc$-anti-linear.
}

\medskip
\state Proof. For any $1<q\le p$ we associate with $D$ an~unbounded
dense defined operator $T_q$
from $X_q\deff L^q(S, E)$ into $Y_q\deff L^q(S,\Lambda^{0, 1}S\otimes E)$
with the~domain of definition $L^{1, q}(S, E)$.
The elliptic regularity of $D$ (see {\sl Lemma 3.2.1} above ) implies that
$$
\Vert \xi \Vert_{L^{1,q}(S, E)} \le C(q)
\left(
\Vert \dbar\xi +R\xi\Vert_{L^q(S, E)} +
\Vert \xi \Vert_{L^q(S, E)} \right).
$$
Consequently, $T_q$ are closed and satisfy the~hypothesis of
{\sl Lemma 7.2.1}. For $q>q_1$ we also have the~natural imbedding $X_q
\hookrightarrow X_{q_1}$ and $Y_q \hookrightarrow Y_{q_1}$ which commutes with
the~operator $D$. Moreover, due to the regularity  of $D$
this imbedding maps $\ker T_q$ {\sl identically} onto $\ker T_{q_1}$.
Thus we can identify $\sfh^0_D(S, E)$ with any $\ker T_q$.

Now note that for $q'\deff q/(q-1)$ we have the~natural isomorphisms
$$
\eqalign{
X_q^*\equiv &(L^q(S, E))^*
\cong L^{q'}(S,\Lambda^{0, 1}S \otimes E^* \otimes K),\cr
Y_q^*\equiv &(L^q(S,\Lambda^{0, 1}S\otimes E))^*
\cong L^{q'}(S,  E^* \otimes K),\cr
}
$$
induced by the pairing of $E$ with $E^*$ and by integration over $S$. One can
easily check that the~dual operator $T_q^*$ is induced by the~differential
operator $-D^*: L^{1,q'}(S, E^* \otimes K)
\to L^{q'}(S,\Lambda^{0, 1} \otimes E^* \otimes K)$
of the~form $D^* =\dbar - R^*$. In fact, for $\xi \in L^{1, q}(S, E)$
and $\eta\in L^{1, q'}(S, E^*\otimes K)$ one has
$$
\langle T_q\xi,\,\eta \rangle=
\int_S \langle \dbar\xi + R\xi,\,\eta \rangle =
\int_S \dbar\langle \xi ,\,\eta \rangle +
\int_S \langle \xi,\, -(\dbar - R^*)\eta \rangle =
\int_S \langle \xi ,\, -D^*\eta \rangle,
$$
since the~integral of any $\dbar$-exact $(1, 1)$-form vanishes. From {\sl
Lemma 7.2.1} we obtain the~natural isomorphisms $\sfh^1_D(S, E)^*\equiv (\coker
T_p)^*\cong \ker T_p^*$ and $\sfh^0_D(S, E)^*\equiv (\ker
T_p)^*\cong \coker T_p^*$,  which yields the~statement of the lemma. \qed

\medskip
\state Corollary 7.2.3. ([G], [H-L-Sk].) {\sl (Vanishing Theorem for
$D$-cohomologies.)} {\it Let $S$ be a~Riemann surface $S$ of the~genus $g$.
Also let $L$ be a~holomorphic {\sl line} bundle over $S$, equipped with
a~differential operator $D=\dbar + R$ with $R\in L^p\bigr(S,
\homr(L,\,  \Lambda^{0, 1}S\otimes L) \bigl)$, $p>2$. If $c_1(L)<0$,
then $\sfh^0_D (S,\, L)=0$. If $c_1(L)>2g-2$, then $\sfh^1_D (S,\, L)=0$.
}

\medskip
\state Proof. Suppose $\xi$ is a~nontrivial $L^{1,p}$-section of $L$ satisfying
$D\xi=0$. Then due to {\sl Lemma 3.1.1}, $\xi$ has only finitely many zeros
$p_i\in S$ with {\sl positive} multiplicities $\mu_i$. One can easily see
that $c_1(L)=\sum \mu_i \ge 0$. Consequently $\sfh^0_D (S,\, L)$ vanishes
if $c_1(L)<0$. The~vanishing result for $\sfh^1_D$ is obtained via the~Serre
duality of {\sl Lemma 7.2.2}.

\qed

\bigskip\bigskip
\noindent \sl 7.3. Tangent Space to the Moduli Space.

\nobreak\medskip \rm
Recall that in (7.1.4) we obtained the following long exact sequence

\smallskip
$$

\def\mapright#1{\smash{\mathop{\longrightarrow}\limits^{#1}}}
\def\mapdown{\Big\downarrow}
\matrix
0& \mapright{}& \sfh^0(S, TS\otimes [A]) &\mapright{} & \sfh^0_D(S, E)
& \mapright{}& \sfh^0_D(S, N_0)   &\mapright{\delta} &\cr
\vphantom{.}&&&&&&&&\cr
& \mapright{}& \sfh^1(S, TS\otimes [A]) &\mapright{} & \sfh^1_D(S, E)
& \mapright{}& \sfh^1_D(S, N_0)         &\mapright{} &0.
\endmatrix
$$
It is most important for us to associate a similar long exact
sequence of $D$-cohomo\-logies to the short exact sequence (7.1.1). Note
that, due to {\sl Lemmas 6.3.1} and 6.3.2, we obtain
the short exact sequence of complexes
$$

\setbox1=\hbox{$\lrar$}
\def\mapright#1{\,\,\smash{\mathop{{\hbox to
\wd1{\hss\hbox{$\displaystyle\longrightarrow$}\hss}}}\limits^{#1}}\,\,}
\def\mapdown#1{\Big\downarrow\rlap{$\vcenter{\hbox{$\scriptstyle#1$}}$}}
\matrix\format\c&\c&\c&\c&\c&\c&\c&\c&\l\\
0&\mapright{}&L^{1,p}(S, TS)&\mapright{du}&
L^{1,p}(S, E)&\mapright{\barr\pr}&
L^{1,p}(S, E){\bigm/}du(L^{1,p}(S, TS))
&\mapright{}&0\cr
& &\mapdown{\dbar_S}& &\mapdown{D}& &\mapdown{\barr D}& &
(7.3.1)%\eqno(7.3.1)
\cr
0&\mapright{}&L^p_{(0,1)}(S, TS)&\mapright{du}&
L^p_{(0,1)}(S, E)&\mapright{}&
L^p_{(0,1)}(S, E){\bigm/}
du(L^p_{(0,1)}(S, TS))&\mapright{}&0\cr\vphantom{.}
\endmatrix
%\eqno(7.3.1)
$$

\smallskip\noindent
where $\overline D$ is induced by $D\equiv D_{u, J}$.

\medskip\noindent
{\bf Theorem 7.3.1.} {\it For $\overline D$ just defined,
$\ker \overline D=\sfh^0_D(S, N_0) \oplus \sfh^0(S, \caln_1)$ and
$\coker \overline D=\sfh^1_D(S, N_0)$.
}

\smallskip\noindent
\bf Proof. \rm Let $\pr_0: \calo(E) \to \calo(N_0)$ and
$\pr_1: \calo(E) \to \caln_1$ denote the natural projections induced
by $\pr: \calo(E) \to \caln\equiv \calo(N_0) \oplus \caln_1$.
Also let $A$, as in (2.2.7), denote the support of $\caln_1$, {\sl \i.e.,} the
finite set of the vanishing points of $du$. Then $\pr_0$ defines maps
$$
\matrix
\pr_0:L^{1,p}(S , E)  &  \longrightarrow
& L^{1,p}(S , N_0)
\cr
\vphantom{.}\cr
\pr_0:L^p_{(0, 1)}(S , E)  &  \longrightarrow
& L^p_{(0, 1)}(S, N_0).
\cr
\endmatrix
\eqno(7.3.2)
$$
Furthermore, in a neighborhood of every point $p\in A$ the sequence (7.1.1)
can be represented in the form
$$
0\longrightarrow \calo
\buildrel \alpha_p \over\longrightarrow \calo^{n-1}
\oplus \calo
\buildrel \beta_p \over\longrightarrow \calo^{n-1}
\oplus\caln_1\vert_p
\longrightarrow 0 \eqno(7.3.4)
$$
with $\alpha_p(\xi)= (0, z^{\nu_p} \xi)$ and $\beta(\xi,\eta)=
(\xi, j^{(\nu-1)}_p\eta)$. Here $z$ denotes local a holomorphic coordinate
on $S$ with $z(p)=0$, $\nu_p$ is the multiplicity of $du$ in $p$,
$j^{(\nu-1)}_p\eta$ is a $(\nu-1)$-jet of $\eta$ in $p$ and $\caln_1
\vert_p$ is a stalk of $\caln_1$ in $p$.

Now let $\overline\xi\in L^{1,p}(S, E)/du(L^{1,p}(S, TS))$ which
satisfies $\overline D(\overline\xi)=0$. This means that $\overline\xi$ is
represented by some $\xi \in L^{1,p}(S, E)$ with $D\xi= du(\eta)$ for some
$\eta \in L^p_{(0, 1)}(S, TS)$. It is obvious that there exists $\zeta
\in L^{1,p}(S, TS)$ such that $\overline\partial \zeta = \eta$ in a
neighborhood of $A$. Consequently, $D(\xi- du(\zeta))=0$ in a neighborhood
of $A$. Denote $\xi_1 := \xi -du(\zeta)$. Due to (6.3.1), in a neighborhood
of $p\in A$ the equation $D\xi_1=0$ is equivalent to
$$
\overline\partial \xi_1+ N(\xi_1, du)=0.
$$
Due to {\sl Lemma 3.1.2} above, $\xi_1= P(z)+o(|z|^{\nu_p} )$ with some
(holomorphic) polynomial $P(z)$. This gives the possibility of  defining
$\varphi_0 := \pr_0(\xi) \in L^{1,p}(S, N_0)$ and $\varphi_1 :=
\pr_1(\xi- du(\zeta)) \in \sfh^0(S, \caln_1)$. Due to (7.1.2) and
(7.1.3), $D_N\varphi_0=0$. If $\overline\partial \zeta' = \eta$
in a neighborhood of $A$ for some other $\zeta'\in L^{1,p}(S, TS)$,
then $\zeta - \zeta'$ is holomorphic in a neighborhood of $A$, and
consequently $\pr_1(du(\zeta - \zeta'))=0$. Thus, the map
$\iota^0: \ker \overline D \to \sfh^0_D(S, N_0) \oplus \sfh^0(S, \caln_1)$,
$\iota^0(\overline\xi)=(\varphi_0,\varphi_1)$ is well-defined.

\smallskip
Assume that $\iota^0(\overline\xi)=0$ for some $\overline\xi\in
\ker\overline D$ and that $\overline\xi$ is represented by $\xi\in
L^{1,p}(S, E)$ with $D\xi= du(\eta)$ for some $\eta \in
L^p_{(0, 1)}(S, TS)$. Let $\zeta \in L^{1,p}(S, TS)$ satisfy
$\overline\partial \zeta = \eta$ in a neighborhood of $A$. The assumption
$\iota^0(\overline\xi)=0$ implies that $\pr_1(\xi- du(\zeta))=0$ and that
$\pr_0(\xi- du(\zeta))=0$ in a neighborhood of $A$. Consequently, $\xi-
du(\zeta)=du(\psi)$ for some $\psi \in L^{1,p}(S, TS)$ and $\xi \in
du(L^{1,p}(S, TS))$. This means that $\overline\xi=0\in \ker\overline D$
and $\iota^0$ is injective.

\smallskip
Let $\varphi_0 \in \sfh^0_D(S, N_0)$ and $\varphi_1 \in \sfh^0(S, {\cal
N}_1)$. For every $p\in A$ fix a neighborhood $U_p$ and representation of
(7.1.1) in the form (7.3.4) over $U_p$. In every $U_p$ we find $\xi=(\xi_0,
\xi_1) \in L^{1,p}(U_p,\cc^{n-1} \times\cc)$ satisfying the following
properties:

\smallskip
{\sl a)} $D\xi=0$;

\smallskip
{\sl b)} $\xi_0$ coincide with $\varphi_0\vert_{U_p}$ under
the identification $\calo(N_0)\vert_{U_p} \cong \calo^{n-1}$;

\smallskip
{\sl c)} $j^{(\nu-1)}_p \xi = \varphi_1\vert_p \in
\caln_1 \vert_p$.

\smallskip
The corresponding $\xi_1 \in L^{1,p}(U_p,\cc)$ can be
constructed as follows. Let $D(\xi_0,\xi_1) = (\eta_0,\eta_1)$. From
$\pr_0(D\xi)= D_N(\pr_0\xi)=0$ one has $\eta_0=0$. In the representation
(7.3.4) the identity $R\scirc du =0$ of {\sl Lemma 6.3.1} means that
$D(0,\xi_1) = (0,\dbar\xi_1)$. Thus, one can find $\xi_1$ with $\dbar\xi_1
=-\eta_1$, which gives $D(\xi_0,\xi_1) =0$. Consequently, $j^{(\nu-1)}_p
\xi$ is well-defined. Adding an appropriate {\sl holomorphic} term to
$\xi_1$ one can satisfy condition {\sl c)}.
Using an appropriate partition of unity, we can construct $\xi'
\in L^{1,p}(S, E)$ such that $\xi'$ coincides with $\xi$ in a
(possibly smaller) neighborhood of every $p\in A$ and such that $\pr_0 \xi'
= \varphi_0$ in $S$. Thus from $D_N\varphi_0=0$ and (2.2.8) one
obtains $D\xi' \in du(L^p_{(0, 1)}(S, TS \otimes [A]))$. But $D\xi'
=0$ in a neighborhood of every $p\in A$, which means that $D\xi' \in du(
L^p_{(0, 1)}(S, TS))$. This shows surjectivity of $\iota^0: \ker
\overline D \to \sfh^0_D(S, N_0) \oplus \sfh^0(S, \caln_1)$.

\medskip
Now consider the case of $\coker \overline D$. Since the operator
$\pr_0$ satisfies the identities $\pr_0\scirc D_E=D_N\scirc \pr_0$ and
$\pr_0\scirc du=0$, the induced map

\smallskip$\displaystyle
\iota^1 : \coker \overline D\equiv L^p_{(0,1)}(S, E)\bigm/
\bigl(D_E(L^{1,p}(S, E) \oplus  du(L^p_{(0,1)}(S, TS)\bigr)
\longrightarrow $

\smallskip \hfill$\displaystyle
\longrightarrow\sfh^1_D(S, N_0)\equiv
L^p_{(0,1)}(S, N_0)/ D_N(L^{1,p}(S, N_0))
\qquad$

\smallskip\noindent
is well-defined. Moreover, the surjectivity of $\pr_0:
L^p_{(0,1)}(S, E) \to L^p_{(0,1)}(S, N_0)$
easily yields the surjectivity of $\iota^1$.

Now assume that $\iota^1 \overline\xi =0$ for some $\overline\xi \in
\coker \overline D$ and  $\overline \xi$ is represented by $\xi
\in  L^p_{(0,1)}(S, E)$. Then the condition $\iota^1
\overline\xi=0$ means that $\pr_0 \xi= D_N \eta$ for some $\eta
\in L^{1,p}(S, N_0)$. Find $\zeta \in L^{1,p}(S, E)$ such that
$\pr_0 \zeta =\eta$. Then $\pr_0 (\xi - D\zeta)=0$, and, due to
(7.1.3), $\xi- D\zeta=du(\varphi)$ for some $\varphi \in
L^p_{(0,1)}(S, TS\otimes [A])$. Find $\psi \in
L^{1,p}(S, TS\otimes [A])$ such that $\dbar \psi =\varphi$
in some neighborhood of $A$. Then
$$
\xi- D\zeta - D(du(\psi))= du(\varphi-\dbar\psi) \in
du(L^p_{(0,1)}(S, TS)).
$$
Consequently, $\overline \xi \in \im\overline D$. This shows the injectivity
of $\iota^1$.\qed

\bigskip
\state Corollary 7.3.2. \it The~short exact sequence $(7.1.1)$ induces
the long exact sequence of $D$-cohomologies
$$

\def\mapright#1{\,\,\smash{\mathop{\longrightarrow}\limits^{#1}}\,\,}
\def\mapdown{\Big\downarrow}
\matrix\format\c&\c&\c&\c&\c&\c&\c&\c&\c\\
0& \mapright{}& \sfh^0(S, TS) &\mapright{} & \sfh^0_D(S, E)
& \mapright{}& \sfh^0_D(S, N_0)\oplus \sfh^0(S, \caln_1)
&\mapright{\delta} &\cr
%
%\vphantom{.}&&&&&&&&\cr
%
& \mapright{}& \sfh^1(S, TS) &\mapright{} & \sfh^1_D(S, E)
& \mapright{}& \sfh^1_D(S, N_0)         &\mapright{} &0.
\endmatrix
$$

\bigskip\noindent\sl
7.4. Reparameterizations.

\smallskip\rm

Let $S$ be an oriented compact real surface without boundary. Consider
the Teichm\"uller space $\ttt_g$  of marked complex structures on $S$.
This is a complex manifold of a dimension
$$
\dimc \ttt_g = \cases  0      &\text{ if $g=0$;}  \\
1      &\text{ if $g=1$;}  \\
3g-3   &\text{ if $g\ge2$}
\endcases
$$
which can be completely characterized in the following way. The product
$S\times \ttt_g$ possesses a complex (\ie holomorphic) structure $J_{S\times
\ttt}$ such that

\sli the natural projection $\pi\mid_{\ttt }: S\times \ttt_g \to \ttt_g$ is
holomorphic, so that for any $\tau\in \ttt_g$ the identification $S \cong S
\times \{\tau\}$ induces the complex structure $J_S(\tau)\deff J_{S\times
\ttt}\ogran_{S\times \{\tau\}}$ on $S$;

\slii for any complex structure $J_S$ on $S$ there exists a uniquely
defined $\tau\in\ttt_g$ and a diffeomorphism $f: S\to S$ such that $J_S =
f^*J_S(\tau)$ (\ie $f: (S, J_S) \to (S, J_S(\tau))$ is holomorphic) and
$f$ is homotopic to the identity map $\id_S: S\to S$.

Denote the automorphism group of $S\times \ttt_g$ by $\bfg$. It is known that
$$
\bfg=\cases {\bold PGl}(2,\cc)       &\text{ for $g=0$,} \\
{\bold Sp}(2,\;\;\zz) \ltimes T^2   &\text{ for $g=1$,} \\
{\bold Sp}(2g,\zz)                  &\text{ for $g\ge2$.}
		\endcases
$$
We shall use the following information about $\ttt_g$.

For $g=0$ the surface $S$ is a Riemann sphere $S^2$ and all complex
structures on $S^2$ are equivalent to the standard one when $S^2\cong
\cc\pp^1$. Thus, $\ttt_0$ consists of one point, and the group $\bfg={\bold
PGl}(2,\cc)$ is the group of biholomorphisms of $\cc\pp^1$.

For $g=1$ the surface $S$ is a torus $T^2$ and $\ttt_1$ is an upper half-plane
$\cc_+=\{\,\tau\in\cc\;:\;\im(\tau)>0\,\}$. The product $S\times \ttt_1$ can
be identified by the quotient $(\cc\times \cc_+)/\zz^2$ with respect to the
holomorphic action
$$
\bigl((n,m),\, (z,\tau)\bigr) \in\zz^2 \times \cc\times \cc_+ \longmapsto
(m,n)\cdot(z,\tau)\deff (z+m +n\,\tau,\tau) \in \cc\times \cc_+.
$$
In this case the subgroup $T^2\equiv \rr^2/\zz^2$ is a connected component of
the identity ${\rom e}\in\bfg$, in particular $T^2$ is normal. The group
$\bfg= {\bold Sp}(2,\;\; \zz) \ltimes T^2$ is a semi-direct product. The
holomorphic action of $T^2$ on the quotient $(\cc\times \cc_+)/\zz^2$ is
given by
$$
\eqalign{
\bigl([t_1,t_2],\, ([z],\tau)\bigr) &\in T^2 \times (\cc\times \cc_+)/\zz^2
\longmapsto
\cr
[t_1,t_2]\cdot([z],\tau)&\deff ([z+t_1 +t_2\tau],\tau )
\in (\cc\times \cc_+)/\zz^2.
}
$$

For any $g\ge0$ the action of $\bfg$ on $S\times \ttt_g$ is effective and
preserves fibers of the projection $\pi\mid_{\ttt }: S\times \ttt_g \to 
\ttt_g$. This induces the action of $\bfg$ on $\ttt_g$. Furthermore, given 
$\tau \in \ttt_g$ and $f\in\bfg$, there exists a unique diffeomorphism of 
$\hat f_\tau:S \to S$ such that
$$
f\cdot(x,\tau)=(\hat f_\tau(x), f\cdot\tau).
\eqno(7.4.1)
$$

For any $\tau\in\ttt_g$ we have natural isomorphisms $T_\tau\ttt_g \cong
\sfh^1(S,\, \calo_\tau(TS))$ and $T_{\rom e}\bfg \cong \sfh^0(S,\,
\calo_\tau(TS))$, where $\calo_\tau(TS)$ denote the sheaf of a section of $TS$
which are holomorphic with respect to the complex structure $J_S(\tau)$.

Later on, we shall denote elements of $\ttt_g$ by $J_S$ and consider them
as corresponding complex structures on $S$.

%%%%%%%%%%%%%%%%%%%%%page120%%%%%%%%%%%%%%%%%%%%%%%%%

\newpage
\noindent{\bigbf Lecture 8 }

\smallskip\noindent{\bigbf  Transversality.}

\medskip\noindent\sl
8.1. Moduli Space of Nonparameterized Curves.

\smallskip\rm
Let $(X,\omega, J\st)$ be a symplectic manifold with some distinguished
$w$-tamed almost-complex structure $J\st$. In our applications $J\st $ will
be integrable, providing, together with $w$, the K\"ahler structure on $X$.
Let $U\Subset X$ be an open relatively
compact subset which can coincide with $X$ if $X$ is compact. Also let $S$ be
a (fixed) compact oriented surface of genus
$g\ge0$, $u_0: S \to X$ a non-constant $C^1$-smooth map such that $u_0(S)\cap
U \not = \emptyset$.

\medskip
Fix $2<p<\infty$. The Banach manifold $L^{1,p}(S,X)$ of all $L^{1,p}
$-smooth maps $u:S \to X$ is smooth with a tangent space $T_uL^{1,p}(S,X)$ equal
to the space $L^{1,p}(S,u^*TX)$ of $L^{1,p}$-smooth sections of a pulled-back
tangent bundle to $X$. Denote by $\cals_U$ a Banach manifold of those
$u\in L^{1,p}(S,X) $, for which $u$ is homotopic to $u_0$
and $u(S) \cap U \not=\emptyset$. Fix an integer $k\ge1$ and denote by
$\calj^k_U$ the set of $C^k$-smooth almost complex structures $J$ in $X$ such
that $\{ x\in X:J(x)\not=J\st(x) \} \comp U$, and which are tamed by $\omega$.
The latter
means that
$J$ is $\omega$-positive, i.e., $\omega(\xi, J\xi)>0$ for every $x\in X$ and
every nonzero $\xi\in T_xX$.

The evaluation map $\ev: S\times \cals_U \times \ttt \times \calj^k_U \to X$,
given by
formula $\ev(x, u, J_S,J)\deff u(x)$, defines a bundle $E\deff \ev^*(TX)$
over $S\times \cals_U \times \ttt \times \calj^k_U$. We equip $E$ with the
natural complex structure, which is equal to $J(u(x))$ in the fiber $E_{(x,
u, J_S,J)}\cong T_{u(x)}X$. We shall denote by $(E_u,J)$ the restriction of
$E$ onto $S\times \{(u,J_S,J)\}$ which is isomorphic to $u^*TX$.

The bundle $E$ with the complex structure $J$ induces complex Banach bundles
$\wh\cale$ and $\wh\cale'$ over a product $\cals_U \times \ttt \times \calj^k_U$ with
fibers
$$
\eqalign{
\wh\cale_{(u,J_S,J)}  & \deff L^{1,p}(S, E_u),\cr
\wh\cale'_{(u,J_S,J)} & \deff L^p(S, E_u\otimes\Lambda^{(0,1)}S).
}
$$
Here $S$ is equipped with the complex structure $J_S$, $\Lambda^{(0,1)}S$ is
a complex bundle of (0,1)-forms on $S$ and $\otimes$ means the tensor
product over $\cc$ of bundles with the corresponding complex structures. The
bundle $\wh\cale$ is simply a pull-back of a tangent bundle $TL^{1,p}(S,X)$
with respect to the projection $(u,J_S,J)\in \cals_U \times \ttt \times \calj^k_U
\mapsto
u \in L^{1,p}(S,X)$.

On the other hand, the bundle $\wh\cale'$ is the target manifold of the
$\dbar$-operator for the map $u\in L^{1,p}(S,X)$. Namely, we have a
distinguished section $\sigma_\dbar$ of $\wh\cale'$,
$$
\sigma_\dbar(u,J_S,J)\deff \dbar_{J_S,J}u \deff
\half\bigl(du + J\scirc du \scirc J_S\bigr).
\eqno(8.1.1)
$$
If $f$ is another section of $\wh\cale'$ (given \eg by some explicit
geometric construction), then one can consider a $\dbar$-equation
$$
\dbar_{J_S,J}u=f(u,J_S,J).
$$
We shall consider only a homogeneous case where $f(u,J_S,J)\equiv
0$. The corresponding set of solutions is closed in $\cals_U \times \calj^k_U$
and will be denoted by
$$
\calp \deff \{\, (u,J_S, J) \in \cals_U \times \ttt \times \calj^k_U \;:\;
\dbar_{J_S,J}u=0 \,\}.  \eqno(8.1.2).
$$
A map $u\in L^{1,p}(S,X)$ such that

$$
\dbar_{J_S,J}u=0\eqno(8.1.3)
$$
\noindent for some $J_S\in\ttt_g$
and $J\in \calj^k_U$ is called $J$-holomorphic, or
$(J_S,J)$-holomorphic, and its image $M\deff u(S)$ is called a $J$-complex
curve.

\smallskip
The regularity theory for $\dbar$-equation says that $\calp$ is a closed
subset of
$$
\calx \deff \{\,(u,J_S,J)\in\cals_U\times \ttt \times \calj^k_U\;:\;
u\in C^1(S,X)\,\}.
$$
Set
$
\calx^* \deff \{\,(u,J_S,J)\in \calx \,:\,
\hbox{$u$ is an imbedding in a neighborhood of some $y\in S$ }\},
$
and let $\calp^*\deff \calp\cap \calx^*$. Then $\calx^*$ is open in $\calx$
and $\calp^*$ is open in $\calp$.
A group $\bfg$ acts on $\cals_U\times \ttt \times \calj^k_U$ and on $\calx$
in a natural way by composition,
$$
(f,u,J_S,J) \in \bfg\times \cals_U\times \calj^k_U
\mapsto f\cdot(u,J_S,J) \deff (u\scirc \hat f\inv, f\cdot J_S,J),
$$
where $\hat f:S\to S$ is a diffeomorphism induced by $f\in\bfg$ and $J_S$,
see (7.4.1). As usual, the inverse  $(\hat f)\inv$ is introduced to preserve
the associative law $(f\cdot g)\cdot(u,J_S,J)=f\cdot(g\cdot(u,J_S,J))$. The
sets $\calx^*$, $\calp$ and $\calp^*$ are invariant with respect to this
action.

Let $\calm\deff\calp^*/\bfg$ be a quotient with respect to this action and
$\pi_\calp: \calp^* \to \calm$ the corresponding projection. Also let 
$\pi_\calj: \calm\to\calj^k_U$ be a natural projection.

\smallskip\noindent
\bf Lemma 8.1.1. \it The projections $\calx^* \longto \calx^*/\bfg$ and
$\pi_\calp:\calp^* \longto \calm$ are principal $\bfg$-bundles.

\smallskip\noindent
\bf Proof. \rm
We consider the case $g\ge2$  first. It is known that in this case the
group $\bfg$ acts proper discontinuously on $\ttt_g$. This implies that the
same is true for the action of $\bfg$ on $\calx$. Moreover, the definition
of $\calx^*$ provides that $\bfg$ acts freely on this set. Consequently, the
map $\calx^*\longto \calx^*/\bfg$ is simply an (unbranched) covering.

\smallskip
Now we consider the case $g=0$. Note that in this case $S=S^2$ and
$\ttt_0=\{J\st\}$. Fix some $(u^0,J\st,J^0)\in \calx^*$.
Let $y_1$, $y_2$, $y_3$ be distinct points
on $S^2$ such that $u^0$ is an imbedding in a neighborhood of every $y_i$,
in particular $du^0$ is non-vanishing in $y_i$. Also let  $Z_i$ be smooth
submanifolds of codimension 2 in $X$, intersecting $u^0(S^2)$ transversally
in each $u^0(x_i)$, respectively.

Let $V\ni(u^0,J\st,J^0)$ be an open set in $\calx^*$, $W$ its
projection on $\calx^*/\bfg$, and $\bfg\cdot V\deff \{\,f\cdot(u,J\st,J):
f\in \bfg, (u,J\st,J)\in V\,\}$ its $\bfg$-saturation. Consider a set
$$
\calz \deff\{\, (u,J\st,J)\in \bfg\cdot V\;:\; u(y_i)\in Z_i\,\}.
$$
One can easily show that if $V$ is chosen sufficiently small, then
$\calz$ is a smooth closed Banach submanifold of $\bfg\cdot V$,
intersecting every
orbit $\bfg \cdot (u,J\st,J)$ transversally at exactly one point.
Moreover, we have a $\bfg$-invariant diffeomorphism
$\bfg\cdot V \cong \bfg \times \calz$, so that $\calz$ is
a local slice of $\bfg$-action at $(u^0,J\st,J^0)$.
This equips a quotient $\calx^*/\bfg$ with
a structure of a smooth Banach manifold such that the projection
$\calx^*\to\calx^*/\bfg$ is a smooth principle $\bfg$-bundle.

\smallskip
The case $g=1$ is a combination of the two above cases. We fix some
$(u^0,J_S^0, J^0)\in \calx^*$ and a point $y$ on $S$ such that $u^0$ is
an imbedding in a neighborhood of $y$;
in particular $du^0$ is non-vanishing in $y$. Also let  $Z$ be a smooth
submanifold of codimension 2 in $X$, intersecting $u^0(S)$ transversally
in $u^0(y)$. As in the case $g=0$, we fix a small open set
$V\ni(u^0,J^0,J^0)$ in $X$ and consider a set
$$
\calz \deff\{\, (u,J_S,J)\in \bfg\cdot V\;:\; u(y)\in Z\,\},
$$
where $\bfg\cdot V \deff \{\,f\cdot(u,J\st,J): f\in \bfg, (u,J\st,J)\in V\,\}$
is $\bfg$-saturation of $V$. If $V$ is chosen sufficiently small, then
$\calz$ is a slice to the action of $T^2\subset \bfg$. Thus, the projection
$\calx^*\longto \calx^*/T^2$ is a principle $T^2$-bundle.

Now we consider the action of ${\bold Sp}(2,\;\; \zz) =\bfg/T^2$ on
$\calx^*/T^2$. The same arguments as in the case $g\ge2$ show that this
action is free and proper discontinuous. This implies that the projection
$\calx^*/T^2\longto \calx^*/\bfg$ is a covering. Consequently,
$\calx^*\longto \calx^*/\bfg$ is a principle $\bfg$-bundle.

\smallskip
The natural inclusion $\calm \hookrightarrow \calx^*/\bfg$ is continuous
and closed. Further, we have a natural $\bfg$-invariant homeomorphism
$\calp^* \cong \calm \,\times_{\calx^*/\bfg} \calx^*$, which gives
a desired structure of a principle $\bfg$-bundle on $\calp^*$ with a base
$\calm$.
\qed

\smallskip\noindent
\bf Remark. \rm One can show that if $(u,J_S,J)\in \calp^*$, then $u$ is
an imbedding in all but finitely many points of $S$, in particular $J_S$
is completely defined by $u$ and $J$. Similarly, every class $\bfg\cdot
(u,J_S,J) \in \calm$ is completely defined by $J\in \calj^k_U$ and
a $J$-complex curve $M\deff u(S)$. Therefore we  shall denote
elements of $\calm$ by $(M,J)$. Our motivation is that the object we
really exploit is a complex curve $M=u(S)$ rather than its
concrete parameterization $u$. We hope that
the reader will not be confused by such a formal incorrectness.

\medskip
Using {\sl Lemma 8.1.1}, one can obtain from $\bfg$-invariant objects
on $\calp^*$ corresponding objects on $\calm$.
For example, on $\calp$ we have a (trivial) $S$-bundle with
a total space $\calp\times S$, natural projection on $\calp$, a complex
structure $J_S$ on a fiber over $(u,J_S,J)$, and
a map $\ev:\calp\times S \to X$ such that $\ev(u,J_S,J;y)\deff u(y)$
which is invariant with respect to the $\bfg$-action.
It gives a $\bfg$-bundle $\pi_\calc:\calc\to \calm$ with the total space
$\calc\deff\calp^* \times_\bfg S$ and fiber $S$. Moreover, $\calc$ is
equipped with a map $\ev: \calc\to X$. We shall imagine
$\calm$ as a moduli space of all  complex  curves in $X$ of
appropriate topological properties, $\pi_\calc:\calc \to\calm$ as a
corresponding universal bundle and $\ev:\calc \to X$ as an evaluation map.
In particular, every fiber $\calc_{(M,J)} \deff \pi_\calc\inv(M,J)$ over
$(M,J)\in \calm$ possesses a canonical complex structure $J_{\calc_{(M,J)}}
=J_S$.

In a similar way we define complex Banach bundles $\cale$ and $\cale'$ over
$\calx^*/\bfg\supset \calm$. Note that we have a natural continuous
$\bfg$-equivariant inclusion $\calx^*\hook \cals_U\times \calj^k_U$.
Let $(u,J_S,J)\in \cals_U\times \calj^k_U$, $f\in\bfg$, and let
$\hat f:S\to S$ be a diffeomorphism defined as in (1.1), so that
$f\cdot(u,J_S,J)= (u\scirc\hat f\inv,f\cdot J_S,J)$. Define the operator
$f_*$ by setting
$$
f_*:s\in\wh\cale_{(u,J_S,J)} \mapsto (\hat f\inv)^*s\in
\wh\cale_{f\cdot(u,J_S,J)}.
$$
This defines a natural lift of a $\bfg$-action on $\cals_U\times \calj^k_U$
to a $\bfg$-actions on $\wh\cale$. In the same way we define the action of
$\bfg$ on $\wh\cale'$. Since the projection $\calx^*\to \calx^*/\bfg$ admits
a local $\bfg$-slice, there exist uniquely defined bundles $\cale$ and
$\cale'$ over $\calx^*/\bfg$, whose lifts onto $\calx^*$ is $\bfg$-equivariant
isomorphic to $\wh\cale$ and $\wh\cale'$, respectively. In particular, for
any $(u,J_S,J)\in \calx^*$ representing $\bfg\cdot(u,J_S,J)\in \calx^*/\bfg$
we have the natural isomorphisms
$$
\eqalign{
\cale_{\bfg\cdot(u,J_S,J)}  \cong  \wh\cale_{(u,J_S,J)}
& = L^{1,p}(S, E_u),\cr
\cale'_{\bfg\cdot(u,J_S,J)} \cong  \wh\cale'_{(u,J_S,J)}
& = L^p(S, E_u\otimes\Lambda^{(0,1)}S).
}
$$

\medskip
The question we are interested in is whether we can deform a given
(compact) $J_0$-holomorphic curve $M_0=u_0(S)$ into a compact complex
curve $M_1$ which is holomorphic with respect to the given (integrable, for
example)
complex structure $J\st$ on $X$. The idea is to use the continuity
method: one finds an appropriate homotopy $h(t)=J_t$, $t\in[0,1]$, of almost
complex structures connecting $J_0$ with $J\st=J_1$, and shows that there
exists a continuous deformation $u_t:S\to X$ of a map $u_0$ into $u_1$ such
that $u_t$ is $J_t$-holomorphic for all $t\in[0,1]$. This can be successfully
done by studying the linearization of the equation $\dbar u=0$.

\smallskip\noindent
\bf Lemma 8.1.2. \it Let $\calx$ be a Banach manifold, $\cale\to \calx$ and
$\cale'\to\calx$ $C^1$-smooth Banach bundles over $\calx$, $\nabla$ and
$\nabla'$ linear connections in $\cale$ and $\cale'$, respectively, $\sigma$
a (local) $C^1$-section of $\cale$ and $D: \cale \to \cale'$ a $C^1$-smooth
bundle homomorphism.

\sli If $\sigma(x)=0$ for some $x\in \calx$, then the map $\nabla\sigma_x:
T_x\calx \to \cale_x$ is independent of the choice of a connection $\nabla$
in $\cale$;

\slii Set $K_x\deff \ker(D_x:\cale_x\to \cale'_x)$ and $Q_x\deff
\coker(D_x:\cale_x\to \cale'_x)$ with the corresponding imbedding
$i_x:K_x \to \cale_x$ and projection $p_x:\cale'_x \to Q_x$. Let
$\nabla^\hom$ be a connection in $\hom(\cale,\cale')$ induced by connections
$\nabla$ and $\nabla'$. Then the map
$$
p_x \scirc(\nabla^\hom\!\! D_x)\scirc i_x:
T_x\calx\to \hom(K_x,Q_x)
$$
is independent of the choice of connections
$\nabla$ and $\nabla'$.

\smallskip\noindent
\bf Remark. \rm Taking into account this lemma, we shall use the following
notation. For $\sigma\in \Gamma(\calx, \cale)$, $D\in \Gamma(\calx,
\hom(\cale, \cale'))$ and $x\in \calx$ as in the hypothesis of the lemma, we
shall denote by $\nabla\sigma_x:T_x\calx \to \cale_x$ and $\barr{\nabla\!\!D}:
T_x\calx\times \ker D_x \to \coker D_x$ the corresponding operators without
pointing out
which connections were used to define them.

\smallskip\noindent
\bf Proof. \rm \sli Let $\wt\nabla$ be another connection in $\cale$. Then
$\wt\nabla$ has a form
$\wt\nabla=\nabla +A$ for some $A\in\Gamma(\calx,\hom(T\calx,\endo(\cale)))$.
Thus, for $\xi\in T_x\calx$ we get $\wt\nabla_\xi\sigma -\nabla_\xi\sigma =
A(\xi,\sigma(x))=0$.

\slii Similarly, let $\wt\nabla'$ be another connection in $\cale'$, and let
$\wt\nabla^\hom$ be a connection in $\hom(\cale,\cale')$ induced by
$\wt\nabla$ and $\wt\nabla'$. Then
$\wt\nabla'$ also has the form $\wt\nabla=\nabla +A'$ for some
$A'\in\Gamma(\calx,\hom(T\calx,\endo(\cale')))$.
Thus, for $\xi\in T_x\calx$ we obtain $\wt\nabla_\xi^\hom D -\nabla_\xi^\hom D=
A'(\xi)\scirc D_x - D_x \scirc A(\xi)$. The statement of the lemma now follows
from the identities $p_x\scirc D_x=0$ and $D_x\scirc i_x=0$.
%\par\nobreak
\qed

\medskip\noindent
{\bigsl 8.2. Transversal Mappings.}

\smallskip
To have the possibility of deforming a complex curve along a given
path of almost-complex structures, it is useful to know in which points
$(u,J_S,J)$ the set $\calp$ of holomorphic maps is a Banach manifold.
Note that by definition the set $\calp$ is essentially an intersection
of the zero section and the $\sigma_\dbar$-section in the total space
of $\cale' $ over $S\times \ttt \times \calj_U^k$. Thus, we are interested in
which points these sections meet transversally.

\state Definition 8.2.1. 
{\it Let $\calx$, $\caly$, and $\calz$ be Banach manifolds
with $C^k$-smooth maps $f:\caly \to \calx$ and $g:\calz\to\calx$, $k\ge1$.
Define the {\sl fiber product} $\caly\times_\calx \calz $ by setting
$\caly\times_\calx \calz \deff \{\,(y,z)\in \caly\times\calz \;:\;
f(y)=g(z)\,\}$. The map $f$ is called {\sl transversal to $g$} at point
$(y,z)\in \caly\times_\calx \calz$ with $x\deff f(y)=g(z)$, and $(y,z)$
is called a {\sl transversality point}, \iff
a map $df_y\oplus -dg_z: T_y\caly \oplus T_z\calz \to T_x\calx$ is
{\sl surjective} and its kernel admits a closed complement.
}

The set of transversality points $(y,z)\in \caly\times_\calx \calz$ we shall
denote by $\caly\times^\trans_\calx \calz$, with $\trans$ symbolizing
the transversality condition.

In the special case, when the map $g:\calz \to \calx$ is a closed imbedding,
the fiber product $\caly\times_\calx \calz$ is simply the pre-image
$f\inv\calz$ of $\calz\subset\calx$. In particular, every point
$(y,z)\in\caly\times_\calx \calz$ is completely defined by the point
$y\in\caly$, $z=f(y)\in\calz\subset \calx$. In this case we simply say that
$f:\caly \to \calx$ is {\sl transversal to $\calz$ in $y\in\caly$},
\iff $(y,f(y))$ is a transversal point of $\caly\times_\calx \calz\cong
f\inv\calz$.

\state Lemma 8.2.1. {\it The set $\caly\times^\trans_\calx \calz$ is open in
$\caly\times_\calx \calz$ and is a $C^k$-smooth Banach manifold with
a tangent space
$$
T_{(y,z)}\caly\times^\trans_\calx \calz =
\ker\bigl( df_y\oplus d(-g_z): T_y\caly \oplus T_z\calz \to T_x\calx \bigr).
$$}\rm

\state Proof. Fix $w_0\deff(y_0,z_0)\in \caly\times^\trans_\calx \calz$
and set $K_0 \deff \ker( df_{y_0}\oplus dg_{z_0}:
T_{y_0}\caly \oplus T_{z_0}\calz \to T_x\calx \bigr)$. Let $Q_0$ be a closed
complement to $K_0$. Then the map $df_{y_0}\oplus dg_{z_0}:Q_0\to T_x\calx$
is an isomorphism.

Due to the choice of $Q_0$, there exists a neighborhood $V\subset
\caly\times\calz$ of $(y_0,z_0)$ and $C^k$-smooth maps $w':V\to K_0$ and
$w'':V\to Q_0$ such that $dw'_{w_0}$ (resp.\ $dw''_{w_0}$) is the projection
from $T_{y_0}\caly \oplus T_{z_0}\calz$ onto $K_0$ (resp.\ onto $Q_0$),
so that $(w',w'')$ are coordinates in some smaller  neighborhood $V_1\subset
\caly\times\calz$ of $w_0=(y_0,z_0)$. It remains to consider the equation
$f(y)=g(z)$ in new coordinates $(w',w'')$ and apply the implicit function
theorem.
\qed

\medskip
We have defined $\calp$ as a pre-image of a zero section $\sigma_0$ of $\cale'$
with respect to the map $\sigma_\dbar:\cals_U\times \calj_S\times \calj^k_U
\to \cale'$.
Due to {\sl Lemma 8.2.1}, the set $\calp$ is a Banach manifold in those points
$(u,J_S,J)\in\calp$, where $\sigma_\dbar$ is transversal to $\sigma_0$.
However, in any point $(u,J_S,J;0)$ on the zero section $\sigma_0$ of
$\cale'$ we have the natural decomposition
$$
T_{(u,J_S,J;0)}\cale' =
d\sigma_0\bigl( T_{(u,J_S,J)}(\cals_U\times \calj_S\times \calj_U^k) \bigr)
\oplus \cale'_{(u,J_S,J)},
$$
where the first component is the tangent space to the zero section of $\cale'$
and the second is the tangent space to the fiber $\cale'_{(u,J_S,J)}$.
Let $p_2$ denote the projection on the second component. Then the
transversality $\sigma_\dbar$ and $\sigma_0$ is equivalent to the surjectivity
of the map $p_2\scirc d\sigma_\dbar: T_{(u,J_S,J)}(\cals_U\times \calj_U^k)\to
\cale'_{(u,J_S,J)}$. But, due to {\sl Lemma 8.1.2 }, this map is a linearization
of $\sigma_\dbar$ at $(u,J_S,J)$ and has a form (8.1.2).

Thus, the transversality of $\sigma_\dbar$
to $\sigma_0$ at $(u,J_S,J)\in \calp$ is equivalent to the surjectivity of
the operator
$$
\eqalign{
\nabla\sigma_\dbar:\;&T_uL^{1,p}(S,X)\oplus T_{J_S}\ttt_g \oplus T_J\calj^k_U
\longto \wh\cale'_{(u,J_S,J)}
\cr
\nabla\sigma_\dbar:\;&(v,\dot J_S,\dot J) \longmapsto
D_{(u,J)}v + J\scirc du \scirc \dot J_S
+ \dot J\scirc du \scirc J_S.
}\eqno(8.2.1)
$$
{\sl Definition 7.1.2} provides that the quotient of $\wh\cale'_{(u,J_S,J)}=
L^p_{(0,1)}(S,E_u)$ by the image of $D_{u,J}$ is $\sfh^1_D(S,E_u)$.
The induced
map $\dot J_S \in T_{J_S}\ttt_g \mapsto J\scirc du \scirc \dot J_S \in
\sfh^1_D(S,E_u)$ is also easy to describe. From equality (8.2.1) and
{\sl Corollary 7.3.2} it follows that its image is equal to the image of the
homomorphism $du\scirc J_S:\sfh^1(S,TS) \to \sfh^1_D(S,E_u)$ and its cokernel
is $\sfh^1_D(S,N_u)$.

It remains to study the image of $T_J\calj^k_U$ in $\sfh^1_D(S,N_u)$. For
$(u,J_S,J)\in \calp$ we define  $\Psi=\Psi_{(u,J)}:T_J\calj^k_U \to
\wh\cale'_{(u,J_S,J)}$ by setting $\Psi_{(u,J)}(\dot J)\deff \dot J\scirc du
\scirc J_S$. Let $\barr\Psi=\barr\Psi_{(u,J)}:T_J\calj^k_U \to
\wh\cale'_{(u,J_S,J)}$ be induced by $\Psi$. Recall that if $(u,J_S,J)\in
\calp$, then $J_S$ is determined by $u$ and $J$.

\state Lemma 8.2.2. \sl(Infinitesimal Transversality). \it Let $(u,J_S,J)\in
\calp^*$. Then the operator $\barr\Psi:T_J\calj^k_U \to \sfh^1_D(S, N_u)$
is surjective.

\state Proof. \rm In {\sl Lecture 3} we proved that for $(u,J_S,J)\in
\calp^*$  mapping $u$ is an imbedding in the neighborhood of all but a finite
number of points $x\in S$. Thus, there exists such a nonempty open set
$V\subset S$ that $u(V)\subset U$ and $u\mid_V$ is an imbedding.

By {\sl Lemma 7.2.2} we have an isomorphism $\sfh^0_D(S,N_u^*\otimes K_S)
\cong \sfh^1_D(S,N_u)^*$. The fact that operators of the type
$D=\dbar + R: L^{1,p}(S,E)\to L^p_{(0,1)}(S,E)$ on the compact Riemann
surface $(S,J_S)$ are Fredholm ensures the existence of the 
finite basis $\xi_1,...,\xi_l$ of
the space $\sfh^0_D(S,N_u^*\otimes K_S)$. By {\sl Lemma 3.1.1}
every $\xi \in \sfh^0_D(S,N_u^*\otimes K_S)$ vanishes not more than
in $c_1(N_u^*\otimes K_S)[S]$ points on $S$ (compare to the proof of
{\sl Corollary 7.2.3}). From here it follows that there exists
$\psi_i\in C_c^k(V, N\otimes \Lambda^{0,1}), i=1,...,l$, generating an
$\rr $-basis of the space $\sfh^1_D(S,N)$.

Take some $\psi_i\in C^k_c(V, N\otimes \Lambda^{0,1})$.
It is a $\cc$-antilinear $C^k$-smooth homomorphism from
$TS\ogran_V$ into $N\ogran_V$ vanishing outside some compact
in $V$. Since $u\ogran_V$ is a $C^k$-imbedding and $u(V)\subset U$,
$\psi_i$ can be represented in the form
$\psi_i =\pr_N \scirc \dot J \scirc du \scirc J_{S}$
for some $J$-anti-linear $C^k$-smooth endomorphism $\dot J$
of $TX$ vanishing outside some compact in $U$. Thus, $\dot J\in
T_J\calj^k_U$ and $\Psi\dot J=\psi_i^{(\nu)}$.
\qed

\state Corollary 8.2.3. \it $\calm$ and $\calp^*$ are $C^k$-smooth Banach
manifolds and $\pi_\calj:\calm \to \calj^k_U$ is a Fredholm map. For
$(M,J)\in\calm$ with $M=u(S)$ there exist natural isomorphisms
$$
\eqalign{
\ker(d\pi_\calj:T_{(M,J)}\calm \to T_J\calj^k_U)&\;\cong\; \sfh^0_D(S,
\caln_M),
\cr
\coker(d\pi_\calj:T_{(M,J)}\calm \to T_J\calj^k_U)&\;\cong\; \sfh^1_D(S,
\caln_M),
}
$$
where $\caln_M=\calo(N_u)\oplus \caln\sing_u$ is a normal sheaf to $M$,
$\sfh^0_D(S, \caln_M)$ denotes $\sfh^0_D(S,N_u)$
$\oplus\sfh^0(S,\caln_u\sing)$
and $\sfh^1_D(S, \caln_M)$ denotes $\sfh^1_D(S, N_u)$. The index of
$\pi_\calj$
$$
\ind_\rr(\pi_\calj)=\ind_\rr(\caln_M)\deff \dimr\sfh^0_D(S, \caln_M)
-\sfh^1_D(S, \caln_M)
$$
is equal to $2(c_1(X)[M]+ (n-3)(1-g))$, where $n\deff\dimc X$.

\smallskip\noindent
\bf Proof. \rm One easily sees that the section $\sigma_{\dbar }$ is
$C^k$-smooth if $J\in \calj^k_U$. For $\calp^*$ the statement follows from
{\sl Lemmas 8.2.1} and
{\sl 8.2.2}. Moreover, $\calp^*$ is a $C^k$-smooth submanifold of $\calx^*$.
Further, local slices of the projection $\pi_\calp:\calp^*\to \calm$
are also $C^k$-smooth. This induces the structure of $C^k$-smooth Banach
manifold on $\calm$.

Consider the natural projection $\pi: \calp^* \to \calj^k_U$.
The tangent space $T_{(u,J_S,J)}\calp^*$ consists of $(v,\dot J_S,\dot J)$
with $D_{u,J}v + \half J\scirc du \scirc \dot J_S + \half\dot J\scirc du
\scirc J_S =0$, and the differential $d\pi: T_{(u,J_S,J)}\calp^* \to
T_J\calj^k_U$ has a form $(v,\dot J_S,\dot J)\in T_{(u,J_S,J)}\calp^*\mapsto
\dot J\in  T_J\calj^k_U$.

The kernel $\ker(d\pi) T_{(u,J_S,J)}\calp^*$
consists of solutions of the equation
$$
D_{u,J}v + \half J\scirc du \scirc \dot J_S =0
%\eqno(2.22)
$$
with $v\in\cale_{(u,J_S,J)}$ and $\dot J_S\in T_{J_S}\ttt_g$.
Since the map $\pi_\calp:\calp^*\to\calm$ is a principal $\bfg$-bundle,
the kernel $\ker(d\pi_\calj:T_{(M,J)}\calm \to T_J\calj^k_U)$ is obtained
from $\ker(d\pi)$ by taking a quotient with respect to the tangent space
to a fiber $\bfg\cdot(u,J_S,J)$, which is equal to $du(\sfh^0(S,TS))$.
Using relations $\sfh^0(S,TS)=
\ker(\dbar_{TS}:L^{1,p}(S,TS)\to L^p(S,TS\otimes\Lambda^{(0,1)}S)$,
$T_{J_S}\ttt_g\cong\sfh^1(S,TS)=\coker(\dbar_{TS})$ and
$du\scirc \dbar_{TS}= D_{(u,J)}\scirc du$, we conclude that
the space $\ker(d\pi_\calj)$ is isomorphic to the quotient
$$
\{v\in L^{1,p}(S,E_u)\;:\; Dv=du(\phi) \hbox{ for some }\phi\in
L^p(S,TS\otimes\Lambda^{(0,1)}S)\}
\!\!\bigm/\!\!
du\bigl(L^{1,p}(S,TS)\bigr).
$$
Hence, $\ker(d\pi_\calj:T_{(M,J)}\calm\to T_J\calj^k_U)\cong
\sfh^0_D(M,\caln_M)$ by {\sl Theorem 7.3.1}. In particular, $\ker(d\pi_\calj)$
is finite dimensional.

\smallskip
Similarly, the image of $d\pi_\calj$ consists of those $\dot J$ for which
the equation
$$
D_{u,J}v + \half J\scirc du \scirc \dot J_S + \half\dot J\scirc
du \scirc J_S =0
$$
has a solution $(v,\dot J_S)$. Hence, $\im(d\pi_\calj)=\ker\barr\Psi$, and
$\coker(d\pi) \cong \sfh^1_D(S,N_u)$. Thus, $d\pi_\calj$ is a Fredholm map.
This implies the Fredholm property for the projection $\pi:\calp^* \to
\calj^k_U$.

\smallskip
Due to {\sl Corollary 7.3.2}, $\ind_\rr(\caln_M)=\ind_\rr(E_u)-\ind_\rr(TS)$.
Using the index and Riemann-Roch theorems, we get from $c_1(E)=c_1(X)[M]$
and $c_1(TS)=2-2g$ the needed formula
$\ind_\rr(\caln)=2\bigl(c_1(X)[M] + n(1-g) - (3-3g)\bigr)=
2(c_1(X)[M] + (n-3)(1-g))$.
\qed

\smallskip
A straightforward application of this statement and the Sard Lemma for 
Fredholm maps gives the following

\state Corollary 8.2.4. {\it If smoothness $k$ of the structures in 
$\calj^k$ is large enough and with $\ind_\rr (\pi_{\calj })=2(c_1(X)[\gamma ] 
+(n-3)(1-g))<0$ for the homology class $\gamma $, then the following holds:

\smallskip
1) For any $(M_0,J_0)\in \calm_{\calj^k}$ with $[M]=[\gamma ]$   
there is a neighborhood $W$ fo $J_0$ in $\calj^k_{\omega }$ and a closed subset $S
\subset W$ of Hausdorff codimension close to $2$ such that for $J\in W
\setminus S$ $\calm_J$ is empty;

\smallskip
2) In the manifold $C^k_{[\gamma ]}(I,\calm_{\calj^k})$ of $k$-smooth paths 
in $\calm_{[\gamma ]} $ there is 
a closed subset $R$ of Hausdorff codimension close to one such that for 
any path $h\in C^k(I,\calm_{\calj^k}\setminus R$ the moduli space 
$\calm_h$ is empty.
}

\medskip\noindent
{\bigsl 8.3. Components of the Moduli Space.}

\smallskip
Before stating the next results, we introduce some new notations.

\nobreak\state
Definition 8.3.1. Let $Y$ be a $C^k$-smooth finite-dimensional manifold,
possibly with $C^k$-smooth boundary $\d Y$, and $h:Y \to \calj^k_U$ a 
$C^k$-smooth map. Define the {\sl relative moduli space}
$$
\calm_h \deff Y\times_{\calj^k_U}\calm \cong
\{\,(u,J_S,y)\in \cals_U\times \ttt_g\times Y\,:\, (u,J_S, h(y))\in \calp^*
\,\}/\bfg
$$
with the natural projection $\pi_h:\calm_h \to Y$. In the special case
$Y=\{J\}\hook \calj^k_U$, we obtain the moduli space of $J$-complex curves
$\calm_J\deff \pi_\calj\inv(J)$. The projection $\pi_h:\calm_h \to Y$
is a fibration with a fiber $\pi_h\inv(y)=\calm_{h(y)}$. We shall denote
elements of $\calm_h$ by $(M,y)$, where $M=u(S)$ with $h(y)$-holomorphic
map $u:S\to X$.

\smallskip\noindent
\bf Lemma 8.3.1. \it Let $Y$ be a $C^k$-smooth finite-dimensional manifold,
and $h:Y \to \calj^k_U$ a $C^k$-smooth map. Suppose that for some
$y_0\in Y$ and $u_0\in \calm_{h(y_0)}$ the map
$\barr\Psi \scirc dh: T_{y_0}Y \to \sfh^1_D(S, N_{u_0})$
is surjective. Then $\calm_h$ is a $C^k$-smooth
manifold in some neighborhood of $(M_0,y_0)\in \calm_h$ with the tangent
space
$$
T_{(M,y)}\calm_h=\ker\bigl(D\;\oplus\; \Psi\scirc dh:
\cale_{u,h(y)}\oplus T_yY \longto \cale'_{u,h(y)} \bigr)
\bigm/ du(\sfh^0(S, TS)).
\eqno(8.3.1)
$$

\smallskip\noindent
\bf Proof. \rm Let $y\in Y$, $(u,J_S,h(y)\in\calp^*$ and $M=u(S)$, so that
$(M,y)\in \calm_h$. From the proof of {\sl Corollary 8.2.3} it follows that the
image of the map $d\pi_\calj:T_{(M,h(y))}\calm\to T_{h(y)}\calj^k_U$ is equal
to $\ker\bigl(\barr\Psi_{(u,h(y))}\bigr)$, and that the cokernel of
$d\pi_\calj$ is mapped by $\barr\Psi$ isomorphically onto $\sfh^1_D(S, N_u)$.
The statement of the lemma now follows from {\sl Lemma 8.2.1}.
\qed

\medskip\noindent
\bf Definition 8.3.2. \rm Let $Y$ be a compact manifold, $h:Y\to \calj^k_U$
a $C^k$-smooth map, $\calm_h\subset \calm \times Y$ a corresponding moduli
space and $(M_0,y_0)\in \calm_h$ a point. A {\sl component
$\calm_h(M_0,y_0)$ of $\calm_h$ through $(M_0,y_0)$} is the set of those
$(M,y)\in\calm_h$, that for every open neighborhood $W$ of $h(Y)\subset
\calj^k_U$ there exists a continuous path $\gamma:[0,1]\to \calm$,
with the following properties:

\smallskip
\item{\sl a)} $\gamma(0)=(M_0,h(y_0))$ and $\gamma(1)=(M,h(y))$,
i.e., $\gamma$ connects $(M_0,y_0)$ and $(M,y)$ in $\calm$;

\smallskip
\item{\sl b)} $J_t\deff\pi_\calj(\gamma(t))\in W\subset\calj^k_U$ for
any $t\in[0,1]$, i.e., the corresponding path of almost complex structures
$J_t$ lies in the given neighborhood $W$ of $h(Y)$.

\medskip\noindent
\bf Lemma 8.3.2 \it Let $h:Y\to \calj^k_U$, $(M_0,y_0)\in \calm_h$ and
$\calm_h(M_0,y_0)$ be as in  Definition 8.3.2. Then

\smallskip
\sli $\calm_h(M_0,y_0)$ is a closed subset of $\calm_h$,

\smallskip
\slii if $\calm_h(M_0,y_0)$ is compact, then there exists a set
$\calm_h^0$ containing $\calm_h(M_0,y_0)$, which is compact and open
in $\calm_h$;

\smallskip
\sliii if $\calm_h(M_0,y_0)$ is not compact, then there exists a continuous
path $\gamma:[0,1)\to \calm_h$ with the following properties:

\smallskip
\item{\sl a)} $\gamma(0)=(M_0,h(y_0))$, i.e., $\gamma$ begins at $(M_0,y_0)$;

\smallskip
\item{\sl b)} there exists a sequence $t_n\nearrow 1$ such that
the sequence $(M_n,J_n)\deff\gamma(t_n)$ lies in $\calm_h$ and is
discrete there, but the sequence $\{J_n\}$ converges to some
$J^*\in \calj_U^k$.

\medskip\noindent
\bf Proof.\ \sli\hskip 1.2em
\rm Let $(M',y')\in \barr{\calm_h(M_0,y_0)}\subset \calm_h$
and $J'=h(y')$.
Let $W$ be any open neighborhood of $h(Y)\subset \calj_U^k$, and let
$\{(M_n,y_n)\}$ be a sequence in $\calm_h$ converging to $(M',y')$. Then
there exists some ball $B\ni(M',J')$ in $\calm$ whose projection on
$\calj^k_U$ lies in $W$. Since some $(M_n,J_n)$ with $J_n=h(y_n)$
lies in $B$, there exists
a path $(M_t,J_t)$ in $\calm$ connecting $(M_0,h(y_0))$ with $(M',J')$
such that $J_t\in W$ for any $t\in [0,1]$. This proves the closedness of
$\calm_h(M_0,y_0)$.

\medskip\noindent
\slii\hskip 1.2em
Let $(M',y')\in \calm_h$ and $J'=h(y')$. Fix a finite dimensional subspace
$F\subset T_{J'}\calj^k_U$ such that the map $D_{u',J'} \oplus \Psi:
\cale_{u',J'} \oplus F \to \cale'_{u',J'}$ is surjective. Let $B\ni0$ be a ball
in $F$. Find a $C^k$-smooth map $H:Y\times B$ such that $H(y,0)\equiv h(y)$
and $dH_{(y',0)}: T_{(y',0)}(Y\times B) \to T_{J'}\calj^k_U$ induces
isomorphism $T_0B\buildrel\cong\over \longto F\subset T_{J'}\calj^k_U$. Then
$\calm_H$ contains a neighborhood $W\ni (M',y',0)$ which is $C^k$-smooth
manifold such that $W\cap \calm_h$ is closed in $W$. This implies that
$\calm_h$ is a locally compact topological space.

Since $\calm_h(M_0,y_0)$ is a compact subset of $\calm_h$, it has an open
neighborhood $V$ whose closure $\barr V \subset\calm_h$ is also compact.
Let $W_i\subset \calj^k_U$ be a fundamental system of neighborhoods
of $h(Y)$, so that $\cap_i W_i =h(Y)$. Let $V_i$ denote the set of those
$(M,y)\in \calm_h$ such that $(M,h(y))$ and $(M_0,h(y_0))$ can be connected by
a continuous path $(M_t, J_t)$ in $\calm$ with $J_t$ lying in $W_i$.
Then $\cap V_i =\calm_h(M_0,y_0)$. The same arguments as in the part \sli
of the proof show that every $V_i$ is both open and closed in $\calm_h$.

We state that there exists a number $N\in \nn$ such that
$\bigl(\cap_{i=1}^N V_i\bigr) \cap\barr V \subset V$. If it is false,
then for any
$n\in \nn$ there would exist $(M_n,y_n)\in \bigl(\cap_{i=1}^n V_i\bigr)
\cap\barr V \bes V$. But, in this case, some subsequence of $\{(M_n,y_n)\}$
would converge to some $(M^*,y^*)\in\bigl(\cap_{i=1}^\infty V_i\bigr)
\cap\barr V \bes V$, which is impossible.
For such $N\in \nn$ the set
$$
\def\capl{\mathop{\cap}}
\calm^0_h\deff \left(\capl_{i=1}^N V_i\right) \cap\barr V
=\left(\capl_{i=1}^N V_i\right) \cap V
$$
satisfies the conditions of the lemma.

\medskip\noindent
\sliii\hskip 1.2em
Suppose that $\calm_h(M_0,y_0)$ is not compact. Then
there exists a discrete sequence $\{(M_n,y_n)\}$ in
$\calm_h(M_0,y_0)$. Since $Y$ is compact, we may assume that $y_n$ converges
to some $y^*$. For any $n\in \nn$ fix a path $\gamma_n:[0,1] \to \calm_h$
connecting $(M_{n-1}, h(y_{n-1}))$ with $(M_n,h(y_n))$. Set $t_n\deff
1-2^{-n}$. For $t\in [t_{n-1},t_n]$ define $\gamma(t)\deff
\gamma_n(2^n(t-t_{n-1}))$. Then $\gamma:[0,1)\to \calm_h$ and
$t_n\nearrow1$ satisfy the conditions of the lemma.
\qed

\smallskip \state Theorem 8.3.3. {\it
Let $(M_0,J_0)\in \calm$ and let $h:[0,1]\to \calj^k_U$ be $C^k$-smooth and
with $h(0)=J_0$. Suppose that there exists a subset $\calm_h^0$ of $\calm_h$
which is compact, open and contains $(M_0,J_0)$. Suppose also that
$\ind(\pi_\calj)=2(c_1(X)[M_0]+ (n-3)(1-g))$ is non-negative.
Then $h$ can be $C^k$-approximated by smooth maps $h_n:[0,1] \to \calj^k_U$
such that

\sli every $\calm_{h_n}$ contains a component $\calm_{h_n}^0$
which is a smooth connected compact manifold of the expected dimension
$\dimr(\calm_{h_n}^0)=\ind(\pi_\calj)+1$;

\slii $\calm_{h_n(0)}^0\deff \pi_{h_n}\inv(0)\cap \calm_{h_n}^0$
and $\calm_{h_n(1)}^0\deff \pi_{h_n}\inv(1)\cap \calm_{h_n}^0$
are also smooth manifolds of the expected dimension $\ind(\pi_\calj)$
($\calm^0_{h_n(1)}$ may be empty),
and $\calm_{h_n}$ is a smooth bordism between $\calm_{h_n(0)}^0$ and
$\calm_{h_n(1)}^0$;

\sliii every $\calm_{h_n(0)}^0$ is connected with $(M_0,J_0)$ by a path in
$\calm$, i.e., there exists a $C^k$-path $\gamma_n: [0,1] \to \calm$ such that
$\gamma_n(0)= (M_0,J_0)$ and $\gamma_n(1)\in \calm_{h_n(0)}^0$; in particular,
every $\calm_{h_n(0)}^0$ is nonempty;

\sliv for every $(M,J)\in \calm_{h_n}$ one has $\dimr\sfh^1_D(S,N_M)\le1$.
}

\state Proof. Denote $\calm_h^0$ by $K$. Let $\cale_K$ and $\cale'_K$ be
the pull-backs onto $K$ of the bundles $\cale\to \calm$ and $\cale'\to
\calm$, respectively. Further, let
$\calt\deff h^*T\calj^k_U$ be the pull-back of the
tangent bundle $T\calj^k_U$ onto $[0,1]$.

Due to {\sl Lemma 8.2.2}, for every $(M,J)\in K$ with $M=u(S)$ we find that
$n_{(M,J)}\in \nn$ and a $C^k$-smooth homomorphism
$P_{(M,J)}: F_{(M,J)} \to \calt$ form a trivial vector bundle $F_{(M,J)}$
over $[0,1]$ of rank $\rank F_{(M,J)}=n_{(M,J)}$ such that the operator
$$
D_{(u,J)}\oplus \Psi_{(u,J)} \scirc P_{(M,J)} :
\cale_{(M,J)}\oplus F_{(M,J)} \to \cale'_{(M,J)}
\eqno(8.3.2)
$$
is surjective.

Since $K$ is compact, we can take an appropriate finite collection
of $(F_{\!(M,J)},$ $P_{\!(M,J)})$ and construct a homomorphism $P: F\to \calt$
from a trivial vector bundle $F\cong [0,1]\times \rr^n$ such that
$D_{(u,J)}\oplus \Psi_{(u,J)} \scirc P : \cale_{(M,J)}
\oplus F \to \cale'_{(M,J)}$ is surjective
for all $(M,J)\in K$.

For $\dot J$ lying in some small ball in $T_J\calj^k_U$ we set
$\exp_J(\dot J)\deff J{(1-J\dot J/2)\over(1+J\dot J/2)}$. Differentiating
the identity $J^2=-1$ gives the relation $J\dot J =- \dot JJ$, and
consequently the equality
$$
\left(J{(1-J\dot J/2)\over(1+J\dot J/2)} \right)^2=-1,
$$
which means that $\exp_J$ takes values in $\calj^k_U$.
Further, it is easy to see that the differential of $\exp_J$ in
$0\in T_J\calj^k_U$ is the identity map of $T_J\calj^k_U$. Thus,
$\exp_J$ is a natural exponential map for $\calj^k_U$.

Take a sufficiently small ball $B=B(0,r)\subset \rr^n$ and
define a map $H^*:[0,1]\times B  \to \calj^k_U$ by setting
$H^*(t, y)\deff \exp_{h(t)}(P(t,y))$. The construction of $H^*$
provides that for all $(M,J)\in K$ with $J=h(t)$ the map
$D_{(u,J)}\oplus \Psi_{(u,J)} \scirc dH^*(t,0) : \cale_{(M,J)}
\oplus T_{(t,0)}([0,1]\times B) \to \cale'_{(M,J)}$ is surjective.
Making an appropriate small perturbation of $H^*$, we obtain
a $C^k$-smooth function $H: [0,1]\times B  \to \calj^k_U$ with the
following properties:

\sli for $H(t,0)=h(t)$, i.e., $H$ is a deformation of $h$ with a parameter
space $B$;

\smallskip
\slii
$
%\displaystyle
D_{(u,J)}\oplus \Psi_{(u,J)} \scirc dH(t,0) : \cale_{(M,J)}
\oplus T_0B \longto \cale'_{(M,J)}
$
is surjective for all $(M,t)\in K$ with $J=h(t)$ and $M=u(S)$;

\sliii in a neighborhood of every $(t,y)\in [0,1]\times B$ with $y\not=0$
the map $H$ is $\calj^\infty_U$-valued and $C^\infty$-smooth.

\smallskip
Let us identify $K$ and $\calm_h$ with the subset of $\calm_H$ using natural
continuous imbedding $(M,t)\in \calm_h \mapsto (M,t,0)\in \calm_H$.
Due to {\sl Lemma 8.2.2}, there exists a neighborhood $V$ of $K$ in $\calm_H$
which is $C^k$-smooth manifold. Taking a smaller neighborhood of $K$ if
it is needed, we may assume that the closure $\barr V$ of $V$ is compact and
does not meet connected components of $\calm_h$ different from $\calm_h^0=K$.
Let $p:V \to B$ be the natural projection such that $(u,t,y)\mapsto y$.

We state that there exists a smaller ball $B_1=B(0,r_1)\subset B$ such that
for every $y\in B_1$ the set $p\inv(y)\subset V$ is compact. Suppose
the contrary is true. Then there would exist a sequence $y_n\in B$ converging
to $0\in B$ such that $p\inv(y_n)$ are not compact. Since $\barr V$ is
compact, there would also exist such $u_n$ and $t_n$ that $(u_n,t_n,y_n)$
lies in the closure $\barr{p\inv(y_n)}\subset \barr V$ but
not in $V$.  Taking an appropriate
subsequence we may assume that $(u_n,t_n,y_n)$ converges to $(u^*,t^*,0)$.
However, in this case $(u^*,t^*,0)\in \calm_h \cap \barr V$ and hence
$(u^*,t^*,0)\in V$. On the other hand $\barr V\bes V$ is compact, and thus
$(u^*,t^*,0)$ must belong to $\barr V\bes V$. The obtained contradiction
shows that the statement is true.

\smallskip
Set $V_1\deff p\inv(B_1)$. Due to the choice of $H$ the resticted
projection $p:V_1\bes K\to \check B_1 \deff B_1\bes\{0\}$ is
$C^\infty$-smooth. Due to the Sard lemma (see, e.g., [Fed], \S {\bf 3.4}),
there exists a dense subset $B_1^* \subset B_1$ such that for any
$y\in B_1^*$ the set $p\inv(y)$ is a $C^\infty$, a smooth compact manifold.
Fix a sequence $y_n\in B_1^*$, converging to $0\in B$ and set $h_n(t)\deff
H(t,y_n)$. Also set $\calm_{h_n}^0 \deff \calm_{h_n} \cap V$, so that
$\calm_{h_n}^0 = p\inv(y_n)$. Then every $\calm_{h_n}^0$ is a
$C^\infty$-smooth nonempty manifold, which is connected with $(M_0,J_0)$
by a path in $\calm$.

\smallskip
Due to {\sl Lemma 8.2.3}, the tangent space to $V_1\subset \calm_H$ at
$(u,t,y)$ is canonically isomorphic to
$$
\ker\!\Bigl(\!D_{u,H(t,y)}\,\oplus\, \Psi\scirc dH:
\cale_{u,H(t,y)}\,\oplus\, T_{(t,y)}\bigl([0,1]{\times} B_1\!\bigr)
\to \cale'_{u,H(t,y)} \!\Bigr)
\!\!\!\Bigm/ \!\!\!du(\sfh^0(S\!\!,\! TS)).
$$
Since $p: V_1 \to B_1$ is a projection of the form $(u,t,y)\in V_1
\mapsto y\in B_1$, the differential $dp_{(u,t,y)}$ maps the tangent vector
of the form $(\dot u,\dot t,\dot y)\in T_{(u,t,y)}V_1$ into $\dot y
\in T_yB_1$. This means that $dp_{(u,t,y)}$ is a restiction on
$\ker(D_{u,H(t,y)}\,\oplus\, \Psi\scirc dH)$ of the linear projection
$p_B:\cale_{u,H(t,y)} \oplus\, T_{(t,y)}\bigl([0,1]{\times} B_1\!\bigr)
\to T_y B_1$ such that $p_B(\dot u,\dot t,\dot y)=\dot y$.
In particular, for $y=y_n$ the map $dp_{(u,t,y)}$ is surjective,
which means the surjectivity of the map
$$
p_B: \ker \bigl(D_{u,H(t,y)}\,\oplus\, \Psi\scirc dH \bigr)
\longto T_{y_n}B_1.
$$
This is equivalent to the surjectivity of
$$
D_{u,H(t,y_n)}\,\oplus\, \Psi\scirc dH\,\oplus\,p_B:
\cale_{u,H(t,y_n)}\,\oplus\, T_{(t,y_n)}\bigl([0,1]{\times} B_1\!\bigr)
\to \cale'_{u,H(t,y_n)} \oplus T_{y_n}B_1
$$
and hence to the the surjectivity of
$$
D_{u,h_n(t)}\,\oplus\, \Psi\scirc dh_n:
\cale_{u,h_n(t)}\,\oplus\, T_t[0,1]
\to \cale'_{u,h_n(t)}.
\eqno(8.3.3)
$$
Consequently, $\dimr\sfh^1_D(S, N_M)\le1$ for any $(M,t)\in \calm_{h_n}\cap
V_1$.
\qed

\state Corollary 8.3.4.
\it Under the conditions of Theorem 8.3.3 suppose additionally that
$S$ is a sphere $S^2$. Then for all points $(M,t)\in \calm^0_{h_n}$ the
associated $D_N$-operator is surjective, i.e., $\sfh^1_{D_N}(S^2,N_{M_t})=0$.

Moreover, $\calm_{h_n}$ is a trivial bordism: $\calm_{h_n(0)}\times [0,1]$.
In particular, for every $h_n(0)$-holomorphic sphere $M_0\in \calm_{h_n(0)}$
there exists a continuous family of $h_n(t)$, holomorphic spheres $M_{n,t}=
u_{n,t}(S^2)$ with $M_{n,0}=M_0$.

\smallskip\noindent\sl Proof. \rm
Suppose that for some $(M,t)\in \calm^0_{h_n} $ we have
$\sfh^1_{D_N}(M,N_M)\not= 0$. Then by (iv) of {\sl Theorem 8.3.3}
\, $\sfh^1_{D_N}(M,N_M)= 1$. But this contradicts {\sl Theorem 7.3.1} for the
case $S=S^2$ and $L=N_M$.

Let for $(M, t)\in \calm^0_{h_n}$ we have $M =u(S)$ and $J=h_n(t)$. Let
also $\dot J \not= 0 \in dh_n(T_t[0,1])$. Then by {\sl Lemma 8.2.1} and  {\sl
Corollary 8.2.3} the tangent space $'_{(M,t)} \calm^0_{h_n}$ is canonically
isomorphic to

%%%%%%%%%%%%%%%%%%%%Question%%%%%%%%%%%%%%%%%

$$
\ker\!\Bigl(\!D_{u,J}\,\oplus\, \Psi:
\cale_{u,J}\,\oplus\, \rr\,\dot J\!\bigr)
\to \cale'_{u,J} \!\Bigr)
\!\!\!\Bigm/ \!\!\!du(\sfh^0(S\!\!,\! TS)),
$$
and the differential of the projection $d\pi_{h_n}: '_{(M,t)} \calm^0_{h_n}
\to T_t[0,1] \cong \rr$ is of the form $d\pi_{h_n}[v,a\dot J] = a$. While in
the case $S=S^2$
the space $\sfh^1(S, TS)$ is trivial, {\sl Corollary 8.3.4} insures the
surjectivity of the operator $D_{u,J}: \cale_{u,J} \to \cale'_{u,J}$. Thus
for $a\not=0 \in \rr$ there exists such a $v\in \cale_{u,J}$  that
$[v,a\dot J]\in '_{(M,t)} \calm^0_{h_n}$. This means that for every
$(M, t)\in \calm^0_{h_n}$ the projection $d\pi_{h_n}: '_{(M,t)}
\calm^0_{h_n} \to T_t[0,1]$ is surjective. While the manifold
$\calm^0_{h_n}$ is compact, there is a diffeomorphism $\calm_{h_n} \cong
\calm_{h_n(0)}\times [0,1]$.
\qed

%But ${\sf ind}D_N$ is even (being equal to the index
%of its complex linear part, \ie $\dbar $). So $\sfh^0_{D_N}(M,N_M)\not= 0$.
%At the same time by the Serre duality for $D$-cohomologies, Lemma 4.7,
%$\sfh^1_{D_N}(M,N_M)^*\cong
%\sfh^0_{D_N}(M,N^*\otimes \Omega )\not= 0$. $\sfh^1_{D_N}(M,N_M)^*\cong
%\sfh^0_{D_N}(M, $
%
%\noindent $N^*\otimes \Omega )\not= 0$.
%
%
%Let $\xi $ be a nonzero section of
% $\sfh^0_{D_N}(M,N_M)$ and $\eta $ a nonzero section of
%
%\noindent $\sfh^0_{D_N}(M,N^*$ $\otimes \Omega )\not= 0$. By (1.11) we have
%$$
%\dbar <\xi ,\eta > = <D_N\xi ,\eta > + <\xi ,D^*_N\eta > = 0.
%$$
%\noindent
%This means that $<\xi ,\eta >$ is a nonzero holomorphic form on $S^2$.
%Contradiction.
%\smallskip
%\qed

\medskip\noindent\sl
8.4. Moduli of Parameterized Curves.

\medskip\rm
Sometimes it may be useful to consider the moduli spaces of parameterised
complex curves. Therefore, in this paragraph $\calp $ will denote the space
of parameterized complex curves on $X$. Suppose that  $(u, J_S, J) \in \calp$, {\sl i.e.,} $u$ belongs to
$C^1(S, X)$ and satisfies the equation $du\scirc J_S = J \scirc
du$. Set
$$
T_{(u, J_S, J)}\calp\deff \bigl\{\, (v,\dot J_S,\dot J) \in T_u\cals
\times T_{J_S} \calj_S \times T_J\calj \, :\,
2D_{u, J}v  +  \dot J \scirc du \scirc J_S +
\dot J_S \scirc du \scirc J =0 \,\bigr\}.
$$
Let $\pr_\calj: \cals \times \calj_S \times \calj \to \calj$
and  $\pr_{(u, J_S, J)} : T_{(u, J_S, J)} \calp \to T_J\calj$ denote
the natural projections.

\smallskip\nobreak
\state Theorem 8.4.1. {\it The map $\pr_{(u, J_S, J)} : T_{(u, J_S, J)}\calp
\to T_J\calj$ is surjective \iff\/  $\sfh^1_D(S, N_0)$ $=0$, and then the
following hold:

\sli  the kernel $\ker(\pr_{(u, J_S, J)})$ admits a closed complementing
space;

\slii for $(\tilde u,\tilde J_S,\tilde J)\in \calp$ close enough to
$(u, J_S, J)$ the projection $\pr_{(\tilde u,\tilde J_S,\tilde J)}$ is
also surjective;

\sliii for some neighborhood $U\subset \cals \times \calj_S
\times \calj$ of $(u, J_S, J)$ the set $\calp \cap U$ is a Banach
submanifold of $U$ with the tangent space $T_{(u, J_S, J)}\calp$
at $(u, J_S, J)$;

\sliv there exists a $C^1$-map $f$ from some neighborhood $V$ of
$J\in \calj$ into $U$ with $f(V)\subset \calp$, $f(J) = (u,J_S,J)$
and $\pr_\calj
\scirc f = \id_V$ such that $\im (df: T_J\calj \to
T_{(u, J_S, J)}\calp)$ is complementing to  $\ker(\pr_{(u, J_S, J)})$.
}

\state Proof. Denote by $\tilde A$ the set of a singular point of 
$M=u(S)$, {\sl i.e.}, the set of cuspidal and self-intersection points of $M$.
Let $\xi \in L^p_{(0, 1)} (S, E)$. Using the local
solvability of the $D$-equation (proved essentially in {\sl Lemma 3.2.1})
and an appropriate partition of unity, we can find $\eta\in L^{1,p}(S, E)$
such that $\xi -2D\eta \in C^1_{(0, 1)}(S, E)$ and $\xi -2D\eta =0$ in a
neighborhood of $\tilde A$. Then we find $\dot J\in T_J\calj \equiv
C^1_{(0, 1)}(X, TX)$ such that $\dot J \scirc du \scirc J_S =\xi-
2D\eta$. The surjectivity of $\pr_{(u, J_S, J)}$ and {\sl Theorem 7.3.1} easily
yield the identity $\sfh^1_D(S, N_0)=0$.

\smallskip\sl
From now on and until the~end of the~proof we suppose that 
$\sfh^1_D(S, N_0)=0$.

\smallskip\rm
Let $ds^2$ be some Hermitian metric on $S$ and let $\calh_{(0, 1)} \subset
C^1_{(0, 1)}(S, TS)$ be a space of $ds^2$-harmonic $TS$-valued
$(0, 1)$-forms on $S$. Then the~natural map $\calh_{(0, 1)} \to
\sfh^1(S, TS)$ is an isomorphism. Furthermore, due to {\sl Corollary 7.3.2},
the~map
$$
g\deff(du, 2D): \calh_{(0, 1)} \oplus L^{1,p}(S, E) \to
L^p_{(0, 1)}(S, E)
$$
is surjective and has a finite-dimensional  kernel. Consequently, there exists
a~closed subspace $Y\subset \calh_{(0, 1)} \oplus L^{1,p}(S, E)$
such that $g : Y \to L^p_{(0, 1)}(S, E)$ is an isomorphism. Let
$$
h=(h_{TS},h_E): L^p_{(0, 1)}(S, E) \to Y\subset \calh_{(0, 1)}
\oplus L^{1,p}(S, E)
$$
be the inversion of $g\vert_Y$.

Take $\dot J\in T_J\calj= \{\, I\in C^1(X, {\sl End}(TX))\, :\,
JI+IJ=0\,\} \equiv C^1_{(0, 1)}(X, TX)$. Then $\dot J \scirc du
\scirc J_S$ lies in $C^0_{(0, 1)}(S, E)$. Let $h(\dot J \scirc du
\scirc J_S) =(\xi,\eta)$ with $\xi\in \calh_{(0, 1)}
\subset C^1_{(0, 1)}(S, TS)$ and $\eta\in L^{1,p}(S, E)$.
Then we obtain 
$$
2D(-\eta) + \dot J \scirc du \scirc J_S +
J \scirc du ( J_S \xi)=0,\eqno(8.4.1)
$$
where we use the identity $J\scirc du \scirc J_S =-du$. Using that 
$T_{J_S}\calj_S =$ $ C^1_{(0, 1)}(S, TS)$,
we conclude that the formula $F(\dot J)= \bigl(J_S h_{TS}(\dot J \scirc du
\scirc J_S), -h_E(\dot J \scirc du \scirc J_S),\dot J \bigr)$ defines a 
{\sl bounded} linear operator $F: T_J\calj \to T_{(u, J_S, J)}\calp$ such that
$\pr_{(u, J_S, J)}\scirc F = \id_{T_J\calj}$. In particular, $\sfh^1_D(S,
N_0)=0$ implies the surjectivity of $\pr_{(u, J_S, J)}$.

\smallskip
The image of the just defined operator $F$ is closed, because the convergence
of $F(\dot J_n)$, $J_n\in T_J\calj$ obviously yields the convergence of
$J_n=\pr_{(u, J_S, J)}\scirc F(\dot J_n)$. One can easily see that $\im(F)$ is
a closed complementing space to $\ker(\pr_{(u, J_S, J)})$.

\smallskip
Further, for $(\tilde u, \tilde J_S,\tilde J) \in \calp$ close enough to 
$(u, J_S, J)$ the~map $\tilde g \deff (d\tilde u, 2D_{\tilde u,\tilde J}): 
Y \to L^p_{(0, 1)}(S,\tilde E)$ is also an~isomorphism. This implies
the~surjectivity $\pr_{(\tilde u,\tilde J_S,\tilde J)}$.

\smallskip
The~statements \sliii and \sliv can easily be obtained from {\sl ii)} and
the implicit function theorem.
\qed

\smallskip\smallskip\noindent\sl
8.5. Gromov Non-squeezing Theorem.

\smallskip\rm
We shall finish this chapter with the proof of the Gromov 
non-squeezing theorem.

Consider an infinite cylinder $\zz (\lambda ):=\Delta (\lambda )
\times \rr^{2n-2}$ in
$\rr^{2n}$. Here $\Delta (\lambda )$ is a disk of radius 
$\lambda $ in the $x_1,y_1$-plane,
and on $\rr^{2n-2}$ we fix coordinates $x_2,y_2,...,x_n,y_n$. By $\bb (R)$
we denote the ball of radius $R$ in $\rr^{2n}$. We consider on $\rr^{2n}$ the
standard symplectic form $\omega =\Sigma_{i=1}^ndx_i\wedge dy_i$.

Let $D\subset \rr^{2n}$ be a domain. Recall that a smooth map
$f:D\to \rr^{2n}$ is called a symplectomorphism if $f^*\omega =\omega $. In
particular, symplectomorphisms preserve the Euclidean volume.

\smallskip
\state Exercise. {\rm Find a volume preserving a linear map from $\bb (R)$ to
$\zz (\lambda )$ for arbitrary  $R$ and $\lambda $.
}

The following theorem shows how far symplectomorphisms are from 
volume preserving maps.

\smallskip
\state Theorem 8.5.1. {\it (Gromov Non-squeezing Theorem). If there exists
a symplectomorphism $f:\bb (R) \to \zz (\lambda )$, then $R\le \lambda$}.

\state Proof. {\rm Taking a smaller ball, if nessessary, and translating
the image, we can suppose that $f(\bb (R))$ is contained in some compact
$\bar\Delta(0,\lambda')\times [0,a]^{2n-2}$ of $\Delta (0,\lambda )\times
\rr^{2n-2}$. Factorizing $\rr^{2n-2}$ by the lattice $a\cdot \zz^{2n-2}$,
we observe that $\overline{f(\bb (R))}\comp \Delta (0,\lambda )\times
T^{2n-2}$. Moreover, our syplectic form $\omega $ descends onto this
factor and still can be decomposed as $\omega = \omega_1 + \omega_2$. Here
$\omega_1 = dx_1\wedge dy_1$ and $\omega_2=\Sigma_{i=2}^ndx_i\wedge dy_i$.

Further, we compactify the disk $\Delta (0,\lambda )$ to a two-dimensional
sphere $S^2$ by an $\eps $-disk $\Delta (0,\eps )$ and extend
$\omega_1$ to a smooth strictly positive $(1,1)$-form, still denoted by
$\omega_1$ on $S^2$, and having
$$
\int_{S^2}\omega_1 = \pi \lambda^2 + \eps . \eqno(8.5.1)
$$
Note that now we have a symplectic imbedding $f:\overline{\bb (R)}
\to S^2\times T^{2n-2}$. Also note that our symplectic form
$\omega $ on $S^2\times T^{2n-2}$  tames the standard complex
structure $J\st$ on this manifold and satisfies
$$
\int_{[\gamma ]}\omega = \pi \lambda^2 + \eps , \eqno(8.5.2)
$$
\noindent where $[\gamma ]=[S^2\times \{ pt\} ]$.
We further observe that $\calm_{[\gamma ],J\st}$ consists of $\{ \cc\pp^1
\times \{ a\} :a\in T^{2n-2}\} $, i.e., this moduli space is diffeomorphic
to a torus $T^{2n-2}$.

\state Lemma 8.5.2. \it (a) For a generic almost-complex structure $J$ on
$S^2\times T^{2n-2}$   tamed by $\omega $ the moduli space
$\calm_{[\gamma ],J}$  is diffeomorphic to $T^{2n-2}$.

\smallskip\noindent
(b) Moreover, for any
(not nessessarily generic!) $J$ and for any point $p\in S^2\times T^{2n-2}$
there exists a $J$-complex curve from $\calm_{[\gamma ],J}$ passing
through $p$.

\smallskip
\state Proof. {\rm (a) For any  $J\in \calj_{\omega }$, any $J$-complex curve
$C$
representing $[\gamma ]$ obviously cannot be decomposed as $C=C_1\cup C_2$,
where $C_k$ are $J_k$-complex for some $J_k\in \calj_{\omega }$. Therefore, we
can apply Theorem 3.2 from the introduction to this chapter. Thus (a) is
proved.

\noindent (b) In fact, for a generic $J$ we have a little bit more. Namely,
for a generic curve $h:[0,1]\to \calj_{\omega }$ with $h(0)=J\st $ consider
an evaluation map $\ev :\calc_h\to S^2\times T^{2n-2}$, where $\calc_h$
is the universal $\cc\pp^1$-bundle over $\calm_{h}$. Recall that
$\calm_h$ is a manifold diffeomorphic to $\calm_{h(0)}\times [0,1]$. This
yields that  $\calc_h$ is diffeomorphic to $\calc_{h(0)}\times [0,1]$.

Now, from the uniqueness properties of complex curves it easily follows 
that for any $t\in [0,1]$ the map $\ev :\calc_{h(t)}\to
S^2\times T^{2n-2}$ is surjective. Thus, we prove 
$\ev :\calc_{J}\to
S^2\times T^{2n-2}$ is surjective for generic $J\in \calj_{\omega }$.

Since generic structures are dense in  $\calj_{\omega }$, we immediately
obtain that the evaluation map $\ev :\calc_{h(t)}\to
S^2\times T^{2n-2}$ is surjective for all $J\in \calj_{\omega }$.
}

\smallskip\rm
This completes the proof of the lemma.

\smallskip
Let us return to the proof of our theorem. Using {\sl Proposition 1.2.1}
from {\sl Lecture 1}, we can extend the structure $f_*J\st $ from
$f(\bb (R))$ onto $S^2\times T^{2n-2}$ to a $\omega $-tamed structure
$J$. By {\sl Lemma 8.5.2} there is a $J$-complex rational curve $C\ni f(0)$.
Remark that
$$
\int_C\omega = \pi \lambda^2 + \eps , \eqno(8.5.3)
$$
\noindent because $[C]=[S^2\times \{ pt\} ]$. While $f$ is a
symplectomorphism,
$$
\int_{f^{-1}(C\cap f(\bb (R))}\omega\st  \le \pi \lambda^2 + \eps .
\eqno(8.5.4)
$$
But $f^{-1}(C\cap f(\bb (R))$ is a complex curve in $\bb (R)\subset \cc^n$
passing through zero. By a well-known estimate due to Alexander-Taylor-Ullman,
see [A-T-U], the area of such a curve is at least $\pi R^2$. So $\pi R^2\le
\pi \lambda^2 +\eps $.
\qed

\medskip\noindent
{\bigsl 8.6. Exceptional Spheres in Symplectic $4$-Manifolds.}

\smallskip\rm
{\sl Corollary 8.3.4} enables us to prove an interesting result about 
exceptional spheres in symplectic $4$-manifolds. Recall that a smooth 
rational curve $C$ in a smooth  complex surface $X$ is called exceptional if 
$[C]^2=-1$. Such curve can be contructed to a point, i.e., there exists a 
smooth complex surface $Y$ and a holomorphic map $h:X\to Y$ such that 
$h(C)=y$ is a point and $h\mid_{X\setminus C}:X\setminus C\to Y\setminus 
\{ y\} $ is a biholomorphism.

Consider now a symplectic $4$-manifold $(X,\omega )$ and let $M$ be a 
symplectic $2$-sphere imbedded to $X$, i.e., there is an imbedding 
$u:\ss^2\to M \hookrightarrow X$ with $u^*\omega$  nowhere zero.

\state Definition 8.6.1. {\it Call $M$ exceptional if $[M]^2=-1$.
}

\smallskip
Using {\sl Lemmas 1.4.2 and 1.4.3} we can construct an 
almost-complex structure 
$J\in \calj_{\omega }$, which is moreover integrable in the neighborhood of 
$M$ and such that $M$ becomes $J$-complex. Now the complex analytic statement 
mentioned above allows us to construct this $M$ to a point.

\state Corollary 8.6.1. {\it Take a maximal, linearly independant in 
$\sfh_2(X,\zz )$ system of 
exceptional symplectic spheres $M_1,...,M_k$ in $X$. Then 
there exists $J\in \calj_{\omega }$ making some symplectic spheres $\tilde 
M_1,...,\tilde M_k$, with $\tilde M_j$ being isotopic to $M_j$ for  
$1\le j\le k$, all $J$-complex.
}
\state Proof. {\rm  
For $k=1$ this is the statement of {\sl Lemma 1.4.2}. Suppose we prove this 
statement for $k-1$ exceptional spheres. Take our spheres $M_j,j=1,...,k$ and 
find a structure $J$ such that $\tilde M_j$ are isotopic to $M_j$ and $J$-
complex, $j=1,...,k-1.$ The normal Gromov operator $D_{J,\tilde M_j}$ is 
surjective on all $\tilde M_j,j=1,...1k-1$ by  {\sl Corollary 7.2.3} and 
thus by {\sl Theorem 7.3.1} and {\sl Corollary 8.2.3} the projection 
$\pi_{\calj }:\calm \to \calj $ is a diffeomorphism in the neighborhood 
of $\tilde M_j$ and $J$ for all $j=1,...,k-1$.

Take a structure $J_1$ such that $M_k$ becomes $J_1$-complex. Again the 
projection $\pi_{\calj }$ is regular in the neighborhood of $M_k$. Using 
{\sl Corollary 8.2.4}, we take a generic path $h$ starting close to $J$, say 
at $J´$  
and ending close to $J_1$ such that all $\calm_h(\tilde M_j´,J´)$ are 
compact and open in corresponding moduli spaces. Here $\tilde M_j´$ 
are $J´$-complex spheres isotopic to $\tilde M_j$.  Applying 
{\sl Corollary 8.3.4}, we obtain an isotopy of 
$\tilde M_j´$ to a $J''$-complex 
spheres $\tilde M_j''$ with $J''$ close to $J_1$. It remains  
to take a $J''$-complex sphere $\tilde M_k''$ isotopic to $M_k$.

\smallskip
\qed
}

\newpage
%%%%%%%%%%% chapter 4 %%%%%%%%%%%%%%

%\magnification=\magstep1
%\input komp.def

% \magnification=\magstep1
%\input lect.def

%plain TeX with amssym files
%version of 12.03.96
%\baselineskip=12.5pt plus 1.7pt
\noindent{\bigbf Chapter IV. Envelopes of Meromorphy of Two-spheres.}

\smallskip
This chapter is devoted to the study of  envelopes of
meromorphy of neighborhoods of two-spheres in  complex algebraic surfaces.
The original question which motivated our studies was asked by A. Vitushkin.

\smallskip\it
Let $M$ be a ``small'' perturbation of the complex line in $\cc\pp^2$. Does
there exist a nonconstant holomorphic function in the neighborhood of such
an $M$?

\smallskip\rm
 It was asked as a test question on the way to
searching for the solution to the Jakobian conjecture, and the answer, as
expected, was negative.

Let us briefly test the possible approaches to the answer to this question.
It is more or less clear that one should try to extend holomorphic 
(or meromorphic via expected nonexistence of holomorphic) functions onto 
the whole $\cc\pp^2$ and then conclude

\smallskip\noindent
{\sl 1.} If $M$ is a complex line itself, then nonconstant 
holomorphic functions do not exist in any neighborhood of $M$ for the 
following reason. First note that in this case  $\cc\pp^2\setminus M=\cc^2$. 
Let $B_N$ denote the closed 
ball of radii $N$ in $\cc^2$. Then $V_N:=\cc\pp^2\setminus B_N$ is a 
fundumental system of strictly pseudoconcave neighborhoods of $M$. Any 
function holomorphic in $V_N$ holomorphically extends onto the ball $B_N$ by 
Hartogs' theorem, and thus becomes holomorphic on the whole $\cc\pp^2$, 
i.e., constant.

One can try to construct such an exhaustion for any $M$. But the generic $M$ 
is totally real outside of three positive elliptic points (this follows from 
the Lai formulae and cancellation theorem of Kharlamov-Eliashberg, 
see Appendix IV)
and one can check that a totally real disk has no small concave tubular 
neighborhoods.

\smallskip\noindent
{\sl 2.} As was just pointed out, a generic perturbation $M$ will have 
exactly three elliptic points with complex tangents. This makes it problematic 
to attempt to construct a family of complex disks (or any other Riemann 
surfaces) with boundary on $M$ starting from our elliptic points. Such an
approach to the construction of the holomorphic envelope of $M$ seems to 
be perfect only for specially imbedded spheres, see [B-G],[B-K],[E], [Sch]
and [F-M].

%%%%%%%%%%%%%%%%%%%%question%%%%%%%%%%%%%%%%%%%%%

\smallskip\noindent
{\sl 3.} The method used here to study  envelopes of meromorphy of 
spheres  is based on Gromov's theory of pseudoholomorphic curves. We 
remark first of all that a small perturbation of a complex sphere is 
symplectic:

\state Definition. \sl 
A $C^1$-smooth immersion $u\:S\to(X,\omega)$ of a real surface $S$ into a 
symplectic manifold $(X,\omega )$ 
is called {\it symplectic} if $u^*\omega$
does not vanish anywhere on~$S$.
\rm

Thus, we shall study envelopes of meromorphy
of neighborhoods of  two-spheres symplectically immersed
in complex surfaces. Here a {\it complex surface\/} means
a (Hausdorff) connected complex two-dimensional manifold~$X$
countable at infinity.

The idea is to perturb the complex structure in the given neighborhood of 
$M$ in such a way that $M$ becomes complex; see {\sl Lemma 1.4.2}, where 
an appropriate family $J_t$ of tamed almost complex structures is 
constructed. Using results from Chapter III, we then construct a family $M_t$ 
of $J_t$-complex spheres and try to extend functions along this family.
In more formal language to extend "along a family" means "onto the envelope 
of meromorphy".

Let $U$ be a domain in~$X$.  Its envelope of 
meromorphy $(\widehat U,\pi)$ is the maximal domain over~$X$ satisfying the 
following conditions:  

(i) there exists a holomorphic 
embedding $i\:U\to\widehat U$ with $\pi\circ i=\id_U$;

(ii)
each meromorphic function $f$ on $U$
extends to a meromorphic function $\widehat f$ on~$\widehat U$,
that is, $\widehat f\circ i=f$.

The envelope of meromorphy exists for each domain $U$.
This can be proved, for example,
by applying the Cartan--Thullen method
to the sheaf of meromorphic functions on $X$, see 
[Iv-1].

In the sequel we shall restrict ourselves to K\"ahler
complex surfaces, that is, we assume that  $X$
carries a strictly positive closed $(1,1)$-form~$\omega$.

The aim  of the present chapter is to prove the following result.

\state Theorem 4.1. \it
Let $u\:S^2\to X$ be a symplectic immersion of the two-sphere $S^2$
in a disk-convex K\"ahler surface~$X$
such that $M:=u(S)$ has only positive double points.
Assume that $c_1(X)[M]>0$.
Then the envelope of meromorphy $(\widehat U,\pi)$
of an arbitrary neighborhood $U$ of $M$
contains a rational curve $C$ with $\pi^*c_1(X)[C]>0$.
\rm

The definition of disk-convexity 
is given in Lecture 10.
At this point, we only observe that all compact manifolds are disk-convex.
As usual,  $c_1(X)$ is the first Chern class of~$X$. In Appendix IV, 
Corollary A4.3.2 we explain why the condition $c_1(X)[M]>0$ is nessessary.

 The exposition includes a construction of a 
complete family of holomorphic deformations of a non-compact complex curve 
in a complex manifold, parameterized by a finite codimension analytic subset 
of  a Banach ball. This will be carried out in Lecture 9.   

  The existence of this family is used to prove a
generalization of Levi's continuity principle, which is applied  to describe
envelopes of meromorphy.

\smallskip\noindent
{\sl 4.}
Another natural approach is due to S. Nemirovski. Suppose that 
for some neighborhood $V\supset M$ there is a nonconstant holomorphic 
function in $V$. Then from the result of Fujita, [Fu], it follows that there 
exists a Stein domain over $\cc\pp^2$ which contains $M$. One can imbed
by the Stout theorem, [St], some neighborhood of $M$ into a compact algebraic 
surface $X$ and observe that $b_+^2(X)>1$. The Seiberg-Witten theory imposes 
an adjunction inequality on $M$:
$$
\vert c_1(X)[M]\vert  + [M]^2\le 0,
$$
see [K-M]. But $c_1(X)[M]=c_1(\cc\pp^2)[M]=3$ and $[M]^2=1$, a contradiction.

\smallskip 
This method seems to work only in some special algebraic 
surfaces like $\cc\pp^2$ or $\cc\pp^1\times \cc\pp^1$ but does not 
require that $M$ be symplectic! We shall give more details in Appendix V.

\newpage
\bigskip\bigskip\noindent
{\bigbf Lecture 9}

\smallskip\noindent
{\bigbf Deformation of Noncompact Curves.}

\medskip\noindent
{\bigsl 9.1. Banach Analytic Sets.}

\smallskip\rm
Now we shall study a process of local deformation of noncompact complex 
curves in complex manifolds. We shall see that the moduli have nice analytic 
structure. Let start from the following

\smallskip
\state Definition 9.1.1. {\it By a {\sl Banach ball } we mean a ball in some
complex linear 
Banach space. A subset $\calm$ of the Banach ball $B$ is said to be a Banach
analytic set of the finite codimension, b.a.s.f.c., if there
exists a holomorphic map $F:B\to \cc^N$ with $N<\infty$ such that  $\calm =
\{ x\in B: F(x)=0\}$.

\smallskip\rm
The importance of this notion consists in the fact that contrary to the
general Banach analytic sets b.a.s.f.c., their properties are similar
to the finite dimensional analytic sets. Namely, the following is true.

\state Theorem 9.1.1. ([Ra])
\it Let $B$ be a ball in a Banach space  $\cal F$, $\calm \subset B$
--- a b.a.s.f.c. \ and  $x_0$ a point in $\calm $. Then there exists a
neighborhood $U\ni x_0$ in $B$ such that $\calm \cap U$ is a union of
finite number of irreducible b.a.s.f.c. \ $\calm_j$.

Moreover, each $\calm_j$ can be represented as a finite ramified covering
over a domain in some closed linear  subspace  ${\cal F}_j
\subset {\cal F}$ of finite codimension.
\rm

\medskip The aim of this paragraph is to prove an existence of a complete
family of holomorphic deformations of a stable curve over a complex manifold
 $X$,
which is parameterized by a b.a.s.f.c. Before stating the result, let us
introduce the following definition.

Let $C$ be a nodal curve, $E$ a holomorphic
vector bundle over $C$, and $C= \cup_{i=1}^l C_i$ a decomposition of $C$
into irreducible components. Suppose that $E$ extends sufficiently
smoothly to the boundary $\d C$.

\state Definition 9.1.2. {\it Let us define an {\sl $L^{1,p}$-section  $v$}
of the bundle $E$ over $C$ as a couple $(v_i)_{i=1}^l$ of $v_i \in
L^{1,p}(C_i,E)$ such that at any nodal point $z \in C_i \cap C_j$
$v_i(z) = v_j(z)$.
Let us define an {\sl $E$-valued $L^p$-integrable $(0,1)$-form  $\xi$ on
$C$} as a couple $(\xi_i)_{i=1}^l$ of $(0,1)$-forms $\xi_i \in L^p(C_i,
E\otimes \Lambda^{(0,1)} )$.
}

\smallskip
Let $L^{1,p}(C,E)$ denote a Banach space of $L^{1,p}$-sections of $E$
over $C$, and  $L^p(C, E\otimes \Lambda ^{(0,1)} )$ a Banach space
of $L^p$-integrable $(0,1)$-forms $C$. Denote by $\calh^{1,p}(C,E)$ a
Banach space of {\sl holomorphic}  $L^{1,p}$-sections of $E$ over $C$.

Analogously, for a complex manifold $X$ by $L^{1,p}(C,X)$ we denote the
set of couples $u=(u_i)_{i=1}^l$ such that  $u_i \in L^{1,p}(C_i,
X)$, with $u_i(z) = u_j(z)$ at any nodal point $z \in C_i \cap C_j$.
One can see that $L^{1,p}(C,X)$ is a Banach manifold with a tangent space
$T_uL^{1,p}(C,X) = L^{1,p}(C,u^*TX)$. Further by $\calh^{1,p}(C,X)$
we shall denote the manifold of {\sl holomorphic} $L^{1,p}$-maps from $C$ to $X$.
Note that for $u\in \calh^{1,p}(C,X)$ one has $u(C) \subset u(\barr C)
\Subset X$, because $L^{1,p}\subset C^{0,1-{2\over p}}$ and therefore 
$u$ is continuous up to the boundary.

\smallskip

\smallskip\noindent\sl
9.2. Solution of a Cousin-type Problem.

\smallskip\rm
An important role plays  the following result on the possibility to 
 solve the following  Cousin-type problem.

\smallskip
\noindent
\bf Lemma 9.2.1. \it
Let $C$ be a nodal curve and $E$ a holomorhpic vector bundle
over $C$, $C^1$-smooth up to the boundary. Let $\{ V_i\}_{i=1}^l$ be
a finite covering of $C$ by Stein domains with piecewise smooth
boundaries. Put $V_{ij} \deff V_i \cap V_j$ and suppose that all triple
intersections $V_i \cap V_j \cap V_k$ á $i \not= j \not= k \not= i$ are
empty.

Then for $2 \le p< \infty$ the \v{C}ech-differential
$$
\matrix
\delta :& \sum_{i=1}^l \calh ^{1,p}(V_i, E)& \lrar&
\sum_{i<j} \calh ^{1,p}(V_{ij}, E)
\cr
\delta :& (v_i)_{i=1}^l        & \longmapsto & (v_i -v_j)
\endmatrix
\eqno(9.1.1)
$$
possesses the following properties:

\sli the image $\im(\delta)$ is closed and has finite codimension; moreover,
$\coker(\delta) = \sfh^1(C, E) = \sfh^1(C_{\sf comp}, E)$, where
$C_{\sf comp}$ denotes the union of compact irreducible components of $C$;

\slii the kernel $\ker(\delta)$ is isomorphic to $\calh ^{1,p}(C, E)$ 
and admits a closed direct complement.
\rm

\state Proof. Before considering the \v{C}ech complex, let us look at the
corresponding $\dbar$-problem. Consider the following operator
$$
\dbar : L^{1,p}(C,E) \lrar L^p(C, E\otimes \Lambda^{(0,1)}_C).
\eqno(9.1.2)
$$
First we shall prove that this operator possesses the same properties as
$(9.1.1)$, i.e., that  $\ker(\dbar)$ admits a closed direct complement,
$\im(\dbar)$ has finite codemension and is closed, and  $\coker(\dbar) =
\sfh^1(C, E) = \sfh^1(C_{\sf comp}, E)$. Moreover, we shall construct a
natural isomorphism between (co)kernels of  $(9.1.1)$ and $(9.1.2)$.

\smallskip
Since the boundary $C$ is smooth, there exist nodal curves $C^+$ and
$C^{++}$ such that $C \Subset C^+ \Subset C^{++}$ and differences
$C^+ \bs \barr C$ (resp.\ $C^{++} \bs \barr C^+$) consist of
an annulai $A^+_\alpha$ (resp.\ $A^{++}_\alpha$) adjacent to the
corresponding components $\gamma_\alpha$ (resp.\ $\gamma^+ _\alpha$) of the
boundary $\d C$ (resp.\ $\d C^+$, see.\ Fig.~17.). Now $E$ extends to a
holomorphic vector bundle over $C^{++}$, again denoted as $E$.

\bigskip
\vbox{\nolineskip\rm\xsize.55\hsize%
\putt[1.05][0]{\advance\hsize-1.05\xsize\parindent=0pt%
\bigskip
\centerline{Fig.~17.\ \ $C\Subset C^+ \Subset C^{++}$.}
\medskip%\def\script{\scriptstyle}%
%\sevenrm
Boundaries of $C^{++}$, $C^+$ and $C$ are marked correspondingly by a bold,
punctured and broker lines.
\smallskip}%
\noindent
\epsfxsize=\xsize\epsfbox{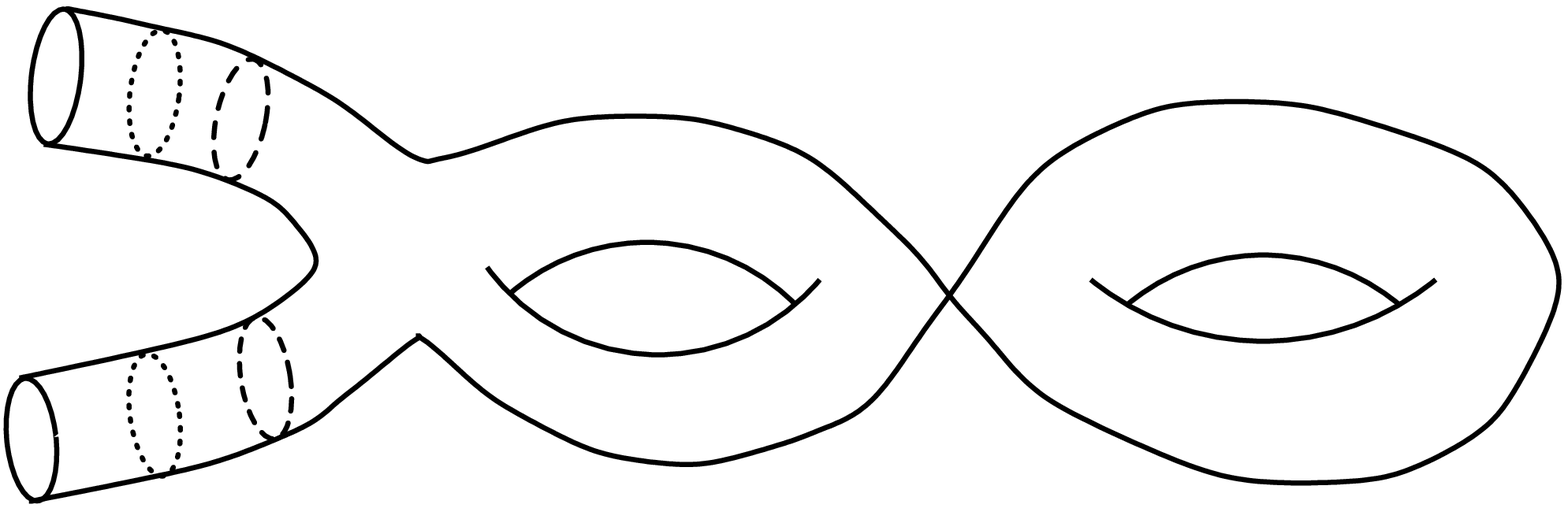}%
}

\bigskip

Consider the following sheaves on $C^{++}$
$$
\matrix \format \r\;\; &\c\; &\; \c\; &\; \l \\
L^{1,p}_\loc(\cdot ,E): &
V & \mapsto &  L^{1,p}_\loc(V,E) \vphantom{\int_3};
\cr
L^p_\loc(\cdot ,E\otimes \Lambda^{(0,1)}_{C^{++}}): &
V & \mapsto &  L^p_\loc(V,E\otimes \Lambda^{(0,1)}_{C^{++}}),
\endmatrix
$$
and a sheaf homomorphism defined by the operator
$$
\dbar: L^{1,p}_\loc(V,E) \lrar L^p_\loc(V,E\otimes \Lambda^{(0,1)}_{C^{++}}).
$$
Then sheaves  $L^{1,p}_\loc(\cdot ,E)$ and $L^p_\loc(\cdot ,E\otimes \barr
\Omega_{C^{++}})$ together with a  $\dbar$-homomorphism form a fine
resolution of the (coherent) sheaf $\calo^E$ of holomorphic sections of
$E$ over $C^{++}$. In smooth points of $C^{++}$ this fact follows from
$L^p$-regularity of the (elliptic) operator $\dbar$, and in nodal points
this can be seen from the following considerations.

Let $z\in C$ be a nodal point which lies on the intersection of
irreducible components $C_i$ and  $C_j$ of the curve $C$. Let $\xi_i$
(resp.\ $\xi_j$) be an $E$-valued $L^p_\loc$-integrable (0,1)-form,
defined in the neighborhood of the point $z$ in $C_i$ (resp.\ $C_j$).
And let $v_i$ (resp $v_j$) be  $L^{1,p}_\loc$-solutions of the equation
$\dbar v_i = \xi_i$ (resp.\ $\dbar v_j= \xi_j$).
Adding a local {\sl holomorphic} section of $E$ over $C_i$, we can achieve
an equality $v_i(z)= v_j(z)$. Then the pair $(v_i, v_j)$ defines
a section of the sheaf $L^{1,p}_\loc(\cdot ,E)$ in the neighborhood of $z$.
This shows that the Dolbeault lemma is valid for the holomorphic bundle $E$
also in the neighborhood of the nodal points.

\smallskip
This gives the natural isomorphisms
$$
\matrix \format \r&\c & \c & \l &\c\\
\ker{\Bigm(}\dbar: L^{1,p}_\loc(C^{++},E) &\lrar &
L^p_\loc(C^{++},E\otimes \Lambda^{(0,1)}_{C^{++}} )
{\Bigm)} &=& \sfh^0(C^{++}, E),
\vphantom{\int_3}
\cr
%\noalign{\rlap{and}}
\coker{\Bigm(}\dbar: L^{1,p}_\loc(C^{++},E) & \lrar&
L^p_\loc(C^{++},E\otimes \Lambda^{(0,1)}_{C^{++}} )
{\Bigm)} &=& \sfh^1(C^{++}, E).
\endmatrix
\eqno(9.1.3)
$$
Note that one also has the same isomorphisms for $C$ and $C^+$. Note also
that there are natural isomorphisms $\sfh^1(C^{++}, E)= \sfh^1(C^+, E) =
\sfh^1(C, E)= \sfh^1(C_{\sf comp}, E)$, induced by the restrictions
$L^p_\loc(C^{++},E\otimes \Lambda^{(0,1)}_{C^{++}}) \to L^p_\loc(C^+,E\otimes
\Lambda^{(0,1)}_{C^+} ) \to L^p_\loc(C,E\otimes \Lambda^{(0,1)}_{C} ) \to
L^p_\loc(C_{\sf comp},E\otimes \Lambda^{(0,1)}_C) $.

Now take some $\xi \in L^2(C^+,E \otimes \Lambda^{(0,1)}_{C^+})$,
which defines a zero cohomology class in $\sfh^1(C^+, E)$. We can extend
$\xi$ by zero to the element $\ti \xi\in L^2_\loc(C^{++}, E\otimes
\Lambda^{(0,1)}_{C^{++}})$. While $[\ti\xi]_\dbar= [\xi]_\dbar =0$,
there exists $\ti v \in L^{1,2}_\loc(C^{++},E)$ such that $\dbar\ti v =
\ti\xi$. The restriction $v \deff \ti v\ogran_{C^+}$ satisfies
$\dbar v = \xi$ and $v \in L^{1,2}(C^+,E)$. This shows that the image of a
continuous (!) operator
$$
\dbar: L^{1,2}(C^+,E) \lrar L^2(C^+,E\otimes
\Lambda^{(0,1)}_{C^+}) \eqno(9.1.4)
$$
is of finite codimension. From here and from the Banach open mapping
theorem it follows that this image is closed. Furthermore,
because $L^{1,2}(C^+,E)$ is a Hilbert space, the kernel of the operator
$(9.1.4)$ admits a direct complement $Q \subset L^{1,2}(C^+,E)$.
Moreover, operator $(9.1.4)$ maps $Q$ isomorphically onto its image.
Thus, operator $(9.1.4)$ splits, i.e., there exists a continuous operator
$$
\def\dash{\raise.2pt\hbox{-}} T^+:L^2(C^+,E\otimes
\Lambda^{(0,1)}_{C^+}) \lrar L^{1,2}(C^+,E)  \eqno(9.1.5)
$$
such that $\im(T^+)=Q$ and for any  $\xi \in L^2(C^+,E\otimes
\Lambda^{(0,1)}_{C^+})$ with $[\xi]_\dbar=0 \in \sfh^1(C^+, E)$ one has
$\dbar(T^+\xi)=\xi$.

\smallskip
Define an operator $T:L^2(C, E\otimes \Lambda^{(0,1)}_C) \lrar
L^{1,2}(C,E)$ in the following way. For $\xi \in
L^2(C, E\otimes \Lambda^{(0,1)}_C)$ extend $\xi$ by zero to $\ti \xi \in
L^2(C^+,E\otimes \Lambda^{(0,1)}_{C^+})$ and put $T(\xi) =
T^+(\ti\xi)\ogran_{C}$. $T$ is obviously continuous; moreover
$$
\norm{T^+(\ti\xi)}_{L^{1,2}(C^+)} \le c\cdot \norm{\xi}_{L^2(C)}
$$
with constant $c$ independent of $\xi$. By the $L^p$-regularity
for the elliptic $\dbar$-operator (see ex., [Mo]) for $2\le p<\infty$ and

%%%%%%%%%%%%%%%%%%%question%%%%%%%%%%%%%%%

$v \in L^{1,p}_\loc(C^+,E)$ there is an interior estimate
$$
\norm{v}_{L^{1,p}(C)} \le c'\cdot
\bigl( \norm{v}_{L^{1,2}(C^+)} + \norm{\dbar v}_{L^p(C^+)} \bigr)
\eqno(9.1.6)
$$
with the constant $c'$ independent on $v$. From this it follows that for
$2\le p<\infty$ and  $\xi \in L^p(C, E\otimes \Lambda^{(0,1)}_C)$ with
$[\xi]_\dbar=0 \in \sfh^1(C, E)$ one has
$$
\norm{T(\xi)}_{L^{1,p}(C)} \le c''\cdot \norm{\xi}_{L^p(C)}
$$
with a constant $c''$ independent on  $\xi$. This means that the operator $T$
is a splitting off the operator $(9.1.2)$. This means that operator
$(9.1.2)$ possesses the properties \sli and  \slii of {\sl Lemma 9.2.1}.

\medskip
Let us return to the \v{C}ech operator $(9.1.1)$. Fix some partition of unity
$\1=\sum_{i=1}^l\phi_i$, subordinate to the covering
$\{V_i\}_{i=1}^l$ of the curve $C$. Take a cocycle
$w=(w_{ij}) \in \sum_{i<j}\calh ^{1,p}(V_{ij}, E)$. For $i>j$ put
$w_{ij} \deff -w_{ji}$.
Define $f_i \deff \sum_j \phi_j w_{ij}$. Then $f_i \in L^{1,p}(V_i, E)$
and $f_i- f_j = w_{ij}$.  Consequently, $\dbar f_i \in L^p(V_i, E\otimes
\Lambda^{(0,1)}_C)$ and $\dbar f_i = \dbar f_j$ ¢ $V_{ij}$. So,
$\dbar f_i = \xi\ogran_{V_i}$ for the correctly defined $\xi \in L^p(C,
E\otimes \Lambda^{(0,1)}_C)$.  Moreover, $(w_{ij})$ and $\xi$ define
the same cohomology class $[w_{ij}]=[\xi]$ in $\sfh^1(C,E)$.

Suppose additionally that the induced cohomology class $[w_{ij}]$ is
trivial. Put $f\deff T(\xi)$ and $v_i= f_i-f$. Then $v_i \in
L^{1,p}(V_i, E)$, $v_i- v_j = w_{ij}$, and $\dbar v_i= \dbar f_i -\dbar f=0$.
Thus, $v\deff(v_i) \in \sum_{i=1}^l \calh ^{1,p}(V_i, E)$ and $\delta(v)= w$.
From here it follows that the formula $T_\delta : w \mapsto v$ defines an
operator $T_\delta$, which is a splitting of $\delta$. The explicit
construction shows that $T_\delta$ is continuous. This completes the proof
of {\sl Lemma 9.2.1}.
\qed

\smallskip\noindent\sl
9.3. Case of a Stein Curve.

\smallskip\rm

\smallskip
\state Lemma 9.3.1. \it Let $C$ be a Stein nodal curve with a piecewise
smooth boundary, and $X$ a complex manifold. Then

\sli $\calh^{1,p}(C,X)$ possesses a natural structure of a complex Banach
manifold with a tangent space $T_u\calh^{1,p}(C,X) =
\calh^{1,p}(C,u^*TX)$ at $u\in \calh^{1,p}(C,X)$.

\slii If  $C'\subset C$ is again a nodal curve, then the restriction map
$\calh^{1,p}(C, X) \to \calh^{1,p}(C', X)$ is holomorphic, and its
differential at $u\in \calh^{1,p}(C,X)$ is again the restriction map
$\calh^{1,p}(C, u^*TX) \to \calh^{1,p}(C', u^*TX)$, $v \mapsto v|_{C'}$.
\rm

\smallskip\noindent
The \bf proof \rm will be given in several steps.

\smallskip\noindent
{\sl Step 1.}\ \
Suppose first that the image $u(C)$ lies in the complex chart $U \subset X$
with complex coordinates $w= (w_1, \ldots, w_n) : U \buildrel \cong
\over \lrar U' \subset \cc^n$. Then the set $\calh^{1,p}(C,U)$ 
can be
naturally identified with the set $\calh^{1,p}(C,U')$, which is an open set
in the Banach space $\calh^{1,p}(C,\cc^n)$. This defines on $\calh^{1,p}(C,U)$
the structure of a complex Banach manifold with the tangent space
$T_u\calh^{1,p}(C,U) \cong \calh^{1,p}(C,\cc^n) \cong \calh^{1,p}(C,u^*TU)$
at the point $u\in \calh^{1,p}(C,U)$.

Note that if $u_t$, $t\in [0,1]$, is a $C^1$-curve in $\calh^{1,p}(C,U)$,
then the tangent vector $v\in \calh^{1,p}(C,u^*TU)$ to $u_t$ at $u_0$
is given by $v(z)= {\d u \over\d t}(z) \in T_{u(z)}U$. This last
expression does not depend on the choice of the complex coordinates
$w= (w_1,\ldots, w_n) : U \to \cc^n$ in $U$. This implies two things.

First, {\sl the complex structure on $\calh^{1,p}(C, U)$ does not depend on
the choice of the complex coordinates $w= (w_1, \ldots, w_n) : U \to \cc^n$
in $U$}. Second, for $C$ property \slii of the lemma is fulfilled.

Thus, if $u(C)$ is contained in a coordinate chart, then {\sl Lemma 9.3.1}
is proved.

\smallskip\noindent
{\sl Step 2.}\ \ Now fix $u_0 \in \calh^{1,p}(C,X)$ and suppose
that a finite covering $\{V_i\}_{i=1}^l$ of the curve $C$ is chosen such
that first, the conditions of {\sl Lemma 9.2.1} are satisfied, and
second, for every $V_i$ {\sl Lemma 9.3.1} holds, e.g., every $u_0(V_i)$
is contained in some coordinate chart $U_i \subset X$.

Set  $V_{ij}\deff V_i\cap V_j$. Choose  balls $B_{ij} \subset\calh^{1,p}
(V_{ij}, u_0^*TX) \cong T_{u_0}\calh^{1,p}(V_{ij},X)$ such that there
exists a biholomorhpism $\psi_{ij} : B_{ij} \buildrel \cong \over \lrar
B'_{ij}
\subset \calh^{1,p}(V_{ij},X)$ with $\psi_{ij}(0) = u_0|_{V_{ij}}$ and
$d\psi_{ij}
(0) = \id: \calh^{1,p}(V_{ij}, u_0^*TX) \to\calh^{1,p}(V_{ij}, u_0^*TX)$.
Then take  an open sets $B_i \subset \calh^{1,p}(V_i, X)$  such that $u_0|_{V_i}
\in B_i$ and for every $u_i \in B_i$ one has $u_i|_{V_{ij}} \in B'_{ij}$.

Now, the holomorphic mappings  $\phi_{ij} :B_i \to B_{ij} 
\subset \calh^{1,p}(V_{ij}, \allowbreak u_0^*TX)$, $\phi_{ij} : u_i
\mapsto \psi_{ij}\inv (u_i|_{V_{ij}})$ are uniquely defined. 
They define a holomorphic map
$$
\matrix
\Phi: \prod_{i=1}^l B_i & \lrar   & \sum_{i<j} \calh^{1,p}(V_{ij},u_0^*TX)
\cr
(u_i)_{i=1}^l    & \mapsto & \phi_{ij}(u_i) - \phi_{ji}(u_j).
\endmatrix
$$
One can easily see that the map $\Phi$ defines the condition of
compatibility of local holomorphic mappings $u_i : V_i \to X$, namely
$(u_i)_{i=1}^l \in \prod_{i=1}^lB_i$ define a holomorphic map $u: C \to X$
if and only if $\Phi(u_i)=0$. Moreover, the differential $d\Phi$ at
point $(u_0|_{V_i})$ coincides with the co-boundary \^Cech operator $(9.1.1)$.
Since  $C$ is Stein, $\sfh^1(C, u_0^*TX) =0$. From {\sl Lemma 9.2.1}
and the implicit function theorem, we see that  the statements  \sli and \slii
of {\sl Lemma 9.3.1} are fulfilled in some neighborhood of the map
$u_0\in \calh^{1,p} (C, X)$.

\smallskip\noindent
{\sl Step 3.}\ \ One can show that for any Stein nodal curve $C$
and any $u \in \calh^{1,p} (C, X)$ the statements \sli and \slii of
{\sl Lemma 9.3.1} are fulfilled in some neighborhood of the mapping $u$,
applying {\sl Step 2} sufficiently many times. If, for example, $C$
is an annulus $A_{r,R}$, then one can cover it by the narrow
annulai $A_{r_i, R_i}$, $0< R_i - r_i <\!\!< 1$, and every $A_{r_i, R_i}$
with small sectors $V_{ij} = \{ z= \rho
e^{\isl \theta} \in \cc \,:\, r_i<\rho < R_i, \alpha_j < \theta < \beta_j\}$
á $0< \beta_j -\alpha_j <\!\!< 1$. Details are left to the reader.
\qed

\medskip\noindent\sl
9.4. Degeneration to a Node.

\smallskip\rm
One of the diffuculties in constructing {\sl holomorphic} families
of stable curves with boundaries comes from the fact that the moduli
space of complex structures on the noncompact surface with boundary
$\Sigma$ does not have a natural complex structure; moreover it can have an odd
real dimension. If, for example, $\Sigma$ is an annulus, then it is
biholomorhic to the standard annulus $A_{r,1}$ with uniquely defined
$r\in (0,1)$; thus, the interval $(0,1)$ is a moduli space in this case.
In the general case, if $\Sigma$ has genus $g$ and $k$ boundary components,
then the real dimension of the moduli space is equal to $d= 6g -6 + 3k$,
except in four cases, where  $\Sigma$ is either a sphere ($g=0$, $k=0$),
a torus ($g=1$, $k=0$), a disk ($g=0$, $k=1$) or an annulus ($g=0$, $k=2$),
see e.g., [Ab]. Note that exactly in this cases the group of biholomorphic
automorphisms of corresponding complex curves $(\Sigma, J)$ has positive
dimension.

%%%%%%%%%%%%%%%%%%%%%question%%%%%%%%%%%%%%%%%%%%%%%

One can correct the situation by introducing additionally  $k$ parameters,
namely, fixing $k$ marked points, one on each boundary component.
Let $A$ be an annulus with boundary circles $\gamma_0$ and $\gamma_1$
and $X$ be a complex manifold.

\state Theorem 9.4.1. \it There exist complex Banach manifolds
$\calm(A, X)$ and $\calc(A, X)$, a holomorphic projection
$\pi_\calc: \calc(A, X)\to \calm(A,X)$ and holomorphic mappings
$\ev: \calc(A, X) \to X$, $z_1 :\calc(A, X) \to \Delta$, $z_2 :\calc(A, X)
\to \Delta$ and $\lambda_\calm :\calm(A, X) \to \Delta$ with the following
properties:

\sli for any $y\in \calm(A,X)$ the fiber of the projection
$C_y\deff \pi_\calc\inv(y)$ is a nodal curve, parameterized by an annulus
$A$; moreover, the mapping $(z_1, z_2) : C_y \to \Delta^2$ defines
a biholomorphism onto the curve $\{ (z_1, z_2)
\in \Delta^2 \;:\; z_1 \cdot z_2 = \lambda_\calm(y) \}$; in particular,
$C_y$ is either a standard node if $\lambda_\calm(y)=0$, or a holomorphic
annulus $\{ |\lambda_\calm(y)| <|z_1| <1 \}$;

\slii the diagram
$$

\def\mapright#1{\!\!\smash{\mathop{\hbox to50pt{\cleaders\lowminus\hfill}
\mkern-13mu\rightarrow}\limits^{#1}}\!\!}
\def\mapdown#1{\Big\downarrow\rlap{$\vcenter{\hbox{$\scriptstyle#1$}}$}}
\matrix
\calc(A, X) & \mapright{ (\ev, z_1, z_2) }& X \times \Delta^2
\cr
\mapdown{\pi_\calc} & &
\hphantom{\hbox{$X\times{.}$}} \mapdown{\lambda= z_1 \cdot z_2}
\cr
\calm(A,X)& \mapright{ \lambda_\calm } &
\hphantom{\hbox{$X\times{.}$}} \Delta
\endmatrix
\eqno(9.4.1)
$$
is commutative; moreover, for every $y\in \calm(A,X)$ the restriction
$\ev|_{C_y}$ is in 

\noindent $\calh^{1,p} (A_a, X)$ with  $a=\lambda_\calm(y)$,
and mappings
$\ev_1: y\in \calm(A,X) \mapsto \ev|_{C_y}(z_1^{-1}(1))$ and
$\ev_2: y\in \calm(A,X) \mapsto \ev|_{C_y}(z_2^{-1}(1))$ are
holomorphic;

\sliii let $C$ be an annulus or a nodal curve with smooth boundary
$\d C =\gamma_1\coprod \gamma_2$, $p_i \in \gamma_i$ be marked points,
and  $u:C \to X$ be a holomorphic $L^{1,p}$-smooth mapping, then there exist
{\sl a unique} $y\in \calm(A,X)$ and a {\sl unique} biholomorphism
$\phi: C \to C_y$ such that $\ev \scirc \phi = u : C\to X$ ¨ $z_i \scirc
\phi (p_i) =1\in \barr \Delta$; in other words, $\calm(A,X)$ parameterizes
holomorphic mappings to $X$ from annulai and nodes with marked points
on the boundary;

\sliv if the commutative diagram with complex spaces $\calw$ and
$\calz$ and holomorphic mappings
$$

\def\mapright#1{\!\!\smash{\mathop{\hbox to50pt{\cleaders\lowminus\hfill}
\mkern-13mu\rightarrow}\limits^{#1}}\!\!}
\def\mapdown#1{\Big\downarrow\rlap{$\vcenter{\hbox{$\scriptstyle#1$}}$}}
\matrix
\calz & \mapright{ (\ev^\calz, \ti z_1,  \ti z_2) }& X \times \Delta^2
\cr
\mapdown{\pi_\calz} & &
\hphantom{\hbox{$X\times{.}$}} \mapdown{\lambda= \ti z_1  \cdot \ti z_2}
\cr
\calw & \mapright{ \lambda_\calw } &
\hphantom{\hbox{$X\times{.}$}} \Delta
\endmatrix
\eqno(9.4.2)
$$
possesses the properties \sli and \slii, in particular, the fibers
$\calz_w \deff\pi_\calz\inv(w)$ should be the nodal curves
with induced maps $f_w\deff\ev^\calz|_{\calz_w} \in \calh^{1,p}(\calz_w, X)$,
then diagrams $(9.4.1)$ and $(9.4.2)$ uniquely fit together to form
the  commutative diagram
$$

\def\mapright#1{\!\!\smash{\mathop{\hbox to50pt{\cleaders\lowminus\hfill}
\mkern-13mu\rightarrow}\limits^{#1}}\!\!}
\def\mapdown#1{\Big\downarrow\rlap{$\vcenter{\hbox{$\scriptstyle#1$}}$}}
\matrix
\calz & \buildrel \wt F \over \lrar &
\calc(A, X) & \mapright{ (\ev, z_1, z_2) }& X \times \Delta^2
\cr
\mapdown{\pi_\calz} & &
\mapdown{\pi_\calc} & &
\hphantom{\hbox{$X\times{.}$}} \mapdown{\lambda= z_1 \cdot z_2}
\cr
\calw & \buildrel  F \over \lrar &
\calm(A,X)& \mapright{ \lambda_\calm } &
\hphantom{\hbox{$X\times{.}$}} \Delta
\endmatrix
\eqno(9.4.3)
$$
here $\lambda_\calm \scirc F = \lambda_\calw$ and $(\ev, z_1, z_2) \scirc
\wt F= (\ev^\calz, \ti z_1, \ti z_2)$;

\slv differential  $d\lambda_\calm: T_y\calm(A,X) \to T_{\lambda_\calm(y)}
\Delta \cong \cc$ is not degenerate in any point $y \in \calm(A,X)$,
and for any $a \in \Delta$ the fiber $\lambda_\calm\inv(a)$ is
naturally isomorphic to the manifold $\calh^{1,p} (A_a, X)$, where $A_a$ 
is a curve $\{ (z_1, z_2)\in
\Delta^2 : z_1 \cdot z_2 = a \}$; in particular, for any $y \in \calm(A,X)$
there is a biholomorphism $C_y \cong A_{\lambda_\calm(y)}$ and a natural
exact sequence
$$
0 \lrar
\calh^{1,p} (C_y, u^*TX) \buildrel \iota_y \over\lrar T_y \calm(A,X)
\buildrel d\lambda_\calm(y) \over\lrar \cc \lrar 0.
$$
\rm

\state Proof. Let $(A, p_1, p_2)$ be a smooth annulus with marked points,
one on each boundary component $\gamma_i \cong S^1$, and let $J$ be a complex
structure on $A$. It is known that $(A, J)$ is biholomorphic to some
annulus $A_{r,1}= \{ r < |z| <1 \}$. It is easy to see that there is
only one isomorphism $\psi: (A, J) \to A_{r,1}$ which smoothly extends
to the diffeomorphism $\psi: \barr A \to \barr A_{r,1}$ with $\phi(p_1) =1$.
Put $a \deff \phi(p_2)$. Now it is obvious that there is a unique
biholomorphism $\phi: (A, J) \to A_a \deff \{ (z_1, z_2)\in
\Delta^2 \;:\; z_1 \cdot z_2 = a \}$ such that $\phi(p_1) =1$ ¨ $\phi(p_2)
=a$.

Thus, the mapping  $\lambda: \Delta^2 \to \Delta$, $\lambda(z_1, z_2)
= z_1 \cdot z_2$ with a fiber  $A_a$ over $a\in \Delta$ represents a
holomorphic moduli space of annulai with marked points on the boundaries,
completed at $a=0$ by a standard node. If $a\not =0$ the coordinates
$z_i$, $i=1, 2$, realize an imbedding of each  $A_a$ as an annulus into
$\cc$ in such a way that the circle $\gamma_i$ becomes an outer unit
circle. When $a \lrar 0$, the annulus $A_a$ degenerates into a standard
node, and each $z_i$ becomes a standard coordinate on the corresponding
component of the node.

\state Remark. In the sequel we denote by $A_a$ an annulus (or corr.\ a node)
{\sl with fixed points $a$ and $1$ on the boundary $\d A_a$}, and also
together with uniquely determined coordinates $z_1$ and $z_2$.

\smallskip
Fix $r$ with $0 <r< 1$. For $|a|<r$ define a mappings $\zeta^a_1,
\zeta^a_2: A_{r,1} \to A_a$ as $\zeta^a_1(z) \deff z$ and $\zeta^a_2(z)
\deff a/z$, so that $\zeta^a_i$ are the inverse of the coordinates $z_i$.
Consider a mapping
$$
\matrix
\Psi_r: \coprod_{|a|<r} \calh^{1,p}(A_a, X) & \lrar &
\calh^{1,p}(A_{r,1}, X) \times \calh^{1,p}(A_{r,1}, X) \times \Delta(r)
\cr
u \in \calh^{1,p}(A_a, X) & \mapsto & (u \scirc \zeta^a_1, u \scirc
\zeta^a_2, a).
\endmatrix
$$
One easily sees that $\Psi_r$ is holomorphic on each  $\calh^{1,p}(A_a, X)$
and the image of $\Psi_r$ consists of such  triples $(u_1,u_2,a)$
that each $u_i \in \calh^{1,p}(A_{r,1}, X)$ extends to $u_i \in \calh^{1,p}
(A_{|a|,1}, X)$ and $u_2(z) = u_1(a/z)$. Thus $\Psi_r$ is injective and
has closed image. On the disjoint union $\calm(A, X) \deff \coprod_{a\in
\Delta}\calh^{1,p} (A_a, X)$ define the topology induced by the mappings
$\Psi_r$.
Clearly it agrees with the topologies on each slice $\calh^{1,p}(A_a, X)$.

Our task now is to construct an appropriate holomorphic structure
on $\calm(A, X)$, which agrees with the above introduced holomorhpic
structures on each slice

\noindent $\calh^{1,p}(A_a, X)$, and also with the topology
on $\calm(A, X)$, defined above.

\smallskip\noindent\sl
Case 1. Let us start with the special case $X= \cc^n$. \rm For
$a\not =0$  every function $f\in \calh^{1,p} (A_a, \cc^n)$ in a unique way
can be decomposed into a Laurent series $f(z_1) =\sum_{i= -\infty} ^\infty c_i
z_1^i$. Set $f^+(z_1) \deff \sum_{i=0}^\infty c_i z_1^i$ and $f^-(z_1)
\deff \sum_{i=-\infty} ^0 c_i z_1^i$. It is convenient to consider the
function  $f^-$ as a function of the variable $z_2 = a/z_1$; then
$f^-(z_2)= \sum_{i=0}^\infty c_{-i} (z_2/a)^i$. One has
$f^+ \in \calh^{1,p}(\{ |z_1| <1\}, \cc^n)$, $f^-\in \calh^{1,p}(\{ |z_2| <
1\}, \cc^n)$, $f^+(0) = f^-(z_2 =0) = c_0$, and $f(z_1) = f^+(z_1) +
f^-(a/z_1) - c_0$, so that a pair $(f^+, f^-)$ defines the holomorhpic
function $\hat f \in \calh^{1,p}(A_0, \cc^n)$. The canonical isomorphisms
$\calh^{1,p}(A_a, \cc^n)\cong \calh^{1,p}(A_0, \cc^n)$ obtained in this was
define
on $\coprod_{|a|<1} \calh^{1,p} (A_a, \cc^n)$ the structure of the trivial
Banach bundle with base $\{ |a| <1\}$ and fiber $\calh^{1,p}(A_0, \cc^n)$,
and thus a structure of a Banach manifold.

One can see that the mapping $\Psi_r : \coprod_{|a|<r} \calh^{1,p}(A_a,
\cc^n) \to
\calh^{1,p}(A_{r,1}, \cc^n)^2 \times \Delta(r)$ is holomorphic.

\smallskip\noindent\sl
Case 2. Let $X = U \subset \cc^n$ be an open subset in $\cc^n$.
\rm Then  $\calm(A, U)$ is also open
in $\calm(A, \cc^n)$ and consequently inherits a holomorphic
structure.
The natural projection $\lambda_\calm : \calm(A, U) \to \Delta$
becomes holomorphic. One easily sees that the differential $d\lambda_\calm$
is not degenerate. Thus, we can define {\sl a universal family of curves}
$\calc(A, U)$ as a fiber product $\calm(A, U) {\times}_\Delta
\Delta^2$ with respect to the mappings $\lambda_\calm : \calm(A, U) \to
\Delta$ and $\lambda: \Delta^2 \to \Delta$, which is a complex Banach
manifold due to the nondegeneracy of $d\lambda_\calm$.

Denote by $\pi_\calc : \calc(A, U) \to \calm(A, U)$ the natural projection.
Then the fiber $C_y$ over $y \in \calm(A, U)$ is biholomorphic to $A_a$ with
$a=\lambda_\calm(y)$. The natural projection from $\calc(A, U)$ onto
$\Delta^2$ induces holomorphic functions $z_1$ and $z_2$ on $\calc(A, U)$
taking values in $\Delta$  such that the property  \sli of
{\sl Theorem 9.4.1} is satisfied.

Now let $a\not=0 \in \Delta$ and $f\in \calh^{1,p} (A_a, \cc^n)$.
Represent $f$ in the form $f(z_1) = f^+(z_1) + f^-(a/z_1) - f_0$, where
$f^\pm \in \calh^{1,p} (\Delta, \cc^n)$ with $f^+(0) = f^-(0) = f_0$.
Analogously, for $f\in \calh^{1,p} (A_0, \cc^n)$ we have $f= (f^+, f^-)$,
where again $f^\pm \in \calh^{1,p} (\Delta, \cc^n)$ á $f^+(0) = f^-(0) = f_0$.
Consider a holomorphic function $\ti f(z_1, z_2) \deff f^+(z_1) + f^-(z_2)-
f_0$, $\ti f\in \calh^{1,p} (\Delta^2, \cc^n)$. Define mappings
$\ti\ev_a:\calh^{1,p} (A_a, \cc^n) \times \Delta^2 \lrar \cc^n$ by the formula
$\ti\ev (f, z_1, z_2) \deff \ti f(z_1, z_2)$. One easily sees that
$\ti\ev_a$ define a holomorphic map $\ti\ev :\calm(A, \cc^n) \times \Delta^2
\to \cc^n$. Let $\ev$ denote the restriction of this map onto
$\calc(A, \cc^n) \subset \calm(A, \cc^n) \times \Delta^2$. It is
straightforward to check that $\calm(A, \cc^n)$, $\calc(A, \cc^n)$,
$\ev: \calc(A, \cc^n) \to \cc^n$ and $z_{1,2} : \calc(A, \cc^n) \to \Delta$
satisfy the statements of {\sl Theorem 9.4.1}.

Thus, for $U \subset \cc^n$ we can define $\ev: \calc(A, U) \to U$ as a
restriction $\ev: \calc(A, \cc^n) \to \cc^n$.
Statements of  {\sl Theorem 9.4.1} again remain true. In particular, if
$G:  U \to U'\subset \cc^n$ is biholomorphic, then the natural bijections
$\calm(A,U) \buildrel \cong \over\lrar \calm(A, U')$ and $
\calc(A, U) \buildrel\cong \over \lrar \calc(A, U')$ are also 
biholomorphisms. This means that the complex structure in  $\calm(A, U)$
does not depend on the imbedding $U \subset \cc^n$.

\smallskip
Let $C= A_{r,1}$ be an annulus and $u: C \to X$ be a holomorphic imbedding.
Then $du : TC \to u^*TX$ is an imbedding of holomorphic bundles over $C$,
and this defines a holomorphic normal bundle as a factor-bundle
$N_C \deff u^*TX/ TC$. Since $C$ is Stein, the bundle $N_C$ is holomorphically
trivial. Fix a holomorphic frame $\sigma_1, \ldots, \sigma_{n-1} \in
\calh^{1,p} (C, N_C)$, $n= \dimc X$, and also its lift $\ti\sigma_1, \ldots,
\ti \sigma_{n-1} \in \calh^{1,p} (C, u^*TX)$. Denote by $B^{n-1}(r)$ the ball
of radius $r$ in $\cc^{n-1}$ with coordinates $w= (w_1, \ldots, w_{n-1})$.
As follows from {\sl Lemma 9.3.1}, there is a holomorphic map
$\Psi : C\times B^{n-1}(r) \to X$ such that $\d \Psi / \d w_i|_{z\in C,w=0}
=\sigma_i(z)$. Thus $\Psi$ is biholomorphic in the neighborhood of
$C\equiv C\times \{0\}$. In particular, for $r$ sufficiently small an
image $U\deff \Psi(C \times B^{n-1}(r))$ is a coordinate chart with
coordinates $(z, w_1, \ldots, w_{n-1})$.

\smallskip\noindent\sl
Step 3. Consider a general case where $C\cong A_a$ is arbitrary, and
$u: C \to X$ is a holomorphic map.

\smallskip\rm If $a=0$, then $C$ is a node; thus, there exists a 
neigborhood $V_0$
of the nodal point such that its image $u(\barr V_0)$ lies in some
coordinate chart of $X$. If $a\not=0$ and the image $u(\barr C)$ is not
contained in any chart of $X$, then $u$ is not constant. Thus, for some
$|a|< r <1$ mapping $u$ will be an imbedding in some neighborhood
of the circle $S^1_r \deff \{|z_1| =r\} \subset C \cong \{|a|< |z_1|<1\}$.

In all cases there is a covering $\{V_0, V_1, V_2\}$ of  $C=\{(z_1,
z_2) \in \Delta^2 : z_1 \cdot z_2 =a\}$ of the form $V_1 =\{ (z_1,z_2)\in C
: r_1 <|z_1| <1 \}$, $V_2 =\{ (z_1,z_2)\in C : r_2 <|z_2| <1 \}$,
$V_0 =\{ (z_1, z_2) \in C : |z_1| < R_1, |z_2|< R_2 \}$, where $0<r_1 <R_1
<1$, $0< r_2 <R_2 <1$, $r_1 \cdot R_2 > |a| <r_2 \cdot R_1$, and $V_0$
is such that the image $u(\barr V_0)$ lies in some chart $U$ of $X$.

\bigskip
%Figure "Covering"
\vbox{\xsize.35\hsize\nolineskip\rm
\putm[.15][0.]{|z_2|}%
\vskip.06\xsize
\putm[.095][.07]{1}%
\putm[.07][.32]{R_2}%
\putm[.09][0.5]{r_2}%
\putm[.345][.32]{V_2}%
\putm[.33][.68]{V_0}%
\putm[.66][.64]{V_1}%
\putm[1.01][.82]{|z_1|}
\putm[.48][.9]{r_1}%
\putm[.66][.9]{R_1}%
\putm[.92][.9]{1}%
\putt[1.2][0]{\advance\hsize-1.4\xsize%
\bigskip
\centerline{Fig.~17.\ \ Covering $\{V_0, V_1, V_2\}$.}
\bigskip%\sevenrm
\parindent=0pt
In this figure curve $A_a$ is drawn in boldface as a piece of the parabola,
the elements of the covering $V_0$, $V_1$ ¨ $V_2$ by broken lines.    
}
\noindent
\epsfxsize=\xsize\epsfbox{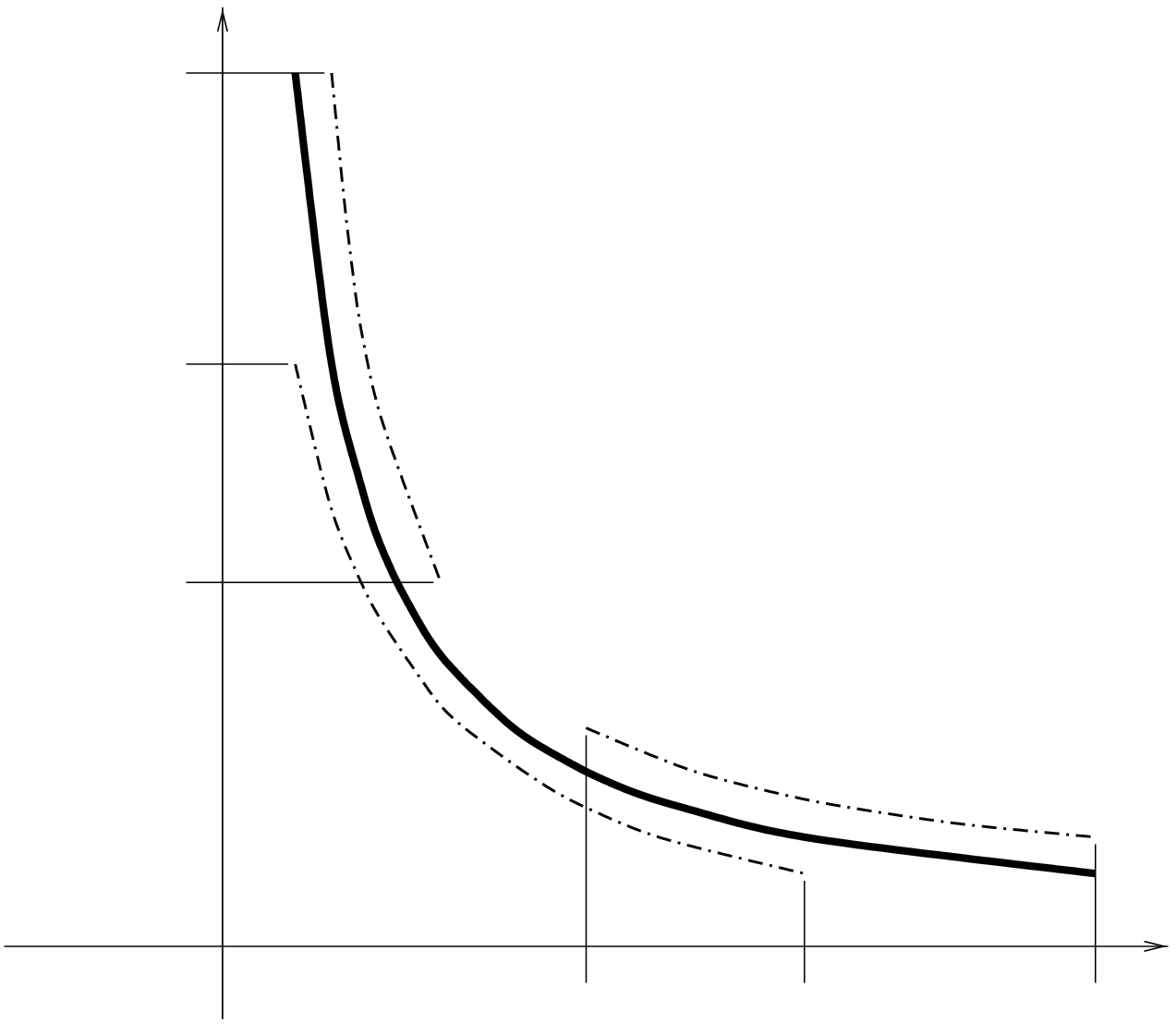}
}

\bigskip

Fix coordinate $z_1$ on $V_1$ and coordinate $z_2$ on $V_2$. On $V_0$ define
new coordinates $\ti z_1 \deff z_1 /R_1$ and $\ti z_2 \deff z_2 /R_2$,
and also the marked points $\ti p_1 \deff R_1$ and $\ti p_2 \deff a/R_2$.
Set  $\ti a \deff a/(R_1R_2)$. Then $V_0 \cong A_{\ti a}= \{(\ti z_1, \ti z_2)
\in \Delta^2 : \ti z_1\cdot \ti z_2 = \ti a\}$ with $\ti z_i(\ti p_1)=1$.
It is easy to see that the change of the complex parameter $a= z_1 z_2$ on
$C$, which parameterizes the holomorphic structures of annulai with marked
points on the boundary, can be reduced to the change of the analogous
parameter in $V_0$. Namely, let $C' \cong A_{a'} = \{ (z_1,z_2) \in \Delta^2:
z_1 \cdot z_2= a'\} $  be obtained as a result of a small deformation
of $a$. Put $\ti a{}'
\deff a'/(R_1R_2)$. Look on $C'$ as the result of patching of complex
curves $V_0'$,
$V_1$ ¨ $V_2$, where $V_0'$, $V_1$ and $V_2$ are defined in  $\Delta^2$
in the  same way $V_0'= \{ (z_1,z_2) \in C': |z_1| <r_1, |z_2| <r_2\}$,
$V_i= \{r_i <|z_i|<1\}$ when $i=1,2$, but now the coordinates  $(z_1,z_2)$
satisfy a new relation $z_1\cdot z_2 = a'$, and $V_0'\cong A_{\ti a{}'}$.

Thus, for a sufficiently small neighborhood $\calw \subset \calm(A, X)$ of
the curve $(C,u)$ over $X$ a holomorphic map is defined as
$$
\matrix \format \c\;&\c\;&\c\;& \c& \c& \c& \c& \c \\
\Theta :& \calw &\lrar& \calm(A,U) &\llap{$\times$}&
\calh^{1,p} (V_1, X) &\times&
\calh^{1,p} (V_2, X);
\cr
\noalign{\vskip2pt}
\Theta :&(C',u') &\mapsto& ((A_{\ti a{}'}, u|_{V_0'}), &&
u|_{V_1},&& u|_{V_2}).
\endmatrix
$$
On the other hand, a couple of maps $u_0 : V_0' \to X$, $u_1 : V_1 \to X$
and $u_2 : V_2 \to X$ define a map $u': C' \to X$, exactly where $u_1$
coincides with $u_0$ on $W_1 \deff V_1 \cap V_0'$, and $u_2$ coincides with
$u_0$ on $W_2 \deff V_2 \cap V_0'$. Note that domains $W_i \subset V_i$
do not change in the process of deformation of the complex structure on $C$.
A construction of gluing from {\sl Step 2} of the proof of {\sl Lemma 9.3.1}
finishes the proof of the theorem.
\qed

\medskip\noindent\sl
9.5. Banach Analytic Structure on the Stable Neighborhood.

\smallskip\rm
We are prepared now to prove the main result of this lecture.
\smallskip
\state Theorem 9.5.1. {\it Let  $(X,J)$ be a complex manifold, and
$(C_0,u_0)$ a stable complex curve over $X$, parameterized by a real
surface $\Sigma$. There exists a b.a.s.f.c. $\calm $ and $\calc$
and holomorphic maps $F: \calc \lrar X$ and  $\pi: \calc \lrar \calm$
%$$
%\matrix \calc & \buildrel F \over \lrar & X \cr \llap{$\pi$}\downarrow \cr
%\calm \endmatrix
%$$
such that

a) for any $\lambda \in \calm $ a fiber  $C_\lambda =\pi^{-1}(\lambda )$
is a nodal curve parameterized by $\Sigma$, and $C_{\lambda_0} \cong C_0$
for some $\lambda_0$;

b) for $F_\lambda \deff F|_{C_\lambda}$ a pair $(C_\lambda ,F_\lambda )$
is a stable curve over $X$; moreover $F_{\lambda_0}=u_0$;

c) if $(C',u')$ is a stable curve over $X$, sufficiently close to
$(C_0,u_0)$ in the Gromov topology, then there exists $\lambda' \in \calm$
such that $(C',u')=(C_{\lambda '},F_{\lambda'} )$;

d) for an appropriate $N{\in}\nn$ and for a small ball $B$ in
Banach space

\noindent $\calh^{1,p}(C_0, u_0^*TX) \allowbreak \oplus \cc^N$
b.a.s.f.c. \ $\calm$ is realized as a zero set of a holomorphic map
$\Phi$ from $B$ to the finite dimensional space $\sfh^1(C,u_0^*TX)$.  \rm
}

\smallskip

\noindent
\bf Proof. \rm The proof is based on the construction of local deformations of stable curves
and on the analysis of the gluing of local models.
 
Let $(C_0, u_0)$ be a stable curve
over a complex manifold $X$ with parameterization $\sigma_0: \Sigma
\to C_0$. Using  {\sl Proposition 3.2}, fix a covering $\{ V_\alpha,
V_{\alpha\beta}\}$ of the surface $\Sigma$, having properties
{\sl i)}--{\sl v$\!$i)} of the above-mentioned proposition. In particular, 
there exist
biholomorphisms $\phi^0 _{\alpha \beta}: \sigma_0 (V_{\alpha \beta}) \lrar
A_{\lambda^0 _{\alpha \beta}}$. From {\sl i)}--\slvi it follows that for every
couple $\bflamb \deff (\lambda _{\alpha \beta})$, sufficiently close to
$\bflamb^0 \deff (\lambda^0 _{\alpha \beta})$, there is a nodal curve
$C$ and parameterization $\sigma: \Sigma \to C$ such that properties
{\sl i}--{\sl v$\!$i} remain true and there exist biholomorphisms
$\phi _{\alpha \beta}: \sigma (V_{\alpha \beta}) \lrar A_{\lambda _{\alpha
\beta}}$. In particular, the complex structures on each $\sigma (W_{\alpha
\beta})$ do not change. Moreover, in $V_{ \alpha\beta}$ there exist
holomorphic coordinates $z_1$ and $z_2$ such that
$z_1 |_{W_{\alpha,\beta}}$ and $z_2|_{W_{\beta, \alpha}}$ do not change
by changing  $\lambda _{\alpha \beta}$, and $z_1 \cdot z_2 \equiv
\lambda _{\alpha \beta}$. This means that a disk $\Delta_{\alpha
\beta}\deff \{ \lambda _{\alpha \beta} \,:\, |\lambda^0 _{\alpha \beta}-
\lambda _{\alpha \beta} | \le \epsi \,\} $ parameterizes {\sl a holomorphic}
family of curves of the type $\sigma(V_{\alpha\beta })$.

One can illustrate the changing of the complex structure in $V_{\alpha\beta}$
by the following picture.

\medskip
\vbox{\nolineskip\xsize.548\hsize%
\putm[.01][.38]{\underbrace{\hskip.35\xsize}_{V_\alpha}}
\putm[.54][.44]{\underbrace{\hskip.46\xsize}_{V_\beta}}
\putm[.27][.16]{\overbrace{\hskip.36\xsize}^{V_{\alpha\beta}}}
\putm[.27][.3]{\underbrace{\hskip.09\xsize}_{W_{\alpha\beta}}}
\putm[.54][.3]{\underbrace{\hskip.09\xsize}_{W_{\beta\alpha}}}
\putt[1.05][0]{\advance\hsize-1.05\xsize\parindent=0pt%
\centerline{Fig.~18. }
\smallskip
In this picture one sees $V_{\alpha\beta}$
\noindent toge\-ther with adjacent sets
$V_\alpha$ and $V_\beta$. One can suppose that by varying 
$\lambda _{\alpha \beta}$, the complex structure varies only in the shaded
domain.
}
\noindent
\epsfxsize=\xsize\epsfbox{twist.eps}
}

\bigskip
Let $N \in \nn$ be the number of the elements of our covering of the type
$V_{\alpha \beta}=V_{\beta \alpha}$. Then for sufficiently small $\epsi>0$
the polydisk
$$
\Delta^N_\lambda \deff \{ \bflamb \deff (\lambda _{\alpha \beta}) \;:\;
|\lambda^0 _{\alpha \beta}- \lambda _{\alpha \beta} | \le \epsi \,\}
$$
parameterizes a holomorphic family of nodal curves $\{C_\bflamb\}_{\bflamb
\in \Delta^N _\lambda }$. Every curve $C_\bflamb = \sigma(\Sigma)$
is obtained by gluing the pieces $\sigma(V_\alpha)$ and
$\sigma(V_{\alpha \beta})$.
Denote $\calm \deff $ 

\noindent $\cup_{\bflamb \in\Delta^N _\lambda} \calh^{1,p}
(C_\bflamb, X)$. Since $\sigma(V_\alpha)$ does not contain nodal points
and the complex structures on $\sigma(V_\alpha)$ are constant and do not
depend on
$\bflamb =(\lambda _{\alpha \beta})$, we can suppose that  $V_\alpha$
are complex curves. We have natural isomorphisms $\sigma(V_
{\alpha \beta}) \cong A_{\lambda _{\alpha \beta}}$. Thus, the following map
is well- defined
$$
\matrix \format \c\;& \c \;& \c& \c& \c&\c \\
\Theta: & \calm &\lrar&
\prod_\alpha \calh^{1,p} (V_\alpha, X)&\times&
\prod_{\alpha \beta} \calm(V_{\alpha \beta}, X)
\cr
\noalign{\smallskip}
\Theta: & (C_\bflamb, u) & \mapsto & \bigl(u|_{V_\alpha}, 
(\lambda _{\alpha \beta}, u|_{V _{\alpha \beta}}) \bigr).
\endmatrix
$$
It is easy to see that a couple $\bigl( u_\alpha, (\lambda_{\alpha \beta},
u_{\alpha
\beta}) \bigr)\in \prod_\alpha \calh^{1,p} (V_\alpha, X)\times \prod_{\alpha
\beta} \calm(V_{\alpha \beta}, X)$ belongs to the image of $\Theta$,
exactly if for all pairs $(\alpha,\beta)$ the gluing conditions
$u_\alpha
|_{W_{\alpha, \beta}} = u_{\alpha \beta}|_{W_{\alpha, \beta}} $ are fulfilled.

Let us repeat the gluing procedure from {\sl Step 2 of Lemma 3.6.}
For this take the balls
$$
\matrix \format \l\,& \l\,& \l\,& \l\,& \l \\
B_\alpha &\subset& \calh^{1,p} (V_\alpha, u_0^*TX)
&\cong& T_{u_0}\calh^{1,p} (V_\alpha,X)
\cr
B_{\alpha \beta} &\subset& \calh^{1,p} (V_{\alpha \beta},u_0^*TX) \oplus \cc
&\cong& T_{u_0}\calm(V_{\alpha \beta}, X)
\cr
B'_{\alpha, \beta} &\subset& \calh^{1,p} (W_{\alpha, \beta}, u_0^*TX)
&\cong& T_{u_0} \calh^{1,p}(W_{\alpha, \beta}, X)
\endmatrix
$$
such that there exist biholomorphisms
$$
\matrix \format \l\,& \l\,& \l\,& \l\,& \l\,& \l \\
\psi_\alpha: & B_\alpha &\buildrel \cong \over \lrar &\psi_\alpha(B_\alpha)&
\subset& \calh^{1,p}(V_\alpha, X)
\cr
\psi_{\alpha \beta}:& B_{\alpha \beta}&
\buildrel \cong \over \lrar& \psi_{\alpha \beta}(B_{\alpha \beta})&
\subset& \calm(V_{\alpha \beta}, X)
\cr
\psi'_{\alpha,\beta}:& B'_{\alpha,\beta}&
\buildrel \cong \over \lrar& \psi'_{\alpha,\beta} (B'_{\alpha,\beta})&
\subset &\calh^{1,p}(W_{\alpha,\beta}, X)
\endmatrix
$$
having the properties
$$
\matrix \format \l\,& \l\,& \l\, \hskip 7pt & \l\,& \l\,& \l \\
\psi_\alpha(0)  &=& u_0|_{V_\alpha}, &
d\psi_\alpha(0) &=&
\id : T_{u_0} \calh^{1,p}(V_\alpha, X) \to T_{u_0} \calh^{1,p}(V_\alpha, X),
\cr
\psi_{\alpha\beta}(0)  &=& u_0|_{V_{\alpha\beta}}, &
d\psi_{\alpha\beta}(0) &=&
\id : T_{(\lambda^0_{\alpha,\beta}, u_0)} \calm(V_\alpha, X) \to
T_{(\lambda^0_{\alpha,\beta}, u_0)} \calm(V_\alpha, X),
\cr
\psi'_{\alpha, \beta}(0)  &=& u_0|_{W_{\alpha, \beta}}, &
d\psi'_{\alpha, \beta}(0) &=&
\id: T_{u_0} \calh^{1,p}(W_{\alpha, \beta}, X)
\to T_{u_0} \calh^{1,p}(W_{\alpha, \beta}, X).
\endmatrix
$$
Shrinking the balls $B_\alpha$ and  $B_{\alpha \beta}$, if nessessary, we
can suppose that for all $\xi_\alpha\in B_\alpha$ and 
$\xi_{\alpha\beta} \in B_{\alpha \beta}$ the restrictions
$\psi_\alpha(\xi_\alpha)|_{W_{\alpha,\beta}}$ and
$\psi_{\alpha \beta}(\xi_{\alpha \beta})|_{ W_{\alpha, \beta}}$
belong to the image $\psi'_{\alpha, \beta} (B'_{\alpha, \beta })$. Consider
a holomorphic map
$$
\matrix\format \c\,& \c& \c & \c\,& \c \,& \l\\
\Psi:& \prod_\alpha B_\alpha &\times&
\prod_{\alpha <\beta} B_{\alpha \beta}
& \lrar   & \sum_{\alpha, \beta} \calh^{1,p}(W_{\alpha, \beta}, u_0^*TX)
\cr
\Psi:& (v_\alpha, && v_{\alpha \beta}) & \mapsto &
\psi^{\prime\;\;-1}_{\alpha, \beta}
\bigl(\psi_\alpha(v_\alpha)|_{W_{\alpha, \beta}}\bigr)
- \psi^{\prime\;\;-1}_{\alpha, \beta}
\bigl(\psi_{\alpha \beta}(v_{\alpha \beta})|_{W_{\alpha, \beta}}\bigr).
\endmatrix
$$
As in analogous situations, which have already appeared in this paper, 
the map $\Psi$
gives the gluing condition of local holomorphic maps
$\psi_\alpha (v_\alpha): V_\alpha \to X$ and
$\psi_{\alpha \beta} (v_{\alpha \beta}): V_ {\alpha \beta}\to X$. Thus
we can identify $\calm \cap \prod_\alpha B_\alpha \times \prod_{\alpha
<\beta} B_{\alpha \beta}$ with the set $\Psi\inv(0)$.

Let us study in detail the behavior of $\Psi$ in the point $y_0 \in \prod
B_\alpha \times \prod B_{\alpha \beta}$, $y_0 = \bigl(\psi_\alpha \inv (u_0
|_{V_\alpha}), \psi_{\alpha \beta} \inv (u_0| _{V_{\alpha \beta}})\bigr)$,
$\Psi(y_0) =0 \in \sum_{\alpha, \beta} \calh^{1,p}(W_{\alpha, \beta},
u_0^*TX)$. One easily sees that the tangent space at $y_0$ is
$$
T_{y_0}\bigl( \prod B_\alpha \times \prod B_{\alpha \beta}\bigr)
=
\sum_\alpha \calh^{1,p} (V_\alpha, u_0^*TX)  \oplus
\sum_{\alpha \beta} \calh^{1,p} (V_{\alpha \beta}, u_0^*TX)
\oplus \cc^N,
$$
and the differential $d\Psi(y_0)$ on the summand  $\sum\calh^{1,p} (V_\alpha,
u_0^*TX)\oplus \sum \calh^{1,p} (V_{\alpha \beta},$

\noindent $ u_0^*TX )$ coincides with
the \v{C}ech  codifferential $(9.1.1)$ with respect to the covering 
$\{V_\alpha, V_{\alpha\beta}\}$ of the curve $C_0$. 
By {\sl Lemma 9.2.1} we can represent
$\sum
\calh^{1,p} (W_{\alpha, \beta}, u_0^*TX)$ as a direct sum $\calw \oplus
\calq$, where $\calw= \im(d\Psi(y_0))$, and $\calq$ is isomorphic to
$\sfh^1(C_0, u_0^*TX)$ and of finite dimension. Let $\Psi_\calw$ and
$\Psi_\calq$ be the components of $\Psi= (\Psi_\calw, \Psi_\calq)$ with
respect to this decomposition, and let $\wt\calm \deff \Psi_\calw\inv(0)$.
By {\sl Lemma 9.2.1} and the implicit function theorem, $\wt\calm$ is a complex
submanifold in $\prod B_\alpha \times\prod B_{\alpha\beta}$ with a tangent
space at $y_0\in \wt\calm$   equal to $\calh^{1,p} (C_0, u_0^*TX) \oplus
\cc^N$, and $\calm$ is defined in $\wt\calm$ as a zero set of a holomorphic
mapping $\Phi\deff \Psi_\calq
|_{\wt\calm}: \wt\calm \to \calq \cong \sfh^1(C, u_0^* TX)$. This defines
on $\calm$ a structure of the Banach analytic set of finite definition.

\medskip
All that remains is to construct the corresponding holomorphic
family of nodal curves $\pi: \calc\to \calm$ and the holomorphic mapping
$F: \calc \to X$. Note that to each ball $B_\alpha$
we can naturally associate a trivial family  $\pi_\alpha: \calc_\alpha
\deff B_\alpha  \times V_\alpha \to B_\alpha$, and to each ball
$B_{\alpha\beta}$ a holomorphic family $\pi_{\alpha\beta}:
\calc_{\alpha\beta} \to B_{\alpha  \beta}$, with fiber $\pi_{\alpha\beta}
\inv( v_{\alpha\beta} )$ equal to $A_{\lambda_{\alpha\beta}}$, where
$\lambda_{ \alpha\beta}$ is uniquely determined by the relation
$\psi_{\alpha\beta} (v_{\alpha\beta}) = (\lambda_{\alpha\beta}, 
u_{\alpha\beta}) \in \calm( V_{\alpha\beta}, X)$.

Extend these families to the families 
$\ti\pi_\alpha: \wt\calc_\alpha \to\prod  B_\alpha \times\prod B_{\alpha\beta}$
and 
$\ti\pi_{\alpha\beta}: \wt\calc_{ \alpha \beta} \to\prod B_\alpha 
\times\prod B_{\alpha\beta}$. It is obvious
that $\wt\calc_\alpha$ and $\wt\calc_{\alpha\beta}$ canonically patch into a
globally defined family of nodal curves $\ti\pi: \wt\calc \to
\prod B_\alpha \times\prod B_{\alpha\beta}$. $\calc$ will be a Banach
manifold. Moreover, we shall obtain correctly defined  holomorphic maps
$F_\alpha: \wt\calc_\alpha \to X$ and $F_{\alpha\beta}: \wt\calc_{
\alpha \beta} \to X$ such that for $z\in V_\alpha$ one has $F(v_\alpha,
v_{\alpha\beta}, z) \deff \psi_\alpha(v_\alpha)[z]$ plus analogous relations
for $F_{\alpha\beta}$.

Let us define $\calc$ as a restriction $\calc \deff \wt\calc|_\calm$. Note
that the restriction of the {\sl trivial} holomorphic family
$\wt\calc_\alpha= V_\alpha \times \prod B_\alpha \times \prod
B_{\alpha\beta}$  ­  $\calm$ is also a trivial holomorphic family. 

%%%%%%%%%%%%%%%%%%%%%question%%%%%%%%%%%%%%%%%%%%%%

Thus,
$\calc$ is a b.a.s.f.d. in the neighborhood of the points $y \in \calc \cap
\wt\calc_\alpha$. Analogously, every holomorphic family of curves
$\wt\calc_{\alpha\beta}$ is defined in the trivial bundle
$\Delta^2 \times \prod B_\alpha \times \prod B_{\alpha \beta} \buildrel\pr
\over\lrar \prod B_\alpha \times\prod B_{\alpha\beta}$
by the condition $z_1\cdot z_2- \lambda_{\alpha\beta} =0$, where
$\lambda_{\alpha\beta}
:B_{\alpha \beta} \to \Delta$ is a holomorphic parameter of the deformation
of the complex structure in $V_{\alpha \beta}$, and $z_1$ with $z_2$ are the
standard coordinates in $\Delta^2$. So $\calc$ is a b.a.s.f.d. also
in the neighborhood of points  $y \in \calc \cap \wt\calc_{\alpha \beta}$.
Since $\calm$ was defined by the condition of the coincidence of local
mappings $F_\alpha$ and $F_{\alpha\beta}$; thus, on $\calc$ a global
holomorphic map $F: \calc \to X$ is defined.

Properties {\it a), b) {\rm and} d) of \sl Theorem 9.5.1} for $\calm$, $\calc$
and $F$ follow directly from the construction, and one obtains 
property {\it c)} by using {\sl Theorem 5.3.2}.
\qed

\medskip\noindent\sl
9.6. Drawing Families of Curves.

\smallskip\rm
In the proof of the Continuity Principle for meromorphic mappings,
which we need in this paper (see\ {\sl Theorem 4.2} of the next paragraph
or {\sl Theorem 5.1.3} from [I-S]), we use the following consequence of the
main results of the present paragraph.

Let $(C_n, u_n)$ be a sequence of irreducible stable curves over $X$,
converging to a stable curve $(C_\infty, u_\infty)$.

\state Lemma 9.6.1.
\it There is an index $N$ and a complex (maybe singular) surface
$Z$ with holomorphic mappings $\pi_Z :Z \to \Delta$ and $F: Z \to
X$, which define together a holomorphic family of stable nodal curves over
$X$, joining $(C_N, u_N)$ with $(C_\infty, u_\infty)$. More precisely, the
following are true:

\smallskip\noindent
$1)$ For every $\lambda\in \Delta$ the fiber $C_\lambda=\pi_Z^{-1}(\lambda)$
is a connected nodal curve with boundary $\d C_\lambda$; a pair
$(C_\lambda, u_\lambda)$ with $u_\lambda \deff F|_{C_\lambda}$ is a stable
curve over $X$.

\noindent
$2)$ All curves $C_\lambda$, except for finite numbers, are connected and
smooth.

\noindent
$3)$ $(C_0,u_0)$ is equal to $(C_\infty, u_\infty)$, and there exists
$\lambda_N\in \Delta$ such that $(C_{\lambda_N},u_{\lambda_N})=(C_N,u_N)$.

\noindent
$4)$ There exist open sets $V_1,\ldots,V_m$ in $Z$ such that each
$V_j$ is biholomorphic to $\Delta \times A_j$, where $A_j$ is an
annulus. Moreover, the  diagram
$$
\matrix
V_j&\buildrel \cong \over \lrar &\Delta\rlap{$\times A_j$}
\cr
\llap{$\pi$}\downarrow & &\downarrow\rlap{$\pi_{\Delta }$}
\cr
\Delta & = & \Delta
\endmatrix
$$
is commutative; each annulus $C_\lambda \cap V_j \cong \{\lambda\}
\times A_j$ is adjacent exactly to one of the components of the
boundary $\d C_\lambda$,
and the number $m$ of domains $V_j$ is equal to the number of the boundary
components  of each curve $C_\lambda$.

%$5)$ on each $V_j$ the map  $F: Z \to X$ lifts to the map
%$\hat F: Z \to \hat U$, \ie there exists a holomorphic maps
%$\hat F_j: V_j\to \hat U$, such that $ \hat\pi \scirc \hat F_j = F|_{V_j}$;
%more over, on $V_j \cap C_{\lambda_N}$ this lifts coincide with $\hat u_N :
%C_N\cong C_{\lambda_N} \to \hat U$.
\rm

\state Remark. This lemma was stated without proof in {\sl Proposition 5.1.1}
from [I-S].

\state Proof. Let $\pi_\calc: \calc \to \calm$ and  $\ev: \calc \to
X$ define a complete family of holomorphic deformations of the stable nodal
curve over $X$  $(C_\infty, u_\infty)$, constructed in {\sl Theorem 9.4.1}.
Let $\lambda^* \in \calm$ parameterize the curve $(C_\infty, u_\infty)$.
Choose a sequence $\lambda_n \lrar \lambda^*$ ¢ $\calm$ such that
$(C_n, u_n) \cong (\pi_\calc\inv(\lambda_n),
\ev|_{\pi_\calc\inv(\lambda_n)})$ for all sufficiently big $n$.

While $(C_\infty, u_\infty) = (\pi_\calc\inv (\lambda^*), 
\ev|_{ \pi_\calc \inv(\lambda^*)})$ lifts in the neighborhood
of the boundary of $\hat U$, thus by shrinking $\calm$, if nessessary, we can
suppose that this is true for all $\lambda \in \calm$.

By Theorem 9.5.1 the space $\calm$ is a b.a.s.f.d. therefore, we can
apply {\sl Theorem 9.1.1}. In particular, there are only finitely many
irreducible components at $\lambda^*$ of $\calm$. Let $\calm_1$ be a
component of $\calm$, which contains infinitely many $\lambda_n$. Represent
$\calm_1$ in the neighborhood of $\lambda^*$ as a proper ramified covering
$\pi_1: \calm_1 \to B_1$ over a Banach ball $B_1$. If $\Delta$ is an
imbedded disk in $B_1$, then $\pi_1\inv(\Delta)$ is a one-dimensional
analytic set, irreducible components which can be also parameterized by the
disks. Thus, there exists a holomorphic mapping $\phi: \Delta \to \calm_1$,
passing through $\lambda^*$ and $\lambda_N$ for some $N>\!>1$. The pre-image
of the family $\pi_\calc: \calc\to \calm$ with respect to $\phi$ defines
a holomorphic family of stable nodal curves over $\Delta•$ with total space
$\pi_Z: Z \to\Delta$ and mapping $F:Z \to X$, containing $(C_\infty,
u_\infty)$ and $(C_N, u_N)$.

%%%%%%%%%%%%%%%%%%%%%%%question%%%%%%%%%%%%%%%%%%%%%

While $C_N$ is smooth, the general fiber $C_\lambda =\pi_Z\inv(\lambda)$, and
again is smooth. Taking if nessessary a smaller disk, we can suppose that
$C_\lambda$ are singular only for a finite number of $\lambda \in \Delta$.
All other properties  1)--5) follow from the construction of the family
$\pi_Z: Z\to \Delta$ and a mapping $F:Z \to X$.
\qed

%%%%%%%%%%%%%%%%%%%%%%%Page 151%%%%%%%%%%%%%%%%%%%%%%%%%%%%

\newpage

\noindent{\bigbf
Lecture 10}

\smallskip\noindent
{\bigbf Envelopes of Meromorphy of Two-Spheres.}

\medskip\noindent\sl
10.1. Continuity Principles Relative to K\"ahler Spaces.

\smallskip
\rm
Our aim in this paragraph is to prove {\sl Theorem 4.1} and give some
corollaries from it. First of all we
need an appropriate form of the so- called ``continuity principle'' for the
extension of meromorphic mappings.

For the notion of meromorphic mapping from a domain $U$ in a complex manifold
into a complex manifold (or space) $Y$ we refer to [Rm]. We only point
out here that meromorphic mappings into $Y = \cc\pp^1$ are exactly meromorphic
functions on $U$, see [Rm].

\smallskip\rm

\state Definition 10.1.1. \it A Hermitian complex space $Y$ is called 
disk-convex if,
for any sequence
$(C_n,u_n)$ of smooth  curves over $Y$ parameterized by the same surface 
$\Sigma $, such that

1) $\area(u_n(C_n))$ are uniformly bounded and

2) $u_n$ $C^1$-converges in the neighborhood of $\d C_n $,

\noindent \it there is a compact $K\subset Y$ which contains all $u_n(C_n)$.

\smallskip\rm
This definition obviously carries over to the case where $Y$ is a symplectic
manifold.
In this case one should consider $(C_n,u_n)$ as $J_n$-holomorphic curves,
with $J_n$
converging to some $J$ (everything in $C^1$-tology) and all structures
being tamed
by a given symplectic form.

Let $U$ now be a domain in the complex manifold $X$, and $Y$ is 
a complex space .

\state Definition 10.1.2. \it An envelope of meromorphy of $U$ relative to $Y$ is the largest
domain $(\hat U_Y,\pi )$ over $X$, which contains $U$ (i.e., there exists an
imbedding
$i:U\to \hat U_Y$ with $\pi \scirc i={\sl Id}$) such that every meromorphic
map
$f:U\to Y$ extends to a meromorphic map $\hat f:\hat U_Y\to Y$.

\smallskip\rm

\smallskip
Using the Cartan-Thullen construction for the germs of meromorphic mappings of
open subsets of $X$ into $Y$ instead of germs of holomorphic functions, one
can prove the existence and uniqueness of the envelope.

\smallskip\noindent
\bf Proposition 10.1.1. \it For any domain $U$ in the complex space 
$X$ and for any
complex space $Y$ there exists a maximal domain $(\hat U_Y,\pi )$ over $X$,
containing $U$ such that every meromorphic mapping $f:U\longrightarrow Y$
extends to a meromorphic mapping $\hat f:\hat U_Y\longrightarrow Y$. Such
a domain is unique up to a natural isomorphism.

\smallskip
\noindent \rm See [Iv-1] for details.

\smallskip

\state Theorem 10.1.2. \tensl (Continuity Principle-I). \it Let $X$ be a
disk-convex complex surface and $Y$ a disk-convex K\"ahler space. Then the
envelope of meromorophy $(\hat U_Y,\pi )$ of $U$
relative to $Y$ is also disk-convex with respect to the pulled-back K\"ahler 
form.

\smallskip\rm
This result can be reformulated in more familar terms as follows.

Let $\{ (C_t,u_t)\}_{t\in [0, 1]}$
be a continuous (in a Gromov topology) family of complex curves over $X$ with boundaries, parameterized by
a unit interval. More precisely, for each $t\in [0, 1[$ a smooth Riemann
surface with
boundary $(C_t, \d C_t)$ is given together with the holomorphic
mapping $u_t: C_t\longrightarrow X$, which is $C^1$-smooth up to the boundary.
Note that $C_1$ is not supposed to be smooth, i.e., it can be a nodal curve!

Suppose that
in the neighborhood $V$ of $u_0(C_0)$ a meromorphic map
$f$ into the complex space $Y$ is given.

\smallskip
\noindent
\bf Definition 10.1.3. \rm We shall say that the map $f$ meromorphically
extends along the family $(C_t,u_t)$ if for every $t\in [0, 1]$ a
neighborhood
$V_t$ of $u_t(C_t)$ is given, and given a meromorphic map $f_t:V_t
\longrightarrow
Y$ such that

a) $V_0 = V$ and $f_0= f$;

b) if $V_{t_1}\cap V_{t_2}\not=\emptyset $ then $f_{t_1}\ogran_{V_{t_1}\cap
V_{t_2}} = f_{t_2}\ogran _{V_{t_1}\cap V_{t_2}}$.

\smallskip\noindent

\smallskip
\noindent

\smallskip
\noindent
\bf Theorem 10.1.3. (\sl Continuity Principle-II) \it Let $U$ be a domain in
the complex surface $X$.
Let $\{ (C_t,u_t)\} _{t\in [0, 1]} $  be a continuous family of complex
curves over $X$ with boundaries in $U_1$, a relatively compact 
subdomain in $U$.
Suppose moreover that $u_0(C_0)\subset U$ and that $C_t$ for $t\in [0,1[$
are smooth. Then every meromorphic mapping
$f$ from $U$ to the disk-convex K\"ahler space $Y$ extends
meromorphically along the family $(C_t,u_t)$.

\smallskip\noindent\rm

Taking an image manifold $Y$ complex line $\cc $ or $\cc \pp ^1$ we obtain the
continuity principles for holomorphic or meromorphic functions. When it is
necessary to emphasize that we are considering  the mappings into a certain
manifold $Y$, we shall refer to the statement above as the \sl continuity
principle relative to $Y$ \rm or \sl c.p.\ for the meromorphic mappings into
$Y$.

\rm

The discussion above leads to the following

\smallskip
\noindent
\bf Corollary 10.1.4.
\it If we have the domain $U$ in a complex surface $X$, a K\"ahler 
space $Y$ and
a family $\{ (C_t,u_t)\} $ satisfying the conditions of the ``continuity
principle'', then the family $\{ (C_t,u_t)\} $ can be lifted onto $\hat
U_Y$, {\sl i.e.,} there exists a continuous family $\{ (C_t,\hat u_t) \} $ of
complex curves over $\hat U$ such that $\pi \scirc \hat u_t = u_t$ for
each $t$.

\smallskip\rm Of course the point here is that the mapping can be extended
to the neighborhood of $u_1(C_1)$, which is a reducible curve having in
general compact components. This makes our situation considerably more
general than the classical one, i.e., when $X$ is supposed to be Stein;
compare with [Ch-St].

\smallskip
\state Remark. Let us explain the meaning of this theorem by an example.
Let $X$ be disk-convex and $U\subset X$ some domain. Furthermore, let
$f:U \to Y$ be a meromorphic map and $\{(C_n, u_n)\}$ a sequence of
stable complex curves over $X$, converging in the Gromov topology to
$(C_\infty, u_\infty)$. Suppose that the images of the boundaries
$u_n (\d C_n)$ and $u_\infty (\d C_\infty)$ are contained in  $U$.
Suppose also that $f$ extends along every curve $u_n (C_n)$. This means
that there exists a complex surface $V_n$, containing $C_n$, and a locally
biholomorphic map $u'_n:V_n \to X$ such that $u'_n|_{C_n} = u_n$ and $f$
meromorphically extend from $u'_n{}\inv U$ onto the whole $V_n$.

The latter is equivalent to the lift of the curves $(C_n, u_n)$ into the
envelope $\hat U$, i.e., to the existence of holomorphic mappings
$\hat u_n: C_n \to \hat U$ such that $\hat \pi \scirc \hat u_n =u_n$. In
other words, one can take as $V_n$ a neighborhood of the lift of  $C_n$
into $\hat U$. One easily sees that the curves  $(C_n,\hat u_n)$ are
stable over $\hat U$, have uniformly bounded areas and converge in the
neighborhood of the boundary $\d C_n$. By {\sl Theorem 10.1.2} and
by the Gromov compactness theorem , some subsequence
$(C_n,\hat u_n)$ converges to the $\hat{U}$-stable curve
$(C_\infty,\hat u_\infty)$, and $\hat\pi \scirc \hat u_\infty
=u_\infty$. This means that $f$ extends along the curve $u_\infty (
C_\infty)$. Thus,  {\sl Theorem 4.2} is a generalization of the
E. Levi continuity principle.

\medskip\noindent
{\sl 10.2. Proof of the Continuity Principle.}

\smallskip
Suppose that there is a subsequnce of our sequence, which we still denote by
$(C_n,u_n)$ and  which is not contained in any compact subset of the envelope
$(\hat U_Y, \pi )$. Put $v_n=\pi \scirc u_n$ and consider our sequence as 
a sequence $(C_n,v_n)$ of stable curves over $X$. While the areas are bounded,
we can suppose that this sequence converges in the Gromov topology by the 
disk-convexity of $X$. Denote by $(C_0,v_0)$ its limit. For $N$ sufficiently 
large take a holomorphic family $(\calc ,\pi ,\calv )$ joining $(C_N,v_N)$ 
with $(C_0,v_0)$ as in {\sl Lemma 9.6.1.}

\smallskip
Let $f:\hat U_Y \to Y$ be some meromorphic mapping of our envelope into a disk-convex
K\"ahler space $Y$. Composing $f$ with $\calv$, we obtain a 
meromorphic map $h:= f\scirc \calv :
\bigcup U_j\cup \pi^{-1}(V(\lambda_N))\to Y$. Here $V(\lambda_N)$ 
is a suffuciently small neighborhood of $\lambda_N \in \Delta $.

Denote by $W$ the maximal connected open subset of $\Delta $ containing
$V(\lambda_N)$
such that $h$ meromorphically
extends onto $\bigcup U_j\cup \pi^{-1}(W)$. We want to prove that $W=\Delta
$. Denote
by $\Gamma_{f_{\lambda }}$ the graph of the restriction $f\mid_{C_{\lambda }}$. Fix
some Hermitian mertic on $\calc $.

\state Lemma 10.2.1. \it For any compact $K\subset \Delta $ there is a constant
$M_K$ such that ${\sl area}(\Gamma_{f_{\lambda }})\le M_K$ for all $\lambda
\in W\cap K$.

\smallskip\noindent\sl
Proof. \rm Shrinking $\calc $ if nessessary, we can suppose that $f$ has 
only a finite number of points of indeterminancy in 
$\bigcup_jU_j$. Denote by $S_W$ the discrete
subset in $W$, which consists of points $s\in W$ such that either the fiber
$C_s$ is singular or  $s$ is the projection of the indeterminancy points of
$f\mid_{\calc \mid_W}$.
Fix a point $\lambda_1\in W\setminus S_W$. Take a path $\gamma :[0,1]\to
W\setminus
S_W$ connecting
$\lambda_1$ with some $\lambda_2\in W\setminus S_W$. Take some relatively compact
in $W\setminus S_W$ neighborhood $V$ of $\gamma ([0,1])$.

Recall that a K\"ahler metric on a complex space $Y$ consists of a locally 
finite
covering $\{ V_{\alpha }\} $ of $Y$ and strictly plurisubharmonic functions
$\phi_{\alpha }$ on $V_{\alpha }$ such that $\phi_{\alpha }-\phi_{\beta }$ are
pluriharmonic on $V_{\alpha }\cap V_{\beta }$. The sets 
$\{ f^{-1}(V_{\alpha })\} $ form a
covering of $\pi^{-1}(V)$ and $dd^cf^*\phi_{\alpha }=w$ is a correctly- defined
semi-positive closed form on $\pi^{-1}(V)$.
We have 
$$
{\sl area}(\Gamma_{f_{\lambda }})={\sl area}(C_{\lambda })+
\int_{C_{\lambda }}w.
$$
\noindent
By the Stokes formula
$$
\int_{C_{\gamma (0)}}w-
\int_{C_{\gamma (1)}}w =
\int_{\cup_tC_{\gamma (t)}}dw -
\int_{\cup_t{\d C_{\gamma (t)}}}w=
$$
$$
=
\int_{\cup_t{\d C_{\gamma (t)}}}w.
$$
\noindent
This is obviously uniformly bounded on $\lambda_2\in K\setminus S_W$, which
implies the uniform bound in $K$ by the Gromov (or Bishop, in this case) 
compactness theorem.

\smallskip
\qed

\state Lemma 10.2.2. \it Mapping $h$ meromorphically extends onto
$\calc $.

\smallskip\noindent\sl
Proof. \rm Let us prove first that $h$ meromorhpically extends onto $\calc $
minus the singular fiber $C_0$, i.e., that $W\supset \Delta \setminus \{ 0\} $.

The same arguments as in [Iv-2] (the only difference being that in [Iv-2] 
curves $C_
{\lambda }$ were disks) show that $(\Delta \cap \d W)\setminus (\{ 0\})
$ is an
analytic variety. Now, using a Thullen-type extension theorem of Siu, 
see [Si-1],
we can extend $h$ onto $\calc \setminus C_0$.

Denote by $\hat C_0$ the union of compact components of $C_0$. We shall show
that $\hat C_0$ contracts to a finite number of normal points. To prove this we
can suppose that $\hat C_0$ is connected. Otherwise, we can apply the same
arguments as below to the connected components. We shall now prove that
$\hat C_0$ contracts to a normal point. Denote by $E_1,...,E_n$ the
irreducible components of $\hat C_0$. All we must prove is that the
matrix $(E_i,E_j)$ is negatively defined, see [Gra]. Denote by $l_i$ the
multiplicity of $E_i$. Thus $\hat C_0=\Sigma_{i=1}^nl_i\cdot E_i$. Denote
by $M$ the $\zz $-module generated by $E_1,...,E_n$ with scalar product defined
by intersection of divisors. Put $D=\hat C_0$. Then we have

\noindent
1) $E_i\cdot E_j\ge 0$ for $i\not= j$;

\noindent
2) $D\cdot E_i\le 0$ for all $i$, because this number is not more than the
intersection of $E_i$ with the nonsingular fiber $C_{\lambda }$.

\smallskip
By Proposition 3 from [Sh], v.1, Appendix I, we have $A\cdot A\le 0$
for all $A\in M$ and $A\cdot A=0$ iff $A$ is proportional to $D$. But
$D\cdot D <\hat C_0\cdot C_{\lambda }=0$, where $C_{\lambda }$ is a smooth
fiber. This proves that $(E_i\cdot E_j)$ is negatively defined.

Therefore, all that is left to prove
is that the normal point is a removable singularity for the 
meromorphic mappings into a disk-convex K\"ahler space.

Let $(\calc ,s)$ be a germ of two-dimensional variety with an isolated normal
singularity $s$. Let a meromorphic map $h:\calc \setminus \{ s\} \to Y$ 
be given.
Realize $\calc $ as a finite proper analytic cover over a bi-disk 
$\Delta^2$ with
$s$ being the only point over zero. Denote by $p:\calc \to \Delta^2$ this
covering.

The composition $h\scirc p^{-1}$ is a multivalued meromorphic map 
from $\Delta^2 \setminus \{ 0\} $
to $Y$. This can be extended to the origin, see [Si-1]. Consider the following
analytic set in $\calc\times Y$:
$$
\Gamma = \{ (x,y)\in \calc \times Y: (p(x),y)\in \Gamma_{h\scirc p^{-1}} \} .
$$
The irreducible component of $\Gamma$, which contains $\Gamma_h$, will be the
graph of the extension of $h$ onto $\calc $.

\smallskip
\qed

Let us turn to the proof of the theorem. We have a holomorphic
family $\calc \to \Delta $ of  stable curves $(C_{\lambda },
v_{\lambda })$
over $X$ such that:

1) $(C_{\lambda_0},v_{\lambda_0})=(C_0,v_0)$;

2) $h=f\scirc \calv $ extends meromorphically onto $\calc $, where
$\calv :\calc \to X$ is an evaluation map.

\smallskip
We want to lift this family to the envelope $(\hat U_Y, \pi )$. Denote
by $E_1,...,E_l$ the set of all irreducible curves in $\calc $
which are contracted
by $\calv :\calc \to X$ to the points. It follows that $\calc \setminus
\bigcup_{j=1}^lE_j$ 
lifts to the envelope by the Cartan-Thullen construction. Note further that
$E_j$
do not intersect $\bigcup_{j=1}^mU_j$. That is why either $E_j\subset \hat C_0$
is a compact component of singular fiber or projects
surjectively
on $\Delta $ under the projection $\pi :\calc \to \Delta $. In
the second case it
intersects $C_0$ and the point $u(C_0 \cap E_j)$ lifts to the
envelope. This proves that $\calc \setminus \hat C_0$ lifts to the
envelope. More precisely, we have shown that $\calc \setminus
\hat C_0'$ lift to the envelope, where $\hat C_0'$
is a union of all compact
components of the singular fiber mapped by $\calv $ into the points.

Denote by $\hat C$ some connected subset of $\hat C_0'$
and choose its
neighborhood $V$  not to intersect other components of $\hat C_0'$.
Then $\calv $ maps $V$ onto the neighborhood of the point $\calv (\hat C)$ in
$X$. Now it is clear that this point lifts to the envelope.

Thus, we have proved that $(C_0,v_0)$ lifts to the envelope $(\hat U_Y,\pi )$
of $U$. This implies that $(C_n,u_n)$ should lie in a compct subset of
$\hat U_Y$. This is the desired contradiction.

\smallskip
\qed

\smallskip\rm
Let us give one corollary of the Continuity Principle just proved.

\state Corollary 10.2.3. \it Let $X$ be a complex surface with one singular
normal point $p$. Let $D$ be a domain in $X$, $\d D\ni p$. Suppose there
is a sequence $(C_n,u_n)$ of stable curves over $X$ converging to
$(C_0,u_0)$ in a Gromov topology and such that

a) there is a compact $K\subset D$ with $u_n(\d C_n)\subset K$ for all
$n$;

b) $p\in u_0(C_0)$.

\noindent
Then every meromorphic function from $D$ extends to the neighborhood of
$p$.

\smallskip\rm
To our knowledge this statement is new also for holomorphic functions.

\smallskip\noindent
\sl 10.3. Construction of Envelopes - I.

\smallskip
\rm

We are now ready to give the proof of {\sl Theorem 4.1}, first under 
some technical assumptions to clarify the main idea. Then in the next 
paragraph the general case will be considered.

Let a K\"ahler surface $(X,\omega )$ be fixed and some symplectic 
immersion $u:\ss^2\to X$, with only positive double intersections be given.
Our technical assumptions for the moment will be the following:

\smallskip\noindent
1. We suppose that $c_1(X)[M]$ is not only positive but moreover  
bigger than the sum $\delta $ of double points of $M:=u(\ss^2)$.

\smallskip\noindent
2. $X$ is supposed to be ``positive'' in the sense that for any almost-complex
structure $J\in \calj_{\omega }$ and $J$-complex sphere $C$ satisfies 
$[C]^2\ge 0$.

\smallskip As we have already mentioned, those conditions will be 
eliminated in the next paragraph. 

\medskip Now let $f$ be a meromorphic function in some neighborhood $U$ of 
$M$. Fix some relatively compact 
$U_1\Subset U$ containing $M$.  
Let $\{ J_t\}_{t\in [0,1]}$ be a family of $C^1$-smooth almost
complex structures on the~$4$-manifold $X$ constructed in {\sl Lemma 1.4.2}.
We lift the structures $J_t$ onto the envelope $(\hat U,\pi )$ in such
a way that those liftings are standard on $\hat U\setminus i(U_1)$,
where $i:U\to \hat U$ is a canonical imbedding.

For $t=1$ our $J_1$-holomorphic sphere $M=M_1$ is given. By assumption
$c_1(X)[M]$ $\delta_1 =p \ge1$. Here $\delta_1=\delta $ is the geometric
self-intersection of $M_1$, i.e., the number of double points. From 
{\sl Theorem 8.4.1} we get a family $M_t=f(J_t)$
of $J_t$-holomorphic
curves for $t$ sufficiently close to $1$, where $f:V\subset \calj \to
\cals \times \calj_S \times \calj$ is
a~local section of $\pr_\calj$ constructed in {\sl Theorem 8.4.1}.

\state Definition 10.3.1. 
\sl We say that the~family $\{M_t\}_{t\in(\hat t,1]}$ of
(possibly reducible) $J_t$-holomorphic curves is semi-continuous if
there exist a (maybe infinite) partition of the~ interval $(\hat t,1]$
of the form $1=t_0 >t_1>\ldots>t_n>\ldots$ with $t_n\searrow \hat t$,
natural numbers $1=N_0 \le N_1 \le\ldots \le N_n \le\ldots$,
$J_t$-holomorphic maps $u_t: \bigsqcup_{j=1}^{N_i} S^j_i\to X$
for $(t_{i+1}, t_i]$ with $M_t\deff u_t(\bigsqcup_{j=1}^{N_i} S^j_i)$
such that

\sli $u_t: (t_{i+1}, t_i] \times \bigsqcup_{j=1}^{N_i} S^j_i\to X$
is a continuous map;

\slii $\area(M_t)$ are uniformly bounded from above;

\sliii $\calh$-$\lim_{t\searrow t_{i+1}} M_t \deff \overline
M_{t_{i+1}} \supset M_{t_{i+1}}$.

The~inclusion $\overline M_{t_{i+1}} \supset M_{t_{i+1}}$ means that
$M_{t_{i+1}}$ has no other components than those of $\overline M_{t_{i+1}}$.

We say that $\{M_t\}$ is a family of spheres if all $S^j_i$ are spheres.

\smallskip\rm
Let $T$ be the~infimum of such $\hat t$, for which there is
a semi-continuous family $\{M_t\}_{t\in (\hat t, 1]}$ of spheres such that
for all irreducible components $M^1_t,\ldots M^{N_i}_t$ of $M_t$,
$t\in (t_{i+1}, t_i]$, one has

{\sl a)} $c_1(X)[M^j_t] -\delta_t^j -\varkappa_t^j=p_t^j \ge1$;

{\sl b)} $\sum_{j=1}^{N_i} p_t^j \ge p$.

Here $\varkappa_t^j$ is the sum of the Milnor numbers of cusp points 
of $M_t^j$, see Appendix 2.  We allow the existence of multiple components, 
{\sl i.e.,} that some of $M_t^j$
can coincide.

\smallskip
The set $T$ is open by {\sl Theorem 8.4.1}. To prove the~closeness of $T$, 
we note that
since all $J_t$ are tamed by the~same form $\omega$, the~areas of a
$J_t$-holomorphic curve are uniformly in $t$ estimated from above and below by
$\int_{M_t} \omega$. Moreover, since $\overline M_{t_{i+1}} \supset
M_{t_{i+1}}$, we obtain $\int_{M_{t_i}} \omega \le \int_{M_{t_{i+1}}}
\omega$. In particular, this implies that the~sequence $\{N_i\}$ is bounded
from above, and hence stabilizes for $i$ big enough. From
Gromov's compactness theorem and disk-convexity of the envelope
(Theorem 2.2.2)we obtain that for every $j={1,\ldots,N_i}$ the
sequence $\{M^j_i\}_{i=1}^\infty$ has a subsequence, still denoted in
the~same way, which converges to a $J_{\hat t}$-holomorphic curve $\overline
M^j\subset \hat U_Y$.

To simplify the~notations, we drop the~upper index $j$ from now on
and write $\overline M$ instead $\overline M^j$, having in mind that we can
do the~same constructions for all $j=1,\ldots,N_i$. Note that all irreducible
components of $\overline M$ are $J_{\hat t}$-holomorphic spheres. We write
$\overline M = \bigcup_{k=1}^d m_k\cdot M_k$, where $M_k$ are distinct
irreducible components of $\overline M$ with the~multiplicities $m_j\ge1$.
The genus formula for $\overline M$ now has (because of multiplicities!)
the~following form:

\smallskip
$$
0=\msmall{[\overline M]^2 - c_1(X)[\overline M] \over2}
+ \sum_{k=1}^d m_k - \sum_{k=1}^d m_k(\delta_k +\varkappa_k)
-\sum_{k < l}m_k\cdot m_l [M_k]\cdot [M_l]
$$
$$
-\sum_{k=1}^d \msmall{m_k(m_k -1) \over2 } [M_k]^2.
\eqno(10.3.1)
$$

\smallskip
The~formula can be obtained by taking the~sum of the genus formulas 
for each $M_k$
and then completing the sum $\sum_{k=1}^d m_k[M_k]^2$ to $[\overline M]^2$.
Here $\delta_k$ and $\varkappa_k$ denote self-intersection and Milnor
numbers
of $M_k$. We also have 
$$
\varkappa^j_i + \delta^j_i = \msmall{ [M^j_i]^2 - c_1(X)[M^j_i] \over2} +1.
\eqno(10.3.2)
$$
Note that $[\overline M]^2 = [M^j_i]^2$ and $c_1(X)[\overline M]=
c_1(X)[M^j_i]$ for $i$ sufficiently big, and that $\sum_{k<l} m_k m_l
\,[M_k][M_l] \ge \sum_{k=1}^d m_k -1$, because the~system of curves
$\bigcup_{k=1}^d M_k$ is connected. From (10.3.1) and (10.3.2) we obtain
$$
\sum_{k=1}^d m_k(\delta_k + \varkappa_k) +
\sum_{k=1}^d \msmall{m_k(m_k-1) \over2 } [M_k]^2
\le \delta^j_i + \varkappa^j_i
$$
and
$$
\sum_{k=1}^d m_k c_1(X) [M_k] = c_1(X)[M^j_i].
$$
Using positivity of $(X,\omega)$, {\sl i.e.,} the~property $[M_k]^2 \ge0$,
we obtain
$$
\sum_{k=1}^d m_k \biggl( c_1(X) [M_k] -\delta_k -\varkappa_k\biggr)
\ge c_1(X)[M^j_i] -\delta^j_i -\varkappa^j_i =p_j.
$$
Thus we can choose among $M_1,\ldots,M_d$ the~subset with properties {\sl a)}
and {\sl b)}.

Returning to the~old notations, we can ``correct'' all $\overline M^j$,
choosing an appropriate subset $M^j \subset \overline M^j$. Now applying {\sl
Theorem 8.4.1} to all $M^j$, we can extend our family $\{M_t\}$ in a~
neighborhood of $\hat t$, and hence define it for all $t\in [0,1]$.

For $t=0$ we get the~next statement.

\smallskip
\rm The next rigidity property of symplectic imbeddings is a straightforward
corollary from this theorem.

\state Corollary 10.3.2. \it Let $M_1$ and $M_2$ be two 
symplectically imbedded spheres
in $\cc \pp^2$. Then any biholomorphism of a neighborhood of $M_1$ onto a
neighborhood of $M_2$ is fractional linear.
\rm

\medskip\noindent
{\bigsl 10.4. Construction of Envelopes - II.}

\smallskip Let us now give the proof in the general case. In fact, the proof 
works not only for the meromorphic functions but also for the meromorphic 
mappings into any disk-convex complex K\"ahler space  $Y$. 

\state Theorem 10.4.1. {\it Let $u:\ss^2\to X$ be a symplectic immersion of 
the sphere $\ss^2$ into a disk-convex K\"ahler surface $X$ such that 
$M=u(\ss^2)$.

Let $u:S^2\to X$ be a symplectic immersion of a sphere,
having also positive self-intersections. Let $U$ be a relatively compact
domain in $X$, which contains $M:=u(S^2)$. Denote by $(\hat U_Y,
\hat\pi_Y)$ its envelope of meromorphy relative to $Y$.}

\smallskip\noindent\sl
Step 1. There exists an ($\omega$-tamed) $J_0\in \calj_U$ such that $M$ is
$J_0$-holomorphic.

\smallskip\rm
This was proved in {\sl Lemma 1.1.2} of [I-S]. Moreover, there exists
a smooth homotopy $h :[0,1]\to \calj_U$ joining
$J_0=h(0) $ and $J\st =h(1)$. Put $M_0:=M$ and denote as in {\sl Lemma 8.3.2}
by
$\calm_h(M_0,J_0)$ a component of $\calm_h$ through the point $(M_0,J_0)$.

\smallskip\noindent\sl
Step 2. Suppose that the component $\calm_h(M_0,J_0)$ is not compact.

\smallskip\rm Then by (iii) of {\sl Lemma 8.3.2} there is a continuous path 
$\gamma :[0,1]\to \calm_h$ starting from $(M_0,J_0)$ with property
(b), see (iii) of {\sl Lemma 2.5}. Consider $J_n$-holomorphic spheres $M_n$
in $\calm_h(M_0,J_0)$, which are discrete there with $J_n\to J^*\in \calj_U$
as in \sliii (b) of {\sl Lemma 2.5}.

If for some $n$ $M_n\cap U=\emptyset$,
then, because $J_n=J\st $ on $\hat U\setminus U$, this $M_n$ will be the 
rational curve we are looking for.

If not, then some sequence, still denoted by $M_n$, will converge in the Gromov
topology to a reducible curve $M^{(1)}$. If $M^{(1)}$ has an irreducible
component $M^{(1)}_0$, lying outside of $U$ and such that
$c_1(X)[M^{(1)}_0]>0$, then $M^{(1)}_0$ is our rational curve.

If not, there exists a component $M^{(1)}_0$ of the limit curve $M^{(1)}$
such that $c_1[M^{(1)}_0]>0$. Repeat Step 2 for $M^{(1)}_0$ instead of
$M_0$. Since the area of a  complex  spheres is bounded
from above (see [G]), after a finite number of steps we obtain a
needed rational curve in the envelope $\hat U_Y$, or arrive at 

\smallskip\noindent\sl
Step 3. $\calm_h(M_0,J_0)$ is compact.

\smallskip\rm
From \sliii of {\sl Theorem 8.3.3} and {\sl Corollary 8.3.4} we immediatly
see that there
are continuous (in fact, piecewise $C^k$) paths $(M^n_t,J^n_t)$ such that

1) $M^n_0=M_0$ for all $n$;

2) $J^n_0=J_0$ for all $n$;

3) $J^n_1\to J\st $.

By the Gromov compactness theorem some subsequence from $M^n_1$ , still denoted
by  $M^n_1$, is converging to a $J\st $-holomorphic nodal curve $C^*$
with $\pi^*c_1(X)[C]>0$. Take some irreducible component $C$ of the curve
$C^*$. Then $C$ is the  rational curve in $\hat U_Y$ we are looking for.

\smallskip
\qed

\bigskip\noindent
\sl 10.5. Examples.

\nobreak\smallskip\rm
Here we discuss a few more examples concerning envelopes of meromorphy.

\state Example 1. Let $(X,\omega) = (\cc\pp^1 \times \cc\pp^1,\omega_{FS}
\oplus \omega_{FS} )$, where $\omega_{FS}$ denote the~Fubini-Studi metric on
$\cc\pp^1$. Note that $c_1(X)= 2[\omega]$. Let $J$ be an~$\omega$-tame
almost complex structure on $X$ and $C$ be a $J$ - complex  curve on $X$.
Denote by $e_1$ and $e_2$ the~standard generators  of $\sfh_2(X,\zz)= \zz^2$
and write $[C]=a\cdot e_1 + b\cdot e_2$. Then we get $a+b=\int_C\omega \ge 1$
and $c_1(X)[C]= 2(a+b)$. Furthermore, by the genus formula $0\le g(C) + \delta
+\varkappa = (2ab -2(a+b))/2 +1 =(a-1)(b-1)$. Thus, we conclude that both $a$
and $b$ are non-negative and $[C]^2=2ab \ge0$. So $\cc\pp^1 \times \cc\pp^1$
is non-negative in our sense.

Let $M$ be an imbedded symplectic sphere in $X$. Then $(a-1)(b-1)=0$ by the    genus formula. Therefore, we can assume that $a=1$ and $b\ge0$. 
Now one concludes that the~following holds:

\state Corollary 10.5.1. \it Let $M$ be a~symplectic sphere in $(\cc\pp^1
\times \cc\pp^1,\omega_{FS} \oplus \omega_{FS})$. Then the~envelope of
meromorphy of any neighborhood of $M$ contains a graph of a rational
map of degree $0\le d \le b$ from $\cc\pp^1$ to $\cc\pp^1$.

\rm

\smallskip
\state Example 2. $(X,\omega) = (\cc\pp^2,\omega_{FS})$. Let $M$ be
a~symplectic surface in $X$ of degree $m\deff \int_M \omega$ with {\sl
positive} self-intersections. Then obviously $c_1(X)[M]=3m$.
Note that we proceed to the~construction of
a~family $\{M_t\}$ if the~condition $c_1(X)[M^j_t] > \varkappa(M^j_t)$ is
satisfied for all irreducible components $M^j_t$ of $M_t$. By the~genus
formula one has $\varkappa(C) \le (d-1)(d-2)/2$ for every  complex 
curve $C$ of degree $d$. So we can proceed if $3m > (m-1)(m-2)/2$, which is
equivalent to $1\le m \le8$. Thus, we have the~following

\state Corollary 10.5.2. \it Let $M$ be a~symplectic surface in $\cc\pp^2$ of
degree
$m\le 8$ with positive self-intersections. Then the~envelope of
meromorphy of any neighborhood of $M$ coincides with $\cc\pp^2$ itself.

\state Remark. \rm Note that the~examples include all imbedded symplectic
surfaces in $\cc\pp^2$ of genus $g\le 21$.

\state Example 3. Let $X$ be a ball in $\cc^2$ and $w$ is a standard Euclidean
form. Blow up the origin in $\cc^2$ and denote by $E$ the exceptional curve.
By $\hat X$ denote the blown-up ball $X$. The blow-up of $\cc^2$ is 
also K\"ahler, and we 
denote by $w_0$ some K\"ahler form there. Consider a sufficiently small
$C^1$-perturbation of $E$. This will be a $w_0$-symplectic sphere in $\hat
X$, which is denoted by $M$. The Chern class of the normal bundle to $M$ 
is equal to that
of $E$ and thus is $-1$. Therefore, $c_1(\hat X)[M]=1$ and 
the proof of Theorem 1
applies. One should only note that in the process of deformation $M_t$ cannot
break into irreducible or multiply covered components in this special case.
The only rational curve in $\hat X$ is $E$. Thus, we see that

\it the envelope of meromorphy of any neighborhood of $M$ contains $E$.

\rm One can then blow down the picture to obtain downstairs a sphere
$M_1$-image of $M$ under the blown-down map. This $M_1$ is homologous
to zero, so cannot be symplectic, and for this $M_1$ our Theorem 1
cannot be applied.

\medskip\noindent\bf
Example 4. \rm Chirka in [C] proved the following ``local version'' of our
Theorem 1, which he called  ``a generalized Hartogs' lemma''. Denote by
$\Gamma $ a graph of continuous function $f:\bar\Delta\to \cc $. Consider
the following ``generalized Hartogs' figure'' in $\cc^2$:
$H_\Gamma := \d \Delta \times \bar\Delta\cup \Gamma $. Chirka showed
that every holomorphic function in the neighborhood of $H_{\Gamma }$
extends holomorphically onto the unit bidisk $\Delta^2$. This is a
corollary of our Theorem 4.1 (as is explained in [C]). Really, denote by 
$(z,w)$ coordinates in $\cc^2$. Any function $f$ holomorphic in the 
neighborhood of $\d \Delta \times \bar\Delta $ is a sum of a function 
holomorphic in the bidisk and a function $f_{-}$ holomorphic in $(\cc\pp^1
\setminus \bar\Delta )\times \Delta $. This $f$ is also holomorphic in 
the neighborhood of $\Gamma $. We need only to extend $f_{-}$. Extending 
$\Gamma $ "inside" $(\cc\pp^1\setminus \bar\Delta )\times \Delta $ we 
obtain $\tilde\Gamma $ - a sphere, homologous to the $\{ pt\} \times \cc\pp^1$
in $\cc\pp^1\times \cc\pp^1$. Rescailing the variable $w$ we can make 
$\tilde\Gamma $ symplectic and find ourselves in the conditions of 
{\bf Example 2}.

The proof of Chirka is roughly
the same (i.e., uses the perturbation of the structure), but somewhat simpler.
It uses, instead of Gromov's techniques, the results of Vekua on so- called
generalized analytic functions.

Answering the question, posed by Chirka, Rosay in [R] constructed an example
showing that a ``generalised Hartogs' lemma'' does not hold in
$\cc^3$.

\smallskip
\smallskip\state
Example 5. This example was explained to us by E. Chirka, and it shows
that our Continuity Principle is not valid when the complex dimension of
the manifold $X$ is more than two. 

Take as $X$ the total space of the bundle $\calo(-1)\oplus \calo (-1)$ 
over $\cc\pp^1$. Denote an affine coordinate on $\cc\pp^1$  by $z$, 
coordinates on the fibers by $\xi_1$, $\xi_2$ and $\eta_1$, $\eta_2$ such
that $\eta_1 = z\xi_1$ and $\eta_2 = z\xi_2$. Identify $\cc\pp^1$ with
the zero section of the bundle. Consider the meromorphic function $f= 
e^{\xi_2 /\xi_1}$. The set of essential singularities of $f$ is 
$\{\xi_1=0 \}$, which contains the zero section $\cc\pp^1$. 

Consider the following sequence of analytic disks $C_n$ in $U\deff X \bs
\{\xi_1=0 \}$, $C_n \deff \{ \xi_2 =0, |z|\le n, \xi_1 = {z\over n} \}$. Then
the function $f$ is defined in a neighborhood of every $C_n$. On the other 
hand, the limit curve of the sequence is $C_0 = \cc\pp^1 \cup \Delta_\infty$, 
where $\Delta_\infty \deff \{ \eta_2=0, z=\infty, |\eta_1| \le 1\}$. 
In particular, $f$ does not extend in a neighborhood of $C_0$.

\newpage\noindent
{\bigbf Appendix IV.}

\noindent
{\bigbf  Complex Points and Stein Neighborhoods.}

\medskip\noindent\bf
 \rm The condition $c_1(X)[M]>0$ (in the case of imbedded surfaces,
i.e., when $\delta =0$) in Theorem 4.1 cannot be dropped. In [N-1]
Nemirovsky, using results of Eliashberg-Kharlamov and Forstneric, showed
that any imbedded complex curve $C$ with $c_1(X)[M]\le 0$ can be perturbed
to an imbedded surface $M$ which has a basis of Stein neighborhoods.

\medskip\noindent\sl
A4.1. Complex Points of Real Surfaces in Complex Surfaces.

\smallskip\rm

Let $S\subset X$ be a real surface embedded in a complex surface.
A~point $p\in S$ is said to be {\it complex\/} 
if the tangent plane $T_p S\subset T_p X$ is a complex line.

If locally $S=\Gamma_f$ is the graph of a smooth complex valued function~$f$,
then a point $(z,f(z))\in S$ is complex
if and only if~${\d f\over \d\bar z}(z)=0$.
It is therefore not difficult to show that for generic embeddings
the complex points are isolated, and near each complex point
there exists a complex system of coordinates~$(z,w)$
such that $S$ is given by the equation 
$$
w=|z|^2+{\gamma\over 2}(z^2+\bar z^2)+\bar o(|z|^2),\quad\gamma\in\R.
\eqno (A4.1.1)
$$
If $|\gamma|<1$, then~$p=(0,0)$ is called an {\it elliptic} complex point,
and if $|\gamma|>1$, then it is called {\it hyperbolic}.
(The `parabolic' case~$|\gamma|=1$ is not generic.)
Note that the property of being elliptic or hyperbolic
is invariant under biholomorphic transformations.

It is easy to check that if $S$ is a graph~$\Gamma_f$,
then elliptic and hyperbolic complex points
correspond to positive and negative zeros
of ${\d f\over \d\bar z}$, respectively.

Assume now that an embedded real surface~$S\subset X$
in the general position is compact and oriented.
Then there are two orientations of $T_pS$ at each complex point~$p\in S$.
The first one is the orientation of~$S$
and the second one is the canonical complex orientation of~$T_pX$.
A~complex point is called positive if these two orientations coincide;
otherwise, it is called negative.

Denote by $e_\pm=e_\pm(S)$ and $h_\pm=h_\pm(S)$ the numbers of
positive and negative elliptic and hyperbolic complex points of~$S$.
Introduce {\it Lai's indices\/} $I_\pm:=e_\pm-h_\pm$.
It turns out that~$I_\pm$ are topological invariants
of the embedding.

\proclaim{Lai's Formulae \cite{L}}
Let $S\subset X$ be a compact oriented real surface
in the general position in a complex surface~$X$.
Then
$$
\matrix
I_+ +I_-&=\chi(S)+S^2,\\
I_+ -I_-&=c_1(X)\cdot[S],
\endmatrix
\eqno(A4.1.2)
$$
where $\chi(S)=2-2g(S)$ is the Euler characteristic of~$S$,
$S^2=[S]\cdot[S]$ is the self-intersection index of~$S$ in~$X$,
and $c_1(X)\in \sfh^2(X, \R)$ is the first Chern class of~$X$.
}

\state Proof.
We shall only sketch the proof. For the first formula,
choose a tangent vector field~$\xi\in TS$, twist it by the complex
structure~$J$ on~$X$ and consider the projection
$\eta=\pi(J\xi)\in NS$ to the normal bundle of~$S$.
The zeros of~$\eta$ are exactly the zeros of~$\xi$
and the complex points of~$S$.
Careful inspection of the signs (using local equations~$(*)$)
yields the first formula.

To prove the second formula, pick a non-vanishing 2-form $\omega$ on~$S$
and observe that for a local frame $\xi_1,\xi_2\in TS$,
the expression $\omega(\xi_1,\xi_2)^{-1}\xi_1\wedge_\C\xi_2$
gives a well-defined section of the linear bundle~$\Lambda^2_\C(TX)|_S$.
The zeros of this section are the complex points of~$S$,
because at these points $\xi_1$ and~$\xi_2$ are linearly dependent over~$\C$.
It remains to calculate the signs of these zeros
and to use the definition of~$c_1(X)$.\qed

Sometimes it is convenient to write Lai's formulae
in the following way:
$$
I_\pm(S)={1\over 2}(\chi(S)+S^2\pm c_1(X)\cdot[S]).\eqno(A4.1.3)
$$
For example, if $C\subset X$ is a non-singular
compact complex curve, then after a small real perturbation
we obtain an embedded real surface~$S$
that obviously has only positive complex points.
Thus, $I_-(S)=0$ and therefore~$-\chi(C)=C^2-c_1(C)\cdot[C]$,
which is precisely the adjunction formula for~$C$.

One can also give the following geometric intepretation
of Lai's indices.
Consider the Grassmanian~$Gr(X,\R,2)$
of oriented real 2-planes in~$TX$.
It has two natural oriented submanifolds ${\frak I}_\pm$
formed by complex lines with the complex and anti-complex
orientation, respectively.
Then $I_\pm=\tau(S)\cdot {\frak I}_\pm$, where $\tau:S\to Gr(X,\R,2)$
is the tangential Gauss map.
This again can be proved by calculating local signs.

\smallskip\noindent\sl
A4.2. Cancellation of Complex Points.
\smallskip\rm
Lai's indices $I_\pm$ are essentially the only topological
invariants of the position of embedded real surfaces
with respect to the complex structure.
The following theorem of Kharlamov and Eliashberg
shows that it is possible to cancel
an elliptic and a hyperbolic complex point with the same sign.
In particular, if $I_+=I_-=0$, then the surface is isotopic
to a totally real one.
This can be regarded as a partial case of Gromov's $h$-principle
for $CR$-embeddings.

\proclaim{Cancellation Theorem}
Let $S\subset X$ be an embedded real surface in a general position,
and assume that $e_+(S)h_+(S)>0$. Then there exists
an isotopic embedded surface~$S'\subset X$
such that $e_+(S')=e_+(S)-1$ and~$h_+(S')=h_+(S)-1$.
}

Join the pair of points that we want to remove
by a smooth path in~$S$ and choose a finite covering of this path
by coordinate neighborhoods such that in each chart
$S$~ there is the graph of a smooth function with natural orientation.
Now to move complex points along this path and to eventually
`cancel' them, it suffices to prove the following local result.

\proclaim{Lemma A4.2.1}
Let $f:\bar\Delta\to\C$ be a smooth function
in a neighborhood of the unit disk.
Assume that ${\d f\over \d\bar z}$ has
two ordinary zeros $a,b\in\Delta$ of opposite sign.
Then there exists a function $g:\bar\Delta\to\C$
that is $C^0$ close to~$f$ in the disk,
coincides with~$f$ near the boundary
and has a non-vanishing~${\d g\over \d\bar z}$.

Conversely, if ${\d f\over \d\bar z}$ does not vanish,
then $g$ can be chosen so that~${\d g\over \d\bar z}$
has two zeros of opposite sign at the given points~$a,b\in\Delta$.
}

\state Proof.
By obvious topological reasons there exists
a smooth function~$\psi$ supported in a thin neighborhood
of the segment~$[a,b]$ such that ${\d f\over \d\bar z}+\psi$
does not vanish in~$\Delta$.
Note also that we may choose~$\psi$
so that $\|\psi\|_{C^0}\le 2\|{\d f\over\d\bar z}\|_{C^0}$.
Hence, if the area of $\supp(\psi)$ is small enough,
then the Cauchy-Green transform $T_{CG}\psi $ is $L^{1,p}$-small by 
Calderon-Zygmund inequality and the graph of~$f+T_{CG}\psi $ is totally real.
It remains to multiply $T_{CG}\psi$ by a cut-off function $\phi $
so that $\phi T_{CG}\psi \equiv 0$ near the boundary of~$\Delta$.
This will not spoil $\d(f+\phi T_{CG}\psi)\over \d\bar z$,
because $T_{CG}\psi$ is holomorphic near the boundary
and $C^0$-small.

The second part of the lemma is proved in exactly the same way.\qed

\smallskip The following statement immediately follows from Lai's formulae 
and Cancellation theorem.

\state Corollary A4.2.2. {\it Let $\ss^2\subset \cc^2$ be a symplectically 
imbedded two-sphere in complex projective plane, which is homologous to the 
projective line. Then there is an imbedded two-sphere $\tilde \ss^2\subset 
\cc^2$, which is $C^0$-close to $\ss^2$ and has precisely three points 
with complex tangents. Moreover, those points are positive and elliptic.
}
  
\smallskip\noindent\sl
A4.3. Neighborhoods of Real Surfaces.
\smallskip\rm
Now we can give a construction of Stein neighborhoods
of deformations of real surfaces satisfying certain
topological conditions.

\proclaim{Theorem A4.3.1}
Let $S\subset X$ be an embedded real compact oriented surface
in a complex surface.
If $S$ satisfies both inequalities
$$
2I_\pm=\chi(S)+S^2\pm c_1(X)\cdot[S]\leq 0, \eqno(A4.3.1)
$$
then there exists an isotopic embedded real surface
with a basis of strictly pseudoconvex Stein neighborhoods.
}

This result was proved in~\cite{F} for~$X=\C^2$.
The case where  $S$ is a smooth deformation of a complex curve
was used in~\cite{N1} and the general case in~\cite{N2}.

\state Proof.
By the Cancellation Theorem there exists an isotopic real surface,
denoted by~$S$, that has only {\it hyperbolic\/}
complex points.

Let $p\in S$ be a totally real or a hyperbolic complex point.
It is not difficult to see from the local equations of~$S$
that there exist coordinate neighborhoods~$U\supset\!\supset V\ni p$
and a non-negative smooth function~$\psi$ in~$U$
with the following properties:
\roster
\item"1)" $S\cap U=\{x\in U\mid \psi(x)=0\}$;
\item"2)" $d\psi\equiv 0$ on $S\cap U$;
\item"3)" $\psi$ is strictly plurisubharmonic in~$U\setminus(S\cap V)$.
\endroster

Choose a finite covering of~$S$ by~$V_j\subset\!\subset U_j$.
Let, furthermore, $\chi_j$~be smooth non-negative functions supported in~$U_j$
such that
\roster
\item"1)" $\chi_j\equiv 1$ on $\bar{V_j}$;
\item"2)" $\displaystyle\sum_{|\alpha|\leq 2} |D^\alpha\chi_j|\leq A\chi_j$
everywhere in~$U_j$ for some constant~$A$.
\endroster
Such functions indeed exist: one may first pick functions
satisfying the first condition and then raise them to the third power. 

Now by the Leibnitz rule we have
$$
L(\chi_j\psi_j)=i\,\d\bar{\d}(\chi_j\psi_j)=
i\bigl(\psi_j\d\bar{\d}\chi_j+\d\psi_j\wedge\dbar\chi_j
+\d\chi_j\wedge\dbar\psi_j+\chi_j\d\bar{\p}\psi_j\bigr).
$$
It readily follows that the Levi form~$L(\chi_j\psi_j)$
is non-negative in some neighborhood~$W_j\supset S\cap U_j$
and positive in~$(W_j\setminus S)\cap V_j$.

Consider the function~$\Phi:=\displaystyle\sum_j\chi_j\psi_j$.
By construction
\roster
\item"1)" $\Phi\equiv 0$ on~$S$ and $\Phi>0$ outside of~$S$;
\item"2)" $\Phi$ is strictly plurisubharmonic in a punctured
neighborhood of~$S$.
\endroster
Thus, the sets $U_\eps=\{\Phi<\eps\}$,
where $\eps>0$ is a regular value of~$\Phi$,
form a basis of strictly pseudoconvex Stein neighborhoods of~$S$.\qed

It is worth observing that an embedded real surface~$S\subset X$
with elliptic complex points cannot have
a Stein neighborhood basis.
Indeed, according to Bishop for an elliptic point~$p\in S$,
there exists a non-trivial continuous family
of holomorphic discs~$f_t:(\Delta,\d\Delta)\to (X,S)$ with~$f_0\equiv p$.
By the continuity principle every  holomorphic function
in a neighborhood of~$S$ extends holomorphically along this family,
and therefore a sufficiently small neighborhood of~$S$
cannot be Stein.

In fact, it can be proved that, at least when~$S^2\geq 0$,
a surface~$S$ that `potentially' has elliptic complex points
(which means that one of the indices~$I_\pm(S)$ is positive)
does not admit a (not necessarily small) Stein neighborhood
(see~\cite{N2} and~\cite{N3}).

\smallskip
Now, consider a complex surface~$X$ and a smooth, compact complex curve 
$C\subset X$. 

\smallskip
\state Corollary A4.3.2. {\it If $c_1(X)\cdot[C]\le 0$ then a generic embedded 
real surface $S$ irotopic to $C$ has a basis of Stein neighborhoods. In 
particular such $S$ has no non trivial envelope of meromorphy.
}

\state Proof. {\rm If $c_1(X)\cdot[C]\le 0$,
then a generic embedded real surface~$S$ isotopic to~$C$
satisfies~$I_\pm(S)\le 0$, because of the adjunction formula
and Lai's formulae. The statement now follows from the Theorem A4.3.1.
}

\smallskip
\qed

\newpage
\noindent\bigbf
Appendix V. Seiberg-Witten Invariants and Envelopes.

\medskip\noindent
{\bigsl A5.1. The Genus Estimate.}
\medskip\rm
Throughout this Appendix $X$ will denote a compact, real, oriented  manifold 
of dimension $4$. Fix some $c\in H^2(X, \rr )$. Let $M$ be a closed, oriented 
surface of genus $g(M)$ imbedded in $X$. For a riemannian metric $g$ on 
$X$ denote by $s_g$ the scalar curvature of $g$, i.e. 
$s_g=g_{ij}R^{ij}={\sl tr}R_g$-trace of the  Ricci curvature. In the case 
when $X=M$ is a compact surface the Gauss-Bonnet theorem tells us that 

$$
{1\over 4\pi }\int_Ms_gd\mu_g = 2-2g(M), \eqno(A5.1.1)
$$
\noindent where $\mu_g$ is the area form associated to $g$.

   The following statement is due to P. Kronheimer. 

\state Theorem A5.1.1. {\it Assume that $[M]^2=0$ and that for every 
riemannian metric $g$ on $X$ one can find a two-form $\omega $ representing 
$c$ such that 

$$
\Vert \omega \Vert_{L^2_g}\le {1\over 4\pi }\Vert s_g\Vert_{L^2_g} + K(X),
\eqno(A5.1.2)
$$
\noindent
where constant $K(X)$ depends only on the topology of $X$. Then

$$ 
\vert c[M]\Vert \le \max \bigl\{ 2g(M)-2,0\bigr\} .\eqno(A5.1.3)
$$
}
\state Proof. {\rm 
\noindent Take some 
neighborhood $U$ of $M$ which is diffeomorphic to $M\times D$, where $D$ is 
a disk in $\rr^2$ with coordinates $x,y$.  Denote by $g_M$ the 
riemannian 
metric on $M$ of constant scalar curvature $s_g$ equal to $4\pi (2-2g(M))$.   

(A5.1.1) tells us that  if $M$ is not a torus, then the ${\sf area}_{g_M}M=1$.
In the case of torus we just take the flat metric with total area one. For 
real positive $R$ consider the metrics $g_{D,R} = R(dx^2 + dy^2)$ on $D$.

Take now the following metric on $X$:

$$
g_R = \cases g_M\oplus g_{D,R} & on \hbox{ } $U$ \cr
              g & away \hbox{ } from \hbox{ } $U$ \cr
\endcases              
$$ 

Let $\omega_R$ be the harmonic representative of 
$c$ with respect to the metric $g_R$. Harmonic form minimases the 
$L^2$-norm in its cohomology class. So by the assumption (A5.1.2) of our 
theorem 

$$
\Vert \omega_R\Vert_{L^2_{g_R}}\le (2g(M)-2)\sqrt{\pi }R + O(1). \eqno(A5.1.4)
$$
Really, $\Vert \omega_R\Vert_{L^2_{g_R}}\le {1\over 4\pi }\vert s_{g_R}\Vert
+ O(1) $ $= \bigl(\int_{M\times D}[4\pi (2g-2)]^2d\mu_{g_R}\bigr)^{{1\over 2}}= 
{1\over 4\pi }\bigl(\pi R^2\cdot 4\pi (2g(M)-2)\bigr)^{{1\over 2}} + O(1) $
$ = \sqrt{\pi } \cdot R\cdot (2g(M)-2) + O(1)$. 

Further 

$$
 \Vert \omega_R\Vert_{L^2_{g_R}}\ge \vert c[M]\vert \sqrt{\pi }R.
 \eqno(A5.1.5)
$$
\noindent
Really, $\Vert \omega_R\Vert^2_{L^2_{g_R}}\ge \int_D\Vert \omega\mid_{M\times 
(x,y)}\Vert^2_{L^2_{g_M}}R^2dxdy$ because $\omega\mid_{M\times (x,y)}$ is 
only one of the coefficients of $\omega $. $\Vert \omega\mid_{M\times 
(x,y)}\Vert^2_{L^2_{g_M}}$ is bounded from below by the $L^2$-norm of the 
corresponding harmonic representative. Note that the harmonic representative 
of $\omega\mid_{M\times(x,y)}=c\mid_M$ is just $c[M]*_{g_M}1$ 
because 
Hodge $*_{g_M}$ commutes with the Laplacian: 
$\Delta_{g_M}: \Delta_{g_M}(*_{g_M}1)=*_{g_M}(\Delta_{g_M}1)=0$. So the 
$L^2$-norm of this harmonic representative and thus of 
$\omega\mid_{M\times (x,y)}$ is bounded from below by $c[M]$. Inequality
(A5.1.5) follows.

Inequalities (A5.1.4) and (A5.1.5) clearly imply inequality (A5.1.2).

%\medskip\noindent{\bigsl Step 2. Case $[M]^2>0$.} Consider now the general 
%case of the surface with positive homological self-intersection.
}             
\qed

This theorem reduces the problem of minimazing the genus of an imbedded 
surfaces in a fixed homology class to the a priori estimate (A5.1.2). This 
estimate can be obtained if the Seiberg-Witten invariant $SW(X,c)$ for our 
cohomology class $c$ is not zero.

\medskip\noindent
{\bigsl A5.2. The use of Seiberg-Witten Invariant.}
\smallskip 

The intersection form on the second cohomology space~$H^2(X,{\Bbb R})$ of a 
compact $4$-manifold $X$ is symmetric, and
hence it has Sylvester indices~$b^\pm(X)$.  By definition, the signature
of~$X$ is the signature $\sigma(X)=b^+ -b^-$ of its intersection form.

If $X$ is a K\"ahler complex surface, then the
maximal positive subspace in~$H^2(X,{\Bbb R})$ is spanned by the class of
the K\"ahler form~$\omega$ and by the classes of the real and imaginary
parts of holomorphic $2$-forms on~$X$. In particular,
$b^+(X)=1+2\dim_\cc\Omega^2(X)$ is an odd positive integer.

Suppose now that for our $4$-manifold $X$ we have that $b^+(X)>1$.  Let us
pick a characteristic class $c\in H^2(X,{\Bbb Z})$, that is, an integral
cohomology class which is congruent~$\mathop{\roman{mod}} 2$ to the second
Stiefel--Whitney class~$w_2(X)\in H^2(X,{\Bbb Z}_2)$.  (For example, the
first Chern class of an almost-complex structure~$J$ on~$X$ is
characteristic.)

In this situation one can define the {\it Seiberg--Witten invariant\/}
$SW(X,c)\in{\Bbb Z}/2{\Bbb Z}$ of the smooth oriented manifold~$X$ (if
$\psi$ is an orientation-preserving diffeomorphism of~$X$,
then $SW(X,\psi^*c) =SW(X,c)$).

Roughly speaking, the invariant~$SW(X,c)$ is non-zero if the moduli space
of gauge equivalence classes of solutions of the Seiberg--Witten
differential equations for a Riemannian metric on~$X$ has non-zero homology
class in the configuration space.  For the precise definition of the
invariant, the reader is referred to the original paper [Wt] and to the
books [Wg] and [Mr].

\state Theorem of Witten - I. {\it 
Let $X$ be a compact K\"ahler surface with~$b^+>1$.  Then the moduli space of 
solutions of the (suitably perturbed) Seiberg--Witten equations can be 
identified with the space of all effective divisors that are Poincar\'e dual 
to the cohomology class $\frac12\bigl(c-c_1(X)\bigr)$.

In particular, $SW(X,c_1(X))=1$.
}

Really, if $c$ is chosen to be equal to $c_1(X)$ there is only one divisor 
Poincar\'e dual to the zero class, namely the zero divisor. Therefore the 
moduli space of the solutions of $SW$-equations in this case is a singleton.

\state Remark. {\rm In the opposite direction, the Seiberg--Witten invariants 
of~$X$
can be used to produce complex curves with given cohomology class,
namely, if $SW(X,c)\ne 0$, then there is a complex curve dual
to~$\frac12\bigl(c-c_1(X)\bigr)$. This statement is due to Taubes, see [Ta].
}

\smallskip
Suppose that, for a our pair $(X,c)$ the invariant~$SW(X,c)$ is
non-zero. Then, by definition, for every metric~$g$ on~$X$, there is a
solution of the Seiberg--Witten equations.  Using the properties of the
equations, one can obtain {\it a priori\/} estimates of the solutions
in terms of the metric~$g$.  We do not write out the equations themselves
and only state the key estimate of this type, which is essentially due to
Witten (see also [Mr] and [Kr]).

\smallskip
\state Theorem of Witten - II. {\it If
$SW(X,c)\ne 0$  then for every  Riemannian metric $g$ on $X$ with there is a
closed two-form~$\eta$ representing this class in the de~Rham cohomology
and such that
$$
\Vert \eta \Vert_{L_g^2}\le {1\over 4\pi }\Vert s_g\Vert_{L^2_g} + K(X),
\eqno(A5.2.1)
$$
where $K(X)$ is a constant depending only on the topology of~$X$.
}
\smallskip
As a consequence we see that

\smallskip
\state Corollary A5.2.1. {\it  The number of classes $c\in H^2(X,\zz )$ 
with non-trivial Seiberg--Witten invariants $SW(X,c)$ is finite.
}
\state Proof. {\rm   Indeed,
for a fixed metric~$g$, the harmonic representatives of these integer
classes are uniformly bounded.

}
\qed 

\smallskip
\state Corollary A5.2.2 {\it Let $X$ be a compact K\"ahler surface with 
$b^{+} >1$.  Let $M\subset X$ be an embedded real surface with non-negative 
self-intersection $M^2\ge0$. Then 
$$
\left| c_1(X)\cdot M \right|+M^2\le \max\{0,2g(M)-2\}. \eqno(A5.2.2)
$$
}
\state Proof. {\rm 
Let us choose an orientation on~$M$ so that~$c_1(X)\cdot M\le 0$. 
We blow up the surface $X$ at $d=M^2$ distinct points and denote the
blown-up surface by $\widetilde X$, the exceptional curves by $E_j$, and
the corresponding blow-down by $\pi:\tilde X\to\nomathbreak X$.
Then $c_1(\tilde X)=\pi^*c_1(X)-\sum_j [E_j]$.

Let $\Sigma$ be the interior connected sum of~$M$ with $\overline{E_j}$
in~$\widetilde X$.  Then
$[\Sigma]=\pi^*[M]-\sum_j [E_j]$, and therefore $\Sigma^2=0$. Observe 
that the genus of $\Sigma$ is equal to that
of $M$.

By Theorem A5.1.1 we have that 
$$
\max\{0,2g(M)\} \ge\bigl|c_1(\tilde X)\cdot\Sigma\bigr|
=\left| c_1(X)\cdot M-M^2\right|=\left| c_1(X)\cdot M \right|+M^2.
$$
}
\qed

\smallskip Corollary A5.2.2, as it is stated, makes no difference between 
tori and spheres, i.e. $\max\{0,2g(M)\} =0$ in both cases. However for 
spheres it can be improved. We shall derive now the following result of
Kronheimer and Mrowka [Kr-M] and Morgan, Szab\'o, and~Taubes [M-S-T].

\smallskip
\state Corollary A5.2.3. {\it  
An embedded oriented real surface~$M\subset X$ with non-negative
self-inter\-sec\-tion index in a compact K\"ahler surface with~$b^+(X)>1$
satisfies the inequality
$$
M^2+\bigr|c_1(X)\cdot M\bigl|\le 2g(M) - 2 \eqno(A5.2.3)
$$
provided that  $M$ is not a $2$-sphere with trivial real homology class.
}
\state Proof. {\rm  The case $g(M)\ge 1$ is covered by (A5.2.2). So let 
$M$ be a sphere. (A5.2.2) tells that $M^2=0$ and~$c_1(X)\cdot M=0$.

Let us blow up the surface~$X$ at one point.  Let $E\subset \tilde X$ be the
exceptional curve.  For every $n\ge 1$, we consider the interior connected
sum~$M_n=E\#nM$.  This is an embedded two-sphere such that $M_n^2=-1$
and~$c_1(\tilde X)\cdot M_n=1$.

Let us use the existence of an orientation-preserving
diffeomorphism $\psi\:\tilde X\to\tilde X$ (with support in a neighbourhood
of~$M_n$) that maps~$M_n$ into itself with opposite orientation.  We have
$\psi^*c_1(\tilde X)=c_1(\tilde X)+2[M_n]$.  If the
class~$[M]\in H^2(X,\rr )$ were non-zero, this would give us an
infinite set of cohomology classes on~$\tilde X$ with non-zero
Seiberg--Witten invariant, which contradicts Corollary A5.2.1.
}
\qed

\state Remark. {\rm 
By the usual adjunction formula, a complex curve $C$ with $C^2\ge 0$
violates our inequality if and only if~$c_1(X)\cdot C>0$.  This situation
occurs in~$\cc\pp^2$ and other complex surfaces with negative
canonical class indeed.  For this reason, it is impossible to drop the
assumption~$b^+>1$ in the above theorem.
}

In A5.3 we shall repeatedly use the fact that complex curves in complex
surfaces with $c_1(X)>0$ violate the ``general'' form of the adjunction
inequality.

\medskip\noindent
{\bigsl A5.3. The Genus Estimate on Stein Surfaces and Envelopes.}

\smallskip We start from the following 

\state Theorem of Stout. {\it 
Let $K\subset Y$ be a compact set in a Stein manifold.  Then there is a 
neighbourhood $V$ of $K$ and  biholomorphic imbedding $h$ of $V$ onto an open 
subset in an affine algebraic variety $X$. Moreover, if $Z\subset Y$ is a 
complex submanifold, then one can find such $h$ that $Z\cap V$ will be  mapped 
onto the intersection of an algebraic subvariety in $X$ with $h(V)$.
}

We shall see now that imbedded real surfaces in
Stein surfaces remind real surfaces in compact K\"ahler surfaces with $b^+>1$.

\smallskip
\state Theorem A5.3.1. {\it 
Let $M\subset Y$ be an imbedded oriented real surface in a Stein surface with 
non-negativ self-intersection.
Suppose that $M$ is not a $2$-sphere with trivial homology class. Then
$$
S^2+\bigr|c_1(X)\cdot S\bigl|\le 2g-2. \eqno(A5.3.1)
$$

}

\state Proof. {\rm 
We may always assume that the homology class of~$M$ is non-zero because
our inequality becomes trivial if~$[M]=0$ and $M$ is~not a sphere.  Hence,
there is a non-singular (non-compact) complex curve~$L\subset Y$ such
that~$L\cdot M\ne 0$.  Indeed, $\sfh_2(X,\zz)$ is torsion-free for
any Stein surface, and hence by the Poincar\'e duality we can find a class
$c\in H^2(Y, \zz )$ with $c\cdot M\ne 0$.  Since $Y$ is a Stein surface,
it follows that $c=c_1({\Cal L})$ for a holomorphic line bundle~$\Cal L$.
Then~$L$ is the divisor of a generic holomorphic section of~$\Cal L$.

Using the algebraic approximation theorem of Stout we obtain an open
imbedding of a neighbourhood of~$M$ into an affine part $X_a$ of a compact 
algebraic surface $X$ such
that~$L$ is mapped to an algebraic curve~$C\subset X$ so that
$C\cdot Y=L\cdot M\ne 0$. Therefore, the real homology class of~$M$
remains non-trivial after the algebraization.

Denote by $H=X\setminus X_a$-"the curve at infinity".  
We have non-zero classes $[M]$ and~$[H]$ in~$H^2(Y,\rr )$ such that
$[M]^2\ge 0$, \ $[H]^2>0$, and $[H]\cdot [M]=0$.  Since the intersection
form is non-degenerate, it immediately follows that $b^+(Y)>1$.  Hence, the
desired inequality is a consequence of Corollary A5.2.3.

}
\qed

\state Remarks. 
{\bf 1. \rm 
It follows from Theorem A4.3.1 that the inequality A5.3.1 is sharp, namely, an
imbedded real surface satisfying this inequality is isotopic to a surface
with a basis of Stein neighbourhoods.  

\smallskip\noindent
{\bf 2.} This theorem is also true for the surfaces $M$ with negative 
self-intersection, see [N-4]. 

\smallskip\noindent
{\bf 3.} Concerning a real surface with non-negative
self-intersection index the Theorem A5.3.1 still holds  for non-Stein strictly 
pseudoconvex domains
as well.  This can proved by the same argument (see [N-2] and [N-3])
by using the Lempert algebraic approximation theorem [Le] for such
domains.

\smallskip\noindent
{\bf 4.} {\sl Immersed Surfaces}
We assume that a surface~$S$ is immersed with $\kappa_+$~positive and
$\kappa_-$~negative ordinary double points.  Then, after embedding a
neighbourhood of~$S$ in an algebraic surface~$Y$, we replace each
positive double point with an embedded handle and perform a blow-up at each
negative double point.  This operation increases the genus by~$\kappa_+$
and does not change~$S^2$ and~$c_1(Y)\cdot S$.  Hence, Theorem~A5.3.1 holds 
with $g$ replaced by~$g+\kappa_+$.
}

\smallskip Now we shall apply the adjunction inequalities (or Genus estimates) 
for imbedded real
surfaces in Stein surfaces to the study of envelopes of holomorphy and 
meromorphy of domains in complex surfaces.  

The main idea is fairly simple.  Suppose that we can solve the Levi problem
for domains over a complex surface~$X$, that is to prove that if $D$ is a 
pseudoconvex  domain over $X$ then $D$ is a Stein domain unless some very 
 special case occurs. I.e. for example $D=X$ (this is the case of $X=\cc\pp^2$)
 or $X$ is foliated by complex curves (case of $X=\cc\pp^1\times \cc\pp^1$.  
 Then Theorem~A5.3.1 yields then topological restrictions on real surfaces in 
 $D$.  The below are due to Nemirovski, see [N-2].

Let $M\subset\cc\pp^2$ be an embedded oriented real surface in the
complex projective plane.  By definition, the degree of~$M$ is the
intersection index of~$M$ with complex lines, $d=M\cdot [\cc\pp^1]$.  We
can choose the orientation of~$M$ so that~$d\ge 0$ and note that $M$ is
homologous to zero if and only if~$d=0$.

For an embedded surface of degree~$d$, we have~$M^2=d^2$ and
$c_1(\cc\pp^2)\cdot M=3d$. In particular, if~$S$ is a smooth complex
curve, then the classical genus formula reads
$g(M)=g_\Bbb C(d)=\frac12(d^2-3d+2)$.

\state Corollary A5.3.2. {\it  
Let $M\subset \cc\pp^2$ be an imbedded oriented real surface of positive
degree~$d>0$ and genus~$g$.  If there is a non-constant holomorphic
function in a neighbourhood of~$M$, then
$$
g\ge\frac12(d^2+3d+2)=g_\Bbb C(d)+3d.
$$

Conversely, if $M$ satisfies this inequality, then it is isotopic to an
embedded surface with a basis of Stein neighbourhoods.
}

\state Corollary: Vitushkin's conjecture. {\it 
If an embedded two-sphere in~$\cc\pp^2$ is not homologous to zero, then
every holomorphic function in a neighbourhood of this sphere is constant.

}

\state Proof. {\rm 
Let us consider a neighbourhood~$U\supset S$.  If there is a non-constant
holomorphic function in~$U$, then, by Theorem~1, the envelope of
holomorphy~$\tilde U\supset U$ is a Stein domain because it cannot coincide
with~$\cc\pp^2$.  Hence, the desired inequality immediately follows
from Theorem~A5.3.1.

Conversely, %%added
if the inequality of the theorem is satisfied, then, by Theorem~4, there
is an isotopic surface with a basis of Stein neighbourhoods.
}
\qed
\smallskip\noindent
\state Corollary A5.3.3. {\it 
If $M\subset \cc\pp^2$ be an imbedded oriented real surface of positive
degree~$d>0$ and genus~$g$ with 

$$
g < \frac12(d^2+3d+2)=g_\Bbb C(d)+3d.
$$
then the envelope of meromorphy of any neighborhood of $M$ consides with 
$\cc\pp^2$.
}

\smallskip
For a real oriented surface $M\subset\cc\pp^1\times\cc\pp^1$, we
consider its ``bidegree''~$d=(d_1,d_2)$, that is, the pair of intersection
indices with horizontal and vertical complex lines.  Then
$c_1(\cc\pp^1\times\cc\pp^1)\cdot S=2(d_1+d_2)$ and $S^2=2d_1d_2$.

A smooth complex curve $C\subset\cc\pp^1\times\cc\pp^1$ has
``positive'' bidegree $d\in{\Bbb Z}^2_+\setminus\{0\}$.  The genus of~$C$
is given by the formula $g_\cc(d)=d_1d_2-d_1-d_2+1$.

\state Corollary A5.3.4. {\it 
Let $S\subset\cc\pp^1\times\cc\pp^1$ be an embedded oriented real
surface of non-zero bidegree~$d=(d_1,d_2)\in{\Bbb Z}^2\setminus\{0\}$ and
genus~$g$.  If there is a non-constant holomorphic function in a
neighbourhood of~$S$, then either
$$
g\ge d_1d_2+|d_1+d_2|+1,\quad d_1d_2>0,
\tag$**$
$$
or $d_1d_2=0$ and the envelope of~$S$ contains a vertical or a horizontal
line.

Conversely, if $S$ satisfies~$(**)$, then it is isotopic to an embedded
surface with a basis of Stein neighbourhoods.

}

\state Remark. {\rm In particular, it is impossible to deform a complex curve
in $\cc\pp^1\times\cc\pp^1$ other than a vertical or a horizontal
line into a real surface with holomorphic functions in its neighbourhoods.
}

\state Proof. {\rm 
By the Theorem of Fujita the envelope of holomorphy of a domain
$U\subset\cc\pp^1\times\cc\pp^1$ either is a Stein domain or
contains one of the fibers.  We also note that if $\tilde U$ contains,
say, a vertical line and intersects the other vertical lines, then
$\tilde U$~coincides with the entire $\cc\pp^1\times\cc\pp^1$ by the usual 
Hartogs Theorem.

It remains to apply Theorems~A5.3.1  and~A4.3.1 in the same way as above.
}
\qed

\smallskip
Let us consider the direct product $X=\cc\pp^1\times Y$, where~$Y$ is a
non-compact Riemann surface.  We recall that $H^2(Y,{\Bbb Z})=0$, and hence
$H^2(X,{\Bbb Z})\cong H^2(\cc\pp^1,{\Bbb Z})\cong{\Bbb Z}$ by the
K\"unneth formula.  For a real oriented surface~$S\subset X$, let
$d=\deg\pi_{\cc\pp^1}$ be the topological degree of its projection to
the $\cc\pp^1$-factor.  Then $S^2=0$ and $c_1(X)\cdot S=2d$.  We also
note that we can orient~$S$ so that~$d\ge 0$.

By the classical Behnke--Stein theorem, $Y$ is a Stein manifold. Hence,
Theorem~A4.3.1 gives the following description of envelopes of holomorphy
over~$X$.  The envelope of a domain $U\subset X$ either is a Stein domain or
contains a fiber~$\cc\pp^1\times\{y\}$.  In the latter case, holomorphic
functions in~$U$ are constant in the $\cc\pp^1$-direction, and therefore
the envelope is of the form $\cc\pp^1\times\pi_Y(U)$.

Using the same argument as above, we obtain the following result.

\state Corollary A5.3.5. {\it 
Let $S\subset\cc\pp^1\times Y$ be an embedded oriented real surface of
genus~$g$ and positive degree~$d>0$.  If $g<d+1$, then every holomorphic
function in a neighbourhood of~$S$ extends to the
set~$\cc\pp^1\times\pi_Y(S)$.

If $g\ge d+1$, then $S$ is isotopic to a surface with a Stein neighbourhood
basis.
}

\smallskip
\state Remark. The results of this section are much more general then 
the {\sl Theorem 4.1} in the case of $X=\cc\pp^2, \cc\pp^1\times \cc\pp^1$. 
For example, if a symplectic sphere $M\subset \cc\pp^2$ is in the homology 
class of the line, then the procedure of the proof of the Theorem 4.1 
gives an isotopy of $M$ to a complex line. However there are knoted 
two-spheres in $\cc\pp^2$ in the homology class of the line (but they are of 
course non symplectic). From the other hand the Theorem 4.1 is valid 
for the spheres not nessessarily in compact K\"ahler surfaces, where the 
envelopes of holomorphy or meromorphy are not nessessarily Stein. And for a  
non-Stein envelope the existence of a compact curve is also a strong 
conclusion. 

To our understanding of the subject  a statement containing the both type of 
results is not avaiable  by neither Gromov no Seiberg-Witten techniques.

\medskip\noindent
{\bigsl A5.4. Two-spheres in $\cc^2$.}
\smallskip
In this section we apply Theorem A5.3.1 to embedded spheres in $\cc^2$ and
its blow-ups.  This enables us to show that an embedded 2-sphere
in~$\cc^2$ is homologous to zero in its envelope of holomorphy. 

It follows from Theorem A5.3.1 that an embedded sphere~$S$ with 
self-intersection
index~$S^2=0$ in a Stein complex surface~$X$ is homologous to zero
in~$\sfh_2(X,\zz)$.  For an embedded two-sphere~$S\subset\cc^2$,
the condition $S^2=0$ holds because $\sfh_2(\cc^2,\zz)=0$.  Hence,
$S$ is homologous to zero in every Stein domain over~$\cc^2$ that
contains a neighbourhood of~$S$.  Since the envelope of holomorphy of
every domain in~$\cc^2$ is a Stein domain, we obtain the following
result [N-3].

\state Corollary A5.4.1. {\it 
An embedded two-sphere $S\subset\cc^2$ is homologous to zero in the
envelope of holomorphy of every neighbourhood~$U\supset S$.
}

One can note that the theorem does not work for {\it immersed\/} spheres.
Indeed, a simple modification of the results in Appendix IV for immersed 
surfaces
(we must replace $S^2$ in Lai's formulas by the Euler number of the normal
bundle) shows that there is an immersed totally real 2-sphere in~$\cc^2$
with a single positive double point.  This sphere has a basis of Stein
neighbourhoods in each of which this sphere is certainly not homologous to
zero.

\medskip\noindent
{\bigsl A5.5.  Attaching complex disks to strictly pseudoconvex domains.}

\smallskip 
Let $U$ be a relatively compact domain with smooth strictly pseudoconvex 
boundary in a Stein surface. We say that an embedded analytic disk 
$D\subset X$,
smooth up to the boundary, is attached to~$U$ (from the outside) if
$D\subset X\setminus U$ and $\partial D\subset\partial U$.  We shall always
assume that~$D$ is transversal to~$\partial U$.

For example, the disk~$\{|z|\le1,\,w=0\}\subset{\Bbb C}^2$ is attached to the
boundary of the $\eps$-neighbourhood of the circle $\{|z|=1+\eps,\,w=0\}$.
We note that the boundary of this disk is not homologous to zero in~$U$.

A more intricate example can be obtained as follows.  Let us consider a
totally real torus~$T\subset{\Bbb C}^2$.  Reversing the proof of 
the Cancellation Theorem,
we create an elliptic and a hyperbolic point in~$T$.  Let us replace a
neighbourhood of the elliptic point by a complex disk and observe that
this disk is attached to a strictly pseudoconvex neighbourhood of the
remaining part of the torus.  The boundary of this disk bounds a real
surface with one handle inside the domain.  An important observation
(due to Forstneri\v{c} [F]) is that, by Lai's formulas, this
construction does not work for a sphere.  In fact, it is impossible
to find an analytic disk and a disk with only hyperbolic complex points
with common boundary and without common interior points.

We generalize these examples following [N-3] and prove that the
boundary of a holomorphic disk attached to a strictly pseudoconvex domain
in~${\Bbb C}^2$ is never sliced by a smooth disk inside the domain.

\state Theorem A5.5.1. {\it 
Suppose that an analytic disk~$D$ is attached to a strictly pseudoconvex
domain~$U$ in a Stein surface~$X$ with $\sfh_2(X,\rr)=0$.  Then there is
no smooth disk inside of~$U$ that has the same boundary.
}

\state Proof. {\rm 
We assume that there is a smooth disk in~$U$ with boundary~$\partial D$.
Smoothing the union of our two disks, we obtain an embedded
sphere~$S\subset X$.  The na\"{\i}ve idea is to paste a tubular
neighbourhood of~$D$ to~$U$ so that the result is a Stein domain.
Then~$S$ is not homologous to zero in this domain, which contradicts
Theorem~9.  The following basic example shows first that such a pasting is
impossible and then provides a remedy for this problem.

Let us consider the standard disk $D=\{w=0,\,|z|\le 1\}$ in~${\Bbb C}^2$
with a strictly pseudoconvex ``collar''
$C=\{ |w|^2+(|z|-a)^2\le (a-1)^2,\,|z|\le 1+\eps\}$ for some~$a>1$
and~$\eps>0$.  If we glue a small neighbourhood of~$D$ to the collar, then
the interior of this set contains a Hartogs figure with appropriate
``walls,'' and therefore this interior is not a Stein domain.  However,
we can construct a Stein neighbourhood ``pinched'' at the origin.

To see this, we note that, in logarithmic coordinates~$\xi=\ln|z|$,
$\eta=\ln|w|$, the domain on the $|z|>1$~side of~$C$ corresponds to the
subgraph of a smooth convex function $\eta=f(\xi)$,
$\xi\in(0,\ln(1+\eps)]$.  The $\eta$-axis is clearly a vertical asymptote
of~$f$.  Let us choose a~point $\xi_0\in(0,\ln(1+\eps)]$ such that the slope
of~$T_{\xi_0}\Gamma_f$ is a positive integer~$n\ge 1$, and replace $f$ by
this tangent for all~$\xi\le\xi_0$.  The subgraph of this new function
defines a neighbourhood of~$D\setminus\{0\}$ that is logarithmically convex,
and hence a Stein domain.  Near the origin, this neighbourhood is of the
form~$|w|\le K|z|^n$, $K>0$.

If we blow up the origin in~${\Bbb C}^2$, that is, make the
change of variables~$z=z'$, $w=w'z'$, then the proper preimage of our
neighbourhood becomes~$|w'|\le K|z'|^{n-1}$.  Thus, after
$n$~blow-ups, we obtain an ordinary neighbourhood of the proper preimage
of~$D\cup C$.  By construction, this domain is locally pseudoconvex and
intersects only the last exceptional curve by a disk.  Hence, it is a
Stein domain by the remark after Takeuchi's.

The situation of the theorem is reduced to this model by a suitable choice
of coordinates in a pinched neighbourhood of~$D$ and by a small perturbation
of~$U$.  Thus, we have the following construction of Stein neighbourhoods.
}

\smallskip
\state Lemma A5.5.2. {\it  
Suppose that an analytic disk~$D$ is attached from the outside to a
strictly pseudoconvex domain~$U$ in a Stein surface $X$.  Then the proper
preimage of $D\cup U$, after several blow-ups at a point~$p\in D$, is
contained in a Stein domain.
}

To complete the proof of the theorem, we recall that
$\sfh_2(X,\rr)=0$, and thus $S^2=0$ and $c_1(X)\cdot S=0$.  Consequently,
for the proper preimage $\tilde S\subset\tilde X$ of the sphere~$S$, after
$n\ge 1$~blow-ups at~$p\in S$ we have ${\tilde S}^2=-n$ and
$|c_1({\tilde X})\cdot{\tilde S}|=n$, which contradicts the inequality of
Theorem A5.3.1.
}
\qed

\smallskip
\state Corollary: Vitushkin's conjecture. {\it 
It is impossible to attach an analytic disk from the outside to a strictly
pseudoconvex domain in~${\Bbb C}^2$ diffeomorphic to the ball.
}

\state Proof. 
The idea is to ``reflect'' the disk inside the domain by inversion.

Let $U\subset{\Bbb C}^2$ be a strictly pseudoconvex domain in~${\Bbb C}^2$
diffeomorphic to the ball. More precisely, we assume that there is a
diffeomorphism $\psi\:\bar U\to\bar B$ of manifolds with boundary. Then
$\psi$ extends to a diffeomorphism $\Psi$ of~${\Bbb C}^2$ that maps~$U$ to
the standard ball~$B$.

If a holomorphic disk~$D$ is attached to $U$ from the outside, then its
image~$\Psi(D)$ is a smooth disk attached to~$B$.  Applying the inversion
with respect to~$\partial B$, we obtain a smooth disk~$D'$ in~$B$ with the
same boundary.  The preimage $\Psi^{-1}(D')$ is a smooth disk in~$U$ with
boundary~$\partial D$, which contradicts the previous theorem.
\qed

\newpage

\spaceskip=4pt plus3.5pt minus 1.5pt
\xspaceskip=5pt plus4pt minus 2pt
\font\csc=cmcsc10
\font\tenmsb=msbm10
%\font\sevenmsb=msbm7
%\font\fivemsb=msbm5
%\newfam\msbfam
%\textfont\msbfam=\tenmsb
%\scriptfont\msbfam=\sevenmsb
%\scriptscriptfont\msbfam=\fivemsb
%\def\Bbb#1{{\fam\msbfam{#1}}}
\def\rr{\hbox{\tenmsb R}}
\def\cc{\hbox{\tenmsb C}}
\newdimen\length
\newdimen\lleftskip
\lleftskip=3.5\parindent%  to be adjusted !!!
\length=\hsize \advance\length-\lleftskip
\def\entry#1#2#3#4\par{\parshape=2  0pt  \hsize%
\lleftskip \length%
\noindent\hbox to \lleftskip%
{\bf[#1]\hfill}{\csc{#2 }}{\sl{#3}}#4%
\medskip% to be adjusted
}
\ifx \twelvebf\undefined \font\twelvebf=cmbx12\fi

\bigskip\bigskip
\centerline{\twelvebf References.}
\bigskip

%{ref}{author}{title}.

\entry{Ab}{Abikoff W.:}{The Real Analytic Theory of Teichm\"uller Space}.
Springer-Verlag (1980).

\entry{Al}{Alexander H.:}{Gromov's method and hulls.} Geometric Complex 
Analysis, Proc. Conf., ed. J. Noguchi et all (1966), World Scientific.

\entry{A-T-U}{Alexander H., Taylor ., Ullman .:}{Areas of projection of 
analytic sets.} Invent. math. {\bf16}, 335-341 (1973).

\entry{Ar}{Aronszajn, N.:}{A unique continuation theorem for elliptic
differential equation or inequalities of the second order.} J.\ Math.\ Pures
Appl., {\bf36}, 235--339, (1957).

\entry{B-K}{Bedford, E., Klingenberg W.:}{On the~envelope of meromorphy of
a 2-sphere in $\cc^2$.}  J.\ Amer.\ Math.\ Soc., {\bf4}, 623-646, (1991).

\entry{B-G}{Bedford, E., Gaveau B.:}{Envelopes of holomorphy of certain 
$2$-spheres in $\cc^2$.} Amer. J. Math. {\bf 105}, 975-1009 (1983).

\entry{B-P-V}{Barth W., Peters, C., Van de Ven, A.:}{Compact complex
surfaces.} Springer Verlag, (1984).

\entry{Bn}{Bennequin, D.:}{Entrelacement et \'equation de Pfaff.}
Ast\'erisque {\bf107--108}, 87--161 (1982).

%\entry{Ch-Sp}{Chern, S.-S., Spanier, E.:}{A theorem on orientable surfaces
%in four-dimensional space.} Comment.\ Math.\ Helv., {\bf25}, 205--209 (1951).

\entry{Ch}{Chirka E.:}{On generalized Hartogs lemma.} To appear in the
volume dedicated to B.\.V.\.Shabat, PHASIS, Moscow.

\entry{Ch-St}{Chirka, E., Stout L.:}{A Kontinuit\"atssatz.} In ``Topics in
Complex Analysis'', Banach Center Publications, {\b31}, 143--150, (1995).

\entry{D-G}{Dethloff G., Grauert H.:}{Deformation of compact Riemann surfaces
$Y$ of genus $p$ with distinguished points $P_1, \ldots P_m \in Y$.}
Int.\ Symp.\ ``Complex geometry and analysis'' in Pisa/Italy, 1988;
Lect.\ Notes in Math., {\bf1422}(1990), 37--44.

\entry{D-M}{Deligne P., Mumford D.:}{The irreducibility of the space of 
curves of a given genus.} IHES Math.\ Publ.,  {\bf36}(1969), 75--109.

\entry{E}{Eliashberg, Y.:}{Filling by holomorphic disks and its applications.}
London Math.\ Soc.\ Lecture Notes, {\bf151}, Geometry of low dimensional
manifolds, (1991).

\entry{F-S}{Fintushel R., Stern R.}{Immersed spheres in $4$-manifolds
and the immersed Thom conjecture}\ Turk. J. Math. v.19, N 2, 145--157 (1995).

\entry{F-M}{Fornaess, J.-E., Ma, D.:}{A $2$-sphere in $\cc^2$ that cannot be 
filled in with analytic disks.} Int. Math. Res. Notes {\bf1}, 17-22 (1995).

\entry{F}{Forstneri\v c, F.:}{Complex tangents of real surfaces in complex
surfaces.} Duke Math.\ J., {\bf67}, 353--376, (1992).

\entry{Fu}{Fujita, R.}{Domaines sans point critique int\'erieur sur l'espace
projectif complexe.} J. Math. Soc. Japan, {\bf15}, 443--473, (1963).

\entry{Ga}{Gahov, F.:}{Boundary values problems.} 3-rd edition, Nauka, 
Moscow (1977).

\entry{Gra}{Grauert, H.:}{\"Uber Modifikationen und exzeptionelle analytische
Mengen.} Math.Ann. {\bf146}, 331--368 (1962).

\entry{Gr-Hr}{Griffiths, P., Harris, J.:}{Principle of algebraic geometry.}
John Wiley \& Sons, N.-Y., (1978).

\entry{G}{Gromov, M.:}{Pseudo- holomorphic curves in symplectic
manifolds.} Invent.\ math. {\bf82}, 307--347 (1985).

\entry{Hv-Po}{Harvey, R., Polking, J.:}{Removable singularities of
solutions of partial differential equations.} Acta Math. {\bf125}, 39-56
(1970).

\entry{Hr-W}{Hartman, P., Winter, A.:}{On the~local behavior of solutions
of non-parabolic partial differential equations.} Amer.\ J.\ Math.,
{\bf75}, 449--476, (1953).

%\entry{Hf}{Hofer, H.:}{Pseudoholomorphic curves in symplectizations with
%applications to the Weinstein conjecture in dimension three.} Invent. Math.
%{\bf114}, 515--563 (1993).

\entry{H-L-S}{Hofer, H., Lizan, V., Sikorav, J.-C.:}{On genericity for
holomorphic curves in four-dimensional almost-complex manifolds.} J.\ of Geom.\
Anal., {\bf7}, 149--159, (1998).

\entry{Hu}{Hummel C.:}{Gromov's Compactness Theorem for Pseudoholomorphic 
Curves.} Birkh\"auser (1997).

\entry{Iv-1}{Ivashkovich, S.:}{Extension of analytic objects by the method
of Cartan-Thullen.} In proceedings of Conference on Complex Analysis and Math.
Phys., 53--61, Krasnojarsk, (1988).

\entry{Iv-2}{Ivashkovich, S.:}{The Hartogs-type extension theorem for
meromorphic maps into compact K\"ahler manifolds.} Invent. Math. {\bf109},
47--54 (1992).

\entry{IS-1}{Ivashkovich, S., Shevchishin, V.:}{Structue of the moduli space 
in the neighborhood of a cusp-curve and meromorphic hulls.} Invent. Math. 
{\bf136}, 571-602 (1999).

\entry{IS-2}{Ivashkovich, S., Shevchishin, V.:}{Deformations of non-compact 
complex curves and envelopes of meromorphy of spheres.} Russ. Math. Sbornik. 
{\bf189}, 1295-1333 (1998).

\entry{IS-3}{Ivashkovich S., Shevchishin V.:}{Gromov Compactness Theorem 
for Stable Curves.} math.DG/9903047 .

\entry{IS-4}{Ivashkovich S., Shevchishin V.:}{Pseudo-holomorphic curves and 
envelopes of meromorphy of two-spheres in~${\Bbb C}P^2$.} Preprint, Bochum 
 (1995), available as e-print math.CV/9804014.
 
\entry{Ko-No}{Kobayashi, S., Nomizu, K.:}{Foundations of differential
geometry.} Vol.II, Interscience Publishers, (1969).

\entry{K}{Kontsevich, M.:}{Enumeration of rational curves via torus actions.}
proc.\ Conf.\ ``The moduli spaces of curves'' on Texel Island, Netherland.
Birkh\"auser prog.\ Math., {\bf129}(1995), 335--368.

\entry{K-M}{Kontsevich, M., Manin Yu.:}{Gromov-Witten classes, quantum 
cohomology, and enumerative geometry.} Comm.\ Math.\ Phys., {\bf164}(1994),
525--562.

\entry{Kr}{Kronheimer, P.:}{Embedded surfaces and gauge theory in three and 
four dimensions.} Preprint, Harvard (1996).

\entry{Kr-M}{Kronheimer, P., Mrowka, T.:}{The genus of embedded surfaces in 
the projective plane.} Math. Res. Lett. {\bf1}, 797--808 (1994).

\entry{La}{Lai H.F.:}{Characteristic classes of real manifolds immersed in 
complex manifolds:} \ Trans. AMS, v.172, 1-33 (1972).

\entry{Le}{Lempert L.:}{Algebraic approximation in analytic geometry.}
\ Invent. Math., v.121, N 2, 335-354 (1995).

\entry{L-M}{Lisca P., Mati\'c G.:}{Tight contact structures and 
Seiberg--Witten invariants.} \ Invent. Math. v.129, N 3, 509-525 (1997).

%\entry{Li}{Lichnerovicz, A:}{Global theory of connections and holonomy
%groups.} Noordhoff Int. Publishing, Leiden, (1976).

\entry{M-W}{Micallef, M., White, B.:}{The structure of branch points in
minimal surfaces and in pseudoholomorphic curves.} Ann. Math. {\bf 139},
35-85 (1994).

\entry{McD-1}{McDuff, D.:}{The local behavior of holomorphic curves  in
almost complex 4-manifolds.} J.Diff.Geom. {\bf 34}, 143-164 (1991).

\entry{McD-2}{McDuff, D.:}{Examples of symplectic structures.} Invent.\ 
math., {\bf 89}, 13--36 (1987).

%\entry{McD-3}{McDuff, D.:}{Singularities and positivity of intersections of
%$J$ - holomorphic  curves.} In \it ``Holomorphic curves in symplectic geometry.''
% \rm Ed. by M.~Audin and J.~Lafontaine, Birkh\"auser, (1994).

\entry{Mg}{Morgan J.W.:}{The Seiberg--Witten Equations and Applications to the
Topology of Smooth Four-Manifolds}  Mathematical Notes 44, princeton 
University press, princeton (1996).

\entry{M-S-T}{ Morgan J.W., Szab\'o Z., Taubes C.H.:}{A product formula for 
Seiberg--Witten invariants and the generalized Thom conjecture} 
J. Differential Geom. l44, N4, 706-788 (1996).

\entry{Mo}{Morrey, C.:}{Multiple integrals in the calculus of variations.}
Springer Verlag, (1966).

\entry{Mr}{Moore, J.:}{Lectures on Seiberg--Witten invariants.} Lecture Notes 
in Math.1629, Springer, Berlin-New York (1996).

\entry{M}{Mumford D.:}{A remark on Mahler's compactness theorem.} proc.\
Amer.\ Math.\ Soc., 28, 289-294 (1971).

\entry{N-1}{Nemirovski S.:}{Stein domains on algebraic surfaces.} Russian
Math.\ Notes, {\bf 60}, 295-298 (1996).

\entry{N-2}{Nemirovski S.:}{Holomorphic functions and embedded real surfaces.} 
Math. Notes {\bf63}, 599-606 (1998).

\entry{N-3}{Nemirovski, S.:}{ On embeddings of the two-sphere in Stein 
surfaces.}\ Doklady RAN , v. 362, (1998).

\entry{N-4}{Nemirovski, S.:}{Complex Analysis and Differential Topology 
on Complex Surfaces.} Russ. Math. Surveys, {\bf 54}, N 4, 47-74 (1999).

\entry{Ra}{Ramis J.-P.:}{Sous-ensembles analytiques d'une vari\'et\'e
banachique complex.} Springer, Berlin (1970).

\entry{Rm}{Remmert, R.:}{Holomorphe und meromorphe Abbildungen komplexer
R\"aume.} Math. Ann. {\bf133}, 328--370 (1957).

\entry{Rf}{Rolfsen, D.:}{Knots and links.}Publish or perish, {\bf N7}, (1976).

\entry{R}{Rosay J.-P.:}{A Counterexample Related to Hartogs Phenomenon (A 
question by Chirka).} Michigan Math. Journal {\bf45}, 529-535 (1998).

\entry{Pa}{Parker T.:}{Bubble tree convergence for harmonic maps} J.\ Diff.\
Geom., {\bf44}(1996), 595--633.

\entry{P-W}{Parker T., Wolfson J.:}{Pseudoholomorphic maps and bubble trees}
J.\ Geom.\ Anal., {\bf3}(1993), 63--98.

\entry{S-U}{Sacks J., Uhlenbeck K.:}{Existence of minimal immersions of
two-spheres} Annal.\ Math., {\bf113}(1981), 1--24.

\entry{S}{Sikorav J.-C.:}{Some properties of holomorphic curves in almost
complex manifolds.} In "Holomorphic curves in symplectic geometry".
Edited by M. Audin, J. Lafontaine. Birkh\"auser, progress in
Mathematics v. 117, Ch.V, 165-189 .

\entry{Sk-2}{Sikorav, J.-C.:}{Singularities of $J$-holomorphic curves.} 
Math.\ Z., {\bf226}, 359--373, (1997).

\entry{Si-1}{Siu, Y.-T.:}{Extension of meromorphic maps into K\"ahler
manifolds.} Annals of Math. {\bf102}, 421--462 (1975).

\entry{Si-2}{Siu, Y.-T.:}{Every Stein subvariety admits a Stein
neighborhood.} Invent.math. {\bf38}, 89--100 (1976).

\entry{Sm}{Smale, S.:}{An infinite dimensional version of Sard's theorem.}
Amer. J. Math. {\bf87}, 861-866 (1965).

\entry{Sh}{Shafarevich I.:}{Basic algebraic geometry. Second edition.}
Springer-Verlag.

\entry{Sch}{Shcherbina, N.:}{On the~polynomial hull of a two-dimensional
sphere in $\cc^2$.} Soviet Math.\ Doklady, {\bf49}, 629--632, (1991).

\entry{St}{Stout L.:}{Algebraic domains in Stein manifolds.} Banach Algebras 
and several complex variables. Proc. Conf. New Haven/Conn. (1983), 
 Contemp. Math. (1984) {\bf32}, 259--266.

\entry{Tk}{Takeuchi, A.:}{Domaines pseudoconvexes infinis et la metrique 
riemannienne dans un espace projectif.} J. Math. Soc. Japan, {\bf16}, 
159--181 (1964).

\entry{Ta}{Taubes, C.:}{$SW \Longrightarrow Gr$: From the Seiberg-Witten 
Equations to Pseudo-Holomorphic Curves.} J. Amer. Math. Soc. {\bf9}, n 3,
845-918 (1996).

%\entry{Wn}{Weinstein, A.:}{Lectures on symplectic manifolds.} CBMS. Reg.
%Conf. Series, {\bf N29}, AMS(1977).

\entry{Wa}{Walker, R.:}{Algebraic curves.} Springer Verlag, (1978).

\entry{Wt}{Witten, E.:}{Monopoles and four-manifolds.} Math. Res. Lett.
{\bf1}, 769--796 (1994).

%\endRefs

\bigskip
\bigskip
\settabs 2 \columns
\+Institute of Applied Problems
&Institute of Applied problems\cr
\+of Mechanics and Mathematics
&of Mechanics and Mathematics\cr
\+Ukrainian Acad. Sci.,
&Ukrainian Acad. Sci.,\cr
\+vul. Naukova 3b, 290053 L'viv
&vul. Naukova 3b, 290053 L'viv \cr
\+Ukraine
&Ukraine \cr
\medskip
\+ U.F.R. de Maht\'ematiques
& Ruhr-Universit\"at Bochum \cr
\+ Universit\'e de Lille-I
& Mathematisches Institut \cr

\+ Villeneuve d'Ascq Cedex
& Universit\"atsstrasse 150 \cr
\+ 59655 France
& NA 4/67  44780 Bochum  Germany \cr
\+ ivachkov\@gat.univ-lille1.fr
& sewa\@cplx.ruhr-uni-bochum.de   \cr

\enddocument
\bye